\documentclass{amsart}
\usepackage[backref]{hyperref}
\usepackage{euscript}
\usepackage{amsmath}
\usepackage{amscd}
\usepackage{amssymb}
\usepackage{enumerate}
\usepackage{amsfonts}
\usepackage{accents}
\usepackage{graphicx}
\usepackage{subcaption}
\usepackage[all]{xy}
\usepackage{mathtools}
\usepackage{setspace}
\usepackage{array}
\usepackage{bbold}
\usepackage{tikz}
\usetikzlibrary{matrix,arrows,calc,intersections,fit,shapes.geometric,shapes.misc}
\usepackage{tikz-cd}

\newsavebox{\forcedcentertikzcdbox}
\usepackage{textcomp}
\makeatletter
\DeclareFontEncoding{LS1}{}{}
\DeclareFontSubstitution{LS1}{stix}{m}{n}
\DeclareMathAlphabet{\mathstscr}{LS1}{stixscr}{m}{n}
\makeatother
\makeatletter
\renewcommand{\maketag@@@}[1]{\hbox{\m@th\normalsize\normalfont#1}}
\makeatother
\usepackage{hyperref}

\DeclareMathOperator{\cMap}{\mathstscr{M\mkern-4mu a\mkern-3.5mu p}}

\newcommand{\ssF}{\mathstscr{F}}
\mathchardef\nhyphen=45

\newcommand{\sma}{\wedge}

\newcommand{\bC}{{\mathbb C}}

\newcommand{\bP}{{\mathbb P}}
\newcommand{\bQ}{{\mathbb Q}}
\newcommand{\bR}{{\mathbb R}}
\newcommand{\bM}{{\mathbb M}}
\newcommand{\bN}{{\mathbb N}}
\newcommand{\bS}{{\mathbb S}}
\newcommand{\bF}{{\mathbb F}}

\newcommand{\bX}{{\mathbb X}}

\newcommand{\bZ}{{\mathbb Z}}

\newcommand{\sI}{\EuScript I}
\newcommand{\sJ}{\EuScript J}
\newcommand{\sV}{\EuScript V}
\newcommand{\sZ}{\EuScript Z}

\newcommand{\bbn}{{\mathbb n}}
\newcommand{\bbm}{{\mathbb m}}
\newcommand{\bbk}{{\mathbb k}}

\newcommand{\bbone}{{\mathbb{1}}}
\newcommand{\bDelta}{{\mathbb \Delta}}

\newcommand{\fC}{{\mathfrak C}}

\newcommand{\cF}{\mathcal F}

\newcommand{\cG}{\mathcal G}

\newcommand{\cA}{\mathcal A}
\newcommand{\aA}{\mathcal A}
\newcommand{\cB}{\mathcal B}
\newcommand{\aB}{\mathcal B}
\newcommand{\cC}{\mathcal C}
\newcommand{\aC}{\mathcal C}
\newcommand{\aM}{\mathcal M}
\newcommand{\cD}{\mathcal D}
\newcommand{\aD}{\mathcal D}
\newcommand{\cE}{\mathcal E}

\newcommand{\cI}{\mathcal I}

\newcommand{\cL}{\mathcal L}
\newcommand{\cM}{\mathcal M}

\newcommand{\fcMR}{\Mbar^{\bR}}
\newcommand{\fcMRreg}{\Mbar^{\bR,\reg}}
\newcommand{\Cbar}{\overline{\cC}}
\newcommand{\fcCR}{\Cbar^{\bR}}

\newcommand{\aN}{\mathcal N}
\newcommand{\cO}{\mathcal O}
\newcommand{\cP}{\mathcal P} 
\newcommand{\cQ}{\mathcal Q} 
\newcommand{\cR}{\mathcal R}
\newcommand{\cS}{\mathcal S}
\newcommand{\cT}{\mathcal T}
\newcommand{\cU}{\mathcal U}
\newcommand{\cV}{\mathcal V}
\newcommand{\cX}{\mathcal X}

\newcommand{\cZ}{\mathcal Z}
\newcommand{\barT}{\overline{\cT}}

\newcommand{\scrS}{\EuScript S}

\newcommand{\scrE}{\EuScript E}
\newcommand{\scrF}{\EuScript F}
\newcommand{\scrT}{\EuScript T}
\newcommand{\scrI}{\EuScript I}

\newcommand{\scrP}{\EuScript P}

\newcommand{\Q}{Q}

\newcommand{\Mbar}{\overline{\cM}}

\newcommand{\Rbar}{\overline{\cR}}

\newcommand{\Tbar}{\overline{\cT}}

\newcommand{\Milnor}{\overline{M}}
\newcommand{\Nash}{\overline{N}}

\newcommand{\rel}{\mathrm{rel}}

\newcommand{\htp}{\simeq}

\newcommand{\op}{\mathrm{op}}

\newcommand{\fib}{\mathrm{fib}}
\newcommand{\base}{\mathrm{base}}
\newcommand{\mfib}{\mathrm{mfib}}
\newcommand{\Ho}{\mathrm{Ho}}

\newcommand{\Filt}{\mathrm{Filt}}

\newcommand{\C}{{\mathrm C}}

\newcommand{\scrJ}{{\EuScript J}}

\newcommand{\id}{\operatorname{id}}

\newcommand{\virdim}{\operatorname{virdim}}
\renewcommand{\mod}{\operatorname{mod}}

\newcommand{\Tw}{\operatorname{Tw}}

\newcommand{\Comp}{\operatorname{Comp}}
\newcommand{\Hom}{\operatorname{Hom}}
\newcommand{\End}{\operatorname{End}}
\newcommand{\Ends}{\operatorname{Ends}}

\newcommand{\Aut}{\operatorname{Aut}}

\newcommand{\Map}{\operatorname{Map}}

\newcommand{\Sp}{\operatorname{Sp}}
\newcommand{\Spcat}{\operatorname{SpCat}}

\newcommand{\dbar}{\overline{\partial}}

\newcommand{\Gr}{\operatorname{Gr}}

\newcommand{\Sub}{\operatorname{Sub}}

\newcommand{\VB}{\operatorname{VB}}
\newcommand{\Vect}{\operatorname{Vect}}

\newcommand{\Diff}{\operatorname{Diff}}
\newcommand{\Id}{\operatorname{Id}}

\renewcommand{\max}{\operatorname{max}}

\def\co{\colon\thinspace}

\newcommand{\hocofib}{\operatorname{hocofib}}

\newcommand{\Coeq}{\operatorname{Coeq}}

\newcommand{\Ev}{\mathrm{Ev}}

\newcommand{\gr}{\mathrm{gr}}
\newcommand{\Fun}{\mathrm{Fun}}

\newcommand{\ob}{\operatorname{ob}}
\newcommand{\mor}{\operatorname{mor}}
\newcommand{\eq}{\mathrm{eq}}

\newcommand{\nd}{\mathrm{nd}}
\newcommand{\hofib}{\mathrm{hofib}}
\newcommand{\Ass}{\mathrm{Ass}}
\newcommand{\Comm}{\mathrm{Comm}}
\newcommand{\disc}{\mathrm{disc}}
\newcommand{\kmod}{\bbk\text{-}\mbox{mod}}
\newcommand{\Span}{\mathrm{Span}}
\newcommand{\surj}{\mathrm{surj}}

\newlength{\dhatheight}
\newcommand{\dhat}[1]{
      \mathbf{\hat{\text{$#1$}}}}
    \newcommand{\that}[1]{
    \settoheight{\dhatheight}{\ensuremath{\hat{#1}}}
    \addtolength{\dhatheight}{-0.35ex}
    \hat{\vphantom{\rule{1pt}{\dhatheight}}
      \smash{\dhat{#1}}}}

\newcommand{\dhatS}[1]{
    \settoheight{\dhatheight}{\ensuremath{\scriptstyle{\hat{#1}}}}
    \addtolength{\dhatheight}{-0.175ex}
    \hat{\vphantom{\rule{1pt}{\dhatheight}}
      \smash{\hat{#1}}}}

  \newcommand{\dhatC}[1]{
    \settoheight{\dhatheight}{\ensuremath{\scriptstyle{\hat{#1}}}}
    \addtolength{\dhatheight}{-0.3ex}
    \hat{\vphantom{\rule{1pt}{\dhatheight}}
      \smash{\hat{#1}}}}

  \newcommand{\thatS}[1]{
    \settoheight{\dhatheight}{\ensuremath{\hat{#1}}}
    \addtolength{\dhatheight}{-0.175ex}
       \hat{\vphantom{\rule{1pt}{1.8ex}}
      \smash{\hat{\vphantom{\rule{1pt}{\dhatheight}}
          \smash{\hat{#1}}}}}}
  
   \newcommand{\thatB}[1]{
    \settoheight{\dhatheight}{\ensuremath{\hat{#1}}}
    \addtolength{\dhatheight}{-0.175ex}
       \hat{\vphantom{\rule{1pt}{2.125ex}}
      \smash{\hat{\vphantom{\rule{1pt}{\dhatheight}}
      \smash{\hat{#1}}}}}}

\numberwithin{equation}{subsection}

\newtheorem{thm}{Theorem}[section]
\newtheorem{cor}[thm]{Corollary}
\newtheorem{lem}[thm]{Lemma}
\newtheorem{lemma}[thm]{Lemma}
\newtheorem{prop}[thm]{Proposition}
\newtheorem{defin}[thm]{Definition}
\newtheorem{def-lem}[thm]{Definition-Lemma}

\theoremstyle{remark}

\newtheorem{rem}[thm]{Remark}
\newtheorem{warning}[thm]{Warning}

\newtheorem{example}[thm]{Example}

\newtheorem{notn}[thm]{Notation}

\newcommand{\superscript}[1]{\ensuremath{^{\textrm{#1}}} }

\renewcommand{\nd}[0]{\superscript{nd}}

\newcommand{\comment}[1]{}
\newcommand{\Top}{\operatorname{Top}}
\newcommand{\GTop}{\operatorname{GTop}}
\newcommand{\Cat}{\mathrm{Cat}}

\newcommand{\Orb}{\operatorname{Chart}_{\cO}}
\newcommand{\wtOrb}{\widetilde{\operatorname{Chart}}_{\cO}}
\newcommand{\Orbo}{\Orb^{\emptyset}}

\newcommand{\reg}{\mathrm{reg}}
\newcommand{\ori}{\mathrm{ori}}

\newcommand{\TOP}{\mathrm{TOP}}

\newcommand{\Ar}{\operatorname{Ar}}

\newcommand{\Kur}{\operatorname{Chart}_\mathcal{K}}

\newcommand{\iso}{\mathrm{iso}}

\DeclareMathOperator*{\hocolim}{hocolim}
\DeclareMathOperator*{\colim}{colim}
\DeclareMathOperator*{\holim}{holim}

\newcommand{\sk}{\mathrm{sk}}
\newcommand{\Set}{\operatorname{Set}}
\newcommand{\Pre}{\operatorname{Pre}}
\newcommand{\Sh}{\operatorname{Sh}}
\newcommand{\Coh}{\operatorname{Coh}}

\title[Arnold Conjecture and Morava $K$-theory]{Arnold Conjecture and Morava $K$-theory}

\author[M.~Abouzaid]{Mohammed Abouzaid}
\author[A.J.~Blumberg]{Andrew J. Blumberg}
\date{\today} 
\thanks{The first author was supported by NSF grants DMS-1308179,  DMS-1609148, and DMS-1564172, as well as the Simons Foundation through its ``Homological Mirror Symmetry'' Collaboration grant.  The second author was supported in part by NSF grants DMS-1812064 and DMS-1564289.}
\begin{document}
\begin{abstract}
  We prove that the rank of the cohomology of a closed symplectic manifold with coefficients in a field of characteristic $p$ is smaller than the number of periodic orbits of any non-degenerate Hamiltonian flow. Following Floer, the proof relies on constructing a homology group associated to each such flow, and comparing it with the homology of the ambient symplectic manifold. The proof does not proceed by constructing a version of Floer's complex with characteristic $p$ coefficients, but uses instead the canonical (stable) complex orientations of moduli spaces of Floer trajectories to construct a version of Floer homology with coefficients in Morava's $K$-theories, and can thus be seen as an implementation of Cohen, Jones, and Segal's vision for a Floer homotopy theory. The key feature of Morava K-theory that allows the construction to be carried out is the fact that the corresponding homology and cohomology groups of classifying spaces of finite groups satisfy Poincar\'e duality.
\end{abstract}
\maketitle
\tableofcontents

\section{Introduction}
\label{sec:introduction}

This introduction is split into two parts: the first explains the consequences of this paper for symplectic topology. The second gives a proof of these applications, providing references to the various parts of the paper where the necessary ingredients are established. 

\subsection{Applications of the main results}
\label{sec:appl-main-results}
Let $(M,\omega)$ be a closed symplectic manifold, and $H \co M \times S^1 \to \bR$ a Hamiltonian function all of whose closed time-$1$ Hamiltonian orbits are non-degenerate (these data will be fixed for the entirety of the paper). Arnol'd conjectured that the number of such orbits is bounded below by the minimal number of critical points of a Morse function on $M$. We prove: 
\begin{thm}\label{thm:bound-orbits-Morava}
  For each natural number $0 \leq n < \infty$, the rank of the 
  generalized homology of $M$ with respect to Morava $K$-theory, $H_*(M;K(n))$, is smaller than the number of time-$1$ closed contractible Hamiltonian orbits of $H$. 
\end{thm}
In the above statement, we operate with the usual conventions in chromatic homotopy theory (we recommend~\cite{Ravenel1992} as a good introductory reference for the chromatic viewpoint and~\cite{Wurgler1991} as a guide to the Morava $K$-theories in particular): the Morava $K$-theories $K(n)$ implicitly depend on a prime $p$ which is not included in the notation, and we adopt the standard convention that the case $n=0$ corresponds to rational homology $H_*(M, \bQ)$.  The above result thus contains as a special case the result of Fukaya and Ono \cite{FukayaOno1999} and Liu and Tian \cite{LiuTian1998}, who extended Floer's construction of Hamiltonian Floer homology \cite{Floer1989b} to arbitrary closed symplectic manifolds.

Our notation for generalized homology groups will however be non-standard, as we write $H_*(M; K(n))$ for the group usually written as $K_*(n)(M)$.  For $n$ strictly positive, the coefficients $H_*(\ast; K(n)) = K(n)_*$ of Morava $K$-theory (i.e., the generalized homology groups of a point) are 
\begin{equation}
  K(n)_* \cong \bF_{p}[v_n^\pm],
\end{equation}
where $v_n$ is a variable of degree of $2(p^n-1)$.  In Theorem~\ref{thm:bound-orbits-Morava}, the rank of $H_*(M;K(n))$ is taken as a module over $\bF_{p}[v_n^\pm]$; since this ring is a graded field and $H_*(M;K(n))$ is a graded module, the rank in Theorem \ref{thm:bound-orbits-Morava} is simply the number of elements of a basis of
$H_*(M;K(n))$ over $K(n)_*$.

The bound of Theorem \ref{thm:bound-orbits-Morava} can be re-expressed in terms of ordinary homology using the Atiyah-Hirzebruch spectral sequence: if we choose $n$ so that $2(p^n-1)$ is larger than the dimension of $M$, the spectral sequence which computes $H_*(M;K(n))$ from 
\begin{equation}
H_*(M,K(n)_*) \cong H_*(M;\bF_p) \otimes K(n)_*
\end{equation}
collapses at the $E_2$ page and so $H_*(M;K(n))$ has the same rank as $H_*(M;\bF_p) \otimes_{\bF_p} K(n)_*$.  That is, all Morava $K$-theories for such $n$ have the same rank which moreover agrees with the rank of $H_*(M;\bF_p)$ over $\bF_p$.  Passing from $\bF_p$ to an arbitrary characteristic $p$ field by the universal coefficients formula, we conclude the statement given in the abstract:
\begin{cor} \label{cor:arnold_f_p}
 The rank of the 
  ordinary homology group $H_*(X;\Bbbk)$ for each characteristic $p$ field $\Bbbk$ is smaller than the number of time-$1$ closed contractible  Hamiltonian orbits of $H$. \qed
\end{cor}

\begin{rem}
In~\cite{FukayaOno2001}, Fukaya and Ono proposed an approach to an integral version of Corollary~\ref{cor:arnold_f_p}. The results of this paper shares with their proposal the use of the natural stable almost complex structure on moduli spaces of Hamiltonian Floer trajectories. In this paper, we combine this stable almost complex structure with orientation theory for generalized cohomology theories, whereas their proposal is to use this stable almost complex structure to choose (virtual) perturbations of the Cauchy-Riemann equation for which all perturbed solutions lie in the locus with trivial isotropy.
\end{rem}

\subsubsection{Hamiltonian Floer $K(n)$-homology}
\label{sec:hamiltonian-floer-kn}

The proof of Theorem \ref{thm:bound-orbits-Morava} follows in broad strokes the previous work on the Arnol'd conjecture initiated by Floer. We begin by considering the lattice $\Pi$ of classes $\beta \in H_2(M; \bZ)$ lying in the image of the map from $\pi_1(\cL M)$, the fundamental group of the free loop space $\cL M$, to the homology of $M$ which associates to each loop of loops the corresponding torus. Recalling that the symplectic form $\omega$ defines an energy map
\begin{equation}
  H_2(M; \bZ) \to \bR
\end{equation}
which assigns to each class $\beta$ the area $\omega(\beta)$ of any representing curve, we define the Novikov ring $\Lambda_*$ with coefficients in $K(n)_*$ to be the $2$-periodic ring whose elements of degree $0$ are infinite sums
\begin{equation} \label{eq:Novikov-ring}
\sum_{\substack{ \beta \in \Pi  \\ \omega(\beta) \to +\infty}} a_{\beta} q^{\beta}
\end{equation}
with $a_\beta \in K(n)_*$.

Here, we use the familiar notation in symplectic topology where such a series consists of an a priori infinite sum indexed by elements of $\Pi$, with the condition that there are only finitely many terms indexed by elements whose energy is smaller than any given integer. This ring can be thought of more abstractly as follows: if $u$ is a variable of degree $2$, we obtain $\Lambda_*$ by completing the group ring of $\Pi$ with coefficients in the $2$-periodic ring $K(n)_*[u,u^{-1}]$, with respect to the filtration associated to the map $\Pi \to \bR$ induced by $\omega$.

We assign to $M$ the $\Lambda_*$ module 
\begin{equation}
H_*(M;\Lambda) \equiv H_*(M;K(n))\otimes_{K(n)_*} \Lambda_*. 
\end{equation}
\begin{rem}
The above formula will end up being an (easy) theorem, rather than a definition. As we shall presently see in Section \ref{sec:floer-kn-homotopy}, there is a ring spectrum $\Lambda$ whose homotopy groups are $\Lambda_*$ and which is a $K(n)$-module.  The left hand side should be thought of as the generalized homology of $M$ with coefficients in $\Lambda$.  Because $K(n)$ is a field, this homology admits the above expression (see e.g.~\cite[3.4]{HoveyStrickland1999}).
\end{rem}

The key goal of this paper is the construction of the \emph{Floer homology of $H$ with coefficients in $K(n)$}, denoted $HF_*(H; \Lambda)$,  as  a $\Lambda_*$-module. For the next result, we shall consider the rank of this module, which we define to be the rank of the associated module over the field of fractions of $\Lambda_*$.  The main application follow immediately from the following result:
\begin{prop} \label{prop:properties_Ham_K_n}
  The Floer homology $HF_*(H; \Lambda)$ satisfies the following properties as a $\Lambda_*$ module:
  \begin{enumerate}
  \item it admits $H_*(M;\Lambda) $ as a summand, and
    \item its rank is bounded above by the number of contractible closed time-$1$ Hamiltonian orbits of $H$.
  \end{enumerate}
\end{prop}

\begin{proof}[Proof of Theorem \ref{thm:bound-orbits-Morava}]
The rank of $HF_*(H; \Lambda)$ is bounded below by the rank of $H_*(M;\Lambda)$ over $\Lambda_*$, which agrees with the rank of $H_*(M;K(n))$ over $K_*(n)$ by base change. We conclude that the number of closed time-$1$ Hamiltonian orbits is bounded as desired.
\end{proof}

\begin{rem}
Though we shall not establish the result in this paper, a slight extension of our constructions implies that the isomorphism type of $HF_*(H; \Lambda) $ as a $\Lambda_*$ module is independent of the Hamiltonian function, and of the auxiliary data required to define it, e.g. of a choice of compatible almost complex structure on $M$. Choosing the Hamiltonian to be $C^2$-small then implies (via a spectral sequence argument), that $HF_*(H; \Lambda)$ is isomorphic to $H_*(M;\Lambda)$. As in the case of ordinary Floer homology, we expect that there will be many alternative ways of proving this result, including using $S^1$-localization as in Fukaya and Ono's original work \cite{FukayaOno1999}, or further pushing the comparison between Morse and Floer theory \cite{Fukaya1997b,PiunikhinSalamonSchwarz1996,FilippenkoWehrheim2018} adopted in this paper.
\end{rem}

\begin{rem}
 As we shall presently discuss, the construction of this paper can be seen as a realization of the vision of Cohen, Jones, and Segal \cite{CohenJonesSegal1995}, that the Floer homology of symplectic manifolds should lift to a (stable) homotopy type, in analogy with the Morse homology of a finite dimensional manifold arising by applying the homology functor to its stable homotopy type. Such a homotopy type would give rise to Floer homology groups with coefficients in generalized cohomology theories. However, it was realized early on that, because the tangent spaces of moduli spaces arising in Floer theory are not in general (stably) trivial, but always have stable almost complex structures, one can only hope for the existence of such Floer homology groups with coefficients in generalized cohomology theories equipped with Thom isomorphism theorems for complex vector bundles, i.e. for complex-oriented cohomology theories. The key point is that such theories satisfy Poincar\'e duality isomorphisms for manifolds equipped with stable almost complex structures. In particular, the Floer (stable) homotopy groups of a Hamiltonian on a general symplectic manifold cannot be defined.

In Section \ref{sec:virt-fund-chains} below, we explain why we have to restrict our coefficients further, and do not define Floer homology with coefficients in arbitrary  complex-oriented cohomology theories. The key problem is that the moduli spaces appearing in Floer theory can sometimes be (locally) described as quotients of manifolds by finite groups actions with non-trivial fixed points, so the construction of generalized Floer homology groups requires considering theories which are equipped with Poincar\'e duality isomorphisms for orbifolds (equipped with stable almost complex structures). The simplest known examples of such theories arise from the Morava K-theories discussed above.
\end{rem}

\subsubsection{The Floer $K(n)$-homotopy type}
\label{sec:floer-kn-homotopy}

We now explain the construction of Hamiltonian Floer $K(n)$-homology, and the proof of Proposition \ref{prop:properties_Ham_K_n}. In Floer theory, the construction of Hamiltonian Floer homology classically proceeds via the construction of a chain complex (the Floer complex). Unfortunately, the only homology theories that are naturally computed by chain complexes agree with ordinary homology \cite{BurdickConnerFloyd1968}.  We must therefore start using the language of stable homotopy theory, and introduce various spectra from which (generalized) homology groups are obtained by considering homotopy groups.  The appearance of homotopy groups becomes less surprising when one is reminded that the homology of a chain complex is the quotient of cycles by boundaries, that cycles correspond to maps from the ground ring considered as free chain complex of rank $1$, and that boundaries correspond to maps which are chain homotopic to zero, so that the quotient corresponds to passing to the equivalence relation of chain homotopy.

\begin{rem}
Unfortunately, there is not a definitive modern exposition of spectra
and the stable category aimed at a general audience.  Adams~\cite[Part
  3]{Adams1974} is the classic introduction to spectra for the reader
unfamiliar with the subject;
Lewis-May-Steinberger~\cite{LewisMaySteinberger1986} is the definitive
treatment of the category of spectra which was the basis for most work
in the subject in the 20th century.  A very nice discussion and
comparison between the modern theories of diagram spectra (which gives a
point-set symmetric monoidal model of the stable category) is
Mandell-May-Schwede-Shipley~\cite{MandellMaySchwedeShipley2001}, and
the equivariant theory is described in Mandell-May~\cite{MandellMay}.

Although in the body of the paper we will require the full strength of
the modern theories, for the purpose of this introduction the reader
should simply have in mind that spectra are collections of based
spaces $\{X_i\}$ indexed by the natural numbers, together with the
datum of maps $\Sigma X_i \to X_{i+1}$ from the suspension of each
space to the next, together with additional structures whose
specification would not help the exposition.  Given a based space $A$,
the associated {\em suspension spectrum} is the collection of spaces
$\{\Sigma^i A\}$, where $\Sigma^i$ denotes the iterated suspension.
Roughly speaking, the category of spectra is a model of the result of
formally inverting the suspension operator $\Sigma$ on the category of
spaces; the inverse of $\Sigma$ (up to homotopy) is then the loops
$\Omega$, which we sometimes write $\Sigma^{-1}$ depending on
context.  The homotopy groups of a spectrum are produced by a process
of stabilization of the homotopy groups of the spaces $\{X_i\}$.

A key point to take from the above vague description is
that operations on based spaces induce operations on spectra.
  At the level of spaces, the operations that we shall use are the
  addition of a disjoint basepoint, which canonically associates to
  each space $X$ a based space $X_+$, the smash product of based
  spaces $(A,\ast)$ and $(B,\ast)$: 
  \begin{equation}
    A \wedge B \equiv A \times B/ A \times \ast \cup \ast \times B,
  \end{equation}
and the construction of the mapping space $F(A,B)$ of (base-point preserving) continuous maps
between based spaces.
  
In order to allow the reader unfamiliar with the subject to understand
the basic ideas, we shall use notation reminiscent of chain
complexes. In particular, we write
\begin{equation}
  C_*(M;\bbk) \equiv M_+ \wedge \bbk \cong  \Sigma^{\infty}_+ M \sma \bbk,
\end{equation}
for the \emph{spectrum of chains} of $M$ with coefficients in a
spectrum $\bbk$. This is the spectrum whose homotopy groups are the
generalized homology groups $H_*(M;\bbk)$, which are usually written
$\bbk_*(M)$. We shall treat the subscript in $C_*$ as entirely
decorative (in the sense that $C_k(M;\bbk)$ is not given any meaning),
though it does allow us to distinguish chains from cochains: we thus write 
\begin{equation}
C^*(M;\bbk) \equiv F(M_+, \bbk)
\end{equation}
for the {\em spectrum of cochains} of $M$ with coefficients in the
spectrum $\bbk$; the homotopy groups of this spectrum are the
generalized cohomology groups $H^*(M;\bbk)$, usually denoted $\bbk^*(M)$.

Our work requires the use of the algebraic structure of the category
of spectra; this is a symmetric monoidal category under the smash
product, for which the sphere spectrum is the unit.  This allows the
definition of ring spectra and categories of module spectra associated
to a ring spectrum in the evident ways.  That is, a ring spectrum
$\bbk$ is a spectrum with a unit map $\bS \to \bbk$ and a
multiplication $\bbk \sma \bbk \to \bbk$ which is appropriately
associative and unital.  Such a ring has a category of \emph{module
  spectra} which we denote $\bbk\nhyphen\mod$, consisting of spectra equipped
with an associative action of $\bbk$; one key example are the chains
with coefficients in $\bbk$, which are equipped with a natural map
\begin{equation}
  C_*(M;\bbk) \wedge \bbk \to C_*(M;\bbk).
\end{equation}

Classical models of spectra (i.e., those described in~\cite{Adams1974} and~\cite{LewisMaySteinberger1986})  only admitted a symmetric monoidal
structure after passage to the homotopy category (for the derived
smash product); the advantage of more modern categories of spectra (e.g. as discussed in \cite{MandellMaySchwedeShipley2001}) is
a point-set symmetric monoidal structure.  Module spectra on the
point-set level are a much more satisfactory technical notion than
modules in the homotopy category; we will be careful to
distinguish in what follows between the two cases.  Nonetheless, the
casual reader does not need to worry extensively about these issues on
a first reading.

We end this aside by warning the reader that the proper development of
the theory requires understanding the modern homotopy theory of ring
spectra and their modules and the methods for computing derived
functors in this context.  For example, we implement much of this by
considering cofibrant ring spectra, and notions of fibrant and
cofibrant modules.  These technical issues will be elided in this
introduction, and are discussed in Appendix
\ref{sec:groups-categ-spectra}.
\end{rem}

Returning to the discussion of the previous section, we fix a prime $p$ and consider a associative ring spectrum $K(n)$ representing Morava $K$-theory; the ring structure was shown to be unique by Angeltveit~\cite{Angeltveit2011}, although there are uncountably many choices of complex orientation (see Appendix \ref{sec:ambid-orient}).  Associated to any monoid $\Pi$, and in particular to the submonoid of $H_2(M,\bZ)$ considered earlier, we have a \emph{monoid ring with coefficients in $K(n)$,} which we denote
\begin{equation}
  K(n)[\Pi] \equiv   K(n) \wedge \Sigma^{\infty}_+ \Pi.
\end{equation}
This is an associative ring spectrum; the product in the monoid ring is defined by combining the operation on $\Pi$ and the product in $K(n)$.  Moreover, the natural inclusion $K(n) \to K(n)[\Pi]$ makes the monoid ring into an algebra under $K(n)$.  In our context, there is additional structure coming from the homomorphism $\omega \co \Pi \to \bR$.  Specifically, $\omega$ induces a grading of this ring by $\bR$; it will be more convenient to pass to the associated decreasing filtration $K(n)[\Pi_{\geq \lambda}] $  indexed by $\lambda \in \bR$, arising from the filtration of $\Pi$ which associates to each real number $\lambda$ all elements whose image under $\omega$ is larger than $\lambda$.

We are interested in a $2$-periodic version of this ring, so we
introduce an associative ring spectrum $P\bS$ which is equivalent to
the infinite wedge $\bigvee_n \Sigma^{-2n} \bS$ of negative and
positive spheres of even dimension, and which we refer to as the
$2$-periodic sphere spectrum. For our purposes it is convenient to have more point-set control, and so
we describe a specific model of this ring in Appendix~\ref{sec:action}.

We define the
$2$-periodic group ring of $\Pi$ with coefficients in $K(n)$ to be the
smash product $ P \bS \wedge K(n)[\Pi] $, which inherits a natural
filtration from $K(n)[\Pi]$, with quotients
\begin{equation}
P\bS \wedge   K(n)[\Pi_{\leq \lambda}] \equiv   P\bS \wedge K(n)[\Pi]/ P\bS \wedge K(n)[\Pi_{\geq \lambda}].
\end{equation}
We define the \emph{spectral Novikov
  ring} $\Lambda$ to be the completion of the group ring with respect
to this filtration, i.e., the corresponding homotopy inverse limit.
\begin{equation}
  \Lambda \equiv \holim_{\lambda} P\bS \wedge   K(n)[\Pi_{\leq \lambda}].
\end{equation}

By construction, the coefficient ring of $\Lambda$ is the graded ring
$\Lambda_*$ from Section~\ref{sec:hamiltonian-floer-kn}, and the
natural map
\begin{equation}
  P \bS \wedge  K(n)[\Pi] \to \Lambda
\end{equation}
induces an equivalence on associated graded spectra with respect to the subfiltration associated to any discrete unbounded subset of $\bR$ (note that all such subsets are countable). 

The main construction of this paper, performed in Section \ref{sec:constr-spectra-floer} below, is of a \emph{spectrum of Floer chains} 
\begin{equation}
 CF_*(H; \Lambda),
\end{equation}
associated to the Hamiltonian $H$, which is a module over $\Lambda$
(in the homotopy category) and is equipped with a complete decreasing
filtration. We denote by 
\begin{equation}
   CF^\lambda_*(H; \Lambda)
 \end{equation}
 the quotient of $CF_*(H; \Lambda)$ by its $\lambda$-filtered part; $CF_*(H;\Lambda)$ is then the homotopy inverse limit $\holim_{\lambda} CF^\lambda_*(H; \Lambda)$.

\begin{defin}\label{def:definition_Floer-homology-quantitative}
  The Floer homology groups $HF_*(H; \Lambda)$ are given by the inverse limit
  \begin{equation}
    HF_k(H; \Lambda) \equiv \lim_{\lambda} \pi_k(  CF^\lambda_*(H; \Lambda) )
  \end{equation}
of the homotopy groups of $ CF^\lambda_*(H; \Lambda)$:
\end{defin}

\begin{rem} \label{rem:we-are-lazy}
A more sensible version of
Definition~\ref{def:definition_Floer-homology-quantitative} would be
to set
  \begin{equation}
     HF_k(H; \Lambda) \equiv \pi_k( CF_*(H; \Lambda) ).
  \end{equation}
However, the ad hoc definition above has the advantage that it allows
us to prove Theorem~\ref{thm:bound-orbits-Morava} without having to
analyze the Mittag-Leffler condition for the inverse system
$\pi_k(CF^\lambda_*(H; \Lambda) ) $ of homotopy groups.  In a more
complete account of Hamiltonian Floer homology with coefficients in
Morava $K$-theories, one would verify this condition (which in fact
does hold).
\end{rem}

We can now explain how the results of the previous section follow from results stated in terms of spectra, which we split into two parts:
\begin{prop} \label{prop:quantitative-Floer-chains-and-residue}
  The filtered spectrum $ CF_*(H; \Lambda)$ admits the structure of a filtered module over the filtered ring $ \Lambda$ in the homotopy category of spectra.
  The associated graded spectrum with respect to a sufficiently fine discrete and unbounded subset of $\bR$ is a free module over $ P\bS \wedge K(n)[\Pi]$, generated by closed time-$1$ periodic orbits of $H$.
\end{prop}

\begin{prop}  \label{prop:retract_CF_ordinary-chains}
The Floer cochains $CF_*(H; \Lambda)$ admit $C_*(M; \Lambda)$ as a
retract in the homotopy category of modules over $\Lambda$.  Moreover,
this retraction is compatible with the filtration by $\Lambda$ in the
following sense: if we regard $C_*(M; \Lambda)$ as a filtered
$\Lambda$-module using the isomorphism $C_*(M; \Lambda) \cong
C_*(M;\bbk) \sma_{\bbk} \Lambda$, then there exists a constant $c$
such that the retraction restricts to a composite
\begin{equation}
C^\lambda_*(M; \Lambda) \to CF^{\lambda+c}_*(H; \Lambda) \to
C^{\lambda+c}_*(M; \Lambda),
\end{equation}
where $c$ does not depend on $\lambda$.
\end{prop}

Both of these results depend on the module structure on the Floer
chains.  In contrast to our use of strict multiplications on
$\Lambda$, notice that here we are making a weaker assertion about
module structures, which are only asserted to exist in the homotopy category of spectra.  This assertion
amounts to the existence of a map of spectra
\begin{equation}
  CF_*(H; \Lambda) \wedge \Lambda \to CF_*(H; \Lambda)
\end{equation}
such that the action of the unit of $\Lambda$ is homotopic to the identity, and the following diagram commutes up to homotopy (the left vertical map is induced by the ring structure on $\Lambda$, and the top horizontal map is the module map):
\begin{equation} \label{eq:module-homotopy-category}
  \begin{tikzcd}
    CF_*(H; \Lambda) \wedge \Lambda \wedge \Lambda \ar[r] \ar[d] &  CF_*(H; \Lambda) \wedge \Lambda \ar[d] \\
   CF_*(H; \Lambda) \wedge \Lambda \ar[r] &  CF_*(H; \Lambda). 
  \end{tikzcd}
\end{equation}
\begin{rem}
As with Remark \ref{rem:we-are-lazy}, implementing the construction of module structures at the level of homotopy categories is done for the sake of technical convenience, and Floer theory provides all the necessary data to produce the higher homotopies in the above diagrams which are required to lift $ CF_*(H; \Lambda)$ to an object in the point-set category of modules over $\Lambda$.
\end{rem}

The assertion in Proposition~\ref{prop:quantitative-Floer-chains-and-residue} that the
module structure is compatible with the filtrations is a consequence of suitable pairings 
\begin{equation}
P\bS \wedge   K(n)[\Pi_{> \lambda_1}]  \sma CF_*^{> \lambda_2}(H; \Lambda) \to CF_*^{> \lambda_1 + \lambda_2}(H; \Lambda),
\end{equation}
where we write $CF_*^{> \lambda}(H; \Lambda)$ for the $\lambda$-filtered part.

Passing to homotopy groups, these structures imply that the Floer homology
$  HF_*(H; \Lambda) $ forms a module over $\Lambda_*$ and that the
module structure is compatible with the induced filtrations on both
sides. This reduces the results of the previous section to
Proposition~\ref{prop:quantitative-Floer-chains-and-residue} and~\ref{prop:retract_CF_ordinary-chains}: 
\begin{proof}[Proof of Proposition~\ref{prop:properties_Ham_K_n}]
By considering the spectral sequence associated to the filtration of $ CF_*(H;
\Lambda)$, we see that the rank of $  HF_*(H; \Lambda) $ over
$\Lambda_*$ is smaller than the number of Hamiltonian orbits. On the other hand, by
Proposition~\ref{prop:retract_CF_ordinary-chains}, considering the spectral
sequences associated to the filtration of $C_*(M; \Lambda)$ and
$CF_*(H;\Lambda)$ implies that the the homology $H_*(M; \Lambda)$ is a submodule of $  HF_*(H; \Lambda)$ and hence has smaller rank. 
\end{proof}
\subsection{Construction of the spectra of Floer chains}
\label{sec:constr-spectra-floer}

\subsubsection{The relative cochains of Floer trajectories}
\label{sec:relat-coch-floer}
Having discussed the consequences of the existence of Floer spectra for symplectic topology, we now turn to their construction. We begin by recalling that, in the construction of ordinary Floer homology, the Floer chain complex is generated by closed time-$1$ Hamiltonian orbits. The differentials in the Floer complex are obtained by (virtual) counts of elements of moduli spaces of Floer trajectories (i.e., Morse gradient flow lines of the $1$-form in the loop space obtained by transgressing $\omega$).

For the spectral generalization, it is convenient to start by introducing the set $\cP$ of \emph{lifts of contractible Hamiltonian orbits} to the cover $\widetilde{\cL M}$ of the free loop space of $M$ associated to the image $\Pi$ of the homomorphism
\begin{equation}
  \pi_1(\cL M) \to H_2(M, \bZ).
\end{equation}
Concretely, a lift of an orbit to an element of $\cP$ corresponds to an equivalence class of choices of \emph{capping discs,} i.e., extensions of the map $S^1 \to M$ to a $2$-disc. Two discs lie in the same equivalence class whenever the corresponding map from a sphere to $M$ represents the trivial homology class. In particular, there is a well-defined action map
\begin{equation}
  \cA \co \cP \to \bR
\end{equation}
which assigns to elements of $\cP$ the integral of $\omega$ over the corresponding capping disc.  As discussed in Section \ref{sec:moduli-spaces-stable}, we can associate to each pair of elements of $\cP$ a moduli space $\fcMR(p,q)$ of (stable) Floer trajectories. The key properties that are relevant at the moment is that this is a compact Hausdorff topological space, which is empty whenever $p$ does not strictly precede $q$ with respect to the above partial order (in practice, we artificially set $\fcMR(p,p) $ to be a point, instead of being empty as geometry dictates).  The description of the topology of the moduli space of Floer trajectories yields natural closed inclusions 
\begin{equation} \label{eq:inclusion-boundary-stratum}
  \fcMR(p,q) \times \fcMR(q,r) \to \fcMR(p,r)
\end{equation}
which are associative and unital. This defines a partial ordering on elements of $\cP$ given by
\begin{equation}
  p < q \textrm{ if and only if } \fcMR(p,q) \textrm{ is non-empty}.
\end{equation}
The action of $\Pi$ on $\cP$ preserves this ordering, and is free with finite quotient given by the set of time-$1$ Hamiltonian orbits. We thus have the principal example of the notion of \emph{flow category} discussed in Section \ref{sec:topol-enrich-categ}:
\begin{defin}
 The category $\fcMR(H)$ of \emph{Floer trajectories for $H$,} has set of objects $\cP$, morphisms the moduli spaces of Floer trajectories $\fcMR(p,q)$, and compositions given by Equation \eqref{eq:inclusion-boundary-stratum}.
\end{defin}
\begin{rem}
  Cohen, Jones, and Segal \cite{CohenJonesSegal1995} introduced the term flow categories to refer to categories whose morphism spaces are manifolds with corners, so that the composition maps give rise to inclusions of boundary strata. Their notion is appropriate for Morse theory, or for Floer-theoretic contexts wherein one can choose perturbations which ensure that all moduli spaces are manifolds, which is not the case in our situation. We adopt the terminology used by Pardon in \cite{Pardon2016}, so that our notion is a generalization of the original notion. If we want to specifically refer to flow categories in which all morphism spaces are manifolds, we shall use the term \emph{manifold flow category}.
\end{rem}

At this stage, we introduce the notation $\partial \fcMR(p,q)$ for the subspace of $\fcMR(p,q)$ consisting of broken Floer trajectories, i.e., those which are images of compositions in $\fcMR(H)$. Letting $ \Omega \bbk$ denote the desuspension of an associative ring spectrum $\bbk$, we introduce a spectrum
\begin{equation}
  C^*_{\rel \partial}( \fcMR(p,q), \Omega \bbk),
\end{equation}
which we call the \emph{spectrum of relative cochains}, and whose homotopy groups compute relative cohomology
\begin{equation}
\pi_{-k}(  C^*_{\rel \partial}( \fcMR(p,q), \Omega \bbk)) \equiv H^{k+1}( \fcMR(p,q),\partial \fcMR(p,q); \bbk).
\end{equation}
As discussed in Section \ref{sec:relat-coch-orbisp}, the specific model for relative cochains that we use admits a functorial map
\begin{equation}
    C^*_{\rel \partial}( \fcMR(p,q), \Omega \bbk) \wedge  C^*_{\rel \partial}( \fcMR(q,r), \Omega \bbk) \to  C^*_{\rel \partial}( \fcMR(p,r), \Omega \bbk),
\end{equation}
which is constructed as a composition of two operations. The first is
the canonical product of cochains, which gives a map 
\begin{multline}
C^*_{\rel \partial}( \fcMR(p,q), \Omega \bbk) \wedge  C^*_{\rel
  \partial}( \fcMR(q,r), \Omega \bbk) \\
\to C^*_{\rel \partial}( \fcMR(p,q) \times \fcMR(q,r), \Omega^2 \bbk).
\end{multline}
The second map is a spectrum-level model
\begin{equation}
 C^*_{\rel \partial}( \fcMR(p,q) \times \fcMR(q,r), \Omega^2 \bbk) \to     C^*_{\rel \partial}( \fcMR(p,r), \Omega \bbk)
\end{equation}
for the boundary homomorphism
 associated to the inclusion of the product $\fcMR(p,q) \times
 \fcMR(q,r)$ as a closed subset of $ \fcMR(p,r) $, which is given by
 Equation \eqref{eq:inclusion-boundary-stratum}.  In particular, we
 obtain a \emph{category of relative cochains,} which we denote
 $C^*_{\rel \partial}( \fcMR(H), \Omega \bbk)$, with objects the
 elements of $\cP$ and morphisms the spectra of relative cochains.
This is our fundamental example of a {\em spectral category}, which is the
 analogue of the notion of a dg category, in the case where the base ring is replaced by a ring spectrum.

The action of $\Pi$ on $\cP$ naturally lifts to the category of relative cochains.  
Moreover, $\Pi$ acts on $\bbk\nhyphen\mod$ via the homomorphism
\begin{equation}
  2 c_1 \co \Pi \to 2\bZ
\end{equation}
associated to the first Chern class, together with the action of
$2\bZ$ on the category of $\bbk$-modules given by suspension; the
precise notion of group actions on spectral categories that we use is
we explained in Section~\ref{sec:groupact}.  We now can describe the
central notion of a virtual fundamental chain (see
Definition~\ref{defn:virtual_fundamental_chain} in the text):
\begin{defin}
  A \emph{virtual fundamental chain} on the moduli spaces of Floer trajectories, with coefficients in $\bbk$, is a $\Pi$-equivariant functor from the category of relative cochains to the category of $\bbk$-modules,
  \begin{equation}
 \delta \co  C^*_{\rel \partial}( \fcMR(H), \Omega \bbk) \to \bbk\nhyphen\mod,
\end{equation}
mapping each element of $\cP$ to a rank-$1$ module, i.e., a $\bbk$-module
weakly equivalent to $\Sigma^n \bbk$ for some $n \in \bZ$.
\end{defin}

In practice we will work with a homotopical representative of a
virtual fundamental chain given by a collection of composable
$\Pi$-equivariant bimodules relating $\Pi$-equivariant spectrally
enriched categories, of which the first is $C^*_{\rel \partial}(\fcMR(H),
\Omega \bbk)$ and the last is $\bbk\nhyphen\mod$.  As discussed in Section
\ref{sec:from-flow-categories}, we can compose bimodules (via the bar
construction), to obtain a bimodule over $C^*_{\rel \partial}(\fcMR;
\Omega \bbk) $ and $\bbk\nhyphen\mod$. Interpreting this bimodule as a
functor with value in $\bbk\nhyphen\mod$, we obtain the desired virtual
fundamental chain.  Alternatively, we can think of a virtual
fundamental chain as a homotopy coherent functor between
$\infty$-categories, but we will not seriously use that viewpoint in
this paper.

\begin{example} \label{ex:2-cell-complex-manifold}
  This is the first of a series of examples in which we consider the following implausibly simple situation: say that $\cP$ consists of only two elements $p< q$, and that $\Pi$ is trivial. In that case, there is a single moduli space $\fcMR(p,q)$ of interest, and the spectral category $C^*_{\rel \partial}( \fcMR, \Omega \bbk)$ consists of a single (non-trivial) morphism space which is the spectrum $C^*(\fcMR(p,q); \Omega \bbk)$.

  To prepare the groundwork for Section \ref{sec:virt-fund-chains}, we assume that we are in a generic situation in which $\fcMR(p,q) $ is manifold whose dimension we denote by $d$. As we shall see in Section \ref{sec:stable-almost-compl-1}, this manifold admits a natural stable almost complex structure, so the spectral version of Poincar\'e duality, which goes under the name Spanier-Whitehead duality \cite{SpanierWhitehead1955}, yields an equivalence
\begin{equation}
  C^*(\fcMR(p,q);  \Sigma^{d} \bbk) \htp C_*(\fcMR(p,q); \bbk),
\end{equation}
whenever $\bbk$ is a complex-oriented cohomology theory.  We shall
spend inordinate effort in
Sections~\ref{sec:comp-coch-kuran},~\ref{sec:local-glob-constr},
and~\ref{sec:coher-comp-tang} specifying such a map at the point set
level in order to establish the necessary functoriality; but the case of a
single manifold is classical (e.g., see~\cite[Part III.10.13]{Adams1974}).

Composing with the map $C_{*}(\fcMR(p,q); \bbk) \to \bbk$ given by the projection to a point gives a map
\begin{equation}
C^*(\fcMR(p,q); \Sigma^{d} \bbk)  \to \bbk
\end{equation}
which after a shift by $d+1$ gives a virtual fundamental chain $\delta$.

In this case, it is straightforward to explain how to extract a homotopy type from a virtual fundamental chain: there is a natural map
\begin{equation}
\bbk \to  C^*(\fcMR(p,q); \bbk) 
\end{equation}
obtained by pullback under the projection to a point. After a shift, the composite with the virtual fundamental chain is a map
\begin{equation}
  \Omega \bbk \to \Omega^{d+1} \bbk.
\end{equation}
We associate to $\delta$ the (homotopy) fibre of this map. If we assume that $\bbk$ is an Eilenberg-Mac Lane spectrum (i.e., we are considering ordinary homology), this would amount to the $2$-term chain complex
\begin{equation}
  \bbk \to \Omega^{d} \bbk,
\end{equation}
with differential prescribed by the fundamental class of $ \fcMR(p,q)$. Of course, this map can only be non-trivial if $d=0$, in which case it precisely corresponds to a (possibly signed) count of elements of $\fcMR(p,q) $. 
\end{example}
\begin{rem}
As should be clear from the above discussion, the homotopy type of $C_{*}(\fcMR(p,q); \bbk)$ serves as the domain of the virtual fundamental chain, and plays no other role in the construction of the homotopy type associated to this simplest case of a flow category. This is related to the idea that the image of the fundamental class of a $\bbk$-oriented manifold in the coefficient ring $\bbk$ depends only on the bordism class of the manifold, in which the key point is that the (relative) fundamental class of a manifold restricts to the fundamental class of its boundary. In our setting, we expect the existence of a notion of bordism of spaces equipped with virtual fundamental chains (and more generally of flow categories equipped with virtual fundamental chains), so that the constructions of this paper only depend on this bordism type. While we do require a result of this nature in our comparison of Morse and Floer theory (c.f. Section \ref{sec:comp-with-morse} below), we do not develop the general notion, and our proof proceeds via ad-hoc methods.
\end{rem}
The results of this paper are summarized by the following:
\begin{thm} \label{thm:Floer-trajectories-have-fundamental-chains}
 The moduli spaces of Floer trajectories admit virtual fundamental chains with coefficients in any Morava $K$-theory $K(n)$.
\end{thm}

We now explain how the above result, together with the construction in
Section~\ref{sec:from-flow-categories} of a homotopy type associated
to each virtual fundamental chain, yields a construction of spectra of
Floer chains:

\begin{proof}[Proof of Proposition~\ref{prop:quantitative-Floer-chains-and-residue}]
Since our references to Section \ref{sec:from-flow-categories} will otherwise be completely opaque, we begin by outlining the construction and properties
of the homotopy type of a virtual fundamental chain.

For each $\lambda$, we can restrict $\cP$ to the
elements of action $\leq \lambda$ and add a new terminal object;
denote this poset by $\cP_{\lambda}$.  We write $C^*_{\rel
  \partial}(\fcMR_{\lambda}; \Omega \bbk)$ for the restriction of the
relative cochains to $\cP_{\lambda}$.  Similarly, we can restrict
$\delta$ to a virtual fundamental chain $\delta_{\lambda}$.  We can
now define the homotopy type of $\delta_{\lambda}$ as the derived 
smash product:
\begin{equation}
|\delta_{\lambda}| = \bS_{\lambda} \sma_{C^*_{\rel \partial}(\fcMR_\lambda;
\Omega \bbk)} \delta_{\lambda},
\end{equation}
where here $\bS_{\lambda}$ is the unique $C^*_{\rel
  \partial}(\fcMR_{\lambda}; \Omega \bbk)$-module specified to be
trivial everywhere except at the new terminal point of
$\cP_{\lambda}$, where it is $\bS$.

The evident collapse maps $\cP_{\lambda_2} \to \cP_{\lambda_1}$ for
$\lambda_2 > \lambda_1$ give rise to maps $|\delta_{\lambda_2}| \to
|\delta_{\lambda_1}|$, and Definition~\ref{defn:invht} then constructs
the Floer chains $CF_*(H; \Lambda)$ as the homotopy limit over these
maps.  By construction, there is an evident decreasing filtration, and
we show in Theorem \ref{thm:Pi-homotopy-action} 
that $CF_*(H; \Lambda)$ is a filtered module spectrum over $\Lambda$.
A computation explained in Proposition~\ref{prop:filt-homot-type}
shows that suitable associated graded spectra are free on the number
of orbits.  Putting this all together proves
Proposition~\ref{prop:quantitative-Floer-chains-and-residue}.
\end{proof}

In the remainder of this introduction, we continue the process of
deriving the proof of our results, which have now been reduced to Theorem \ref{thm:Floer-trajectories-have-fundamental-chains}, and to Proposition
\ref{prop:retract_CF_ordinary-chains} whose proof is postponed until Section
\ref{sec:comp-with-morse}.

\subsubsection{Virtual fundamental chains from global Kuranishi charts}
\label{sec:virt-fund-chains}

The basic idea behind the proof of Theorem \ref{thm:Floer-trajectories-have-fundamental-chains} follows the outline in the toy case of Example \ref{ex:2-cell-complex-manifold}, but there are two problems: at a fundamental level, the elements of the moduli spaces of Floer trajectories may have non-trivial groups of automorphisms, so that one requires a notion of orientations for orbifolds. At a technical level, the moduli spaces we encounter are not themselves orbifolds, but have natural local presentations as quotients of topological spaces equipped with thickenings to topological manifolds. To be precise, we need the following notion:
\begin{defin} \label{def:Kuranishi-chart}
  A \emph{Kuranishi chart} is a quadruple $(X, V, s, G)$, consisting of the following data:
\begin{enumerate}
\item (Symmetry group) a finite group $G$,
\item (Thickened chart) a $G$-manifold $X$ (paracompact and Hausdorff, and possibly with boundary), for which the action of $G$ is assumed to be locally modeled after a linear representation, 
\item (Obstruction space) a finite dimensional $G$ representation $V$ equipped with an invariant inner product,
  \item (Defining section) and a $G$-equivariant map $s \co X \to V$.
  \end{enumerate}
  We write $Z = s^{-1}(0)$ for the \emph{zero locus}, and define the \emph{footprint} to be the quotient space $Z/G$.
\end{defin}
\begin{rem}
Working with smooth instead of topological manifolds, the appearance of the above definition in symplectic topology goes back to Fukaya and Ono \cite{FukayaOno1999}; they explain that it can be traced back to Kuranishi's work on deformation theory of complex manifolds. The case of topological manifolds was considered by Pardon \cite{Pardon2016}, where the condition of local linearity was not required.  The terminology of footprints is due to McDuff and Wehrheim \cite{McDuff2017}.
\end{rem}

As a toy case, the reader should have in mind a situation in which the moduli spaces $\fcMR(p,q)$ of Floer trajectories admit a \emph{global Kuranishi chart,} i.e., are homeomorphic to a footprint of a Kuranishi chart as above, where we still restrict attention to finite groups (a more general notion of a global Kuranishi chart would include the case in which $G$ is a compact Lie group, $V$ is replaced by a $G$-vector bundle over a $G$-manifold $X$, and $s$ is a section of $V$).

Our results will rely on considering Kuranishi charts with the property that the underlying manifold $X$ admits an orientation with respect to a given cohomology theory. Concretely, the orientations we consider arise from the data of stable complex structures, so we need to define this notion for topological manifolds. There are in fact two such notions; we refer the reader to \cite{MadsenMilgram1979} for the definitions used below. In order to formulate them, we recall that work of Milnor \cite{Milnor1964} and Kister \cite{Kister1964} associates to each topological manifold a fibre bundle whose structure group $\TOP(d)$ is the group of homeomorphisms of $\bR^d$ fixing the origin. This is usually called the \emph{tangent microbundle} of a topological manifold. We have natural maps $\TOP(d) \to \TOP(d+1)$, and define the direct limit to be the group $\TOP$ of \emph{stable homeomorphisms of Euclidean space.} The stable tangent space of each topological manifold $X$ is thus classified by a map
   \begin{equation}
 X \to B \TOP.
  \end{equation}
  There is a natural map $U(d) \to \TOP(2d)$ given by the action of the unitary group on $\bC^d \cong \bR^{2d}$, which induces a map $U \to \TOP$ by taking direct limits.  A \emph{stable complex lift of the tangent microbundle} of $X$ is a lift of the classifying map to $BU$. 

While we expect that all the examples we consider admit stable complex structures in the above sense, it shall turn out that the following strictly weaker notion is sufficient for our purpose: consider the inclusion $\TOP(d) \to F(d)$ into the (group-like) monoid of based self-homotopy equivalences of  $S^d$. We also have natural maps $F(d) \to F(d+1)$, and the direct limit is the monoid $F$ of stable self-homotopy equivalences of spheres, which is equipped with a natural map $\TOP \to F$. Associated to each manifold  $X$ is the composite map
   \begin{equation}
 X \to B F,
  \end{equation}
  which classifies the \emph{stable tangent spherical fibration} of $X$.  A \emph{stable complex lift of the spherical fibration} of $X$ is then a lift of this map to $BU$ (this is strictly less information than a lift of the tangent microbundle, because the composite with the map $BU \to B\TOP$ need not be homotopic to the original map). This is the notion that we will use, as the additional flexibility of working with arbitrary spherical fibrations will allow us to provide extremely concrete constructions of duality isomorphisms; we illustrate this by continuing Example \ref{ex:2-cell-complex-manifold}:
\begin{example}\label{ex:2-cell-complex-Kuranishi-M-chart}
  Assuming still that we only have two Hamiltonian orbits, we relax the condition that $\fcMR(p,q)$ is a manifold, and require only that it admit a global Kuranishi chart for which the group $G$ is trivial. In that case, we have a closed inclusion
  \begin{equation}
    \fcMR(p,q) \to X,
  \end{equation}
  defined as the zero-locus of a map $s \co X \to V$. Whenever the spherical fibration of $X$ has a stable complex lift, and assuming that $\bbk$ admits a complex orientation (see Section \ref{sec:orient-kuran-charts} below), we may apply a generalized version of Alexander duality discussed in Section \ref{sec:span-whit-dual-5} to obtain an equivalence
  \begin{equation}
      C^*(\fcMR(p,q);  \Sigma^{2n} \bbk) \htp C_*(X, X \setminus \fcMR(p,q);  \bbk),
    \end{equation}
    where the right hand side is the spectrum of relative chains of the pair $(X, X \setminus \fcMR(p,q))$ with coefficients in $\bbk$, whose homotopy groups compute relative homology
    \begin{equation}
      \pi_k (   C_*(X, X \setminus \fcMR(p,q);  \bbk)) \equiv H_k(X, X \setminus \fcMR(p,q);  \bbk),     
    \end{equation}
and which can be defined as the smash product of $\bbk$ with the mapping cone of the inclusion of $X \setminus \fcMR(p,q) $   in $X$ (c.f. Definition \ref{def:mapping_cone}).    We call this spectrum the \emph{virtual cochains} following Pardon \cite{Pardon2016}. The defining map $s$ yields a map
    \begin{equation}
    C_*(X, X \setminus \fcMR(p,q);  \bbk) \to C_*(V, V \setminus 0; \bbk), 
     \end{equation}
     of relative chains. The right hand side is a rank-$1$ module over $\bbk$ (corresponding to the fact that the reduced homology of a sphere has rank $1$), so that the composite
     \begin{equation}
       C^*(\fcMR(p,q);  \Sigma^{2n} \bbk) \to   C_*(V, V \setminus 0; \bbk), 
     \end{equation}
defines the desired virtual fundamental class after desuspension.
   \end{example}

   The next step is to consider global charts with non-trivial group actions. The solution for formulating Poincar\'e duality for orbifolds with coefficients in ordinary cohomology is well-understood: it suffices to invert the order of the isotropy in the coefficient ring. One way to understand why this works is to consider the (Borel) equivariant cohomology of a space $Z$ equipped with the action of a finite group $G$. This cohomology group is the ordinary cohomology of the space
   \begin{equation}
     BZ \equiv E G \times_{G} Z
   \end{equation}
obtained by applying the Borel construction to $Z$. If we work with a coefficient ring $R$ in which the order of $G$ is inverted, the natural maps
   \begin{equation} \label{eq:isomorphism_cohomology-quotient-invariant-cohomology}
     H^*(Z/G; R) \to H^*(BZ; R) \to (H^*(Z; R))^G
   \end{equation}
   are isomorphisms. If $Z$ is a manifold, this allows one to deduce Poincar\'e duality for $Z/G$ from Poincar\'e duality for $Z$ (this is essentially the point of view used by Fukaya and Ono \cite{FukayaOno1999}), while if $Z$ is a closed subset of a manifold $X$, we can use the above to relate the cohomology of the quotient to the equivariant homology of the pair $(X, X \setminus Z)$ (this is the point of view used in Pardon's work \cite{Pardon2016}). Of course, both of these points of view need to be globalised in order to be of much use.

Since we cannot control the order of $G$, the above strategy
ultimately leads to working with rational coefficients. Studying
generalized cohomology theories therefore requires a new idea, and the
key point is to focus on the middle term of Equation
\eqref{eq:isomorphism_cohomology-quotient-invariant-cohomology}. Since
this involves a Borel construction, we are led to employ the
techniques of equivariant stable homotopy theory.  More specifically,
Equation~\eqref{eq:isomorphism_cohomology-quotient-invariant-cohomology}
is a version of the {\em norm map} from the homotopy orbits to the
homotopy fixed-points, the study of which is the subject of the theory
of Tate cohomology and Tate spectra.

The appearance of the Morava $K$-theories in our work is now suggested
by the foundational observation of Ravenel~\cite{Ravenel82} that the
Morava $K$-theory $H^*(BG, K(n))$ for a finite group $G$ has finite rank
over $K(n)_*$.  Greenlees and Sadofsky~\cite{GreenleesSadofsky1996}
interpreted this to show that $BG$ is self-dual with respect to Morava
$K$-theories in the sense that there is a corresponding isomorphism
\begin{equation}
H_*(BG, K(n)) \cong  H^*(BG, K(n)) 
\end{equation}
which arises from a comparison of spectra.

It is natural to expect now that there should be some kind of
Poincar\'e duality for orbifolds with coefficients in Morava
$K$-theory, and indeed we have the following result of Cheng \cite{Cheng2013}:
   \begin{thm}
     If a finite group $G$ acts on a closed smooth manifold $X$ of dimension $2d$, then the datum of an almost complex structure on $X$ that is preserved by $G$ induces an isomorphism
     \begin{equation}
       H^*(BX; K(n)) \cong H_{2d-*}(BX; K(n)),
     \end{equation}
     between the $K(n)$ homology and cohomology of the classifying space $BX \equiv X \times_{G} EG$.
   \end{thm}
We reprove this result in greater generality, using the same basic building blocks: this duality follows from 
Spanier-Whitehead duality, the vanishing of Tate cohomology for $K(n)$-local theories~\cite{GreenleesSadofsky1996}, and the Adams isomorphism \cite{Adams1984}. These matters are discussed at great length in Section \ref{sec:comp-coch-kuran} where we refer to the appendices for some of the technical details, but the following continuation of Example \ref{ex:2-cell-complex-Kuranishi-M-chart} provides a basic summary.
   \begin{example}\label{ex:2-cell-complex-Kuranishi-P-chart}
     Assume that $ \fcMR(p,q)$ admits a global Kuranishi chart such that $V$ is a complex representation, and $X$ admits a stable almost complex structure.  Following Greenlees and May's treatment of Tate cohomology \cite{GreenleesMay1995}, we consider the composition of $G$-equivariant maps
     \begin{equation} \label{eq:transfer-map}
      E G_{+}  \wedge C^*(Z; \bbk) \to C^*(Z; \bbk)   \to C^*(EG ; C^*( Z; \bbk)). 
    \end{equation}
    induced by the projection $E G_+ \to S^0$.    The key result of Greenlees and Sadofsky~\cite{GreenleesSadofsky1996} that we use is that, whenever $\bbk$ is a Morava $K$-theory, the above map is a $G$-equivariant equivalence, so that it can be inverted up to homotopy to yield a map
     \begin{equation}
       C^*(B Z; K(n)) \htp C^*(E G; C^*(Z; K(n)))^G  \dashrightarrow  \left(E G_{+}  \wedge C^*(Z; K(n))\right)^G,
     \end{equation}
after passing to fixed points, where we use $\dashrightarrow $ to
indicate the fact that we invert a canonically given map (more
generally, if we are representing a map in the homotopy category by a
zig-zag with backwards maps equivalences).  Next, we use the assumption that $X$ has a stable almost complex structure and that $K(n)$ is complex oriented to construct an equivalence
     \begin{equation}
    C^*(Z; \Sigma^{d} K(n))  \dashrightarrow   C_*(X, X \setminus Z; K(n)),
    \end{equation}
which is again given by Spanier-Whitehead duality, where $d$ is the dimension of $X$.  This equivalence can be factored through $G$-equivariant maps, so that we obtain a homotopy class of maps
    \begin{equation}
      \left(E G_{+}  \wedge C^*(Z; \Sigma^{d} K(n))\right)^G \dashrightarrow \left(E G_{+}  \wedge   C_*(X, X \setminus Z; K(n))  \right)^G.
    \end{equation}
    Finally, we use the Adams isomorphism, a deep result in
    equivariant stable homotopy theory (discussed in Appendix~\ref{sec:adams-isom-ambid}) which asserts that the fixed points and orbits of free spectra are equivalent, to get a map
    \begin{equation}
\left(E G_{+}  \wedge C_*(X, X \setminus Z; K(n))   \right)^G \dashrightarrow  C_*(B X, B (X \setminus Z); K(n)),
\end{equation}
where $B X$ and $B(X \setminus Z)$ again refer to the Borel constructions on the $G$-spaces $X$ and $X \setminus Z$. We have, as before, a natural map
\begin{equation}
C_*(B X,B (X \setminus Z); K(n))  \to C_*(B V, B (V \setminus 0); K(n))  
\end{equation}
and the fact that $V$ is a complex representation and $K(n)$ is
complex oriented trivialize the action of $G$ on the representation
sphere $S^V$ in the category of $K(n)$-modules (see
Section~\ref{sec:eqMU} and in particular
Corollary~\ref{cor:thom-zigzag} for discussion of this untwisting),
which yields a map
\begin{equation}
  C_*(B V, B (V \setminus 0); K(n))  \to  C_*(V, V \setminus 0; K(n))
\end{equation}
splitting the natural homotopy class of maps in the other direction. Applying the appropriate desuspensions provides a map
\begin{equation}
  C^*(B Z; \Omega K(n)) \to \Omega^{d+1} C_*(V, V \setminus 0; K(n)),
\end{equation}
whose composition with the pullback of cochains from $Z$ to $BZ$ is the desired fundamental chain.
\end{example}

\subsubsection{Kuranishi flow categories}

In the finite-dimensional approach to virtual fundamental chains, one starts by observing that the moduli spaces of Floer trajectories admit covers by footprints of Kuranishi charts. The difficulty in formulating a global notion lies in stating the required data along overlaps of charts. Our solution is to introduce,  in Definition \ref{def:stratified-Kuranishi-charts}, a monoidal category $\Kur$ of \emph{Kuranishi charts}. We believe that this brings substantial clarity to various constructions, as much of the content of this paper can be interpreted in terms of constructions of lax monoidal functors from $\Kur$ (and related categories) to categories of spectra.
\begin{rem}
Since the theory of Kuranishi charts contains the theory of orbifold charts, the reader may be concerned by the fact that a complete theory of maps of orbifolds needs to be formulated in terms of a $2$-category (see e.g. Joyce \cite{Joyce2019}). We avoid these $2$-categorical subtleties because the explicit geometric constructions of Floer theory can be performed at the $1$-categorical level.
\end{rem}

We begin by formulating the variant of the notion of Kuranishi structure \cite{FukayaOno1999}, Kuranishi space \cite{Joyce2019},  implicit atlas \cite{Pardon2016}, or Kuranishi atlas \cite{McDuff2017} that we use. We define a \emph{Kuranishi presentation} of $\fcMR(p,q)$ to be a diagram in $\Kur$ (i.e., a functor from a category $A(p,q)$ to $\Kur$), equipped with a homeomorphism between the colimit of the corresponding diagram in the category of topological spaces (under the footprint functor) with $\fcMR(p,q)$. We impose some mild technical conditions on our diagrams, the most important of which is the requirement that the collection of transition functions relating charts covering each point in $\fcMR(p,q)$ is contractible. The construction of such a Kuranishi presentation is given in Section \ref{sec:kuran-pres-moduli}; as with all the other finite dimensional approaches to the construction of virtual fundamental chains, the thickenings are moduli spaces of Floer trajectories with additional marked points, which solve perturbed pseudo-holomorphic curve equations. The additional choice of marked points gives rise to the group $G$ of symmetries, and the perturbations are chosen from the vector space $V$. 

Recalling that the spaces $\fcMR(p,q)$ are the morphisms of a category, we need to formulate the multiplicativity of Kuranishi presentations. This takes the form of functors $A(p,q) \times A(q,r) \to A(p,r)$, which are associative in the sense that they satisfy the axioms of a bicategory, and of coherent natural transformations in the diagram
\begin{equation}
  \begin{tikzcd}
    A(p,q) \times A(q,r) \ar[r] \ar[d] &  A(p,r) \ar[d] \\
    \Kur \times \Kur \ar[r] & \Kur.
  \end{tikzcd}
\end{equation}
Such structure is encoded by the notion of a lift of a flow category to a Kuranishi flow category (see Section \ref{sec:kuran-orbisp-flow}), and the construction of such a lift in the context of Hamiltonian Floer theory is given in Section \ref{sec:kuran-pres-prod}.
\begin{rem}
  The reader who is familiar with Floer theory will observe that the space $\fcMR(p,q)$ has a virtual dimension $\virdim \fcMR(p,q) \in \bZ$, which should agree with the quantity $\dim X - \dim V$ for each Kuranishi chart of this space. The maps that arise in the Kuranishi presentation of $\fcMR(p,q)$ are thus maps between charts of the same virtual dimension.  On the other hand, we have
  \begin{equation}
    \virdim \fcMR(p,q) +  \virdim \fcMR(q,r) = \virdim \fcMR(p,r) - 1,
  \end{equation}
so that the product of charts for $ \fcMR(p,q) $ and $ \fcMR(q,r) $ has virtual dimension one smaller than the virtual dimension of charts for $\fcMR(p,r) $. Instead of considering arbitrary maps between charts of different virtual dimension, we shall restrict attention to those that arise by composing equi-dimensional maps with inclusions of boundary strata. This is why the definition of the category $\Kur$ in Section \ref{sec:strat-kuran-charts} incorporates the data of stratifications.
\end{rem}

In order to discuss the relevance of Kuranishi flow categories to the
construction of virtual fundamental chains, we shall introduce the
notion of a \emph{tangentially twisted fundamental chain.} It is extremely convenient for this purpose to consider the (Milnor) model $MX$ for the spherical tangent fibration of a topological manifold $X$, whose fibre at a point $x \in X$ is the cone of the inclusion of $X \setminus x$ in $X$ (see Section \ref{sec:span-whit-dual-5}). As with every spherical fibration over $X$, we obtain a spectrum of tangentially twisted cochains
\begin{equation}
  C^*(Z; MX \wedge \bbk)
\end{equation}
associated to every map $Z \to X$ and every spectrum $\bbk$, whose homotopy groups compute cohomology with twisted coefficients.

If $Z$ is a compact subset of $X$, the main advantage of this model is the existence of a natural equivalence
\begin{equation}
  C_*(X, X \setminus Z) \to C^*(Z; MX)
\end{equation}
which realizes Spanier-Whitehead duality. If $X$ is not assumed to be closed, then the right hand side should incorporate a condition of compact support as in Poincar\'e duality for non-compact manifolds, and if $Z$ intersects the boundary of $X$, we should consider the relative cochains of the pair $(Z, \partial X \cap Z)$ as in Lefschetz duality.

In the setting of Kuranishi charts, it is more natural to consider the \emph{virtual tangent bundle,} which is the desuspension $MX^{-V}$ by the $G$-representation $V$. To keep track of equivariance, we pull back $MX^{-V}$ to the Borel construction $BZ$, and consider, as in Example \ref{ex:2-cell-complex-Kuranishi-P-chart}, the Borel equivariant cochains
\begin{equation}
   C^{*,c}_{\rel \partial}(BZ; MX^{-V} \wedge \bbk),
 \end{equation}
 which are compactly supported (along the $Z$ direction), relative the intersection of $Z$ with the boundary of $X$. Specializing to the case $\bbk$ represents a Morava $K$-theory $ K(n)$, we again refer to Section \ref{sec:span-whit-dual-5} for the construction of a zig-zag of equivalences
   \begin{equation}
     C^{*,c}_{\rel \partial}(BZ; MX^{-V} \wedge \bbk) \leftarrow \cdots \to  C_*(X, X \setminus Z; \Sigma^{-V} \bbk)
   \end{equation}
   in the case of a single chart. The results of Section \ref{sec:local-glob-constr} then show that this zig-zag can be made homotopy coherent in an appropriate sense. Using the map of pairs $(X, X \setminus Z) \to (V, V \setminus 0)$ as in Example \ref{ex:2-cell-complex-Kuranishi-P-chart}, it is then relatively straightforward to map $ C_*(X, X \setminus Z; \Sigma^{-V} \bbk)$ to a rank-$1$ module over $\bbk$.

   The outcome of this discussion is that every Kuranishi flow category, and in particular the lift of $\fcMR(H)$ constructed in Section \ref{sec:kuran-flow-categ-1}, admits a twisted fundamental chain. In order to formulate this notion, we shall consider, for each pair $(p,q)$ of Hamiltonian orbits, a spectrum denoted
   \begin{equation} \label{eq:Borel-equivariant-twisted-cochains}
     C^{*,c}_{\rel \partial}(B\cZ; M\cX^{-V} \wedge \bbk)(p,q)
   \end{equation}
which is obtained by gluing together the Borel equivariant spectra of twisted cochains $ C^{*,c}_{\rel \partial}(BZ; MX^{-V} \wedge \bbk)$ over all Kuranishi charts of $ \fcMR(p,q)$.
   \begin{rem}
   The precise construction of these twisted cochains will incorporate an additional shift by the difference $d_p - d_q$ of integers associated to each orbit, which record the degrees of the corresponding cells. We shall delay discussing this datum in the introduction until its importance to the existence of complex orientation becomes clear below.   
   \end{rem}
   \begin{prop} \label{prop:existence-twisted-fundamental-chain}
     If $\bbk$ represents a Morava $K$-theory, there is a $\Pi$-equivariant spectral category $ C^{*,c}_{\rel \partial}(B\cZ; M\cX^{-V} \wedge \bbk) $ with objects $p \in \cP$ and morphisms given by Equation \eqref{eq:Borel-equivariant-twisted-cochains}, which admits a $\Pi$-equivariant functor
     \begin{equation} \label{eq:twisted-fundamental-chain}
      C^{*,c}_{\rel \partial}(B\cZ; M\cX^{-V} \wedge \bbk) \to \bbk\nhyphen\mod 
    \end{equation}
assigning to each element of $\cP$ a rank-$1$ module over $\bbk$.
  \end{prop}
  \begin{proof}
    The construction of the spectral category of tangentially twisted cochains is given in Section \ref{sec:span-whit-dual-3}. In Section \ref{sec:span-whit-dual-1}, we compare this to a spectral category built from the virtual cochains, which is proved, in Section \ref{sec:impl-atlas-with} to admit the desired functor to the category of modules.
  \end{proof}
  \begin{rem}
    As should be clear from the above discussion, Proposition
    \ref{prop:existence-twisted-fundamental-chain} is a general result
    about lifts of flow categories to Kuranishi presentations, and its
    proof uses no specific feature of Hamiltonian Floer theory. The
    fact that we restrict its statement to Morava $K$-theory is thus
    slightly surprising: the issue lies in Equation \eqref{eq:transfer-map},
    which is the step comparing the twisted cochains, which are built
    from cochains of sections of spectra over classifying spaces of
    finite groups, with the chain theory.  This
    comparison holds more generally for \emph{$K(n)$-local} cohomology theories,
    which in particular have the property that they vanish on spaces which are $K(n)$-acyclic. 
\end{rem}

\subsubsection{Oriented Kuranishi charts}
\label{sec:orient-kuran-charts}

At this stage, it remains to compare twisted and ordinary cochains in order to construct virtual fundamental chains from twisted ones. As discussed in Section \ref{sec:virt-fund-chains}, this amounts to a lifting problem for the tangent spherical fibration of Kuranishi charts from the space $F$ of stable self homotopy equivalences of spheres to the unitary group $U$. Our lifting will factor through the stable orthogonal group $O = \colim_d O(d)$, in the sense that we exhibit an equivalence between the stable spherical fibration $MX$ of the Kuranishi charts that we consider, and the sphere bundle $S^{TX}$ of a vector bundle $TX$.

\begin{rem} \label{rem:how-to-deal-with-smooth-charts}
If we worked with smooth Kuranishi charts, the vector bundle $TX$
would, as the notation indicates, be the ordinary tangent bundle. In
that case, it is standard to use an exponential map to compare
$S^{TX}$ to the model $MX$ for the spherical fibration which we find
convenient for proving Spanier-Whitehead duality. Of course, the
exponential map uses a choice (of Riemannian metric), which one must
keep track of while trying to prove the functoriality and multiplicativity of the comparison. While this is a slightly technical point to discuss in the introduction, our solution to this problem may be of interest to readers who would appreciate a preview of the relevant part of the paper (Sections \ref{sec:smoothly-fiber-kuran-1}, \ref{sec:tang-spac-whitn} and \ref{sec:nash-spher-fibr-1}).

Consider, as Nash did in \cite{Nash1955}, the space $NX$ of Moore paths in $X$ which do not return to their starting point so that, in particular, the only constant path is the one parameterized by an interval of length $0$. If $X$ is smooth, we impose in addition the condition that all paths are differentiable at the origin. The derivatives at the origin and the evaluation at the other endpoint yield natural maps
\begin{equation}
 TX  \leftarrow NX \to X \times X
\end{equation}
of fibre bundles over $X$. Denoting by $0$ the section of $NX$ consisting of constant paths, we find that the above maps restrict to maps
\begin{equation}
   TX \setminus 0  \leftarrow NX \setminus 0 \to X \times X \setminus X.
\end{equation}
Defining $NX|0$ and $TX|0$ to be the fibre bundles obtained by taking the cone over the complement of $0$ at each point, we have induced maps
\begin{equation}
   TX|0  \leftarrow NX|0 \to MX,
\end{equation}
which are equivalences of spherical fibration over $X$. It is then easy to compare the spherical fibration $TX|0$ with $S^{TX}$ via another zig-zag.

It should be clear to the reader that the construction above is free of choices, functorial, and multiplicative. Choices would enter in picking inverse homotopy equivalences, but we systematically avoid choosing inverses in this paper.
\end{rem}

One the main technical points of this paper is to avoid using the smooth structure on the moduli spaces of pseudo-holomorphic maps. There are potential approaches to proving this smoothness result both from the finite dimensional point of view \cite{FukayaOhOhtaOno2016} or using the theory of polyfolds \cite{HoferWysockiZehnder2009}. Instead, we introduce the notion of a \emph{flag smooth} manifold, which consists of the data of a topological manifold $X$ equipped with a topological submersion  over a smooth manifold, together with a smooth structure on all fibers (a related notion appeared already in \cite{Siebert1999}). Note that such a structure yields a vector bundle on $X$, which we denote $TX$, defined as the direct sum of the tangent spaces of the fibers with the pullback of the tangent space of the base. It is quite easy to adapt the Nash argument comparing different models of tangent spherical fibrations to this setting. 

\begin{rem}
  A more careful construction than the one given in this paper would show that a flag smooth structure with total space $X$ gives a lift of the classifying map of the topological microbundle $X$  to a vector bundle (rather than simply a lift of the corresponding spherical fibration).  In dimension greater than $5$, smoothing theory \cite{KirbySiebenmann1977} implies that a lift of the topological microbundle of a manifold determines an underlying smooth structure.
  
The results that we need are not so delicate as to require smoothing theory. One way to think about this is that the homotopy type associated to moduli spaces depends only on their structured bordism class. Standard arguments in Pontryagin-Thom theory imply that the relevant smooth and topological $G$ bordism groups are isomorphic (the case of dimension $4$ is more delicate than that of higher dimensional manifolds, but was worked out by Freedman and Quinn \cite{FreedmanQuinn1990}).
\end{rem}

By simply requiring the datum of a flag smooth structure on the thickening, it is straightforward to extend the notion of a flag smooth manifold to that of a flag smooth Kuranishi chart. It is less straightforward to formulate the functoriality of the notion. This is done in Section \ref{sec:fiber-kuran-pres}, where we construct a monoidal category $\Kur^{fs}$ of flag smooth charts, using a notion that is entirely motivated by the geometry of moduli spaces of Floer trajectories (c.f. Definition \ref{def:smoothly_fibred_complete}). This category admits a forgetful functor to $\Kur$, so that it makes sense to define a \emph{flag smooth Kuranishi flow category} to be a lift of the diagram underlying a Kuranishi flow category from $\Kur$ to $\Kur^{fs}$.

The relevance of this notion to Floer theory is given by the fact that the thickenings that we construct admit forgetful maps to abstract moduli spaces of curves with marked points. These abstract moduli spaces have natural smooth structures essentially coming from their description in terms of submanifolds of complex algebraic varieties. The fibers of these forgetful maps are naturally embedded inside smooth Banach manifolds as zero loci of smooth Fredholm sections of Banach bundles, and hence are smooth whenever the section is transverse to the origin. Putting these ingredients together yields the following result which is proved in Section \ref{sec:fibr-smooth-pres}.
\begin{prop}
The flow category $\fcMR(H)$ lifts to a flag smooth Kuranishi flow category. \qed
\end{prop}

Having constructed a lift of the spherical fibration to a real vector bundle, we now consider the problem of lifting the corresponding stable bundle to a stable complex vector bundle. The relevant notion is that of a \emph{relative stable complex orientation,} i.e. we have for each orbit $p$ a stable complex vector space
\begin{equation}
  V_p = (V_p^+, V_p^-)
\end{equation}
whose formal difference is the dimension of the corresponding cell of the Floer chains. This is one of the additional complications of Floer theory when thought of as infinite dimensional Morse theory, as we can assume $V_p^- = 0$ in the finite dimensional situation. A complex orientation of a Kuranishi chart relative the pair $V_p$ and $V_q$ thus consists of a $G$-equivariant complex vector bundle $I$ over each chart, together with a $G$-representation $W$ and an isomorphism
\begin{equation}
  V_p^+ \oplus \bR  \oplus W \oplus TX \oplus V_q^- \cong  V_p^- \oplus W  \oplus I \oplus V_q^+. 
\end{equation}
There is a natural topology on the space of such isomorphisms, so we construct a category $\Kur^{\ori}(V_{p}, V_{q})$ of Kuranishi charts with stable complex structures relative the pairs $V_p$ and $V_q$, as a category internal to topological spaces (i.e. with both objects and morphisms equipped with a natural topology, and continuous compositions). This accounts for the functoriality of orientations. The multiplicativity of orientations then arises from composition functors
\begin{equation}
  \Kur^{\ori}(V_{p}, V_{q}) \times \Kur^{\ori}(V_{q}, V_{r}) \to \Kur^{\ori}(V_{p}, V_{r})
\end{equation}
which equip the collection of such categories with the structure of a bicategory which we denote $\Kur^\ori$. Considering the monoidal category $\Kur$ as a bicategory with a single $0$-cell, we have a forgetful 2-functor
\begin{equation}
  \Kur^\ori \to \Kur,
\end{equation}
so it again makes sense to define a \emph{complex oriented Kuranishi flow category} as a lift of a Kuranishi flow category to $ \Kur^\ori$.

Returning to Floer theory, we have the following summary of the results of Section \ref{sec:stable-almost-compl-1}:
\begin{thm} \label{thm:Floer-trajectory-category-lifts}
The flow category $\fcMR(H)$ admits a lift to a complex-oriented Kuranishi flow category. \qed
\end{thm}

\begin{proof}[Proof of Theorem \ref{thm:Floer-trajectories-have-fundamental-chains}]
Given Theorem \ref{thm:Floer-trajectory-category-lifts}, this is essentially a consequence of the long work of Sections \ref{sec:impl-atlas-with}, \ref{sec:local-glob-constr}, and \ref{sec:coher-comp-tang}, which proves that every complex oriented Kuranishi flow category is equipped with a virtual fundamental chain with coefficients in Morava $K$-theory. The following diagram may help the reader trace through the various maps; each entry is a spectrally enriched $\Pi$-equivariant category associated to the complex oriented Kuranishi presentation of $\fcMR(H)$ \begin{equation}
    \begin{tikzcd}
      C^{*,c}_{\rel \partial}(\cB {\cZ}; T {X}|0^{- V} ) &  C^{*,c}_{\rel \partial}(\cB {\cZ}; N {X}|0^{- V} ) 
     \arrow{l}[swap]{\textrm{Evaluation}}  \ar[d,"\textrm{Evaluation}"] \\
  C^{*,c}_{\rel \partial}(\cB {\cZ}; S^{T{X} - V} \wedge \bbk)  \ar[d,dashed,swap,"\textrm{Index}"] \ar[u]  & C^{*,c}_{\rel \partial}(\cB {\cZ}; M {X}^{- V} ) \\
  C^{*,c}_{\rel \partial}(\cB {\cZ}; S^{I - V - \ell} \wedge \bbk)  \ar[d,dashed,swap,"\textrm{Orientation}"] &   (E G_{+}  \wedge  C^{*,c}_{\rel \partial}( {\cZ}; M {X}^{- V} ))^{G} \ar[u,dashed,swap,"\textrm{Ambidexterity}"] \\
 C^{*,c}_{\rel \partial}(\cB {\cZ}; \Omega \bbk)  \ar[d,dashed,swap,"\textrm{Compactness}"] & (E G_{+}  \wedge \cX|\cZ^{- V} )^{G} \ar[u,dashed,swap,"\textrm{Duality}"] \ar[d,dashed,"\textrm{Adams}"]   \\
C^*_{\rel \partial}(\cB {\cZ}; \Omega \bbk) & B \cX|\cZ^{- V} \ar[d,"\textrm{Augmentation}"]    \\
C^*_{\rel \partial}(\fcMR(H); \Omega \bbk) \ar[u,"\textrm{Pullback}"] & \bbk\nhyphen\mod
    \end{tikzcd}
  \end{equation}
  More precisely, each solid morphism is a $\Pi$-equivariant functor, and each dashed arrow is a $\Pi$-equivariant bimodule representing an equivalence. As discussed in Section \ref{sec:from-flow-categories}, we can compose bimodules (via the bar construction), to obtain a bimodule over $C^*_{\rel \partial}(\fcMR(H); \Omega \bbk) $ and $\bbk\nhyphen\mod$. Interpreting this bimodule as a functor with value in $\bbk\nhyphen\mod$, we obtain the desired virtual fundamental chain.
\end{proof}

\subsubsection{Comparison with Morse theory}
\label{sec:comp-with-morse}

We now return to discuss Proposition \ref{prop:retract_CF_ordinary-chains}: in order to compare Floer homology with ordinary homology, we use Morse theory as an intermediate step, as in \cite{Fukaya1997b,PiunikhinSalamonSchwarz1996}. The starting point is a Morse-Smale function $f$ on $M$. 
Consider the set $\cX$ of lifts of critical points of $f$ to the cover $\widetilde{\cL M}$, partially ordered by the values of the action functional and by those of $f$. These are the objects of a flow category $\Tbar$, whose morphisms $ \Tbar(x,y) $ are compactified moduli spaces of gradient flow lines connecting critical points $x$ and $y$. The constructions of this paper, implemented in a much simpler setting with trivial isotropy groups and obstruction spaces, yield a virtual fundamental chain on this flow category, with coefficients in $K(n)$. We denote by $CM_*(f; \Lambda)$ the corresponding homotopy type constructed using the methods of Section \ref{sec:from-flow-categories}. The following result, which is essentially due to Cohen, Jones, and Segal \cite{CohenJonesSegal1995}, follows from the results proved in  Appendix \ref{sec:morse-theor-setup}:
\begin{thm}
  There is a natural equivalence $CM_*(f; \Lambda) \cong C_*(M; \Lambda)$.
\end{thm}
\begin{proof}
In Appendix \ref{sec:morse-theor-setup} we construct a homotopy type $CM_*(f; \bbk)$ from the flow category with objects critical points of $f$.  The inclusion of constant filling discs gives a distinguished map $M \to \widetilde{\cL M}$, so that we can write each element of $\cX$ canonically as a pair $([x],\pi)$ with $[x]$ a critical point of $f$ and $\pi$ an element of $\Pi$. The  key point is that the space of morphisms between objects $([x],\pi)$ and $([x'],\pi')$ is empty whenever $\pi \neq \pi'$. 
 This implies that we have an equivalence
  \begin{equation}
    CM_*(f; \Lambda) \cong CM_*(f; \bbk) \wedge_{\bbk} \Lambda. 
  \end{equation}
  In Proposition \ref{prop:computation-Morse}, we prove the equivalence of $ CM_*(f; \bbk)$ with $C_*(M; \bbk)$. Since
  \begin{equation}
    C_*(M; \Lambda) \cong C_*(M; \bbk) \wedge_{\bbk} \Lambda,
  \end{equation}
the desired result follows.  
\end{proof}

The comparison of Morse and Floer theory can then be implemented as follows: consider the partial ordering on $\cX \amalg \cP$ extending the partial order of these two sets by
\begin{equation} \label{eq:Hofer_energy-constant}
  x < p \textrm{ if and only if } \cA(x) < \cA(p) - c
\end{equation}
for some constant $c$ that will depend on the Hamiltonian $H$ (i.e., its Hofer energy). In particular, a Hamiltonian orbit of $H$ never precedes a critical point of $f$. In Section \ref{sec:hamilt-floer-theory-1}, we construct a flow category $\fcMR(f,H)$ with such objects, with morphisms between critical points given by moduli spaces of gradient trajectories, those between Hamiltonian orbits by moduli spaces of Floer trajectories, and those between a critical point $x$ and an orbit $p$ given by a mixed moduli space $\fcMR(x,p)$ consisting of a pseudo-holomorphic plane with a marked point, which is asymptotic to $p$ along the end, and a gradient flow line from $x$ to the image of the marked point, at illustrated in Figure \ref{fig:moduli-spaces-mixed}.
\begin{figure}[h]
  \centering
  \begin{tikzpicture}
    \coordinate[label= left:$x$] (x) at (-4,0);
    \draw [thick,->] (x) -- (-2,0);
    \draw [thick](x) -- (0,0);
    \draw  (4,1) -- (3,1) .. controls (2, 1) and (0,1) .. (0,0) ;
    \draw  (4,-1) -- (3,-1) .. controls (2,-1) and (0,-1) .. (0,0) ;
    \draw (4,1) arc (90:-90:.2 and 1);
    \draw  [dashed] (4,1) arc (90:270:.2 and 1);
     \draw [fill,color=black] (0,0) circle (.05);
   \coordinate[label= right:$p$] (p) at (4.2,0);
  \end{tikzpicture}
  \caption{A representation of an element of the moduli space $\fcMR(x,p) $.}
  \label{fig:moduli-spaces-mixed}
\end{figure}
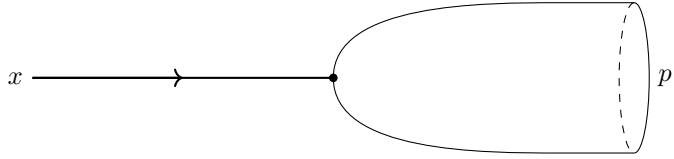
\begin{prop}
The flow category $\fcMR(f,H)$ admits a virtual fundamental chain with coefficients in $K(n)$. 
\end{prop}
\begin{proof}
  The construction of a lift to a complex oriented Kuranishi flow category is done in Section~\ref{sec:stable-almost-compl-1}. The existence of the fundamental chain then follows from the results of Section~\ref{sec:impl-atlas-with},~\ref{sec:local-glob-constr}, and~\ref{sec:coher-comp-tang}.
\end{proof}

At this stage, we can appeal to Section~\ref{sec:maps-homotopy-types},
in which we prove that the homotopy type associated to the fundamental
chain on $\fcMR(f,H)$ is the cofiber of a map  
\begin{equation}
  CM_*(f; \Lambda)  \to CF_*(H; \Lambda), 
\end{equation}
induced from maps
\begin{equation}
CM_*^{\lambda}(f; \Lambda) \to CF_*^{\lambda+c}(f; \Lambda)
\end{equation}
where $c$ is the constant appearing in Equation \eqref{eq:Hofer_energy-constant}.

We complete the reduction of the results stated in the introduction to
those proved in the paper with:

\begin{proof}[Proof of Proposition  \ref{prop:retract_CF_ordinary-chains}]
Having constructed a map of homotopy types, it remains to show that it
splits. In Section \ref{sec:morse-theory}, we construct a flow
category $\fcMR(f,H,f)$ with objects indexed by three copies of the
set $\cX$ of critical points of the Morse function $f$, which are
ordered as $\cX_- < \cX_0 <\cX_+$, and one copy of the set $\cP$ of
Hamiltonian orbits. The key properties of this category are that it
admits $\fcMR(f,H)$ as the full subcategory with objects $\cX_-$ and
$\cP$, that $\cX_0$ and $\cP$ consist of incomparable objects, and
that the full subcategory with objects  $\cP \amalg \cX_+ $ has the
property that there are no morphisms from $p$ to $x_+$ unless
$\cA(x) \leq \cA(y) - c$ and $f(x) \leq f(y)$, with this space of
morphisms consisting of a point whenever $x=y$.

In Section \ref{sec:stable-almost-compl-1}, we prove that this
category admits a lift to a complex oriented Kuranishi flow category,
and hence determines a homotopy type over $\Lambda$.  It is clear from
the construction of $\fcMR(f,H,f)$ that the hypotheses of
Proposition~\ref{prop:contractible} hold for the subcategories with
objects $\cX_- \amalg \cX_0 $ and $\cX_+ \amalg \cX_0 $; i.e., that
the associated homotopy types are acyclic.  We can then apply
Proposition~\ref{prop:retract} to conclude that this data specifies a
retraction as a $\Lambda$-module of the map represented by
$\fcMR(f,H)$ from Morse to Floer chains.
\end{proof}

\subsection{Outline of the paper}
\label{sec:outline-paper}

We end this introduction with Figure \ref{fig:Dependency}, which
provides a dependency diagram for the contents of this paper: dashed
arrows indicate minimal dependencies, i.e. that the essential point of
the target section can be understood without reference to the source: 

\begin{figure*}[h]
  \centering
  \begin{tikzpicture}[node distance=1.5cm]
    \tikzstyle{box} = [rectangle, rounded corners, minimum width=2cm, minimum height=1cm,text centered, draw=black]
    \tikzstyle{part} = [minimum width=2cm, minimum height=1cm,align=center]

    \node (Kur) [part] {Kuranishi \\ presentations};
    \node (floer) [part, left of=Kur, xshift=-1.5cm] {Floer \\ theory};
\node (AlgTop) [part, right of=Kur, xshift=1.5cm] {Algebraic \\ Topology};
    
  \node (flow) [box, below of=Kur] {Section \ref{sec:topol-orfif-flow}};
  \node (kuranishi) [box, below of=flow] {Section \ref{sec:kuran-flow-categ-2}};
  \node (kursmooth) [box, below of=kuranishi] {Section \ref{sec:fiber-kuran-pres}};
    \node (kuror) [box, below of =kursmooth] {Section \ref{sec:compl-orient-kuran-1}};

  \node (homotopy) [box, right of=flow, xshift=1.5cm] {Section \ref{sec:from-flow-categories}};
\node (virtual) [box, below of =homotopy] { Section \ref{sec:impl-atlas-with}};
    \node (compar) [box, below of =virtual] {Section \ref{sec:comp-coch-kuran}};
    \node (cohI) [box, below of =compar] {Section \ref{sec:local-glob-constr}};
      \node (cohII) [box, below of =cohI] {Section \ref{sec:coher-comp-tang}};

      \node (PSS) [box, left of=kuranishi, xshift=-1.5cm] {Section \ref{sec:hamilt-floer-theory-1}};
      \node (kurHam) [box, below of =PSS] {Section \ref{sec:kuran-flow-categ-1}};
      \node (Hamor) [box, below of =kurHam] {Section \ref{sec:stable-almost-compl-1}};
      \node (Morse) [box, below of =Hamor] {Appendix \ref{sec:morse-theory-and-homotopy-types}};

\node (Spectra) [box, right of=homotopy, xshift=1.5cm] {Appendix \ref{sec:groups-categ-spectra}};
      \draw [dashed,->] (Spectra) .. controls  +(up:1cm) and +(up:1cm) .. (flow);

 \node (Morava) [box, right of=compar, xshift=1.5cm] {Appendix \ref{sec:orientations} and \ref{sec:adams-isom-ambid}};
      \draw [dashed,->] (Morava) -- (compar);

  \draw [dashed, ->] (flow) -- (homotopy);
  \draw [->] (homotopy) -- (virtual);
  \draw [->] (flow) -- (kuranishi);
  \draw [->] (kuranishi) -- (virtual);
  \draw [->] (kuranishi) -- (kursmooth);
  \draw [->] (kursmooth) -- (kuror);
  \draw [->] (virtual) -- (compar);
  \draw [->] (compar) -- (cohI);
  \draw [->] (cohI) -- (cohII);
  \draw [->] (kuror) -- (cohII);
  \draw [->]  (cohII) -- (Morse);
   \draw [dashed, ->] (kursmooth) -- (kurHam);
  
  \draw [dashed, ->] (flow) -- (PSS);
  \draw [->] (PSS) -- (kurHam);
  \draw [->] (kurHam) -- (Hamor);
  \draw [->] (kuranishi) -- (kurHam);
  \draw [->] (kuror) -- (Hamor);
  \draw [->, dashed] (Morse) .. controls  +(left:2cm) and +(left:2cm) .. (PSS);
\end{tikzpicture}
  
 \caption{Dependency diagram: dashed arrows indicate minimal dependencies.}
  \label{fig:Dependency}
\end{figure*}
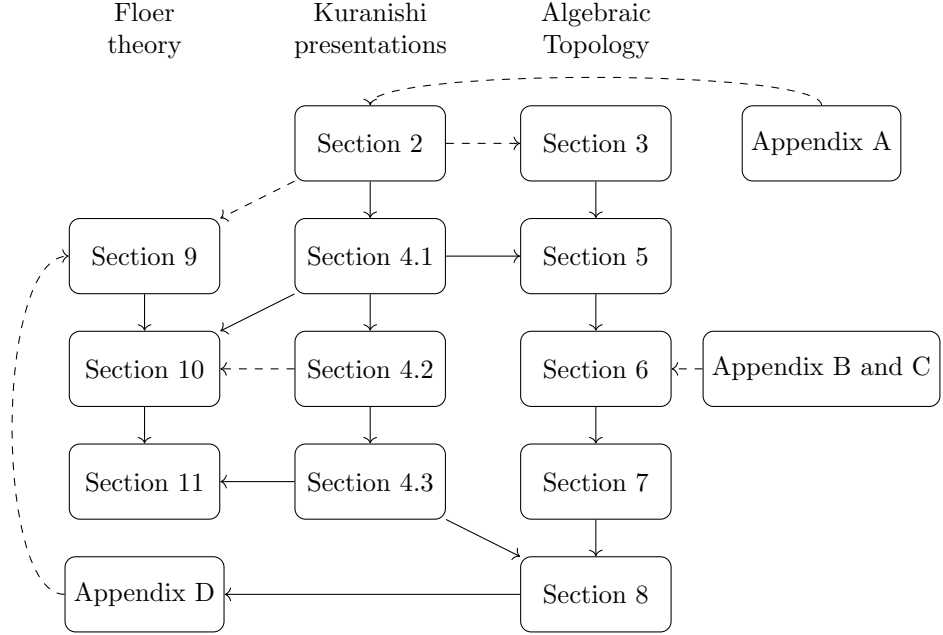

\subsection*{Acknowledgements}
We would like to thank John Berman, Kenji Fukaya, Mike Hill, Mike
Hopkins, Thomas Kragh, Alexander Kupers, Robert Lipshitz, Mike
Mandell, Dusa McDuff, Kate Ponto, Frank Quinn, Dietmar Salamon, Katrin
Wehrheim, and Shmuel Weinberger, for their insights on various issues
related to this project, as well as Tyler Lawson, Michael Miller,
Semon Rezchikov, and Ivan Smith for extensive comments on a
preliminary version of this paper. The first author would particularly
like to thank Jacob Lurie for explaining to him the ambidexterity of
Morava $K$-theory, John Pardon for early discussions about using the
formalism of \cite{Pardon2016} for the purpose of constructing the
Floer homotopy type, as well as Paul Seidel for an early conversation
about potential applications of this paper's main construction.

\part{The homotopy type of flow categories}
\label{part:kuran-pres}

\section{Topological and orbispace flow categories}
\label{sec:topol-orfif-flow}

As observed by Cohen-Jones-Segal \cite{CohenJonesSegal1995}, filtered
homotopy types arise in Floer theory from appropriately oriented flow
categories.  The purpose of this section is to introduce the formalism of such categories, to define the corresponding spectrally enriched categories of \emph{relative cochains,} and to introduce the notions of  \emph{orbispace and Kuranishi flow categories} which will be relevant to Floer theoretic applications.

\subsection{Topologically enriched flow categories}
\label{sec:topol-enrich-categ}

\subsubsection{Equivariant flow categories}
\label{sec:equiv-flow-categ}

Let $\cP$ be a partially ordered set which is \emph{locally finite
dimensional} in the sense that, for any pair of comparable elements
$p < q$, there is a bound on the number of elements of a totally
ordered subset with minimum $p$ and maximum $q$; this technical condition will not be relevant in this section, but will be essential in Section \ref{sec:from-flow-categories}, 
\begin{rem}
In Floer theory, the set $\cP$ will be the set of lifts of Hamiltonian
orbits to the universal cover of the free loop space, or to an
intermediate cover on which the action functional is
well-defined. This set is naturally ordered by action, but this
ordering need not be locally finite dimensional.  Nonetheless, Gromov
compactness implies that one may in some sense discretize this
ordering so that local finiteness holds, allowing the constructions of
this section to directly apply.
\end{rem}

We treat $\cP$ as a category in the usual way, with objects the
elements of $\cP$ and morphisms given by an arrow from $p$ to $q$ if
and only if $p \leq q$. We shall also fix a discrete group $\Pi$ acting freely on $\cP$ as a
category; explicitly, $\Pi$ has a free action on the elements of $\cP$
that is compatible with the order in the sense that $p < q$ then
$\pi(p) < \pi(q)$ for $\pi \in \Pi$.

We will extend the action of
$\Pi$ to enriched categories over $\cP$; we are primarily interested
in topologically or spectrally enriched categories.  
In order to specify actions of groups on categories, we use the
language of $2$-categories, i.e., categories equipped with a category
of morphisms for each pair of objects.  Associated to $\Pi$ is a
$2$-category $B\Pi$ with a single object, $1$-cells the elements of
$\Pi$, and $2$-cells encoding the product on $\Pi$.  An action of
$\Pi$ on a category is now some kind of $2$-functor from $B\Pi$ to the
$2$-category of enriched categories, enriched functors, and enriched
natural transformations.  Explicitly, the data of the action on an
enriched category $\aC$ consists of:
\begin{enumerate}
\item For each $\pi \in \Pi$ an enriched functor
\begin{equation}
\gamma_{\pi} \colon \aC \to \aC,
\end{equation}
\item for each pair $\pi_1, \pi_2 \in \Pi$ a natural comparison
transformation 
\begin{equation}\label{eqn:compose}
\gamma_{\pi_1} \circ \gamma_{\pi_2} \to \gamma_{\pi_2 \pi_1},
\end{equation}  
\end{enumerate}
and this data satisfies unitality and associativity properties that we
describe in more detail in Appendix~\ref{sec:groupact}.  When the
comparison transformation (Equation~\eqref{eqn:compose}) is the
identity, we refer to this as a strict action.  A strict action can
also be described as a functor from $G$ (regarded as a category with a
single object) to the category of categories.  Many of the actions we
consider are strict.  However, we will also consider pseudo-actions
(where the transformation is an isomorphism) and homotopical actions
(where the transformation is merely a weak equivalence).

\begin{defin} \label{def:top_flow_category}
A \emph{$\Pi$-equivariant flow category over $\cP$} is a topologically
enriched category with object set $\cP$, equipped with a strict action
of $\Pi$ extending the action on the set of objects, and such that
morphisms from $p$ to $q$ are given by a point if $p=q$, and are otherwise empty
unless $p < q$.
\end{defin}

We note at this stage that there is tautological example of such a
flow category, given by considering $\cP$ as a topological category
with discrete spaces of morphisms; this will be the key example to
consider in Section \ref{sec:from-flow-categories}. 

\begin{rem}
The example to keep in mind for Floer theory is the following: on a symplectic manifold equipped with a non-degenerate Hamiltonian we consider the set $\cP$ of lifts of time-$1$ orbits to the universal cover of the free loop space (or an intermediate cover on which the action functional is singly valued), with partial order induced by action. The group $\Pi$ is then given by the deck transformations of this cover, and the morphisms in the flow category are given by the Gromov-Floer compactification of the moduli spaces of Floer trajectories (see Section \ref{sec:hamilt-floer-theory-1}).
\end{rem}

\begin{rem} \label{rem:pathological_spaces}
We warn the reader that we do not assume that the morphism spaces are
CW complexes or Euclidean Neighborhood Retracts (ENRs), nor that they
have the homotopy types of such spaces. Instead, the main properties
implied by the geometric constructions we consider is that the
morphism spaces are compact and are locally homeomorphic to closed
subsets of Euclidean space.  In particular, they are Hausdorff.
\end{rem}

\begin{example} The simplest non-trivial example arises in the case $\cP$ is 
  a totally ordered set consisting of a triple of elements $ \{ p < q < r\} $ and
the group $\Pi$ is trivial.  Consider the data of topological spaces
  $\cM(p,q)$, $\cM(q,r)$ and $\cM(p,r)$, and a closed embedding 
  \begin{equation}
     \cM(p,q) \times \cM(q,r) \to \cM(p,r),
  \end{equation}
whose image will be denoted $\partial \cM(p,r)$.  This composition
  map suffices to specify a topologically enriched category with
  object set $\cP$, where we define 
\begin{equation}
\cM(p,p) = \cM(q,q) = \cM(r,r) = \ast
\end{equation}
and composition maps with $\ast$ to be the identity.
\end{example}

\subsubsection{Collared completion of flow categories}
\label{sec:coll-compl-flow}

In Section~\ref{sec:relat-coch-orbisp}, we shall consider a certain
spectrally enriched category associated to a flow category, whose
homotopy groups compute the relative cohomology of morphism spaces. We
will use a lift to spectra of the boundary homomorphism in cohomology
with coefficients in a spectrum $\bbk$ 
\begin{equation}
 H^*(Y; \bbk) \to H^{*+1}(X,Y; \bbk) = H^*(X, Y; \Omega \bbk) 
\end{equation}
associated to an inclusion $Y \subset X$. In order for such a map to
be functorial, it is convenient for $Y$ to be equipped with a collar,
corresponding to the variable appearing in the delooping of the right
hand side. Our purpose in this section is to replace each flow
category by a collared completion, which will be essential for the
functoriality of later constructions. 

\begin{defin}
For each pair $p,q \in \cP$, we define the partially ordered set
$2^{\cP}(p,q)$ to consist of totally ordered subsets of $\cP$ all of
whose elements lie strictly between $p$ and $q$, with ordering given
by inclusion.  
\end{defin}

We assign to each element $Q \in 2^{\cP}(p,q)$ the cube
\begin{equation}
    \kappa^Q \equiv \prod_{i \in Q} \kappa^i
\end{equation}
of dimension $|Q|$, where each $\kappa^i$ is an 
interval $[0,1]$. 
We adopt the
convention that $\kappa^\emptyset= \ast$ when $Q$ is empty.
Given an inclusion $Q \subset P$ of elements of $ 2^{\cP}(p,q)$, we
have a natural map 
\begin{equation}
    \kappa^Q \to \kappa^P
\end{equation}
associated to setting all coordinates not in $Q$ to equal $0$; this
construction defines a functor from $2^{\cP}(p,q)$ to the category of
topological spaces.

If $\cM$ is a flow category, we can think of each element $Q = (q_1,
\ldots, q_n)$ of $2^{\cP}(p,q)$ as a composable sequence in $\cP$, so
we associate to it the space
\begin{equation} \label{eq:space_associated_to-collar}
  \cM(Q) \equiv \cM(p,q_1) \times \cdots \times \cM(q_n, q).
\end{equation}
This construction is contravariantly functorial in $Q$, in the sense
that an inclusion $Q \to P$ of elements of $2^{\cP}(p,q)$ induces a
natural map $\cM(P) \to \cM(Q)$ by composition in $\cM$.

Since the cube $\kappa^Q$ on $Q$ is covariantly functorial in $Q$, we
can define the \emph{collared completion} of $\cM(p,q)$ to be the union of
the values of the functor $\cM(-) \times \kappa^{(-)}$ quotiented by
the equivalence relation that glues these spaces along the maps
\begin{equation}
  \cM(Q) \times \kappa^Q \leftarrow \cM(P) \times \kappa^Q \to \cM(P)
  \times \kappa^P
\end{equation}
for each inclusion $Q \to P$:
\begin{defin}
The collared completion of $\cM(p,q)$ is the coend
\begin{equation}
\hat{\cM}(p,q) \equiv   \int^{Q \in 2^{\cP}(p,q)}  \cM(Q) \times \kappa^Q,
\end{equation}
of the functors $\cM$ and $\kappa$ on $2^{\cP}(p,q)$.
\end{defin}

\begin{rem}
The space $\hat{\cM}(p,q)$ is homeomorphic to the homotopy colimit of
the functor $ 2^{\cP}(p,q) \to \Top$, which assigns $\cM(Q)$ to
$Q$. That construction would naturally yield a decomposition into
simplices, whereas later constructions will rely on the cubical
decomposition that we highlight.
\end{rem}
 The collared completion is naturally equipped with a projection
 map
 \begin{equation}
\hat{\cM}(p,q) \to \cM(p,q)      
 \end{equation}
which is induced by the collapse maps $\kappa^Q \to \ast$. It also has a natural notion of boundary:
\begin{defin}\label{defn:collared-completion}
The \emph{boundary} $\partial \hat{\cM}(p,q)$  of $\hat{\cM}(p,q) $ is the
subset where one or more collar
coordinate equals $1$. 
\end{defin}
Observe that, by construction, the inclusion $\partial
\hat{\cM}(p,q) \to \hat{\cM}(p,q)$ is a Hurewicz cofibration (i.e. a
map satisfying the homotopy extension property, which in particular
implies that it is a closed inclusion); it is in fact locally modeled
after the canonical example of such a cofibration, which is the
inclusion of a space in its product with the interval $[0,1)$.

\begin{rem}
Our use of the term boundary is justified by considering the case
studied by Cohen, Jones, and Segal~\cite{CohenJonesSegal1995}: if each
space $\cM(p,q)$ is a smooth manifold with corners, such that the
composition maps  
\begin{equation}
\cM(p,q) \times \cM(q,r) \to \cM(p,r)
\end{equation}
enumerate the codimension-$1$ boundary strata, then $\hat{\cM}(p,q)$
is again a manifold with corners, with boundary given as above.
\end{rem}

\begin{figure}[h]
  \centering
 \begin{subfigure}[c]{0.49\textwidth}
   \centering
   \begin{tikzpicture}
      \node (p) at (-1,0) {$p$};
      \node (q1) at (0,1) {$q_1$} ;
      \node (q2) at (0, .5) {$q_2$};
      \node (vdots) at (0,0) {$\vdots$};
      \node  (qk) at (0,-1) {$q_k$};
      \node  (q) at (1,0) {$q$};
      \draw[->] (p) -- (q1);
      \draw[->] (q1) -- (q);
      \draw[->] (p) -- (q2);
      \draw[->] (q2) -- (q);
      \draw[->] (p) -- (qk);
      \draw[->] (qk) -- (q);
    \end{tikzpicture}
  \end{subfigure}
\begin{subfigure}[c]{0.49\textwidth}
  \centering
  \begin{tikzpicture}
      \node (p) at (0,0) {$p$};
      \node (q1) at (1,0) {$q_1$} ;
      \node (q2) at (2,0) {$q_2$};
      \node  (q) at (3,0) {$q$};
      \draw[->] (p) -- (q1);
      \draw[->] (q1) -- (q2);
      \draw[->] (q2) -- (q);
    \end{tikzpicture}
  \end{subfigure}
 \begin{subfigure}[c]{0.49\textwidth}
   \centering
   \begin{tikzpicture}
     \coordinate (p) at (0,0);
     \coordinate (q1) at (0,1);
     \coordinate (q2) at (60:1);
     \coordinate (dots) at (30:1);
      \coordinate (qk) at (1,0);
     \filldraw (p) circle (2pt) node[align=left,   left] {$\emptyset$};
     \filldraw (q1) circle (2pt) node[align=left,   above] {$q_1$};
     \filldraw (q2) circle (2pt) node[align=left,   right] {$q_2$};
     \filldraw (dots) circle (2pt) node[align=left,   right] {$\cdots$};
      \filldraw (qk) circle (2pt) node[align=left,   right] {$q_k$};
      \draw (p) -- (q1);
      \draw (p) -- (q2);
      \draw (p) -- (dots);
      \draw (p) -- (qk);
    \end{tikzpicture}
  \end{subfigure}
\begin{subfigure}[c]{0.49\textwidth}
  \centering
  \begin{tikzpicture}
 \coordinate (p) at (0,0);
     \coordinate (q1) at (0,1);
     \coordinate (q2) at (1,0);
     \coordinate (q1q2) at (1,1);
     \filldraw (p) circle (2pt) node[align=left,   left] {$\emptyset$};
     \filldraw (q1) circle (2pt) node[align=left,   above] {$q_1$};
     \filldraw (q2) circle (2pt) node[align=left,   right] {$q_2$};
     \filldraw (q1q2) circle (2pt) node[align=left,   right] {$q_1q_2$};
     \draw (p) -- (q1) -- (q1q2) -- (q2) -- cycle;
    \end{tikzpicture}
   \end{subfigure}
  \caption{The cubical complex of morphisms from $p$ to $q$ associated to small posets.}  \label{fig:cubical-sets}
\end{figure}
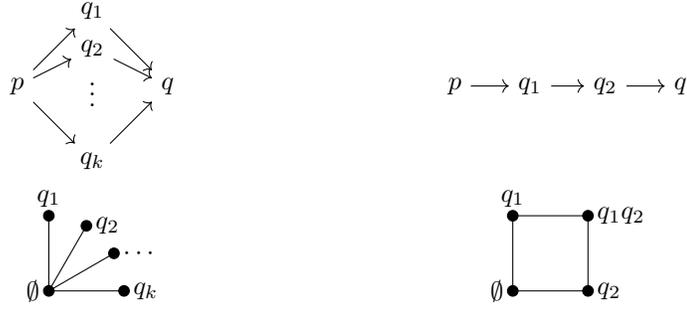

\begin{example}
In the case $\cM = \cP$, it is straightforward to compute that
\begin{equation}
\hat{\cP}(p,q) \equiv \colim_{Q \in 2^{\cP}(p,q) } \kappa^Q
\end{equation}
is a cubical complex obtained by gluing along these inclusions. By
the assumption that $\cP$ is locally finite dimensional, $\hat{\cP}(p,q)$
is a finite dimensional cubical complex, which deformation retracts to
the inclusion of the point associated to $Q = \emptyset$.  We consider
the following special cases: 
\begin{enumerate}
\item If $p = q$, $\hat{\cP}(p,p)$ is a point and $\partial \hat{\cP}(p,p) = \emptyset$.
\item If $p$ and $q$ are successive elements, 
\begin{equation}
\hat{\cP}(p,q) = \ast \quad\textrm{and}\quad \partial \hat{\cP}(p,q) =
\emptyset.
\end{equation}
\item If there is a unique element $q_1$ such that $p < q_1 < q$, then
\begin{equation}
\hat{\cP}(p,q) = [0,1] \quad\textrm{and}\quad \partial \hat{\cP}(p,q)
= \{ 1\},
\end{equation}
 with the endpoint $0$ the image of the inclusion of
 $\kappa^{\emptyset}$, and the endpoint $1$ the image of the inclusion
 of $\partial \hat{\cP}(p,q)$. 
\item If there are two incomparable elements between $p$ and $q$,
  i.e., $p < q_1 < q$ and $p < q_2 < q$, then 
\begin{equation}
\hat{\cP}(p,q) = [0,1] \vee_{0} [0,1] \quad \textrm{and} \quad 
 \partial \hat{\cP}(p,q) = \{1\} \cup \{1\},  
\end{equation}
where the union means the disjoint union of the endpoints.
More generally, as shown on the left side of Figure
\ref{fig:cubical-sets}, if there are $k$ incomparable 
elements $\{q_i\}$ between $p$ and $q$, then 
\begin{equation}
\hat{\cP}(p,q) = \underbrace{[0,1] \vee_0 [0,1] \vee_0 \ldots \vee_0
  [0,1]}_{k} \quad\textrm{and}\quad 
 \partial \hat{\cP}(p,q) = \bigcup_k \{1\}
\end{equation}
\item Given the poset $p < q_1 < q_2 < q$, then 
\begin{equation}
\begin{gathered}
\hat{\cP}(p,q) \cong [0,1] \times [0,1] \\ 
\partial \hat{\cP}(p,q) \cong \{(1,t) \, \mid \, t \in [0,1]\} \vee
\{(s,1) \, \mid \, s \in [0,1]\}
\end{gathered}
\end{equation}
 with basepoint $(1,1)$ (see the right side of Figure \ref{fig:cubical-sets}).  More generally, if $\cP$ is
  a totally ordered finite poset with $k$ elements (and bottom element
  $p$ and top element $q$), then 
\begin{equation}
\hat{\cP}(p,q) = [0,1]^{|\cP|-2}
\end{equation}
with boundary the faces where at least one coordinate is $1$.

\item Given the poset $p < q_1 < q < z$ and $p < q_2 < q$, we have
  that $\hat{\cP}(p,q)$ consists of the product $[0,1] \times [0,1]$
  corresponding to the totally ordered subset $p < q_1 < q < z$ and
  the product $[0,1] \times [0,1]$ corresponding to the totally
  ordered subset $p < q_2 < q < z$, glued along the face corresponding
  to $q$.  The boundary is the set of edges where the coordinate
  associated to either $q$, $q_2$, or $q_1$ is $1$.

\end{enumerate}
\end{example}

Given totally ordered subsets $P \in 2^{\cP}(p,q)$ and $Q \subset 2^{\cP}(q,r)$, there is a natural map
\begin{equation}
\theta_{p,q,r} \colon  \kappa^P\times \kappa^Q \to \kappa^{P \amalg \{ q \} \amalg Q}
\end{equation}
specified by the inclusion on coordinates in $P$ and $Q$ and setting
the coordinate $q$ to $1$.  These maps are associative:

\begin{lem}
Given a well a totally ordered subset $R \subset  2^{\cP}(r,s) $, the
following diagram commutes: 
\begin{equation}
  \begin{tikzcd}
\kappa^P \times \kappa^Q \times \kappa^R \ar[rr,"\theta_{p,q,r} \times \id"] \ar[d,swap,"\id \times \theta_{q,r,s}"] && \kappa^{P \amalg \{q \} \amalg Q} \times \kappa^R \ar[d,"\theta_{p,r,s}"] \\
\kappa^P \times \kappa^{Q \amalg \{r\} \amalg R} \ar[rr,swap,"\theta_{p,q,s}"] && \kappa^{P \amalg \{q \} \amalg Q \amalg \{r\} \amalg R}.
  \end{tikzcd}  
\end{equation} \qed
\end{lem}

Furthermore, the map $\theta_{p,q,r}$ is compatible with inclusions of subsets in $2^{\cP}(p,q)$ and $2^{\cP}(q,r)$. Combining this with the natural homeomorphism
\begin{equation}
  \cM(P) \times\cM(Q) \to  \cM(P \cup  \{ q \} \cup Q)
\end{equation}
for $P \subset 2^{\cP}(p,q)$ and $Q \subset 2^{\cP}(q,r)$, which is contravariantly functorial in both variables, we obtain an associative product map 
\begin{equation} \label{eq:product_cubical_collar_category}
\hat{\cM}(p,q) \times \hat{\cM}(q,r) \to \hat{\cM}(p,r),
\end{equation}
whose image is contained in $\partial  \hat{\cM}(p,r)$. By construction, the space of morphisms $ \hat{\cM}(p,p) $ is a single point, which acts as the unit.

\begin{defin}
The \emph{collared completion of $\cM$} is the topological category
$\hat{\cM}$ with object set $\cP$ and morphism spaces from $p$ to $r$
given by  $\hat{\cM}(p,r)$. 
\end{defin}

Applying this construction to $\cM = \cP$ yields the \emph{collar category} $\hat{\cP}$ with morphisms the cubical complexes $\kappa(p,q)$ (where we take
the unique point in $\kappa(p,p)$ to be the identity map), and the
associative and unital composition specified by Equation
\eqref{eq:product_cubical_collar_category}.
\begin{rem}
The collared category $\hat{\cP}$ is constructed to have the universal
property that an enriched functor from $\hat{\cP}$ to spaces or
spectra is the same data as a homotopy coherent diagram over $\cP$.
More generally, the collared completion produces a fattened indexing
category that allows us to encode coherent composition homotopies
explicitly.  See~\cite{Vogt1973} and~\cite{CordierPorter1997} for a
early treatments of homotopy coherent category theory,
and~\cite{Leitch1974} for an appearance of this particular
construction (and also~\cite{Segal1974}).
\end{rem}

\begin{lem}
The collared completion $\hat{\cM}$ is equipped with an enriched 
projection functor $\hat{\cM} \to \hat{\cP}$ such that the following
diagram commutes 
\begin{equation} \label{eq:forget_collars_project_to_P}
  \begin{tikzcd}
    \hat{\cM} \ar[d] \ar[r] & \cM \ar[d] \\
    \hat{\cP} \ar[r] & \cP,
  \end{tikzcd}
\end{equation}
where the horizontal maps are induced by the projection maps that
collapse the collar directions. \qed
\end{lem}

Since the construction of collars did not involve any choices, the
group $\Pi$ has a natural induced action on the collared completions.
More precisely, for $\pi \in \Pi$, let $\pi(Q)$
denote the evident subset of $2^{\cP}(\pi(p),\pi(q))$.  Then $\pi$ takes
$\kappa^P$ to $\kappa^{\pi(P)}$ (where the action on the coordinates
is trivial) and the inclusion $\kappa^P \to \kappa^Q$ induces an
inclusion $\kappa^{\pi(P)} \to \kappa^{\pi(Q)}$.  As a consequence,
$\pi$ acts on morphisms via a natural map
\begin{equation}
  \hat{\cM}(p,q) \to
\hat{\cM}(\pi(p),\pi(q))
\end{equation}
which is the identity.
Moreover, $\pi$ is evidently compatible with
the composition maps $\theta$, i.e.,  
\begin{equation}
\pi(\theta_{p,q,r}(f,g)) = \theta_{\pi(p), \pi(q),\pi(r)}(\pi(f), \pi(g)),
\end{equation}
and clearly preserves the unit.  This discussion establishes that
$\hat{\cM}$ is a $\Pi$-equivariant flow category.

Moreover, it is clear from the construction that the functors in
Diagram~\eqref{eq:forget_collars_project_to_P} are strictly
$\Pi$-equivariant.  Here a $\Pi$-equivariant functor between
equivariant categories $\aC$ and $\aD$ is a $2$-natural transformation
between the $2$-functors $B\Pi \to \aC$ and $B\Pi \to \aD$.
Explicitly, this amounts to an enriched functor $F \colon \aC \to \aD$
and natural transformations $\gamma_\pi \circ F \to F \circ
\gamma_\pi$ (which in this case are through identities) that satisfy
compatibility diagrams listed in Appendix~\ref{sec:equifunc}.

Summarizing, we have:
\begin{lemma}
The collared completion $\hat{\cM}$ of a $\Pi$-equivariant flow
category $\cM$ is a $\Pi$-equivariant flow category, and all the
arrows in Diagram \eqref{eq:forget_collars_project_to_P} are strictly
$\Pi$-equivariant. \qed
\end{lemma}

\subsection{Relative cochains of topological flow categories}
\label{sec:relat-coch-orbisp}

Given a topological flow category $\cM$, we can construct a spectral
category (i.e., a category enriched in spectra) by applying the
functor $\Sigma^{\infty}_+$ that adds a disjoint basepoint and takes
the suspension spectrum to the mapping spaces.  It will turn out that,
for the purpose of associating homotopy types to flow categories, we
need to consider instead a category of \emph{relative cochains} constructed in
this section.

\subsubsection{Background on spectra}
\label{sec:background-spectra}

Throughout this section (and the paper), we take as our point-set
models of equivariant spectra the closed symmetric monoidal categories
of orthogonal $G$-spectra for various finite groups $G$ (although only
the case of trivial $G$ is used in this section).  See
Appendix~\ref{sec:backgr-categ-g} for a rapid review of the
foundations of the theory of equivariant orthogonal spectra.

We will write $\Sp_G$ to denote the category of orthogonal
$G$-spectra, which is symmetric monoidal with respect to the smash
product $\sma$ and has unit the sphere spectrum $\bS$.  We will denote
by $F(X,Y)$ the internal mapping object, i.e., the spectrum of maps
from $X$ to $Y$.  We will frequently make use of the natural smash
product map 
\begin{equation}
  \label{eq:smash_product_map_functions}
F(X,Y) \sma F(W,Z) \to F(X \sma W, Y \sma Z)  
\end{equation}
which is the adjoint of the smash of the evaluation maps.  This map is
associative and unital; see
Proposition~\ref{prop:associative-smash-function-spectra}.  Note
however that this map is not usually a weak equivalence except in the
presence of dualizability hypotheses (e.g., see~\cite[\S
  III]{LewisMaySteinberger1986} for a comprehensive treatment of
formal duality theory).

\begin{defin}
For a space $Z$ and an orthogonal spectrum $X$, the \emph{spectrum of cochains of $Z$ with coefficients in $X$} is given by
\begin{equation} \label{eq:cochain_notation}
C^*(Z; X) = F(\Sigma^{\infty} Z_+, X^{\mfib}),
\end{equation}
where $\Sigma^{\infty}$ denotes a particular strong monoidal model for
the suspension spectrum functor (given by the functor $F_{0}(-) \colon
\GTop \to \Sp_G$ discussed in Appendix~\ref{sec:ortho}) and
$(-)^{\mfib}$ is the lax monoidal fibrant replacement functor on
orthogonal spectra described in Appendix~\ref{sec:fibrant}.
\end{defin}

We will need the compactly-supported version of the cochains.

\begin{defin}
The \emph{spectrum of compactly supported cochains on $Z$ with coefficients in $X$} is  the orthogonal
spectrum 
  \begin{equation}\label{eq:compactly-supported-maps}
C^{*,c}(Z; X) = F(\Sigma^{\infty} Z^+, X^{\mfib}),
\end{equation} 
where $Z^+$ denotes the one-point compactification.  Note that our assumptions on
$Z$ suffice to ensure that the one-point compactification is
Hausdorff.
\end{defin}

We also have relative versions.

\begin{defin}
For a Hurewicz cofibration $Y
\to Z$, the \emph{spectrum of relative cochains of the pair $(Y,Z)$ with coefficients in $X$} is the orthogonal spectrum
\begin{equation}\label{eq:relative-cochains}
C^*(Z,Y;X) = F(\Sigma^{\infty} Z/Y, X^{\mfib}),
\end{equation}
and $C^{*,c}(Z,Y;X)$ is the compactly
supported analogue. 
\end{defin}

With the familiar convention in algebraic topology that the quotient
of a space by the empty set corresponds to adding a disjoint
basepoint, the absolute case in Equations \eqref{eq:cochain_notation}
and ~\eqref{eq:compactly-supported-maps} corresponds to the case $Y =
\emptyset$. 

Since $(-)^{\mfib}$ is lax monoidal, we have a natural smash product
map 
\begin{equation}
C^*(Z_1; X_1) \sma C^*(Z_2; X_2) \to C^*(Z_1 \sma Z_2;
X_1 \sma X_2)
\end{equation}
induced from the smash product map on mapping spectra; this is
associative and unital.  There are analogous maps on $C^{*,c}(-;-)$
and the relative versions, using the canonical isomorphism $X_1^+ \sma
X_2^+ \cong (X_1 \times X_2)^+$ (which requires that the spaces
involved be locally compact Hausdorff).

The fibrant replacement is required to ensure that our cochains have
the right homotopy type; we now explain what we mean by this.  As
alluded to in Remark \ref{rem:pathological_spaces}, one subtlety that
arises is that we cannot in general assume that $Z$ has the homotopy
type of a CW complex, as we shall assume only that it is locally compact, paracompact, 
and Hausdorff.  Therefore, we are not computing the usual 
derived mapping space, as we do not expect $C^*(-;-)$ to preserve
weak equivalences in the first variable, although it does preserve
homotopy equivalences.  

We are however correctly computing the cohomology.  Recall that 
 by the main theorem of~\cite{Huber1961}, the \v{C}ech 
cohomology groups can be computed using the space of maps from a
$k$-space into Eilenberg-Mac Lane spaces.  As a consequence, for a
fibrant model of $HR$ the homotopy type $C^*(Z; HR)$ has homotopy
groups that recover \v{C}ech cohomology.
More generally, analogous results hold for $C^*(Z; X)$ for any fibrant
orthogonal spectrum $X$, where  we have in mind~\cite{Brown73,
  Jardine97} for \v{C}ech theory in the context of generalized
cohomology; e.g., see~\cite[7.1.0.1]{Lurie2009} (note that in our context the \v{C}ech and hypercover
localizations agree, as the covering dimension of the spaces $Z$ that
arise is finite; e.g., see ~\cite[7.2.1.12]{Lurie2009}).

\begin{rem}
If we were studying the \emph{homology} of Kuranishi presentations,
the potential pathologies of the spaces we consider would force us to
consider the associated pro-spaces constructed either from nerves of
\v{C}ech covers, or neighbourhoods in the ambient manifold (this is
the subject of Shape theory~\cite{Mardevsic2000}),
see~\cite[Appendix A.9]{Pardon2016} for the corresponding discussion
at the level of ordinary homology. 
\end{rem}

\subsubsection{The relative cochains of a flow category}
\label{sec:relat-coch-coll}

\begin{notn}
The notion of relative cochains will require the use of certain
\emph{desuspension functors.}  It is convenient to label the spheres
appearing in these desuspensions by elements $q$ of the partially
ordered set $\cP$, so we begin by introducing a real line
\begin{equation}
    \ell_q \equiv \bR^{\{[q]\}},
\end{equation}
where we recall that $[q]$ is the equivalence class of $q$ under the
action of $\Pi$. By construction, this choice is $\Pi$-invariant in
the sense that $\ell_q = \ell_{\pi \cdot q}$ for each $\pi \in \Pi$. 

The main advantage of working with the line $\ell_q$ rather than the
real line $\bR$ is that it breaks the symmetry on the direct sum
$\ell_q \oplus \ell_p $, when considering more than one orbit, which
makes is possible to unambiguously write down maps on this direct sum
that depend on its decomposition into factors. 

We write $S^{\ell_{q} }$ for the associated
copy of $S^1$ obtained by one-point compactification, and, given a
spectrum $\bbk$, we write  
\begin{equation}
  \Omega^{\ell_q} \bbk \equiv F(\Sigma^\infty S^{\ell_{q} }, (\bbk)^{\mfib})
\end{equation}
 for the desuspension.
\end{notn}

Let $\cM$ be a $\Pi$-equivariant flow category over $\cP$, and $\bbk$ a cofibrant ring spectrum:
\begin{defin}
For each pair of distinct
elements $(p,q) \in \cP$,
the \emph{relative cochains of $\cM$ with coefficients in $\bbk$} are the $\Omega^{\ell_{q}}
\bbk$-valued cochains on $\hat{\cM}(p,q)$ relative to $\partial
\hat{\cM}(p,q)$, i.e., 
\begin{equation}
 C^*_{\rel \partial} ({\cM}; \Omega \bbk)(p,q) \equiv   C^* (\hat{\cM}(p,q), \partial \hat{\cM}(p,q);\Omega^{\ell_{q}} \bbk).
\end{equation}

In the special case $p = q$, we define
\begin{equation}
C^*_{\rel \partial}({\cM}; \Omega \bbk)(p,p) \equiv  \bbk.
\end{equation}
\end{defin}
\begin{rem}
In Floer and Morse theory, it is in some sense more natural to introduce a line $\ell_{pq}$ associated to a pair of objects, which corresponds to translation in the moduli space of flow lines connecting them. The group of translation is canonically identified with $\bR$, so that the labelling by orbits is again only a matter of convenience to record correspondence between factors. Our notation breaks the symmetry between input and output, and formally associates this line of translations to the output. One can proceed without breaking symmetry, at the cost of replacing the composition formulae which we will presently define by more complicated ones.
\end{rem}

We have observed earlier that the inclusion $\partial \hat{\cM}(p,q)
\to \hat{\cM}(p,q)$ is always a Hurewicz cofibration, so that the
relative cochains $C^*_{\rel \partial} ({\cM}; \Omega \bbk)(p,q)$
compute the relative (shifted) cohomology groups.

\begin{rem}
We may identify the spaces in the spectrum $C^*_{\rel
  \partial} ({\cM}; \Omega \bbk)(p,q)$ as follows.  By adjunction, for
any based space $Z$ we have natural equivalences
\begin{equation}
F(\Sigma^{\infty} Z, (\Omega^{\ell_{q}} \bbk)^{\mfib}) \htp
F(\Sigma^{\infty} Z, \Omega^{\ell_{q}} \bbk^{\mfib}) \cong
F(\Sigma^{\infty} S^{\ell_{q}} \sma Z, \bbk^{\mfib})   
\end{equation}
of orthogonal spectra, and so for a finite-dimensional representation
$V$ we have an equivalence
\begin{multline} 
\left(C^*_{\rel \partial}({\cM}; \Omega \bbk)(p,q)\right)(V) \htp \\
\Map\left(S^{\ell_{q}} \sma \hat{\cM}(p,q)/ \partial \hat{\cM}(p,q), \Omega^{\infty} \Omega^V \bbk^{\mfib}\right),  
\end{multline}
i.e., we may identify the value of this spectrum at $V$ with maps whose domain is $\hat{\cM}(p,q) \times
S^{\ell_{q}}$, subject to the condition that the boundary in the
first factor and the basepoint in the second factor map to the
basepoint in the target.
\end{rem}

\begin{example} \label{rem:acyclic_morphism_spaces}
Let $\cP$ be the natural numbers with their usual ordering, and let
$\bbk$ be the sphere spectrum for specificity. In that case,
$\hat{\cP}(p,q) = [0,1]^{q-p-1}$, and the inclusion of $\partial
\hat{\cP}(p,q) \subset \hat{\cP}(p,q) $ is a homotopy equivalence
unless $p$ and $q$ are successive elements, in which case the boundary
is empty.  This implies that the spectrum $C^*_{\rel \partial}
(\hat{\cP}; \Omega \bS)(p,q) $ is contractible (i.e., has
\emph{trivial} homotopy groups) except for successive elements, where
it is equivalent to $\Omega \bS$.  The reason for having non-trivial
models for acyclic spectra is that we shall need the composition maps,
which are about to be defined, to encode specified null homotopies.
\end{example}

Given a triple $(p,q,r)$, we now define a composition map
\begin{equation} \label{eq:product_relative_cochains_partial_ordered_set}
\Psi_{p,q,r} \colon C^*_{\rel \partial}( {\cM}; \Omega \bbk)(p,q) \sma C^*_{\rel \partial}( {\cM}; \Omega \bbk)(q,r) \to   C^*_{\rel \partial}( {\cM}; \Omega \bbk)(p,r).
\end{equation}
The starting point is to fix a homeomorphism of the interior of the collar $\mathring{\kappa} = (0,1)$ with $\bR$, and obtain an identification
\begin{equation}\label{eq:map_collar_lines-coordinates}
\mathring{\kappa}^q \cong \ell_{q} .  
\end{equation}

Given a pair $P$ and $Q$ of elements of $2^{\cP}(p,q)$ and
$2^{\cP}(q,r)$, the projection away from the $q$ coordinate
induces a map
\begin{equation}
\kappa^{Q \amalg \{q\} \amalg P} \to \kappa^Q \times \kappa^{P},
\end{equation}
and the inverse to the homeomorphism of
Equation~\eqref{eq:map_collar_lines-coordinates} induces a map
\begin{equation}
(\kappa^q)_+ \sma S^{\ell_{r}} \to S^{\ell_{q} \oplus \ell_{r}}.
\end{equation}
Combining these, we define the map
\begin{equation} \label{eq:map_collar_lines}
\Delta_{P,Q} \colon \kappa^{Q  \amalg \{q\} \amalg  P}_+  \sma
S^{\ell_{r}} \to 
(\kappa^{Q} \times \kappa^{P})_+ \sma \kappa^q_+ \sma S^{\ell_r} 
\to
(\kappa^{Q} \times \kappa^{P})_+ \sma S^{\ell_{q} \oplus \ell_{r}}.
\end{equation}

The key point is that this map takes both endpoints of
the $\kappa^q$ factor in $\kappa^{Q \amalg \{q\} \amalg P} $, to the
basepoint of the sphere $S^{\ell_{q} \oplus \ell_{r}} \cong S^{\ell_{q}}
\wedge S^{\ell_{r}} $.  This means, in particular, that any
basepoint-preserving map from the left-hand side extends (as a
constant map) to any cube labeled by a set that does not contain $q$.

Using the inverse of the homeomorphism 
\begin{equation}
\cM(P) \times\cM(Q) \to  \cM(P \cup  \{ q \} \cup Q),
\end{equation}
we have a map $\Delta_{P,Q}$
\begin{equation} \label{eq:map_collar_lines_with_cat}
\begin{tikzcd}
(\cM(P \cup  \{ q \} \cup Q) \times \kappa^{Q  \amalg \{q\} \amalg  P})_+  \sma
S^{\ell_{r}} \ar[d] \\ 
(\cM(P) \times\cM(Q) \times \kappa^{Q} \times \kappa^{P})_+ \sma
\kappa^q_+ \sma S^{\ell_r}
\ar[d] \\ 
(\cM(P) \times\cM(Q) \times \kappa^{Q} \times \kappa^{P})_+ \sma
S^{\ell_{q} \oplus \ell_{r}}.
\end{tikzcd}
\end{equation}

Thus, since the maps $\Delta_{P,Q}$ are compatible with inclusions,
they extend over the colimit to give rise to maps 
\begin{multline}
\Delta_{p,q,r} \colon \hat{\cM}(p,r)_+ \sma S^{\ell_{r}} \to
(\hat{\cM}(q,r) \times \hat{\cM}(p,q))_+ \sma S^{\ell_{q} \oplus
  \ell_{r}}  \\ \to  (\hat{\cM}(q,r)_+ \sma S^{\ell_r}) \sma (\hat{\cM}(p,q)_+ \sma
S^{\ell_{q}}).
\end{multline}
When we take $p=q$ or $q=r$, these maps are the identity.  

\begin{lem}\label{lem:Delta-coassociative}
The maps $\Delta_{p,q,r}$ are coassociative, in the sense that for $p
\leq q \leq r \leq s$, the following diagram commutes
\begin{equation}
\begin{tikzcd}
\hat{\cM}(p,s)_+ \sma S^{\ell_{s}} \ar[r,"\Delta_{p,q,s}"] \ar[d,swap,"\Delta_{p,r,s}"] &
(\hat{\cM}(q,s)_+ 
\sma S^{\ell_s})
\sma (\hat{\cM}(p,q)_+ \sma S^{\ell_q}) \ar[d,"\Delta_{q,r,s} \sma \id"] \\
\begin{gathered}
(\hat{\cM}(r,s)_+ \sma S^{\ell_s}) \sma \\ (\hat{\cM}(p,r)_+ \sma
S^{\ell_r})  
\end{gathered}
 \ar[r,swap,"id \sma \Delta_{p,r,s}"] &
 \begin{gathered}
(\hat{\cM}(r,s)_+ \sma S^{\ell_s}) \sma (\hat{\cM}(q,r)_+ \sma
S^{\ell_r}) \\ \sma (\hat{\cM}(p,q)_+ \sma S^{\ell_q}).   
 \end{gathered}
\end{tikzcd}
\end{equation}
\end{lem}

\begin{proof}
It is straightforward to chase elements around these diagrams using
the formulas above.  A key aspect to note is that going down and over
we permute $S^{\ell_s}$ past $S^{\ell_r}$ and then $S^{\ell_r}$ past
$S^{\ell_q}$.  Going the other way, we permute $S^{\ell_s}$ past
$S^{\ell_q}$ and then $S^{\ell_s}$ past $S^{\ell_r}$.
\end{proof}

We now use these maps to make the relative cochains
into a spectral category.

\begin{defin}
For $p \neq q \neq r$, the map $\Psi_{p,q,r}$ is the composition
\begin{equation}\label{eq:comp}
  \begin{tikzcd}
    C^*_{\rel \partial}( {\cM}; \Omega \bbk)(q,r) \ar[d] \sma C^*_{\rel \partial}( {\cM}; \Omega \bbk)(p,q) \\
    C^*_{\rel \partial} (({\cM}(q,r), \partial {\cM}(q,r)) \times ({\cM}(p,q), \partial {\cM}(p,q))  ;\Omega^{\ell_{q} \oplus \ell_{r}} (\bbk \sma \bbk))  \ar[d] \\
    C^*_{\rel \partial}( {\cM}; \Omega \bbk)(p,r) 
  \end{tikzcd}
\end{equation}
where the top arrow is induced by the smash product of mapping spectra
and the bottom map is induced by the dual of the map $\Delta_{p,q,r}$
and the multiplication $\bbk \sma \bbk \to \bbk$.  (In the middle
spectrum of Diagram \eqref{eq:comp}, we are using the standard
notation where for pairs $(X,Y)$ and $(W,Z)$, we write $(X,Y) \times
(W,Z)$ to denote the pair $(X \times W, X \times Z \cup Y \times W)$.)

When $p=q$ or $q=r$, we define $\Psi_{p,p,q}$ and $\Psi_{p,q,q}$ in
terms of the composite 
\begin{equation}
\bbk \sma (\Omega^{\ell_{q}} \bbk)^{\mfib} \to 
\bbk^{\mfib} \sma (\Omega^{\ell_{q}} \bbk)^{\mfib} \to (\bbk \sma
\Omega^{\ell_{q}} \bbk)^{\mfib} \to (\Omega^{\ell_{q}} \bbk)^{\mfib},
\end{equation}
using the canonical homeomorphism $\bbk^{\mfib} \cong
F(\bS,\bbk^{\mfib})$. 
\end{defin}

Since the unit map $\bS \to C^*_{\rel \partial}({\cM}; \Omega
\bbk)(p,p)$ is the unit map $\bS \to \bbk$, we can conclude
that the composition is unital:

\begin{lem}
The diagram
\begin{equation}
\begin{tikzcd}
\bS \sma C^*_{\rel \partial}( {\cM}; \Omega \bbk)(p,q) \ar[r] \ar[dr] & C^*_{\rel \partial}({\cM}; \Omega \bbk)(p,p) \sma C^*_{\rel \partial}( {\cM}; \Omega \bbk)(p,q) \ar[d] \\
& C^*_{\rel \partial}({\cM}; \Omega \bbk)(p,q)
\end{tikzcd}
\end{equation}
commutes, as does the analogous diagram on the other side.
\end{lem}

\begin{proof}
The unitality diagram can be written as 
\begin{equation}\label{eqn:unitality}
\begin{tikzcd}
\bS \sma C^*_{\rel \partial}( {\cM}; \Omega \bbk)(p,q) \ar[r] \ar[dr] & \bbk^{\mfib} \sma C^*_{\rel \partial}( {\cM}; \Omega \bbk)(p,q) \ar[d] \\
& C^*_{\rel \partial}({\cM}; \Omega \bbk)(p,q).
\end{tikzcd}
\end{equation}
Since the map $\bS \sma F(X,\bS) \to F(X, \bS \sma \bS)$ is the unit
$\bS \to F(\bS, \bS)$ followed by the composition pairing $F(\bS, \bS)
\sma F(X,\bS) \to F(\bS \sma X, \bS \sma \bS)$, it is formal that the
diagram 
\begin{equation}
\begin{tikzcd}
\bS \sma F(X, \bS) \ar[r] \ar[dr] & F(X, \bS \sma \bS) \ar[d] \\
& F(X, \bS)
\end{tikzcd}
\end{equation}
commutes (e.g., see the discussion
before~\cite[6.12]{ElmendorfKrizMandellMay1997}).  Since the right hand
vertical map is induced by the product on $\bbk$, the unitality
diagram in Equation~\eqref{eqn:unitality} commutes.
\end{proof}

Moreover, it is straightforward to check that the product maps are
associative, using Lemma~\ref{lem:Delta-coassociative}:
\begin{lem}
Given a quadruple $p \leq q \leq r \leq s$ of elements of $\cP$, the product maps defined above fit in a commutative diagram of spectra
\begin{equation}
\begin{tikzcd}  \begin{gathered}
    C^*_{\rel \partial}( {\cM}; \Omega \bbk)(p,q) \wedge  C^*_{\rel \partial}( {\cM}; \Omega \bbk)(q,r) \\ \wedge  C^*_{\rel \partial}( {\cM}; \Omega \bbk)(r,s)
  \end{gathered}
     \ar[d] \ar[r] &
     \begin{gathered}
C^*_{\rel \partial}( {\cM}; \Omega \bbk)(p,r) \\  \wedge  C^*_{\rel \partial}( {\cM}; \Omega \bbk)(r,s)       
     \end{gathered}
 \ar[d]   \\
    C^*_{\rel \partial}( {\cM}; \Omega \bbk)(p,q) \wedge  C^*_{\rel \partial}( {\cM}; \Omega \bbk)(q,s) \ar[r] & C^*_{\rel \partial}( {\cM}; \Omega \bbk)(p,s). 
\end{tikzcd}
\end{equation} \qed
\end{lem}

For the next definition, recall that a category $\aC$ enriched in
spectra consists of a category $\aC$ equipped with mapping spectra
$\aC(x,y)$ for each pair of objects $x,y$ and suitably associative and
unital composition maps.  This is the stable homotopy theory
generalization of a dg category.  See Appendix~\ref{sec:spectral} for
a precise definition and some of the technical properties we need.

\begin{defin}\label{defn:relcoch}
For a $\Pi$-equivariant flow category $\cM$, the \emph{category of relative
  cochains with coefficients in $\bbk$}, denoted  $C^*_{\rel \partial}({\cM}; \Omega
\bbk)$, is the spectral category with 
\begin{itemize}
\item objects the elements of $\cP$, and 
\item morphism spectra for pairs $(p,q)$ given by $C^*_{\rel \partial}({\cM}; \Omega
\bbk)(p,q)$.  
\end{itemize}
The composition is specified by $\Psi_{p,q,r}$ and the unit is induced
by the unit for $\bbk$.
\end{defin}

Moreover, the only choice made in the construction is the
identification in Equation \eqref{eq:map_collar_lines-coordinates},
which we fix, so that the resulting category acquires a natural action
by the group $\Pi$ that acts on $\cP$.  Specifically, $\Pi$ acts on morphisms via the natural identity
map 
\begin{equation}
 C^*_{\rel \partial}({\cM}; \Omega \bbk)(p,q) \to C^*_{\rel
  \partial}({\cM}; \Omega \bbk)(\pi p, \pi q)  
\end{equation}
for each element $\pi \in \Pi$.  It is straightforward to verify that
this action on the morphism spectra is compatible with composition and
preserves the unit. 

\begin{lemma}\label{lem:strictactrelcoch}
 The category of relative cochains $C^*_{\rel \partial}({\cM};\Omega \bbk)$
 has a strict action of $\Pi$. \qed
\end{lemma}

There are pullback maps
\begin{equation}
C^*_{\rel \partial}({\cP};\Omega \bbk)(p,q) \to C^*_{\rel
  \partial}({\cM};\Omega \bbk)(p,q)
\end{equation}
induced by the projection $\hat{\cM}(p,q) \to \hat{\cP}(p,q)$, which
are strictly compatible with the $\Pi$-action and with composition and
thus induce a strict $\Pi$-equivariant spectral functor.  

\begin{lemma}
There is a strict $\Pi$-equivariant functor
\begin{equation}
C^*_{\rel \partial}({\cP};\Omega \bbk) \to C^*_{\rel
  \partial}({\cM};\Omega \bbk).  
\end{equation} \qed
\end{lemma}

\begin{example}
The basic case of the category $C^*_{\rel \partial}({\cP};\Omega \bS)$
will be particularly important in
Section~\ref{sec:from-flow-categories}.
\begin{enumerate}
\item If $\cP$ consists of a single element $p$, then
\begin{equation}
C^*_{\rel \partial}({\cP}; \Omega \bS)(p,p) = \bS
\end{equation}
and the composition is induced by the natural map
\begin{equation}
\bS \sma \bS \cong \bS
\end{equation}

\item If $\cP$ consists of two elements, then 
\begin{equation}
C^*_{\rel \partial}({\cP}; \Omega \bS)(p,q) \cong (\Omega^{\ell_{q}}
\bS)^{\mfib} \htp \Omega \bS.
\end{equation}
The two compositions are both essentially the identity map.

\item If $\cP = \{ p < q < r \}$, then 
\begin{equation}
\begin{aligned}
C_{\rel \partial}^*({\cP};\Omega \bS)(p,q) \htp \Omega \bS,\, C_{\rel
  \partial}^*({\cP}; \Omega \bS)(q,r) \htp \Omega \bS \\
\textrm{and}\quad C_{\rel \partial}^*({\cP}; \Omega
  \bS)(p,r) \htp \ast.
\end{aligned}
\end{equation}
  In terms of these descriptions, the
  non-identity composition $\Omega \bS \sma \Omega \bS \to \ast$ is
  the terminal map.  (In other words, the actual composition is
  null-homotopic.)
  
\item If $\cP$ consists of elements $p < q_1 < q$ and $p < q_2 < q$
  with $q_1$ and $q_2$ incomparable, then for the pairs $(p,q_1)$,
  $(p,q_2)$, $(q_1, q)$, and $(q_2, q)$ the mapping spectrum is
  equivalent to $\Omega \bS$.  For $(p,q)$, we have that 
\begin{equation}
C^*_{\rel \partial}({\cP}; \Omega \bS)(p,q) \htp F(S^1, \Omega \bS) \htp
\Omega^2 \bS.
\end{equation}
This computation uses the identification
\begin{equation}
\kappa(p,q) / \partial \kappa(p,q) \cong \left([0,1] \vee [0,1]\right)
/ \{\{1\} \times [0,1] \cup [0,1] \times \{1\}\} \cong S^1
\end{equation}
In this case, the composition $\Omega \bS \sma \Omega \bS \to \Omega^2
\bS$ is homotopic to the identity.
\end{enumerate}

\end{example}

\subsection{Orbispace flow categories and the Borel construction}
\label{sec:orbisp-flow-categ}

The formalism of topological flow categories is designed to reflect
structures present in finite-dimensional Morse theory. In the
infinite-dimensional setting, a new geometric phenomenon (bubbling)
causes the moduli spaces of flow lines to potentially admit additional
symmetries, so that they acquire the structure of orbispaces.  The
purpose of this section is to provide a generalization of the flow
category formalism that encodes this situation.

Since the geometry provides us with natural charts, we do not appeal to
the general formalisms of orbispaces recently developed in homotopy theory
\cite{HenriquesGepner2007,Korschgen2018,Schwede2018}, preferring to
give a direct construction.  The outcome of this section is that we
can associate to categories enriched in orbispaces a topological
category via the Borel construction.

\begin{notn}
We will frequently use constructions of homotopy colimits and limits,
which we do using explicit models in terms of the (co)bar
construction.  See Appendix~\ref{sec:spectralbar} for a review of
the bar construction in the context of spectra, and Appendix \ref{sec:hoco} for a review of homotopy
(co)limits and for various technical
results that we need. 
\end{notn}

\subsubsection{Orbispace presentations}
\label{sec:orbisp-pres}

We begin by defining the category $\Orbo$ of \emph{charts of
  orbispaces}, which is sometimes referred to as the category of
equivariant spaces (the subscript $\emptyset$ is required for
consistency with later notation, where we will introduce stratified
orbispaces): 

\begin{defin}
The category $\Orbo$ has objects pairs $(Z, G)$, where $G$ is a finite
group and $Z$ is a $G$-space.  A morphism $(Z_0, G_0) \to (Z_1, G_1)$
is specified by a homomorphism $p\colon G_0 \to G_1$ and a
$G_0$-equivariant map $Z_0 \to p^* Z_1$.  Composition of group
homomorphisms and of maps of spaces defines the composition of
morphisms in $\Orbo$.
\end{defin}

In the above definition, $p^* Z_1$ refers to considering $Z_1$ as a
$G_0$ space with the action induced by the map from $G_0$ to $G_1$. We
shall often abuse notation and write $Z_1$ for $p^* Z_1$.  Note that
the composition depends on the fact that for $p_1 \colon G_0 \to G_1$
and $p_2 \colon G_1 \to G_2$, there is an identification 
\begin{equation}
(p_2 \circ p_1)^* Z_2 = p_2^* (p_1^* Z_2).
\end{equation}

\begin{rem}
Although we do not need this formalism in the paper, $\Orbo$ is a
enriched indexed category~\cite{Shulman2013} --- the enrichments on
mapping spaces vary with the domain and these enrichments are
compatible in a precise sense.
\end{rem}

We have a natural functor
\begin{align} \label{eq:quotient_functor}
 \Orbo & \to \Top \\
  (Z,G) & \mapsto Z/G.
\end{align}

In particular, associated to any functor 
\begin{equation}
(Z_\bullet,G_\bullet)  \co  A \to \Orbo
\end{equation}
with domain a small category $A$, we can construct a space
\begin{equation}
\colim_{\alpha \in A}  Z_\alpha/G_\alpha.
\end{equation}
as the colimit of the composite $A \to \Orbo \to \Top$.

There is a distinguished class of morphisms in $\Orbo$:
\begin{defin}
An \emph{open embedding} $(Z_0, G_0) \to (Z_1, G_1)$ is a morphism of
charts of orbispaces such that  
\begin{enumerate}
\item $\ker(G_0 \to G_1)$ acts freely on $Z_0$, and 
\item the induced map $Z_0 \times_{G_0} G_1 \to Z_1$ is an open
  embedding.  
\end{enumerate}
  
\end{defin}

The freeness of the action by the kernel readily implies:

\begin{lem} \label{lem:isotropy_group_invariant}
If $f \co (Z_0, G_0) \to (Z_1, G_1)$ is an open embedding,
then the map $G_0 \to G_1$ induces an isomorphism
$G_{z_0} \to G_{f(z_0)}$ of stabilizer groups for each $z
\in Z_0$. \qed
\end{lem}

The cartesian product induces a symmetric monoidal structure on
$\Orbo$.

\begin{lem}
The category of charts of orbispaces is symmetric monoidal, with
product  
\begin{equation}
    \Orbo \times \Orbo \to \Orbo
\end{equation}
specified on objects by the assignment
\begin{equation}
\left((Z_{1}, G_1), (Z_2, G_{2})\right) \mapsto  (Z_{1}
\times Z_{2}, G_1 \times G_2).
\end{equation}
On morphisms, we define the map
\begin{equation}
    f_1 \times f_2 \co (Z_{1} \times Z_{2}, G_1 \times G_2)   \to (Z'_1 \times Z'_2, G'_1 \times G'_2)
\end{equation}
induced by arrows $f_i \co (Z_i, G_i) \to (Z'_i,G'_i)$ to be
the product of the maps of spaces and of groups.  \qed
\end{lem}

We now define the version of an orbispace structure on a topological
space that we use.

\begin{defin}
An \emph{orbispace cover} of a topological space $\cM$ consists of the
following data: 
\begin{enumerate}
\item A diagram $(Z_\bullet, G_\bullet) \co A \to \Orbo$ factoring
  through the subcategory of open embeddings.
\item A homeomorphism $ \colim_{\alpha \in A}  Z_\alpha/G_\alpha \to  \cM$.
\end{enumerate}
The \emph{isotropy group} $G_{[z]}$ of a point $[z] \in \cM$ (with
respect to an orbispace cover) is the stabilizer group
$(G_{\alpha})_z$ for any lift of $[z]$ to $z \in Z_{\alpha}$.
\end{defin} 

Notice that in fact the data of the topological space $\cM$ is
redundant; it suffices to simply provide the diagram $(Z_\bullet,
G_\bullet) \co A \to \Orbo$ and define $\cM$ to be the colimit.  As a
consequence, we will often describe an orbispace presentation in terms
of the diagram and suppress $\cM$.

\begin{defin}\label{defn:orbipres}
An orbispace cover of $\cM$ is an \emph{orbispace presentation} if for
each point $[z] \in Z$, the nerve of $A_{[z]}$ is contractible, where
$A_{[z]} \subset A$ denote the full subcategory consisting of objects
$\alpha$ of $A$ such that the image of the map $Z_\alpha \to \cM$
contains $[z]$.
\end{defin}

This condition is essential to ensure that the algebraic constructions
that we will study, and that are expressed in terms of the
presentation $(Z_\bullet, G_\bullet) \co A \to \Orbo $, accurately
reflect the intrinsic geometric features of the space $\cM$.  The
first instance of the importance of this constraint is the following
consistency check.

\begin{lem}
If $(Z_\bullet,
G_\bullet) \co A \to \Orbo$ is an orbispace presentation of $\cM$, the natural map
\begin{equation}
  \hocolim_{\alpha \in A} Z_\alpha/G_\alpha \to \cM
\end{equation}
is a quasifibration with contractible homotopy fibers (in particular,
a weak equivalence).
\end{lem}

\begin{proof}
Since the homotopy colimit is over a functor to $\Top$ taking values
in open embeddings of subspaces of $\cM$, the nerve condition allows
us to apply the homotopical Siefert-van Kampen
theorem~\cite[A.3.1]{Lurie2014} and conclude that the natural map is a
weak equivalence.  The nerve condition also implies that the actual
fibers are contractible and therefore that the inclusion of the fiber
in the homotopy fiber is an equivalence; the map is a quasifibration.
\end{proof}

Given a pair $(\cM, \partial \cM)$, and an orbispace cover $(Z_\bullet,
G_\bullet) \co A \to \Orbo$ of $\cM$, we define an orbispace cover
\begin{equation}
(\partial Z_\bullet, G_\bullet) \co A \to \Orbo
\end{equation}
of $\partial \cM$, which we call the boundary cover, as follows: for each $\alpha$ in $A$, projection
induces a map 
\begin{equation}
Z_\alpha \to Z_\alpha / G_\alpha \to \cM.
\end{equation}
This permits us to define $\partial Z_\alpha$ as the inverse image of
the boundary:
\begin{equation}
 \partial Z_\alpha \equiv Z_\alpha \times_{\cM} \partial \cM.
\end{equation}
Notice that $\partial Z_\alpha$ inherits a $G_\alpha$ action since the
pullback factors through the orbits $Z_\alpha / G_\alpha$ and is
evidently functorial in $A$.

The following result makes it possible to use the above notion in our constructions:
\begin{lemma}
Let $(\cM, \partial \cM)$ be a pair.  Then the diagram $(\partial Z_\bullet,
G_\bullet)$ is an orbispace cover of $\partial \cM$.  If $(Z_\bullet,
G_\bullet)$ is an orbispace presentation of $\cM$, then this diagram
is an orbispace presentation of $\partial \cM$. \qed
\end{lemma}

\begin{proof}
Since colimits in spaces are stable under base change, 
\begin{equation}
\colim_{\alpha \in A} \partial Z_\alpha / G_\alpha \cong 
(\colim_{\alpha \in A} Z_\alpha / G_\alpha) \times_{\cM} \partial \cM
\cong \partial \cM.
\end{equation}
It is straightforward to check that the open embedding condition is
satisfied.  Moreover, it is clear that the contractibility of the
nerves of $A_{[z]}$ for $(Z_\bullet, G_\bullet)$ implies the same for
$(\partial Z_\bullet, G_\bullet)$.
\end{proof}

\subsubsection{Orbispace flow categories}
\label{sec:orbisp-flow}

We now return to the setting of Section~\ref{sec:topol-enrich-categ}.
We assume that $\cP$ is a partially ordered set with a free action of
a discrete group $\Pi$.  In order to proceed, we will need to employ
the notions of bicategory (a weak $2$-category) and various kinds of
$2$-functors between bicategories.  We review these definitions in
Appendix~\ref{sec:2categories}, but we also take pains to spell out in
some detail the structure we are working with.

First, note that any monoidal category $\cC$ gives rise to an example
of a bicategory sometimes denoted $B\cC$ in which there is a single
object, the category of endomorphisms of that object is given by
$\cC$, and the composition is given by the monoidal product.  This is
a bicategory and not a strict $2$-category because the product in
a monoidal category is only associative up to coherent isomorphism.
We are particularly interested in this construction applied to
$\Orbo$; we will often abusively denote the resulting bicategory by
$\Orbo$.  Note that if we replace $\cC$ by an equivalent strict
monoidal category (where the associativity isomorphisms are the
identity), then this construction yields a strict $2$-category.

For strict $2$-categories $\cC$ and $\cD$, a strict functor $F \colon
\cC \to \cD$ is specified by a function $\ob(\cC) \to \ob(\cD)$,
functors $F_{x,y} \colon \cC(x,y) \to \cD(Fx,Fy)$ for each pair of
objects $x$, $y$, and natural identities $Ff \circ Fg \to F(f \circ
g)$ along with associativity and unitality data expressed by natural
transformations.  But for bicategories $\cC$ and $\cD$, it is useful
to work with functors where the composition transformation is not the
identity.  A lax functor $F \colon \cC \to \cD$ is specified by a
function $\ob(\cC) \to \ob(\cD)$, functors $F_{x,y} \colon \cC(x,y)
\to \cD(Fx,Fy)$ for each pair of objects $x$, $y$, and natural
transformations $Ff \circ Fg \to F(f \circ g)$ along with
associativity and unitality data expressed by natural transformations
(see Definition~\ref{defin:lax-functor}).

\begin{example}\label{ex:enriched-cat-as-lax}
One of the original motivations for the definition of bicategories and
lax functors was the example of enriched categories.  Let $\cV$ be a
monoidal category regarded as a bicategory as above and let $S$ be a
set.  We will regard $S$ as a bicategory by letting each category of
morphisms be the terminal category with a single object and single
morphism.  Then a lax functor $S \to \cV$ is equivalent to a category
enriched in $\cV$ with object set $S$.
\end{example}

Next, we need to introduce the (somewhat less standard) notion of a
$\Pi$-equivariant bicategory $\cC$.  We explain this in gory detail in
Appendix~\ref{sec:equibicat}, but for now we point out that we have in
mind an extremely strict notion insofar as we mean an action of $\Pi$
on the set of objects ($0$-cells) such that for $\pi \in \Pi$ we have
identity maps
\begin{equation}
\cC(x,y) \to \cC(\pi x, \pi y)
\end{equation}
which are compatible with the horizontal composition.  In fact, in our
main examples, we will work with $\Pi$-equivariant $2$-categories
(strict bicategories).  A $\Pi$-equivariant $2$-functor $\cC \to \cD$
is a $2$-functor which is strictly compatible with the $\Pi$ action on
the $1$-cells; this compatibility is expressed by natural identities
interchanging $F$ and the action.

We are now ready to define an orbispace flow category.  In light of
Example~\ref{ex:enriched-cat-as-lax}, an orbispace flow category
should be thought of as a way of presenting a category enriched in
orbispaces.  In the following, we regard $\Orbo$ as a bicategory
having trivial $\Pi$-action.

\begin{defin} \label{def:orbispace_flow_category}
  An \emph{orbispace flow category} consists of the following data:
  \begin{enumerate}

  \item  A $\Pi$-equivariant $2$-category $A$, with object set $\cP$,
    and such that $A(p,p) = \ast$ (i.e., the category with a single
    object and morphism) and $A(p,q)$ is empty unless $p \leq q$.

  \item A strictly $\Pi$-equivariant normal lax $2$-functor $(Z_\bullet,
    G_\bullet) \co A \to \Orbo$, such that the value for any pair
    $(p,q)$ is an orbispace presentation.  Here normal means that the
    $2$-functor strictly preserves identities.

\end{enumerate}
\end{defin}

More explicitly, such a flow category consists of categories $A(p,q)$
for each pair of comparable elements of $\cP$, together with
composition functors 
\begin{equation} \label{eq:map_product_indexing}
 A(q,r) \times A(p,q)  \to A(p,r),
  \end{equation}
that are strictly associative.  The action of $\Pi$ permutes the
categories $A(p,q)$ in that we require that $A(p,q)$ be equal to
$A(\pi p, \pi q)$ for $\pi \in \Pi$.

Moreover, we have an orbispace presentation $(Z_\bullet,
G_\bullet)(p,q)  \co A(p,q) \to \Orbo$ for each ordered pair $p<q$,
and natural transformations
\begin{equation}
(Z_\bullet, G_\bullet)(p,q) \times  (Z_\bullet, G_\bullet)(q,r)
  \Rightarrow (Z_\bullet, G_\bullet)(p,r)    
\end{equation}
such that the two induced natural transformations
\begin{equation}
(Z_\bullet, G_\bullet)(p,q) \times  (Z_\bullet, G_\bullet)(q,r) \times
  (Z_\bullet, G_\bullet)(r,s) \Rightarrow (Z_\bullet, G_\bullet)(r,s)  
\end{equation}
on the product $ A(p,q) \times A(q,r) \times A(r,s) $ coincide (there
are similar diagrams and conditions for units which we do not describe
here).

The $\Pi$-equivariance of the $2$-functor to $\Orbo$ means that the
orbispace presentations 
\begin{equation}
(Z_\bullet, G_\bullet)(p,q)  \co A(p,q) \to \Orbo
\end{equation}
and
\begin{equation}
(Z_\bullet, G_\bullet)(\pi p, \pi q)  \co A(\pi p,\pi q) \to \Orbo
\end{equation}
coincide.

We can associate to each orbispace flow category $(A, Z_\bullet,
G_\bullet)$ a $\Pi$-equivariant topological flow category $\cM_A$,
i.e., a lax functor from the chaotic category of $\cP$ to the
bicategory associated to the monoidal category of spaces. 

\begin{prop}\label{prop:orbispace-to-flow-category}
Let $(A, Z_\bullet, G_\bullet)$ be an orbispace flow category.  Then 
there is a $\Pi$-equivariant topological flow category with object set 
$\cP$ and morphisms from $p$ to $q$ given by 
\begin{equation}
\cM_A  (p,q) = \colim_{\alpha \in A(p,q)} Z_\alpha/G_\alpha. 
\end{equation}
\end{prop}

\begin{proof}
First, observe that the composition 
\begin{multline}
\colim_{\alpha \in A(p,q)} Z_\alpha/G_\alpha \times \colim_{\beta \in
  A(q,r)} Z_\beta/G_\beta   \to \\
\colim_{(\alpha,\beta) \in A(p,q) \times
  A(q,r)} \left(Z_\alpha/G_\alpha \times Z_\beta/G_\beta\right)  \to \\
 \colim_{(\alpha,\beta) \in A(p,q) \times
  A(q,r)} (Z_\alpha \times Z_\beta)/(G_\alpha \times G_\beta) \to \colim_{\gamma \in A(p,r)} Z_\gamma/G_\gamma
\end{multline}
is associative and unital.  The assignment $Z_\bullet/G_\bullet$
specifies a monoidal functor $\Orbo \to \Top$ and therefore induces a
$2$-functor on the associated bicategories; this gives rise to a
normal lax $2$-functor from $A$ to $\Top$.  Since $A$ is a
$2$-category and $\colim$ is a monoidal functor in the indexing
category, $\cM_A$ as specified is a topological category.  Put more
succinctly, applying the colimit functor to the hom categories of a
normal lax $2$-functor from $A$ to $\Top$ yields a lax functor from
the chaotic category on $\cP$ to $\Top$.

Next, we need to describe the action of $\Pi$.  For each $\pi \in
\Pi$, the strict $2$-functor $\gamma_{\pi}$ encoding the action
permutes the colimits above and so gives rise to a functor
$\gamma^{\cM}_{\pi} \colon \cM_A \to \cM_A$; it is straightforward to
check that this is associative and unital.
\end{proof}

\begin{rem}
Our example of interest fits into the framework of
Definition~\ref{def:orbispace_flow_category}, and imposing additional
strictness reduces technical complexity.  However, we can in fact take
as input for our construction a bicategory $A$ rather than a strict
$2$-category without much trouble.  Rectifying $A$ to a biequivalent
strict $2$-category (see Theorem~\ref{thm:2catrect} for a precise
statement), we can proceed as above in
Proposition~\ref{prop:orbispace-to-flow-category}, although the resulting
action of $\Pi$ satisfies the associativity composition only up to natural
isomorphism; the collection of functors specify a pseudo-action of
$\Pi$ on $\cM_A$.  Rectifying this action to a strict action
of $\Pi$ yields the desired flow category structure on $\cM_A$.

More generally, our definition can be phrased in a homotopy-coherent
setup; ultimately, this will be the correct form of these
constructions, and in this context the distinctions we are making
above disappear.
\end{rem}

We can now formulate the notion of a flow category that has an
orbispace structure.

\begin{defin}
If $\cM$ is a topological flow category, a \emph{lift of $\cM$ to an
  orbispace category} is an orbispace category $(A, Z_\bullet, G_\bullet)$ equipped with an
isomorphism of topological categories
\begin{equation}
  \cM_A \to \cM.
\end{equation}
\end{defin}

In particular, note that a lift to an orbispace category means we have
homeomorphisms $\cM_A(p,q) \cong \cM(p,q)$ for all $(p,q)$.

\subsubsection{The Borel construction}
\label{sec:borel-construction-1}

The Borel construction can be applied to an orbispace flow category.
We will write $E G$ for the two-sided bar construction $B(G, G,
\ast)$, providing a functorial model for the universal space for $G$. Using this model, we write
\begin{equation}
E \colon \Orbo \to \Orbo  
\end{equation}
for the strong monoidal functor
specified on objects by the assignment
\begin{equation}
(Z, G) \mapsto (Z \times EG, G).
\end{equation}
The assertion that this is a strong monoidal functor amounts to the fact that the natural map
\begin{equation}
(Z_1 \times EG_1) \times (Z_2 \times EG_2) \cong (Z_1 \times Z_2)
  \times (E(G_1 \times G_2))
\end{equation} 
is a homeomorphism for our chosen model of the universal space.

Formally, the functor $E$ has the effect of cofibrantly replacing each
space in the Borel model structure for its group of equivariance.
Composing $E$ and the quotient yields the Borel construction; we can
alternatively describe this as the bar construction $B(\ast, G, X)$.  

\begin{lem}\label{lem:borel-orbo-functor}
There is a strong monoidal Borel construction functor 
\begin{equation}
B \colon \Orbo \to \Top
\end{equation}
specified on objects by the assignment 
\begin{equation}
(Z, G) \mapsto (Z \times EG / G) = Z \times_G EG
\end{equation}
and on morphisms by assigning to $f\colon (Z,G) \to (Z',G')$ the
composite 
\begin{equation}
Z \times_G EG \to (f^* Z' \times f^* EG') / G \to Z' \times_{G'} EG' 
\end{equation}
The monoidal structure is induced by the homeomorphism
\begin{equation}
(Z_1 \times_{G_1} EG_1) \times (Z_2 \times_{G_2} EG_2) \cong
(Z_1 \times Z_2) \times_{G_1 \times G_2} (E(G_1 \times G_2)).
\end{equation}
\end{lem}

Since the Borel construction is functorial, we can apply it to
orbispace presentations.

\begin{defin}
The \emph{Borel
  construction} of an orbispace presentation
$(Z_\bullet, G_\bullet) \colon A \to \Orbo$ is the functor
\begin{equation}
  BZ_{\bullet} \co A \to \Top 
\end{equation}
produced as the composite
\begin{equation}
\begin{tikzcd}
A \ar[r] & \Orbo \ar[r] & \Top
\end{tikzcd}
\end{equation}
which assigns to each object $\alpha \in A$ the space
\begin{equation}\label{eq:Borel_construction}
    BZ_{\alpha} \equiv Z_\alpha \times_{G_\alpha} E G_\alpha.
  \end{equation}
\end{defin}

Passing to the homotopy colimit, we obtain a Borel homotopy type for
an orbispace presentation.
  
\begin{defin}
The
\emph{Borel homotopy type} of an orbispace presentation $(A, Z_\bullet, G_\bullet)$,  is the space  
\begin{equation}
BZ(A) \equiv \hocolim_{\alpha \in A} BZ_\alpha.
\end{equation}
\end{defin}

Note that the notation is potentially ambiguous since the Borel homotopy
type depends not only on the indexing category $A$, but also on the
functor $A \to \Orbo$.  In practice, no confusion shall arise, and we
shall further abuse notation by dropping $A$ from the notation when
this is warranted.

By construction, for an orbispace presentation with associated space
$\cM$, we have a natural map
  \begin{equation}
     BZ  \to \cM   
   \end{equation}
induced by the canonical natural transformation $Z_\alpha
\times_{G_\alpha} EG_\alpha \to Z_\alpha / G_\alpha$
determined by the terminal map $EG_\alpha \to \ast$ and the
natural map 
\begin{equation}
\hocolim_{\alpha \in A} Z_\alpha/G_\alpha \to \colim_{\alpha \in A}
Z_\alpha/G_\alpha.
\end{equation}

As a sanity check for our definitions, we identify the fiber of this
map; this is another consequence of the acyclic nerve condition in the
definition of an orbispace presentation.

\begin{lem}
Given an orbispace presentation, the fibre $B_{[z]} Z$ of the map $BZ
\to \cM$ at a point $[z] \in \cM$ is equivalent to the classifying
space $B G_{[z]}$.
\end{lem}

\begin{proof}
Let $G_{\alpha_z}$ denote the stabilizer group of a point $z \in
Z_\alpha$ mapping to $[z]$; this subgroup of $G_\alpha$ is
well-defined up to conjugacy, and the fibre of
\begin{equation}
BZ_\alpha = Z_\alpha \times_{G \alpha} EG_{\alpha} \to Z_\alpha /
G_\alpha \to \cM
\end{equation}
can be identified with $B G_{\alpha_z}$.  The fibre of the map $BZ$ to
$\cM$ is then given by the homotopy colimit of $B G_{\alpha_z}$ over
$A$, as can be seen by factoring this map through the homotopy colimit
of $Z_\alpha/G_\alpha$:
\begin{equation}
BZ_\bullet \otimes_A N(- \downarrow A) \to (Z_\bullet / G_\bullet)
\otimes_A N(- \downarrow A) \to (Z_\bullet / G_\bullet) \otimes_A *.
\end{equation}
where the second arrow is induced by the map from the functor $A \to
\Top$ that is the nerve of the overcategory at an object of $A$ to the
constant functor on a point.  Here, we are using the ``tensor product
of functors'' description of the coend computing the (homotopy)
colimit; see Appendix~\ref{sec:hoco} for a review.

Each map of $f \co \alpha \to \beta$ of charts whose image contains
$[z]$ induces a map $B G_{\alpha_z} \to B G_{\beta_z}$ which is an
equivalence by Lemma~\ref{lem:isotropy_group_invariant}, since we
assumed that $G_f$ acts freely on the fibres.  Since the nerve of
$A_{[z]}$ is contractible by the assumption that we start with an
orbispace presentation, we conclude that the homotopy colimit
\begin{equation}
  \hocolim_{A} B G_{\alpha_z} \cong \hocolim_{A_{[z]}} B G_{\alpha_z}
\end{equation}
is equivalent to $B G_{[z]}$.
\end{proof}

Given an orbispace lift $(A, Z_\bullet, G_\bullet)$ of a strictly
$\Pi$-equivariant flow category $\cM$ over a partially ordered set
$\cP$, we will define an associated strictly $\Pi$-equivariant
topological category $BZ$.  The argument for the following proposition
is the same as the proof of
Proposition~\ref{prop:orbispace-to-flow-category}, relying on
Lemma~\ref{lem:borel-orbo-functor} to produce a lax $2$-functor
from $A$ to $\Top$ and the fact that the homotopy colimit has the same
monoidal properties as the colimit.

\begin{prop}\label{prop:borel-category}
Let $(A, Z_\bullet, G_\bullet)$ be an orbispace lift of $\cM$.  Then
there is a strictly $\Pi$-equivariant topological category $BZ$ with
objects those of $\cP$ and morphism spaces determined as
\begin{equation}
  BZ(p,q) \equiv BZ(A(p,q))
\end{equation}
via the Borel construction.  The composition 
\begin{equation}\label{eq:topcomp}
BZ(A(p,q)) \times BZ(A(q,r)) \to BZ(A(p,q) \times A(q,r)) \to BZ(A(p,r))
\end{equation}
is determined by the horizontal composition in $A$ and
Lemma~\ref{lem:borel-orbo-functor}. \qed
\end{prop}

By construction, $BZ$ is compatible with $\cM$ in the following sense.

\begin{lem}
There is a natural $\Pi$-equivariant functor
 \begin{equation}
    B Z \to \cM, 
 \end{equation}
which is the identity on objects. \qed
\end{lem}

We may think of $B Z$ as the topological flow category associated by
the Borel construction to an orbispace flow category. The relative
cochains of the morphism spaces in $BZ$ yield a strictly
$\Pi$-equivariant category 
\begin{equation}
  C^*_{\rel \partial}( B Z; \Omega \bbk)  
\end{equation}
by applying the construction from Section~\ref{sec:relat-coch-coll},
with morphism spectra 
\begin{equation} \label{eq:definition_relative_cochains_Borel-naive}
  C^*( \widehat{B Z}(p,q), \partial \widehat{B Z}(p,q); \Omega \bbk).
\end{equation}
The functoriality of this construction readily implies: 
\begin{lem}\label{lem:underlying-cochain-comparison}
There is a natural $\Pi$-equivariant functor
\begin{equation}
C^*_{\rel \partial}( {\cM}; \Omega \bbk) \to   C^*_{\rel \partial}( B Z; \Omega \bbk) .
\end{equation} \qed
\end{lem}

While Lemma~\ref{lem:underlying-cochain-comparison} is conceptually
satisfying, it shall play no direct role in our construction: we shall
henceforth work exclusively with $ C^*_{\rel \partial}( B Z; \Omega
\bbk) $ as our model for the relative cochains of an orbispace flow
category.  The intuition here is that this spectral category reflects
some of the features of the given orbispace category which should
not be forgotten, as one would do by considering the ordinary relative
cochains of the underlying topological category.

The comparison map that will play a key role in our construction is
the spectral functor 
\begin{equation}
 C^*_{\rel \partial}( {\cP}; \Omega \bbk) \to   C^*_{\rel \partial}( B
 Z; \Omega \bbk)   
\end{equation}
induced by the fact that $BZ$ is a category over $\cP$.

\section{Partial orders and homotopy types}
\label{sec:from-flow-categories}

The purpose of this section is to explain how, given a locally finite
dimensional partially ordered set $\cP$ as in the previous section, we
may associate a homotopy type filtered by $\cP$ to each spectral
functor from $C^{*}_{\rel \partial}(\cP; \Omega \bbk)$ to the category
of $\bbk$-modules. For specificity, we assume that $\cP$ is
additionally equipped with an order preserving map to $\bR$, which we
refer to as the action (see Section \ref{sec:homotopy-types-from}
below).  Assuming that $\cP$ carries a free action (over $\bR$) of a
group $\Pi$ with a homomorphism to $\bR$, and that the spectral
functor is equivariant, we show that the associated filtered homotopy type has
an action by a Novikov ring associated to $\Pi$.

\subsection{The spectral Novikov ring}
\label{sec:spectr-novik-ring}

In this section, we will construct a spectral Novikov ring associated
to an associative ring spectrum $\bbk$ and a discrete group $\Pi$
equipped with action and index homomorphisms  
\begin{equation} \label{eq:action_and_degree_homomorphisms}
  \cA \times \deg \co \Pi \to \bR \times \bZ  .
\end{equation}
These Novikov rings arise as certain
completions of the ``spectral group ring'' $\bbk[\Pi]$.  

We begin by recalling the algebraic version of the Novikov ring that
we construct.  We define
\begin{equation}
\begin{tikzcd}
\Pi_{> \lambda} := \{\pi \in \Pi \, \mid \, \cA(\pi) > \lambda\} &&
\Pi_{\geq \lambda} := \{\pi \in \Pi \, \mid \, \cA(\pi) \geq \lambda\}
\end{tikzcd}
\end{equation}
and $\Pi_{< \lambda}$ and $\Pi_{\leq \lambda}$ analogously.  

Let $R$ be a commutative ring.  We want to define the Novikov ring
$\Lambda_{R,\Pi}$ to be the completion of the group ring $R[\Pi]$ with
respect to the decreasing filtration induced by $\cA$, i.e., by
$\Pi_{> \lambda}$.  As an $R$-module, this can be described as the
inverse limit 
\begin{equation}
\Lambda_{R,\Pi} = \lim_{\lambda} \frac{R[\Pi]}{R[\Pi_{> \lambda}]}.
\end{equation}  
However, the preceding candidate definition of the Novikov ring
$\Lambda_{R,\Pi}$ obscures a description of ring structure, as
$R[\Pi]/\Pi_{\lambda}$ is not itself a ring.  For this purpose, we can
consider the group ring $R[\Pi_{\geq 0}]$ and consider the completion
of this at the decreasing filtration specified by $\Pi_{> \lambda}$,
for $\lambda \geq 0$.  This can be described as the inverse limit 
\begin{equation}
\Lambda^+_{R,\Pi} = \lim_{\lambda} \frac{R[\Pi_{\geq 0}]}{R[\Pi_{\geq 
    \lambda}]}.
\end{equation} 
Here since $R[\Pi_{\geq \lambda}]$ is an ideal in $R[\Pi_{\geq 0}]$,
this inverse limit is evidently itself a ring.  To produce the Novikov
ring itself, we invert the elements $\Pi_{> 0}$.  Note that we have
not assumed that $\Pi$ is commutative, but this localization makes
sense because $\Pi_{\geq 0}$ and $\Pi_{> 0}$ are closed under
conjugation:
\begin{equation}
\cA(\pi_1 \pi_2 \pi_1^{-1}) = \cA(\pi_1) + \cA(\pi_2) - \cA(\pi_1) =
\cA(\pi_2).
\end{equation}
Then for any $x \in \Pi^0_+$ and $y \in \Pi^0$, we have $(yxy^{-1})y =
yx$.  That is, the (left) Ore condition is satisfied.

We now turn to produce a version of these completions in the context
of the ring spectrum $\bbk$.  Note that in our application $\bbk$ is
not a commutative ring spectrum.  Nonetheless, we can retain control
on the completions we study because $\pi_* (\bbk)$ is graded
commutative and because the existence of the action map ensures that
the Ore condition will hold.

\begin{defin}
The spectral group ring $\bbk[\Pi]$ is the associative ring spectrum
with underlying spectrum the free $\bbk$-module $\bbk \sma
\Sigma^\infty_+ \Pi$, with multiplication map
\begin{equation}
(\bbk \sma \Sigma^\infty_+ \Pi) \sma (\bbk \sma \Sigma^\infty_+ \Pi) \cong (\bbk \sma \bbk) \sma
\Sigma^\infty_+ (\Pi \times \Pi) \to \bbk \sma \Sigma^\infty_+ \Pi.
\end{equation}
Here the first map is induced by the transposition maps and the second
map is induced by the multiplications on $\bbk$ and $\Pi$ respectively.
\end{defin}

We will now introduce a graded variant of the group ring $\bbk[\Pi]$.
Specifically, we want to twist the wedge summand associated to $\pi
\in \Pi$ by the sphere of dimension $\deg(\pi$).  To describe the
multiplicative structure, we find it convenient to use the strictly
multiplicative system of even spheres $\{\bS[n]\}_{n \in \bZ}$
described in Appendix~\ref{sec:action}; these have the property that
$\bS[n] \htp S^{-2n}$ and there are associative multiplication maps
\begin{equation}
\bS[n] \sma \bS[m] \to \bS[n+m]
\end{equation}
for all $m,n \in \bZ$.  We build the group ring using the category
of graded ring spectra; we review this in detail in
Appendix~\ref{sec:action} as well, but for now we simply note that for
a (possibly non-unital) monoid $M$, an $M$-graded spectrum is a
collection of spectra $\{E_m\}$ for $m \in M$.  An $M$-graded ring
spectrum is an $M$-graded spectrum equipped with associative
multiplication maps $E_m \sma E_n \to E_{n+m}$ for $m,n \in M$; we
require these to be unital if $M$ is.  Associated to an $M$-graded
spectrum $E$ is an underlying spectrum that is equivalent to the wedge
$\bigvee_{m \in M} E_m$; when $E$ is an $M$-graded ring spectrum, the
underlying spectrum is a ring spectrum.

\begin{defin}
The {\em graded spectral group ring} $\Sigma^{\deg} \bbk[\Pi]$ is the
underlying spectrum associated to the $\Pi$-graded spectrum whose
value at $\pi$ is
\begin{equation} \label{eq:shift_k_by-deg}
\bS[-\deg(\pi)] \sma \bbk,
\end{equation}
equipped with multiplication maps 
\begin{equation} \label{eq:product_graded_twisted_group_ring}
\left(   \bS[-\deg(\pi_1)] \sma \bbk \right) \sma \left( \bS[-\deg(\pi_2)] \sma \bbk \right) \to \bS[-\deg(\pi_1 \cdot \pi_2)] \sma \bbk
\end{equation}
induced by the products on $\{\bS[n]\}$, $\bbk$, and $\Pi$.
\end{defin}

We now define a decreasing filtration on $\Sigma^{\deg} \bbk[\Pi]$
induced by the action homomorphism.  For this purpose, we use the
homotopical version of the filtered derived category (see
Appendix~\ref{appendix:complete} for a quick review or~\cite[\S
  5.1]{BhattMorrowScholze2019} and~\cite{GwilliamPavlov2018} for more
detailed treatments).  By a decreasing filtration on a spectrum $X$ we
mean a functor $F \colon \bR^{\op} \to \Sp$ such that each map
$F(\lambda) \to F(\lambda')$ for $\lambda > \lambda'$ and $F(-\infty)
= \hocolim_F \htp X$.  

We let $\Sigma^{\deg} \bbk[\Pi_{> \lambda}]$ denote the underlying
spectrum associated to the $\Pi_{> \lambda}$-graded spectrum whose
value at $\pi \in \Pi_{> \lambda}$ is given by
Equation~\eqref{eq:shift_k_by-deg}.  For $\lambda > \lambda'$ there is
a natural inclusion $\Pi_{> \lambda} \to \Pi_{> \lambda'}$ that
induces commutative diagrams 
\begin{equation}
\begin{tikzcd}
\Pi_{> \lambda} \ar[d] \ar[r] & \Pi && \Sigma^{\deg} \bbk[\Pi_{> \lambda}] \ar[d] \ar[r] & \Sigma^{\deg}
\bbk[\Pi] \\
\Pi_{> \lambda'} \ar[ur] & && \Sigma^{\deg} \bbk[\Pi_{> \lambda'}] \ar[ur]  
\end{tikzcd}
\end{equation}
Thus, this specifies a decreasing filtration on $\Sigma^{\deg}
\bbk[\Pi]$.

We can complete a filtered spectrum by taking the Bousfield
localization with respect to the maps which induce equivalences of
associated graded spectra; see Definition~\ref{defn:completion} (and
the surrounding discussion in Appendix~\ref{appendix:complete}).  The
colimit of the resulting filtered spectrum is the completion of
$\Sigma^{\deg} \bbk[\Pi]$.  To describe the completion explicitly,
notice that there is an induced map of cofibers
\begin{equation}
\bbk[\Pi] / \bbk[\Pi_{> \lambda}] \to \bbk[\Pi] /
\bbk[\Pi_{> \lambda'}]
\end{equation}
and hence a functor $\bR^{\op} \to \Sp$ specified on objects by the
assignment 
\begin{equation}
\lambda \mapsto \bbk[\Pi] / \bbk[\Pi_{> \lambda}].
\end{equation}
Since colimits commute with the smash product, there is a natural
isomorphism  
\begin{equation}
\bbk[\Pi] / \bbk[\Pi_{> \lambda}] \cong \bbk[\Pi_{\leq \lambda}],
\end{equation}
where $\Pi_{\leq \lambda}$ denotes the subset of $\Pi$ of elements of
action less than or equal to $\lambda$.

\begin{defin}
The {\em completed graded group ring} is the completion
of $\Sigma^{\deg} \bbk[\Pi]$ with respect to the decreasing filtration by
action.  This is specified by the formula
\begin{equation}
\Sigma^{\deg} \bbk((\Pi)) = \holim_{\lambda \in \bZ^{\op}} \left(\Sigma^{\deg}
\bbk[\Pi] / \Sigma^{\deg} \bbk[\Pi_{> \lambda}]\right).
\end{equation}
\end{defin}

\begin{rem}
Note that we have passed to $\bZ \subset \bR$ for the purpose of
taking the (homotopy) limit above. The only property that we use is
that $\bZ$ is discrete, closed under addition, and terminal with
respect to the usual ordering (hence initial when considered in
$\bR^{\op}$). 
\end{rem}

By construction, the completion $\Sigma^{\deg} \bbk((\Pi))$ is the
homotopy colimit of a decreasing filtration defined via
\begin{equation}
\lambda \mapsto \Sigma^{\deg} \bbk[\Pi_{> \lambda}] / \holim_\lambda
\bbk[\Pi_{> \lambda}], 
\end{equation}
which can be rewritten as
\begin{equation}
\lambda \mapsto \hofib\left(\Sigma^{\deg} \bbk((\Pi)) \to
\Sigma^{\deg} \bbk[\Pi] / \Sigma^{\deg} \bbk[\Pi_{> \lambda}]\right),
\end{equation}
and with the natural induced maps on homotopy fibers.

The completion map induces an equivalence on associated graded objects
in the following sense.

\begin{lem}
For each $\lambda \in \bR$, there is an induced equivalence
\begin{equation}
\bbk[\Pi] / \bbk[\Pi_{> \lambda}]  \htp
\bbk((\Pi)) / \hofib \left(\bbk((\Pi)) \to \bbk[\Pi] /
\bbk[\Pi_{> \lambda}] \right).
\end{equation}
For $\lambda' < \lambda$, there is an induced equivalence
\begin{equation}
\bbk[\Pi_{> \lambda}] / \bbk[\Pi_{> \lambda'}] \htp 
\frac{\bbk((\Pi)) / \hofib \left(\bbk((\Pi)) \to \bbk[\Pi_{> \lambda}]\right)}{
\bbk((\Pi)) / \hofib \left(\bbk((\Pi)) \to \bbk[\Pi_{>
    \lambda'}]\right)}.
\end{equation} \qed
\end{lem}

Next, we turn to discussion of the multiplicative structure.  The
filtration is compatible with the product: for $\lambda_1 >
\lambda_1'$ and $\lambda_2 > \lambda_2'$, there is a commutative
diagram  
\begin{equation}
\begin{tikzcd}
\Sigma^{\deg} \bbk[\Pi_{> \lambda_1}] \sma \Sigma^{\deg}
\bbk[\Pi_{> \lambda_2}] \ar[r] \ar[d] & \ar[d] \Sigma^{\deg} \bbk[\Pi_{> \lambda_1
     + \lambda_2}]. \\
\Sigma^{\deg} \bbk[\Pi_{> \lambda'_1}] \sma \Sigma^{\deg}
\bbk[\Pi_{> \lambda'_2}] \ar[r] & \Sigma^{\deg} \bbk[\Pi_{> \lambda'_1
     + \lambda'_2}].
\end{tikzcd} 
\end{equation}
induced by the fact that the action map $\aA$ is a homomorphism.
These maps are associative, and therefore imply that the decreasing
filtration on $\Sigma^{\deg} \bbk[\Pi]$ makes it into a monoid in the
category of filtered spectra.  Since the completion functor is lax monoidal
(see Theorem~\ref{thm:compmonoid}), there is an induced multiplicative
structure on the completion as a filtration.  Moreover, the passage to
the colimit is symmetric monoidal, and so we we have an induced ring
structure on $\Sigma^{\deg} \bbk[\Pi]$.

\begin{thm}
The Laurent series ring $\Sigma^{\deg} \bbk((\Pi))$ is an associative
ring spectrum. \qed
\end{thm}

\begin{notn}
Since $\Pi$ and $\bbk$ are fixed in the paper, we shall denote this completed 
group ring by $\Lambda$, and call it \emph{the Novikov ring spectrum}.
\end{notn}

On homotopy groups, the Mittag-Leffler condition is satisfied and so
we have the following computation.

\begin{lem}
There is a canonical isomorphism of graded rings
\begin{equation}
\pi_* (\Sigma^{\deg} \bbk((\Pi))) \cong \Lambda_{\pi_*(\bbk), \Pi}
\end{equation}
where each indeterminate $\pi \in \Pi$ has degree $-\deg(\pi)$. \qed
\end{lem}

\begin{proof}
The assertion about the additive structure follows by direct
calculation, using the Mittag-Leffler condition.  To see that the
multiplicative structure is correct, we can give an alternate
construction of $\Sigma^{\deg} \bbk((\Pi))$ by constructing an
analogue of $\Lambda^+$, completing $\Sigma^{\deg}\bbk[\Pi_{\geq 0}]$
using the decreasing filtration given by $\Sigma^{\deg}\bbk[\Pi_{>
    \lambda}]$ and inverting positive degree elements of $\Pi$.
The fact that the Ore condition holds implies that we can describe
this localization in terms of a calculus of fractions~\cite[\S
  7.2.3]{Lurie2014}, and in particular has the expected universal
property with respect to inverting elements.  Now the assertion
follows by a direct comparison of homotopy groups.
\end{proof}

\subsection{Homotopy types from cellular diagrams}
\label{sec:homotopy-types-from}

For the remainder of this paper, we assume that we have an action map
\begin{equation}
\cA \co  \cP \to \bR,
\end{equation}
that is compatible with the partial order on $\cP$ in the sense that
if there is arrow from $p$ to $q$, then $\cA(p) \leq \cA(q)$.
Moreover, we strengthen the condition of local finiteness by requiring
that
\begin{equation}
  \label{eq:finite_bounded_arrows}
  \parbox{31em}{for every $\lambda \in \bR$, and $p \in \cP$, there are only finitely many elements $q$ of $\cP$ receiving an arrow
from $p$, and so that  $ \cA(q) < \lambda$.}
\end{equation}

Furthermore, we assume that the action maps on $\Pi$ and $\cP$ are
compatible with the $\Pi$-action on $\cP$ in the sense that the
diagram 
\begin{equation}
  \begin{tikzcd}
   \Pi \times \cP \ar[r,"\cA"] \ar[d] & \bR \times \bR \ar[d,"+"] \\
\cP \ar[r,swap,"\cA"] & \bR 
  \end{tikzcd}
\end{equation}
commutes, i.e., the formula
\begin{equation}\label{eq:action-homomorphism}
  \cA(\pi \cdot p) = \cA(\pi) + \cA(p) 
\end{equation}
holds.

Our goal in this section is to construct a spectrum $|\delta_\lambda|$
for each action level $\lambda$ from the data of a $\Pi$-equivariant
spectral functor 
\begin{equation}
\delta \colon C^*_{\rel \partial}(\cP; \Omega \bS) \to \kmod
\end{equation}
In the next section, we shall explain how to assemble these spectra to
a inverse system. 

First, we need to explain what we mean by a $\Pi$-equivariant spectral
functor to $\bbk$-modules.  We begin by describing the action of $\Pi$
on the category $\kmod$ of (left) $\bbk$-modules.  Roughly
speaking, we want $\Pi$ to act on $\kmod$ via 
\begin{equation}
\pi \mapsto \Sigma^{\deg \pi} (-).
\end{equation}
We do this by smashing with the coherent system of models for spheres
$\{\bS[n]\}$.  Writing $\gamma_{\pi}$ for the functor
\begin{equation}
\bS[-\deg \pi] \sma (-) \colon \kmod \to \kmod,
\end{equation}
the action of $\Pi$ on $\kmod$ is specified by the assignment
\begin{equation}
\pi \mapsto \gamma_{\pi}.
\end{equation}
The multiplicative structure of the system $\{\bS[n]\}$ induces an
associative collection of composition natural transformations
$\gamma_{\pi_1} \circ \gamma_{\pi_2} \to \gamma_{\pi_1 \pi_2}$
specified by the formula
\begin{equation}
\bS[-\deg \pi_1] \sma \left(\bS[-\deg \pi_2] \sma (-)\right) \to
\bS[-\deg (\pi_1 + \pi_2)] \sma (-).
\end{equation}
However, although these transformations are natural equivalences, they
are not the identity as would be the case for a strict action.  We
refer to this structure as a homotopy action of $\Pi$; this is a lax
functor from the $2$-category generated by $\Pi$ to the $2$-category
of categories, where the associativity transformations are through
weak equivalences.  See Appendix~\ref{sec:group-acti-spectr} and in
particular Theorem~\ref{thm:kmodact} for more discussion. Equivariant functors $F \colon \cC \to \cD$ between categories with
homotopy action of $\Pi$ then come in two varieties, depending on the
direction of the natural transformations expressing the action.  We
will consider those specified by a spectral functor $F$ along with a
family of natural transformations  
\begin{equation}
\gamma_{\pi}^{\cD} \circ F \to F \circ \gamma_{\pi}^{\cC} 
\end{equation}
that are suitably compatible with the associativity and unit
transformations.  For instance, specializing to the case where the
domain category $\cC$ has a strict action of $\Pi$ and the range is
$\kmod$ with the action of $\Pi$ specified above, a $\Pi$-equivariant
functor is specified by natural transformations
\begin{equation}
\bS[-\deg(\pi)] \sma F(-) \to F(\pi -) 
\end{equation}
such that the diagram
\begin{equation}
\begin{tikzcd}
\bS[-\deg(\pi_1 )] \sma \bS[-\deg(\pi_2)] \sma F(c) \ar[r] \ar[d] 
& \bS[-\deg(\pi_1 + \pi_2)] \sma F(c) \ar[d] 
\\
\bS[-\deg(\pi_1)] \sma F(\pi_2 c) \ar[r] & F((\pi_1 + \pi_2) c)  
\end{tikzcd}
\end{equation}
commutes, along other coherence diagrams which we discuss in
Definition~\ref{defn:equivariant-functors}.

We now give the main definition of this section, for which we fix  a partially ordered set $\cP$ with a free action of $\Pi$, and a $\Pi$-equivariant map $\cA$ to $\bR$, satisfying Condition \eqref{eq:finite_bounded_arrows}:
\begin{defin}
  A \emph{$\Pi$-equivariant $\cP$-cellular
  diagram} is a $\Pi$-equivariant spectrally enriched functor  
\begin{equation}
\delta \colon C^*_{\rel \partial}(\cP; \Omega \bS) \to \kmod.
\end{equation}
That is, $\delta$ is a $\Pi$-equivariant
$C^*_{\rel \partial}(\cP; \Omega \bS)$-module over $\bbk$.
\end{defin}

In order to relate the notion of cellular diagram to flow categories,
we consider the analogous generalization:
\begin{defin}\label{defn:virtual_fundamental_chain}
  A \emph{virtual fundamental chain} for a topological flow category $\cM$ is a
  $\Pi$-equivariant spectrally enriched  functor
  \begin{equation} \label{eq:virtual_fundamental_chain}
    \delta \co C^*_{\rel \partial}(\cM; \Omega \bbk) \to
    \kmod
  \end{equation}
  with the property that the image of each object is equivalent to a
  graded line, i.e., a $\bbk$-module equivalent to $\Sigma^n \bbk$ for
  some $n \in \bZ$.
\end{defin}
\begin{rem}
  It is tempting to require that the functor in Equation \eqref{eq:virtual_fundamental_chain} be enriched over $\bbk$, but this in fact does not make sense in the setting we are considering: the ring spectrum $\bbk$ is only an associative spectrum, so that it does not make sense to consider categories enriched over $\bbk$. 
\end{rem}

As a special case, we define a virtual fundamental chain for a
Kuranishi flow category $\bX$ to be a $\Pi$-equivariant functor
\begin{equation}
\delta \co C^*_{\rel \partial}(B\cZ; \Omega \bbk) \to \kmod
\end{equation}
where $B \cZ$ is the topological flow category from
Proposition~\ref{prop:borel-category}.
  
\begin{rem}
Our formulation of the notion of virtual fundamental chain is
analogous to the one used by Pardon \cite{Pardon2016}.  The
constructions of virtual fundamental chains in Floer theory that exist
in the literature
\cite{FukayaOno1999,LiuTian1998,Pardon2016,FilippenkoWehrheim2018} amount to such a construction in the category of $H
\bQ$ modules.
\end{rem}

As observed earlier, the projection map to the collar and the unit of $\bbk$ induce a
natural $\Pi$-equivariant functor 
\begin{equation}
C^*_{\rel \partial}(\cP; \Omega \bS) \to  C^*_{\rel \partial}(\cM; \Omega \bbk).
\end{equation}
Therefore, we have the following
immediate corollary.

\begin{cor}
A virtual fundamental chain determines a $\Pi$-equivariant
$\cP$-cellular diagram.
\end{cor}

\begin{rem}
In this paper, the value of $\delta$ will always be a graded line.
However, this condition is not required for the general theory, and
applications to symplectic topology which go beyond the ones
considered here naturally lead one to consider more general values.
\end{rem}

Our construction depends on the homotopy class of $\Pi$-equivariant
spectral functors represented by the virtual fundamental chain.  In
fact, in our application, a $\Pi$-equivariant $\cP$-cellular diagram
$\delta$ is presented as a sequence of zig-zags of $\Pi$-equivariant
spectral functors between spectral categories with $\Pi$-actions
\begin{equation}
\begin{tikzcd}
C^*_{\rel \partial}(\cP; \Omega \bS) \ar[r] & \cC_1 & \ar[l,"\htp"] \cC_2
\ar[r] & \cdots & \ar[l,"\htp"] \cC_n \ar[r] & \kmod
\end{tikzcd}
\end{equation}
such that all of the functors pointing left are equivalences of spectral
categories, i.e., functors which are homotopically fully-faithful and
essentially surjective.  (Such functors are typically referred to as
DK-equivalences; see Definition~\ref{defn:DK-equivalence} for a
review.)  A simplifying aspect of our situation is that all of the
categories $\aC_i$ have strict actions of $\Pi$ and the functors are
strictly equivariant, with the exception of $\kmod$ (as discussed
above) and of hence the last functor $\cC_n \to \kmod$.  

It is possible to formally rectify such a zig-zag in order to work
with a representative that is an actual functor.  Since our
constructions use this functor in a bar construction, we can do this
using the interpretation of spectral functors as bimodules.  (We 
explain this perspective in detail in Appendix~\ref{sec:spectral}, but
we give an indication of the strategy here.)

Suppose that we have a zig-zag
\begin{equation}
\begin{tikzcd}
\cC_1 & \cC_2 \ar[r,"G"] \ar[l,swap,"F"] & \cC_3.
\end{tikzcd}
\end{equation}
Associated to the spectral functor $F$ is the $\cC_1 \sma
\cC_2^{\op}$-module specified by the assignment 
\begin{equation}
(x,y) \mapsto \cC_1(Fy,x)
\end{equation}
and associated to the spectral functor $G$ is the $\cC_2 \sma
\cC_3^{\op}$-module specified by the assignment 
\begin{equation}
(x,y) \mapsto \cC_3(y,Gx)
\end{equation}
Note that since $F$ is a DK-equivalence, there is an equivalence 
\begin{equation}
\cC_1(Fy,x) \htp \cC_1(Fy, Fz) \htp \cC_2(y,z)
\end{equation}
for $Gz \htp x$.  Then the tensor product of these bimodules (which we
can compute using the bar construction) yields a $\cC_1 \sma
\cC^{\op}_3$-module that represents the composite functor.

For the remainder of this section, we will tacitly suppress the issue
of rectifying a zig-zag and write in terms of a representative that is
an honest functor.

\subsubsection{Homotopy type of sub-level sets}
\label{sec:homotopy-type-sub}

We now begin to construct the homotopy type associated to a
$\Pi$-equivariant cellular diagram.  We do this in terms of the
filtration of $\cP$ by the action map $\cA$.

\begin{defin}
For each real number $\lambda$, let $\cP_\lambda$ denote the quotient
of $\cP$ by the elements $p$ such that $\cA(p) > \lambda$.  That is, 
\begin{equation}
\cP_\lambda = \{p \in \cP \,\mid\, \cA(p) \leq \lambda\} \coprod \{\infty\},
\end{equation}
where the partial order is inherited from $\cP$ and $\infty$ is a new
terminal object.
\end{defin}

Note that $\cP_{\lambda}$ no longer admits an action of $\Pi$;
instead, each element $\pi \in \Pi$ induces a functor from
$\cP_{\lambda}$ to $\cP_{\lambda + \cA(\pi)}$; we explain this in more
detail below.

\begin{rem}
The reader should have in mind the following analogy: the homotopy
type we seek to define corresponds to a Laurent series ring
$\bbk((z))$, with degree of monomials corresponding to action. In
order to construct this, one may start with $z^{n} \bbk[z^{-1}]$, and
take an inverse limit ($n$ plays the role of $\lambda$).
\end{rem}

The filtration on $\cP$ induces one on $\Pi$-equivariant
$\cP$-cellular diagrams.

\begin{defin}
Given a $\Pi$-equivariant $\cP$-cellular diagram
$\delta$, we define the spectral functor 
\begin{equation}
    \delta_\lambda \co C^*_{\rel \partial}(\cP_\lambda; \Omega \bS)
    \to \kmod
  \end{equation} 
to be the restriction of $\delta$ to $\cP_\lambda$, where we stipulate
that $\delta_\lambda(\infty) = \ast$.
\end{defin}

Next, we will define the homotopy type associated to $\delta_\lambda$
in terms of a suitable derived tensor product, which we compute
explicitly via a two-sided bar construction.  For each $\lambda$, we
have the following distinguished module over $C^*_{\rel
  \partial}(\cP_\lambda; \Omega \bS)$.

\begin{defin}
Let $\bS_\lambda$ denote the (contravariant) functor
\begin{equation}
 C^*_{\rel \partial}(\cP_\lambda; \Omega \bS) \to \Sp
\end{equation}
which is uniquely specified by the assignment of $\bS$ to the
terminal object, and the zero-object $\ast$ to every other object.
The structure map 
\begin{equation}
C^*_{\rel \partial}(\cP_\lambda; \Omega \bS)(\ast,\ast) = \bS
\to F(\bS, \bS)
\end{equation}  
is the canonical map given as the adjoint of the multiplication map
$\bS \sma \bS \to \bS$.
\end{defin}

We briefly recall the details of the two-sided bar construction for
spectral categories and functors; some technical results we require
are discussed in Appendix~\ref{sec:spectralbar}.  

\begin{defin}
For a spectral category $\cC$, a left $\cC$-module (i.e., covariant
functor) $F$, and a right $\cC$-module (i.e., contravariant functor)
$G$, the two-sided bar construction $B(G,\cC,F)$ is the spectrum given
as the geometric realization of the simplicial spectrum $B_\bullet(G,
\cC, F)$ defined as
\begin{equation}
[k] \mapsto \bigvee_{c_0, c_1, \ldots, c_k} G(c_k) \sma
\cC(c_{k-1},c_k) \sma \ldots \sma \cC(c_0,c_1) \sma F(c_0),
\end{equation}
with face maps induced by the compositions and module actions and
degeneracies by the units of $\cC$.
\end{defin}

In order for this to have the correct homotopy type, we need to ensure
that the smash products compute the derived smash product; for this
purpose, we use the notion of a {\em pointwise-cofibrant} spectral
category, which is a spectral category $\aC$ where the mapping spectra
$\aC(x,y)$ are cofibrant spectra for each pair of objects $x$ and $y$.
As we review in Appendix~\ref{sec:spectral}, we can replace $\cC$ by a
pointwise-cofibrant category $\cC'$ equipped with a DK-equivalence
$\aC' \to \aC$ and at least one of the pullbacks of $G$ and $F$ to
$\aC'$-modules by a pointwise-cofibrant module.  There is a further
subtlety insofar as we need these cofibrant replacements to preserve
the action of $\Pi$; we discuss how to accomplish this in
Appendix~\ref{sec:equicof} when the action of $\Pi$ on $\aC$ is
strict, which suffices for our applications.

\begin{defin}\label{defn:homotopytype}
Given $\lambda$, we define the homotopy type of $\delta_\lambda$ to be
  \begin{equation}
    |\delta_\lambda| = B(\bS_\lambda,C^*_{\rel \partial}(\cP_\lambda;
    \Omega \bS), \delta_\lambda),
  \end{equation}
  which is a model for the (derived) tensor product of functors 
\begin{equation}
\bS_\lambda \sma^L_{C^*_{\rel \partial}(\cP_\lambda; \Omega \bS)}
\delta_\lambda.
\end{equation}
We will also denote this bar construction by $B(\bS_\lambda,
\delta_\lambda)$, and refer to it as the geometric realization of
$\delta_\lambda$.
\end{defin}

\begin{rem}
Definition~\ref{defn:homotopytype} should be thought of as an analogue
of the iterated cone construction that arises in the Cohen-Jones-Segal
approach to the Floer homotopy type.  In their setting, they construct
an explicit resolution of the module corresponding to
$\delta_{\lambda}$ in order to derive the tensor product of functors.
There are also certain differences since they fix $\cP = \bZ$ and work
with chains rather than cochains.  Roughly speaking, in the context of
chains, we would consider the stabilized flow category
$\Sigma^\infty_+ \hat{\cP}^+_\lambda$, where $\cP^+_\lambda$ is the
full subcategory of $\cP$ consisting of elements of action $< \lambda$
and an adjoined ``terminal'' object which receives maps from all
elements except the extremal ones (i.e., those that have no outgoing
maps).  We can then assign a homotopy type to functors $\delta \colon
\Sigma^\infty_+ \hat{\cP} \to \Sp$ by restricting to a functor
$\delta_\lambda \colon \Sigma^{\infty}_+ \hat{\cP}^+_\lambda \to \Sp$ (defined
to be $\ast$ on the new terminal object) by passage to the homotopy
colimit.  This formulation recovers the Cohen-Jones-Segal homotopy
type when $\cP = \bZ$. 
\end{rem}

The action of $\Pi$ on $\cP$, $C^*_{\rel \partial}(\cP_\lambda; \Omega
\bS)$, and $\delta_\lambda$ assemble to equip a suitable inverse limit
over $\lambda$ of the spectra $|\delta_\lambda|$ with an action by the
completed twisted group ring $\Sigma^{\deg} \bbk((\pi))$.  We will
describe this action in detail in Section \ref{sec:homotopy-type-as}
below; for the time being, we will consider a fixed $\lambda$ and
ignore the equivariant structure.

We now explain how to compute the geometric realization
$|\delta_\lambda|$ in simple examples, which ultimately build up to
the computation of the associated graded spectrum in the general
case.  Since $\bS_{\lambda}$ is trivial except at the terminal object,
we can rewrite the bar construction as follows.

\begin{lem}
The bar construction $B(\bS_\lambda,C^*_{\rel \partial}(\cP_\lambda;
    \Omega \bS), \delta_\lambda)$ is the simplicial spectrum whose $k$-simplices are
    \begin{equation}
       \bigvee_{\substack{p_0 \leq p_1 \leq \ldots p_{k-1} \\ \cA(p_{k-1})
  \leq \lambda}}
C^*_{\rel \partial}(\infty, p_k) \sma 
C^*_{\rel \partial}(p_{k-1}, p_k) \sma \ldots
\sma C^*_{\rel \partial}(p_0, p_1) \sma \delta_\lambda(p_0),
\end{equation}
where the $0$-simplices are the trivial spectrum $\ast$.  The face map
$d_0$ is the trivial map, the face maps $d_1, \ldots, d_{n-1}$ are
given by the composition in $C^*_{\rel \partial}(\cP_\lambda; \Omega
\bS)$, and the face map $d_n$ is determined by the module structure on
$\delta$.
\end{lem}

\begin{proof}
The identification of the $0$-simplices follows from the fact that 
\begin{equation}
B_0(\bS_\lambda, C^*_{\rel \partial}(\cP_\lambda; \Omega \bS),
\delta_\lambda) = \bigvee_{p \in \cP_\lambda} \bS_\lambda(p) \sma
\delta_\lambda(p).
\end{equation}
Since $\bS_\lambda(p) = \ast$ unless $p = \infty$ in which case it is
equivalent to $\bS^{\mfib}$, this coproduct collapses to
$\delta_\lambda(\infty)$, which is also $\ast$ by definition.
Therefore, the $0$-simplices of this bar construction are always
homeomorphic to $\ast$.  The description of the higher simplices
follows from the definitions in the same fashion.
\end{proof}

\begin{lem}
The resulting bar construction is split (see
Section~\ref{sec:spectralbar}) with the splitting given by the 
subspectra 
\begin{equation}
  \bar{B}_\bullet(\cP,\lambda) = \bar{B}_\bullet(\bS_\lambda, C^*_{\rel
  \partial}(\cP_\lambda; \Omega \bS), \delta_\lambda),
\end{equation}
obtained by restricting the coproducts in the expression for the
two-sided bar construction to sequences of objects in $\cC$ without
repeated elements. \qed
\end{lem}

This property implies that the geometric realization of the simplicial
spectrum $B_\bullet(\cP,\lambda)$ can then be computed in terms of the
filtered colimit of the skeleta computed via the pushouts
\begin{equation}\label{eq:skpush}
  \begin{tikzcd}
    \bar{B}_{n}(\cP,\lambda) \times \partial \Delta_n \ar[d] \ar[r] &
\sk_{n-1} B(\cP,\lambda) \ar[d] \\
\bar{B}_{n}(\cP,\lambda) \times \Delta_n \ar[r] & \sk_n B(\cP,\lambda).
  \end{tikzcd}
\end{equation}
Here the top horizontal map is induced by the face maps.

We now compute a series of simple examples.  In these cases, since
$\cP_\lambda$ is finite, there are no nondegenerate simplices in
degree larger than the maximal length of a  totally ordered subset of
$\cP_\lambda$.

\begin{example} \label{ex:compute_filtered_homotopy_type}
\hspace{5 pt}
\begin{enumerate}
\item If $\cP_\lambda$ consists only of the object $\infty$ (e.g.,
  $\cP$ has no objects of action $\leq \lambda$), then
  $|\delta_\lambda|$ is homeomorphic to a point.

\item \label{item:filtered_homotopy_one_object} If $\cP_\lambda$ has a single object $p < \infty$, then the
  $1$-simplices of $|\delta_\lambda|$ are given by
\begin{equation}
\bS_{\lambda}(\infty) \sma C^*_{\rel \partial}(\cP_\lambda; \Omega
\bS)(p,\infty) \sma \delta_{\lambda}(p) \htp \bS \sma \Omega \bS \sma
\delta(p) \htp \Omega \delta(p),
\end{equation}
as all other terms in the wedge are $\ast$.  Moreover, there are no
higher simplices in the reduced complex.  The geometric realisation
thus is equivalent to the quotient of $\Delta^1_+ \sma (\Omega
\delta(p))$ which identifies $\{0,1\} = \partial \Delta^1$ with the
basepoint in $\Omega \delta(p)$; it is therefore equivalent to the
suspension $\Sigma \Omega \delta(p) \htp \delta(p)$.

\item If $\cP_\lambda$ consists of objects $p < q < \infty$, then the
  $1$-simplices have nontrivial terms corresponding to $(p,\infty)$
  and $(q,\infty)$.  Furthermore, by
  Remark~\ref{rem:acyclic_morphism_spaces}, the term corresponding to
  $(p,\infty)$ is contractible.  As a consequence the $1$-simplices are
  equivalent to  
\begin{equation}
\bS_{\lambda}(\infty) \sma C^*_{\rel \partial}(\cP_\lambda; \Omega
\bS)(q,\infty) \sma \delta_{\lambda}(q) \htp \bS \sma \Omega \bS \sma
\delta(q) \htp \Omega \delta(q).
\end{equation}
The nondegenerate $2$-simplices are only non-vanishing for the triple
$(p,q,\infty)$; calculating, we obtain
\begin{align*}
\bS_{\lambda}(\infty) \sma C^*_{\rel \partial}(\cP_\lambda; \Omega
\bS)(q,\infty) \sma C^*_{\rel \partial}(\cP_\lambda; \Omega \bS)(p,q)
\sma \delta(p) \\
\htp \bS \sma \Omega \bS \sma \Omega \bS \sma \delta(p)
\htp \Omega^2 \delta(p).
\end{align*}
The first face map lands in $\ast$, the second in the contractible
nondegenerate $1$-simplex, and the third is induced by structure map 
\begin{equation}
C^*_{\rel \partial}(\cP_\lambda)(p,q) \sma \delta(p) \to \delta(q).
\end{equation}
Equation~\eqref{eq:skpush} shows that the geometric realisation is
therefore equivalent to the cofiber of the structure map $\Omega
\delta(p) \to \delta(q)$.  To be more precise, the geometric
realization can be described as the spectrum obtained by the quotient of
$\Delta^2_+ \sma (\Omega^2 \delta(p))$ which identifies the first face
with $\Delta^1_+ \sma (\Omega \delta(p))$ via the (loop of) the
structure map and the other two faces with $\ast$, and identifies the
two vertices of the first face with $\ast$.

If the structure map was the identity, then we are looking at the
based tensor with the quotient of $\Delta^2$ that identifies the horn
$\Lambda^2_0$ to $\ast$; this is contractible, and so the result is
$\ast$.  When the structure map is null-homotopic, we conclude that
there is an equivalence
\begin{equation}
  \label{eq:homotopy_type_two_cells_wedge}
|\delta_\lambda| \htp \delta(p) \vee \delta(q).  
\end{equation}

\item If $\cP_\lambda$ consists of objects $(p,q_1, q_2, \infty)$ with
  $p < q_1 < \infty$ and $p < q_2 < \infty$, then the
  $1$-simplices have non-vanishing contributions potentially from the
  terms corresponding to $(q_1,\infty)$, $(q_2,\infty)$, and $(p,\infty)$.
  Computing as above, the contributions to the $1$-simplices from the
  first two terms are thus the wedge 
\begin{equation}
\Omega \delta(q_1) \vee \Omega \delta(q_2).
\end{equation}
For $(p,\infty)$, since $\kappa(p,\infty)$ is the wedge $[0,1] \vee
[0,1]$ (where the intervals have basepoint $0$) and $\partial
\kappa(p,\infty)$ is the union of the two endpoints $\{1\}$, we find
that 
\begin{equation}
C^*_{\rel \partial}(\cP_\lambda; \Omega \bS)(p, \infty) \sma \delta(p)
\htp \Omega^2 \delta(p).
\end{equation}
That is, in total the one-simplices are the wedge
\begin{equation}
\Omega \delta(q_1) \vee \Omega \delta(q_2) \vee \Omega^2 \delta(p).
\end{equation}

The $2$-simplices have contributions from the triples $(p,q_1, \infty)$
and $(p, q_2, \infty)$; thus the $2$-simplices are the wedge
\begin{equation}
\Omega^2 \delta(p) \vee \Omega^2 \delta(p).
\end{equation}

Therefore, the homotopy type of the geometric realization is given by
the quotient of the wedge
\begin{equation}
\begin{aligned}
(\Omega \delta(q_1) \sma \Delta^1_+) &\vee (\Omega \delta(q_2) \sma
\Delta^1_+) \vee (\Omega^2 \delta(p) \sma \Delta^1_+) \\ &\vee 
(\Omega^2 \delta(p) \sma \Delta^2_+) \vee (\Omega^2 \delta(p) \sma
\Delta^2_+)
\end{aligned}
\end{equation}
where the two copies of $\Omega^2 \delta(p) \sma \Delta^2_+$ are
attached to each other along weak equivalences on a face to $\Omega^2
\delta(p)$ in the $1$-simplices and to $\Omega \delta(q_1) \sma
\Delta^1_+$ and $\Omega \delta(q_2) \sma \Delta^1_+$ via the structure
maps along the other faces.  This is equivalent to the homotopy
pushout of the diagram   
\begin{equation}
\begin{tikzcd}
\Sigma \Omega \delta(q_1) & \Sigma \Omega^2 \delta(p) \ar[l] \ar[r] &
\Sigma \Omega \delta(q_2),
\end{tikzcd}
\end{equation}
where the maps in the pushout are determined by the structure maps of
$\delta$.

Now assume that the structure maps are null-homotopic.  Since the
connecting map in the cofiber sequence 
\begin{equation}
\begin{tikzcd}
\Sigma \Omega^2 \delta(p) \ar[r] &
\Sigma \Omega \delta(q_1) \vee \Sigma \Omega \delta(q_2) \ar[r] &
\relax|\delta\relax| \ar[r] & \Sigma^2 \Omega^2 \delta(p)
\end{tikzcd}
\end{equation}
is null homotopic, the homotopy type is equivalent to 
\begin{equation}
\Sigma \Omega \delta(q_1) \vee \Sigma \Omega \delta(q_2) \vee \Sigma^2
\Omega^2 \delta(p) \htp \delta(q_1) \vee \delta(q_2) \vee \delta(p).
\end{equation}
We can see this more explicitly as follows.  For each
$\Omega^2 \delta(p)$, the face maps are trivial except for the map to
$\Omega^2 \delta(p)$ which is a weak equivalence.  Therefore, the
contribution to the realization from the $2$-simplices is the based
tensor of $\Omega^2 \delta(p)$ with the gluing of two copies of
$\Delta^2$ along a face, with the boundary collapsed to $\ast$.
\end{enumerate}
\end{example}

More generally, the descriptions given in
Equation~\eqref{eq:homotopy_type_two_cells_wedge} can be extended by
considering the filtration on $|\delta_\lambda|$ induced by $\cA$. The
key computation required to analyze this filtration is given by the
following: 

\begin{lem}\label{lem:justp}
Let $\delta_{\lambda,p}$ denote the module on $\cP_{\lambda}$ which
has value $\delta(p)$ at $p$ and $\ast$ everywhere else.  Then the bar
construction 
\begin{equation}
B(\bS_{\lambda}, C^*_{\rel \partial}(\cP_\lambda; \Omega \bS), \delta_{\lambda, p}) 
\end{equation}
is naturally equivalent to $\delta(p)$.
\end{lem}

\begin{proof}
The bar construction in question has $k$-simplices
\begin{equation}
\bigvee_{p, q_0, q_1, \ldots, q_{k-1}, \infty}
\delta(p) \sma C^*_{\rel \partial}(p,q_0)
\sma C^*_{\rel \partial}(q_0,q_1) \sma \ldots
\sma C^*_{\rel \partial}(q_{k-1},\infty),
\end{equation}
for a totally ordered subset $Q = (p, q_0, \ldots, q_{k-1})$ of $\cP$.
For the nondegenerate simplices, we have $p < q_0 < \ldots < q_{k-2}$.
The $0$-simplices are $\ast$, and the $1$-simplices are $C^*_{\rel
  \partial}(\cP_{\lambda}, \Omega \bbk)(p,\infty) \sma \delta(p)$.

To analyze this bar construction, observe that since the action maps
\begin{equation}
\delta_{\lambda, p}(p) \sma C^*_{\rel \partial}(\cP_\lambda;
\bS)(p,q) \to \delta_{\lambda, p}(q)
\end{equation}
are all the trivial map and hence the face map $d_0$ is trivial, there
is a natural equivalence
\begin{equation}
B(\bS_{\lambda}, C^*_{\rel \partial}(\cP_\lambda; \Omega \bS),
\delta_{\lambda, p}) \cong \delta_p \sma B(\bS_{\lambda}, C^*_{\rel
  \partial}(\cP_\lambda; \Omega \bS), \bS_p),
\end{equation}
where here $\bS_p$ abusively denotes the covariant analogue of
$\bS_\lambda$.

First, we analyze the case when the restriction of $\cP_{\lambda}$ to
elements over $p$ is isomorphic to the finite set $\{0,1,\ldots,m\}$
with the standard order, which we will write as $\{p,q_0, q_1, \ldots,
q_{m-2}, \infty\}$.  In this case, the $k$-simplices of the bar
construction are contractible except when $k = m$, where 
\begin{equation}
C^*_{\rel \partial}(p,q_0)
\sma C^*_{\rel \partial}(q_0,q_1) \sma \ldots
\sma C^*_{\rel \partial}(q_{m},\infty) \htp \Omega^{m} \bS.
\end{equation}
The evident map to the simplicial object obtained by setting the
$\ell$-simplices to $\ast$ for $\ell < m$ clearly induces a weak
equivalence on geometric realizations, and so the geometric
realization of the bar construction is $\Sigma^m \Omega^{m} \bS \htp 
\bS$.

We can interpret this calculation as follows, which is useful for
understanding the general case: the projection induces the canonical
map 
\begin{equation}\label{eq:cochain-projection}
C^*_{\rel \partial}(\cP_{\lambda}; \Omega \bS)(p,q) \to
C^*(\cP_{\lambda}; \Omega \bS)(p,q).
\end{equation}
Since the cubical complex $\kappa(p,q)$ contracts to the cone point
for any $p$ and $q$, the latter is equivalent to $\bS$, and the
projection map is a natural equivalence.

We now consider general posets $\cP_{\lambda}$.  Consider the
simplicial spectrum $C^*(\cP_\lambda)_{\bullet}$ with $k$-simplices 
\begin{equation}
\bigvee_{p, q_0, q_1, \ldots, q_{k-1}, \infty}
C^*(p,q_0) \sma C^*(q_0,q_1) \sma \ldots
\sma C^*(q_{k-1},\infty),
\end{equation}
for a totally ordered subset $Q = (p, q_0, \ldots, q_{k-1})$ of
$\cP_\lambda$.  The $0$-simplices are $\ast$, and the $1$-simplices
are $C^*(\cP_{\lambda}, \Omega \bbk)(p,\infty)$.  The structure maps
are induced by the composition.  The geometric realization of
$C^*(\cP_\lambda)$ is equivalent to $\bS$; since the cubes are
contractible, the $k$-simplices are equivalent to $\Omega^k \bS$, and
the structure maps induce an equivalence to the bar construction on
the cochains of $S^1$, i.e., $\bS$.

Since the projection maps of Equation~\eqref{eq:cochain-projection}
are compatible with the composition, they induce a simplicial map
\begin{equation}
B_\bullet(\bS_{\lambda}, C^*_{\rel
  \partial}(\cP_\lambda; \Omega \bS), \bS_p) \to
C^*(\cP_\lambda)_{\bullet}.
\end{equation}
We will now argue that this map is an equivalence.  By construction,
each cubical complex $\kappa(p,q)$ can be written as 
\begin{equation}
\kappa(p,q) = \hocolim_{Q \in 2^{p,q}} \kappa(Q),
\end{equation}
where $Q$ is a totally ordered subset and the inclusions are
determined by setting coordinates to $0$.  The boundary $\partial
\kappa(p,q)$ can analogously be written as a homotopy colimit over the
same diagram
\begin{equation}
\partial \kappa(p,q) = \hocolim_{Q \in 2^{p,q}} \partial \kappa(Q),
\end{equation}
This implies that we can decompose the cochains and relative cochains
as homotopy limits and express the projection map as
\begin{multline}
C^*_{\rel \partial}(\kappa(p,q); \Omega \bS) \cong
\holim_{Q \in 2^{p,q}} C^*_{\rel \partial}(Q; \Omega \bS)
\\ \to
C^*(\kappa(p,q); \Omega \bS) \cong \holim_{Q \in 2^{p,q}} C^*(Q;
\Omega \bS).
\end{multline}
These decompositions are compatible with the simplicial structure,
since $\kappa(p,q)$ is a subcomplex of $\kappa(p,q')$ for $q' > q$.
By the finiteness hypothesis on $\cP$, we can pull these homotopy
limits outside the geometric realization; the result now follows from
the calculation for a cube.

\end{proof}

We now return to study the filtration on $|\delta_\lambda|$: observe
that if we restrict to objects with action $> n$ or $\geq n$ for $n
\leq \lambda$, we obtain full subcategories of $\cP_{\lambda}$ which
we will denote $\cP_{\lambda,  >n}$ and $\cP_{\lambda, \geq n}$.

\begin{defin}\label{defn:bounded-homotopy-type}
Let $|\delta_{\lambda, >n}|$ be the bar construction 
\begin{equation}
B(\bS_\lambda,C^*_{\rel \partial}(\cP_{\lambda,>n};
    \Omega \bS), \delta_{\lambda, >n}),
\end{equation}
where we are restricting to the poset $\cP_{\lambda, >n}$.  We define  
$|\delta_{\lambda, \geq n}|$ analogously to be the bar construction
\begin{equation}
B(\bS_\lambda,C^*_{\rel \partial}(\cP_{\lambda, \geq n};
    \Omega \bS), \delta_{\lambda, \geq n}),
\end{equation}
specified by restricting to the poset $\cP_{\lambda, \geq n}$.
\end{defin}

There is a natural inclusion $|\delta_{\lambda, >n}| \to
|\delta_{\lambda, \geq n}|$.  Using this, we can identify the
``associated graded'' piece corresponding to a fixed $n \leq
\lambda$, under a hypothesis on the action filtration on $\cP$ that
holds in our examples.

\begin{prop} \label{prop:computed_associated_graded}
If the subset of $\cP$ of elements of action $n$ is discrete (i.e.\ $p
\leq q$ and $\cA(p) = \cA(q) = n$ implies that $p = q$), then there is
an equivalence
\begin{equation}
|\delta_{\lambda, \geq n}| / |\delta_{\lambda, >n}| \htp
\bigvee_{\cA(p) = n} \delta(p). 
\end{equation}
\end{prop}

\begin{proof}
Writing 
\begin{align*}
|\delta_{\lambda, \geq n}| &= B(\delta_{\lambda, \geq n}, C^*_{\rel
  \partial}(\cP_{\lambda, \geq n}, \Omega \bS), \bS_{\lambda}) \\&\cong
B(\delta_{\lambda, \geq n}, C^*_{\rel \partial}(\cP_{\lambda, \geq n},
\Omega \bS), C^*_{\rel \partial}(\cP_{\lambda, \geq n}, \Omega \bS))
\sma_{C^*_{\rel \partial}(\cP_{\lambda, \geq n}, \Omega \bS)}
\bS_{\lambda}
\end{align*}
and analogously for $|\delta_{\lambda, >n}|$.  Next, observe that 
\begin{equation}
B(\delta_{\lambda, > n}, C^*_{\rel \partial}(\cP_{\lambda, > n}, \Omega \bS), \bS_{\lambda}) 
=
B(\delta_{\lambda, > n}, C^*_{\rel \partial}(\cP_{\lambda, \geq n},
\Omega \bS), \bS_{\lambda}),
\end{equation}
since any chain in $\cP$ which passes through $q$ such that $\cA(q) =
n$ will give rise to a trivial contribution in the righthand side, as
$\delta_{\lambda > n}$ evaluated on the minimal element of that chain
will have action $\leq n$ and so be $\ast$.

Therefore using the fact that colimits commute, the quotient in
question can be computed in terms of the induced quotient of functors
\begin{equation}
\begin{tikzcd}
B(\delta_{\lambda, > n}, C^*_{\rel \partial}(\cP_{\lambda, > n},
\Omega \bS), C^*_{\rel \partial}(\cP_{\lambda, > n}, \Omega \bS))
\ar[d] \\
B(\delta_{\lambda, \geq n}, C^*_{\rel \partial}(\cP_{\lambda, \geq n},
\Omega \bS), C^*_{\rel \partial}(\cP_{\lambda, \geq n}, \Omega \bS)).
\end{tikzcd}
\end{equation}
Now, as functors, for $p$ such that $\cA(p) > n$, we have that
\begin{equation}
B(\delta_{\lambda, \geq n}, C^*_{\rel \partial}(\cP_{\lambda, \geq n},
\Omega \bS), C^*_{\rel \partial}(\cP_{\lambda, \geq n}, \Omega \bS))
\end{equation}
and
\begin{equation}
B(\delta_{\lambda, > n}, C^*_{\rel \partial}(\cP_{\lambda, > n},
\Omega \bS), C^*_{\rel \partial}(\cP_{\lambda, > n}, \Omega \bS))
\end{equation}
coincide, and so the quotient vanishes for $p$ such that $\cA(p) > n$.
As a consequence, we can reduce to computing the bar construction with
respect to the module $\delta_{\lambda = \cA(p)}$ which coincides with
$\delta$ on $p$ such that $\cA(p) = n$ and is $\ast$ elsewhere.

Under our assumption, there are no arrows between elements $p$ and $q$
of action equal to $n$ if they are distinct. This implies that the bar
construction for this module splits as a wedge over $p$ such that
$\cA(p) = n$.  The result then follows from Lemma~\ref{lem:justp}.
\end{proof}

\begin{rem}
In fact, we can drop the assumption of Proposition
\ref{prop:computed_associated_graded}, at the cost of considering a
filtration indexed by $\cP$ itself; since the applications we consider
do not require this generality, we work with the above simplified
setting. 
\end{rem}

The layers of this filtration on $|\delta_{\lambda}|$ are attached by
the module structure maps of $\delta_{\lambda}$; when these are
null-homotopic, the filtration splits.  To see this, we appeal to
the notion of the Kan suspension of a
simplicial spectrum~\cite{Kan1963}, which we will denote by
$\tilde{\Sigma}$.  

\begin{defin}
For a simplicial spectrum $X_\bullet$, the $k$-simplices of the Kan
suspension $\tilde{\Sigma} X_\bullet$ are determined by the formula 
\begin{equation} \label{eq:Kan_suspension_formula}
(\tilde{\Sigma} X)_k = X_0 \vee X_1 \vee \ldots \vee X_{k-1}.
\end{equation}
On the summand $X_\ell$ in $(\tilde{\Sigma} X)_k$, we have
degeneracies specified by the formulas
\begin{equation}
s_i = 
\begin{cases}
s_{i-(k-\ell)} \colon X_\ell \to X_{\ell+1} \qquad  &i \geq k-\ell \\
\id \colon X_\ell \to X_{\ell} \qquad &i<k-\ell
\end{cases}
\end{equation}
and face maps
\begin{equation}
d_i =
\begin{cases}
\id \colon X_\ell \to X_\ell \qquad &i<k-\ell \\
d_{i-(k-\ell)} \colon X_\ell \to X_{\ell-1} \qquad &i \geq k-\ell,
\ell>0 
\end{cases}
\end{equation}
and as the trivial map for $i=k$, $\ell=0$.
\end{defin}

The Kan suspension has the effect of shifting the nondegenerate
simplices up a simplicial degree.  Moreover, it models the ordinary
suspension, in the following sense: for cofibrant simplicial spectra
$X$, there is a natural equivalence
\begin{equation}
\Sigma |X| \to |\tilde{\Sigma} X|
\end{equation}
of spectra.

\begin{prop}
Assume that the subset of $\cP$ of elements of action $n$ is discrete
(i.e.\ $p \leq q$ and $\cA(p) = \cA(q) = n$ implies that $p = q$).
When the action maps $C^*_{\rel \partial}(p,q) \sma
\delta_{\lambda}(p) \to \delta_{\lambda}(q)$ are null-homotopic for
all $p$ such that $\cA(p) = n$, the filtration splits and we find
\begin{equation}
|\delta_{\lambda, \geq n}| \htp |\delta_{\lambda, >n}| \vee
\bigvee_{\cA(p) = n} \delta(p).
\end{equation}
\end{prop}

\begin{proof}
The connecting map in the cofiber sequence
\begin{equation}
|\delta_{\lambda, >n}| \to |\delta_{\lambda, \geq n}| \to
\bigvee_{\cA(p) = n} \delta(p) \to \Sigma|\delta_{\lambda, >n}|
\end{equation}
induces the composite
\begin{equation}
|\delta_{\lambda, \cA(p) = n}| \to \Sigma |\delta_{\lambda, >n}| \to
|\tilde{\Sigma} \delta_{\lambda, >n}| 
\end{equation}
to the Kan suspension.  The composite is the geometric realization of
a simplicial map, which can be described as follows.  For $k \geq 1$,
each nondegenerate $k$-simplex in $|\delta_{\lambda,
  \cA(p)=n}|$ corresponds to a totally-ordered subset $\{p, q_1,
\ldots, q_{k-1}, \infty\}$ where $\cA(q_1) > \cA(p) = n$.  For $k \geq
2$, each nondegenerate simplex in $\tilde{\Sigma} \delta_{\lambda,
  >n}$ corresponds to a totally-ordered subset $(q_1, \ldots, q_{k-1},
\infty)$ with $\cA(q_1) > n$.  The map on the nondegenerate
$k$-simplices is then the map induced by the action maps $C^*_{\rel
  \partial}(p,q_1) \sma \delta(p) \to \delta(q_1)$ along with the
smash product of the identity maps, except on the $1$-simplices where
it is the collapse map.  This assignment is clearly compatible with
the simplicial structure.  The map on the degenerate simplices is
analogous.  Although we have not assumed that the null-homotopies of
the action maps are compatible, because the face map corresponding to
the module action in $|\delta_{\lambda, \cA(p)=n}|$ is trivial, the
homotopies assemble to produce a null homotopy of the connecting map.
\end{proof}

Using the fact that 
\begin{equation}
|\delta_{\lambda, > n}| = \colim_{m < n} |\delta_{\lambda \geq m}|,
\end{equation}
we can conclude the following proposition.

\begin{prop}\label{prop:nullhomtype}
Suppose that all the composition maps in $\delta_\lambda$ other than
the identities are null homotopic.  Then we have an equivalence
 \begin{equation} \label{eq:homotopy_type_big_wedge}
  |\delta_\lambda| \htp \bigvee_{\cA(q) \leq \lambda} \delta(q).
\end{equation} \qed
\end{prop}

\subsubsection{Signpost: Morse and Floer theory}
\label{sec:signp-morse-theory}

It is useful at this stage to note that the constructions of the above section suffice to produce the Morse homotopy type of a function on a closed manifold, as discussed in  Appendix \ref{sec:morse-theory-and-homotopy-types}.  In this setting, discussed already by Cohen-Jones-Segal \cite{CohenJonesSegal1995}, the group $\Pi$ is trivial, the category $\cP$ is finite, and the morphisms in the flow category $\cM$ are topological manifolds with stratified boundary, arising as moduli spaces of flow lines of a Morse-Smale function. 

To briefly summarize the construction, recall that moduli spaces of
flow lines in Morse theory are oriented relative the positive-definite
subspace of the Hessian matrix at each critical point. We write
$\delta_p$ for the corresponding orientation line, which is a (graded)
rank-$1$ module over $\bbk$, and introduce the  spectral category $
C_*( \cM, \delta)$, whose objects are those of $\cM$, and whose
morphisms assign to a pair $(p,q)$: 
  \begin{equation}
 C_*( \cM, \delta)(p,q)  \equiv \cM(p,q)_+ \wedge \Hom_\bbk(\delta_p,\delta_q) .
  \end{equation}
There is a canonical functor
  \begin{equation} \label{eq:projection_map_chains_moduli_morse}
 C_*( \cM, \delta)(p,q)  \to \kmod,
  \end{equation}
given by projecting $\cM(p,q)$ to a point.
  
In Appendix \ref{sec:morse-theory-and-homotopy-types}, we explain (in a variant of the work of \cite{CohenJonesSegal1995}), that there is a zig-zag of equivalences
\begin{equation} \label{eq:zig-zag-chains_cochains_Morse}
C^*_{\rel \partial}(\cM; \Omega \bbk)   \leftarrow \cdots \to  C_*( \cM, \delta)
\end{equation}
which arises from an appropriate application of Poincar\'e duality
with coefficients in $\bbk$ (i.e., Spanier-Whitehead duality). In
principle, this requires that the flow category be oriented with
respect to $\bbk$, i.e., that the underlying manifolds admit
compatible orientations, but this is true for the sphere spectrum, and
hence holds for every spectrum.  Nonetheless, we limit the discussion
in Appendix \ref{sec:morse-theory-and-homotopy-types} to complex
oriented theories, since it allows us to appeal to the results of the
main part of the paper without modification (the only place they
appear is in the proof of
Lemma~\ref{lem:bimodule_equivalence_Morse_theory}).

Using the $2$-sided bar construction to compose the zig-zag in
Equation \eqref{eq:zig-zag-chains_cochains_Morse} with the functor in
Equation \eqref{eq:projection_map_chains_moduli_morse}, and the
pullback morphism $C^*_{\rel \partial}(\cP; \Omega \bbk)  \to
C^*_{\rel \partial}(\cM; \Omega \bbk)  $, where $\cP$ is the partially
ordered set of critical points, we obtain a functor 
  \begin{equation}
        C^*_{\rel \partial}(\cP; \Omega \bbk) \to \kmod.
      \end{equation}
Since $\cP$ is finite in this context, the truncation $\cP_\lambda$
does not depend on the choice of constant $\lambda$ as long as it is
sufficiently large. Applying the results of Section
\ref{sec:homotopy-type-sub} thus produces the \emph{Morse homotopy
  type}. In order to compare this to the classical homotopy type of
the underlying manifold, we develop special tools in Appendix
\ref{sec:morse-theory-and-homotopy-types}, that are variants of the
tools which we introduce in Section \ref{sec:maps-homotopy-types}
below, in order to compare Morse and Floer theory.

\subsection{The homotopy type as an inverse limit}
\label{sec:homotopy-type-as}

Returning to the general theory, the construction of $\cP_\lambda$ is evidently contravariantly
functorial in $\lambda$, i.e., for $\lambda_0 < \lambda_1$ there is a
canonical functor
\begin{equation}
\cP_{\lambda_1} \to \cP_{\lambda_0}
\end{equation}
given by sending all of the elements $p \in \cP_{\lambda_1}$ with
$\cA(p) > \lambda_0$ to the terminal object.  However, the description
of the functoriality of the associated homotopy types requires more
work.  In the bar construction, the attachment of cells for 
$\delta_p$ is controlled by totally ordered subsets of $\cP_{\lambda}$
with top entry $\infty$; changing the length of these totally ordered
subsets by collapsing shifts the simplicial degree in which these
attaching maps are represented.  To handle this, we need to consider a
suitable intermediate construction that interpolates the effects of
shifts of degree.

To explain our construction, we begin by working through the most
basic example:

\begin{example}
Let $\cP$ be a partially ordered set with two objects $p < q$ such
that $A(p) = 0$ and $A(q) = 1$.  We assume we are given a functor
\begin{equation}
\delta \co C^{*}_{\rel \partial}(\cP; \Omega \bS) \to \kmod.
\end{equation}
Then $\cP_{0}$ is the poset $p \to \infty$ and $\cP_1$ is the
poset $p \to q \to \infty$.  We will construct a representative of a
homotopy class of maps $|\delta_1| \to |\delta_0|$ in a way that
will extend to general $\cP$.

To do this, we consider the poset category $\cP_{0,1}$ with
non-identity arrows generated by the following diagram

\begin{equation}
  \begin{tikzcd}
    \infty_0 \ar[r] & \infty_1 \\
    p \ar[r] \ar[u]  & q. \ar[u]
  \end{tikzcd}
\end{equation}

We define a functor $\delta_{01}$ on the relative cochains of
$\cP_{0, 1}$, which is specified on objects as
\begin{equation}
\begin{cases}
\delta_{01}(p) = \delta_p, \\
\delta_{01}(q) = \delta_q, \\
\delta_{01}(\infty_0) = \delta_{01}(\infty_1) = \ast,
\end{cases}
\end{equation}
and on morphisms, $\delta_{01}$ is determined by $\delta$.  We define the
functor $\bS_{0,1}$ to vanish on all objects but $\infty_1$, where it
is $\bS$.  The homotopy type associated to $\delta_{01}$ is then given
by two-sided bar construction:

\begin{equation}
|\delta_{01}| \equiv B(\bS_{0,1},   C^{*}_{\rel \partial}(\cP_{0,1};
\Omega \bS) , \delta_{01}).
\end{equation}

There is an evident functor $\cP_{0,1} \to \cP_{1}$ specified on
objects by sending both $\infty_0$ and $\infty_1$ to $\infty$.  We
specify a map of bar constructions 
\begin{equation}
B_\bullet(\bS_{0,1},   C^{*}_{\rel \partial}(\cP_{0,1};
\Omega \bS) , \delta_{01})
\to
B_\bullet(\bS_1,   C^{*}_{\rel \partial}(\cP_1;
\Omega \bS) , \delta_1)
\end{equation}
as follows.
The pullbacks of the functors $\delta_1$ and $\bS_1$ to
$C^{*}_{\rel \partial}(\cP_{0,1}; \Omega \bS)$ induce natural
transformations 
\begin{equation}
\delta_{01} \Rightarrow \delta_1 \qquad\textrm{and}\qquad \bS_{01}
\Rightarrow \bS_1
\end{equation}
as functors on $\cP_{0,1}$.  There is a functor
\begin{equation}
C^*_{\rel \partial}(\cP_{0,1}; \Omega \bS)(p,\infty_1) \to C^*_{\rel
  \partial}(\cP_{1}; \Omega \bS)(p,\infty)
\end{equation}
determined by the restriction map
\begin{equation}
 C^*_{\rel \partial}(\kappa^{p\infty_0\infty_1} \vee
 \kappa^{pq\infty_1}; \Omega \bS) \to  C^*_{\rel \partial}(
 \kappa^{pq\infty}; \Omega \bS). 
\end{equation}
We map $C^*_{\rel \partial}(p,\infty_2)$ to $C^*_{\rel \partial}(p,
\infty)$ using the restriction map and we map the composite $C^*_{\rel
  \partial}(p,\infty_1) \sma C^*_{\rel \partial}(\infty_1, \infty_2)$
to the degenerate simplex $C^*_{\rel \partial}(p,\infty) \sma
C^*_{\rel \partial}(\infty, \infty)$ using the composite of the
composition and the restriction map.  This is clearly a simplicial
map, and therefore we have a map of geometric realizations
\begin{equation}\label{eq:zig}
  |\delta_{01}| \to |\delta_{1}|.  
\end{equation}
We claim that the map in Equation~\eqref{eq:zig} is an equivalence.
Observe that by construction it is compatible with the filtration by
action, and so we have a commutative diagram of cofiber sequences
\begin{equation}
\begin{tikzcd}
\relax |\delta_{01,>0}| \ar[r] \ar[d] & \relax|\delta_{01, \geq 0}| = |\delta_{01}|
\ar[r] \ar[d] & \ar[d] \delta_p \\
\relax|\delta_{1, > 0}| \ar[r] & \relax|\delta_{1, \geq 0}| = |\delta_1| \ar[r] &
\relax\delta_p  
\end{tikzcd}
\end{equation}
where the vertical maps are induced by the comparison map.  Since
the outer maps are equivalences, the middle map must be as well.

Let $\cP'_0$ denote the subcategory of $\cP_{0,1}$ with objects $p \to
\infty_0 \to \infty_1$. We write $\delta_{01}^0$ for the functor
induced by $\delta$, and $\bS_1$ for the functor induced by
$\bS_{0,1}$.  We define
\begin{equation}
|\delta_{01}^0| = B(\bS_1, C^*_{\rel \partial}(\cP'_0, \Omega \bS),
\delta_{01}^0).
\end{equation}

Collapsing $q$ to $\infty_1$ induces a functor $\cP_{0,1} \to
\cP'_{0}$; this produces a comparison map 
\begin{equation}\label{eq:proj}
|\delta_{01}| \to |\delta_{01}^0|
\end{equation}
as follows.  For simplices that correspond to chains that do not
contain $q$, we use the identity map and the restriction map.
Otherwise, the comparison map of Equation~\eqref{eq:proj} is induced
by the composite of the map of bar constructions arising from the
composite of the map $\delta_{01} \to \delta_{01}^0$ and composition
and restriction as above.  Note that this map is usually not an
equivalence; it is our model for the projection.

We complete the construction by comparing $|\delta^0_{01}|$ to
$|\delta_0|$.  By inspection, the filtration by action
(Proposition~\ref{prop:computed_associated_graded}) implies 
that $B_{\bullet}(\bS_1, C^*_{\rel \partial}(\cP'_0, \Omega \bS),
\delta_{01}^0)$ and $B_{\bullet}(\bS_0, C^*_{\rel \partial}(\cP_0,
\Omega_{\bS}), \delta_0)$ are abstractly equivalent.

We construct an explicit zig-zag exhibiting this equivalence using the
Kan suspension.  Specifically, there is a map
\begin{equation} \label{eq:map_of_homotopy_types_after_suspend_loop}
|\delta_{01}^0| \to |\tilde{\Sigma} \Omega B_{\bullet}(\bS_0,
C^*_{\rel \partial}(\cP_0, \Omega \bS), \delta_0)|
\end{equation}
induced by a map of simplicial spectra $\delta_{01}^{0} \to
\tilde{\Sigma} \Omega \delta_0$.  Here $\Omega \delta_0$ is the
cotensor with $S^1$ in the category of simplicial spectra; this is
computed levelwise, and amounts simply to looping each wedge summand
in the $k$-simplices.  Since there is a natural equivalence
\begin{equation}
|\tilde{\Sigma} \Omega B_{\bullet}(\bS_0,
C^*_{\rel \partial}(\cP_0, \Omega \bS), \delta_0)| \longleftarrow
\Sigma |\Omega B_{\bullet}(\bS_0,
C^*_{\rel \partial}(\cP_0, \Omega \bS), \delta_0)|
\end{equation}
and natural equivalences
\begin{equation}
     \begin{tikzcd}
     \Sigma {|}\Omega B_{\bullet}(\bS_0,
C^*_{\rel \partial}(\cP_0, \Omega \bS), \delta_0){|} \ar[r] & \Sigma \Omega {|}\delta_0{|} \ar[r] &
{|}\delta_0{|},
    \end{tikzcd}
\end{equation}
given this claim we can conclude the desired comparison result.  (Here
we are using the fact that the loop functor commutes with geometric
realization of simplicial spectra.)

The map of simplicial spectra $|\delta_{01}^{0}| \to
\tilde{\Sigma} \Omega |\delta_0|$ is constructed as follows.  The
$1$-simplices of $|\delta_{01}^{0}|$ are specified by the contractible
spectrum 
\begin{equation}
\C^*_{\rel \partial}(\cP_{0,1}; \Omega \bS)(p, \infty_1) \sma \delta_p
\end{equation}
and the non-degenerate $2$-simplices are 
\begin{equation}
C^*_{\rel \partial}(\cP_{0,1}; \Omega \bS)(p, \infty_0) \sma 
C^*_{\rel \partial}(\cP_{0,1}; \Omega \bS)(\infty_0, \infty_1) \sma \delta_p
\htp \Omega^2 \delta_p.
\end{equation}

The $2$-simplices of $\tilde{\Sigma} \Omega \delta_0$ are 
\begin{equation}
\Omega \left(C^*_{\rel \partial}(\cP_0; \Omega \bS)(p, \infty) \sma \delta_p\right) \htp \Omega^2\delta_p,
\end{equation}
and the $1$-simplices and the $0$-simplices are $\ast$.  The map in
this case is specified on the $2$-simplices as the equivalence
smashing together the identity map on $\delta_p$, the
equivalence 
\begin{equation}
C^*_{\rel \partial}(\cP_{0,1}; \Omega \bS)(p, \infty_0) \to C^*_{\rel
  \partial}(\cP_1; \Omega \bS)(p, \infty)
\end{equation}
induced by the functor $\cP_{0,1} \to \cP_{1} $, and the identification
\begin{equation}
C^*_{\rel \partial}(\cP_{0,1}; \Omega \bS)(\infty_0, \infty_1) \to
\Omega \bS.
\end{equation}
The idea is that the extra $\Omega \bS$ term arising from looping and
the suspension serves as the receptacle of the mapping spectra in the
domain that do not involve $p$.  On the $1$-simplices, that map is
just the collapse map; this is a weak equivalence.  A straightforward
check verifies that these levelwise maps are compatible with the
simplicial identities.

Therefore, we end up with a zig-zag of maps, in which all but one
arrow is an equivalence: 
\begin{equation}
  \begin{tikzcd}
  {|}\delta_{1}{|}  & \arrow{l}[swap]{\htp} | \delta_{01}{|}  \arrow[r]  \ar[r] & {|}\delta_{01}^0{|} & \ar[l,swap,"\htp"]   \cdots \ar[r,"\htp"]  & {|}\delta_0{|}.     
  \end{tikzcd}
\end{equation}
\end{example}

We now explain a generalization of this procedure for constructing
representatives of a homotopy class of maps from
$|\delta_{\lambda_2}|$ to $|\delta_{\lambda_1}|$ that collapse the
elements in $\cP_{\lambda_2}$ of action greater than $\lambda_1$.

\begin{defin}
Given a pair $\lambda_1 < \lambda_2$ of real numbers, the partially
ordered set $\cP_{\lambda_2, \lambda_1}$ is specified as having
elements the union of those of $\cP_{\lambda_2}$, together with two
additional elements $\{\infty_1, \infty_2\} $.

In mild abuse of notation, we will write $\cP \cap \cP_{\lambda_2,
  \lambda_1}$ to denote the subset of $\cP_{\lambda_2,
  \lambda_1}$ coming from $\cP$.  

The morphisms of $\cP_{\lambda_2,\lambda_1}$ are uniquely determined
by the following requirements:

\begin{enumerate}

\item The poset $\cP_{\lambda_2}$ is a full subcategory of
  $\cP_{\lambda_2, \lambda_1}$.

\item There is a unique map $\infty_1 \to \infty_2$; that is the
  only map with domain $\infty_1$.

\item For each $p \in \cP \cap \cP_{\lambda_2, \lambda_1}$, there is a
  unique map $p \to \infty_1$ when $\cA(p) \leq \lambda_1$ and no map
  from $p \to \infty_1$ otherwise.

\item There is a unique map $p \to \infty_2$ for each $p \in \cP \cap
  \cP_{\lambda_2, \lambda_1}$.

\end{enumerate}
\end{defin}
Note that the object $\infty_2$ is a terminal object and the
subcategory $\infty_1 \to \infty_2$ is a terminal spine, in the sense
that all maps in $\cP_{\lambda_2, \lambda_1}$ with source $\infty_1$
or $\infty_2$ have target in the subcategory.  Moreover, the
undercategory of $\infty_1$ in $\cP_{\lambda_2, \lambda_1}$ is
isomorphic to $\cP_{\lambda_1}$.

We now explain how to associate a homotopy type to $\cP_{\lambda_2,
  \lambda_1}$ and to a $\cP$-cellular diagram $\delta$.

\begin{defin}
Define $\bS_{\lambda_2, \lambda_1}$ to be the right
$C^*_{\rel \partial}(\cP_{\lambda_2, \lambda_1}; \Omega \bS)$-module
specified on objects $p \in \cP_{\lambda_2, \lambda_1}$ via
\begin{equation}
\begin{cases}
\bS_{\lambda_2, \lambda_1}(p) = \bS \qquad &p = \infty_{\lambda_2}
\\ 
\bS_{\lambda_2, \lambda_1}(p) = \ast \qquad &\textrm{otherwise} 
\end{cases}
\end{equation}
and with the evident structure maps.  Define $\delta_{\lambda_2,
  \lambda_1}$ to be the left $C^*_{\rel \partial}(\cP_{\lambda_2,
  \lambda_1}; \Omega \bS)$-module specified on objects $p \in
\cP_{\lambda_2, \lambda_1}$ as  
\begin{equation}
\begin{cases}
\delta_{\lambda_2, \lambda_1}(p) = \delta_p \qquad &p \in \cP \cap \cP_{\lambda_2, \lambda_1} \\
\delta_{\lambda_2, \lambda_1}(p) = \ast \qquad &\textrm{otherwise}
\end{cases}
\end{equation}
and with structure maps inherited from $\delta$.
\end{defin}

We make the following definition in analogy with
Definition~\ref{defn:homotopytype}.

\begin{defin}\label{defn:spinehomotopytype}
The homotopy type $|\delta_{\lambda_2, \lambda_1}|$ is the geometric
realization of the simplicial spectrum given as the bar construction
\begin{equation}
B(\bS_{\lambda_2, \lambda_1}, C^*_{\rel \partial}(\cP_{\lambda_2,
  \lambda_1}; \Omega \bS), \delta_{\lambda_2, \lambda_1}),
\end{equation}
which is a model for the (derived) tensor product of functors 
\begin{equation}
\bS_{\lambda_2, \lambda_1} \sma^L_{C^*_{\rel \partial}(\cP_{\lambda_2,
    \lambda_1}; \Omega \bS)} \delta_{\lambda_2, \lambda_1}. 
\end{equation}
\end{defin}

There is an evident functor $\cP_{\lambda_2, \lambda_1} \to
\cP_{\lambda_2}$ defined to be the identity on $\cP \subset
\cP_{\lambda_2, \lambda_1}$ and to take $\infty_1$ and $\infty_2$ to
the terminal object of $\cP_{\lambda_2}$.  We can use this to define a
comparison map on homotopy types as follows.

\begin{lem}\label{lem:zig}
There exists a natural map
\begin{equation}
B(\bS_{\lambda_2, \lambda_1},C^*_{\rel \partial}(\cP_{\lambda_2,
  \lambda_1}; \Omega \bS), \delta_{\lambda_2, \lambda_1}) 
\to
B(\bS_{\lambda_2}, C^*_{\rel
  \partial}(\cP_{\lambda_2}; \Omega \bS), \delta_{\lambda_2})
\end{equation}
specified by restriction and collapsing.
\end{lem}

\begin{proof}
A nontrivial $k$-simplex of $B_{\bullet}(\bS_{\lambda_2},C^*_{\rel
  \partial}(\cP_{\lambda_2}; \Omega \bS), \delta_{\lambda_2})$
corresponds to a totally ordered length $k$ subset $Q$ of
$\cP_{\lambda_2}$ of the form  
\begin{equation}
p = q_0 < q_1 < q_2 < \ldots < q_{k-1} < \infty.
\end{equation}
Consider a subset $\tilde{Q}$ of $\cP_{\Lambda}$ that is in the
inverse image of $Q$ under the projection $\cP_{\Lambda} \to \cP_{\max
  \Lambda}$; this labels a simplex of the domain.  If $Q = \{p, q_1,
\ldots, q_{k-1}, \infty\}$ and $\cA(q_{k-1}) > \lambda_1$, then
$\tilde{Q} = \{p, q_1, \ldots, q_{k-1}, \infty_2\}$ and we map the
associated $k$-simplex by the identity map.  Otherwise, $\tilde{Q}$
either has the form 
\begin{equation}
p = q_0 < q_1 < q_2 < \ldots < q_{k-1} < \infty_{2}.
\end{equation}
or   
\begin{equation}
p = q_0 < q_1 < q_2 < \ldots < q_{k-1} < \infty_{1} < \infty_{2},
\end{equation}
with $\cA(q_{k-1}) \leq \lambda_1$.  Since the cubical complex
$\hat{\cP}_{\lambda_2}(q_{k-1}, \infty)$ is a subcomplex of
$\hat{\cP}_{\lambda_2, \lambda_1}(q_{k-1}, \infty_2)$ and analogously 
\begin{equation}
\partial \hat{\cP}_{\lambda_2}(q_{k-1}, \infty) \subseteq
\partial \hat{\cP}_{\lambda_2, \lambda_1}(q_{k-1}, \infty_2),
\end{equation}
there are restriction maps
\begin{equation}
C^*_{\rel \partial}(\cP_{\lambda_2, \lambda_1}; \Omega \bS)(q_{k-1},
\infty_2) \to 
C^*_{\rel \partial}(\cP_{\lambda_2}; \Omega \bS)(q_{k-1}, \infty)
\end{equation}
which are compatible with the composition and module structure maps.
Therefore, in the first case we have a map 
\begin{equation}
\begin{tikzcd}
C^*_{\rel \partial}(p, q_1) \sma C^*_{\rel \partial}(q_1, q_2) \sma
\ldots \sma C^*_{\rel \partial}(q_{k-1}, \infty_2) \ar[d] \\
C^*_{\rel \partial}(p, q_1) \sma C^*_{\rel \partial}(q_1, q_2) \sma
\ldots \sma C^*_{\rel \partial}(q_{k-1}, \infty) 
\end{tikzcd}
\end{equation}
induced by the restriction and in the second case we have a map
\begin{equation}
\begin{tikzcd}
C^*_{\rel \partial}(p, q_1) \sma C^*_{\rel \partial}(q_1, q_2) \sma
\ldots \sma C^*_{\rel \partial}(q_{k-1}, \infty_1) \sma C_*^{\rel
  \partial}(\infty_1, \infty_2) \ar[d] \\
C^*_{\rel \partial}(p, q_1) \sma C^*_{\rel \partial}(q_1, q_2) \sma
\ldots \sma C^*_{\rel \partial}(q_{k-1}, \infty) \sma C^*_{\rel
  \partial}(\infty, \infty) 
\end{tikzcd}
\end{equation}
determined by the identity maps, the composition, and the
restriction.  This specifies a simplicial map and hence on realization
a map of spectra.
\end{proof}

Under our standing hypotheses listed at the beginning of Section
\ref{sec:homotopy-types-from}, this map is an equivalence: 

\begin{prop}\label{prop:spine}
The comparison map defined in Lemma~\ref{lem:zig} is a weak
equivalence. 
\end{prop}

\begin{proof}
By construction, the comparison map is compatible with the action
filtration and therefore induces maps
\begin{multline}
B(\bS_{\lambda_2, \lambda_1},C^*_{\rel \partial}(\cP_{\lambda_2,
  \lambda_1, \geq n}; \Omega \bS), \delta_{\lambda_2, \lambda_1, \geq
  n})  \\
\to
B(\bS_{\lambda_2}, C^*_{\rel
  \partial}(\cP_{\lambda_2, \geq n}; \Omega \bS), \delta_{\lambda_2,
  \geq n}) 
\end{multline}
and
\begin{multline}
B(\bS_{\lambda,\lambda_1},C^*_{\rel \partial}(\cP_{\lambda_2,
  \lambda_1, > n}; \Omega \bS), \delta_{\lambda_2,\lambda_1, > n})  \\
\to
B(\bS_{\lambda_2}, C^*_{\rel
  \partial}(\cP_{\lambda_2, > n}; \Omega \bS), \delta_{\lambda_2, >
  n})
\end{multline}
for each $n$.  Since $\delta_{\Lambda}(\infty_i) = \ast$ for all
$\infty_i$, these maps induce an equivalence on associated graded
spectra with respect to the action filtration.  By our hypothesis on
$\cP$, the filtration is finite and thus equivalences on associated
graded spectra induce by induction an equivalence on the geometric
realizations. 
\end{proof}

The point of $|\delta_{\lambda_2, \lambda_1}|$ is to provide an
intermediate object that maps to both $\delta_{\lambda_1}$ and
$\delta_{\lambda_2}$.  Lemma~\ref{lem:zig} constructed the comparison
map to $\delta_{\lambda_2}$; we now turn to explain the map to
$\delta_{\lambda_1}$.

We will write $\cP^{\lambda_1}_{\lambda_2, \lambda_1}$ to denote the
full subcategory of $\cP_{\lambda_2, \lambda_1}$ spanned by the objects of
$\cP$ of action smaller than or equal to $\lambda_1$, $\infty_1$, and
$\infty_2$.  There is a natural functor $\cP_{\lambda_2, \lambda_1}
\to \cP^{\lambda_1}_{\lambda_2, \lambda_1}$ that 
collapses $q \in \cP \cap \cP_{\lambda_2, \lambda_1}$ such that
$\cA(q) > \lambda_1$ to $\infty_2$.  We will construct a corresponding
collapse map on homotopy types.

\begin{lem}\label{lem:collapse}
There is a map of spectra
\begin{equation}
B(\bS_{\lambda_2,\lambda_1}, C^*_{\rel \partial}(\cP_{\lambda_2,
  \lambda_1}; \Omega \bS), \delta_{\lambda_2, \lambda_1}) \to
B(\bS^{\lambda_1}_{\lambda_2, \lambda_1}, C^*_{\rel
  \partial}(\cP^{\lambda_1}_{\lambda_2, \lambda_1}; \Omega \bS),
\delta^{\lambda_1}_{\lambda_2, \lambda_1}) 
\end{equation}
induced by restriction and collapse.
\end{lem}

\begin{proof}
Simplices in the domain can be labeled by totally ordered subsets $Q =
\{q_1, q_2, \ldots, \infty_2\}$ of $\cP_{\Lambda}$.  The map is in question
is defined levelwise on wedge summands as follows:
\begin{enumerate}
\item Simplices labeled by subsets $Q$ such that $\cA(q_1) >
  \lambda_1$ are taken to the basepoint in the range.  
\item All other simplices are mapped using the identity on mapping
  spectra for $q_i, q_{i+1} \in \cP \cap \cP_{\lambda_2, \lambda_1}$,
  the restriction on pairs $q_i, \infty_{2}$ with $q_i \in
  \cP_{\lambda_1}$ (as in the construction of the map in
  Lemma~\ref{lem:zig}), and the identity on mapping spectra for the
  pair $\infty_{1}, \infty_{2}$.
\end{enumerate}
It is straightforward to check that these assignments assemble into a
simplicial map.
\end{proof}

Finally, we can compare the range to $B(\bS_{\lambda_1}, C^*_{\rel
  \partial}(\cP_{\lambda_1}; \Omega \bS), \delta_{\lambda_1})$, as
follows.

\begin{prop}\label{prop:shiftcompare}
There is a zig-zag of equivalences 
\begin{equation}
  \begin{tikzcd}[row sep=small]
    B(\bS^{\lambda_1}_{\lambda_2, \lambda_1}, C^*_{\rel
  \partial}(\cP^{\lambda_1}_{\lambda_2, \lambda_1}; \Omega \bS), 
\delta^{\lambda_1}_{\lambda_2, \lambda_1})  \\
\ar[u, "\htp"] \ar[d, swap, "\htp"] \cdots \\ 
B(\bS_{\lambda_1},
C^*_{\rel \partial}(\cP_{\lambda_1}; \Omega \bS),
\delta_{\lambda_1}).
  \end{tikzcd}
\end{equation}
\end{prop}

\begin{proof}
The point is that the Kan suspension of $\Omega \delta_{\lambda_1}$
has no $1$-simplices and the $k$-simplices corresponds
bijectively to the $k-1$-simplices of $\cP^{\lambda_1}_{\lambda_2,
  \lambda_1}$.  Specifically, we define a map of simplicial spectra
\begin{equation}
B(S^{\lambda_1}_{\lambda_2, \lambda_1}, \cP_{\lambda_2, \lambda_1},
\delta_{\lambda_2, \lambda_1}) \to \tilde{\Sigma} \Omega
B(S_{\lambda_1}, \cP_{\lambda_1}, \delta_{\lambda_1})
\end{equation}
as the identity on terms coming from $\delta_p$, the identity on
mapping spectra corresponding to objects in $\cP \cap \cP_{\Lambda}$,
and the natural comparison map to $\Omega \bS$ for the morphism
spectrum of $\infty_1, \infty_2$.  It is straightforward to see that
this is a weak equivalence by inspecting the levelwise cofiber.
\end{proof}

We now will use the preceding work to define a zig-zag representing
the projection map $\delta_{\lambda_2} \to \delta_{\lambda_1}$.  

\begin{defin}
We denote by $\Delta_{\lambda_1}^{\lambda_2}$ the zig-zag representing
a homotopy class  
\begin{equation}
\delta_{\lambda_2} \to \delta_{\lambda_1}
\end{equation}
constructed as the composite zig-zag 
\begin{equation}
\begin{tikzcd}
{|}\delta_{\lambda_2}{|} & \ar[l,"\htp"]
|\delta_{\lambda_2, \lambda_1}{|} \ar[r] & {|}
\delta^{\lambda_1}_{\lambda_2, \lambda_1} {|} & 
\ar[r,swap,"\htp"] \ar[l,"\htp"] \ldots & {|} \delta_{\lambda_1} {|},
\end{tikzcd}
\end{equation}
where the leftmost equivalence is defined in Lemma~\ref{lem:zig}, the
righthand middle map in Lemma~\ref{lem:collapse}, and the
comparison zig-zag of equivalences in
Proposition~\ref{prop:shiftcompare}.
\end{defin}

The following sanity check verifies that these projection maps behave
the way we expect.  Here recall that from
Definition~\ref{defn:bounded-homotopy-type} that $|\delta_{\lambda_2,
  > \lambda_1}|$ is the homotopy type obtained by restricting to the
subcategory of $\cP_{\lambda_2}$ of objects with action $> \lambda_1$.

\begin{prop}\label{prop:sanity}
For $\lambda_2 > \lambda_1$, there is an equivalence
\begin{equation}
|\delta_{\lambda_2, > \lambda_1}| \htp \hofib(|\delta_2| \to |\delta_1|).
\end{equation}
\end{prop}

\begin{proof}
It suffices to consider the homotopy fiber of the map
\begin{equation}\label{eq:projection}
B(\bS_{\lambda_2, \lambda_1}, C^*_{\rel \partial}(\cP_{\lambda_2,
  \lambda_1}; \Omega \bS), \delta_{\lambda_2, \lambda_1}) \to
B(\bS^{\lambda_1}_{\lambda_2, \lambda_1}, C^*_{\rel
  \partial}(\cP^{\lambda_1}_{\lambda_2, \lambda_1}; \Omega \bS),
\delta^{\lambda_1}_{\lambda_2, \lambda_1})
\end{equation}
from Lemma~\ref{lem:collapse}.  We do this by factoring the map, as
follows.  First, we have the map 
\begin{equation}
B(\bS_{\lambda_2, \lambda_1}, C^*_{\rel \partial}(\cP_{\lambda_2,
  \lambda_1}; \Omega \bS), \delta_{\lambda_2, \lambda_1}) \to
B(\bS_{\lambda_2, \lambda_1}, C^*_{\rel \partial}(\cP_{\lambda_2,
  \lambda_1}; \Omega \bS), \delta^{\lambda_1}_{\lambda_2, \lambda_1}), 
\end{equation}
Then we compose with the map
\begin{equation}
B(\bS_{\lambda_2, \lambda_1}, C^*_{\rel \partial}(\cP_{\lambda_2,
  \lambda_1} \Omega \bS), \delta^{\lambda_1}_{\lambda_2, \lambda_1})
\to B(\bS^{\lambda_1}_{\lambda_2, \lambda_1}, C^*_{\rel
  \partial}(\cP^{\lambda_1}_{\lambda_2, \lambda_1}; \Omega \bS),
\delta^{\lambda_1}_{\lambda_2, \lambda_1})
\end{equation}
induced by the restrictions.  The composite is precisely the map from
Equation~\eqref{eq:projection}.  The homotopy fiber of the first of
these is clearly equivalent to $|\delta_{\lambda_2, > \lambda_1}|$;
the simplices in the fiber are precisely those corresponding to chains
that are contained in $\cP_{> \lambda_1}$.  The second map is a weak
equivalence by the filtration argument, since it is action-preserving. 
\end{proof}

Next, we have the following associativity property for these
zig-zags.

\begin{lem}\label{lem:zigassoc}
For $\lambda_0 < \lambda_1 < \lambda_2$, the zig-zags
$\Delta_{\lambda_1}^{\lambda_2} \circ \Delta_{\lambda_0}^{\lambda_1}$
and $\Delta_{\lambda_0}^{\lambda_2}$ represent the same map in the
homotopy category.
\end{lem}

\begin{proof}
For this purpose it is convenient to introduce another construction
akin to $\cP_{\lambda', \lambda}$ with an additional intermediate
object.  Specifically, we consider the category $\cP_{\lambda_3,
  \lambda_2, \lambda_1}$, which has two auxiliary objects $\infty_1$
and $\infty_2$ in addition to the terminal object $\infty$, and
morphisms specified such that:
\begin{enumerate}
\item $\cP_{\lambda_1}$, $\cP_{\lambda_2}$, and $\cP_{\lambda_3}$ are
  full subcategories, 
\item  the unique non-terminal map with domain $\infty_1$ has codomain $\infty_2$,
  and
\item for $p$ such that $\cA(p) \leq \lambda_1$, $p$ has a unique map
  to $\infty_1$ and $\infty_2$ and for $p$ such that $\cA(p) \leq
  \lambda_2$, $p$ has a unique map to $\infty_2$.
\end{enumerate}
For convenience, we write $|\delta_{\lambda_3, \lambda_2, \lambda_1}|$
for the associated homotopy type.

Now, there are collapse maps
\begin{equation}
\begin{tikzcd}
\relax|\delta_{\lambda_3, \lambda_2, \lambda_1}| \ar[r] &
\relax|\delta^{\lambda_2}_{\lambda_3, \lambda_2, \lambda_1}| \ar[r] & 
\relax|\delta_{\lambda_2, \lambda_1}|
\end{tikzcd}
\end{equation}
and 
\begin{equation}
\begin{tikzcd}
\relax|\delta_{\lambda_3, \lambda_2, \lambda_1}| \ar[r] & 
\relax|\delta^{\lambda_1}_{\lambda_3, \lambda_2, \lambda_1}| & \ldots
\ar[l] \ar[r] & \relax|\delta_{\lambda_1}|.
\end{tikzcd}  
\end{equation}
In addition, there are collapse maps
\begin{equation}
\begin{tikzcd}
\relax |\delta_{\lambda_3, \lambda_2, \lambda_1}| \ar[r] & \relax
|\delta_{\lambda_3, 
  \lambda_1}| &&
\relax |\delta^{\lambda_1}_{\lambda_3, \lambda_2, \lambda_1}| \ar[r] &
\relax |\delta^{\lambda_1}_{\lambda_3, \lambda_1}|
\end{tikzcd}
\end{equation} 
which are weak equivalences, by the argument for
Proposition~\ref{prop:spine}.

These comparison weak equivalences are compatible with the collapse
zigzags to $|\delta_{\lambda_1}|$.  Specifically, the diagram
\begin{equation}
\begin{tikzcd}
\relax|\delta_{\lambda_3, \lambda_2, \lambda_1}| \ar[d] \ar[r] &
\relax|\delta_{\lambda_3, \lambda_1}| \ar[d] \\
\relax|\delta^{\lambda_1}_{\lambda_3, \lambda_2, \lambda_1}| \ar[r] &
\relax|\delta^{\lambda_1}_{\lambda_3, \lambda_1}| 
\end{tikzcd}
\end{equation}
commutes.  This implies that the collapse zigzag $|\delta_{\lambda_3,
  \lambda_2, \lambda_1}| \leftrightarrow |\delta_{\lambda_1}|$ is
homotopic to the zig-zag $|\delta_{\lambda_3, \lambda_1}|
\leftrightarrow |\delta_{\lambda_1}|$.

Next, the collapse map $|\delta_{\lambda_3, \lambda_2, \lambda_1}| \to
|\delta^{\lambda_1}_{\lambda_3, \lambda_2, \lambda_1}|$ factors
through $|\delta^{\lambda_2}_{\lambda_3, \lambda_2, \lambda_1}|$, 
and the diagram 
\begin{equation}
\begin{tikzcd}
\relax |\delta_{\lambda_3, \lambda_2}| \ar[d] & 
\relax |\delta_{\lambda_3, \lambda_2, \lambda_1}| \ar[l] \ar[d] & \\
\relax |\delta^{\lambda_2}_{\lambda_3, \lambda_2}| & \ar[l]
\relax |\delta^{\lambda_2}_{\lambda_3, \lambda_2, \lambda_1}| \ar[d]
\ar[r] & \relax |\delta_{\lambda_2, \lambda_1}| \ar[d] \\
& \relax |\delta^{\lambda_1}_{\lambda_3, \lambda_2, \lambda_1}| \ar[r] & 
\relax |\delta^{\lambda_1}_{\lambda_2, \lambda_1}|
\end{tikzcd}
\end{equation}
commutes.  The result now follows from the fact that all of the
horizontal maps are weak equivalences, again by the argument for
Proposition~\ref{prop:spine}.
\end{proof}

Ideally, we would now define the homotopy type of a $\Pi$-equivariant
$\cP$-cellular diagram $\delta \co C^{*}_{\rel \partial}(\cP; \Omega
\bS) \to \kmod$ as an analogous homotopy limit of a diagram
\begin{equation}
\holim_\lambda |\delta_\lambda|
\end{equation}
over the zig-zags $\Delta_{\lambda_1}^{\lambda_2}$ for each $\lambda_2
> \lambda_1$ in $\bR$, with a fixed choice of $\Lambda_{\lambda_2,
  \lambda_1}$ for each pair $(\lambda_1, \lambda_2)$.  However, this
definition is unmanageable as written; indexing over $\bR$ raises
various technical problems.

Instead, we may consider a countable indexing set $I \subset \bR$
which is final (i.e., admitting a subsequence going to $+\infty$), and
which we give the induced (total) order from $\bR$.  We have a diagram
constructed from the spectra $\{|\delta_i|\}$ and the zig-zags
connecting $|\delta_i|$ and $|\delta_{i-1}|$; we denote the indexing
category for this diagram by $\tilde{I}$.  We can compute
$\holim_{\tilde{I}} |\delta_{-}|$ as a model for the filtered homotopy
type.  Lemma~\ref{lem:zigassoc} implies that any choice of diagram
$\tilde{I}$ will result in an equivalent homotopy limit.  (See
Appendix~\ref{sec:rect} for a discussion of homotopy limits over
zig-zags.)

Since we shall later study the action of the Novikov ring on this
homotopy type, it is convenient to assume that $I$ is closed under
addition. We thus specialise to consider the set $\{k \epsilon\}_{k
  \in \bZ} $.  The next lemma records the fact that the homotopy type
is independent of the choice of $\epsilon$. 

\begin{lem}\label{lem:repindependent}
For any constants $\epsilon_1, \epsilon_2 \in \bR$, there is a natural 
zig-zag of equivalences
\begin{equation}
\holim_{k \in \bZ} |\delta_{k \epsilon_1}| \htp
\holim_{k \in \bZ} |\delta_{k\epsilon_2}|.
\end{equation} 
\qed
\end{lem}

This invariance result justifies the following definition.

\begin{defin}\label{defn:invht}
The homotopy type $|\delta|$ of a $\Pi$-equivariant $\cP$-cellular diagram
$\delta$ is the object of $\Ho(\Sp)$ represented by 
\begin{equation}
\holim_{k \in \bZ} |\delta_{k\epsilon}|.
\end{equation}
for any choice of positive non-zero constant $\epsilon$.
\end{defin}

When studying the module structure, we will need to work with various
specific representatives of $|\delta|$.

\begin{defin}
Let $I \subset \bR$ be a countable discrete subset.  We write
$\delta(I)$ to denote the diagram of zig-zags indexed by $I$. 
\end{defin}

We will typically work with systems $\delta(\ldots \to k\epsilon
\to \ldots)$ and translations thereof.  Given a countable discrete
diagram $I \subset \bR$, we can consider the restricted diagrams
$I_{\leq \lambda} = \{i \in I \, \mid \, i \leq \lambda\}$; the
inclusions $I_{\leq \lambda_1} \to I_{\leq \lambda_2}$ for $\lambda_2
> \lambda_1$ induce compatible maps of homotopy limits
\begin{equation} \label{eq:filtration_full_homotopy_type}
\holim_{I_{\leq \lambda_2}} |\delta_{i}| \to \holim_{I_{\leq
    \lambda_1}} |\delta_i|.
\end{equation}
Therefore, we regard the homotopy type of a $\Pi$-equivariant
$\cP$-cellular diagram as equipped with a natural filtration by
action. We record the following variant of
Proposition~\ref{prop:computed_associated_graded}:

\begin{prop}\label{prop:filt-homot-type}
Let $I$ be $\epsilon$-dense and suppose there are no morphisms from
$p$ to $q$ whenever $\cA(q) \leq \cA(p) + \epsilon$.  Then the
associated graded spectra for the filtration on $|\delta|$ given by
Equation~\eqref{eq:filtration_full_homotopy_type} can be described for
$\lambda_2 = n\epsilon$ and $\lambda_1 = (n-1)\epsilon$ as 
\begin{equation}
\hofib(\delta_{\lambda_2} \to \delta_{\lambda_1}) \htp
\bigvee_{\cA(p) \in (\lambda_1, \lambda_2]} \delta(p).
\end{equation}
\qed
\end{prop}

Although we did not construct it directly this way, we view $|\delta|$
as a completion with respect to a decreasing filtration.  As a
consequence, $|\delta|$ is itself equipped with the decreasing
filtration given by the homotopy fibers of the natural maps $|\delta|
\to |\delta_{\lambda}|$, induced as above.  Here by a filtration we
again mean a zig-zag with upwards arrows weak equivalences, induced by
passage to homotopy fibers from the diagrams
\begin{equation}\label{eqn:delta_dec_filt}
\begin{tikzcd}
{|}\delta{|} \ar[r] \ar[dr, dashrightarrow] \ar[ddr] & {|}\delta_{\lambda_2}{|} \\
& \vdots \ar[u] \ar[d] \\
& {|}\delta_{\lambda_1}{|}
\end{tikzcd}
\end{equation}
for $\lambda_2 > \lambda_1$.  Therefore, choosing representatives for
the zig-zags, we get a filtration
\begin{equation}\label{eqn:unused-filtration}
\ldots \to \hofib(|\delta \to |\delta_{\lambda_2}|) \to \hofib(|\delta
\to |\delta_{\lambda_1}|) \to \hofib(|\delta \to |\delta_{\lambda_0}|)
\to \ldots 
\end{equation}

\begin{prop}
This decreasing filtration is complete in the sense that
\begin{equation}
|\delta| \htp \hocolim_{\lambda \to -\infty} \left(\hofib(\delta \to
\delta_{\lambda})\right) 
\end{equation}
and
\begin{equation}
\holim_{\lambda \to \infty} \left( \hofib(\delta \to
\delta_{\lambda})\right) \htp \ast. 
\end{equation} \qed
\end{prop}

It will be useful later on to have the following alternate
characterization of the terms in this filtration.

\begin{lem}\label{lem:deltafilt}
Suppose that $|\delta|$ is represented by the homotopy limit over $I
\subset \bR$.  Then for any $\lambda_0$, there is an equivalence
\begin{equation}
\hofib(\delta \to \delta_{\lambda_0}) \htp \holim_{i \in I}
|\delta_{i, > \lambda_0}|.
\end{equation}
\end{lem}

\begin{proof}
Since homotopy limits commute up to weak equivalence, we can compute
\begin{equation}
\hofib(\delta \to \delta_{\lambda_0}) \htp \holim_i \left(\hofib(\delta_i
\to \delta_{\lambda_0})\right).
\end{equation}
The result now follows from Lemma~\ref{prop:sanity}.
\end{proof}

Here notice that the maps in the homotopy limit system
$\holim_{\lambda} |\delta_{\lambda, > \lambda_0}|$ are precisely the
zig-zags $\Delta_{\lambda_1}^{\lambda_2}$ restricted to represent a
homotopy class of maps $|\delta_{\lambda_2, > \lambda_0}| \to
|\delta_{\lambda_1, > \lambda_0}|$.  We refer to these terms as the
{\em truncated homotopy limits}.

The maps in the filtration are induced by the natural inclusions.

\begin{lem}
For $\lambda_2 > \lambda_1 > \lambda_0$, there are natural inclusion maps
\begin{equation}
|\delta_{\lambda_2, > \lambda_1}| \to |\delta_{\lambda_2, > \lambda_0}|
\end{equation}
which are strictly associative.
\end{lem}

\begin{proof}
The maps are induced by the evident levelwise simplicial inclusions;
the category $C^*_{\rel \partial}(\cP_{\lambda_2, > \lambda_1}; \Omega
\bS)$ is a full subcategory of $C^*_{\rel \partial}(\cP_{\lambda_2, >
  \lambda_0}; \Omega \bS)$.  Associativity is evident.
\end{proof}

Assembling these, we have the following induced maps on homotopy
limits arising as the composite of the inclusions levelwise. 

\begin{lem}
Fix $I \subset \bR$.  Then for each $\lambda_0, \lambda_1 \in I$ such
that $\lambda_1 > \lambda_0$, there are natural maps 
\begin{equation}
\holim_{i \in I} |\delta_{i,> \lambda_1}| \to 
\holim_{i \in I} |\delta_{i, > \lambda_0}|
\end{equation}
which are strictly associative. \qed
\end{lem}

In light of this, we take the system
\begin{equation}\label{eqn:filtration}
\ldots \to \holim_{i \in I} |\delta_{i, > \lambda_2}| \to 
\holim_{i \in I} |\delta_{i, > \lambda_1}| \to 
\holim_{i \in I} |\delta_{i, > \lambda_0}| \to \ldots
\end{equation}
as our model of the decreasing filtration on $|\delta|$.

\begin{lem}\label{lem:delta-complete}
The decreasing filtration of Equation~\eqref{eqn:filtration} is
complete; the homotopy limit is trivial and the homotopy colimit is
equivalent to $\delta$.
\end{lem}

\begin{rem}
Note that we have not in fact constructed a map of filtered spectra
comparing this decreasing filtration with the one induced by
Equation~\eqref{eqn:unused-filtration}; the natural comparison maps 
are compatible with the structure maps only up to homotopy.  With a
little more work it is possible to rigidify these comparison maps, but
since for our purposes it suffices to simply take the filtration
described in Equation~\eqref{eqn:filtration} as the definition, we do
not carry this out here.
\end{rem}

\begin{example} \label{ex:colimit_infinite_chain}
Say that $\cP = \bZ$, and $\delta$ is the functor mapping all objects
to $\bbk$, and all morphisms to $0$. According to
Proposition~\ref{prop:nullhomtype} and
Equation~\eqref{eq:homotopy_type_big_wedge}, the homotopy type of 
$|\delta_\lambda|$ for each real number $\lambda$ is a wedge of copies
of $\bbk$ indexed by integers smaller than $\lambda$.  Therefore, by
definition the homotopy type of $\delta$ is the homotopy limit
\begin{equation}
\holim_\lambda \bigvee_{k \leq \lambda, k \in \bZ} \bbk,
\end{equation}
where the structure maps in the inverse system are given by projecting
away from the missing summands.  Since the system of homotopy groups
is evidently Mittag-Leffler, there is an additive identification
\begin{equation}
\pi_* |\delta| \cong \bbk_*((t)) \equiv \bbk_*[[t]][t^{-1}],
\end{equation} 
where here $R((t))$ denotes Laurent series (i.e. series in $t$ and
$t^{-1}$ with finitely many monomials of negative exponent) and
$\bbk_*$ denotes the homotopy groups of $\bbk$.  In the formula
above, the formal variable $t$ is in degree $0$.  In order to put the
variable in a different degree, we can consider suspensions of $\bbk$.
This identification can be promoted to a multiplicative equivalence,
although this is a bit complicated.  We explain both aspects of the
general situation in more detail below.

On the other hand, one could try to set $\lambda = -\infty$ in the
definition of the homotopy type, and avoid taking the inverse limit
over $\lambda$.  However, in this example, the resulting two-sided bar
construction is acyclic; there does not exist an element $z \in \cP$
such that $(z,\infty)$ is a pair of successive elements in $\cP_{-\infty}$.
\end{example}

\subsection{The equivariant structure of the filtered homotopy type}

We now introduce the action of $\Pi$ on the homotopy limit of the
spectra $\{|\delta_\lambda|\}$.  We begin by looking at the
interaction of the action with the filtrations by $\lambda$ on $\Pi$
and $\cP$.

First, note that Equation~\eqref{eq:action-homomorphism} implies that
the collection of posets $\{\cP_\lambda\}$ is compatible with the
action of $\Pi$ on $\cP$, in the following sense.  

\begin{lemma}
For each $\pi \in \Pi$ and $\lambda \in \bR$, we have induced isomorphisms
\begin{align}
\alpha_{\pi} \colon \cP_{\cA(p) = \lambda} &\to \cP_{\cA(p) = \cA(\pi) + \lambda} \\
\alpha_{\pi} \colon \cP_{\lambda, > k} &\to \cP_{\lambda + \cA(\pi), >
  k + \cA(\pi)} \\
\alpha_{\pi} \colon \cP_{\lambda} &\to \cP_{\lambda + \cA(\pi)},
\end{align}
where we stipulate that $\pi$ acts trivially on the terminal object.
These maps are evidently strictly associative in the sense that
$\alpha_{\pi_1} \circ \alpha_{\pi_2} = \alpha_{\pi_1 \pi_2}$ and
unital in that $\alpha_{e}$ is the identity map.  \qed
\end{lemma}

Moreover, the proof of Lemma~\ref{lem:strictactrelcoch} establishing a
strict action of $\Pi$ on $C^*_{\rel \partial}(\cP; \Omega \bS)$
lifting the action on $\cP$ induces an analogous action on the
collection $\{C^*_{\rel \partial}(\cP_\lambda; \Omega \bS)\}$ which
is compatible with the categorical structure. 

\begin{prop}
For each $\pi \in \Pi$ and $\lambda \in \bR$, we have induced homeomorphisms
\begin{equation}
\alpha_{\pi} \colon C^*_{\rel \partial}(\cP_{\lambda}; \Omega \bS)(p,q) \to C^*_{\rel
  \partial}(\cP_{\lambda + \cA(\pi)}; \Omega \bS) (\pi p, \pi q) 
\end{equation}
for each $\pi \in \Pi$.  These maps are strictly associative in the
sense that $\alpha_{\pi_1} \circ \alpha_{\pi_2} = \alpha_{\pi_1
  \pi_2}$ and unital in that $\alpha_{e}$ is the identity map.

The maps $\alpha_{\pi}$ are compatible with the composition and
categorical unit in the sense that we have commutative diagrams 
\begin{equation}
  \begin{tikzcd}
    \begin{gathered}
             C^*_{\rel \partial}(\cP_{\lambda})(q,r)  \\ \sma C^*_{\rel
  \partial}(\cP_{\lambda})(p,q)
    \end{gathered}
    \ar[r] \ar[d] &
    \begin{gathered}
      C^*_{\rel 
  \partial}(\cP_{\lambda + \cA(\pi)}) (\pi q, \pi r)  \\
\sma C^*_{\rel \partial}(\cP_{\lambda + \cA(\pi)}) (\pi p, \pi
q)      \end{gathered}
\ar[d] \\
C^*_{\rel \partial}(\cP_{\lambda})(p,r)  \ar[r]
& C^*_{\rel 
  \partial}(\cP_{\lambda + \cA(\pi)}) (\pi p, \pi r) 
  \end{tikzcd}
\end{equation}
(where here the coefficients are understood to be $\Omega \bS$)
and
\begin{equation}
  \begin{tikzcd}
   \bS \ar[r] \ar[dr] & C^*_{\rel \partial}(\cP_{\lambda}; \Omega
\bS)(p,p) \ar[d] \\
& C^*_{\rel 
  \partial}(\cP_{\lambda + cA(\pi)}; \Omega \bS) (\pi p, \pi p). 
  \end{tikzcd}
\end{equation} \qed
\end{prop}

We now turn to describe the interaction of $\Pi$ with the left and
right $C_*^{\rel \partial}(\cP_{\lambda}; \Omega \bS)$-modules that
are used in the bar construction that defines $|\delta_{\lambda}|$.

\begin{lemma}
For $\pi \in \Pi$, there are maps
\begin{equation}
\beta_{\pi} \colon \bS_{\lambda}(p) \to \bS_{\lambda + \cA(\pi)}(\pi p)
\end{equation}
specified by the identity maps.  These maps are strictly associative
and unital.  They are also clearly compatible with the module action
map.\qed
\end{lemma}

The action of $\Pi$ on $\delta$ is lax, which complicates the
description of the structure that arises on the terms of the
filtration.  In the following proposition, for $\pi \in \Pi$ we will
write $\bS[\pi]$ to denote $\bS[- \deg \pi]$.

\begin{prop}
For $\pi \in \Pi$ and $\lambda \in \Lambda$, there are maps of spectra 
\begin{equation}
\gamma_\pi \colon \delta_{\lambda}(p) \sma \bS[\pi] \to
\delta_{\lambda + \cA(\pi)}(\pi p)
\end{equation}
The diagrams
\begin{equation}
\begin{tikzcd}
\delta_\lambda(p) \sma \bS[\pi_1] \sma \bS[\pi_2] \ar[d,swap,
  "\gamma_{\pi_1}"] \ar[r] & 
\delta_{\lambda}(p) \sma \bS[\pi_1 \pi_2] \ar[d,"\gamma_{\pi_1 \pi_2}"]
\\
\delta_{\lambda + \cA(\pi_1)}(\pi_1 p) \sma \bS[\pi_2] \ar[r] &
\delta_{\lambda + \cA(\pi_1) + \cA( \pi_2)} ((\pi_1 \pi_2) p)
\end{tikzcd}
\end{equation}
commute, and analogous diagrams express the associativity of these
maps.

The action maps are also compatible with the module structure on
$\delta$, in that the following diagrams commute:
\begin{equation}
  \begin{tikzcd}
    C^*_{\rel \partial}(\cP_{\lambda}; \Omega \bS)(p,q) \sma \delta_{\lambda}(p)
\ar[r] \ar[d,swap,"\alpha_\pi \sma \gamma_\pi"] \sma \bS[\pi] &
\delta_{\lambda}(q) \ar[d,"\gamma_\pi"] \sma \bS[\pi] \\
C^*_{\rel \partial}(\cP_{\cA(\pi) + \lambda}; \Omega \bS)(\pi p,\pi q)
\sma \delta_{\lambda + \cA(\pi)}(\pi p) 
\ar[r] & \delta_{\lambda + \cA(\pi)}(\pi q) 
  \end{tikzcd}
\end{equation}
and the analogous associativity diagrams also commute.
 \qed
\end{prop}

\begin{prop}\label{prop:action}
For each $\pi \in \Pi$ and $\lambda \in \bR$, there is a natural map
of spectra 
\begin{multline}
\alpha_{\pi} \colon B(S_{\lambda}, C^*_{\rel \partial}(\cP_{\lambda}; \Omega \bS),
\delta_{\lambda}) \sma \bS[\pi] \\ \to B(S_{\lambda+\cA(\pi)}, C^*_{\rel
  \partial}(\cP_{\lambda+\cA(\pi)}; \Omega \bS), \delta_{\lambda+\cA(\pi)}).
\end{multline}
The diagrams
\begin{equation}
\begin{tikzcd}
B(S_{\lambda}, \delta_{\lambda}) \sma \bS[\pi_1] \sma \bS[\pi_2] \ar[r,"\alpha_{\pi_1}"]
\ar[d,swap,"\alpha_{\pi_1 \pi_2}"] &
   B(S_{\lambda}, \delta_{\lambda}) \sma \bS[\pi_1 \pi_2]
   \ar[d] \\
B(S_{\lambda+\cA(\pi_1)}, \delta_{\lambda+\cA(\pi_1)}) \sma \bS[\pi_2] \ar[r]
&
B(S_{\lambda+\cA(\pi_1 + \pi_2)}, \delta_{\lambda+\cA(\pi_1 + \pi_2)})
\end{tikzcd}
\end{equation}
commute, as do the analogous diagrams expressing the associativity of these
maps.  This structure is weakly unital in the sense that, for $0 \in \Pi$, the
map $\alpha_0$ is related to the identity via the unit weak
equivalence $\bS \to \bS[e]$.
\end{prop}

\begin{proof}
For each $k$, the smash product $\beta_\pi \sma \underbrace{\alpha_\pi \sma \ldots \sma
  \alpha_{\pi}}_{(k-1)} \sma \gamma_{\pi}$ specifies a map on $k$-simplices
\begin{multline}
B_k(S_{\lambda},
C^*_{\rel \partial}(\cP_{\lambda}; \Omega \bS), \delta_{\lambda}) \sma
\bS[\pi] \\ \to
B_k(S_{\lambda+\cA(\pi)}, C^*_{\rel
  \partial}(\cP_{\lambda+\cA(\pi)}; \Omega \bS), \delta_{\lambda+\cA(\pi)}).
\end{multline}
The compatibility of the constituent maps with the unit, composition,
and module structure imply that these maps assemble into a simplicial
map 
\begin{multline}
B_\bullet(S_{\lambda},
C^*_{\rel \partial}(\cP_{\lambda}; \Omega \bS), \delta_{\lambda}) \sma
\bS[\pi] \\ \to
B_\bullet(S_{\lambda+\cA(\pi)}, C^*_{\rel
  \partial}(\cP_{\lambda+\cA(\pi)}; \Omega \bS), \delta_{\lambda+\cA(\pi)}).
\end{multline}
The associativity and unitality follow immediately from the discussion above.
\end{proof}

We now turn to describe the compatibility of the $\Pi$ action with the
zigzags appearing the homotopy limit system of
Definition~\ref{defn:invht}.  It is clear that $\pi \in \Pi$ takes a
representative of the projection for $\lambda_1 < \lambda_2$ to a
representative of the projection for $\cA(\pi) + \lambda_1 < \cA(\pi)
+ \lambda_2$.

\begin{prop}\label{prop:equizig}\label{prop:act_zigzag}
For $\pi \in \Pi$, the following diagram of zigzags commutes:
\begin{equation}
\begin{tikzcd}[row sep=small, column sep=tiny, transform shape, nodes={scale=.95}]
{|}\delta_{\lambda_2}{|} \sma \bS[\pi] \ar[d] & \ar[l,"\htp"] \ar[d]
|\delta_{\lambda_2, \lambda_1}{|} \ar[r] \sma \bS[\pi] & {|}
\delta^{\lambda_1}_{\lambda_2, \lambda_1} {|} \sma \bS[\pi]
\ar[d] &
\ar[r,swap,"\htp"] \ar[l,"\htp"] \ar[d] \ldots & {|} \delta_{\lambda_1}
{|} \sma \bS[\pi] \ar[d] 
\\
{|}\delta_{\cA(\pi)+\lambda_2}{|} & \ar[l,"\htp"]
|\delta_{\cA(\pi)+\lambda_2, \cA(\pi)+\lambda_1}{|} \ar[r] & {|}
\delta^{\cA(\pi) + \lambda_1}_{\cA(\pi)+\lambda_2, \cA(\pi) + \lambda_1} {|} & 
\ar[r,swap,"\htp"] \ar[l,"\htp"] \ldots & {|} \delta_{\cA(\pi)+\lambda_1} {|} 
\end{tikzcd}
\end{equation} \qed
\end{prop}

As a consequence, we obtain an action on systems of the following
form.

\begin{prop}\label{prop:action-systems}
For $\pi \in \Pi$ and $I \subset \bR$ a countable discrete subset,
there are natural maps 
\begin{equation}
\delta(I) \sma \bS[\pi] \to \delta(\cA(\pi) + I)
\end{equation}
and therefore there are maps of spectra
\begin{equation}
\alpha_{\pi} \colon \left(\holim_{k \in \bZ} |\delta_{k\epsilon}|\right)
\sma \bS[\pi] \to
\holim_{k \in \bZ} \left(|\delta_{k\epsilon}| \sma \bS[\pi]\right) 
\to
\holim_{k \in \bZ} |\delta_{\cA(\pi) + k\epsilon}|.
\end{equation}  \qed
\end{prop}

We also have the analogous result for the truncated homotopy limits.
Specifically, using the evident analogues of
Proposition~\ref{prop:action} and Proposition~\ref{prop:act_zigzag},
we obtain the following action on the terms of the filtration:

\begin{cor}\label{cor:action-on-filtration} 
For an element $\pi \in \Pi$ and $\lambda_0 \in \bR$, there are maps of spectra
\begin{multline}
\alpha_{\pi} \colon \left(\holim_{k \in \bZ} |\delta_{k \epsilon, >
  \lambda_0}|\right) \sma \bS[\pi] \to \\
\holim_{k \in \bZ} \left(|\delta_{k \epsilon, >
  \lambda_0}| \sma \bS[\pi]\right) \to
\holim_{k \in \bZ} |\delta_{k \epsilon + \cA(\pi), > \lambda_0 +
  \cA(\pi)}|.
\end{multline} \qed
\end{cor}

The following proposition records the way these maps assemble into
coherent actions maps for $\Pi$.  (We write this out for the absolute
case; analogous results hold for the truncated homotopy limits.)

\begin{prop}\label{prop:action-assoc}
For $\pi_1, \pi_2 \in \Pi$, the composite
\begin{equation}
\begin{tikzcd}[row sep=small, column sep=small, transform shape, nodes={scale=.9}]
{|}\delta_{\lambda_2}{|} \sma \bS[\pi_1] \sma \bS[\pi_2] \ar[r] &
{|}\delta_{\cA(\pi_1)+\lambda_2}{|} \sma \bS[\pi_2]
\ar[r] & {|}\delta_{\cA(\pi_1)+\cA(\pi_2)+\lambda_2}{|} \\
{|}\delta_{\lambda_2, \lambda_1}{|} \sma \bS[\pi_1] \sma \bS[\pi_2]
\ar[u,"\htp"] \ar[d] \ar[r] & {|}\delta_{\lambda_2+\cA(\pi_1),
  \lambda_1+\cA(\pi_1))}{|} \sma \bS[\pi_2] 
\ar[u,"\htp"] \ar[d] \ar[r] & {|}\delta_{\cA(\pi_1) + \cA(\pi_2) +
  \lambda_2, \cA(\pi_1) + \cA(\pi_2) +\lambda_1}{|}
\ar[u,"\htp"] \ar[d] \\
{|}\delta^{\lambda_1}_{\lambda_2, \lambda_1}{|} \sma \bS[\pi_1] \sma \bS[\pi_2]
\ar[r] & {|}\delta^{\lambda_1 + \cA(\pi_1)}_{\lambda_2+\cA(\pi_1),
  \lambda_1+\cA(\pi_1))}{|} \sma \bS[\pi_2] 
\ar[r] & {|}\delta^{\cA(\pi_1) + \cA(\pi_2) + \lambda_1}_{\cA(\pi_1) +
  \cA(\pi_2) + \lambda_2, \cA(\pi_1) + \cA(\pi_2) +\lambda_1}{|} 
\\
\vdots \ar[u,"\htp"] \ar[d,swap,"\htp"] & \vdots \ar[u,"\htp"]
\ar[d,swap,"\htp"] & \vdots \ar[u,"\htp"] \ar[d,swap,"\htp"] \\
{|} \delta_{\lambda_1} {|} \sma \bS[\pi_1] \sma \bS[\pi_2] \ar[r] &  {|}
\delta_{\cA(\pi_1)+\lambda_1} {|} \sma \bS[\pi_2] \ar[r] & {|}
\delta_{\cA(\pi_1)+\cA(\pi_2)+\lambda_1} {|}
\end{tikzcd}
\end{equation}
coincides with the zig-zag
\begin{equation}
\begin{tikzcd}
{|}\delta_{\lambda_2}{|} \sma \bS[\pi_1 \pi_2] \ar[r] & 
{|}\delta_{\cA(\pi_1+\pi_2)+\lambda_2}{|} \\ 
{|}\delta_{\lambda_2, \lambda_1}{|} \sma \bS[\pi_1 \pi_2]
\ar[u,"\htp"] \ar[r] \ar[d] & 
{|}\delta_{(\pi_1+\pi_2)(\lambda_2, \lambda_1)}{|} \ar[u,"\htp"] \ar[d] \\
{|} \delta^{\lambda_1}_{\lambda_2, \lambda_1} {|} \ar[r] & {|} \delta_{\cA(\pi_1+\pi_2)+\lambda_1, (\pi_2 + \pi_1)(\Lambda)} {|} \\
\vdots \ar[u,"\htp"] \ar[d,swap,"\htp"] \ar[r] & \vdots \ar[u,"\htp"] \ar[d,swap,"\htp"] \\ 
{|} \delta_{\lambda_1} {|} \sma \bS[\pi_1 \pi_2] \ar[r] & {|}
\delta_{\cA(\pi_1 + \pi_2)+\lambda_1} {|}  
\end{tikzcd}
\end{equation}
up to the product on $\{\bS[-n]\}$.  

The analogous diagrams which express the associativity of the
action of $\Pi$ on these homotopy types, commute. 
\qed
\end{prop}

We now want to assemble these results into a description of
multiplicative structures on $|\delta|$.  To be precise, since $\Pi$
is discrete, to produce a map $\bbk[\Pi] \sma Z \to Z$ in the homotopy
category, it suffices to produce a map of monoids $\Pi \to
\Map_{\Ho(\Sp)}(Z,Z)$.  For this purpose, it is enough to produce
action maps $a_{\pi} \colon Z \to Z$ for each $\pi \in \Pi$ and show
that $a_{\pi_1} \circ a_{\pi_2} = a_{\pi_1 \pi_2}$ for all $\pi_1,
\pi_2 \in \Pi$; note that these action maps are maps in the homotopy
category.  Producing a map $\Sigma^{\deg} \bbk[\Pi] \sma Z \to Z$ in the
homotopy category is not substantially more difficult; considering the
construction of $\Sigma^{\deg} \bbk[\Pi]$ in terms of $\Pi$-graded
spectra, it suffices to produce shifted action maps $a_{\pi} \colon Z
\sma \bS[\pi] \to Z$ such that the diagram
\begin{equation}
\begin{tikzcd}
Z \sma \bS[\pi_1] \sma \bS[\pi_2] \ar[r,"\id \sma \mu"] \ar[d,swap,"a_{\pi_1} \sma \id"] & Z \sma \bS[\pi_1 \pi_2] \ar[d,"a_{\pi_1 \pi_2}"] \\
Z \sma \bS[\pi_2] \ar[r,swap,"a_{\pi_2}"] & Z 
\end{tikzcd}
\end{equation}
homotopy commutes.

Using this technique, we have the following result about the $\Pi$
action on $|\delta|$.

\begin{thm}\label{thm:Pi-homotopy-action}
The homotopy type $|\delta|$ as an orthogonal spectrum is a homotopy
module over $\Sigma^{\deg} \bbk[\Pi]$.  That is, there is an action
map 
\begin{equation}
\Sigma^{\deg} \bbk[\Pi] \sma |\delta| \to |\delta|
\end{equation}
which is associative and unital in the homotopy category. 
\end{thm}

\begin{proof}
For an element $\pi_1 \in \Pi$, we begin by constructing the action map.
Writing $I = \{k\epsilon \, | \, k \in \bZ\}$ and $I_1 = \cA(\pi_1) + I$,
we have a zig-zag 
\begin{equation}
\begin{tikzcd}
\displaystyle \holim_{k \in I} |\delta_k| \sma \bS[\pi] \ar[r,"\alpha_{\pi_1}"] & \displaystyle \holim_{k \in I_1 } |\delta_k| &
\displaystyle \holim_{k \in I \cup I_1} |\delta_k| \ar[l,swap,"\htp"] \ar[r,"\htp"] & \displaystyle \holim_{k \in I} |\delta_k|,
\end{tikzcd}
\end{equation}
where here the homotopy limit over $I \cup I_1$ denotes the diagram
with morphism zig-zags induced by the order on $I \cup I_1$ as well as
morphism zig-zags coming from $I$ and $I_1$.  Lemma~\ref{lem:zigassoc}
(along with the discussion of homotopy limits over zigzags in
Section~\ref{sec:rect}) shows that the two unlabeled maps in the
zig-zag are weak equivalences.

To see that the collection of these maps induce an action of
$\Sigma^{\deg} \Pi$ in the homotopy category, we fix elements $\pi_1$
and $\pi_2$ and consider the following diagram, where we write $I_2 =
\cA(\pi_2) + I$ and $I_{12} = \cA(\pi_2) + I_1 = \cA(\pi_1 \pi_2) + I$.
\begin{equation}
\begin{tikzcd}[row sep=small, column sep=small, transform shape, nodes={scale=.9}]
\displaystyle \holim_{k \in I} |\delta_k| \sma \bS[\pi_1] \sma \bS[\pi_2] \ar[r,"\alpha_{\pi_1}"]
\ar[dr,swap,"\alpha_{\pi_1 \pi_2}"] & \ar[d,"\alpha_{\pi_2}"]
\displaystyle \holim_{k \in I_1 } |\delta_k| \sma \bS[\pi_2] &
\displaystyle \holim_{k \in I \cup I_1} \ar[d,"\alpha_{\pi_2}"]
|\delta_k| \ar[l,swap,"\htp"] \ar[r,"\htp"] \sma \bS[\pi_2] &
\displaystyle \holim_{k \in I} |\delta_k| \ar[d, "\alpha_{\pi_2}"]
\sma \bS[\pi_2] \\
&\displaystyle  \holim_{k \in I_{12}} |\delta_k| & \displaystyle \holim_{k \in I_2 \cup I_{12}}
|\delta_k| \ar[l,swap,"\htp"]
\ar[r,"\htp"] & \displaystyle \holim_{k \in I_2} |\delta_k| \\
& & \displaystyle \holim_{k \in I \cup I_{12} \cup I_2} |\delta_k|
\ar[ul,"\htp"] \ar[r] 
\ar[dr,swap,"\htp"] \ar[u,swap,"\htp"]& \displaystyle \holim_{k
  \in I \cup I_2} |\delta_k| \ar[u,swap,"\htp"] \ar[d,"\htp"]\\
&&& \displaystyle \holim_{k \in I} |\delta_k|. 
\end{tikzcd}
\end{equation}
All of the subdiagrams commute: the top left triangle commutes by
Proposition~\ref{prop:action-systems}, the two squares commute by
Proposition~\ref{prop:action-assoc}, and the remaining subdiagrams
commute by the fact that restriction of diagrams is associative on
homotopy limits.  Composing with the natural equivalence $\bS[\pi_1]
\sma \bS[\pi_2] \to \bS[\pi_1 \pi_2]$, we conclude that $\alpha_{\pi_1
  \pi_2} = \alpha_{\pi_1} \alpha_{\pi_2}$ in the homotopy category.
Unitality is clear from Proposition~\ref{prop:action-systems} and
Proposition~\ref{prop:act_zigzag}.
\end{proof}

\begin{rem}\label{rem:cohact}
We can rectify these homotopy coherent actions to genuine actions by
considering a more elaborate generalization of
Lemma~\ref{lem:zigassoc}; since we do not need this for the main
applications of this paper, we leave this refinement for future work.
\end{rem}

Next, we want to extend this action to an action by the Novikov ring
$\Lambda$, the completion of $\Sigma^{\deg} \bbk[\Pi]$.  To do this,
we need to describe the interaction of the decreasing filtration on
$|\delta|$ with the decreasing filtration on $\Sigma^{\deg}
\bbk[\Pi]$.  Recall that the filtration on $\Sigma^{\deg} \bbk[\Pi]$
is given by
\begin{equation}
\lambda \mapsto \Sigma^{\deg} \bbk[\Pi_{> \lambda}]
\end{equation}
and the filtration on $|\delta|$ is defined in
Equation~\eqref{eqn:filtration} as
\begin{equation}
\lambda \mapsto \holim_{i \in I} |\delta_{i, >\lambda}|.
\end{equation}

The next lemma, combined with Corollary~\ref{cor:action-on-filtration}
and Proposition~\ref{prop:action-assoc}, show that these filtrations
are compatible with the action of $\Pi$ on $|\delta|$.  

\begin{lem}\label{lem:raw-module-structure}
For an element $\pi \in \Pi$ and $\lambda_1 > \lambda_0$, the
following diagram commutes
\begin{equation}
\begin{tikzcd}
\displaystyle{\holim_{i \in I_{\lambda_1}}}
|\delta_{\lambda_i, > \lambda_1}| \sma \bS[\pi] \ar[r,"\alpha_\pi"] \ar[d] &
\displaystyle{\holim_{i \in I_{\lambda_1+\cA(\pi)}}} |\delta_{\lambda_i, > \lambda_1 + \cA(\pi)}| \ar[d] \\
\displaystyle{\holim_{i \in I_{\lambda_0}}}
|\delta_{\lambda_i, > \lambda_0}| \sma \bS[\pi] \ar[r,"\alpha_{\pi}"] &
\displaystyle{\holim_{i \in I_{\lambda_0 + \cA(\pi)}}} |\delta_{\lambda_i, > \lambda_0 + \cA(\pi)}|, 
\end{tikzcd}
\end{equation}
where the vertical maps are the natural inclusions.\qed
\end{lem}

The argument for Theorem~\ref{thm:Pi-homotopy-action} generalizes to
establish that the action maps described in
Lemma~\ref{lem:raw-module-structure} yields a multiplication that is
compatible with the filtrations.

\begin{prop}\label{prop:filt-module-structure}
There are action maps in the homotopy category 
\begin{equation}
\Sigma^{\deg} \bbk[\Pi_{> \lambda_2}] \sma \holim_{i \in I}
|\delta_{i, > \lambda_1}| \to \holim_{i \in I} |\delta_{i, > \lambda_1
  + \lambda_2}| 
\end{equation} \qed
\end{prop}

The fact that the diagram in Lemma~\ref{lem:raw-module-structure}
strictly commutes implies that the action maps of
Proposition~\ref{prop:filt-module-structure} in fact assemble into an
action in the homotopy category of filtered spectra.  Note that is
essential here that we have constructed $\Sigma^{\deg} \bbk[\Pi]$ as a
filtered ring prior to passage to the homotopy category.

\begin{thm}
The action maps above give $|\delta|$ the structure of a filtered
module in the homotopy category of filtered spectra over the filtered
ring $\Sigma^{\deg} \bbk[\Pi]$. \qed
\end{thm}

Passing to completions, we can conclude the following.

\begin{thm}\label{thm:filtered-module}
The action maps above induce on $|\delta|$ the structure of a module
over the filtered ring $\Lambda$ in the homotopy category of filtered
spectra.
\end{thm}

\begin{proof}
As discussed in Section~\ref{appendix:complete} (e.g., see
Theorem~\ref{thm:compmonoid}), completion is a lax symmetric monoidal
functor on the homotopy category of filtered spectra, where the
category of complete filtered spectra is endowed with the completed
smash product of filtered spectra.  Since $|\delta|$ is complete for
the decreasing filtration as shown in Lemma~\ref{lem:delta-complete}, the
induced module structure in the complete category gives rise to one in
the category of filtered spectra.
\end{proof}

Since passage to the underlying object (i.e., the homotopy colimit of
the filtration) is a lax monoidal functor,
Theorem~\ref{thm:filtered-module} yields the following corollary. 

\begin{cor}\label{cor:Lambda-module}
The spectrum $|\delta|$ is a $\Lambda$-module in the homotopy
category.
\qed
\end{cor}

\subsection{Maps of homotopy types}
\label{sec:maps-homotopy-types}

This section discusses the methods that we shall use to compare the
filtered homotopy type, obtained by applying the methods of the
previous section in Floer theory, to the classical homotopy type
associated to the underlying symplectic manifold.  The same method
can be used to prove independence of the choice of auxiliary data,
e.g., the Hamiltonian and the almost complex structure, though we
shall not establish this independence result in this paper.

Comparison maps arise from the following typical construction in Floer
theory: we have a group $\Pi$, a pair of partially ordered sets
$\{\cP^i\}_{i=0,1}$ with free actions of the $\Pi$, equipped with action maps compatible with
the ordering, and associated flow cellular diagrams 
\begin{equation}
\delta^i \co C^*_{\rel \partial}(\cP^i; \Omega \bS) \to \kmod.
\end{equation}
The first step in understanding invariance is to be able to construct
maps of the associated homotopy types, so we begin with the data which
we will use to map $|\delta^0|$ to $|\delta^1|$.  Fix a constant $c
\in \bR$.

\begin{defin}
Let $\cP^{01} = \cP^0 \amalg \cP^1$
denote the partially ordered set where the ordering between elements
of $\cP^0$ and $\cP^1$ is unchanged, and the only new relations are
specified as follows:  
\begin{equation} \label{eq:bound_action_categorical_continuation}
  \parbox{30em}{ For $p_0 \in \cP^0, p_1 \in \cP^{1}$, $p_0 <
    p_1$ if and only if $\cA(p_0) \leq \cA(p_1) - c$.  }
\end{equation}
\end{defin}

We can define the relative cochains $C^*_{\rel \partial}(\cP^{01};
\Omega \bS)$ as above in Definition~\ref{defn:relcoch}.  Although the
action is no longer strictly compatible with the partial order, the
constructions work without modification.  (In addition, note that we
can always reindex $\cP_1$ in order to arrange for $c$ to be $0$.)

We now explain the class of modules $\delta$ we consider in this context.
Recall that for any spectral category $\cC$, given a $\cC$-module $M$
we can form the module $\Omega M$, which pointwise is given by the
formula $(\Omega M)(c) = \Omega M(c)$.

\begin{defin}\label{defn:categorical-continuation}
Let $\delta^0$ and $\delta^1$ be modules over $\cP^0$ and $\cP^1$
respectively.  We define \emph{categorical continuation data} from
$\delta^0$ to $\delta^1$ to consist of a cellular diagram
$\delta^{01}$ over $\cP^{01} $ whose restriction to $\cP^0$ coincides
with $\delta^0$, and whose restriction to $\cP^1$ coincides with
$\Omega \delta^1$. 
\end{defin}

A particularly interesting case is where we are given a
fixed partially ordered set $\cP$, and identifications $\cP^0 \cong \cP
\cong \cP^1$.  In this case, we assume that the constant $c$ in the
definition of the ordering on $\cP^{01}$ vanishes, and we consider a
fundamental chain $\delta$ on $\cP$, together with isomorphisms
$\delta^0 \cong \delta$ and $\delta^1 \cong \Omega \delta$.

\begin{defin}
A choice of categorical continuation data is \emph{unitriangular} if
we are in the situation above and for each pair $(p,q) \in \cP$ with
corresponding elements $q_0 \in \cP^0$ and $p_1 \in \cP^1$ such that 
$\cA(q_0) = \cA(p_1)$, the structure map 
\begin{equation}
C^*_{\rel \partial}(\cP^{01}, \bS)(q_0,p_1) \sma \delta^{0}(q_0) \to \Omega
\delta^1(p_1)
\end{equation}
vanishes unless $p=q$, in which case it is given by the equivalence
\begin{equation}
\Omega \bS \sma \delta(p) \to \Omega \delta(p)
\end{equation}
induced by the canonical identification
\begin{equation}
C^*_{\rel \partial}(\cP^{01}, \bS)(q_0,p_1) \htp \Omega \bS.
\end{equation}
\end{defin}

\begin{rem}
As suggested by the terminology, a unitriangular continuation datum
should be thought of as an upper triangular matrix with identities
along the diagonal.  In particular, note that the condition we impose
does not uniquely determine the categorical continuation data,
although there are of course a family of constraints imposed by the
composition maps.
\end{rem}

Given categorical continuation data $\delta^{01}$, we define the
homotopy type $|\delta^{01}|$ in analogy with
Definition~\ref{defn:invht}.

\begin{defin}
For each real number $\lambda$, $\cP^{01}_{\lambda}$ is the quotient
of $\cP^{01}$ by elements $p \in \cP_0$ such that $\cA(p) > \lambda$
and $q \in \cP_1$ such that $\cA(q) > \lambda + c$.  We define
$\delta^{01}_{\lambda}$ to be the functor induced by $\delta^{01}$ by
restriction, as in Section \ref{sec:homotopy-type-sub}.  The geometric
realization is then the homotopy limit
\begin{equation}
|\delta^{01}| = \holim_{k \in \bZ} |\delta^{01}_{k \epsilon}|,
\end{equation}
Corollary~\ref{cor:Lambda-module} implies that
$|\delta^{01}|$ is a $\Lambda$-module in the homotopy
category.
\end{defin}

For each $\lambda$, there is an inclusion $\cP^{1}_{\lambda} \to
\cP^{01}_{\lambda}$, which induces a comparison map on homotopy
types.

\begin{prop}
Let $\delta^{01}$ be categorical continuation data for $\delta^0$ and
$\delta^1$.  Then for each $\lambda$, there are maps
\begin{equation}
\iota_{1,\lambda} \colon \Omega |\delta^1_{\lambda+c}| \to
|\delta^{01}_{\lambda}| 
\end{equation}
induced by the inclusion of $\cP^{1}_\lambda$ in $\cP^{01}_\lambda$.
On passage to homotopy limits, these maps induce a map of spectra 
\begin{equation}
\iota_1 \colon \Omega |\delta^1| \to |\delta^{01}|.
\end{equation}
There are action maps in the homotopy category
\begin{equation}
\Omega |\delta^1_{\lambda + c}| \sma \bS[\pi] \to
|\delta^{01}_{\lambda + \cA(\pi)}|
\end{equation}
which induce a map of filtered $\Lambda$-modules
\begin{equation}
\Sigma^{\deg} \bbk[\Pi_{> \lambda_1}] \sma \Omega |\delta^1_{\lambda_2
  + c}| \to
|\delta^{01}_{> \lambda_1 + \lambda_2}|.
\end{equation}
As a consequence, the map $\iota_1$ is a map of $\Lambda$-modules.
\end{prop}

\begin{proof}
Since there are no maps in $\cP^{01}$ from objects in $\cP_1$ to
objects in $\cP_0$, the inclusion of posets $\cP^{1}_{\lambda+c}
\to \cP^{01}_{\lambda}$ induces a map of simplicial spectra
\begin{equation}
\iota_1 \colon B_{\bullet}(\bS^{1}_{\lambda+c}, C^*_{\rel \partial}(\cP^{1}_{\lambda+c},
\Omega \bS), \Omega \delta^1_{\lambda+c}) \to B_{\bullet}(\bS^{01}_{\lambda}, 
C^*_{\rel \partial}(\cP^{01}_{\lambda}, \Omega \bS), \delta^{01}_{\lambda}).
\end{equation}
These maps are clearly compatible with the zig-zags in the homotopy
limit system and therefore induce a map $\Omega |\delta^1| \to
|\delta^{01}|$ on passage to the homotopy limits.  Moreover, a
straightforward check using the constructions of the previous section
shows that these maps are compatible with the filtration (in the
homotopy category) and induce a map of $\Lambda$-modules in the
homotopy category.
\end{proof}

We can identify the cofiber of the map $\iota_1$ in terms of
$|\delta^0|$.  A standard difficulty with working with modules in the
homotopy category is that cofibers do not automatically inherit module
structures.  In our case, we exploit the fact that we can directly
obtain the module structures on the cofiber of $\iota_1$ by working
with a concrete construction.

\begin{thm}\label{thm:cofiber}
There are homotopy cofiber sequences
\begin{equation}
\Omega |\delta^1_{> \lambda+c}| \to |\delta^{01}_{> \lambda}| \to 
|\delta^0_{> \lambda}|,
\end{equation}
and on passage to homotopy limits a homotopy cofiber sequence of
homotopy $\Lambda$-modules 
\begin{equation}
\Omega |\delta^1| \to |\delta^{01}| \to |\delta^0|.
\end{equation}
\end{thm}

\begin{proof}
Since homotopy (co)fibers commute with homotopy inverse limits, it
essentially suffices to consider what happens at a fixed $\lambda$.
Moreover, as homotopy cofibers commute with geometric realization, it
suffices to analyze the geometric realization of the simplicial object
produced by the levelwise cofiber.  The cofiber is easy to compute
levelwise; any simplex corresponding to a path that lies completely in
the image of $\cP^1_{\lambda}$ in $\cP^{01}_{\lambda}$ is contracted
to $\ast$.  As a consequence, we can describe the homotopy cofiber as
the geometric realization of the bar construction
\begin{equation}
B_\bullet(\bS^{01}_{\lambda}, C^*_{\rel
  \partial}(\cP^{01}_{\lambda}, \Omega \bS), \widetilde{\delta^{0}_{\lambda}}),
\end{equation}
where $\widetilde{\delta^0_{\lambda}}$ denotes the module on $C^*_{\rel
  \partial}(\cP^{01}_{\lambda}, \Omega \bS)$ which is the restriction
of $\delta^0$ on $p \in \cP^{0}_{\lambda}$ and $\ast$ otherwise.  This
description of the homotopy cofiber is clearly equipped with a natural
action of $\bbk[\Pi_{> \lambda}]$ that is induced from the maps
\begin{multline}
B_\bullet(\bS^{01}_{\lambda}, C^*_{\rel
  \partial}(\cP^{01}_{\lambda}, \Omega \bS),
\widetilde{\delta^{0}_{\lambda}}) \sma \bS[\pi] \to \\ 
B_\bullet(\bS^{01}_{\lambda+\cA(\pi)}, C^*_{\rel
  \partial}(\cP^{01}_{\lambda+\cA(\pi)}, \Omega \bS),
\widetilde{\delta^{0}_{\lambda+\cA(\pi)}}),
\end{multline}
and which is compatible with the maps $\iota_1$ as maps of modules.
Therefore there is an induced $\Lambda$-action after completion.

To identify this homotopy cofiber, we now proceed as in the proof of
Proposition~\ref{prop:spine}.  Specifically, the analogue of the
construction of Lemma~\ref{lem:zig} produces a simplicial map from the
homotopy cofiber to $|\delta^0|$.  Since the action of $\bbk[\Pi_{>
    \lambda}]$ is compatible with composition and with the restriction map,
this simplicial map is compatible with the action of $\bbk[\Pi_{>
    \lambda}]$ and induces a map of $\Lambda$-modules on passage to
completions.  On the other hand, the argument for
Proposition~\ref{prop:spine} shows that this map is a weak
equivalence.  We conclude that there is a natural equivalence from the
homotopy cofiber to $|\delta^0|$ which is a filtered map of homotopy
$\Lambda$-modules.
\end{proof}

As an immediate corollary, the connecting map yields the desired
comparison maps.

\begin{cor}\label{cor:contmap}
A choice of categorical continuation data determines maps 
\begin{equation}
|\delta^0_{\lambda}| \to |\delta^1_{\lambda + c}|
\end{equation}
and 
\begin{equation}
|\delta^0_{> \lambda}| \to |\delta^1_{> \lambda + c}|
\end{equation}
that are compatible with the filtered action of $\Sigma^{\deg}
\bbk[\Pi]$ and on passage to completion maps of homotopy
$\Lambda$-modules
\begin{equation} \label{eq:induced_map_continuation}
|\delta^0| \to |\delta^1|
\end{equation} 
with homotopy cofiber a homotopy $\Lambda$-module equivalent to
$\Sigma |\delta^{01}|$.  \qed
\end{cor}

Next, when the categorical continuation data is unitriangular, the
homotopy type is contractible.

\begin{prop}\label{prop:contractible}
Let $\delta_{01}$ be unitriangular categorical continuation data for
$\delta_0$.  Then $|\delta_{01}|$ is contractible. In particular, the
induced map in Equation~\eqref{eq:induced_map_continuation} is an
equivalence.  
\end{prop} 

\begin{proof}
This follows from a variation of the filtration argument used to
establish Proposition~\ref{prop:nullhomtype}.  Specifically, we
consider the same basic argument, but now we additionally filter the
objects with equal action and inductively reduce to the case of the bar
construction $B(\bS_{\lambda}^{01}, C^*_{\rel
  \partial}(\cP^{01}_\lambda; \Omega \bS), \delta^{01}_{\lambda, p})$
where $\delta^{01}_{\lambda,p}$ is the module which is nontrivial on
the two copies of $p$ (in $\cP^0$ and $\cP^1$) and $\ast$ otherwise.
In this case the argument for Lemma~\ref{lem:justp} shows that the
homotopy type reduces to the cofiber of the action map $\Omega \delta(p)
\to \Omega \delta(p)$.  Since this homotopy cofiber is contractible by
definition, the contribution to the associated graded is trivial.  
\end{proof}

Next, we shall describe additional categorical data that suffices to
prove that the map in Equation~\eqref{eq:induced_map_continuation}
splits.  Consider a pair 
$(\cP^{00},\delta^{00})$ and $(\cP^{01}, \delta^{01})$ of cellular
diagrams, and categorical continuation data $\cP^{01}$ representing a
map $|\delta^{00}| \to |\delta^{01}|$.

\begin{defin}
\emph{Categorical retraction data} for the map $\delta^{00} \to
\delta^{01}$ consists of:
\begin{enumerate}
\item A category $\cP^{\Box}$
with objects the disjoint union $\cP^{00} \coprod \cP^{01} \coprod
\cP^{10} \coprod \cP^{11}$, where $\cP^{10} = \cP^{00}$ and $\cP^{11}
= \cP^{10}$, and morphisms indicated by the diagram
\begin{equation}\label{eq:retract-box}
\begin{tikzcd}
\cP^{01} \ar[r] & \cP^{11} \\
\cP^{00} \ar[r] \ar[u] & \cP^{10}, \ar[u] 
\end{tikzcd}
\end{equation}
where each arrow specifies morphisms as in the definition of
categorical continuation data for the domain and codomain, such that
the analogue of
Equation~\eqref{eq:bound_action_categorical_continuation} holds (i.e.,
there is a constant $c$ so that there is no morphism in $\cP^{\Box}$ from $p$
to $q$ unless $\cA(p) \leq \cA(p) - c$).
\item A module $\delta_\Box$ on $C^*_{\rel \partial}(\cP^{\Box},
  \Omega \bS)$ such that:
\begin{enumerate}
\item The restriction to any adjacent pair of summands (i.e., terms
  connected by arrows in Equation~\eqref{eq:retract-box}) in the
  disjoint union specifies categorical continuation data as do the
  pairs $(\cP^{00} \coprod \cP^{01}, \cP^{10} \coprod \cP^{11})$, 
  $(\cP^{00} \coprod \cP^{10}, \cP^{01} \coprod \cP^{11})$, and
  $(\cP^{00} \coprod \cP^{10} \coprod \cP^{01}, \cP^{11})$.
\item The restriction to the pairs $\cP^{00} \to \cP^{10}$ and $\cP^{10} \to
  \cP^{11}$ specifies unitriangular categorical continuation data.
\end{enumerate}
\end{enumerate}
\end{defin}

In the following, we will write $\delta^{00,01}$ for the restriction
of $\delta^{\Box}$ to the subdiagram spanned by $\cP^{00} \coprod
\cP^{01}$, and analogously for other subsets of the sum.  The
existence of categorical retraction data imposes fairly stringent
conditions on the modules obtained by restricting to various
subdiagrams.  In particular, if we regard $\delta^{00,01}$ as
categorical continuation data for $\delta^0$ and $\delta^1$, then 
\begin{enumerate}
\item $\delta^{00,10}$ is unitriangular
categorical continuation data for $\delta^0$ and $\delta^0$,
\item $\delta^{10,11}$ is unitriangular categorical continuation data for
$\Omega \delta^0$ and $\Omega \delta^0$, 
\item and $\delta^{01,11}$ is categorical continuation data for
  $\Omega \delta^1$ and $\Omega \delta^0$.
\end{enumerate}

Therefore, we can think of the up and over direction in the square as
representing a composite $\delta^0 \to \delta^1 \to \delta^0$ and the
over and up direction in the square as representing a composite
$\delta^0 \to \delta^0 \to \delta^0$ of identity maps.  
The idea for the next proposition is that this data represents a
commutative diagram.

\begin{prop}\label{prop:retract} 
Given categorical retraction data, the map of homotopy
$\Lambda$-modules $f \colon |\delta^{0}| \to
|\delta^{1}|$ represented by the categorical continuation
data $\delta^{00,01}$ is split.  For fixed $\lambda$, this splitting
restricts to a composite
\begin{equation}
|\delta^{0}_{\lambda}| \to |\delta^{1}_{\lambda+c}| \to |\delta^{1}_{\lambda+c}|
\end{equation} 
and a composite 
\begin{equation}
|\delta^{0}_{> \lambda}| \to |\delta^{1}_{> \lambda+c}| \to
|\delta^{1}_{> \lambda+c}|
\end{equation} 
on the restricted homotopy limits.
\end{prop}

\begin{proof}
Recall that a cofiber sequence $X \to Y \to C \to \Sigma X$ in any
stable category splits if the map $C \to \Sigma X$ is null-homotopic.
(In fact, this is true in any triangulated category.)

There is a cofiber sequence
\begin{equation}\label{eqn:box2}
\Omega^2 |\delta^0| \to |\delta^{\Box}| \to |\delta^{00,10,01}|,
\end{equation}
and the inclusion $\Omega^2 |\delta^0| \to |\delta^{\Box}|$ factors as both the
composites $\Omega^2 |\delta^0| \to |\delta^{01,11}| \to |\delta^{\Box}|$ and
$\Omega^2 |\delta^0| \to |\delta^{10,11}| \to |\delta^{\Box}|$.
Since $|\delta^{10,11}| \htp \ast$ by hypothesis, we see that this map
is null-homotopic.  

Moreover, since the cofiber of the inclusion
$|\delta^{01,11}| \to |\delta_{\Box}|$ is $|\delta^{00,10}| \htp
\ast$, we see that the map $\Omega^2 |\delta_0| \to |\delta^{01,11}|$
is itself null-homotopic.  First, the cofiber
sequence 
\begin{equation}
\Omega^2 |\delta^1| \to \Omega^2 |\delta_0| \to |\delta^{01,11}| 
\end{equation}
implies that the map $\Omega^2 |\delta^1| \to \Omega^2 |\delta_0|$ is a
retraction, i.e., that the cofiber sequence
\begin{equation}
\Omega |\delta^{01,11}| \to \Omega^2 |\delta^1| \to \Omega^2 |\delta_0|
\end{equation}
is split.  Rewriting, this splitting exhibits a weak equivalence
\begin{equation}
|\delta^1| \htp |\delta_0| \vee \hocofib(f).
\end{equation}

To see that $f$ is split by $g$, we observe that we have a map 
\begin{equation}
|\delta_{01,11}| \to |\delta_{\Box}| \to |\delta_{00,01}|,
\end{equation}
where the composite map can be described as the canonical map
\begin{equation}
\Omega \hocofib(g) \to \hocofib(f).
\end{equation}
As this map is an equivalence, we can conclude that $f$ is split with
retraction $g$.

The assertions about the restriction to fixed $\lambda$ are
straightforward, since the maps are compatible with $\lambda$ and with
the zig-zags in the restricted homotopy limit.

\end{proof}

\part{Virtual fundamental chains from Kuranishi flow categories}
\label{part:orient-kuran-flow}
\section{Kuranishi flow categories and orientations}
\label{sec:kuran-flow-categ}

\subsection{Kuranishi flow categories}
\label{sec:kuran-flow-categ-2}

The purpose of this section is to formalise the idea of a lift of a
topological flow category to a flow category with morphism spaces
equipped with a variant of the notion of Kuranishi structure
introduced by Fukaya and Ono~\cite{FukayaOno1999}.  From our point of
view, such a lift is a further refinement of the definition of an
orbispace flow category introduced in
Section~\ref{sec:orbisp-flow-categ}, and we shall make use of that
notion in this section. As discussed in the introduction, we find the
various existing formalisms of Kuranishi structures ill-adapted for
the formal constructions which are required in this paper;
Definition~\ref{def:Kuranishi-flow-category} implements our desired
notion of a lift.

\subsubsection{Kuranishi Charts}
\label{sec:anoth-categ-charts}

We begin by elaborating on the notion of a Kuranishi chart given in
the introduction (Definition~\ref{def:Kuranishi-chart}), and
introducing a category for which these are the objects.

\begin{defin}\label{defn:Kuranishi-chart}
A \emph{Kuranishi chart} is a quadruple $(X, V, s, G)$ consisting of
the following data: 
\begin{enumerate}
\item (Symmetry group) a finite group $G$,
\item (Thickened chart) a $G$-manifold $X$ (paracompact and Hausdorff, and possibly with boundary),
\item (Obstruction space) a finite dimensional $G$-representation $V$ equipped with an invariant inner product, and
\item (Defining section) a $G$-equivariant map $s \co X \to V$.
\end{enumerate}
We write $Z = s^{-1}(0)$ for the \emph{zero locus} (which is locally
compact), and define the \emph{footprint} of the chart to be the
quotient space $Z/G$.
\end{defin}

We define the \emph{boundary of a Kuranishi chart} to be the chart
\begin{equation}
\partial \bX \equiv  (\partial X,V,s, G)
\end{equation}
We say that a Kuranishi chart is \emph{without boundary} if the
boundary of $X$ is empty.  By convention, we require that $X$ have
pure dimension, so that we can assign to each Kuranishi chart a
\emph{virtual dimension} 
\begin{equation}
  \dim \bX = \dim X - \dim V.  
\end{equation}

\begin{rem}
In our definition of $G$-manifold $X$, we require the condition that
the $G$-action be \emph{locally Euclidean,} i.e., that there is a
$G$-invariant neighbourhood of each orbit which is $G$-equivariantly
identified with a neighbourhood of an orbit in a
$G$-representation. This implies that $X$ is a $G$-ENR (see
e.g.,~\cite{tomDieck1987}), and in particular that $X/G$ has the
homotopy type of a $CW$-complex. Note that an example of
Bing~\cite{Bing1952} shows that there are $\bZ/2\bZ$ actions on $S^3$
which do not satisfy this property, and results of
Quinn~\cite[Proposition 2.1.4]{Quinn82} show that there are action of
finite groups on discs which are locally Euclidean, but are not
conjugate to smooth actions.
\end{rem}

Definition~\ref{defn:Kuranishi-chart} specifies the objects of the
category of Kuranishi charts.  We now define the morphisms.

\begin{defin}
A map $f \co \bX \to \bX'$ of Kuranishi charts is given by the following data:
\begin{enumerate}
\item a homomorphism $G \to G'$,
\item an isometric embedding $V \to V'$ which is $G$-equivariant, and 
\item a $G$-equivariant map $X \to X'$ preserving the boundary, which commutes with the defining sections.
\end{enumerate}
Denoting by  $G_f^\perp$ the kernel of the map of groups, and by
$V_f^\perp$ the quotient $V'/V$, we require the following properties
to hold: 
\begin{enumerate}
\item the action of  $G_f^\perp$ on $X$ is free and the map from the quotient to $X'$ is an open embedding in the inverse image of $V \subset V'$, and
     \item  near each point in $X'$ lying in the image of $X$ under
       $f$, there is a product chart $U(X)/G_f^\perp \times
       U(V_f^\perp )$ with $U(X)$ an open subset in $X$ and
       $U(V_f^\perp)$ an open subset in $V_f^\perp$, such that the
       following diagram commutes:
\begin{equation}\label{eq:transverse_cut_out}
  \begin{tikzcd}
    U(X) \times U(V_f^\perp) \arrow{r} \arrow{d} & X'  \arrow{d} \\
X \times V_f^\perp \arrow{r} & V_f^\perp.  
  \end{tikzcd}        
\end{equation}
  \end{enumerate}
\end{defin}
The last condition above amounts to the requirement that the projection $X' \to V_f^\perp$ be a topological submersion near the image of $X$, with fibre containing $X/G_f^\perp$ as an open subset.
\begin{rem} \label{rem:why-we-need-underline_V}
The requirement of topological submersion is stronger than that of topological transversality used e.g., in~\cite{Pardon2016}. In the smooth case, it can be replaced by a condition on tangent spaces. 
\end{rem}
\begin{rem} \label{rem:surjection_groups_Kuranishi-charts}
  Starting in Section \ref{sec:impl-atlas-with}, we shall assume that the group homomorphism $G \to G'$ is surjective, as this simplifies various constructions, and the outcome of Part \ref{part:floer-theor-constr} is that the output of Floer theory are maps of Kuranishi charts satisfying this property. The first instance where this condition is used is discussed in Remark \ref{rem:surjection-implies-bundle-trivial} below.
\end{rem}
Composition is defined in the obvious way; the only point to check is that the two properties we impose are preserved.
\begin{lem}
Given maps $\bX \overset{f}{\to} \bX' \overset{g}{\to} \bX''$, the kernel $G_{g \circ f}^\perp$ acts freely on $X$, and the composition $X'' \to V_{g \circ f}^\perp$ is a topological submersion near the image of $X$.
\end{lem}
\begin{proof}
If a non-trivial element of $G_{g \circ f}^\perp$ lies in $G_f^\perp$, it acts freely by the condition on the map $\bX \to \bX'$; otherwise, it maps to a non-trivial element of $G_g^\perp$, which must act freely on $X'$ by the assumption on $\bX' \to \bX''$. Since the map $X \to X'$ is equivariant, we conclude that there are no fixed points in $X$.

To check the second property, observe that the product $U(X) \times (U(V_f^\perp) \times U(V_g^\perp))$ gives the desired chart.
\end{proof}

From the above definition, it is apparent that the virtual dimensions
of $\bX$ and $\bX'$ must agree if there is a map between them, and so
the category of Kuranishi charts decomposes as the disjoint union of
categories indexed by the integers.

\begin{rem}
One natural generalisation is to allow for $V$ to be a vector bundle, and our constructions can be carried out with only minor modifications in this context. It is theoretically possible to allow $V$ to be a $\TOP$-microbundle \cite{Milnor1964}, but the authors are not aware of a context where such generality arises naturally, and we expect that a substantial modification of our methods would be required in this setting to account for appropriate formulations of transversality.  
\end{rem}

\begin{prop}
The category of Kuranishi charts has a monoidal structure, given by
the natural product of Kuranishi charts, which assigns to a pair $\bX
= (X, V,s, G) $ and $\bX' = (X', V' ,s', G') $ the chart 
\begin{equation}
 \bX \times \bX' \equiv  (X \times X', V \oplus V' ,s \oplus s', G \times G'). 
\end{equation}
This product is naturally compatible with maps of Kuranishi charts, and
the unit is the chart $(\ast, \{0\}, s, \{e\})$, where $s$ is the
unique map taking $\ast$ to $0$.
\end{prop}

\begin{rem}
It is straightforward to see that $\bX \times \bX' $ is naturally
isomorphic to $\bX' \times \bX$, i.e., that the monoidal structure on
the category of Kuranishi charts is in fact symmetric.  We shall never
appeal to such symmetries in this paper, so we omit the corresponding
discussion.
\end{rem}

\begin{defin}
We have a \emph{zero locus} functor from the category of Kuranishi
charts to $\Orbo$ given by 
\begin{equation} \label{eq:pass_to_zero}
 (X, V ,s, G)    \mapsto  (s^{-1}(0),G).
\end{equation}
Composing the zero locus functor with the quotient functor $\Orbo \to
\Top$ given by the assignment $(Z,G) \mapsto Z/G$, we obtain the
\emph{footprint functor} from the category of Kuranishi charts to
$\Top$ given by 
\begin{equation}
(X, V, s, G) \mapsto s^{-1}(0)/G.
\end{equation}
(Here we use the terminology introduced by McDuff and
Wehrheim~\cite{McDuff2017}.)

The compatibility of boundaries with maps of Kuranishi charts yields an
endofunctor 
\begin{equation}
\bX \mapsto \partial \bX = (\partial X,V,s, G)
\end{equation}
which assigns to a chart its boundary.
\end{defin}

\begin{rem} \label{rem:hint_at_geometry-Kuranishi}
In applications, $V$ is a choice of inhomogeneous data (obstruction bundle) for a Cauchy-Riemann equation on a family of Riemann surfaces, and  $X$ is an associated moduli space of maps from a Riemann surface with marked points satisfying geometric constraints. The geometric setup allows for increasing the number of marked points (with additional constraints), and enlarging the obstruction space. While enlarging the obstruction space from $V_0$ to $V_1$ yields a map of moduli spaces, considering a setup where the moduli spaces $X_0$ and $X_1$ essentially differ only in the number of marked points yields a correspondence $ X_0  \leftarrow X \to X_1  $. In order to arrive at the abstract setting we are considering, we shall use the fact that the correspondence $X$ is also a Kuranishi chart. This procedure explains some of the combinatorial complexity in Part \ref{part:floer-theor-constr}, which we have traded for the straightforward functoriality of our notion of maps of Kuranishi charts.
\end{rem}

\subsubsection{Stratified orbispace presentations}
\label{sec:strat-orbisp-pres}

Let $\cS$ be a partially ordered set. The following notion will be
essential for the construction of this paper: 
\begin{defin}\label{defn:stratified-space}
An \emph{$\langle \cS \rangle$-stratification} of a topological space
$\cM$ is an assignment of a closed subset $\partial^q \cM$ to each
element $q \in \cS$, such that 
\begin{equation}
\partial^q \cM \cap \partial^r \cM = \emptyset    
\end{equation}
whenever $q$ and $r$ are not comparable.

A map $f \colon \cM \to \cB$ of $\langle \cS \rangle$-stratified spaces is
a continuous map $\cM \to \cB$ such that $f(\partial^q \cM) \subseteq
\partial^q \cB$. 
\end{defin}
The above notion leads to a stratification in the usual sense, with
the stratum $\partial^Q \cM$ associated to each subset $Q$ of $\cS$
being empty unless $Q$ is totally ordered, in which case we set   
\begin{equation}
\partial^Q \cM \equiv \bigcap_{q \in Q} \partial^q \cM.
\end{equation}
When $Q = \emptyset$, we will interpret $\partial^Q \cM$ to be $\cM$.  It
is clear that for a stratified map $f \colon \cM \to \cB$ we have an
inclusion $f(\partial^Q \cM) \subseteq \partial^Q \cB$ for all $Q$. 

In this setting, we shall find it convenient to define a notion of
orbispace presentation that involves choices of charts for each
possible stratum: 

\begin{defin}
The category $\Orb \langle \cS \rangle$ of \emph{$\langle \cS
  \rangle$-stratified orbispace charts} is the category with
\begin{enumerate}
\item objects the pairs $(Z, G)$,
where $Z$ is an $\langle \cS \rangle$-stratified topological space,
and $G$ a finite group acting on $Z$ via stratification-preserving
maps, and 
\item the morphisms are specified by a homomorphism $p\colon G_0 \to G_1$
and a $G_0$-equivariant stratification-preserving map $Z_0 \to p^*
Z_1$.
\end{enumerate}
\end{defin} 

Note that the case $S = \emptyset$ yields the category
$\Orb^\emptyset$ discussed in Section~\ref{sec:orbisp-pres}, and that
we have a functor 
\begin{equation}\label{eq:functor_forget_stratification}
\Orb \langle \cS \rangle \to \Orb^\emptyset
\end{equation}
which forgets the stratification.

\begin{notn}
Given a totally ordered subset $Q$ of $\cS$, we obtain a new partially
ordered subset 
\begin{equation}
  \partial^Q \cS \subset \cS
\end{equation}
consisting of all elements of $\cS$ that are comparable to every
element of $Q$, but do not lie in $Q$.  We understand
$\partial^{\emptyset} \cS$ to be $\cS$.  
\end{notn}

This notation is chosen to make the following lemma hold.

\begin{lem}
If $\cM$ is an $\langle \cS \rangle$-stratified space, then the stratum
$\partial^Q \cM$ is naturally a $\langle \partial^Q \cS
\rangle$-statified space.  In addition, for $Q' \subseteq Q$
totally-ordered subsets of $\cS$, we have
\begin{align}
  \partial^{Q} \cS & = \partial^{Q \setminus Q'} \left(\partial^{Q'} \cS\right) \\
  \partial^{Q} \cM & = \partial^{Q \setminus Q'} \left(\partial^{Q'} \cM\right),
\end{align}
respectively considered as subsets of $\cS$ and $\cM$. In particular, when $\cM'$ is a $\langle \partial^{Q'} \cS
\rangle$-stratified space, we can regard $\partial^{Q \setminus Q'} \cM'$
as a $\langle \partial^{Q} \cS \rangle$-stratified space by restriction.
\end{lem}
\begin{proof}
For $q \in \partial^Q \cS$, we set $\partial^q (\partial^Q \cM) =
\partial^q \cM \cap \partial^Q \cM$; the first assertion is clear.  For the second one, we note that an element of $\cS$ is comparable to every element of $Q$ if and only if it is comparable to every element of $Q'$ and $Q \setminus Q'$, which implies the remaining statements.
\end{proof}

We now define the category of orbispace charts in this context.

\begin{defin}\label{defin:stratorbicharts}
The category $\Orb$ of \emph{stratified orbispace charts} is the category
whose objects are triples $(\cS, Z, G)$, where $\cS$ is a
partially-ordered set,  and $(Z,G)$ is an object of $\Orb \langle \cS
\rangle$.

A morphism $(\cS, Z, G)$ to $(\cS', Z', G')$ is determined by an
order-preserving isomorphism $\rho \colon \cS \cong \partial^{Q} \cS'$
for some totally ordered subset $Q$ of $\cS'$ and a 
morphism
\begin{equation}
  f \colon (Z,G) \to (\partial^{Q} Z', G'),
\end{equation}
in $\Orb \langle \cS \rangle$, where we regard $\partial^{Q} Z'$ as a
$\langle \cS \rangle$-stratified space via $\rho$.

Given morphisms 
\begin{align}
(f,\rho) \co (\cS, Z,G) & \to (\cS', \partial^{Q} Z',G') \\
(g, \rho') \co (\cS', Z',G') & \to (\cS'', \partial^{Q'} Z'', G''),
\end{align}
the composite morphism $(g \circ f, \rho' \circ \rho)$ is specified as follows: we
let 
\begin{equation}
Q^{g \circ f} = \rho'(Q) \amalg Q',
\end{equation}
where by $\rho'(Q)$ we abusively mean the composition $Q \to \cS'
\cong \partial^{Q'} \cS'' \to \cS''$.  Since $\partial^{Q'} \cS'' $
consists of elements of $\cS''$ which are comparable to all elements
of $Q'$, the union $\rho'(Q) \amalg Q' $ is totally ordered.  Then the
isomorphism $\cS \cong \partial^{\rho'(Q) \amalg Q'} \cS''$ is
determined by $\rho$ and $\rho'$, and the map $(Z,G) \to (Z'', G'')$
is the composite of of $f$ with the restriction of $g$.
\end{defin}

The category of stratified orbispace charts inherits a monoidal
structure from its constituent components.

\begin{prop}
The category $\Orb$ has a monoidal structure, where the product of
$(\cS, Z, G)$ and $(\cS', Z', G')$ is given by $(\cS \amalg \cS', Z
\times Z', G \times G')$.  The unit is given by $(\emptyset, \ast,
\{e\})$.
\end{prop}

In many of our geometric applications, the stratified
orbispace charts we deal with live in the simpler subcategory of
$\Orb$ where morphisms are determined by actual equalities $\cS =
\partial^{Q} \cS'$.

\begin{defin}
The category $\Orb^{\iso}$ is the subcategory of $\Orb$ with the same
objects but where the morphisms are specified by identities $\cS =
\partial^{Q} \cS'$.  For a given $\cS$, the category $\Orb^{\cS}$ is the
full subcategory of $\Orb^{\iso}$ spanned by the objects $(\partial^Q
\cS, Z, G)$ as $Q$ varies over the totally-ordered subsets of $\cS$.
\end{defin}

We can now define the notion of an orbispace presentation for a
stratified space:

\begin{defin}
A \emph{stratified orbispace presentation} of an $\langle \cS
\rangle$-stratified space $\cM$ consists of
  \begin{enumerate}
  \item A small category $A$, equipped with a functor $(Z,G) \co A \to
    \Orb^{\cS}$, and
    \item a continuous map
\begin{equation}
\colim_{\alpha \in A} Z_\alpha/G_\alpha \to \cM
\end{equation}
which is a stratified homeomorphism of $\langle \cS
\rangle$-stratified spaces, where the colimit is taken in the category of topological
spaces, which is equipped with the induced stratification.
\end{enumerate}
We require that the restriction of the functor specified in the first
condition to the subcategory of $A$ with image in the category of
$\langle \partial^Q \cS \rangle$-charts be an orbispace presentation
of $\partial^Q \cM$ (note that the second condition requires that each
object of $A$ maps to such a chart).
\end{defin}
We remind the reader that the condition of being an orbispace
presentation is a condition of contractibility of the nerve of the
subcategory of charts whose image contain any given point of the
underlying space. We can formulate this notion in exactly the same way
in the stratified context, using the forgetful functor from
Equation~\eqref{eq:functor_forget_stratification}.

One should think of a stratified orbispace presentation of $\cM$ as a
collection of orbispace presentations of all of the strata of $\cM$,
functorially depending on the choice of stratum.

\subsubsection{$\langle \cS \rangle$-Kuranishi charts}
\label{sec:langle-cs-rangle}

Except for the fact that we allow partially ordered sets more general than the natural numbers, and that we work with topological manifolds rather than smooth ones, the following notion is essentially equivalent to the notion considered by J\"anich \cite{Janich1968}:
\begin{defin}
An \emph{$\langle \cS \rangle$-manifold} is a manifold $X$ with boundary, which is $\langle \cS \rangle$-stratified in the sense of the previous section,  such
that each stratum $\partial^Q X$ is a manifold with boundary admitting a
neighbourhood which is homeomorphic to 
\begin{equation} \label{eq:local_stratification_manifold}
  \partial^Q X \times (-\infty,0]^Q,
\end{equation}
via a homeomorphism that preserves strata.
\end{defin}

To clarify the above definition, if $P$ is a subset of $Q$, the
intersection of $\partial^P X$ with the image of the chart in Equation~\eqref{eq:local_stratification_manifold} is assumed
to be the product of $ \partial^Q X \times (-\infty,0]^{Q \setminus
    P}$.

\begin{defin} 
We denote by $\Kur\langle \cS \rangle $ the category of $\langle \cS
\rangle$-Kuranishi charts: the objects $\bX = (X,V,s, G)$ are as
before, except that $X$ is a $\langle \cS \rangle$-manifold, and the
$G$-action preserves the stratification.  Morphisms are required to
preserve the strata and the witness to the topological submersion in
Equation \eqref{eq:transverse_cut_out} is required to be a stratified
map.
\end{defin}

The boundary functor of Kuranishi charts has a stratified analogue.

\begin{lem}
Let $S$ be a partially-ordered set and let $Q \subset \cS$ be a totally
ordered subset.  If $X$ is an $\langle \cS \rangle$-manifold, then
$\partial^Q X$ is an $\langle \partial^Q \cS \rangle$-manifold.
\end{lem}

\begin{proof}
If $Q \subset P$, then a neighbourhood of $\partial^P X$ in
$\partial^Q X$ is locally homeomorphic to $\partial^P X \times
[0,\infty)^{P \setminus Q}$. 
\end{proof}

If $X$ admits an action of $G$ that preserves the stratification, then
$\partial^Q X$ inherits an action.
\begin{defin}
For each totally ordered subset $Q \subseteq \cS$, there is a functor 
\begin{equation}
\partial^Q \co  \Kur\langle \cS \rangle \to \Kur\langle \partial^Q \cS
\rangle 
\end{equation}
specified by mapping $\bX = (X, V, s, G)  $  to
\begin{equation}
\partial^Q \bX \equiv (\partial^Q X, V,  s|\partial^Q X, G).
\end{equation}   
\end{defin}

The monoidal structure on the category of Kuranishi charts
corresponds, in the stratified setting, to the existence of natural
functors
\begin{equation}
\Kur\langle \cS_1 \rangle \times \Kur\langle \cS_2 \rangle \to
\Kur\langle \cS_1 \amalg \cS_2 \rangle
\end{equation}
 whenever $\cS_1 \amalg \cS_2$ is ordered in such a way that $p_1 <
 p_2$ whenever $p_i \in \cS_i$. This makes the disjoint union of the
 categories $\Kur \langle \cS \rangle$ into a monoidal category. In
 Section \ref{sec:strat-kuran-charts} below, we shall consider a
 category of \emph{stratified charts} with objects the disjoint union
 of the categories $\Kur \langle \cS \rangle$, but with some
 additional morphisms corresponding to inclusions of boundary strata.

\begin{rem}
  We note that a stratification was not assumed in Section \ref{sec:topol-enrich-categ}, when topological and orbispace flow categories were considered. In that context, the stratification on each morphism space is induced by the structure maps. Thus, a Kuranishi presentation of a morphism space in a flow category inherits a stratification of its zero locus. We shall require a stratification of the thickening as well, which  justifies the introduction of the above notion.
 \end{rem}

\subsubsection{Stratified Kuranishi charts}
\label{sec:strat-kuran-charts}

We now introduce the main category of Kuranishi charts which we use: 

\begin{defin} \label{def:stratified-Kuranishi-charts}
The category of \emph{stratified Kuranishi charts}, denoted $\Kur$,
has
\begin{enumerate}
\item objects consisting of a partially ordered set $\cS$ and an
  object $\bX$ of $\Kur\langle \cS \rangle $, and
\item morphisms $f$ from  $ (\cS, \bX)$ to $(\cS',\bX')$ specified by an
order-preserving isomorphism $\rho \colon \cS \to \partial^Q \cS'$ for
a totally ordered subset $Q$ of $\cS'$ and a map $\bX \to \partial^{Q}
\bX'$ of $\Kur\langle \cS \rangle$ charts, where we regard
$\partial^{Q} \bX'$ as a $\Kur\langle \cS \rangle$ chart via $\rho$.
\end{enumerate}

Composition of morphisms is defined as in
Definition~\ref{defin:stratorbicharts}.
\end{defin}

By construction, we have a fully faithful embedding $\Kur \langle \cS
\rangle \subset \Kur$ which lands in the subcategory of morphisms
where $S = S'$ (i.e. $\rho$ is the identity map).  We have a functor
\begin{equation}
    \Kur \times \Kur \to \Kur
\end{equation}
which assigns to each pair $\bX_{1} \in \Kur \langle \cS_{1} \rangle$
and $\bX_{2} \in \Kur \langle \cS_{2} \rangle$ the product Kuranishi
chart $\bX_{1} \times \bX_{2} \in  \Kur\langle \cS_1 \amalg \cS_2
\rangle$. On morphisms, we define the map 
\begin{equation}
f_1 \times f_2 \colon \bX_1 \times \bX_2 \to \partial^{Q_1 \amalg Q_2}
(\bX'_1 \times \bX'_2)
\end{equation}
using the isomorphism $\partial^{Q_1} \bX'_1 \times \partial^{Q_2}
\bX'_2 \cong \partial^{Q_1 \amalg Q_2} (\bX'_1 \times \bX'_2)$.

The functoriality and associativity of the above product are recorded
in the following lemma. 

\begin{lem}
The product of Kuranishi charts equips $\Kur$ with a monoidal
structure, with unit $(\emptyset, (\ast, \{0\}, s, \{e\}))$. The
zero-locus functor $\Kur \to \Orb$ is strongly monoidal. \qed  
\end{lem}

We again often want to restrict attention to the subcategory of $\Kur$
in which the morphisms of partially ordered sets are isomorphisms.

\begin{defin}\label{defn:kuriso}
The category $\Kur^{\iso}$ is the subcategory of $\Kur$ with the same
objects but where the morphisms are specified by identities $\cS =
\partial^{Q} \cS'$.  For a given $\cS$, the category $\Kur^{\cS}$ is the
full subcategory of $\Kur^{\iso}$ spanned by the objects $(\partial^Q
\cS, \bX)$ as $Q$ varies over the totally-ordered subsets of $\cS$.
\end{defin}

Notice that given an isomorphism $\cS \cong \cS'$, there is an induced
isomorphism of categories $\Kur \langle \cS \rangle \to \Kur \langle
\cS' \rangle$. In the context of Section \ref{sec:topol-orfif-flow}, we are given a partially ordered set $\cP$, equipped with an action of a group $\Pi$, and this induced isomorphism yields:
\begin{lem}\label{lem:kurequi}
For $\pi \in \Pi$ and a subset $\cQ \subset \cP$, there is a natural
isomorphism of categories 
\begin{equation}
\alpha_{\pi} \colon \Kur \langle \cQ \rangle \to \Kur \langle \pi \cQ
\rangle.
\end{equation}
This assignment is strictly associative since $\alpha_{\pi} \circ
\alpha_{\pi'} = \alpha_{\pi \pi'}$. \qed
\end{lem}

Given an object $X$ of $\Kur^{\cS}$, we can produce an object of $\Kur
\langle \cS \rangle$ by collaring, i.e., attaching a cube $\kappa^Q
\times \partial^Q X$ to $X$ for each totally-ordered subset $Q$,
exactly as in Definition~\ref{defn:collared-completion}.

\begin{lem}\label{lem:collar-functor}
There is a collar functor 
\begin{equation}
(\widehat{-}) \colon \Kur^{\cS} \to \Kur \langle \cS \rangle.
\end{equation} \qed
\end{lem}

We will use the collaring functor extensively in
Section~\ref{sec:impl-atlas-with}.  In particular, we will rely on the
fact that the collaring functor is compatible with products in the
following sense.

\begin{prop}
Let $\cS_1$, $\cS_2$, and $\cS_{12} $ be partially ordered sets, and
$Q$ a totally ordered subset of $\cS_{12}$. An isomorphism $\cS_1
\amalg \cS_2 \cong \partial^Q \cS_{12}$ determines a commutative
diagram: 
\begin{equation}
\begin{tikzcd}
\Kur^{\cS_1 \amalg \cS_2} \ar[d,swap,"(\widehat{-})"] & \Kur^{\cS_1} \times \Kur^{\cS_2}
\ar[l] \ar[r] & \Kur^{\cS_{12}} \ar[d, "(\widehat{-})"] \\
\Kur \langle \cS_1 \amalg \cS_2 \rangle &\Kur \langle \partial^Q \cS_{12}
\rangle  \ar[l,swap,"\cong"] & \Kur \langle \cS_{12}
\rangle \ar[l,swap, "\partial^Q"].
\end{tikzcd}
\end{equation} \qed
\end{prop}

\subsubsection{Kuranishi presentations}
\label{sec:anoth-categ-charts-pres}
Let $\cM$ be a compact Hausdorff space that is $\langle \cS
\rangle$-stratified, for $\cS$ a partially ordered set.

\begin{defin}
A $d$-dimensional \emph{$\langle \cS \rangle$-Kuranishi presentation}
of $\cM$ consists of the following data: 
\begin{enumerate}
\item An indexing category $A$.
\item A stratified orbispace presentation $(\partial^{Q_\bullet}\cS_\bullet, Z_\bullet, G_\bullet) \co A \to \Orb^{\cS}$ of $\cM$.
\item A lift $\bX \co A \to  \Kur^{\cS}  \to \Orb^{\cS}$ of the
  orbispace presentation via the zero locus functor that factors
  through the subcategory of $\Kur^{\cS}$ consisting of
  charts of virtual dimension $d - |Q_\bullet|$.
\end{enumerate}
\end{defin}
Note that the value of the partially ordered set $\cS_\bullet$ is in fact constant, and given by $\cS$, because we imposed this condition in the definition of an orbispace presentation, and the footprint functor preserves the datum of the partially ordered set.

Whenever $\cS = \emptyset$, we say that this is a \emph{closed
  Kuranishi presentation} of $\cM$. Whenever $A$ is a singleton, we
say that $\cM$ is equipped with a \emph{global Kuranishi chart}. We
shall often abuse notation and write $\bX \co A \to \Kur$ for a Kuranishi presentation, neglecting to write down the
homeomorphism to $\cM$, and even the stratification data; 
indeed, the space $\cM$ and its $\langle \cS \rangle$-stratification
are determined by the functor $\bX$, as is it homeomorphic to the
colimit of the footprints $Z_\alpha/G_\alpha$ over all charts.

\subsubsection{Kuranishi flow categories}
\label{sec:kuran-orbisp-flow}
We now return to the setting of Sections \ref{sec:topol-enrich-categ}
and \ref{sec:orbisp-flow}: $\cP$ is a partially ordered set with a
free action of $\Pi$. Additionally, we fix an assignment of integers
$d_p \in \bZ $ for each element $p \in \cP$, which is $\Pi$-equivariant
in the sense that 
\begin{equation}
  d_{\pi \cdot p} \equiv  d_p  + 2 \deg \pi.  
\end{equation}
\begin{rem}
The reader familiar with Morse theory may want to keep in mind that,
in this context, the integer $d_p$ assigned to a critical point is the
dimension of the positive-definite subspace of the Hessian matrix.
\end{rem}

We write $\cP(p,q)$ for the partially ordered set of elements of $\cP$
which are strictly between $p$ and $q$. Given a triple $p < q < r$, we
have an equality 
\begin{equation} \label{eq:orbits-inclusion-product-stratum}
  \cP(p,q) \amalg \cP(q,r) = \partial^{q} \cP(p,r),
\end{equation}
where $\partial^{q} \cP(p,r)$ denotes as above the complement of
$\{q\}$ in the subset of elements of $\cP(p,r)$ that are comparable to $q$.
This observation plays a key role in the following, which should be
compared with Definition \ref{def:orbispace_flow_category}: 
\begin{defin} \label{def:Kuranishi-flow-category}
  A $\Pi$-equivariant \emph{Kuranishi flow category} $\bX$ with objects $\cP$ consists of:
  \begin{enumerate}
   \item A strictly $\Pi$-equivariant $2$-category $A$, with object set
     $\cP$, and such that $A(p,p) = \ast$, and $A(p,q)$ is empty
     unless $p \leq q$.
 \item A strict $2$-functor $\bX \co A \to \Kur$, which assigns to
   $A(p,q)$ a $ \langle \cP(p,q) \rangle $-Kuranishi presentation of
   dimension $d_p - d_q-1$ (except if $p=q$, in which case we require
   this presentation to have dimension $0$, and to be a point), and
   which is strictly $\Pi$-equivariant.  (Here we are regarding $\Kur$
   as a bicategory with a single object.)
\end{enumerate} 

\end{defin}
Note the analogy with Definition \ref{def:orbispace_flow_category},
with the proviso that we now require $\bX$ to be a strict rather than
a lax functor.  In particular, the $2$-functor $A \to
\Kur$ consists of a functor $A(p,q) \to \Kur$ for each pair $(p,q)$,
so that the following diagram commutes 
\begin{equation}\label{eq:Kuranishi-flow-category-isomorphism_triple}
   \begin{tikzcd}
      A(p,q) \times A(q,r) \ar[r] \ar[d] &  A(p,r) \ar[d]  \\
      \Kur  \times \Kur \ar[r] 
      &   \Kur .
    \end{tikzcd}
  \end{equation}
The equivariance of the functor in Definition
\ref{def:Kuranishi-flow-category} comes from the trivial $\Pi$-action
on $\Kur$, and the natural isomorphisms of Lemma~\ref{lem:kurequi}
corresponding to relabelling a $ \langle \cP(p,q) \rangle $-Kuranishi
presentation to obtain a $ \langle \cP(\pi \cdot p, \pi \cdot q)
\rangle $ presentation.

For later constructions, we shall also fix a $\Pi$-invariant
assignment $V_{p} \equiv (V^+_{p}, V^-_{p})$ of a pair of vector
spaces for each element $p \in \cP$, so that $V^-_{p}$ is a complex
vector space, with the property that each orbit has an element $p_0$
such that
\begin{equation} \label{eq:degree_basepoint_lifted_orbit}
\dim V^+_{p_0} - \dim V^-_{p_0} =  d_{p_0}.
\end{equation}

\begin{rem} \label{rem:choose-slice-stable-index}
In Morse theory, $V^+_{p}$ is the positive eigenspace of the Hessian
at the critical point associated to $p$ and we shall set $V^-_{p}$ to
vanish,  while in Floer theory the stable vector space $V_p$
represents the virtual index of an operator on a plane with asymptotic
ends associated to a Hamiltonian orbits.
\end{rem}

\subsection{Flag smooth Kuranishi presentations}
\label{sec:fiber-kuran-pres}

The purpose of this subsection and the next is to introduce the
refinements of the notion of Kuranishi presentation that are required
to construct virtual fundamental chains in Floer theory, with
coefficients in complex oriented cohomology theories. The reader who
is mostly interested in the formal aspects of the theory may want to
postpone reading this section until after Section
\ref{sec:smoothly-fiber-kuran-1}, while the reader who is mostly
interested in Hamiltonian Floer theory should probably first read
Section \ref{sec:hamilt-floer-theory-1} and refer to this section
while reading Section \ref{sec:kuran-flow-categ-1}. 

\subsubsection{Fibered Kuranishi charts}
\label{sec:fiber-kuran-pres-2}

We begin by considering the structure we are facing before any choice
of smooth structure.  Since all constructions require stratifications,
we shall first fix a partially ordered set $\cS $. If $X$ and $B$ are
$\langle \cS \rangle$-manifolds, a stratified submersion
\begin{equation}
  X \to B
\end{equation}
is a stratified map which is locally homeomorphic to a projection   
\begin{equation}
  \bR^{n_0 + n_1} \times [0,\infty)^{Q} \to \bR^{n_1} \times [0,\infty)^{Q},  
\end{equation}
in a neighbourhood of each point lying in $\partial^Q X$. 

\begin{defin}
A \emph{fibered $\langle  \cS \rangle$-Kuranishi chart} is a  triple
$(\bX, B,  \pi)$, consisting of an $\langle \cS \rangle$-Kuranishi
chart, a $G$-equivariant $\langle \cS \rangle$-manifold $B$, and a
$G$-equivariant stratified topological submersion $\pi \co X \to B$. 
\end{defin}

We form a category of fibered stratified Kuranishi charts as follows.

\begin{defin}
The category of fibered $\langle \cS \rangle$-Kuranishi charts has
objects as above and a morphism $f \co \alpha \to \beta$ is given by a
map of the corresponding Kuranishi charts, and an equivariant
topological stratified submersion $B_\alpha \to B_\beta$, which fit in
a commutative diagram 
\begin{equation} \label{eq:commutative_diagram_fibered_map}
  \begin{tikzcd}
    X_\alpha \ar[r]  \ar[d] & X_\beta \ar[d] \\
    B_\alpha  \ar[r] & B_\beta,
  \end{tikzcd}
\end{equation}
such that the map $ X_\beta \to  B_\beta \times V_\beta/V_\alpha$ is a
topological submersion near the image of $X_\alpha$, and contains this
image as an open subset. It is straightforward to define composition
as in the construction of the category of Kuranishi charts.
\end{defin}

\begin{rem}
The condition on morphisms should be thought of as a strengthening of
the transversality condition for Kuranishi charts: writing $X_{\beta,
  p}$ for the fibre of $X_\beta$ over a point $p$ in $B_\beta$, we
obtain a Kuranishi chart $\beta_p$ given by this total space, an
action of $G_\beta$ inherited from $X_\beta$, and the projection to
$V_\beta$. Writing $\alpha_p$ for the corresponding Kuranishi chart
obtained as the fibre of the composition $X_\alpha \to B_\alpha \to
B_\beta$, the transversality in Diagram
\eqref{eq:commutative_diagram_fibered_map} is thus the assertion that
we have a map $\alpha_p \to \beta_p$ of Kuranishi charts for each $p
\in B_\beta$. 
\end{rem}

\subsubsection{Flag smooth structures}
\label{sec:fibr-smooth-struct}

In the language of the previous section, a flag smooth Kuranishi chart
is a fibred Kuranishi space together with a choice of smooth structure
on the base $B$ and a fibrewise smooth structure on the topological
submersion $X \to B$, in the following sense: we have a choice of
atlas for $X$ consisting of product charts $\bR^{n} \times U \to X$
over charts $U \to B$, with transition functions which are
continuously differentiable in the fibre direction (depending
continuously on the base). Before implementing this idea, we explain
why a naive approach fails: 

\begin{rem}
  Assuming that the projections $X_\alpha \to B_\alpha$ and $X_\beta
  \to B_\beta$ in Diagram \eqref{eq:commutative_diagram_fibered_map}
  are both equipped with fibrewise smooth structures, one natural
  condition to impose is that these two smooth structures be
  compatible in the sense that the fibres of $X_\alpha \to B_\beta$
  are also equipped with a smooth structure such that the composition
  with the map to $B_\alpha$ is smooth. We could try to define a
  morphism to be such a choice of smooth structure, but it is not
  clear how to define compositions in this context.  We believe that
  it might be possible to resolve this problem by working in an
  $\infty$-categorical context (i.e., by introducing a space of
  compositions), but we did not explore such a solution.
\end{rem}

Our way to handle this issue is to record an additional projection to
smooth manifolds with maps going in the other direction:  
\begin{defin} \label{def:equivariant-submersion-category} 
The category of \emph{equivariant submersions} of smooth $\langle \cS
\rangle$-manifolds has objects $(G, B \to B')$ consisting of a finite
group $G$, and a $G$-equivariant submersion $B \to B'$ of $\langle \cS
\rangle$-manifolds. A morphism $(G_0, B_0 \to B'_0) \to (G_1, B_1 \to
B'_1)$ consists of a homomorphism $G_0 \to G_1$, and a commutative
diagram 
\begin{equation}
\begin{tikzcd}
B_0 \ar[r] \ar[d] & B_1 \ar[d] \\
B'_0 & \ar[l] B'_1
\end{tikzcd}
\end{equation}
of smooth $G_0$-equivariant maps.
\end{defin}
Note that it follows immediately that the arrow $B_0 \to B_1$ is a
submersion, but that $B'_1 \to B'_0$ need only be submersive along the
image of the composition $B_0 \to B_1 \to B'_1$. 

\begin{defin} \label{def:smoothly_fibred_complete}
  A \emph{flag smooth $\langle  \cS \rangle$-Kuranishi chart} $\alpha$ consists of:
  \begin{enumerate}
  \item an $\langle \cS  \rangle$-Kuranishi chart $\bX_\alpha$,
  \item a smooth $G_\alpha$-equivariant submersion of $\langle \cS  \rangle $-manifolds $B_\alpha \to B'_\alpha$.
  \item a $G_\alpha$-equivariant stratified topological submersion  $X_\alpha \to B_\alpha$, and
    \item a  $G_\alpha$-invariant fiberwise smooth structure on the composition $X_\alpha \to B'_{\alpha}$. 
  \end{enumerate}
  We require these data to be compatible in the following sense:
  \begin{enumerate} 
  \item The restriction of the map $s_\alpha \co X_\alpha \to V_\alpha$ to a fiber of $X_\alpha \to B'_{\alpha}$ is smooth near $Z_\alpha$, and
    \item the restriction of the map $X_\alpha \to B_\alpha$ to a fiber of $X_\alpha \to B'_{\alpha}$ is a smooth submersion near $Z_\alpha$ onto the corresponding fiber of the map $B_\alpha \to B'_{\alpha}$.
  \end{enumerate}
\end{defin}

We define a category of flag smooth charts as follows.

\begin{defin}
The category $\Kur^{fs}\langle  \cS \rangle$ of {\em flag smooth
  $\langle \cS \rangle$-Kuranishi charts} has objects charts as in
Definition \ref{def:smoothly_fibred_complete}.  A morphism from
$\alpha$ to $\beta$ consists of a map $f$ of Kuranishi charts as
before and a commutative diagram of stratified maps  
\begin{equation}
  \begin{tikzcd}
    X_\alpha \ar[r] \ar[d,"f"] & B_\alpha \ar[r, two heads] \ar[d, two heads] & B'_{\alpha} \\
     X_\beta \ar[r] & B_\beta \ar[r, two heads]  & B'_{\beta} \ar[u]
  \end{tikzcd}
\end{equation}
in which each arrow labelled $\twoheadrightarrow$ is an equivariant smooth submersion, and such that, for each $p \in B_\beta $, the map
\begin{equation}
  X_{\beta,p} \to V_\beta/V_\alpha
\end{equation}
is a smooth submersion near the image of $X_{\alpha,p}$. We shall
impose an additional condition as follows.  The fibres of $X_\alpha
\to B'_{\beta}$ acquire smooth structures in two different ways:
\begin{enumerate}
\item from their inclusions in the fibres of $X_\alpha \to
B'_{\alpha}$ and the fact that $X_\alpha \to B_\alpha$ is smooth on
fibres, and 
\item from the inclusion of their free quotients in the
fibres of $X_\beta \to B'_\beta$ and the fact that these fibres
submerse to $V_\beta/ V_\alpha$. 
\end{enumerate}
We requires that these two smooth structures agree.
\end{defin}

We define composition in terms of the composition of morphisms of
Kuranishi charts, as well as the two right squares in the following
diagram: 
\begin{equation}
  \begin{tikzcd}
    X_\alpha \ar[r] \ar[d] & B_\alpha \ar[r, two heads] \ar[d, two heads] & B'_{\alpha} \\
    X_\beta \ar[r]   \ar[d]& B_\beta \ar[r, two heads]   \ar[d, two heads] &  B'_{\beta} \ar[u] \\
    X_\gamma \ar[r] & B_\gamma \ar[r, two heads]  & B'_{\gamma} \ar[u] .
  \end{tikzcd}
\end{equation}
\begin{lem}
  \label{lem:composition-flag-smooth}
  The morphisms in $\Kur^{fs}\langle  \cS \rangle$ are closed under composition.
\end{lem}
\begin{proof}
  The key point to check is that the two induced smooth structures on
  the fibres of the map $X_\alpha \to B'_\gamma$ agree. Denoting by
  $X_{\alpha, p'} $ such a fibre for $p' \in B'_\gamma$,  this follows
  by considering the following diagram 
\begin{equation}
  \begin{tikzcd}
X_{\alpha, \pi_{\alpha}(p')} & X_{\alpha, p'} \ar[r]  \ar[dr] \ar[l] \ar[ld]  &  X_{\gamma,p'} \\
X_{\alpha, \pi_{\beta}(p')} \ar[u]   \ar[r] & X_{\beta, \pi_{\beta}(p')}  & X_{\beta,p'} \ar[l] \ar[u].
  \end{tikzcd}
\end{equation}
\end{proof}

\subsubsection{Flag smooth Kuranishi flow categories}
\label{sec:fibr-smooth-kuran}

To proceed further, we combine the categories associated to different
choices of partially-ordered sets labelling strata into a single
category, taking Definition~\ref{def:stratified-Kuranishi-charts} as
our model: 

\begin{defin} \label{def:equivariant_submersion_monoidal_category}
The category of \emph{equivariant submersions of stratified manifolds}
is the category with
\begin{enumerate}
\item Objects given by $(\cS,  G, B \to B')$, with
$\cS$ a partially ordered set, and $(G, B \to B') $ an object of the
category of smooth equivariant submersions of $\langle \cS
\rangle$-manifolds. 

\item Morphisms $f \colon (\cS_0, G_0, B_0 \to B'_0) \to (\cS_1,  G_1, B_1
\to B'_1)$ specified by an order-preserving 
isomorphism $\rho \colon \cS_0 \cong \partial^Q \cS_1$
for some totally ordered subset $Q$ of $\cS_1$ and a map
  \begin{equation}
  (G_0, B_0 \to B'_0) \to (G_1, \partial^{Q} B_1 \to \partial^{Q} B'_1)  
  \end{equation}
  of equivariant submersions of $\langle \cS_0 \rangle$-manifolds, 
where we regard $\partial^{Q} B_1$ and $\partial^Q B'_1$ as $\langle
\cS_0 \rangle$-manifolds via $\rho$.
\end{enumerate}

Composition is defined as in
Definition~\ref{def:stratified-Kuranishi-charts}.  This category is
monoidal with the evident structure maps.
\end{defin}

We may proceed in exactly the same way to pass to Kuranishi charts:

\begin{defin}
The category $\Kur^{fs} $ of \emph{stratified flag smooth Kuranishi
  charts} has objects $( \cS, \alpha)$ consisting of a partially
ordered set $\cS$ and an object $\alpha$ of $\Kur^{fs}\langle \cS \rangle$.  A morphism $f \co ( \cS_0,\alpha_0) \to ( \cS_1,  \alpha_1)$  consists of a morphism in $\Kur$, and a lift of the corresponding morphism in $ \Kur \langle \cS_0 \rangle $ to  $\Kur^{fs}\langle   \cS_0 \rangle$. 
\end{defin}

The following lemma helps make sense of the above definition.

\begin{lem}
The compatibility of the flag smooth structure
with the stratification yields a lift of the restriction functor 
\begin{equation}
\Kur \langle \cS \rangle \to  \Kur \langle \partial^{Q} \cS \rangle
\end{equation}
associated to a totally ordered subset $Q$ of $\cS$ to a functor
\begin{equation}
 \Kur^{fs}\langle  \cS \rangle \to \Kur^{fs}\langle  \partial^{Q} \cS \rangle.
\end{equation}
\end{lem}

In particular and as before, there are no morphisms from $ (\cS_0,
\alpha_0)$ to $(\cS_1, \alpha_1)$ unless $ \cS_0 \cong \partial^Q
\cS_1$ for some $Q$.  We will again sometimes want to restrict to
subcategories where this isomorphism is in fact the identity:
\begin{defin}\label{defn:kurisofs}
The subcategory $\Kur^{fs,\iso}$ consists of those morphisms for which
the isomorphism $\cS_0 \cong \partial^Q \cS_1$ is the identity.  The
full subcategory $\Kur^{fs,\cS} \subset \Kur^{fs,\iso}$ is spanned by
the objects of the form $\partial^Q \cS$ as $Q$ varies over the
totally-ordered subsets of $\cS$.
\end{defin}

There is a natural product functor
\begin{equation} \label{eq:product-fibrewise-smooth-charts}
    \Kur^{fs} \langle  \cS_1 \rangle \times   \Kur^{fs} \langle  \cS_2 \rangle   \to \Kur^{fs} \langle \cS_1 \times \cS_2 \rangle
  \end{equation}
  lifting the corresponding functor for Kuranishi spaces. Explicitly,
  given flag smooth charts $\alpha_1$ and $\alpha_2$, we define the
  product $\alpha_1 \times \alpha_2$ to have underlying Kuranishi
  chart given by the product $\bX_1 \times \bX_2$ and use the
  submersion of smooth $\langle \cS_1 \amalg \cS_2\rangle$-manifolds 
  \begin{equation}
    B_{\alpha_1 \times \alpha_2} \equiv B_{\alpha_1} \times B_{\alpha_2} \to B'_{\alpha_1} \times B'_{\alpha_2} \equiv B'_{\alpha_1 \times \alpha_2}.
  \end{equation}
  There is an induced projection
  \begin{equation}
   X_{\alpha_1 \times \alpha_2} \to B'_{\alpha_1 \times \alpha_2},
 \end{equation}
and the product induces a fibrewise smooth structure over
$B'_{\alpha_1 \times \alpha_2}$. Verifying the remaining properties is
straightforward, as is the functoriality of this construction. 

The above product functor is naturally coherently associative, and is
compatible with maps in $\Kur^{fs}$, so we conclude: 
\begin{lem}
The product of flag smooth charts defines a monoidal structure on
$\Kur^{fs}$, for which the forgetful functor to $\Kur$ is strictly
monoidal. \qed 
\end{lem}

We now lift the notion of flow categories to the flag smooth setting. 

\begin{defin} \label{def:flag-smooth-presentation}
  A $\Pi$-equivariant \emph{flag smooth Kuranishi flow category} with
  object set $ \cP$ is a Kuranishi flow category equipped with
  compatible lifts of all diagrams from $\Kur$ to $\Kur^{fs}$.  
\end{defin}
Note that unpacking this definition involves the categories
$\Kur^{fs,\cS}$ from Definition \ref{defn:kurisofs}. 

\subsection{Complex-oriented Kuranishi charts}
\label{sec:compl-orient-kuran-1}

We continue the development of an abstract framework for formulating the notion of (stable) complex orientations of Kuranishi charts and presentations. We begin by considering a refinement of the notion of a flag smooth chart, which is equipped with a canonical tangent space, and then discuss the data required to fix a stable isomorphism between this tangent space and a complex vector bundle. We note that, throughout this section, we only explicitly use the projection map $X_\alpha \to B_\alpha$, and the map $B_\alpha \to B'_\alpha$ will play no role in our construction.

\subsubsection{Charts with tangent bundles}
\label{sec:splitt-exact-sequ}

 We may assign to each flag smooth chart $\alpha$ a vector bundle on $X_\alpha$ given by the direct sum of $TB_\alpha$ and the fibrewise tangent bundle  $   T^\alpha X_\alpha $ (see Section \ref{sec:smoothly-fiber-kuran-1} below for an extended discussion). Our goal is to compare these tangent bundles for the source and target of maps of flag smooth charts. The key result of this section is Lemma \ref{lem:tangent-space-split-morphisms-well-defined}.

We start by associating to a partially ordered set $\cS$ a category $\Kur^{\scrT}\langle  \cS
\rangle$  of \emph{$\langle \cS \rangle$-Kuranishi charts with tangent bundles:}  this is an internal category in spaces, i.e., a category which is equipped with a topological space of objects and
a topological space of morphisms (see Appendix~\ref{sec:internal} for
a quick summary of the definition and properties of internal categories):
\begin{defin}
 An object of $\Kur^{\scrT}\langle  \cS
\rangle$ consists of 
\begin{enumerate}
\item an object $\alpha$ of $\Kur^{fs}\langle   \cS
  \rangle$,
\item an open $G_\alpha$-invariant subspace $Z^\ori_\alpha$ of $Z_\alpha$, 
\item inner products on $T^\alpha X_\alpha$ and $T B_\alpha$
  (as vector bundles on $Z_\alpha$ and $B_\alpha$), which are
  $G_\alpha$-invariant and are compatible with the stratification in
  the sense that the restriction of $T B_\alpha$ to the stratum
  labelled by $Q \subset \cS$ splits as an othogonal direct sum of $ T
  \partial^{Q} B_\alpha$ with the direction normal to each codimension
  $1$ boundary stratum $\partial^{q} B_\alpha$ for $q \in Q$: 
\begin{equation}
TB_{\alpha} \cong \bR^{Q} \oplus T \partial^Q B_\alpha.
\end{equation}
\end{enumerate}
The topology on the space of inner products on a fixed vector space
determines the topology on the space of objects. 

A morphism in $\Kur^{\scrT}\langle  \cS \rangle  $ from a lift of a
flag smooth Kuranishi chart $\alpha$ to a lift of $\beta$ consists of
a morphism $f$ in $  \Kur^{fs}\langle   \cS \rangle$ from $\alpha$ to
$\beta$ mapping $Z^\ori_\alpha$ to $Z^\ori_\beta$, such that  the
short exact sequences 
\begin{align}
0 \to  T^\beta B_\alpha \to   & TB_\alpha   \to  TB_\beta \to 0 \\
0 \to   T^\alpha X_\alpha  \to &  T^\beta X_\beta \to V_f^\perp \oplus    T^\beta B_\alpha \to 0
\end{align}
induce the same inner product on $V_f^\perp \oplus    T^\beta
B_\alpha$ (as a vector bundle on $Z^\ori_\alpha$). The topology on the
space of morphisms is given by pulling back the subspace topology on
the product of the spaces of objects associated to the source and the
target. 
\end{defin}

To be more explicit about the space of morphisms, we require that
$TB_\alpha \to  TB_\beta $ be an orthogonal projection (with kernel
that we denote $  T^\beta B_\alpha $), and that $T^\alpha X_\alpha
\to  T^\beta X_\beta $ be an isometric embedding. Recall that the map
$V_\alpha \to V_\beta$ is an isometric embedding, which equips
$V_f^\perp  $ with an inner product. The definition of a morphism of
flag smooth charts implies that the cokernel of $ T^\alpha X_\alpha
\to  T^\beta X_\alpha$ is identified with the direct sum $V_f^\perp
\oplus    T^\beta B_\alpha $, and we require that the inner product
induced from its description as a quotient splits as a direct sum of
the inner product on $V_f^\perp $ with the restriction of the inner
product on $TB_\alpha $ to $T^\beta B_\alpha$.

Note that the space of morphisms can be alternatively described as a subspace of the space $G_\alpha$-equivariant inner products on $   T^\beta X_\beta $, $T B_\alpha$, and $T B_\beta $ (the last datum is required because the map $B_\alpha \to B_\beta$ is not in general surjective). Equipping it with the induced topology, the natural map from the space of morphisms lifting $f$ to the spaces of objects lifting $\alpha$ and $\beta$ is continuous; i.e., the source and target maps are continuous. The next result asserts that composition is continuous and well-defined:
  \begin{lem}
    Given lifts of composable arrows $f \co \alpha \to \beta$ and $g \co \beta \to \gamma$ to $  \Kur^{\scrT} \langle   \cS \rangle $ whose restrictions to the space of lifts of $\beta$ agree, the inner products on $   T^\gamma X_\gamma $ and $T B_\alpha$ induce the same inner product on $V_{f \circ g}^\perp \oplus    T^\gamma B_\alpha $.
\end{lem}
    \begin{proof}
      Since the two inner products on $TB_\beta$ associated to $f$ and $g$ agree, the restriction to the subspace $T^\gamma B_\alpha$ of the inner product on $TB_\alpha$ associated to $f$ splits as a direct sum $ T^\beta B_\alpha \oplus  T^\gamma B_\beta$, with the two summands equipped with the inner products associated to the lifts of morphisms $f$ and $g$. Similarly, the inner product on $T^\gamma X_\gamma $ induces an inner product on its quotient by $T^\alpha X_\alpha $ which splits as an orthogonal direct sum of $  V_f^\perp \oplus    T^\beta B_\alpha$ and $ V_g^\perp \oplus    T^\gamma B_\beta$. The result follows by collecting factors.
    \end{proof}

    For the next result, we write $  T X_\alpha $ for the direct sum $T^\alpha X_\alpha \oplus TB_\alpha$.
  
\begin{lem} \label{lem:tangent-space-split-morphisms-well-defined}
  Each morphism  in $\Kur^{\scrT}\langle   \cS \rangle $ lifting an arrow $f$ induces an isomorphism
  \begin{equation}
    T X_\alpha \oplus  V_f^\perp \cong T X_\beta 
  \end{equation}
  of equivariant vector bundles over $Z^\ori_\alpha$.    A lift of composable arrows $f$ and $g$ induces a commutative diagram
   \begin{equation}
     \begin{tikzcd}
       TX_\alpha \oplus V_{f}^\perp \oplus V_g^\perp \ar[r] \ar[d] &   T X_\beta \oplus V_g^\perp\ar[d] \\
  TX_\alpha \oplus V_{g \circ f}^\perp \ar[r] &   TX_\gamma  .  
     \end{tikzcd}
   \end{equation}
\end{lem}
\begin{proof}
Using the orthogonal decomposition associated to an inner product, the datum of a morphism in $\Kur^{\scrT}\langle   \cS \rangle $ induces isomorphisms
  \begin{align}
    T^\beta X_\beta  & \cong   T^\alpha X_\alpha   \oplus  T^\beta B_\alpha  \oplus V_f^\perp  \\
     T B_\alpha & \cong  T B_\beta  \oplus T^\beta B_\alpha .
  \end{align}
   Combining these two isomorphisms, we obtain
   \begin{align}
     T X_\beta & \cong  T^\beta X_\beta \oplus T B_\beta \\
     & \cong  T^\alpha X_\alpha   \oplus  T^\beta B_\alpha  \oplus T B_\beta   \oplus V_f^\perp \\
     & \cong T^\alpha X_\alpha \oplus T B_\alpha \oplus V_f^\perp \equiv T X_\alpha \oplus V_f^\perp.
   \end{align}
   The compatibility with compositions is straightforward.
 \end{proof}

In the remainder of this section, we discuss the multiplicativity of
these constructions.  As before, the starting point is the
construction of a category $\Kur^{\scrT}$ with objects consisting
of the union of the categories $\Kur^{\scrT}\langle \cS \rangle$
as $\cS$ varies.

\begin{defin}
The \emph{category of Kuranishi spaces with tangent spaces} is the
internal category $\Kur^{\scrT}$ in topological spaces with 
\begin{enumerate}
\item space of objects given by the union of the object spaces of
  $\Kur^{\scrT}\langle \cS \rangle$ as $\cS$ varies over all
  partially-ordered sets, equipped with the topology as a disjoint union, and

\item space of morphisms specified by the stipulation that a
  morphism between lifts of objects in $\Kur^{fs}\langle \cS_0
  \rangle$ and $\Kur^{fs}\langle \cS_1 \rangle$ consists of a morphism
  $f$ in $\Kur$ and a lift of the induced arrow in $ \Kur^{fs} \langle
  \cS_0 \rangle $ to $\Kur^{\scrT}\langle \cS_0 \rangle$.
\end{enumerate}
\end{defin}

The key point, as in Section~\ref{sec:fibr-smooth-kuran}, is the
restriction functor 
\begin{equation} \label{eq:restriction_flag_smooth_inner_product}
 \Kur^{\scrT}\langle  \cS \rangle \to \Kur^{\scrT}\langle  \partial^Q \cS \rangle
\end{equation}
induced by the compatibility of inner products with stratifications.
Note that restriction induces canonical isomorphisms 
\begin{align}
  T^{\alpha} X & \cong T^{\alpha} \partial^Q X \\
  T B_{\alpha} & \cong \bR^{Q} \oplus T \partial^Q B_\alpha,
\end{align}
where the second relies on the choice of inner product. Taking the
direct sum, we conclude:
\begin{lem} \label{lem:tangent_space_chart_splits_at_boundary}
The restriction functor in Equation
\eqref{eq:restriction_flag_smooth_inner_product} induces an
isomorphism 
\begin{equation}
T X_{\alpha}  \cong \bR^{Q} \oplus T \partial^Q X_\alpha.
\end{equation} \qed
\end{lem}

We have a natural product map
\begin{equation}
\Kur^{\scrT} \langle  \cS_1 \rangle \times   \Kur^{\scrT} \langle
\cS_2 \rangle   \to \Kur^{\scrT} \langle \cS_1 \times \cS_2 \rangle, 
\end{equation}
lifting the product in Equation
\eqref{eq:product-fibrewise-smooth-charts} by taking the direct sum of
the inner products on the (fibrewise) tangent spaces.  The naturality
and continuity of this construction is summarised as follows: 

\begin{lem}
The product of Kuranishi charts with tangent spaces defines a monoidal
structure on $\Kur^{\scrT}$, for which the forgetful functor to $\Kur$
is strictly monoidal. \qed 
\end{lem}

We again will consider the restricted subcategories of
$\Kur^{\scrT}$ relative to a fixed partially-ordered subset
$\cS$, as in Definition~\ref{defn:kuriso} and ~\ref{defn:kurisofs}.

\begin{defin}
The subcategory $\Kur^{\scrT,\iso}$ consists of those morphisms
for which the isomorphism $\cS_0 \cong \partial^Q \cS_1$ is the
identity, and  the full subcategories $\Kur^{\scrT,\cS} \subset
\Kur^{\scrT,\iso}$ are spanned by the objects of the form $\partial^Q
\cS$ as $Q$ varies over the totally-ordered subsets of $\cS$.
\end{defin}

Returning to the context of Definition
\ref{def:flag-smooth-presentation}, we have the following lift: 
\begin{defin} \label{def:flag-smooth-split-presentation}
A \emph{Kuranishi flow category with tangent bundles} consists of a
flag smooth Kuranishi presentation $\bX \co A \to \Kur^{fs}$, a
$\Pi$-equivariant topological $2$-category $A^{\scrT}$ over $A$, and a
commutative diagram
\begin{equation}
\begin{tikzcd}
A^{\scrT} \ar[r] \ar[d] & \Kur^{\scrT} \ar[d] \\
A \ar[r] & \Kur^{fs}.
\end{tikzcd}
\end{equation}
We require that:
\begin{enumerate}
\item the spaces of objects and of morphisms of $
A^{\scrT}$ lifting each object and each morphism in $A$ be
contractible,
\item that the target maps from morphisms to objects in
$A^{\scrT}$ be a fibration,  
\item and that the induced functor from
$A$ to $\Orb$ given by $Z^\ori$ be a stratified orbispace
presentation.
\end{enumerate}
\end{defin}

As described in Appendix~\ref{sec:internal}, a topological bicategory
is a bicategory where the categories of morphisms are internal
categories in spaces.  A $\Pi$-equivariant topological bicategory is
defined completely analogously as in the non-enriched setting,
requiring that the relevant action maps are given by internal
functors.  Note also that for each category of morphisms, the fact
that the target of the functor $A^{\scrT}(p,q) \to A(p,q)$ is discrete
means that the induced maps on objects and morphisms are locally
constant.

\subsubsection{Categories of Kuranishi charts with relative orientations}
\label{sec:categ-kuran-charts}

The constructions of the previous sections lead to a well-defined and functorial notion of Kuranishi charts equipped with tangent spaces. In this section, we introduce a notion of stable almost complex structure for such tangent spaces. More precisely, the situation we shall encounter when studying Floer theory in Section \ref{sec:stable-almost-compl-1} is that of (relative) stable complex structures (see Appendix \ref{sec:orientations} for some basic theory of relative orientations)

To this end, we consider a pair $V_0 = (V^+_0,V^-_0)$ and $V_1 = (V^+_0,V^-_1)$ of stable vector spaces, with complex structures on $V^-_0$ and $V^-_1$. 
\begin{defin} \label{def:Complex-oriented-Kuranishi-chart}
  The category $\Kur^{\ori}(V_0, V_1)$ of \emph{Kuranishi charts equipped with stable complex structures relative $V_0$ and $V_1$,} is the internal category in spaces with objects consisting of
  \begin{enumerate}
  \item an object $\alpha$ of $\Kur^{\scrT}$,
  \item a complex $G_\alpha$-equivariant vector bundle $I^\bC_\alpha$ over $Z^\ori_\alpha$,
  \item a complex $G_\alpha$-representation $W$,
  \item a finite set $O_\alpha$, and  
  \item a $G_\alpha$-equivariant isomorphism 
\begin{equation} \label{eq:stable-isomorphism-tangent-with-complex}
  V_1^+ \oplus \bR^{O_\alpha}  \oplus W_\alpha \oplus TX_\alpha \oplus V_0^- \cong  V_1^- \oplus I^\bC_\alpha \oplus W_\alpha \oplus V_0^+
\end{equation}
of vector bundles over  $Z^\ori_\alpha$.
  \end{enumerate}
The topology on the space of objects is given by the topology on the space of such isomorphisms (and the topology on objects of $\Kur^{\scrT}$).

A morphism in $\Kur^{\ori}(V_0, V_1)$ consists of
\begin{enumerate}
\item a morphism $f \co \alpha \to \beta$ in $\Kur^{\scrT}$ (recall that this entails a choice of a totally ordered subset $Q_f$ of $\cS_\beta$),
\item a bijection $O_\alpha \cong  O_\beta \amalg Q_f$,
\item an isomorphism
\begin{equation} \label{eq:splitting-complex-bundle-morphism}
  I^\bC_\beta \cong I^\bC_\alpha \oplus V_f^{\perp}
\end{equation}
of vector bundles on $Z_\alpha$, and
\item a $G_\alpha$-equivariant embedding $W_\alpha \to W_\beta$, whose quotient we denote $W_f^\perp$ such that the following diagram commutes:
\begin{equation}
\begin{tikzcd}[column sep=small]
  \begin{gathered}
W_f^\perp  \oplus V_f^{\perp}  \oplus \\  V_1^+ \oplus  \bR^{O_\alpha} \oplus  TX_\alpha \oplus W_\alpha \oplus V_0^-   \end{gathered}\ar[r] \ar[d] &
\begin{gathered}
W_f^\perp  \oplus V_f^{\perp}  \oplus \\ V_1^- \oplus I^\bC_\alpha \oplus W_\alpha  \oplus   V_0^+ \end{gathered} \ar[d] \\
V_1^+ \oplus \bR^{O_\beta}  \oplus W_\beta  \oplus  TX_\beta \oplus  V_0^- \ar[r] &  V_1^- \oplus I^\bC_\beta \oplus W_\beta \oplus V_0^+,
  \end{tikzcd}
\end{equation}
where we use the isomorphism $  V_f^{\perp} \oplus  TX_\alpha \cong \bR^{Q_f} \oplus TX_\beta$ obtained by combining Lemma \ref{lem:tangent-space-split-morphisms-well-defined} and Lemma \ref{lem:tangent_space_chart_splits_at_boundary}.
\end{enumerate}

We topologise the space of morphisms by taking the product of the topology of morphisms in $\Kur^{\scrT}$ with the topology on the space of splittings in Equation \eqref{eq:splitting-complex-bundle-morphism} and the topology on the space of objects.
\end{defin}

We define $\Kur^{\ori} \langle  \cS \rangle(V_0,V_1)$ to be the subcategory of objects lifting $\Kur^{\scrT} \langle  \cS \rangle $. This is the essential case to consider, and the general case is only included to facilitate the discussion of products in Section \ref{sec:bicat-kuran-charts} below.

\subsubsection{The bicategory of Kuranishi charts with complex structures}
\label{sec:bicat-kuran-charts}
We now study the multiplicativity of the construction. We begin by
constructing a bicategory with $1$-cells given by the categories from
Section \ref{sec:categ-kuran-charts}.

The key construction is the following:
\begin{defin}
There is a functor
  \begin{equation} \label{eq:product-oriented-Kuranishi-charts}
    \Kur^{\ori}(V_0, V_1) \times \Kur^{\ori}(V_1, V_2) \to \Kur^{\ori}(V_0, V_2)
  \end{equation}
which assigns to a pair $\alpha \in \Kur^{\ori}(V_0, V_1)$ and $\beta
\in \Kur^{\ori}(V_1, V_2)$ a product $\alpha \times \beta \in
\Kur^{\ori}(V_0, V_2)$ given by
\begin{enumerate}
\item the product of the underlying objects in $ \Kur^{\scrT}$, 
\item the union $O_{\alpha \times \beta} \equiv O_\alpha \amalg
  O_\beta$, 
\item the product vector bundle $I^\bC_{\alpha \times \beta} \equiv  I^\bC_\alpha \times I^\bC_\beta $
  over $Z^{\ori}_{\alpha \times \beta}$, 
\item the $G_{\alpha \times \beta}$ complex representation
  \begin{equation}
 W_{\alpha \times \beta} \equiv W_\alpha \oplus V_1^- \oplus W_\beta,
  \end{equation}
\item and the composite isomorphism 
\begin{align}
  V_2^+ \oplus \bR^{O_{\alpha \times \beta}}  &\oplus W_{\alpha \times
    \beta} \oplus TX_{\alpha \times \beta} \oplus V_0^- 
\\
& \cong V_2^+ \oplus W_\beta \oplus TX_{\beta}  \oplus  \bR^{O_{\beta}}  \\
& \qquad \oplus V_1^- \oplus \bR^{O_{\alpha}} \oplus  W_{\alpha} \oplus TX_{\alpha}  \oplus V_0^-  \\
& \cong V_2^- \oplus  W_{\beta} \oplus I^\bC_\beta \oplus V_1^+ \oplus W_\alpha \oplus TX_{\alpha}  \oplus  \bR^{O_{\alpha}}  \oplus V_0^+ \\
&  \cong  V_2^- \oplus  W_{\beta} \oplus I^\bC_\beta \oplus  V_1^- \oplus W_\alpha \oplus I^\bC_{\alpha} \oplus V_0^+ \\
& \cong V_2^- \oplus  W_{\alpha \times \beta} \oplus I^\bC_{\alpha \times \beta} \oplus V_0^+.
\end{align}
\end{enumerate}
\end{defin}

Since the only operations used in the construction of the product are
disjoint unions of sets, direct sum of vector spaces, and products of
topological spaces, we can immediately deduce the following lemma.

\begin{lem}
There are natural homeomorphism between the composites
\begin{equation}
\begin{tikzcd}
  \begin{gathered}
    \Kur^{\ori}(V_0, V_1) \times \Kur^{\ori}(V_1, V_2) \\ \times \Kur^{\ori}(V_2, V_3)  \end{gathered}  \ar[r] \ar[d] &
  \begin{gathered}
\Kur^{\ori}(V_0, V_2) \\ \times \Kur^{\ori}(V_2, V_3)  \end{gathered} \ar[d] \\
    \Kur^{\ori}(V_0, V_1) \times \Kur^{\ori}(V_1, V_3) \ar[r] &  \Kur^{\ori}(V_0, V_3),
  \end{tikzcd}
\end{equation}
and these associativity diagrams are coherent. \qed
\end{lem}

The above discussion allows us to construct a bicategory as follows:

\begin{defin}
  The topological bicategory $\Kur^{\ori}$ of \emph{Kuranishi
    charts with relative orientations} has $0$-cells given by pairs $V
  = (V^+,V^-)$ consisting of a real vector space $V^+$ and a complex
  vector space $V^-$, $1$-cells given by the categories
  $\Kur^{\ori}(V_0, V_1)$ and composition of $1$-cells given by
  Equation \eqref{eq:product-oriented-Kuranishi-charts}.
\end{defin}

We now have all the structural results required to define the key
notion of this paper.

\begin{defin} \label{def:complex-oriented-flow-category}
A \emph{complex oriented Kuranishi flow category} consists of a
Kuranishi flow category $\bX \co A \to \Kur$ together with a
$\Pi$-equivariant topological $2$-category $A^{\ori}$ with $0$-cells
$\cP$, equipped with a strictly $\Pi$-equivariant $2$-functor
$A^{\ori} \to A$ and a strictly $\Pi$-equivariant $2$-functor
$A^{\ori} \to \Kur^{\ori}$.  Here we use the natural isomorphisms of
Lemma~\ref{lem:kurequi} and the assignment described in
Equation~\eqref{eq:degree_basepoint_lifted_orbit} to express the
compatibility with the $\Pi$-action.

We require that:
\begin{enumerate}
\item the diagram
\begin{equation}
\begin{tikzcd}
A^{\ori} \ar[r] \ar[d] &  A \ar[d] \\
\Kur^{\ori} \ar[r] & \Kur
\end{tikzcd}
\end{equation}
commutes, 
\item the $2$-functor $A^{\ori} \to A$ induces an acyclic
  fibration of nerves for each $A^{\ori}(p,q) \to A(p,q)$, 
\item and the functors $A(p,q) \to \Orb$ induced by $Z^\ori$ are
  stratified orbispace presentations. 
\end{enumerate}
\end{defin}

\section{The virtual cochains of Kuranishi presentations}
\label{sec:impl-atlas-with}

\subsection{The virtual cochains of a Kuranishi chart}
\label{sec:first-non-trivial-example-cochains}

In this section, we introduce our second model for cochains, whose
existence depends on a choice of Kuranishi presentation.

\begin{notn}
We will use the notational convention that for a pair of spaces $A
\subseteq B$, the symbol $B|A$ denotes the homotopy cofiber (mapping
cone) usually denoted $C(B, B\setminus A)$.  Explicitly, this is the
union $B \cup C(B \setminus A)$, where the basepoint of $C(B \setminus
A)$ is given by the cone point $1$.  If $A$ is a based space, with
basepoint disjoint from $B$, we use the same notation for the cofiber
in the category of based spaces, i.e. we collapse the cone on the
basepoint.  Note that the homotopy cofiber is a functor of pairs.
\end{notn}

Associated to each $\langle \cS \rangle$-Kuranishi chart $\bX$ is the
homotopy cofiber $X|Z$ of the inclusion of the complement of the
zero-locus $Z = s^{-1}(0)$ into the domain $X$ of the chart.  Because
the section $s \colon X \to V$ is $G$-equivariant, $X | Z$ is a
$G$-space with basepoint given by the cone point.  Applying the Borel
construction and appropriate shifts, we have the following definition,
where we adopt Pardon's terminology~\cite{Pardon2016}.

\begin{defin} \label{def:virtual_cochains_chart}
The \emph{spectrum of virtual cochains} of a Kuranishi chart is the
Borel construction
 \begin{equation}
 B X|Z ^{-V} \equiv C_*(EG; F(S^{V} , ( X | Z)^{\mfib} ))_{G} 
= EG_+  \sma_G  F(S^{V}, ( X | Z)^{\mfib})
\end{equation}
where here $(-)^{\mfib}$ denotes the monoidal fibrant replacement
functor (see Definition~\ref{def:lax-fibrant}) in the category of
orthogonal $G$-spectra on the trivial universe.
\end{defin}

\begin{rem}
Although it might appear that we would want to use the fibrant
replacement functor on the complete universe, since we are considering
the Borel homotopy type of $F(S^V, (X | Z)^{\mfib})$ it suffices to
work over the trivial universe.
\end{rem}

Another way of thinking of the spectrum of virtual cochains is to
view it as the ``total spectrum'' of a fiberwise spectrum over $BG$
with fiber the desuspension of $X | Z$ by $V$.

It is convenient to break up the construction of virtual cochains into
two steps: the first step associates to a Kuranishi chart $\bX$ the
triple $(X|Z, V, G)$.  The second step is the Borel functor from such
triples to spectra.  The essential problem is that the first step is
not functorial: a map $f \co \alpha \to \beta$ of Kuranishi charts
does not \emph{canonically} induce a map
\begin{equation}
 V_f^\perp \times X_\alpha|Z_\alpha \to X_\beta|Z_\beta,
\end{equation}
though such a map arises from a choice of splitting of the normal
bundle of the image of $X_\alpha$ in $X_\beta$ (appropriately
interpreted, since the underlying manifolds and maps need not be
smooth).  Fixing such a choice, we would obtain a canonical map
\begin{equation}
  V_f^\perp|0 \wedge X_\alpha | Z_\alpha \to X_{\beta} | Z_{\beta},
\end{equation}
but the desired map of desuspensions would require a comparison between
$V_f^\perp|0 $ and the standard sphere $S^{V_f^\perp}$.  Keeping track
of the equivariance, coherence, and multiplicativity of such choices
of map is technically demanding. Instead, our solution will produce a
zig-zag, obtained from Pardon's degeneration to the normal
cone.

\subsection{Degeneration along an interval}
\label{sec:degen-along-an}

Our goal in this section is to associate to an arrow $f \co \alpha \to
\beta$ of Kuranishi charts in $\Kur \langle \cS \rangle$ a canonical
map of virtual cochains in the homotopy category, which is realised by
a zig-zag of maps that are defined independently of any choices. To
simplify our discussion, we shall assume, as discussed in Remark
\ref{rem:surjection_groups_Kuranishi-charts}, that the map $G_\alpha
\to G_\beta$ is a surjection.  This implies that $V_\alpha$ and hence
$V_f^\perp$ and $X_f = X_\alpha \times_{G_\alpha} G_\beta$ are
$G_\beta$-spaces.  In particular, the direct sum decomposition
$V_\alpha \oplus V_f^\perp \cong V_\beta$ and the inclusion $X_f
\subset X_\beta$ are $G_\beta$-equivariant.

Consider the set of triples 
\begin{equation} \label{eq:degeneration_normal_cone_one_dimension}
\{(t,x,v) \subset [0,1] \times X_{\beta} \times V_\beta \, \mid \,
s_{\beta}(x) = t \pi_f( v)\},
\end{equation}
where $\pi_f$ is the projection from $V_\beta$ to its subspace
$V_\beta/V_\alpha = V_f^\perp$.  The transversality assumption for
maps of Kuranishi charts implies that the fiber over $t = 0$ contains
the product $X_f \times V_f^\perp$ as an open subset.
\begin{defin}
The \emph{degeneration to the normal cone of $f$}, denoted  $X_{\alpha,\beta}$, is  the space obtained from Equation~\eqref{eq:degeneration_normal_cone_one_dimension} by removing the
complement of $X_f \times V_f^\perp$ from the fibre over $0$. 
\end{defin}
Note that the projection induces a map $X_{\alpha, \beta} \to [0,1]$, such that the
fiber at $0$ is canonically identified with $X_f \times V_f^\perp$,
and the fiber away from $0$ with $X_\beta$.
\begin{rem} \label{rem:surjection-implies-bundle-trivial}
If we drop the assumption that $G_\alpha \to G_\beta$ is a surjection,
then $V_f^\perp$ fails to be a $G_\beta$ representation, but the
construction above can be replaced by observing that $(X_\alpha \times
V_f^{\perp}) \times_{G_\alpha} G_\beta$ is naturally a
$G_\beta$-equivariant vector bundle over $X_f$. The surjectivity
assumption thus allows us to avoid considering general vector bundles
in our notions of morphisms of Kuranishi charts.
\end{rem}

For the next result, we write $\bbone$ for the partially ordered set $0 < 1$. 
\begin{lem} \label{lem:degeneration_gives_Kuranishi-chart}
The $G_\beta$ action on $X_\beta$ extends naturally to
$X_{\alpha,\beta}$ in such a way that the projection map
$s_{\alpha,\beta} \co X_{\alpha,\beta} \to V_\beta$ is
$G_\beta$-equivariant. The quadruple $\bX_{\alpha,\beta} \equiv
(X_{\alpha,\beta},V_\beta, s_\beta , G_\beta)$ is a $\langle \cS
\times \bbone \rangle$-Kuranishi chart, equipped with natural maps 
\begin{equation}
(\bX_\alpha \times V_f^\perp, \cS ) \to ( \bX_{\alpha,\beta} , \cS
  \times \bbone) \leftarrow (\bX_\beta, \cS)
\end{equation}
of stratified Kuranishi charts, where the required isomorphisms of
partially ordered sets are given by the evident identifications 
\begin{equation}
\cS \times \{0\} \cong \partial^{\cS \times \{0\}} (\cS \times \bbone) \quad\textrm{and}\quad
\cS \times \{1\} \cong \partial^{\cS \times \{1\}} (\cS \times \bbone).
\end{equation}
\end{lem}

\begin{proof}
The projection $X_\beta \to V_f^\perp$ is locally trivial near each
point in $X_\alpha/G^\perp_f$ by assumption.  This implies that the
map $X_{\alpha,\beta} \to V_f^\perp$ is locally modelled, near each
point in $X_f$, after the product of $X_f $ with the subset of $[0,1]
\times V_f^\perp \times V_\beta$ given by elements $(t,u,v)$ such that
$u = t v$.  We may trivialise the second factor as $[0,1] \times
V_f^\perp$, which implies that $X_{\alpha,\beta}$ is a manifold.  We
conclude that $\bX_{\alpha,\beta}$ is a stratified Kuranishi chart,
and the existence of the desired maps associated to the boundary of
the interval follow by inspection.
\end{proof}

In particular, we obtain the $G_\beta$-space $Z_{\alpha,\beta}$  as
the inverse image of $0$ under $s_{\alpha,\beta}$. We also denote by
$Z_f$ the quotient of $Z_\alpha$ by $G_f^\perp$. 

\begin{rem}
An alternative point of view on the degeneration to the normal cone
is that the map $X_{\alpha,\beta} \to [0,1]$ is a (topological)
submersion whose fibre is an $\langle \cS \rangle$-Kuranishi
chart. One could thus formulate a notion of \emph{Kuranishi charts
  over a base}. Since the only such charts we encounter arise from
explicit constructions, the additional formalism would just complicate
matters. 
\end{rem}

Lemma~\ref{lem:degeneration_gives_Kuranishi-chart} induces a zig-zag 
\begin{equation}
V_f^\perp|0  \wedge X_\alpha | Z_\alpha \to  X_{\alpha,\beta} |
Z_{\alpha,\beta} \leftarrow   X_{\beta} | Z_{\beta}, 
\end{equation}
which suggests using $V_\alpha|0$ as a model for the sphere.  However, the
fact that the natural map 
\begin{equation}
  V_\alpha|0 \wedge V_f^{\perp}|0 \to V_\beta|0  
\end{equation}
fails to be a homeomorphism leads us to prefer using the standard
sphere $S^{V_\alpha}$.  We thus enlarge $X_{\alpha,\beta}$ to admit an
embedding of the product of $X_f \times S^{V_f^\perp}$, by considering
its closure $\bar{X}_{\alpha,\beta}$ under the inclusion 
\begin{equation}
[0,1] \times X_{\beta} \times V_\beta \subset \left( [0,1] \times
X_{\beta} \right)_+ \wedge S^{V_\beta} 
\end{equation}
where we recall that $S^{V_\beta}$ is the $1$-point compactification
of $V_\beta$ with basepoint at infinity.

The map $s_{\alpha,\beta}$ extends to a map
\begin{equation}
  \bar{X}_{\alpha,\beta} \to S^{V_\beta},
\end{equation}
given by the composition
\begin{equation}
X_{\alpha,\beta} \to  \left( [0,1] \times  X_{\beta}\right)_+ \sma
S^{V_\beta} \to S^{V_\beta}.
\end{equation}

At this stage, we note that we have a natural inclusion
\begin{equation} \label{eq:inclusion-wedge-cone-to-cone-wedge}
X_f|Z_f \wedge  S^{V_f^\perp} \to  \left( (X_{f})_+ \sma S^{V_f^\perp}
\right) | Z_{f}  
\end{equation}
as can be seen via the identification
\begin{equation}
X_f|Z_f \wedge  S^{V_f^\perp} \cong C((X_{f})_+
\sma S^{V_f^\perp} , (X_f \setminus Z_f)_+ \sma S^{V_f^\perp}),
\end{equation}
with the cofiber taken in the category of based spaces, and using the
inclusion  
\begin{equation} 
(X_f \setminus Z_f)_+ \sma S^{V_f^\perp} \subset ((X_{f})_+ \sma
  S^{V_f^\perp}) \setminus \left( Z_f \times \{0\} \right).
\end{equation}

The natural embeddings:
\begin{equation}
(X_{f})_+ \sma S^{V_f^\perp} \to \bar{X}_{\alpha, \beta} \leftarrow X_{\beta}
\end{equation}
which are respectively over the endpoints $0$ and $1$, thus induce
maps of $G_\beta$-spaces: 
\begin{equation}
S^{V_f^\perp} \wedge X_f | Z_f \to  \bar{X}_{\alpha, \beta}|
Z_{\alpha,\beta} \leftarrow X_{\beta} | Z_{\beta}. 
\end{equation}
      
Composing further with the map induced by projection from $X_\alpha$
to $X_f$ yields 
\begin{equation}
S^{V_f^\perp} \wedge X_\alpha | Z_\alpha \to  \bar{X}_{\alpha, \beta}| Z_{\alpha,\beta} \leftarrow X_{\beta} | Z_{\beta}.
\end{equation}

\begin{lem} 
The map $X_{\beta} | Z_{\beta}  \to \bar{X}_{\alpha, \beta}|
Z_{\alpha,\beta} $ is a $G_\beta$-equivariant homotopy equivalence. 
\end{lem}

\begin{proof}
This map factors through the inclusion $X_{\alpha, \beta}|
Z_{\alpha,\beta} \subset \bar{X}_{\alpha, \beta}| Z_{\alpha,\beta} $,
which is a homotopy equivalence because the inclusion $V|0 \subset
S^V|0$ is a homotopy equivalence.  The inclusion of the interior of a
manifold with boundary is a homotopy equivalence, and the interior of
$X_{\alpha,\beta}$ is the product of $X_\beta$ with an interval. 
\end{proof}

Desuspending by $V_\beta$, we obtain a diagram of $G_\beta$-spectra: 
\begin{equation} \label{eq:maps_desuspended_virtual_cochains}
F( S^{V_\alpha}, (X_f | Z_f)^{\mfib}) \to  F(S^{V_\beta},
(\bar{X}_{\alpha, \beta}| Z_{\alpha,\beta})^{\mfib}) \leftarrow
F(S^{V_\beta}, (X_{\beta} | Z_{\beta})^{\mfib}). 
\end{equation}
where the first map is obtained by decomposing $S^{V_\beta}$ as
$S^{V_\alpha} \wedge S^{V_f^\perp}$, and applying the identity on
the $S^{V_f^\perp}  $ factor. 
 
Applying the (desuspended) Borel construction, we obtain a diagram of
non-equivariant spectra 
\begin{equation}
BX_\alpha | Z_\alpha^{-V_\alpha}  \to  B \bar{X}_{\alpha, \beta}|
Z_{\alpha,\beta}^{-V_\beta}  \leftarrow   B X_{\beta} |
Z_{\beta}^{-V_\beta} 
\end{equation}
with the property that the left-pointing arrow is an equivalence. 

At this stage, we note that, since there is a $G_\beta$-homeomorphism
$X_f | Z_f \cong X_\alpha | Z_\alpha \ltimes_{G_\alpha} G_\beta$, the
fact that the $G_f^\perp$ action on $X_\alpha | Z_\alpha$ is free implies
the following corollary of Lemma \ref{lem:equivalence_borel_free_action}:

\begin{lem}
The projection from $X_\alpha | Z_\alpha$ to $ X_f | Z_f$ induces a
natural equivalence 
\begin{equation}
B X_\alpha | Z_\alpha^{-V_\alpha}  \cong 
B X_f | Z_f^{-V_\alpha}.
\end{equation} \qed
\end{lem}

\subsection{Simplicial degeneration to the normal cone}
\label{sec:simpl-degen-norm}

While our ultimate goal will be to construct a degeneration to the
normal cone for cubical diagrams in $\Kur \langle \cS \rangle$, we
begin by considering the simplicial case. Let $\bDelta \Kur \langle
\cS \rangle$ denote the category of simplices of $\Kur \langle \cS
\rangle$, i.e., objects are functors $\bbn \to \Kur \langle \cS
\rangle$, where $\bbn$ is the ordered set $(0, \ldots, n)$ and
morphisms are commutative diagrams.
  
We assign to each simplex $\sigma \co \bbn \to \Kur \langle \cS
\rangle$ an $\langle \cS \times \bbn \rangle$-Kuranishi chart as
follows: define 
\begin{align}
  V_\sigma & \equiv V_{\sigma(n)} \\
  G_\sigma & \equiv G_{\sigma(n)},
\end{align}
and consider the submanifold
\begin{equation}
X_{\sigma} \subset \Delta^n \times X_{\sigma(n)} \times V_{\sigma(n)}
\end{equation}
consisting of elements $(t_0, \ldots, t_{n}, x, v)$ with $ \sum t_i = 1$, such that the following two conditions hold:
\begin{enumerate}
\item The projection of $x$ to $V_\sigma$ satisfies 
\begin{equation}
s_{\sigma(n)}(x) = \sum_{i=0}^{n} T_i \cdot \pi_i v,
\end{equation}
where $\pi_i$ is the projection from $ V_{\sigma(n)}$ to the orthogonal complement of $V_{\sigma(i)}$ in $V_{\sigma(i+1)}$, and $T_i = \sum_{i \leq j} t_j$.
\item Whenever $t$ lies the image of $\Delta^m \to \Delta^n$, for a map $\iota \co \bbm \to \bbn$, then $x$ lies in the image of $X_{\sigma(\iota(m))}$  in $X_{\sigma(n)}$.
\end{enumerate}
The action of $G_\sigma$ on $X_{\sigma(n)} \times V_{\sigma(n)}$ yields an action on $ X_{\sigma}$, for which the projection map $s_\sigma$ to $V_\sigma$ is equivariant. As in Lemma \ref{lem:degeneration_gives_Kuranishi-chart}, the quadruple $(X_\sigma, V_\sigma, s_\sigma, G_\sigma) $ is an $\langle \cS \times \bbk \rangle$-Kuranishi chart.

We now consider the partial compactification
\begin{equation}
   X_\sigma \subset  \bar{X}_\sigma
\end{equation}
obtained as the closure of $X_\sigma$ in $\left( \Delta^n \times
X_{\sigma(n)} \right)_+ \sma S^{V_{\sigma(n)}} $. The group $G_\sigma$
acts on $\bar{X}_\sigma$, and we have a natural inclusion $Z_\sigma
\equiv s^{-1}_{\sigma}(0) \subset \bar{X}_\sigma$, given as the
inverse image of the origin in $S^{V_\sigma}$.  It is useful to give a
slightly more formal definition at this stage:
\begin{lem}\label{lem:pushout_description_compactified_degeneration}
There is a natural pushout diagram of based spaces
\begin{equation}
\begin{tikzcd}
\displaystyle{\coprod_{\rho \overset{g}{\to} \tau \overset{f}{\to}
    \sigma}  (X_{f \circ g} \times V_g^\perp)_+ \sma S^{V_f^\perp} }
\ar[r] \ar[d] &      \displaystyle{\coprod_{\rho \overset{f \circ
      g}{\to} \sigma}     (X_{f \circ g})_+ \sma S^{V_{f \circ
      g}^\perp} }\ar[d] \\ 
\displaystyle{\coprod_{ \tau \overset{f}{\to} \sigma}  (X_{f})_+ \sma 
  S^{V_f^\perp} }\ar[r] & \bar{X}_\sigma, 
\end{tikzcd}
\end{equation}
where the coproducts are taken over injective maps. \qed
\end{lem}

Given a map $f \co \tau \to \sigma$ of simplices in $\Kur \langle \cS
\rangle$, we write $X_f$ for  the image of $X_\tau$ in $X_\sigma$. We
write $f$ as well for the corresponding map $\Delta^\tau \to
\Delta^\sigma$ of simplices, and $\Delta^f$ for its image. Let
$V_f^\perp$ denote the orthogonal complement of $V_\tau$ in
$V_\sigma$.  

\begin{lem} \label{lem:orthogonal_splitting_Pardon_degeneration}
The fiber of $X_\sigma$ over the interior of $ \Delta^f$ projects
homeomorphically to the product of $V_f^\perp$  with $X_{f}$. The
closure of this fiber is homeomorphic to $V_f^\perp \wedge
\bar{X}_f$. 
\end{lem}

\begin{proof}
In the expression $\sum T_i \pi_i v$ defining $X_\sigma$, all
coefficients $T_i$ vanish along $\Delta^f$ for integers $i$ which
larger than the label of each vertex in the image of $f$, while all
other coefficients are either non-zero, or agree with the
corresponding coefficient defining $X_f$. This implies that the fiber 
of $X_\sigma$ over $\Delta^f$ splits as the product of $V_f^\perp$
with $X_f$. 
\end{proof}

\begin{cor} \label{cor:simplicial_decomposition_restricts}
If the image of $f$ contains the maximal vertex of  $\Delta^\sigma$,
the induced map 
\begin{equation}
X_f \to X_\sigma \times_{\Delta^\tau} \Delta^\sigma
\end{equation}
is a $G_{\sigma}$-equivariant homeomorphism. \qed
\end{cor}

\begin{rem}
At this stage, we could define the \emph{simplicial virtual cochains}
of a functor $A \to \Kur \langle \cS \rangle$ to be the homotopy
colimit of the desuspensions 
\begin{equation}
\hocolim_{\sigma \in \bDelta A}  B\bar{X}_\sigma | Z_\sigma^{-V_\sigma}.    
\end{equation}
In such a setting, it would be possible to formulate the
multiplicativity of the degeneration of the normal cone in terms of
bi-simplicial objects.  We shall find it much more convenient to use
cubical objects. The reader should have in mind the analogous problem
for proving multiplicativity of the Serre spectral sequence, which can
be established by a bi-simplicial argument \cite{Dress1967} replacing
the original cubical proof. 
\end{rem}

\subsection{Bi-simplicial degenerations}
\label{sec:bi-simpl-degen}
In this section, we compare the product of degenerations with the degeneration over a bi-simplex associated to its natural (prismatic) simplicial subdivision.

Consider simplices $\sigma \co \bbn \to \Kur \langle \cS_1 \rangle$ and $\tau \co \bbm \to \Kur \langle \cS_1 \rangle$. The product of the corresponding degenerations is a space
\begin{equation}
   X_\sigma \times  X_\tau  \to \Delta^{n} \times \Delta^m.  
\end{equation}
We would like to express this spaces as a union (colimit) of degenerations over the standard triangulation of the product of simplices. To this end, consider a map $ f \co \bbk \to \bbn \times \bbm$, yielding a map
\begin{equation} \label{eq:map_simplex_product}
  f \co \Delta^k \to   \Delta^{n} \times \Delta^m,
\end{equation}
and let
\begin{equation}
   (\sigma \times  \tau) \circ f \co   \bbk \to \Kur \langle \cS_1 \cup \cS_2 \rangle
 \end{equation}
 denote the composition of $f$ with the product of $\sigma$ and $\tau$.
\begin{lem}\label{lem:product_degenerations_agree}
  If  $f(k) = (n,m)$, there is a natural homeomorphism
\begin{equation}
  f^*(  X_\sigma  \times  X_\tau) \to  X_{(\sigma \times  \tau) \circ f}
\end{equation}
over $\Delta^k$.
\end{lem}
\begin{proof}
The fibers of both spaces over a point in $\Delta^k$ are contained in
     \begin{equation}
      X_{\sigma(n)} \times  X_{\tau(m)}  \times V_{\sigma(n)} \times V_{\tau(m)},       
    \end{equation}
    and we shall show that they are equal. To this end, note that the projection $\pi_i $ associated to $i \in \bbk$ is the direct sum of the projections $\pi_{f_{\bbn}(i)} \oplus \pi_{f_{\bbm}(i)} $ on $V_{\sigma(n)} \times V_{\tau(m)}$, which are given by
    \begin{align}
      \pi_{f_{\bbn}(i)} & = \sum_{\pi_{f_{\bbn}(i-1)} <  j \leq \pi_{f_{\bbn}(i)} } \pi_j \\
      \pi_{f_{\bbm}(i)} & = \sum_{\pi_{f_{\bbm}(i-1)} <  j \leq \pi_{f_{\bbm}(i)} } \pi_j.
    \end{align}
A collection of elements $(x,y) \in  X_{\sigma(n)} \times  X_{\tau(m)}  $ and $(u,v) \in V_{\sigma(n)} \times V_{\tau(m)}$ lies in the degeneration associated to the composition if we have
    \begin{align}
      \left(s_{\sigma}(x), s_{\tau}(y) \right) & = \sum_{i=0}^{k} T_i \left( \sum_{\pi_{f_{\bbn}(i-1)} <  j \leq \pi_{f_{\bbn}(i)} } \pi_j    u,  \sum_{\pi_{f_{\bbm}(i-1)} <  j \leq \pi_{f_{\bbm}(i)} } \pi_j  v \right) \\
      &  = \left( \sum_{j=0}^{n} \left(\sum_{   j \leq \pi_{f_{\bbn}(i)} }t_i\right) \pi_{j} u, \sum_{j=0}^{m}  \left(\sum_{   j \leq \pi_{f_{\bbm}(i)} }t_i\right)   \pi_{j}  v \right) \\
     &  =  \left( \sum_{j=0}^{n} T_j \pi_{j} u, \sum_{j=0}^{m}  T_j \pi_{j}  v \right).  
\end{align}
In the last step, we use the fact that the map in Equation \eqref{eq:map_simplex_product} is given in coordinates by $t_{f_{\bbn}(i)} =  t_{f_{\bbm}(i)} = t_i$, and all other coordinates vanish. The last expression above is exactly the formula for the pullback of degenerations.
\end{proof}

\subsection{Cubical degeneration to the normal cone}
\label{sec:cubic-degen-norm}

Denote by $\bbone$ the category $0 \to 1$ as before, and $\bbone^n$
the product category, which is equipped with canonical morphisms
called \emph{edges} which are the product of identities in all but one
of the factors.  We refer to a functor $\sigma \co \bbone^{n} \to \Kur
\langle \cS \rangle$ as an $n$-\emph{cube} in $\Kur \langle \cS
\rangle$.  Cubes in $\Kur \langle \cS \rangle$ form a category
$\Box \Kur \langle \cS \rangle$ where a morphism $\sigma \to \tau$
is a factorisation of $\sigma$ through a \emph{map of cubes}
$\bbone^{n} \to \bbone^{m}$.  We shall only consider the subcategory
$\Box_{\nd}\Kur \langle \cS \rangle$ consisting of \emph{non-degenerate}
cubes, i.e., those which do not factor through a degeneracy.  
Since we exclude degenerate cubes, a map of cubes is given by the
inclusion of a facet, which can be combinatorially recorded as a map
from a subset of $\{1, \ldots, m\}$ to $\{0,1\}$.

As explained in Theorem~\ref{thm:boxnd}, in the context in which we
work, the evident inclusion functor $\Box_{\nd} \Kur \to \Box \Kur$
induces an equivalence on homotopy colimits.  As a consequence, our
choice to ignore degenerate cubes is harmless.  Since we never work
with the full category of cubes in the main text, in abuse of notation
we will henceforth write $\Box \aC$ to denote $\Box_{\nd} \aC$ (for
any category $\aC$).

Let $\square^{n}$ denote the product $[0,1]^n$.  

\begin{defin}
Let $\sigma \co \bbone^{n} \to \Kur \langle \cS \rangle$ be a cube. 
The {\em cubical degeneration to the normal cone} is the union of the
simplicial degeneration of the normal cones 
\begin{equation}
X_{\sigma}\equiv \colim_{\iota \co \bbk \to \bbone^n} X_{\sigma \circ \iota},
\end{equation}
where the colimit is over the standard simplicial subdivision of the
square whose simplices are labelled by functors $\bbk \to \bbone^n$.
\end{defin}

This colimit is given by the union of the spaces associated to
top-dimensional simplices (injective order preserving maps $\bbn  \to
\bbone^n$), glued along codimension $1$ interior simplices. The
interior condition on facets is equivalent to the condition that $0^n$
and $1^n$ are respectively the minimal and maximal element of the
corresponding map  $\bbn \to \bbone^{n}$. Inductively applying Lemma
\ref{lem:product_degenerations_agree}, we conclude that $X_\sigma$ is
a manifold with boundary equipped with a projection map 
\begin{equation}
  X_\sigma \to \square^n,
\end{equation}
which is a
topological submersion, so that the inverse image of every boundary
stratum of the cube is a manifold (with boundary).

The above construction is equivariant with respect to the action of
$G_\sigma \equiv G_{\sigma(1^n)}$ and admits a natural equivariant
projection map to $V_\sigma \equiv V_{\sigma(1^n)}$, yielding the
Kuranishi chart 
\begin{equation}
(X_\sigma, V_\sigma,  s_\sigma, G_\sigma) \in \Kur \langle \cS \times
\bbone^k \rangle.
\end{equation}

We denote by $Z_\sigma \subset X_\sigma$ the inverse image of the
origin in $V_\sigma$. The closure of $X_\sigma$ in $\left( \square^{n}
\times X_{\sigma(1^n)} \right)_+ \sma S^{V_{\sigma(1^n)}}$ will be
denoted $\bar{X}_{\sigma}$ as before. 

Using the Borel construction, we obtain the spectrum of
\emph{cubical virtual cochains} by applying the construction of
Definition~\ref{def:virtual_cochains_chart} to a cube $\sigma$ in order to
obtain $B\bar{X}_\sigma | Z_\sigma^{-V_\sigma}$.

We now turn to explain the functoriality of this construction.
Given a morphism $f \co \tau \to \sigma$ of cubes in $\Kur \langle \cS
\rangle $, we write $\bX_f$ for the free quotient of $\bX_\tau$ of
$G_f^\perp$, and define $V_f^\perp$ as before to be the orthogonal
complement of $V_\tau$ in $V_\sigma$.
Lemma~\ref{lem:orthogonal_splitting_Pardon_degeneration} allows us to
associate to each such morphism a map of spaces  
\begin{equation} \label{eq:map-suspension-cube-face-degeneration}
 S^{V_f^\perp} \wedge \bar{X}_\tau | Z_\tau\to \bar{X}_\sigma | Z_\sigma,
\end{equation}
which is equivariant with respect to the homomorphism $G_\tau \to
G_\sigma$. 

\begin{rem}\label{rem:change-topology}
While we have an inclusion $X_\sigma \subset \bar{X}_\sigma$, the
induced map $X_\sigma | Z_\sigma \to \bar{X}_\sigma | Z_\sigma$ is in
fact bijective. To see this, note that $\bar{X}_\sigma$ is obtained
from $X_\sigma$ by adding a single point, but that this point is
identified with the cone point in $\bar{X}_\sigma | Z_\sigma $. Thus,
the map $X_\sigma | Z_\sigma \to \bar{X}_\sigma | Z_\sigma $ can be
alternatively described as changing the topology on $X_\sigma |
Z_\sigma$ in such a way that a sequence in $X_f \times V_f^\perp$
whose second coordinate goes to infinity converges to the cone point. 
\end{rem}

Combining the maps constructed above, each arrow $f \co \tau \to
\sigma$ yields a composition: 
\begin{equation}\label{eq:cube-Kur-compose}
\begin{aligned}
    B\bar{X}_\tau | Z_\tau^{-V_\tau} & \equiv    (E G_\tau)_+
    \sma_{G_\tau} F(S^{V_\tau} , ( \bar{X}_\tau | Z_\tau)^{\mfib} )
    \\ 
    & \to (EG_\sigma)_+ \sma_{G_\sigma}  F(S^{V_f} , ( \bar{X}_f | Z_f)^{\mfib} ) \\
    & \to (EG_\sigma)_+ \sma_{G_\sigma}  F(S^{V_f^\perp} \wedge S^{V_f} , S^{V_f^\perp} \wedge ( \bar{X}_f | Z_f)^{\mfib} )  \\
    & \to (EG_\sigma)_+ \sma_{G_\sigma}  F(S^{V_\sigma},  ( S^{V_f^\perp} \wedge \bar{X}_f | Z_f)^{\mfib} )  \\
    & \to  (EG_\sigma)_+ \sma_{G_\sigma}  F(S^{V_\sigma}, ( \bar{X}_\sigma | Z_\sigma)^{\mfib} ) \equiv B\bar{X}_\sigma | Z_\sigma^{-V_\sigma},
\end{aligned}
\end{equation}
where we used the quotient by $G_f^\perp$ in the first step, the
identity on $S^{V_f^\perp}$ in the second, the multiplicativity of the functor
$( \_)^\mfib$ in the third step, and
Equation~\eqref{eq:map-suspension-cube-face-degeneration} in the last
step.  Since all these maps are associative, we conclude:  

\begin{lem}
The assignment
\begin{equation}
\sigma \mapsto B\bar{X}_\sigma | Z_\sigma^{-V_\sigma}
\end{equation} extends to a functor $\square \Kur \langle \cS \rangle
\to \Sp$ which is specified on morphisms by the composition Equation~\eqref{eq:cube-Kur-compose}. \qed
\end{lem}
We shall call this functor the \emph{virtual cochains} functor.

Appealing to Lemma \ref{lem:product_degenerations_agree}, and using
the external multiplicativity of the fibrant replacement functor
$( \_)^\mfib$ (see Proposition~\ref{prop:Q-lax-mon}), we find that this
construction is multiplicative in the sense that the product of
Kuranishi charts induces a map 
\begin{equation}
B\bar{X}_\sigma | Z_\sigma^{-V_\sigma}  \wedge  B\bar{X}_\tau |
Z_\tau^{-V_\tau} \to  B\bar{X}_{\sigma \times \tau} | Z_{\sigma \times
  \tau}^{-V_{\sigma \times \tau}}.
\end{equation}
A straightforward check shows that this product is functorial in
$\sigma$ and $\tau$.  Moreover, the product map is itself naturally
associative for triples, in the sense that the following diagram
commutes:  
\begin{equation}
\begin{tikzcd}[column sep=small]
B\bar{X}_\sigma | Z_\sigma^{-V_\sigma}  \wedge  B\bar{X}_\tau |
Z_\tau^{-V_\tau} \wedge B\bar{X}_\rho | Z_\rho^{-V_\rho} \ar[r] \ar[d]
&  B\bar{X}_{\sigma \times \tau} | Z_{\sigma \times \tau}^{-V_{\sigma
    \times \tau}} \wedge B\bar{X}_\rho | Z_\rho^{-V_\rho}  \ar[d] \\ 
B\bar{X}_\sigma | Z_\sigma^{-V_\sigma}  \wedge B\bar{X}_{ \tau \times
  \rho} | Z_{ \tau\times \rho}^{-V_{\tau\times \rho}}
\ar[r]  &  B\bar{X}_{\sigma \times \tau \times \rho} | Z_{\sigma
  \times \tau\times \rho}^{-V_{\sigma \times \tau\times \rho}} 
\end{tikzcd}
\end{equation}
Moreover, it is evidently unital.

Putting this all together, we conclude:

\begin{lem}
The virtual cochains functor is lax monoidal. \qed
\end{lem}

\subsection{Virtual cochains of Kuranish flow categories}
\label{sec:virt-coch-homot}

We now have the tools at hand to extend the definition of virtual
cochains from charts to presentations, and from there to flow
categories. For the next definition, recall from
Lemma~\ref{lem:collar-functor} that we have a collar functor $\bX \to
\hat{\bX}$ from $\Kur^{\cS}$ to $\Kur \langle \cS \rangle $.

\begin{defin}
The  \emph{virtual cochains} of a Kuranishi presentation $\bX \co A
\to \Kur^{\cS}$  is the homotopy colimit
\begin{equation}
B\bar{X}|Z^{-V}(A) \equiv  \hocolim_{\sigma \in \square A}  B\bar{\hat{X}}_\sigma|\hat{Z}_\sigma ^{-V_\sigma}.
\end{equation}
\end{defin}
We often write $B\bar{X}|Z^{-V} $ for the virtual cochains, omitting the domain category from the notation.

This construction is functorial with respect to natural
transformations of Kuranishi presentations, and passing from a
presentation to a stratum labelled by an element $q$ yields a natural
morphism 
\begin{equation} \label{eq:map_virtual_cochains_stratum}
{ B  \partial^q\bar{X}| \partial^qZ^{-V}} \to  B\bar{X}|Z^{-V}.
\end{equation}

Moreover, we have the following multiplicativity property:

\begin{lem} \label{lem:map_virtual_cochains_product}
Given Kuranishi presentations $\bX(i) \co A(i) \to \Kur^{\cS_i}$ for
$i = \{1,2\}$  we have a natural equivalence 
\begin{equation}
B\bar{X}|Z^{-V}(A(1)) \sma B\bar{X}|Z^{-V}(A(2)) \to
B\bar{X}|Z^{-V}(A(1) \times A(2)),
\end{equation}
which is associative in the sense that, for a triple of Kuranishi
presentations with domain $A(i)$ for $i \in \{1,2,3\}$, the following
diagram commutes: 
\begin{equation}
\begin{tikzcd}[row sep=small, column sep=tiny]
  \begin{gathered}
    B\bar{X}|Z^{-V}(A(1))  \wedge   B\bar{X}|Z^{-V}(A(2)) \\
     \wedge   B\bar{X}|Z^{-V}(A(3))
  \end{gathered}
 \ar[d] \ar[r] &     \begin{gathered}
  B\bar{X}|Z^{-V}(A(1) \times A(2)) \\
     \wedge   B\bar{X}|Z^{-V}(A(3))
  \end{gathered}  \ar[d] \\
  \begin{gathered}
B\bar{X}|Z^{-V}(A(1))  \wedge  \\   B\bar{X}|Z^{-V}(A(2) \times A(3))    
  \end{gathered}
 \ar[r] &   B\bar{X}|Z^{-V}(A(1) \times A(2) \times A(3)),
\end{tikzcd}
\end{equation}
\end{lem}

\begin{proof}
There is a natural functor $\square A_1 \times \square A_2 \to
\square( A_1 \times A_2)$ given by the product of squares. Letting
$\bX_{12} \co A_1 \times A_2 \to \Kur^{ \cS_{1} \amalg \cS_{2}}$
denote the product presentation, and pulling back to $\square A_1
\times \square A_2 $ we have a natural transformation 
\begin{equation}
 B\bar{\hat{X}}_1|\hat{Z}_1^{-V_1} \wedge
 B\bar{\hat{X}}_2|\hat{Z}_2^{-V_2} \to B\bar{\hat{X}}_{12}
 |\hat{Z}_{12}^{-V_{12}} 
\end{equation}
of functors on $\square A_1 \times \square A_2 $.  This natural 
transformation is induced similarly to the discussion of the
multiplicativity of charts above, using the map
\begin{equation}
\begin{tikzcd}
(EG_1)_+  \sma_{G_1}  F(S^{V_1} , ( \bar{X}_1 \mid
Z_1)^{\mfib}) \sma
(EG_2)_+ \sma_{G_2}  F(S^{V_2} , ( \bar{X}_2 \mid
Z_2)^{\mfib} ) \ar[d] \\
E(G_1 \times G_2)_+  \sma_{G_1 \times G_2}  F(S^{V_{12}} , (\bar{X}_{12} \mid
Z_{12})^{\mfib} ).
\end{tikzcd}
\end{equation}

Composing all these maps with the covariance of $\hocolim$, we get
\begin{equation}
\hocolim_{\square A_1 \times \square A_2}    B\bar{\hat{X}}_1|\hat{Z}_1^{-V_1} \wedge  B\bar{\hat{X}}_2|\hat{Z}_2^{-V_2} \to  \hocolim_{\square \left( A_1 \times A_2\right)}    B\bar{\hat{X}}_{12} |\hat{Z}_{12}^{-V_{12}}.       
\end{equation}
Finally, as explained in Equation~\eqref{eq:hocosmash}, we have a natural map
\begin{multline}
\hocolim_{\square A_1}B\bar{\hat{X}}_1|\hat{Z}_1^{-V_1} \wedge    \hocolim_{\square A_2}    B\bar{\hat{X}}_2|\hat{Z}_2^{-V_2}    \\ \to   \hocolim_{\square A_1 \times \square A_2}    B\bar{\hat{X}}_1|\hat{Z}_1^{-V_1} \wedge  B\bar{\hat{X}}_2|\hat{Z}_2^{-V_2} 
  \end{multline}
whose composition with the above gives the desired result.
Associativity follows from the associativity of the constituents of
these transformations.
\end{proof}

We now return to the setting of Section \ref{sec:kuran-orbisp-flow}:
given a Kuranishi flow category $\bX$ with objects $\cP$ with a choice
of vector spaces as $V_p^\pm$ for each element $p \in \cP$ as in the
discussion preceding Remark \ref{rem:choose-slice-stable-index},
and a basepoint $p_0$ in each orbit of the action of $\Pi$ in $\cP$
which satisfies Equation \eqref{eq:degree_basepoint_lifted_orbit}, we
begin by defining a reduced degree
\begin{equation}
    \deg p \equiv \frac{d_{p_0} - d_p}{2}.
  \end{equation}
Recalling that our assumption is $\dim \hat{X} - \dim V = d_p - d_q - 1$, we then have the following equality.
  \begin{multline} \label{eq:dimension-moduli-space-degree-difference}
  1 + \dim \hat{X}_\alpha - \dim V_\alpha -  \dim V^+_{p} + \dim V^+_{q} + \dim V^-_{p} - \dim V^-_{q} \\ = 2 (\deg q - \deg p)
  \end{multline}
  for each Kuranishi chart corresponding to an object of $A(p,q)$.

  Given a cube $\sigma \in \square A(p,q) $, we define
  \begin{multline}
    B\bar{\hat{X}}_\sigma|\hat{Z}^{-V - d}_\sigma \equiv  \\
    C_*\left(E G_\sigma;  F(S^{V_\sigma}, (\hat{X}_\sigma | \hat{Z}_\sigma)^{\mfib}) \wedge S^{V^-_p + V^+_q - (V^+_p + V^-_q)} \right)_{G_\sigma} [ \deg p  - \deg q],
  \end{multline}
where the shift $\deg p  - \deg q $ corresponds to smashing with the spheres from Section \ref{sec:action}, and
 \begin{equation}
    S^{V^-_p + V^+_q - (V^+_p + V^-_q)} \equiv  S^{V^-_p + V^+_q} \wedge  S^{- (V^+_p + V^-_q)},
  \end{equation}
  using the Mandell-May model for negative spheres from Appendix \ref{sec:ortho}.
\begin{rem}
It might seem simpler to replace Equation
\eqref{eq:morphisms_category_virtual_cochains} by the Borel
equivariant chains of $F(S^{V_\sigma + V^+_p + V^-_q}, (S^{V^+_q +
  V^-_p} \wedge \bar{\hat{X}}_\sigma |
\hat{Z}_\sigma)^{\mfib})$. However, it is not clear to the authors how
to compose morphisms with this simpler model. On the other hand, it
would be possible at this stage to use the spheres $ F(S^{V_\sigma+
  V^+_p + V^-_q}, (X_\sigma | Z_\sigma)^{\mfib}) \wedge S^{V^-_p +
  V^+_q}$ as a model for desuspension instead of the Mandell-May
spheres. This will have the desired functoriality and multiplicativity
property, but does not allow us to cleanly separate the part of the
construction that are required for orientations from the rest.
\end{rem}

Next, for each pair $(p,q)$ of objects in $\cP$, we define the
spectrum of morphisms in the category of virtual cochains as the
homotopy colimit 
 \begin{equation} \label{eq:morphisms_category_virtual_cochains}
  B\bar{\cX}|\cZ^{-V - d}(p,q) \equiv \hocolim_{ \sigma \in \square A(p,q) }     B\bar{\hat{X}}_\sigma|\hat{Z}^{-V - d}_\sigma 
 \end{equation}
 over $\sigma \in \square A(p,q) $, 

\begin{lem}\label{lem:category_virtual_cochains_composition}
For a Kuranishi flow category $\bX$ and objects $p,q,r \in \cP$, there
is a natural associative composition map:
\begin{equation}
B\bar{\cX}|\cZ ^{-V-d}(p,q) \wedge  B\bar{\cX}|\cZ^{-V-d}(q,r) \to
B\bar{\cX}|\cZ^{-V-d}(p,r). 
\end{equation}
\end{lem}

\begin{proof}
Before desuspensions, the map is induced by Equations \eqref{eq:map_virtual_cochains_stratum} and Lemma \ref{lem:map_virtual_cochains_product}.  To incorporate shifts, it suffices to construct a functorial map for products of charts. This is obtained from the multiplicativity of the model of spheres  from Section \ref{sec:action} and from the multiplicativity of the fibrant replacement functor. Explicitly, given triples $p < q < r$, we use the composition
\begin{equation}
\begin{tikzcd}[row sep=small, column sep=tiny]
  \begin{gathered}
    F(S^{V_{\sigma_1} }, (\bar{\hat{X}}_{\sigma_1} |
  \hat{Z}_{\sigma_1})^{\mfib}) \wedge S^{V^-_p + V^+_q - (V^+_p +
    V^-_q)} \\
  \wedge    F(S^{V_{\sigma_2}},
  (\bar{\hat{X}}_{\sigma_2} | \hat{Z}_{\sigma_2})^{\mfib}) \wedge
  S^{V^-_r + V^+_p - (V^+_r + V^-_p)}
\end{gathered}
 \ar[d] \\ 
 \begin{gathered}
  F(S^{V_{\sigma_1} + V_{\sigma_2} }, (\bar{\hat{X}}_{\sigma_1} |
\hat{Z}_{\sigma_1} \wedge \bar{\hat{X}}_{\sigma_2} |
\hat{Z}_{\sigma_2})^{\mfib}) \wedge  S^{V^-_p - V^+_p} \wedge  S^{V^+_p  -  V^-_p}) \\
  \wedge S^{V^-_r + V^+_q- (V^-_q + V^+_r)} 
 \end{gathered}
 \ar[d] \\ 
F(S^{V_{\sigma_1} + V_{\sigma_2} }, (\bar{\hat{X}}_{\sigma_1} |
\hat{Z}_{\sigma_1} \wedge \bar{\hat{X}}_{\sigma_2} |
\hat{Z}_{\sigma_2})^{\mfib})  \wedge S^{V^-_r + V^+_q - (V^-_q +
  V^+_r)}. 
\end{tikzcd}
\end{equation}
Associativity follows from the fact that $A$ is a $2$-category.
\end{proof}

The strict action of $\Pi$ on $\bX$ implies the following result.  

\begin{prop}
For a Kuranishi flow category $\bX$, the virtual cochains form a
$\Pi$-equivariant spectrally enriched category with object set $\cP$
and morphism spectra $B\bar{\cX}|\cZ^{-V - d}(p,q)$ for objects $p$
and $q$. \qed
\end{prop}

Finally, given an associative ring spectrum $\bbk$ we produce a new
spectral category by smashing each
morphism spectrum with $\bbk$.

\begin{defin}\label{def:virtual-cochains}
The $\Pi$-equivariant spectral category $B\bar{\cX}|\cZ^{-V-d} \wedge
\bbk$ has object set $\cP$, with morphism spectra given by the 
$\bbk$-modules 
\begin{equation}
\left(B\bar{\cX}|\cZ^{-V-d} \sma \bbk\right)(p,q) = 
B\bar{\cX}|\cZ^{-V-d}(p,q) \sma \bbk.
\end{equation}
\end{defin}

\begin{rem}
It is tempting to think of $ B\bar{\cX}|\cZ^{-V-d} \wedge \bbk$ as a
category enriched in $\bbk$-modules, but this only makes sense if
$\bbk$ is commutative (at least $E_2$).  In the applications of
interest in this paper, $\bbk$ is typically only an associative ring
spectrum.
\end{rem}

\subsection{An augmentation of the virtual cochains}
\label{sec:an-augm-dirt}

The collapse map from a based space to $S^0$ induces an augmentation
(i.e., a map to the ground ring) on reduced homology.  Our goal in this
section is to obtain a similar construction on virtual cochains. 

Given a Kuranishi chart $\alpha$, the starting point is the map
$X_\alpha|Z_\alpha \to S^{V_\alpha}|0$ induced by the projection from
$X_\alpha$  to $S^{V_\alpha}$.  This map induces a map of $G$-spectra 
\begin{equation}
  F( S^{V_\alpha},  (X_\alpha | Z_\alpha)^{\mfib})   \to F( S^{V_\alpha}, (S^{V_\alpha}| 0)^{\mfib}).\end{equation}
Morally speaking, the augmentation arises by passing to the Borel
construction; $F(S^{V_\alpha}, (S^{V_\alpha}| 0)^{\mfib}) \htp \bS$
and $(EG_\alpha)_+ \sma_{G_\alpha} \bS \htp \Sigma^{\infty}_+
  BG_\alpha$, and so the collapse map $\Sigma^{\infty}_+ BG_\alpha \to
\bS$ now induces the augmentation.  However, in order to study the
multiplicative properties of the augmentation, we need to be precise
about the zig-zag.  We have the following diagram of comparison maps:
\begin{equation}
\begin{tikzcd}[column sep= small]
F( S^{V_\alpha}, (S^{V_\alpha}| 0)^{\mfib}) &
F( S^{V_\alpha}, (S^{V_\alpha})^{\mfib}) \ar[l,swap, "\htp"] &
F( S^{V_\alpha}, S^{V_\alpha}) \ar[l] &
\bS \ar[l] \ar[bend right=20, swap, "\htp", ll]
\end{tikzcd}
\end{equation}
Applying the Borel construction, we obtain a natural weak equivalence
\begin{equation}
(EG_{\alpha})_+ \sma_{G_\alpha} \bS \to (EG_{\alpha})_+
  \sma_{G_\alpha} F( S^{V_\alpha}, (S^{V_\alpha}| 0)^{\mfib}).
\end{equation}

Next, we have the chain of homeomorphisms
\begin{equation}
\begin{aligned}
(EG_{\alpha}) \sma_{G_\alpha} \bS = \Sigma^\infty_+ B(G_\alpha,
  G_\alpha, \ast) \sma_{G_\alpha} \bS &\cong 
B(\Sigma^\infty_+  G_\alpha,
  \Sigma^\infty_+ G_\alpha, \bS) \sma_{G_\alpha} \bS \\
&\cong
B(\bS, \Sigma^{\infty}_+ G_\alpha, \bS) \cong \Sigma^{\infty}_+
BG_\alpha.
\end{aligned}
\end{equation}

We now obtain a natural zig-zag 
\begin{equation}\label{eq:reduction_chart}
\begin{tikzcd}[column sep=small]
  B S^{V_\alpha}|0^{-V_\alpha} \equiv C_*(E G_\alpha;F( S^{V_\alpha},
  (S^{V_\alpha}| 0))^{\mfib})_{G_\alpha} & \ar[l,swap,"\htp"]
  \Sigma^{\infty}_+ BG_\alpha \ar[r] & \bS,
\end{tikzcd}  
\end{equation}
and composing with the map 
\begin{equation}
\begin{tikzcd}
B\bar{X}_\alpha|Z_\alpha^{-V_\alpha} \ar[r] & B
S^{V_\alpha}|0^{-V_\alpha}
\end{tikzcd}
\end{equation}
yields the \emph{augmentation associated to a chart}.

More generally, given a cube $\sigma$ of Kuranishi charts, we have a natural map
\begin{equation}  \label{eq:reduction_cube}
  B\bar{\hat{X}}_\sigma|\hat{Z}_\sigma^{-V_\sigma} \to B  S^{V_\sigma}|0^{-V_\sigma}.
\end{equation}
In order to discuss the functoriality of this map, we note that an inclusion $V_\alpha \to V_\beta$ induces a natural inclusion
\begin{equation} \label{eq:standard_spheres_act_on_cones}
  S^{V_f^\perp} \wedge \left( S^{V_\alpha}|0\right) \to   S^{V_\beta}|0, 
\end{equation}
arising from the fact that the left hand side is equivariantly
homeomorphic to $S^{V_\beta} | (S^{V_f} \times 0$).  Desuspending and
passing to the Borel construction then yields an equivalence  
\begin{equation}
B S^{V_\alpha}|0^{-V_\alpha}  \to   B S^{V_\beta}| 0^{-V_\beta}.
\end{equation}

This construction is multiplicative in the following sense.  Recall
from Equation~\eqref{eqn:cone-sphere-smash} that there is a natural
map
\begin{equation}
S^{V_1}|0 \sma S^{V_2}|0 \to S^{V_1 \oplus V_2}|0.
\end{equation}

\begin{prop}\label{prop:sphere-comparison}
Given a pair of inclusions $f_i \co V_{\alpha_i} \to V_{\beta_i}$, for
$i \in \{0,1\}$, we have an identification $V_{f_0 \times f_1}^\perp
\cong V_{f_0}^\perp  \oplus V_{f_1}^\perp$, which, together with
Equation \eqref{eq:standard_spheres_act_on_cones} gives a commutative
diagram  
\begin{equation}
  \begin{tikzcd}
    S^{V_{f_0}^\perp} \wedge \left( S^{V_{\alpha_0}}|0\right) \wedge  S^{V_{f_1}^\perp} \wedge \left( S^{V_{\alpha_1}}|0\right) \ar[r] \ar[d] & S^{V_{\beta_0}}|0 \wedge  S^{V_{\beta_1}}|0 \ar[d] \\
    S^{V_{f_0 \times f_1}^\perp} \wedge \left( S^{V_{\alpha_0} \oplus V_{\alpha_1}}|0\right) \ar[r] &  S^{V_{\beta_0} \oplus V_{\beta_1}}|0.   
  \end{tikzcd}
\end{equation}
Moreover, these identifications are associative for inclusions $f_i
\co V_{\alpha_i} \to V_{\beta_i}$ for $i \in \{0,1,2\}$.
\end{prop}

\begin{proof}
To see that the compatibility diagram commutes is a straightforward
check of the formulas.  The claim about associativity follows from the
commutative diagrams 
\begin{equation}
  \begin{tikzcd}
S^{V_{\beta_0}}|0 \wedge  S^{V_{\beta_1}}|0 \sma S^{V_{\beta_2}}|0
\ar[d] \ar[r] &  \left( S^{V_{\beta_0} \oplus V_{\beta_1}}|0 \right) \wedge
S^{V_{\beta_2}}|0 \ar[d] \\
S^{V_{\beta_0}}|0 \wedge  \left(S^{V_{\beta_1} \oplus V_{\beta_2}}|0\right) \ar[r] &
S^{V_{\beta_0} \oplus V_{\beta_1} \oplus V_{\beta_2}} \\
  \end{tikzcd}
\end{equation}
and the associativity of the smash product.
\end{proof}

Passing to desuspensions and applying the Borel construction, we
obtain a commutative diagram of spectra: 
\begin{equation}
  \begin{tikzcd}
B S^{V_{\alpha_0}}|0^{-V_{\alpha_0}} \wedge  B S^{V_{\alpha_1}}|0^{-V_{\alpha_1}} \ar[r] \ar[d] &   B S^{V_{\beta_0}}|0^{-V_{\beta_0}} \wedge  B S^{V_{\beta_1}}|0^{-V_{\beta_1}}   \ar[d] \\
    B S^{V_{\alpha_0} \oplus V_{\alpha_1}}|0^{-V_{\alpha_0} \oplus V_{\alpha_1}} \ar[r] & B S^{V_{\beta_0} \oplus V_{\beta_1}}|0^{-V_{\beta_0} \oplus V_{\beta_1} } .   
  \end{tikzcd}
\end{equation}
Proposition~\ref{prop:sphere-comparison} implies that these
compositions are associative as well.  Moreover, we observe that there
is a natural map $\bS \to B S^0|0$ induced by the homeomorphism $S^0 \to
S^0|0$, and so we can conclude the following proposition.

\begin{lem}
The assignment $\bX \mapsto B S^V|0^{-V}$ lifts to a lax monoidal functor
\begin{equation}
B S^V|0^{-V} \co  \Kur \to \Sp.
\end{equation} \qed
\end{lem}

Pulling this functor back to $\square \Kur$ under the evaluation map
given by the vertex $1^n$ of the $n$-dimensional cube, the following
result is straightforward to check: 

\begin{lem}\label{lem:augment}
Equation~\eqref{eq:reduction_cube} defines a lax monoidal
transformation
\begin{equation}
B\bar{X}_\sigma | Z_\sigma^{-V_\sigma} \to BS^V|0^{-V}
\end{equation}
of functors 
$\square \Kur \to \Sp$.
\qed
\end{lem}

We now turn to defining the augmentation of the virtual cochains of a
Kuranishi presentation.

\begin{defin}
Given a Kuranishi presentation $\bX \co A \to \Kur$, we define
\begin{equation}
BS^{V}|0^{-V}(A) \equiv \hocolim_{\alpha \in A} B S^{V_\alpha}| 0^{-V_\alpha}.
\end{equation}
\end{defin}

Lemma~\ref{lem:augment} has the following immediate corollary.

\begin{cor}\label{cor:augmentation}
There is an augmentation of
virtual cochains 
\begin{equation}
   B\bar{X} |Z^{-V}(A) \to B S^{V}|0^{-V}(A)
\end{equation}
associated to each Kuranishi presentation. \qed
\end{cor}

This comparison map is functorial in $A$ and multiplicative.

\begin{prop}\label{prop:multiplicative-augmentation}
Given a pair of Kuranishi presentations $\bX(i) \co A(i) \to
\Kur^{\cS_i}$ for $i \in \{1,2\}$, there is a map of homotopy
colimits  
\begin{equation}
B S^{V}|0^{-V}(A(1)) \wedge B S^{V}|0^{-V}(A(2)) \to B
S^{V}|0^{-V}(A(1) \times A(2))
\end{equation}
which fits into a commutative diagram
\begin{equation}
  \begin{tikzcd}
    B\bar{X}|Z^{-V}(A(1))  \wedge   B\bar{X}|Z^{-V}(A(2)) \ar[d]  \ar[r] &        B\bar{X}|Z^{-V}(A(1) \times A(2)) \ar[d] \\
   B S^{V}|0^{-V}(A(1)) \wedge B S^{V}|0^{-V}(A(2)) \ar[r] &  B S^{V}|0^{-V}(A(1) \times A(2)).
  \end{tikzcd}
\end{equation}
for Kuranishi presentations $\bX_1, \bX_2$, and $\bX_3$, the evident
associativity diagram commutes. \qed
\end{prop}

We apply the above discussion to a Kuranishi flow category $\bX$ with
object set $\cP$.  

\begin{defin}
For a Kuranishi flow category $\bX$, the $\Pi$-equivariant spectral
category $B S^{\cV}|0^{-\cV-d}$ has object set $\cP$ and morphism 
spectra 
\begin{equation}
S^{\cV}|0^{-\cV-d}(p,q) \equiv  \hocolim_{\alpha \in A(p,q)}  B
S^{V_\alpha}| 0^{-V_\alpha} \wedge S^{V^-_p + V^+_{q} - (V^+_{p}
+ V^-_{q})}[\deg p  - \deg q],
\end{equation}
in analogy with
Equation~\eqref{eq:morphisms_category_virtual_cochains}. The
compositions in this spectrally enriched category are induced by the
lax monoidal structure on the functor $B S^V|0^{-V}$, the multiplicative
structure on the spheres, and the action of $\Pi$ by the action on
$\cP$ and $A$.
\end{defin}

Corollary~\ref{cor:augmentation} and
Proposition~\ref{prop:multiplicative-augmentation} yield a
$\Pi$-equivariant spectral functor
\begin{equation}
B\bar{\cX}|\cZ^{-V-d} \to BS^{\cV}|0^{-\cV-d}.
\end{equation}

In order to analyse the outcome of this construction, we consider the map
\begin{equation}
  \bS \to  S^{V_\alpha}|0^{-V_\alpha} = F(S^{V_\alpha},(S^{V_\alpha}|0)^{\mfib})
\end{equation}
induced by the inclusion $ S^{V_\alpha} \subset S^{V_\alpha}| 0$.
This is an equivalence, and the diagram 
\begin{equation}
  \begin{tikzcd}
    \bS \ar[r] \ar[dr] &   S^{V_\alpha}|0^{-V_\alpha} \ar[d] \\
    & S^{V_\beta}|0^{-V_\beta}
  \end{tikzcd}
\end{equation}
commutes for each arrow of Kuranishi charts. This transformation is
compatible with products, so we conclude: 

\begin{lem} \label{lem:constant-to-Borel-on-canceling-spheres} 
There is a lax monoidal equivalence from the constant functor $\bS$ on
$\Kur$ to $S^{V}|0^{-V}$. \qed 
\end{lem}

We now construct a zig-zag of $\Pi$-equivariant spectral functors that
represent a homotopy class of functors
\begin{equation}
B S^{\cV}|0^{-\cV-d}  \to \Sp,   
\end{equation}
mapping $p$ to $S^{ V^+_{p} - V^-_p } [-\deg p]$.

We consider the $\Pi$-equivariant spectral category $\sZ$ with object
set $\cP$, and morphisms from $p$ to $q$ given by  
\begin{equation}
\sZ(p,q) = \hocolim_{\alpha \in A(p,q)} (BG_{\alpha})_+ \sma S^{V^-_p + V^+_{q} -
(V^+_{p} + V^-_{q}) } [ \deg p - \deg q]. 
 \end{equation}

Lemma \ref{lem:constant-to-Borel-on-canceling-spheres} supplies a
$\Pi$-equivariant spectral functor $\sZ \to B S^{\cV}|0^{-\cV-d}$.  On
the other hand, there is a $\Pi$-equivariant spectral functor from
$\sZ$ to the category of spectra induced by the projection map $B G_{\alpha} \to
\ast$ and the map
\begin{multline}
S^{V^-_p + V^+_{q} - (V^+_{p} + V^-_{q}) }  \wedge S^{ V^+_{p} - V^-_p } \cong   \\
S^{ V^+_{q} - V^-_{q} } \wedge  S^{V^-_p - V^+_{p} }  \wedge S^{ V^+_{p} - V^-_p }      \to  S^{ V^+_{q} - V^-_{q} }
\end{multline}
arising by permuting the spheres and using the natural map $S^V \sma
S^{-V} \to S^0$.  The resulting zig-zag 
\begin{equation}
\begin{tikzcd}
B\bar{\cX}|\cZ^{-V-d} \ar[r] & B S^{\cV}|0^{-\cV-d} & \ar[l] \sZ \ar[r] & \Sp 
\end{tikzcd}
\end{equation}
of $\Pi$-equivariant spectral functors provides our model for the
augmentation.

Using the bar construction to produce a formal composition rectifying
this zig-zag (as discussed in Section~\ref{sec:equifunc}), we conclude
the following proposition by smashing with a ring spectrum $\bbk$.

\begin{prop}
If $\bbk$ is a cofibrant associative ring spectrum, the augmentation
yields a $\Pi$-equivariant spectral functor 
\begin{equation}
B\bar{\cX}|\cZ^{-V-d} \wedge \bbk \to \bbk-\mod,   
\end{equation}
category, mapping $p$ to a spectrum which is equivalent to
$\bS^{V_{p}}[\deg p] \wedge \bbk$.  \qed
\end{prop}

\subsection{Signpost: basic examples}
\label{sec:signpost}

We end this section by outlining how a closed Kuranishi presentations
gives rise to a two-cell homotopy type.

First, let $\bX \co A \to \Kur$ be a closed Kuranishi presentation of
virtual dimension $d-1$. We shall construct a diagram 
\begin{equation} \label{eq:sign_post_Kuranishi-spaces-2cell}
 \Omega \bbk \to  C^*(B Z ; \Omega \bbk) \dashleftarrow \cdots \dashrightarrow  B X|Z^{-V-d} \wedge \bbk   \to S^V|0^{-V-d} \wedge \bbk   \cong \Omega^{d} \bbk
\end{equation}
where the first arrow comes from the unit, and the last arrow from
smashing the augmentation from Section \ref{sec:an-augm-dirt} with
$\bbk$.  The construction of the middle zig-zag will rely on the
assumption that (i) the charts $X_\alpha$ admit $G_\alpha$-equivariant
orientations with respect to $\bbk$, which are moreover compatible
with maps of charts, and (ii) the spectrum $\bbk$ is ambidextrous (see
Section \ref{sec:ambidextrous-spectra}). 

Inverting the arrows in the zig-zag, the outcome is a map  $\Omega
\bbk \to \Omega^{d}\bbk$. The \emph{homotopy type} associated to $\bX$
is the cofibre of the corresponding map. As a reality check, we note
that, if the map vanishes, then this cofibre is equivalent to the sum
$\bbk \oplus \Omega^{d} \bbk$. In Morse theory, whenever the index
difference between two critical points is $d$, then the moduli space
of flow lines between them is $d-1$-dimensional, which justifies this
shift. 

The construction of a homotopy type with more cells labelled by a
totally ordered set $\cP $ is associated to a diagram 
\begin{equation} \label{eq:sign_post_Kuranishi-spaces}
C^*_{\rel \partial}(\cP; \Omega \bbk) \to C^*_{\rel \partial}(B \cZ;
\Omega  \bbk) \dashleftarrow \cdots \dashrightarrow  B
\bar{\cX}|\cZ^{-V-d} \wedge \bbk \to \bbk-\mod, 
\end{equation}
of functors between spectrally enriched categories. The first arrow
goes back to the constructions of Section \ref{sec:relat-coch-coll},
and the last one is given by the multiplicativity of the
augmentation. The middle zig-zag will again follow from certain
orientability assumptions on the Kuranishi presentations that we
consider.

 \section{Comparison of cochains on a Kuranishi chart}
\label{sec:comp-coch-kuran}

The goal of this section is to relate compactly supported and virtual
cochains on a single Kuranishi chart. This will involve the use of
several geometric and homotopy-theoretic ideas: Spanier-Whitehead
duality, ambidexterity, the Adams isomorphism, flag smooth structures,
models of the tangent spherical fibration, and universal orientations.
We also begin to explain how this comparison extends to spectral
categories by discussing the compatibility of the maps we construct
with inclusions of boundary strata, and products of charts.  

To start, we recall that Spanier-Whitehead duality is a
generalization, at the level of stable homotopy types, of Poincar\'e
duality. Since our charts are not necessarily compact, and virtual
cochains in some sense correspond to homology, it is thus natural that
the first step in the comparison with ordinary cochains is the passage
to a compactly supported theory.

\subsection{Compactly supported cochains}
\label{sec:comp-supp-coch-4}

Given a orbispace chart $(Z,G)$, we let $Z^+$ denote the $1$-point
compactification of $Z$ regarded as based at the point at infinity,
and consider the space $BZ^+$ obtained as the reduced Borel
construction on $Z^+$, i.e., the quotient of the map $EG \times_G \ast
\to EG \times_G Z^+$ that includes the basepoint.

\begin{defin}[c.f. Equation~\eqref{eq:compactly-supported-maps}]
Given a vector space $E$, we define the compactly supported cochains
of $Z$ as the spectrum of maps 
\begin{equation}
  C^{*,c}( BZ; \Omega^E \bbk) = F(BZ^+, \Omega^E \bbk^{\mfib}).
\end{equation}
\end{defin}

Note that this is an abuse of terminology because there is no compact
support condition in the direction of the Borel construction, but only
with respect to the base $Z$.  The natural map $BZ \to BZ^+$ induces a map  
\begin{equation}\label{eq:compact-comparison}
   C^{*,c}(BZ; \Omega^E \bbk) \to   C^{*}(BZ; \Omega^E \bbk)
\end{equation}
which is an isomorphism whenever $Z$ is compact. We shall mostly be
interested in the case $E$ is a (trivialized) line which we denote
$\ell$.  

We will assume that $Z$ is an $\langle \cS \rangle$-stratified space
for a poset $\cS$, in the sense of
Definition~\ref{defn:stratified-space}.  Recall that this means that
we have a closed subset $\partial Z \subset Z$ which is a union of
closed subsets
\begin{equation} \label{eq:boundary_decomposition_orbichart}
   \partial Z = \bigcup_{q \in \cS} \partial^q Z 
\end{equation}
with $\partial^q Z \cap \partial^pZ = \emptyset$ whenever $q$ and $p$
are not comparable.  We denote by $\hat{Z}$ the collared completion of
$Z$ as in Definition~\ref{defn:collared-completion} (see also
Lemma~\ref{lem:collar-functor}).

We now define the relative compactly supported cochains to be the
spectrum of cochains on $\hat{Z}$ which vanish at the
boundary; note that the collaring function ensures that the inclusion
of $B\partial\hat{Z}^+$ into $B\hat{Z}^+$ is a Hurewicz cofibration.

\begin{defin}[c.f. Equation~\eqref{eq:relative-cochains}]\label{defin:relative-compactly-supported}
The relative compactly supported cochains of an $\langle \cS
\rangle$-stratified space are given by the spectrum
\begin{equation} \label{eq:relative_cochains_stratified_orbichart}
   C^{*,c}_{\rel \partial}(B\hat{Z}; \Omega^\ell \bbk) \equiv
   C^{*,c}(B\hat{Z}, B\partial \hat{Z}; \Omega^\ell \bbk) =
   F(B\hat{Z}^+/B\partial\hat{Z}^+, \Omega^\ell \bbk^{\mfib}).
\end{equation}
\end{defin}

The map in Equation~\eqref{eq:compact-comparison} extends to a natural
comparison map 
\begin{equation}
C^{*,c}_{\rel \partial}(B\hat{Z}; \Omega^{\ell}\bbk) \to C^{*}_{\rel \partial}(B\hat{Z}; \Omega^{\ell} \bbk).
\end{equation}
This comparison map is multiplicative, since the inclusion $Z \to Z^+$
is compatible with the homeomorphism $Z_1^+ \sma Z_2^+ \cong (Z_1 \times
Z_2)^+$.

\begin{lem}\label{lem:compactly-supported-product}
Given a pair $(Z_{1},G_{1})$ and $(Z_{2},G_{2})$ of orbispace charts,
and real lines $\ell_1$ and $\ell_2$, we have a commutative diagram
\begin{equation}
  \begin{tikzcd}
    C^{*,c}_{\rel \partial}(B \hat{Z}_{1}; \Omega^{\ell_1} \bbk)
    \wedge C^{*,c}_{\rel \partial}(B \hat{Z}_{2}; \Omega^{\ell_2}
    \bbk) \ar[r] \ar[d] & C^{*,c}_{\rel \partial}(B (\hat{Z}_{1}
    \times \hat{Z}_{2}); \Omega^{\ell_1 \oplus \ell_2} \bbk) \ar[d]
    \\ C^{*}_{\rel \partial}(B \hat{Z}_{1}; \Omega^{\ell_1} \bbk)
    \wedge C^{*}_{\rel \partial}(B \hat{Z}_{2}; \Omega^{\ell_2} \bbk)
    \ar[r] & C^{*}_{\rel \partial}(B (\hat{Z}_{1} \times \hat{Z}_{2});
    \Omega^{\ell_1 \oplus \ell_2} \bbk).
    \end{tikzcd}
\end{equation}
The evident associativity diagrams strictly commute. \qed
\end{lem}

For a codimension $1$-stratum associated to an element $q \in Q$, we
have a \emph{boundary map} 
\begin{equation}
C^{*,c}_{\rel \partial}(B\partial^q \hat{Z}; \Omega^{\kappa \oplus
      \ell} \bbk) \to 
C^{*,c}_{\rel \partial}(B\hat{Z}; \Omega^{\ell}\bbk),
\end{equation}
where $\kappa$ is the real line associated to the collar.
The boundary map fits into the commutative diagram 
\begin{equation}\label{eqn:compactly-supported-boundary-map}
  \begin{tikzcd}
    C^{*,c}_{\rel \partial}(B\partial^q \hat{Z}; \Omega^{\kappa \oplus
      \ell} \bbk) \ar[r] \ar[d] &   C^{*}_{\rel \partial}(B\partial^q
    \hat{Z}; \Omega^{\kappa \oplus \ell}  \bbk) \ar[d] \\
     C^{*,c}_{\rel \partial}(B\hat{Z}; \Omega^{\ell}\bbk) \ar[r] &
     C^{*}_{\rel \partial}(B\hat{Z}; \Omega^{\ell} \bbk).
  \end{tikzcd}
\end{equation}

\subsection{Spanier-Whitehead duality}
\label{sec:span-whit-dual-5}

The main goal of this subsection is to construct what we call the
\emph{Milnor model} for the tangential twist of (compactly supported
relative) cochains associated to Kuranishi charts, and compare them,
via Spanier-Whitehead duality, to the virtual cochains.  This will
require a choice of coefficients for which the classifying spaces of
finite groups are dualizable (i.e., behave as Poincar\'e duality
spaces); it is at this point in our work that restriction to the
Morava $K$-theories becomes relevant.  We conclude the section with
the use of the Adams isomorphism and the norm map from homotopy
orbits to homotopy fixed-points to pass from chains to cochains.

\subsubsection{Spanier-Whitehead duality for manifolds}
\label{sec:span-whit-dual}

Let $X$ be a topological manifold.  We regard $X \times X$ as a space
over $X$ via the projection $\pi \colon X \times X \to X$ to the first
factor.

\begin{defin}
Let $\Milnor X$ be the spherical fibre bundle over $X$ obtained from
$X \times X$ by taking the fiberwise cone of $(X \times X \setminus
\Delta)$, where $\Delta$ denotes the diagonal.
\end{defin}

As a set, this is  
\begin{equation}
  \Milnor X \equiv \coprod_{x \in X} X|x  
\end{equation}
where $X|x = C(X, X \setminus x)$ as in Section
\ref{sec:first-non-trivial-example-cochains}. Since the homeomorphism
group of connected manifolds acts transitively, the first coordinate
projection map $\Milnor X \to X$ is a fibre bundle. There is a
canonical section 
\begin{equation}
  X \to   \Milnor X
\end{equation}
given by the cone point in each fibre, and the image of this section
is contained in a contractible subbundle 
\begin{equation}
  \Milnor_0 X \equiv   \coprod_{x \in X} C(X \setminus x) 
\end{equation}
given by the union of all the cones.

\begin{rem}
We call $\Milnor X$ the \emph{Milnor model} of the tangent spherical
fibration of $X$, because of its close relationship with the tangent
microbundle introduced by Milnor in~\cite{Milnor1964}, and our
inability to find a result in the literature that uses this
construction preceding Milnor's work. 
\end{rem}

Let $i \colon Z \subset X$ be a closed subset.  We shall abuse
notation and write $\Milnor X$ for the restriction of the Milnor model
to $Z$, i.e., the pullback along the inclusion $i$.  We now explain a 
construction of Atiyah-Spanier-Whitehead duality as an equivalence 
\begin{equation} \label{eq:SW-without-boundary}
\Sigma^{\infty} X | Z = \Sigma^{\infty} C(X, X \setminus Z) \to
C^{*,c_0}(Z ;  \Milnor X ), 
\end{equation}
where the righthand side is the spectrum of compactly
supported sections relative to $\Milnor_0 X$.

Recall from Section~\ref{sec:parametrized-spectra} that the compactly
supported sections are defined in terms of the spectrum $C^{*}(Z ;
\Milnor X) = \Gamma_Z(\Milnor X)$ of sections of the parametrized
suspension spectrum of the Milnor fibration.  Specifically, in  
Definitions~\ref{defin:space-of-sections}
and~\ref{defin:spectrum-of-sections}, we define the spectrum of compactly
supported sections of $\Milnor X$ relative to $\Milnor_0 X$ as the
levelwise subspectrum 
\begin{equation}
C^{*,c_0}(Z; \Milnor X) = \Gamma^c_Z(\Milnor X, \Milnor_0 X)
\end{equation}
of $C^{*}(Z ; \Milnor X)$ consisting of sections whose values, outside
a compact subset of $Z$, lie in $\Milnor_0 X$.  There is an evident
inclusion map 
\begin{equation}
C^{*,c}(Z ; \Milnor X \wedge \bbk) \to C^{*,c_0}(Z ; \Milnor X \wedge
\bbk ). 
\end{equation}

The duality map is induced by the space-level map
\begin{equation} \label{eq:constant_sections}
X | Z \to \Map_Z(Z; \Milnor X) 
\end{equation}
that includes $X|Z$ as the ``constant'' sections, where here
$\Map_Z(Z; \Milnor X) $ denotes the space of sections.  More precisely,
we start with the map 
\begin{equation}
X \to \Map_Z(Z; \Milnor X)
\end{equation}
which assigns to each point $x \in X$ the section of $\Milnor X$ that
maps every point $z \in Z$ to $x$ in the fiber.  These sections are
clearly compactly supported relative to $\Milnor_0 X$.  Moreover, if $x$ lies in
$X \setminus Z$, this section lands in the cone part for every $z \in
Z$, so that this map extends tautologically to a map
\begin{equation}\label{eq:SW-duality-CS}
X|Z \to \Gamma^c_{Z}(\Milnor X, \Milnor_0 X),
\end{equation}
where the target is the space of compactly supported sections of
$\Milnor X$ relative to $\Milnor_0 X$ over $Z$.

\begin{prop} \label{prop:comparison_models_compactly_supported_Milnor}
Let $X$ be a topological manifold and $Z \subset X$ a closed subset.
If $\bbk$ is a cofibrant spectrum, the Spanier-Whitehead duality map
of Equation~\eqref{eq:SW-duality-CS} and the inclusion induce natural
equivalences:
\begin{equation} \label{eq:comparison_compactly_supported_Milnor-without-boundary}
\Sigma^\infty X|Z \wedge \bbk \to  C^{*,c_0}(Z ; \Milnor X \wedge \bbk
) \leftarrow  C^{*,c}(Z ; \Milnor X \wedge \bbk).
\end{equation}
\end{prop}

\begin{proof}
The proof proceeds by inducting on the number of elements of a
contractible cover and uses the continuity properties of these
invariants.  When $Z$ is a contractible closed set, the theorem is
tautologically true.  Since each of the functors in the zig-zag
satisfies the Mayor-Vietoris property, we can conclude that the
theorem holds for closed subsets which admit a locally finite
contractible cover.  Every closed subset of a manifold is an
intersection of such sets, and the result now follows since each of
the terms in the zig-zag represent cohomology theories (as functors of
closed subsets of $X$) that are continuous in the sense of~\cite[\S
  2]{Spanier87}, i.e., have the property that the colimit over closed
neighborhoods of $A$ is equal to the value at $A$.  For $\Sigma^\infty
X|Z$, continuity is a consequence of the fact that $X|Z_i$ becomes a
filtered system and the homotopy cofiber commutes with filtered
homotopy colimits.  For the spectra of sections, this follows as
in~\cite[Thm. 5]{LeeRaymond1968}; see also
Section~\ref{sec:parametrized-spectra}.
\end{proof}

We shall refer to the middle and right hand side of Equation
\eqref{eq:comparison_compactly_supported_Milnor-without-boundary} as
the two models for the \emph{spectrum of compactly supported sections}
with coefficients in $\bbk$.

\begin{rem}
In the classical point of view on Poincar\'e duality, the left hand
side is the homology of the complement of $Z$ (which we call virtual
cochains of $Z$ following Pardon), and the right hand side is the
(twisted) cohomology of $Z$.
\end{rem}

The next result asserts that the construction of the Milnor fibration
is compatible with products.  For a pair of sectioned parametrized
spaces (refered to as ex-spaces) $f_1 \colon E_1 \to X_1$ and $f_2
\colon E_2 \to X_2$, there is an external smash product $E_1 \barwedge
E_2$ over $X_1 \times X_2$ that on fibers is the smash product, and which we review in Definition \ref{def:external_smash_para}.

\begin{lem}
If $X_{1}$ and $X_{2}$ are manifolds, there is a natural map
\begin{equation} \label{eq:product_map_milnor}
  \Milnor X_{1} \barwedge \Milnor X_{2} \to \Milnor ( X_{1} \times
  X_{2} )
\end{equation}
of sectioned spaces over $X_{1} \times X_{2}$ induced by the product
on homotopy cofibers 
\begin{equation}
C(X_1, X_1\setminus\{x_1\}) \sma C(X_2, X_2\setminus \{x_2\}) \to
C(X_1 \times X_2, (X_1 \times X_2) \setminus \{(x_1,x_2)\})
\end{equation}
as in Equation~\eqref{eqn:cone-difference-prod}. This product map is
associative in the sense that the diagram of ex-spaces over $X_1
\times X_2 \times X_3$
\begin{equation}
\begin{tikzcd}
\Milnor X_{1} \barwedge \Milnor X_2 \barwedge \Milnor X_3 \ar[r]
\ar[d] & \Milnor X_1 \barwedge \Milnor (X_2 \times X_3) \ar[d] \\
\Milnor (X_1 \times X_2) \barwedge \Milnor X_3 \ar[r] & \Milnor (X_1
\times X_2 \times X_3) 
\end{tikzcd}
\end{equation}
commutes. \qed
\end{lem}

We now state the compatibility of our models for the Spanier-Whitehead
duality maps with the product structure. For the formulation of the
next result, we use the natural maps
\begin{equation}
\begin{tikzcd}
C(X_1, X_1 \setminus Z_1) \sma C(X_2, X_2 \setminus Z_2) \ar[d] \\
C(X_1 \times X_2, (X_1 \times X_2) \setminus (Z_1 \times Z_2))
\end{tikzcd}
\end{equation}
of Equation~\eqref{eqn:cone-difference-prod} throughout the 
diagram and note that the product maps on spectra of compactly
supported cochains are defined using the monoidal structure on spaces of
sections (see
Proposition~\ref{prop:multiplicative-compact-sections-spectra}) and
the multiplicativity of the fibrant replacement functor.
   
\begin{lem}\label{lem:SW-map-multiplicative}
The product map induces a commutative diagram
\begin{equation} \label{eq:multiplicativity_SW_maps}
  \begin{tikzcd}[column sep=small]
    \Sigma^\infty X|Z_{1} \wedge \Sigma^\infty X|Z_{2}  \ar[r] \ar[d] & \Sigma^\infty \left(X_{1} \times X_{2} |Z_{1} \times Z_{2}\right) \ar[d] \\
C^{*,c_0}(Z_{1} ; \Milnor X_{1}) \wedge C^{*,c_0}(Z_{2} ; \Milnor X_{2})\ar[r] & C^{*,c_0}(Z_{1} \times Z_{2} ; \Milnor \left(X_{1} \times X_{2}\right)) \\
C^{*,c}(Z_{1} ; \Milnor X_{1}) \wedge C^{*,c}(Z_{2} ; \Milnor X_{2}) \ar[u] \ar[r]&   \ar[u] C^{*,c}(Z_{1} \times Z_{2} ; \Milnor \left(X_{1} \times X_{2}\right))
  \end{tikzcd}
\end{equation}
for any pair $Z_{1}$ and $Z_{2}$ of closed subsets and analogous
commuting associativity diagrams.  Finally, for any associative ring
spectrum $\bbk$ we have a corresponding diagram with coefficients in
$\bbk$. \qed 
\end{lem}

\begin{rem}
Note the fact that $\Milnor X$ is defined without assuming that $X$ is
a manifold, and that the Spanier-Whitehead duality map is also defined
without this assumption.  The manifold condition only enters in the
proof that the map is an equivalence, and in fact that condition is
only required near $Z$. We shall use this flexibility in later
constructions.
\end{rem}

\begin{rem} \label{rem:lack-of-functoriality-Milnor}
A fundamental disadvantage of the model $\Milnor X$ is the absence of
adequate maps associated to covering maps $\tilde{X} \to X$. The other
models of the tangent spherical fibration which we will later use are
equipped with a canonical map from the pullback of the fibration
associated to $X$ to the fibration associated to $\tilde{X}$. For the
Milnor model, such a map only exists if the covering space is
inessential in the sense that it is a product of $X$ with a discrete
set. We shall return to this issue in Section
\ref{sec:inessential-charts}. 
\end{rem}

\subsubsection{Spanier-Whitehead duality for manifolds with corners}
\label{sec:span-whit-dual-2}

For a manifold $X$ with boundary, the inclusion map in Equation
\eqref{eq:constant_sections} is still well-defined, but is
unfortunately not well-adapted to encoding the relationship between
duality for $Z$ and for its intersection with the boundary. Instead,
we shall work with cochains relative the boundary, after attaching a
collar.  This allows us both to avoid questions of cofibrancy of the
inclusion of the boundary and to be able to define boundary
homomorphisms at the level of cochains.

We work in the setting of $\langle \cS \rangle$-manifolds as in
Section \ref{sec:langle-cs-rangle}: applying the construction of
Section \ref{sec:topol-enrich-categ} in this context, we obtain a
collared space $\hat{X}$ which is covered by the union of the product
of strata $\partial^Q X$ with the corresponding cube collar
$\kappa^Q$, which we find convenient to identify with $[0,1]^Q$ for
the beginning of this section. 
\begin{lem}
If $X$ is an $\langle \cS \rangle$-manifold with boundary, so is
$\hat{X} $, and there is a natural homeomorphism 
\begin{equation}
\partial^Q\hat{X} \cong \widehat{\partial^QX}.
\end{equation}
\end{lem}
\begin{proof}
The proof that $\hat{X}$ is a manifold is local: the assumption that
$X$ is an $\langle \cS \rangle$-manifold amounts to the assertion
that the stratification near $\partial^Q X$ is given by the product
of $\partial^Q X$ with $(-\infty,0]^Q$.  We conclude that a
  neighbourhood of $\partial^Q X \times \kappa^Q$ in $\hat{X}$ is
  homeomorphic to $\partial^Q X \times (-\infty,1]^Q$. 
\end{proof}

We have a version of the Milnor model in the stratified setting.

\begin{defin}
We define the \emph{completed Milnor spherical fibration}
\begin{equation}
 \Milnor \hat{X} \to \hat{X}   
\end{equation}
as the restriction of the Milnor tangent bundle of the topological
manifold without boundary obtained as the union of $X$ with infinite
collars $ \partial^Q X \times [0,\infty)^Q$ to the collared manifold
  $\hat{X}$. 
\end{defin}

Given a closed subset $Z \subset X$, which we equip with the
stratification inherited from $X$, we obtain an inclusion $\hat{Z}
\subset \hat{X}$. As before, we also abusively write $\Milnor \hat{X}
$ for the restriction of the completed Milnor spherical fibration to
$\hat{Z}$, which inherits a stratification from $\hat{X}$.  For each
totally ordered subset $Q$ of $\cS$, we obtain an inclusion of pairs
\begin{equation} \label{eq:inclusion_pairs_closed_subset_manifold_boundary}
(\partial^Q \hat{X}, \partial^Q \hat{Z}) \subset (\hat{X}, \hat{Z})  .
\end{equation}

The restriction of the Milnor spherical bundle to the boundary admits
a subfibration 
\begin{equation}
  \Milnor^{\mathrm{in}} \hat{X} \subset \Milnor \hat{X}  
\end{equation}
over the boundary, consisting of ``inward pointing'' vectors, i.e.
points lying in $\hat{X}$ rather than in the infinite completion.  
\begin{lem}
The fibre of $\Milnor^{\mathrm{in}} \hat{X}$ over every (boundary)
point is contractible.   
\end{lem}
\begin{proof}
  Since the result is local, it suffices to show that the fibre of $\Milnor^{\mathrm{in}} (\bR^{n-k} \times  (0, 1]^k)$ at a boundary point is contractible at the deepest stratum.  There is a natural equivalence
  \begin{equation}
        \Milnor_x \bR^{n-k} \wedge  \Milnor^{\mathrm{in}}_1 (0,\infty)^{\wedge k}  \to \Milnor^{\mathrm{in}}_{(x,1^k)} (\bR^{n-k} \times  (0, 1]^k)
  \end{equation}
  induced by the product in Equation \eqref{eq:product_map_milnor}.  The space $   \Milnor^{\mathrm{in}}_{1} (0,\infty)$ is contractible, so the result follows.
\end{proof}

At this stage, we introduce the compactly supported relative cochains.

\begin{defin}
The compactly supported relative cochains is the spectrum
\begin{equation}
C^{*,c_0}_{\rel \partial}(\hat{Z} ; \Milnor \hat{X} ) \subset
C^{*}(\hat{Z} ; \Milnor \hat{X} )
\end{equation}
consisting of sections whose values lie in $\Milnor_0 \hat{X} $
outside a compact set, and in $\Milnor^{\mathrm{in}} \hat{X}$ at the
boundary.
\end{defin}

As in the previous section, this is ultimately a model for the
compactly supported cochains relative the boundary: 
\begin{lem}
  The inclusion of sections with value the base section outside a compact set and on the boundary induces an equivalence
  \begin{equation}
    C^{*,c}_{\rel \partial}(\hat{Z} ;  \Milnor \hat{X} )   \to C^{*,c_0}_{\rel \partial}(\hat{Z} ;  \Milnor \hat{X})
  \end{equation}
\end{lem}
\begin{proof}
  This follows from the contractibility of the fibres of $\Milnor^{\mathrm{in}} \hat{X}$.
\end{proof}

Our goal now is to construct a commutative diagram
\begin{equation} \label{eq:boundary_map_commute_in_chain_and_cochains}
  \begin{tikzcd}
\Sigma^\infty \partial^Q \hat{X}| \partial^Q \hat{Z}  \ar[r] \ar[d] & C^{*,c_0}(\partial^Q \hat{Z}; \Milnor  \partial^Q \hat{X})  \ar[d, dashed] \\
\Sigma^\infty \hat{X}| \hat{Z}  \ar[r, dashed] & C^{*,c_0}_{\rel \partial}(\hat{Z}; \Milnor \hat{X}) 
 \end{tikzcd}
\end{equation}
in which the solid horizontal map is defined as in the previous section and the solid vertical map is the map induced by the inclusion
of pairs in Equation \eqref{eq:inclusion_pairs_closed_subset_manifold_boundary}.

It is convenient at this stage to write $\Milnor \kappa^Q \to \kappa^Q$ for the space whose fibre at a point $t \in \kappa$ is the cone of the inclusion of the complement of $t$ in $\kappa$. Note that this definition involves some abuse of terminology as the fibres over points lying in the boundary are contractible because we do not complete to a manifold. With this model, there is a natural diagram
\begin{equation} \label{eq:collared_fibration_Milnor_tangent_manifold_boundary}
  \begin{tikzcd}
    \Milnor  \partial^Q \hat{X} \barwedge  \Milnor \kappa^Q \ar[r] \ar[d] & \Milnor \hat{X} \ar[d] \\
    \partial^Q \hat{X} \times \kappa^Q \ar[r] & \hat{X}
  \end{tikzcd}
\end{equation}
where the top horizontal map is induced by the product $\partial^Q \hat{X} \times \kappa^Q \to \hat{X}$ after passing to the product, and taking the cone of the complement of the diagonal.

We now fix the section $t \mapsto (1, |1-t|)$ of the Milnor fibration
of the one dimensional cube, which induces, by taking the product and
the maximum of collar coordinates, a section of $\Milnor \kappa^Q
$. By definition, the value of this section along the boundary facets
where any coordinate vanishes agrees with the cone point, and the
value at the corner $1^Q$ lies in the inner pointing part. This
induces a map of section spaces 
\begin{multline} \label{eq:map_sections_boundary_to_interior_Milnor_spaces}
  \Map_{\partial^Q \hat{Z}}(\partial^Q \hat{Z}; \Milnor \partial^Q \hat{X}) \to \\\Map_{\partial^Q \hat{Z}}(\partial^Q \hat{Z} \times \kappa^Q; \Milnor \partial^Q \hat{X} \wedge \Milnor \kappa^Q)
 \to  \Map_{ \hat{Z}}(\hat{Z};   \Milnor \hat{X}) 
\end{multline}
where the first map is the product with the given section, and the
second is the composition of the top horizontal map in Equation
\eqref{eq:collared_fibration_Milnor_tangent_manifold_boundary} with
the extension away from the collar by the constant section with value
the cone point. These maps preserve the condition that a compact
subset be mapped to the cone part of the Milnor model. Passing to
spectra, we obtain the right vertical map in Diagram
\eqref{eq:boundary_map_commute_in_chain_and_cochains}.  It thus remains
to define the bottom horizontal map. 

The most obvious such map arises from Equation
\eqref{eq:constant_sections}, noting the fact the image of such
constant sections takes the desired values at the boundary. However,
this choice would not yield a commutative diagram, so we change it in
the collar to obtain a map 
\begin{equation}
\hat{X} | \hat{Z} \to \Map_{\hat{Z}}(\hat{Z} ; \Milnor \hat{X})
\end{equation}
which assigns to a point $x$ in the collar with collar coordinate $\{t_q(x)\}_{q \in Q}$ the section
\begin{equation}  \label{eq:constant_sections_modified_in_collar}
z \mapsto   \begin{cases}
  (x, \max_{q} t_q(x)) & \textrm{ if } z \in Z \\
  (x, \max_{q} |t_q(x) - t_q(z)|)& \textrm{ if } z \in \partial^Q Z \times \kappa^Q,
  \end{cases}
\end{equation}
and which we extend to the cone direction in $ \hat{X} | \hat{Z} $ by
using the maximum of the cone coordinates.  This is clearly
continuous.

The diagram
\begin{equation}
  \begin{tikzcd}
\partial^Q \hat{X}| \partial^Q \hat{Z}  \ar[r]    \ar[d] &   \Map_{\partial^Q \hat{Z}}(\partial^Q \hat{Z}; \Milnor \partial^Q \hat{X}) \ar[d] \\
\hat{X} | \hat{Z} \ar[r] &  \Map_{ \hat{Z}}(\hat{Z};   \Milnor \hat{X}) 
  \end{tikzcd}
\end{equation}
then commutes by construction, which yields Diagram
\eqref{eq:boundary_map_commute_in_chain_and_cochains} after passing to
spectra.

We now state the analogue of Proposition
\ref{prop:comparison_models_compactly_supported_Milnor}, whose proof
is entirely similar: 
\begin{prop}
The Spanier-Whitehead maps, and the inclusion of the base-section
induce equivalences 
\begin{equation}
\Sigma^\infty \hat{X}| \hat{Z} \to  C^{*,c_0}_{\rel \partial}(\hat{Z}; \Milnor \hat{X})  \leftarrow C^{*,c}_{\rel \partial}(\hat{Z}; \Milnor \hat{X}). 
\end{equation} \qed
\end{prop}

We note a final infelicity of our construction: while we have arranged
for Diagram \eqref{eq:boundary_map_commute_in_chain_and_cochains} to
commute, there seems to be no natural way to construct a boundary map 
\begin{equation}
 C^{*,c}(\partial^Q \hat{Z}; \Milnor \partial^Q \hat{X})  \to  C^{*,c}_{\rel \partial}(\hat{Z}; \Milnor \hat{X})
\end{equation}
in such a way that the comparison maps above give a commutative diagram. The problem is that the section of $\Milnor \kappa$ used to construct Equation \eqref{eq:map_sections_boundary_to_interior_Milnor_spaces} does not vanish at the boundary. However, it is straightforward to see that the map
\begin{equation} \label{eq:section_Milnor_collar_vanish_boundary}
 t \mapsto (1/2, 2|1/2 - t|) 
\end{equation}
does satisfy this property. 

In order to compare these two constructions, we introduce an
intermediate spectrum of sections constructed using
Equation~\eqref{eq:section_Milnor_collar_vanish_boundary}, which we
denote by $C^{*,c'_0}_{\rel \partial}(\hat{Z}; \Milnor \hat{X})$.
There is then a zig-zag of comparisons
\begin{equation} \label{eq:interpolate_notions_compact_cochains}
 C^{*,c}_{\rel \partial}(\hat{Z}; \Milnor \hat{X}) \leftarrow C^{*,c'_0}_{\rel \partial}(\hat{Z}; \Milnor \hat{X}) \to  C^{*,c_0}_{\rel \partial}(\hat{Z}; \Milnor \hat{X}),
\end{equation}
where the left pointing arrow is the inclusion of sections which are
inward pointing at $1/2$ and take value in the cone region over
$[1/2,1]$ and the right pointing arrow is the composition of the
restriction to the complement of $[1/2,1] \times \partial \hat{Z}$,
together with the identification of this complement with $\hat{Z}$
given by rescaling the collar coordinate. 

\subsubsection{Spanier Whitehead duality for Kuranishi charts}
\label{sec:span-whit-dual-4}

Let $X$ be a $G$-equivariant $\langle \cS \rangle$-manifold, and let
$Z \subset X$ be a closed $G$-invariant subset. We assume throughout
that the stratification is preserved by $G$.  The constructions of the
previous sections only involved the choice of a section of the
fibration associated to the collar $\kappa$, and hence are clearly
$G$-equivariant.  We conclude: 

\begin{prop}\label{prop:SW-start}
There is a commutative diagram of $G$-spectra
\begin{equation} \label{eq:SW-equivalence-Gamma-spectra}
\begin{tikzcd}
  \Sigma^\infty \partial^Q \hat{X}|\partial^Q \hat{Z} \ar[r] \ar[d] & \Sigma^\infty \hat{X}| \hat{Z}  \ar[d] \\
  C^{*,c_0}_{\rel \partial}(\partial^Q \hat{Z}; \Milnor \partial^Q \hat{X})   \ar[r] & C^{*,c_0}_{\rel
    \partial}(\hat{Z}; \Milnor \hat{X}) \\
  C^{*,c'_0}_{\rel \partial}(\partial^Q \hat{Z}; \Milnor \partial^Q \hat{X})   \ar[d] \ar[u] \ar[r] & C^{*,c'_0}_{\rel
  \partial}(\hat{Z}; \Milnor \hat{X}) \ar[u] \ar[d] \\
  C^{*,c}_{\rel \partial}(\partial^Q \hat{Z}; \Milnor \partial^Q \hat{X})  \ar[r]& C^{*,c}_{\rel
  \partial}(\hat{Z}; \Milnor \hat{X}) 
\end{tikzcd}
\end{equation}
in which all vertical arrows are equivalences. \qed
\end{prop}

We now will pass to homotopy fixed points for the action of $G$.  For
convenience, we introduce the following notation.

\begin{notn}
We will write 
\begin{align}
C^*(BG; Y) &= C^*(EG; Y)^{G} \\
C^{*,c_0}_{\rel \partial}(B \hat{Z}; \Milnor \hat{X}) & = C^*(EG;
C^{*,c_0}_{\rel \partial}(\hat{Z}; \Milnor \hat{X}))^G
\\ C^{*,c}_{\rel \partial}(B \hat{Z}; \Milnor \hat{X}) & = C^*(EG;
C^{*,c}_{\rel \partial}(\hat{Z}; \Milnor \hat{X}))^G
\end{align}
for the homotopy fixed-point spectra.

As explained in Appendix~\ref{sec:parametrized-spectra} (see
Proposition~\ref{prop:sections-BG}), for a $G$-spectrum $Y$ the
spectrum $C^*(BG ; Y) $ is equivalent to the spectrum of sections of
a parametrized spectrum over $BG$ induced by the Borel construction.  
\end{notn}

Applying the functor $C^*(EG;-)^{G}$ to
Diagram~\eqref{eq:SW-equivalence-Gamma-spectra} in
Proposition~\ref{prop:SW-start}, we obtain the diagram of spectra: 
\begin{equation} \label{eq:SW-milnor-Borel}
  \begin{tikzcd}
    C^*(B G, \partial^Q \hat{X}|\partial^Q \hat{Z}) \ar[r] \ar[d] &  C^*(B G,\hat{X}|\hat{Z})  \ar[d]\\
    C^{*,c_0}_{\rel \partial}(B \partial^Q \hat{Z}; \Milnor \partial^Q \hat{X}) \ar[r] & C^{*,c_0}_{\rel \partial}(B\hat{Z}; \Milnor \hat{X}) \\
    C^{*,c'_0}_{\rel \partial}(B\partial^Q \hat{Z}; \Milnor \partial^Q \hat{X}) \ar[u]  \ar[r] \ar[d] &  C^{*,c'_0}_{\rel \partial}(B \hat{Z}; \Milnor \hat{X}) \ar[u]  \ar[d] \\
    C^{*,c}_{\rel \partial}(B\partial^Q \hat{Z}; \Milnor \partial^Q \hat{X})  \ar[r]&  C^{*,c}_{\rel \partial}(B \hat{Z}; \Milnor \hat{X}) 
 \end{tikzcd}
\end{equation}
in which all horizontal maps are equivalences.

The homeomorphism $E G_{1} \times E G_{2} \cong E G_{1} \times E
G_{2}$ and the external multiplicative structure on fixed points
induce an associative pairing
\begin{equation}
C^*(EG_1; Y_1) \sma C^*(EG_2; Y_2) \to C^*(E(G_1 \times G_2); Y_1
\sma Y_2).
\end{equation}

\begin{lem}
Given $\langle \cS_i \rangle$-manifolds $X_i$ for $i \in \{1,2\}$, and
closed inclusions $Z_i \subset X_i$ of $G_i$-invariant subsets, 
the map $\hat{X}|\hat{Z}_{1} \wedge \hat{X}|\hat{Z}_{2} \to
\hat{X}_{1} \times \hat{X}_{2} |\hat{Z}_{1}  \times \hat{Z}_{2}$
induces a map 
\begin{equation}
   C^*(B G_{1}, \hat{X}|\hat{Z}_{1}) \wedge  C^*(B G_{2}, \hat{X}|\hat{Z}_{2}) \to C^*(B ( G_{1} \times G_{2}) ,  \hat{X}_{1} \times \hat{X}_{2} |\hat{Z}_{1}  \times \hat{Z}_{2}).
\end{equation}
These maps satisfy the evident associativity diagram for a triple of pairs $(X_i, Z_i)$. \qed
\end{lem}

If the space we consider arises from a Kuranishi chart, it is natural
at this stage to also desuspend by the $G$-representation $V$.

\begin{defin}
For a $G$-representation $V$, we define 
\begin{equation}
C^{*,c}(B \hat{Z}; \Milnor  \hat{X}^{-V}) \equiv C^*(EG,F(S^V,
C^{*,c}(\hat{Z}; \Milnor  \hat{X})) )^G.
\end{equation}
\end{defin}

Writing $F_{\hat{Z}}(S^V, (\Milnor \hat{X})^{\mfib})$ for the
spectrum of compactly supported sections of $(\Milnor 
\hat{X})^{\mfib})$ (denoted in
Appendix~\ref{sec:parametrized-spectra} by $\Gamma^c_{\hat{Z}}(S^V,
(\Milnor \hat{X})^{\mfib})$), we can use the adjunction homeomorphism  
\begin{equation}
F(S^V, C^{*,c}(\hat{Z}; \Milnor  \hat{X})) \cong C^{*,c}(\hat{Z};
F_{\hat{Z}}(S^V, (\Milnor  \hat{X})^{\mfib}))
\end{equation}
described in Lemma~\ref{lem:cotensor-inside} to see that this construction is homeomorphic to
the spectrum of sections that are compactly supported in the direction
of $\hat{Z}$.

\begin{lem}
An isomorphism from a product of Kuranishi charts $\bX_1 \times \bX_2$
to a boundary stratum of a Kuranishi chart $\bX$, determines a
commutative diagram 
\begin{equation} \label{eq:SW-milnor-twisted-desuspend-coefficients}
  \begin{tikzcd}[ column sep=tiny, transform shape, nodes={scale=1}]
    C^*(B G_{1}, \hat{X}|\hat{Z}^{-V_1}_{1} ) \wedge  C^*(B G_{2},
    \hat{X}|\hat{Z}^{-V_2}_{2} ) \ar[r] \ar[d] &  C^*(B
    G,\hat{X}|\hat{Z}^{-V_1\oplus V_2} )  \ar[d] \\
    C^{*,c_0}_{\rel \partial}(B \hat{Z}_{1}; \Milnor
    \hat{X}^{-V_1}_{1}  ) \wedge C^{*,c_0}_{\rel \partial}(B
    \hat{Z}_{2}; \Milnor  \hat{X}^{-V_2}_{2}  ) \ar[r] &
    C^{*,c_0}_{\rel \partial}(B \hat{Z}; \Milnor \hat{X}^{-V_1\oplus V_2} )  \\
    C^{*,c'_0}_{\rel \partial}(B \hat{Z}_{1}; \Milnor
    \hat{X}^{-V_1}_{1}  ) \wedge C^{*,c'_0}_{\rel \partial}(B
    \hat{Z}_{2}; \Milnor  \hat{X}^{-V_2}_{2}  ) \ar[r] \ar[u] \ar[d] &
    C^{*,c'_0}_{\rel \partial}(B \hat{Z}; \Milnor \hat{X}^{-V_1\oplus V_2} ) \ar[u] \ar[d] \\
    C^{*,c}_{\rel \partial}(B \hat{Z}_{1}; \Milnor  \hat{X}^{-V_1}_{1}
    ) \wedge C^{*,c}_{\rel \partial}(B \hat{Z}_{2}; \Milnor
    \hat{X}^{-V_2}_{2}  )  \ar[r] & C^{*,c}_{\rel \partial}(B \hat{Z};
    \Milnor \hat{X}^{-V_1 \oplus V_2} )
 \end{tikzcd}
\end{equation}
in which the vertical arrows on the right are equivalences and the
vertical arrows on the left induce equivalences on the derived smash
product.
\end{lem}

\begin{proof}
The vertical maps are compositions of product maps with boundary maps:
commutativity of the boundary maps follows by desuspending Diagram
\eqref{eq:SW-equivalence-Gamma-spectra} with the negative spheres,
applying the Borel cochains, and using the compatibility of function
spectra with smash products as in Equation
\eqref{eq:smash_product_map_functions} in the top right corner, and
similarly in the other entries of the middle and right column.
\end{proof}

Note that we can collapse the collar to compare the first row in
Diagram \eqref{eq:SW-milnor-twisted-desuspend-coefficients} with the
uncollared analogue. Moreover, smashing $\Milnor \hat{X}$ with a
spectrum $\bbk$ yields a corresponding commutative diagram with
coefficients in any spectrum.

\begin{rem} \label{rem:pushforward_pullback_not_commute}
At this stage, we point out the following difficulty: if $\bX \to
\bX'$ is an inessential map of $\langle \cS \rangle$-Kuranishi charts
(i.e. the corresponding covering map is a product), we obtain a
diagram  
\begin{equation}
\begin{tikzcd}
X|Z \ar[d] \ar[r] & \Map_{X}(Z, MX) \\
X'|Z'  \ar[r] & \Map_{X'}(Z', MX') \ar[u]
\end{tikzcd}
\end{equation}
which does not commute: the problem is that the section induced by an
element of $X|Z$ under the top horizontal map vanishes on the other
points lying in the same orbit, while the map induced by the
composition is invariant under the action. We shall return to this
discussion in Section~\ref{sec:local-glob-constr}, where we introduce
a map in the stable category represented by a zig-zag to resolve this
issue.
\end{rem}

\subsubsection{The Adams and norm isomorphisms for charts}
\label{sec:ambid-adams-isom}

There are two remaining steps in the comparison between virtual and
compactly supported cochains.  First, we use a natural map 
\begin{equation}
(-)_{hG} \to (-)^{hG}
\end{equation}
from the homotopy orbits to the homotopy fixed points, known as the
``norm map,'' which we review in Appendix~\ref{sec:norm-map}.

For any finite group $G$, the composition
\begin{equation}
 E G_+ \to S^0 \to \Map(EG_+, S^0) 
\end{equation}
induced by the projection $EG \to \ast$ yields for any $G$-spectrum
$Z$ a map 
\begin{equation}
\Sigma^\infty E G_+ \wedge Z \to F(\Sigma^\infty E G_+, Z).
\end{equation}
Passing to $G$-fixed points, we obtain the norm map
\begin{equation}
(\Sigma^\infty E G_+ \wedge Z)^G \to F(\Sigma^\infty E G_+, Z)^G.
\end{equation}
The righthand side is the homotopy fixed
point spectrum of $Z$ and the Adams isomorphism connects the lefthand
side via a natural zig-zag of weak equivalences to the homotopy orbits
$Z_{hG} = EG_+ \sma_G Z$.  The cofiber of the norm map is by
definition the Tate fixed points of $Z$.

Specializing to $X = Z \sma \bbk$ for a $G$-spectrum $Z$ and a Morava
$K$-theory spectrum $\bbk$, the norm map takes the form  
\begin{equation}\label{eq:norm}
(\Sigma^\infty E G_+ \sma X \sma \bbk)^G \to F(\Sigma^\infty E G_+,
  \bbk \sma X)^G,
\end{equation}
Since in this case the Tate fixed point spectrum vanishes,
Equation~\eqref{eq:norm} is an equivalence, as discussed in
Appendix~\ref{sec:ambidextrous-spectra}.  For the next statement, we
fix for each universe $U$ a fibrant replacement functor $\cQ_{U}(-)$ for
orthogonal $G$-spectra such that the assignment is externally
multiplicative in the sense that there are associative natural
transformations
\begin{equation}
\cQ_{U_1} \times \cQ_{U_2} \to \cQ_{U_1 \oplus U_2}
\end{equation}
of functors on $\Sp_{G_1} \times \Sp_{G_2}$ induced by the external
smash product (see Section~\ref{sec:fibrant} for the explicit
construction we use and a discussion of its
properties).  Note that in general, even when $U_1$ and $U_2$ are
complete universes, $U_1 \oplus U_2$ will not be a complete $G_1
\times G_2$ universe.  We ultimately handle this issue by working with
specific universes given in terms of countable direct sums of the
regular representation and using the natural map 
\begin{equation}
\rho_{G_1} \oplus \rho_{G_2} \to \rho_{G_1} \otimes \rho_{G_2} \cong
\rho_{G_1 \times G_2}
\end{equation}
specified by the linear extension of the maps specified by the
formulas 
\begin{equation}
\begin{aligned}
g_1 \mapsto \frac{1}{\sqrt{\# G_2}} \sum_{g \in G_2} g_1 \otimes g \\
g_2 \mapsto \frac{1}{\sqrt{\# G_1}} \sum_{g \in G_1} g \otimes g_2.
\end{aligned}
\end{equation}
We discuss this point further in the paragraphs surrounding  Equation~\ref{eq:regembed}.

\begin{lem}
When $\bbk$ is a Morava $K$-theory, the norm map induces a zig-zag of
equivalences  
\begin{equation}
\begin{tikzcd}
\left(\cQ_U(E G_+ \sma  X| Z^{-V}  \sma \bbk)\right)^{G} \ar[r,"\htp"] & (\cQ_U( C^*(E G; X|
  Z^{-V} \sma \bbk))^{G} \\
& \ar[u,"\htp"] C^*(B G,  X| Z^{-V}  \sma \bbk).
\end{tikzcd}
\end{equation}
Moreover, an isomorphism from a product of Kuranishi charts $\bX_1
\times \bX_2$ to a boundary stratum of a Kuranishi chart $\bX$
determines a commutative diagram  
\begin{equation}
  \begin{tikzcd}[column sep=tiny]
    \begin{aligned}
    (\cQ_{U_1}(E G_{1,+} \wedge X| Z^{-V_1}_1  \wedge \bbk))^{G_1}  \wedge 
  \\  (\cQ_{U_2}(E G_{2,+}  \wedge X| Z^{-V_2}_2  \wedge
    \bbk))^{G_2}  
    \end{aligned}
    \ar[r] \ar[d] &
    \begin{aligned}
    (\cQ_{U_1 \oplus U_2}(EG_{12,+} \wedge \\ X|
    Z^{-V_1\oplus V_2}_{12}  \wedge \bbk))^{G_{12}}      
    \end{aligned}
     \ar[d] \\
    \begin{aligned}
(\cQ_{U_1}( C^*(E G_1; X| Z^{-V_1}_1    \wedge \bbk))^{G_1}  \wedge \\
 (\cQ_{U_2}( C^*(E G_2; X| Z^{-V_2}_2    \wedge \bbk))^{G_2}      
    \end{aligned}
 \ar[r] &
 \begin{aligned}
    (\cQ_{U_1\oplus U_2}( C^*(E G_{12}; \\ X| Z^{-V_1\oplus V_2}   \wedge \bbk))^{G_{12}} 
 \end{aligned}
   \\
    \begin{aligned}
  C^*(B G_{1}; {X}|{Z}^{-V_1}_{1} ) \wedge \\  C^*(B G_{2};
  {X}|{Z}^{-V_2}_{2} )     
    \end{aligned}
  \ar[r] \ar[u] &  C^*(B
  G_{12};{X}|{Z}^{-V_1\oplus V_2}_{12} ) \ar[u] .
 \end{tikzcd}
\end{equation}
\end{lem}

\begin{proof}
The zig-zag of equivalences is comprised of the norm map, which is an
equivalence by hypothesis on $\bbk$, and the composite
\begin{multline}
C^*(BG; X|Z^{-V}) \cong F(EG_+,(X|Z^{-V})^{\mfib})^G \\ \to (\cQ_{U}
F(EG_+, (X|Z^{-V})^{\mfib}))^G.
\end{multline}
The commutativity of the product diagrams are a consequence of the externally
multiplicative properties of the functor $\cQ_{-}$ (see
Proposition~\ref{prop:Q-lax-mon} and
Lemma~\ref{lem:mult-universe-embedding}) and of the norm map (see
Proposition~\ref{prop:norm-mult}).
\end{proof}

On the other hand, the spectrum $EG_+ \wedge X| Z^{-V} \wedge \bbk$ is
(by construction) $G$-free.  For $G$-free spectra, the Adams
isomorphism provides a natural weak equivalence between the
$G$-homotopy orbits and the $G$-fixed points, see
Appendix~\ref{sec:adams-isom-ambid}.  To be more precise,
Definition~\ref{def:Adams-pointset-biequivariant} and
Proposition~\ref{prop:Adams-external-mult} establish the following:
\begin{lem}
There is a zig-zag of equivalences
\begin{equation}
\begin{tikzcd}
B X|Z^{-V} \sma \bbk & \ar[l,swap,"\htp"] \ldots \ar[r,"\htp"] & (\cQ_U (EG_+ \wedge X| Z^{-V}  \wedge
\bbk))^{G}.
\end{tikzcd}
\end{equation}
Moreover, an isomorphism from a product of Kuranishi charts $\bX_1
\times \bX_2$ to a boundary stratum of a Kuranishi chart $\bX$,
induces a commutative diagram 
\begin{equation}
\begin{tikzcd}
  \begin{aligned}
(\cQ_{U_1} (E G_{1,+} \wedge X| Z^{-V_1}_1  \wedge \bbk))^{G_1}
   \wedge  \\ (\cQ_{U_2} (E G_{2,+}  \wedge X| Z^{-V_2}_2  \wedge
    \bbk))^{G_2}    
  \end{aligned}
  &
  \begin{aligned}
    (B{X}|{Z}^{-V_1}_{1} \wedge  \bbk) \wedge \\
    (B {X}|{Z}^{-V_2}_{2} \wedge
    \bbk) 
  \end{aligned}
  \ar[l]  \\
  \ar[u] \ldots \ar[d]  & \ar[l] \ar[d] \ldots \ar[u]\\
  (\cQ_{U_1 \oplus U_2}
    (EG_+ \wedge X| Z^{-V_1\oplus V_2}  \wedge \bbk))^{G}   
  & B
        {X}|{Z}^{-V_1 \oplus V_2} \wedge \bbk \ar[l] .   
\end{tikzcd}
\end{equation}
\end{lem}

\begin{rem}
Our eventual treatment of the Adams isomorphism will be substantially
more complicated than might be gleaned from the above discussion,
because of the need to establish suitable functoriality with respect
to change of group. See Section~\ref{sec:using-adams-isom} and
Appendix~\ref{sec:adams-isom-ambid}.
\end{rem}

\subsubsection{Signpost: Spanier-Whitehead duality and Ambidexterity
  isomorphisms for charts} 
\label{sec:signp-span-whit}

To simplify the situation, let us consider a closed Kuranishi chart of
dimension $d-1$, and choose $\bbk$ to be a Morava $K$-theory spectrum:
We have constructed a zig-zag of equivalences   
\begin{equation}
  \begin{tikzcd}
  B X| Z^{-V-d} \wedge \bbk & \ar[l] \cdots \ar[r] &  (E G_+ \wedge X| Z^{-V-d} \wedge  \bbk)^G \ar[d] \\
 C^{*,c}(B Z; \Milnor  X^{-V-d}\wedge \bbk)&  \ar[l] \cdots \ar[r] &  C^*(B G,  X| Z^{-V-d} \wedge  \bbk ).  
  \end{tikzcd}  
\end{equation}
There are two essential remaining points at this stage: 
\begin{enumerate}
\item Studying the compatibility of these constructions with maps and
  products of Kuranishi charts, and 
\item comparing the final term with an untwisted cochain spectrum. 
\end{enumerate}
We carry out the first part in
Section~\ref{sec:local-glob-constr}. For the second, we note that, if
the manifold $X$ is $\bbk$-oriented, we obtain an equivalence of spectra
over $Z$:
\begin{equation}
\Milnor X  \wedge \bbk \cong Z \barwedge \Sigma^{d-1} \bbk,
\end{equation}
where $ Z \barwedge \Sigma^{d-1} \bbk $  is the trivial parametrized
spectrum with fibre $S^{d-1} \wedge \bbk $.  If this orientation is
$G$-equivariant, we conclude that we have an equivalence of orthogonal
$G$-spectra  
\begin{equation}
C^{*,c}(Z; \Milnor  X^{-V-d} \wedge \bbk)  \to C^{*,c}(Z;  \Omega \bbk).
\end{equation}
Assuming that we have started with a global chart, we end up with the
desired zig-zag 
\begin{equation}
\Omega \bbk \to   C^{*}(BZ;  \Omega \bbk) \leftarrow \cdots \to  B X|
Z^{-V-d} \wedge \bbk \to \Omega^{d} \bbk, 
\end{equation}
thus producing a homotopy type.

Because the problem of writing down orientations is non-trivial in the
geometric contexts that motivate our applications, we shall next turn
to the problem of finding conditions which guarantee the existence of
orientations. 

\subsection{Flag smooth Kuranishi charts and tangential spherical fibrations}
\label{sec:smoothly-fiber-kuran-1}

As we review in Section~\ref{sec:ambid-orient}, the Morava
$K$-theories are complex oriented cohomology theories.  As a
consequence, our strategy to produce $K(n)$-orientations of Kuranishi
presentations is to lift the spherical tangent fibration arising from
the Milnor construction to a (stably complex) vector bundle, then
appeal to the orientability of complex vector bundles with respect to
the generalized cohomology theories which we consider.  The most
straightforward way to do so would be to work in the context of smooth
Kuranishi charts; however, applying this to Floer theory requires more
analytic work than we are willing to carefully implement. Instead, our
strategy is to use the fact that the Kuranishi charts which appear in
our examples submerse over smooth manifolds with smooth fibres. This
notion was axiomatised in Section~\ref{sec:fiber-kuran-pres}, but we
recall all the basic notions here (since we do not discuss
functoriality):

If $X$ and $B$ are (topological) $\langle \cS \rangle$-manifolds, we
say that a map 
\begin{equation}
  X \to B  
\end{equation}
is a topological submersion if it is locally homeomorphic to a projection
\begin{equation}
  \bR^{n_1+n_2} \times [0,\infty)^{Q} \to \bR^{n_1} \times [0,\infty)^{Q},  
\end{equation}
in a neighbourhood of each point lying in a stratum of $X$ labelled by $Q$, via a stratum-preserving map.  We also define a fibrewise smooth structure on a topological submersion $X \to B$ as a choice of atlas for $X$ consisting of product charts $\bR^{n_2} \times U \to X$ over charts $U \to B$, with transition functions which are continuously differentiable in the fibre direction (depending continuously on the base). 

With the above in mind, we recall from Definition
\ref{def:smoothly_fibred_complete}, that a flag smooth $\langle \cS
\rangle$-Kuranishi chart consists of
\begin{enumerate}
\item a  Kuranishi chart $\bX \in
\Kur\langle \cS \rangle $, 
\item a smooth $\langle \cS
\rangle$-manifold $B$ equipped with a $G$ action,
\item a $G$ equivariant topological stratified submersion $\pi \co X \to B$, and
\item a fibrewise smooth structure on $\pi$.
\end{enumerate}
We require that the map $s \co X \to V$ be smooth on each fibre of
$\pi$.  

\begin{rem}
In applications, $X$ will be a moduli space of maps from a family of Riemann surfaces with marked points, and $B$ will be an abstract moduli space of Riemann surfaces. In order to arrange for $B$ to also be stratified by $\cS$, we shall impose constraints on the  degenerations in $X$.
\end{rem}

We shall presently see that such a flag smooth structure on a Kuranishi chart induces a lift of its Milnor spherical fibration to a vector bundle. The starting point is to consider the fibrewise tangent space
\begin{equation}
  T^\pi X \to X
\end{equation}
of the projection $X \to B$. 

\begin{defin}
We define the tangent space of $X$ to be the direct sum of the fibrewise tangent space with the
pullback of the tangent space of $B$, 
\begin{equation}
  TX \equiv   T^\pi X  \oplus \pi^* T B.
\end{equation}
We write $S^{TX}$ for the sphere bundle of $TX$, and $S^{TX}|0$
for the spherical fibration over $X$ obtained as the fibrewise cone
of the complement of the origin in $S^{TX}$. 
\end{defin}
For our later discussion, we note that the inclusion
$S^{TX} \to S^{TX}|0 $ is a fibrewise homotopy equivalence.

There is a natural notion of a product of flag smooth Kuranishi
charts, given by taking the product of the submersions, and the
construction of the tangent bundle of the chart is multiplicative, in
the following sense:
\begin{lem}
For flag smooth Kuranishi charts $\bX_1$ and $\bX_2$, there is a
commutative diagram 
\begin{equation}
  \begin{tikzcd}
S^{TX_1}  \barwedge S^{TX_2}  \ar[r] \ar[d] & S^{T (X_1 \times X_2)} \ar[d]\\  
    S^{TX_1}|0 \barwedge S^{TX_2}|0 \ar[r] & S^{T (X_1 \times X_2)}|0
  \end{tikzcd}
\end{equation}
of homotopy equivalences of spherical fibrations over $X_1 \times X_2$. \qed 
\end{lem}

Unfortunately, there is no direct map from $S^{TX}|0$ to $\Milnor
X$. We shall instead construct a correspondence comparing them via a
variant of Nash's model for the tangent space.  

\subsubsection{The Nash spherical fibration}
\label{sec:nash-spher-fibr}

\begin{defin}[c.f. \cite{Nash1955}]
If $B$ is a topological manifold, a \emph{Nash path} $\gamma \co
[0,\infty) \to B$ is a path which is either constant or satisfies the
  property that $\gamma(t) \neq \gamma(0)$ for $t \neq 0$. 
\end{defin}

We now assume that $B$ is a smooth manifold.

\begin{defin}
When $B$ is smooth, we define $N B \to B$ to be the fiber bundle of
Nash paths which are differentiable at the origin.  We let $\Nash B$ be
the mapping cone (homotopy cofibre) of the inclusion of paths with non-zero derivative.
\end{defin}

The evaluation maps at any non-zero point (say $1$) and of the
derivative at the origin yield a zig-zag: 
\begin{equation}
  \Milnor B \leftarrow \Nash B \to S^{TB}|0.
\end{equation}
These maps can be shown to be fibrewise equivalences over $B$ by a
direct application of Nash's argument
(c.f. Proposition~\ref{prop:Nash-equivalence-square} below).  This
zig-zag is compatible with products: 

\begin{lem}
The product of paths induces a fibrewise equivalence $\Nash B_1
\barwedge \Nash B_2 \to \Nash \left(B_1 \times B_2 \right)$, which
fits in a commutative diagram  
\begin{equation} \label{eq:evaluation_maps_spherical_fibration}
\begin{tikzcd}
\Milnor B_1 \barwedge \Milnor B_2 \ar[d] & \ar[l]  \Nash B_1\barwedge \Nash B_2  \ar[r] \ar[d] & S^{TB_1}|0  \barwedge S^{TB_2}|0  \ar[d] \\
\Milnor  \left( B_1 \times B_2 \right) & \ar[l]  \Nash \left(  B_1 \times B_2\right) \ar[r] & S^{T \left(B_1 \times B_2 \right)}|0 .
\end{tikzcd}
\end{equation} \qed
\end{lem}

Assuming that $\pi \co X \to B$ is a topological submersion over a
smooth manifold $B$ equipped with a  flag smooth structure, and that
both $X$ and $B$ are without boundary, we have: 
\begin{defin}\label{defin:nash-fibration} 
The \emph{Nash tangent space} $N X$ is the space of maps $\gamma
\colon [0,\infty)^2 \to X$ whose restriction to the diagonal is a Nash
  path, and such that the following properties hold near the origin:  
\begin{enumerate}
\item The composition $\pi \circ \gamma$ is differentiable, and is independent of the second coordinate.
\item The family of curves in the fibres of $\pi$ parametrised by the first coordinate are differentiable at the origin, and the corresponding path of derivatives is continuous.
\end{enumerate}

The \emph{Nash spherical fibration} $\Nash X$ is the cone of the
complement of the locus $ \{ \gamma \mid d \gamma = 0\}  \subset N X $
given by the vanishing of both directional derivatives at the origin. 
\end{defin}

The differentiability conditions we have imposed imply that the Nash path obtained by restriction to the diagonal can only be constant whenever the two directional derivatives vanish.

\begin{prop} \label{prop:Nash-equivalence-square}
  The map $\Nash X \to X$ is a sectioned fibre bundle, and the evaluation maps
  \begin{equation}
    (X \times X, X \times X \setminus X) \leftarrow (N X, N X
    \setminus \{ \gamma \mid  d \gamma = 0\}  )  \rightarrow (TX, TX \setminus X)  
  \end{equation}
are fibrewise homotopy equivalences of pairs. In particular, the
induced maps 
\begin{equation}
\Milnor X \leftarrow \Nash X \to S^{TX}|0
\end{equation}
are equivalences of spherical fibrations, which are compatible with
products as in
Diagram~\eqref{eq:evaluation_maps_spherical_fibration}.
\end{prop}

\begin{proof}
Consider the group $\Diff^\pi(X)$ of homeomorphisms of $X$ which lift
a diffeomorphism of $B$ and whose restriction to fibers is given by a
continous family of smooth maps. The action of $\Diff^\pi(X) $ is
locally transitive, and lifts to the spaces $\Milnor X \leftarrow
\Nash X  \rightarrow S^{TX}|0 $. It thus suffices to show that the map
of fibers is a homotopy equivalence. 

To this end, we choose a local trivialisation of $X$ as a product $
\bR^{n+k} \cong \bR^n \times \bR^k$, with the map to $B$ given by the
projection to $\bR^k$. Letting $\Nash_0 \bR^{n+k}$ denote the fiber
at the origin, the inclusion of linear maps gives a splitting of the map 
\begin{equation}
\Nash_0 \bR^{n+k} \to     \bR^{n+k}| \bR^{n+k} \setminus 0
\end{equation}
where the right hand side can be identified either with the fiber of
$\Milnor \bR^{n+k}$ at the origin, or with the fiber of
$S^{n+k}|0$. It thus remains to prove that $\Nash_0 \bR^{n+k}$
deformation retracts onto its subspace of linear maps.

The argument given by Nash applies verbatim: there are two steps to
the homotopy, the second of which is the straight-line homotopy, but
the problem is that this may violate the condition that $0$ be the
only point along the diagonal mapping to the origin. So the first step
is to run the homotopy which ``pulls back'' under the obvious
continuous splitting of $ [0,\epsilon)^2 \to [0,\infty)^2$. The key
point is to pick $\epsilon$ continously varying in all the
parameters so that the resulting straight-line homotopy avoids the
basepoint. This is ensured by the $C^1$ conditions near the
origin, which provide a function $\epsilon_\gamma$ such that the
origin is the only intersection point of the image of the diagonal
embedding $[0,\epsilon_\gamma] \to [0,\infty)^2$ under $\gamma$
with the ray generated by the negative of the sum of the
directional derivatives of $\gamma$ at $0$. Since all
constructions take place before taking cones, this construction
gives the desired fiberwise homotopy equivalence, which
preserves the sections by construction. 
\end{proof}

The preceding proposition justifies the terminology \emph{spherical
  fibration} for $\Nash X$, since it proves that its fiber is homotopy
equivalent to a sphere of dimension equal to that of $X$.

If $X$ admits an action of a group $G$ preserving the flag smooth
structure, then our three models for the tangential spherical
fibration are $G$-equivariant fiber bundles over $X$, whose underlying
non-equivariant fiber bundles are homotopy equivalent.

\begin{rem}
In fact, these fiber bundles are equivariantly homotopy equivalent,
but this shall not be required for our purpose. 
\end{rem}

\subsubsection{Boundary maps for tangential fibrations}
\label{sec:bound-maps-tang}

Given a manifold $X$ with corners, we noted in Section
\ref{sec:span-whit-dual-2} that there is a natural map $\Milnor
\partial X \barwedge \Milnor \kappa \to \Milnor \hat{X}$ over the
collar, induced by a section of the Milnor fibration of the collar
$\Milnor \kappa \to \kappa$. Our goal in this section is to extend
this boundary map to the other models of the tangent spherical
fibration.

If $\pi\co X \to B$ is a flag smooth $\langle \cS \rangle$-Kuranishi chart, then the induced map of collared completions
\begin{equation} \label{eq:fibrewise_smooth_submersion-completions}
\hat{\pi} \co  \hat{X} \to \hat{B}  
\end{equation}
admits the structure of a flag smooth topological submersion, which is canonical up to contractible choice: the smooth structure on $\hat{B}$ is determined by the germ of a collar along the boundary strata of the manifold $B$. We fix such a choice in this section.

Define $N\hat{X}$ to be the restriction to $\hat{X}$ of the Nash
tangent space of the union of $X$ with infinite collars as a flag
smooth submersion over the union of $B$ with infinite collars. 

The product decomposition along the collar yields a commutative diagram
\begin{equation} \label{eq:map_product_tangents_boundary_interior}
  \begin{tikzcd}
    \Milnor \partial^Q \hat{X} \barwedge \Milnor \kappa^Q \ar[d] &    \Nash \partial^Q \hat{X} \barwedge \Nash \kappa^Q \ar[d] \ar[l] \ar[r] & S^{T \partial^Q \hat{X}}|0 \barwedge S^{T\kappa^Q}|0 \ar[d]  \\
    \Milnor \hat{X} &  \Nash \hat{X} \ar[l] \ar[r] & S^{T \hat{X}}|0
  \end{tikzcd}
\end{equation}
where the second vertical arrow is given by the composition
\begin{equation}
  [0,\infty)^2 \to [0,\infty)^2 \times [0,\infty) \to \partial^Q \hat{X} \times \kappa^Q \to \hat{X},  
\end{equation}
in which the first arrow is $(t_1,t_2) \to (t_1, t_2, t_2)$.

We shall use this diagram to construct a map from sections of the tangent space over a boundary stratum to global sections. We start by considering the interval: given a pairs of points $(t_0,t_1) \in \kappa$ consider the path $\gamma(t_0,t_1) \co [0,\infty) \to \kappa$ whose restriction to $[0,1]$ is the straight path from $t_0$ to $t_1$, and which has constant value $t_1$ outside the unit interval. This construction gives an embedding
\begin{equation} 
    \Milnor \kappa  \to \Nash \kappa
\end{equation}
which is a splitting of the evaluation map, and can be composed with the derivative at the origin to give a map $\Milnor \kappa \to  S^{T \kappa}|0$. Taking products, we obtain maps
\begin{equation}
\Milnor \kappa^Q  \to \Nash \kappa^Q \to  S^{T\kappa^Q}|0 .
\end{equation}

Fixing the product section of $\Milnor \kappa^Q$ chosen in Equation \eqref{eq:section_Milnor_collar_vanish_boundary}, which vanishes at the endpoints of the collar, we can thus assign to a section of each of these spherical fibrations over $ \partial^Q \hat{X}$ a corresponding section of the spherical fibration over $\hat{X}$ which vanishes at the boundary, and thus obtain a commutative diagram:
\begin{equation}
      \begin{tikzcd}
    \Milnor \partial^Q \hat{X}  \ar[d] &    \Nash \partial^Q \hat{X}  \ar[d] \ar[l] \ar[r] & S^{T\partial^Q \hat{X}}|0  \ar[d] \\
    \Milnor \hat{X} &  \Nash \hat{X} \ar[l] \ar[r] & S^{T \hat{X}}|0.
  \end{tikzcd}
\end{equation}

\subsubsection{Tangentially twisted cochains}
\label{sec:tang-twist-coch}

The outcome of the previous section is that we can assign to a flag
smooth Kuranishi chart $\bX = (X,V, s,G,\pi,B)$ a $G$-equivariant
spherical fibration $\Nash X$ over $Z$, with maps to the Milnor and
vectorial tangent spaces. Passing to the spectra of section over the
Borel construction, and to the collared completion, we obtain a
diagram of maps
\begin{equation}
\begin{tikzcd}
C^{*,c}(B\hat{Z}; \Milnor  \hat{X}^{-V} )  &  C^{*,c}(B\hat{Z};
\Nash  \hat{X}^{-V} ) \ar[l] \ar[d] \\  \ar[r] C^{*,c}(B\hat{Z}; S^{T  \hat{X} -V} ) &  C^{*,c}(B\hat{Z}; S^{T
  \hat{X} -V}|0 ) .
\end{tikzcd}
\end{equation}

It is straightforward to check that these maps are compatible with
products, as long as germs of smooth collars of the bases are
appropriately chosen: 

\begin{defin} \label{def:smooth_collars_compatible}
A \emph{compatible choice of smooth collars} for the strata of an
$\langle \cS \rangle$-stratified smooth manifold $B$ consists of a
germ of smooth embedding for each  
  \begin{equation}
    \partial^{Q} B     \times [0,1)^{Q \setminus Q'} \to \partial^{Q'} B,
  \end{equation}
  for each pair $Q' \subset Q$ of totally ordered subset $Q$, which extends the inclusion $ \partial^Q B  \times \{0\}^{Q \setminus Q'} \to B$. We require that, for each triple $Q'' \subset Q' \subset Q$, the following diagram commutes:
  \begin{equation}
    \begin{tikzcd}
      \partial^{Q} B   \times [0,1)^{Q \setminus Q'} \times [0,1)^{Q' \setminus Q''} \ar[r] \ar[d] &  \partial^{Q'} B \times [0,1)^{Q' \setminus Q''} \ar[d] \\
      \partial^{Q} B   \times [0,1)^{Q \setminus Q''} \ar[r] & \partial^{Q''} B  .
    \end{tikzcd}
  \end{equation}
\end{defin}

Compatible choices of smooth collars on $\langle \cS_i
\rangle$-stratified smooth manifolds $B_i$ for $i \in \{1,2\}$ induce
such a choice for their product, and a choice of smooth collars on $B$
induces one on each boundary stratum.   

\begin{lem}
An isomorphism of the product $\bX_1 \times \bX_2$ of a pair of flag
smooth Kuranishi charts with a stratum of a Kuranishi chart $\bX$, and
a choice of smooth collars for the strata of $B$ which restricts to
product collars for the strata of $B_1 \times B_2$, determine a
commutative diagram  
\begin{equation}
\begin{tikzcd}[column sep=tiny, transform shape, nodes={scale=1}]
C^{*,c}_{\rel \partial}(B\hat{Z}_{1}; \Milnor \hat{X}^{-V}_{1} ) \wedge C^{*,c}_{\rel \partial}(B\hat{Z}_{2}; \Milnor \hat{X}^{-V}_{2} )   \ar[r] &   C^{*,c}_{\rel \partial}(B\hat{Z}; \Milnor \hat{X}^{-V} )   \\
C^{*,c}_{\rel \partial}(B\hat{Z}_{1}; \Nash \hat{X}^{-V}_{1} ) \wedge C^{*,c}_{\rel \partial}(B\hat{Z}_{2}; \Nash \hat{X}^{-V}_{2} ) \ar[u] \ar[d] \ar[r] &   C^{*,c}_{\rel \partial}(B\hat{Z}; \Nash \hat{X}^{-V} ) \ar[u] \ar[d]  \\
C^{*,c}_{\rel \partial}(B\hat{Z}_{1}; S^{T \hat{X} -V}_{1}|0 ) \wedge C^{*,c}_{\rel \partial}(B\hat{Z}_{2}; S^{T\hat{X}-V}_{2}|0 )  \ar[r]
&    C^{*,c}_{\rel \partial}(B\hat{Z}; S^{T \hat{X} -V}|0 ) 
\end{tikzcd}
\end{equation}
in which the righthand vertical arrows are equivalences and the
lefthand vertical arrows induce equivalences of derived smash
products.   \qed 
\end{lem}

We now provide a functorial and multiplicative comparison between the
models $S^{TX}|0$ and $S^{TX}$. The key problem is that the section of
$S^{T \kappa}|0$ chosen to map sections on the boundary to sections on
$X$ does not lie in the ordinary tangent bundle. On the other hand,
any identification $(0,1) \cong \bR$ gives a section of $S^{T\kappa}$
which vanishes at the boundary, and hence a section of $S^{T
  \kappa}|0$ by inclusion.  Thus, it will suffice to compare the
constructions associated to these two sections.

Let $\hat{\kappa}$ denote the interval $[0,2]$, and consider the
section of $S^{T \hat{\kappa}}|0$ which is given by the section
induced by  Equation \eqref{eq:section_Milnor_collar_vanish_boundary}
in the interval $[0,1]$, and by (the closure of) an identification of
the interval  $(1,2)$ with $ \bR$ along the second part. These two
sections both vanish (i.e., have value the basepoint) at $1$, so the
construction is well-defined. Taking products, we obtain a section of
$S^{T \hat{\kappa}^Q}|0$ for every finite set $Q$. The key observation
is that we have two collapse maps 
\begin{equation}
\kappa^Q \leftarrow  \hat{\kappa}^Q \to \kappa^Q
\end{equation}
corresponding to the first and the second interval, and these map the
``concatenated sections'' to the first and the second section,
respectively.

Let $\hat{\hat{X}}$ be the completion of a flag smooth $\langle \cS
\rangle$-manifold $X$ obtained by attaching the collar $\hat{\kappa}^Q
\times \partial^Q X$ along the stratum labelled by $Q$.  We then have
natural maps 
\begin{equation}\label{eq:collar-collapse-Milnor}
\hat{X} \leftarrow \hat{\hat{X}} \to \hat{X}
\end{equation}
obtained respectively by collapsing the first and the second sets of
collars. 

Whenever $Q$ labels a boundary stratum of $X$, we obtain a commutative diagram
\begin{equation}
\begin{tikzcd}
S^{T \partial^Q \hat{X}}|0  \ar[d] &    S^{T \partial^Q
  \hat{\hat{X}}}|0  \ar[d] \ar[l] \ar[r] & S^{T \partial^Q \hat{X}}|0
\ar[d] &  S^{ T\partial^Q \hat{X}} \ar[d]  \ar[l]\\ 
S^{T \hat{X}}|0 &  S^{T \hat{\hat{X}}}|0 \ar[l] \ar[r] & S^{T
  \hat{X}}|0&   S^{ T\hat{X}}  \ar[l]
\end{tikzcd}
\end{equation}
where the first vertical arrow corresponds to the section that moves
in the cone direction, and the third vertical arrow uses the
identification of the collar with $\bR$.  

We thus obtain a diagram of maps of compactly supported relative
cochains.  This comparison diagram is multiplicative:

\begin{lem}\label{lem:multiplicativity-comparison-Whit-tangent}
An isomorphism of the product $\bX_1 \times \bX_2$ of a pair of flag
smooth Kuranishi charts with a stratum of a Kuranishi chart $\bX$, together with a compatible choice of collars, 
determines a commutative diagram 
\begin{equation}
\begin{tikzcd}[column sep=tiny, transform shape, nodes={scale=1}]
C^{*,c}_{\rel \partial}(B\hat{Z}_{1}; S^{T \hat{X}-V}_{1}|0 ) \wedge C^{*,c}_{\rel \partial}(B\hat{Z}_{2}; S^{T \hat{X} -V}_{2}|0 )   \ar[r] &   C^{*,c}_{\rel \partial}(B\hat{Z}; S^{T \hat{X}-V}|0 )   \\
C^{*,c}_{\rel \partial}(B\hat{\hat{Z}}_{1}; S^{T \hat{\hat{X}}-V}_{1}|0 ) \wedge C^{*,c}_{\rel \partial}(B\hat{\hat{Z}}_{2}; S^{T \hat{\hat{X}}-V}_{2}|0 ) \ar[u] \ar[d] \ar[r] &   C^{*,c}_{\rel \partial}(B\hat{\hat{Z}}; S^{T \hat{\hat{X}}-V}|0 ) \ar[u] \ar[d]  \\
C^{*,c}_{\rel \partial}(B\hat{Z}_{1};S^{T \hat{X}-V}_{1}|0 ) \wedge C^{*,c}_{\rel \partial}(B\hat{Z}_{2}; S^{T \hat{X} -V}_{2}|0 )  \ar[r] 
 &    C^{*,c}_{\rel \partial}(B\hat{Z}; S^{T \hat{X} -V}|0 ) \ \\
C^{*,c}_{\rel \partial}(B\hat{Z}_{1}; S^{T \hat{X}_{1} -V} ) \wedge C^{*,c}_{\rel \partial}(B\hat{Z}_{2};S^{T \hat{X}_{2} -V} )  \ar[r]\ar[u]
&    C^{*,c}_{\rel \partial}(B\hat{Z}; S^{T\hat{X} - V})\ar[u]
\end{tikzcd}
\end{equation}
in which the righthand vertical arrows are equivalences and the
lefthand vertical arrows induce equivalences of derived smash
products.  \qed 
\end{lem}

\subsection{Complex-oriented charts}
\label{sec:compl-lifts-tang}

In this section, we shall complete the comparison between different
models for the cochains on a chart.  We prove that a choice of stable
almost complex structure on the virtual tangent bundle induces a
trivialization of the tangentially twisted cochains whenever our
choice of spectral coefficients are complex-oriented.  The formalism
that we adopt for doing this is designed to allow us to globalise
these constructions.

\subsubsection{Stable vector bundles over orbispace charts}
\label{sec:stable-vect-bundl}

We begin by considering vector bundles attached to charts.

\begin{defin}
We define a \emph{stable vector bundle} over an orbispace chart $(Z,
G)$ to be a pair $(I, V)$, where $I$ is a $G$-equivariant vector
bundle over $Z$ and $V$ is a finite-dimensional real
$G$-representation.  We write $d_I$ and $d_V$ for the dimensions of
$I$ and $V$ respectively, and define the \emph{virtual dimension} of
such a stable bundle to be the difference $d_I - d_V$.
\end{defin} 

Let $\ell$ denote a trivialized real line; i.e., as above, we have a
fixed identification of $\ell$ with $\bR$.  Both $I$ and $V$ give rise
to vector bundles over the Borel construction $BZ = E G \times_{G} Z
$.  In a mild abuse of notation, we denote the associated spherical
fibrations by $S^I$ and $S^V$.

\begin{defin}
Given a stable vector bundle over an orbispace chart, we have the
parametrized spectrum
\begin{equation}
  S^{I -V - \ell} \equiv F_{BZ}(S^{V + \ell}, (S^I)^{\mfib}).
\end{equation}
For an associative ring spectrum $\bbk$, we associate to this setup
the space of compactly supported relative cochains from
Definition~\ref{defin:relative-compactly-supported}:
\begin{equation}
  C^{*,c}_{\rel \partial}\left(B Z; S^{I -V - \ell} \sma \bbk
  \right)
\equiv F(B\hat{Z}^+/B\partial\hat{Z}^+, S^{I-V-\ell} \sma \bbk^{\mfib}).
\end{equation}
\end{defin}

\begin{rem}
If $Z$ is the zero locus of a Kuranishi chart, then the tangent bundle
of $X$, together with the obstruction space $V$, yields a stable
vector bundle over $Z$. Note as well that we are considering the case
of an $\langle \cS \rangle$-Kuranishi chart here; in the case of a
$\langle \partial^Q \cS \rangle$-chart, considered as an object of
$\Kur^{\cS}$, we shall more generally take the direct sum of $I$
with $\bR^Q$.
\end{rem}

Consider the following situation: 
\begin{enumerate}
\item $\bX_1$ and $\bX_2$ are a pair of Kuranishi charts, and $\bX_1
  \times \bX_2$ is a codimension $1$ boundary stratum of a chart
  $\bX$, and
\item $\bX_i$ and $\bX$ are respectively equipped with stable vector
  bundles $(I_{i}, V_i)$ and $(I,V)$ and lines  $\ell_i$ and $\ell$. 
\end{enumerate}

\begin{defin}
We define a \emph{compatibility isomorphism} to be $G$-equivariant
isomorphisms of vector bundles over $Z_{1} \times Z_{2}$ and inner
product spaces 
\begin{align} \label{eq:isomorphism_product_bundles}
I & \cong I_1 \times I_2 \\ \label{eq:isomorphism_product_vector_spaces}
V & \cong V_1 \times V_2.
\end{align}
\end{defin}

Pulling back $I$ under the projection $\hat{Z} \to Z$ yields a
collared vector bundle over $\hat{Z}$, i.e., a canonical
identification over each part of the collar with the pullback of the
restriction of $I$ to the corresponding stratum. 

\begin{lem}\label{lem:comp-products}
Given charts as above and a compatibility isomorphism, the product of
compactly supported relative cochains discussed in Lemma~\ref{lem:compactly-supported-product} followed by the boundary
inclusion of Equation~\eqref{eqn:compactly-supported-boundary-map}
yield a map   
\begin{equation}\label{eq:product_map_cochains_vector_bundle}
\begin{tikzcd}
C^{*,c}_{\rel \partial}\left(B \hat{Z}_{1}; S^{I_1-V_1 - \ell_1} \wedge \bbk
\right) \wedge C^{*,c}_{\rel \partial}\left(B \hat{Z}_{2}; S^{I_2-V_2 -
  \ell_2} \wedge \bbk  \right) \ar[d]  \\  C^{*,c}_{\rel
  \partial}\left(B \hat{Z}; S^{I-V - \ell} \wedge \bbk \right). 
\end{tikzcd}
\end{equation} \qed
\end{lem}

The compatibility isomorphisms in our setting arise as part of
classifying map data that is compatible the monoidal structure on the
categories of Kuranishi charts equipped with tangent bundles that we
consider; this implies coherent associativity on the product maps
described in Lemma~\ref{lem:comp-products}.  We will return to this
point in Section~\ref{sec:funct-mult-triv}.

\subsubsection{Complex vector bundles and orientations}
\label{sec:compl-vect-bundl}

We now study the consequences of a complex structure on the vector
bundles associated to charts.

\begin{defin}
A \emph{complex stable vector bundle} over an orbispace chart $(Z,G)$
is a stable vector bundle over the chart such that $I$ and $V$ are
equipped with complex structures, which are preserved by the
$G$-actions and compatible with the inner product.  
\end{defin}
Passing to the
Borel construction, a complex stable vector bundle gives rise to a stable complex vector bundle over $BZ$.
\begin{rem}
In our application, the stable complex vector bundle will arise from a 
complex-linear $\dbar$ operator over the moduli space of Floer
solutions. In order to define such an operator, one must impose
growth/decay conditions at infinity which are best formulated by
choosing a complex trivialisation over each orbit of the tangent
bundle of the ambient symplectic manifold. With such a choice, we
shall consider solutions which converge to a constant along one end,
and vanish along another. This asymmetric choice is required for
multiplicativity to hold in the form discussed in the next section.  
\end{rem}

It is illuminating to express the notion of a complex orientation in
terms of parametrized spectra.  For a $G$-representation $V$ over a
$G$-space $Z$, denote by $S^V_{BZ}$ the parametrized space over $BZ$
with total space $S^V \times BZ$ and fiber $S^V$; as usual, when $V =
\bR^n$, we write $S^n_{BZ}$ for this parametrized space.  Let $\bbk$
denote an associative ring spectrum.  Then the discussion in
Section~\ref{sec:multiplicative-orientations} shows that a
$\bbk$-orientation of a spherical fibration $f \colon E \to BZ$ with
fiber $S^V$ gives rise to a trivialization in the form of an
equivalence of parametrized spectra
  \begin{equation}
E \sma \bbk \to S^{|V|}_{BZ} \sma \bbk
  \end{equation}
over $BZ$, where we are abusively denoting the fibrewise smash product of
parametrized spectra and the tensor of parametrized spectra over
spectra as ``$\sma$''.

If we assume that $\bbk$ admits an equivariant complex
orientation, returning to the setting of a complex stable vector
bundle over an orbispace chart, we have an equivalence
\begin{equation}
S^{I} \wedge  \bbk \to S_{BZ}^{d_I} \sma \bbk.
\end{equation}
More generally, we have the following zig-zag of equivalences
\begin{equation}
  \begin{tikzcd}
C^{*,c}_{\rel \partial}\left(B Z; S^{I-V - \ell} \wedge \bbk  \right)
\ar[r] &  C^{*,c}_{\rel \partial}(B Z; S^{d_I - V - \ell} \sma \bbk)  \\  & C^{*,c}_{\rel \partial}(B Z; S^{d_I -d_V - \ell} \wedge \bbk) \ar[u].   
  \end{tikzcd}
\end{equation}
Finally, desuspending by the models for spheres discussed in
Appendix~\ref{sec:action} and applying the multplicative comparisons
of Proposition~\ref{prop:lurie-comparison}, we obtain a zig-zag of
equivalences  
\begin{equation}\label{eq:cancel-spheres}
C^{*,c}_{\rel \partial}(B Z;  S^{I - V - \ell} \wedge \bbk)
\left[\frac{d_V - d_I}{2}\right] \to \cdots \leftarrow  C^{*,c}_{\rel
  \partial}(B Z; \Omega^\ell \bbk). 
\end{equation}

\subsubsection{Stable complex vector bundles over orbispace charts}
\label{sec:stable-compl-vect}

Returning to the setting of Section \ref{sec:stable-vect-bundl}, we
say that a compatibility isomorphism for Kuranishi charts equipped
with stable complex bundles consists of maps as in Equations
\eqref{eq:isomorphism_product_bundles} and
\eqref{eq:isomorphism_product_vector_spaces} which respect the
complex structures.

As discussed in Section~\ref{sec:multiplicative-orientations},
assuming that $\bbk$ admits a multiplicative complex orientation, given
vector bundles $I_1$ and $I_2$ which are respectively defined over $BZ_1$ and $BZ_2$,  the
diagram 
\begin{equation}
\begin{tikzcd}
(S^{I_1} \sma  \bbk) \barwedge (S^{I_2} \sma \bbk) \ar[r] \ar[d] & (S_{BZ_1}^{d_{I_1}}
  \sma \bbk) \barwedge (S_{BZ_2}^{d_{I_2}} \sma \bbk) \ar[d] \\ 
S^{I_1 \bar{\oplus} I_2} \sma \bbk \ar[r] & S_{BZ_1 \times
  BZ_2}^{d_{I_1} + d_{I_2}} \sma \bbk
\end{tikzcd}
\end{equation} 
of spectra over $BZ_1 \times BZ_2$ commutes, where here the fiberwise
smash product is the external smash product from spaces over $BZ_1$ and $BZ_2$
to $BZ_1 \times BZ_2$ and analogously with the direct sum
$\bar{\oplus}$.  Then a compatibility isomorphism for Kurniahsi charts
equipped with stable complex bundles determines a commutative diagram  
\begin{equation}
    \begin{tikzcd}
      \begin{aligned}
        C^{*,c}_{\rel \partial}\left(B \hat{Z}_{1}; S^{I_1-V_1 - \ell_1} \wedge  \bbk \right)  \wedge  \\  C^{*,c}_{\rel \partial}\left(B \hat{Z}_{2};  S^{I_2-V_2 -\ell_2} \wedge \bbk \right)
      \end{aligned}
       \ar[r] \ar[d] &  C^{*,c}_{\rel \partial}\left(B \hat{Z};  S^{I-V - \ell} \wedge \bbk  \right) \ar[d] \\
       \begin{aligned}
      C^{*,c}_{\rel \partial}\left(B \hat{Z}_{1}; S^{d_{I_1}-V_1 - \ell_1} \wedge  \bbk \right)  \wedge  \\  C^{*,c}_{\rel \partial}\left(B \hat{Z}_{2};  S^{d_{I_2}-V_2 -\ell_2} \wedge \bbk \right)   
       \end{aligned}
      \ar[r]  &  C^{*,c}_{\rel \partial}\left(B \hat{Z};  S^{d_I -V - \ell} \wedge \bbk  \right) \\
      \begin{aligned}
      C^{*,c}_{\rel \partial}(B\hat{Z}_{1}; S^{d_{I_1} - d_{V_1} -
        \ell_1} \wedge \bbk) \wedge \\  C^{*,c}_{\rel
        \partial}(B\hat{Z}_{2};  S^{d_{I_2} - d_{V_2} -\ell_2} \wedge
      \bbk)  
      \end{aligned}
        \ar[r] \ar[u] & C^{*,c}_{\rel \partial}(B\hat{Z};  S^{d_I - d_V - \ell}  \wedge \bbk) \ar[u]
    \end{tikzcd}
\end{equation}
in which the righthand vertical arrows are equivalences and the
lefthand vertical arrows induce equivalences of derived smash
products.  Using the
multiplicativity of the spheres $\bS[n]$ and the comparisons of
Equation~\eqref{eq:cancel-spheres}, we obtain a commutative diagram 
\begin{equation} 
    \begin{tikzcd}[column sep=tiny, transform shape, nodes={scale=.95}]
      \begin{gathered}
      C^{*,c}_{\rel \partial}\left(B \hat{Z}_{1}; S^{I_1-V_1 - \ell_1} \wedge  \bbk \right)\textstyle{[\frac{d_{V_1} - d_{I_1}}{2}]}  \wedge \\   C^{*,c}_{\rel \partial}\left(B \hat{Z}_{2};  S^{I_2-V_2 -\ell_2} \wedge \bbk \right)\textstyle{[\frac{d_{V_2} - d_{I_2}}{2}] } 
      \end{gathered}
       \ar[r] \ar[d] &  C^{*,c}_{\rel \partial}\left(B \hat{Z};  S^{I-V - \ell} \wedge \bbk  \right)[\frac{d_V - d_I}{2}] \ar[d] \\
      \cdots \ar[r] & \cdots \\      
      \begin{aligned}
      C^{*,c}_{\rel \partial}(B\hat{Z}_{1}; \Omega^{ \ell_1} \bbk)
      \\ \wedge  C^{*,c}_{\rel \partial}(B\hat{Z}_{2};  \Omega^{\ell_2}
      \bbk)  
      \end{aligned}
        \ar[r] \ar[u] & C^{*,c}_{\rel \partial}(B\hat{Z};  \Omega^\ell \bbk).  \ar[u]
    \end{tikzcd}
\end{equation}

We will study the associativity of this diagram (and rely on coherent
compatibility conditions coming from a monoidal structure on the
category of charts) in Section \ref{sec:compl-orient-flow-1} below.

\subsubsection{Charts equipped with (relative) complex orientations}
\label{sec:charts-equipped-with}

In order to describe relative complex orientations, we start with a 
pair of stable vector spaces 
\begin{equation}
V_p = (V^+_p, V^-_p) \quad\textrm{and}\quad V_q = (V^+_q, V^-_q)
\end{equation}
such that $V^-_p$ and $V^-_q$ are equipped with complex structures
(the reason for this requirement will be clear once we discuss the
multiplicative notion later in this section).

\begin{defin}
A \emph{stable complex lift} of the tangent space of a flag smooth
$\langle \cS \rangle$-Kuranishi chart $(X,V,s,G, B, \pi)$ relative to
$V_p$ and $V_q$ and a real line $\ell_q$,  consists of the following
data:
\begin{enumerate}
\item A complex vector bundle on $I$ over $Z$ and a complex
$G$-representation $W$, 
\item a complex structure on $V$, compatible
with the $G$-action and the inner product, 
\item and a $G$-equivariant real isomorphism of vector bundles over
  $Z$: 

\begin{equation}\label{eq:stable-isomorphism-relative-complex}
V^-_p \oplus  \ell_{q} \oplus  T X \oplus    W  \oplus  V^+_q  \cong
V^+_p \oplus I \oplus    W \oplus V^-_q.
\end{equation}
\end{enumerate}
\end{defin}

\begin{rem}
The definition is modeled after the intended application in Floer
theory, see in particular
Equation~\eqref{eq:oriented-charts-family-iso}
(and also
c.f. Definition~\ref{def:Complex-oriented-Kuranishi-chart}). The point
is that the tangent space of moduli spaces is naturally isomorphic to
the linearisation of a Cauchy-Riemann operator with inhomogenous terms
that are not complex linear. Gluing operators associated to the
cylindrical ends yields an operator that admits a deformation to a
complex linear operator, canonically up to contractible choice
\cite{FloerHofer1993,WehrheimWoodward2015,Seidel2008a}. We have to
work with stable vector spaces $V_p$ and $V_q$ (rather than ordinary
ones) because the operators associated to the cylindrical ends may not
be surjective. 
\end{rem}

We introduce the notation:  
\begin{align} \label{eq:morphisms-categories-TX-I-and-bimodule}
 \scrT(\bX) & \equiv  F_{BZ}( S^{V}, (S^{TX})^{\mfib}) \wedge S^{V_{p}^{-} +  V^+_q - (  V^+_p + V_q^-) }  \\
  \scrP(\bX) & \equiv  F_{BZ}(S^{V^+_p  + \ell_{q} + W + V + V_q^-}, (S^{V^-_p + \ell_{q} + T X + W + V^+_q})^{\mfib})\\
  & \nonumber \cong  F_{BZ}(S^{ V^+_p  + \ell_{q} + W + V + V_q^- }, (S^{V^+_p + I + W + V^-_q })^{\mfib}) \\
  \scrI(\bX) & \equiv  F_{BZ}(S^{V+ \ell_q}, (S^{I})^{\mfib}),
\end{align}
for the induced parametrized spectra over $BZ$, where we omit the
additional data required to formulate a stable lift from the
notation.  Here also note that the notation $S^{V_{p}^{-} +  V^+_q - (
  V^+_p + V_q^-) }$ (and subsequent manipulations of these terms)
involves the standard negative spheres $S^{-V} = F_V S^0$ in
orthogonal $G$-spectra (see Section~\ref{sec:negative-MM-spheres}). 
In this context, $(-)^{\mfib}$ denotes the fiberwise multiplicative
fibrant replacement functor (see
Definition~\ref{def:lax-fibrant-param}).

We have natural $G$-equivariant fiberwise equivalences
\begin{equation}\label{eq:fiberwise-comparisons}
  \begin{tikzcd}
    \scrT(\bX)    \ar[r] &  F_{BZ}(S^{V + V^+_p + V_q^-}, (S^{V_{p}^{-} + TX + V^+_q})^{\mfib})  \ar[d, "\wedge S^{\ell_{q} + W}  "]  \\
\scrP(\bX) \ar[r,"="] &  F_{BZ}(S^{V^+_p  + \ell_{q} + W + V + V_q^-}, (S^{V^-_p + \ell_{q} + T X + W + V^+_q})^{\mfib})  \\
 \scrI(\bX)     \ar[r, "\wedge S^{  V^+_p + V_q^-  }   "] &  F_{BZ}(S^{ V^+_p  + \ell_{q} + W + V + V_q^- }, (S^{V^+_p + I + W + V^-_q })^{\mfib})   \ar[u,"\cong"]
  \end{tikzcd}
\end{equation}
of parametrized spectra over $BZ$.

Passing to spectra of compactly supported sections, the equivalences
of Equation~\eqref{eq:fiberwise-comparisons} give rise to a zig-zag of
equivalences  
\begin{equation}
  \begin{tikzcd}
 C^{*,c}_{\rel \partial}\left(B Z; \scrI(\bX) \right)  \ar[r] &
 C^{*,c}_{\rel \partial}\left(B Z; \scrP(\bX) \right)   &  
 C^{*,c}_{\rel \partial}\left(B Z; \scrT(\bX) \right)   \ar[l]
  \end{tikzcd}
\end{equation}

Given a stratum of $\partial \bX$ of $\bX$,  with normal bundle
$\kappa$, we set 
\begin{align}
 \scrT(\partial \bX)  & \equiv  F_{B \partial
\hat{Z}}(S^{  V }, (S^{T \partial X})^{\mfib})  \wedge S^{V_{p}^{-} + V^+_q - (V^+_p + V_q^-)}    \\
  \scrP(\partial \bX) & \equiv   F_{B \partial
\hat{Z}}(S^{V^+_p  + \ell_{q} + \kappa +  W + V + V_q^-}, (S^{V^-_p + \ell_{q} + T X + W + V^+_q})^{\mfib}) \\
  & \nonumber \cong   F_{B \partial
\hat{Z}}(S^{ V^+_p  + \ell_{q} + \kappa + W + V + V_q^- }, (S^{V^+_p + I + W + V^-_q })^{\mfib}) \\
  \scrI(\partial \bX) & \equiv F_{B \partial
\hat{Z}}(S^{V + \ell_q + \kappa}, (S^I)^\mfib), 
\end{align}
where we are implicitly using the isomorphism $S^{T \partial \hat{X} +
  \kappa} \to S^{T\hat{X}}$ of parametrized spectra over $B \partial
\hat{Z}$. Note that these definitions are slightly awkward because the
definition of these spectra uses the ambient chart $\bX$. 

\begin{lem}\label{lem:boundary-pre-bimodule-maps}
We have a commutative diagram
\begin{equation}
  \begin{tikzcd}
     \scrT(\partial \bX) \ar[r] \ar[d] &  \scrP(\partial \bX) \ar[d] &   \scrI(\partial \bX) \ar[l] \ar[d]  \\
  F(S^\kappa ,  \scrT( \bX))  \ar[r] &    F(S^\kappa , \scrP( \bX)) &    F(S^\kappa ,  \scrI( \bX)) \ar[l],
  \end{tikzcd}
\end{equation}
where the vertical maps are induced by smashing with $S^\kappa$.\qed
\end{lem}

Passing to spectra of compactly supported sections yields the
following comparison diagram.

\begin{lem}
The maps of Lemma~\ref{lem:boundary-pre-bimodule-maps} induce
commutative diagrams:
\begin{equation}
\begin{tikzcd}
C^{*,c}_{\rel \partial}\left(B \partial \hat{Z};   \scrT(\partial \bX)\right)  \ar[r] \ar[d] & C^{*,c}_{\rel \partial}\left(B \hat{Z};  \scrT( \bX) \right)   \ar[d] \\
C^{*,c}_{\rel \partial}\left(B \partial \hat{Z}; \scrP(\partial \bX) \right)  \ar[r]  &    C^{*,c}_{\rel \partial}\left(B \hat{Z};  \scrP( \bX)  \right) \\
C^{*,c}_{\rel \partial}\left(B\partial \hat{Z}; \scrI(\partial \bX) \right) \ar[u]  \ar[r] & C^{*,c}_{\rel \partial}\left(B \hat{Z}; \scrI( \bX) \right) \ar[u].
\end{tikzcd}
\end{equation} \qed
\end{lem}

\subsubsection{Multiplicativity of relative orientations}
\label{sec:mult-relat-orient}  

We now begin to discuss multiplicativity of relative orientations.  In
contrast to the previous discussion, where we constructed comparison
zig-zags of multiplicative functors, we shall ultimately perform the
multiplicative comparison using a bimodule which represents an
equivalence.  Here we set up the action maps without considering the
full bimodule structure.

To start, let 
\begin{equation}
V_p = (V^+_p,V^-_p) \quad\textrm{and}\quad V_q = (V^+_q,V^-_q)
\quad\textrm{and}\quad V_r = (V^+_r,V^-_r)
\end{equation}
be a triple of stable vector spaces, and consider a triple of data
\begin{equation}
\{(\bX_i, I_i, W_i, \ell_i)\}_{i = 0,1} \quad\textrm{and}\quad (\bX,
I, W, \ell)
\end{equation}
consisting of flag smooth Kuranishi charts equipped 
with stable complex vector bundles, and associated lifts of the
tangent space relative $V_p $ and  $V_q $  for $\bX_1$,   $V_q$ and
$V_r$ for $\bX_2$, and  $V_p$ and  $V_r$  for $\bX$.  We assume that
$\bX_1 \times \bX_2$ is a boundary stratum of $\bX$ and further that
we are given complex isomorphisms 
\begin{align}
I_1 \oplus I_2 & \cong I \\
W_1 \oplus V^-_q \oplus W_2 & \cong W.
\end{align}

We express the necessary compatibility of the lifts as follows.

\begin{defin}
The complex lifts of the tangent spaces of $\bX_i$ are
\emph{compatible with} the complex lift of the tangent space of $\bX$
if the following diagram commutes: 
\begin{equation} \label{eq:compatible_complex-orientation}
\begin{tikzcd}[ column sep=small]
  \begin{gathered}
    V^-_p \oplus  \ell_{q} \oplus  T X_1 \oplus    W_1  \oplus \\
    V^-_q \oplus  \ell_{r} \oplus  T X_2 \oplus    W_2  \oplus  V^+_r 
  \end{gathered}
 \ar[r] \ar[d] &  V^-_p \oplus   T X \oplus    W   \oplus  \ell_{r} \oplus  V^+_r  \ar[d] \\
 \begin{gathered}
   V^-_p \oplus  \ell_{q} \oplus  T X_1 \oplus    W_1  \oplus  \\
   V^+_q \oplus  I_2 \oplus W_2  \oplus  V^-_r
 \end{gathered}
  \ar[d]   & V^+_p \oplus I \oplus    W  \oplus  V^-_r. \\
     V^+_p \oplus I_1 \oplus    W_1 \oplus V^-_q \oplus  I_2 \oplus W_2  \oplus  V^-_r \ar[ur]        & 
\end{tikzcd}\end{equation}
\end{defin}

In order to state the consequence of this compatibility at the level
of spectra, we use the notation from Equation
\eqref{eq:morphisms-categories-TX-I-and-bimodule}: 
\begin{lem} \label{lem:toy-case-bimodule-comparison-tangent}
Assume that the complex lifts of the tangent spaces of $\bX_i$ are
compatible with the complex lift of the tangent space of $\bX$.  Then
there are natural maps 
\begin{align}
   \scrT(\bX_1) \wedge \scrT(\bX_2) & \to F(S^{\kappa}, \scrT(\bX)) \\
   \scrI(\bX_1) \wedge \scrI(\bX_2) & \to F(S^{\kappa}, \scrI(\bX))  \\
  \scrP(\bX_1) \wedge  \scrT(\bX_2) & \to  F(S^{\kappa}, \scrP(\bX)) \\
  \scrI(\bX_1) \wedge \scrP(\bX_2) & \to F(S^{\kappa}, \scrP(\bX)),
\end{align}
of parametrized spectra over $B \hat{Z}_1 \times B \hat{Z}_2$ such
that the following three diagrams commute: 
\begin{equation} \label{eq:left-module-complex-orientation}
  \begin{tikzcd}
    \scrI(\bX_1) \wedge \scrI(\bX_2) \ar[d] \ar[r]  &  \scrI(\bX_1) \wedge \scrP(\bX_2) \ar[d]  \\
    F(S^{\kappa}, \scrI(\bX)) \ar[r] &  F(S^{\kappa}, \scrP(\bX))
  \end{tikzcd}
\end{equation}
\begin{equation} \label{eq:right-module-complex-orientation}
  \begin{tikzcd}
 \scrT(\bX_1) \wedge \scrT(\bX_2) \ar[d]  \ar[r]  &      \scrP(\bX_1) \wedge \scrT(\bX_2)  \ar[d] \\
 F(S^{\kappa}, \scrT(\bX)) \ar[r]  & F(S^{\kappa}, \scrP(\bX)) 
  \end{tikzcd}
\end{equation}
\begin{equation} \label{eq:bimodule-complex-orientation}
  \begin{tikzcd}
    \scrI(\bX_1) \wedge \scrT(\bX_2) \ar[d] \ar[r]  &  \scrI(\bX_1) \wedge \scrP(\bX_2) \ar[d]  \\
   \scrP(\bX_1) \wedge \scrT(\bX_2)  \ar[r] &  F(S^{\kappa}, \scrP(\bX))
  \end{tikzcd}
\end{equation}
\end{lem}
\begin{proof}
The first two maps are entirely straightforward to construct. The
third is defined as the composition 
\begin{equation}
\begin{tikzcd}
  \begin{aligned}
    F_{B \hat{Z}_1}(S^{V^+_p  + \ell_{q} + W_1 + V_1 + V_q^-}, (S^{V^-_p + \ell_{q} + T X_1 + W_1 + V^+_q})^{\mfib}) \\
    \wedge  F_{B \hat{Z}_2}( S^{V_2}, (S^{TX_2})^{\mfib}) \wedge S^{V_{r}^{-} +  V^+_p - (V^+_r + V_p^-) } 
  \end{aligned}
    \ar[d] \\
    \begin{aligned}
      F_{B \hat{Z}_1 \times B \hat{Z}_2}(S^{V^+_r  + \ell_{q} + W_1 + V_1 + V_2 + V_q^-}, \\
      (S^{V^-_p + \ell_{q} + T X_1 + TX_2+ W_1 + V^+_q}  \wedge S^{V_{r}^{-} - V^-_p} )^{\mfib}) 
    \end{aligned}
   \ar[d] \\
F_{B \hat{Z}_1 \times B \hat{Z}_2}(S^{V^+_r  + \ell_{q} + W_1 + V_{12} + V_q^-}, (S^{V^-_r + \ell_{q} + T (X_{1} \times X_{2}) + W_1 + V^+_q} )^{\mfib})   \ar[d,"\wedge S^{W_2 + \kappa} "] \\
F_{B \hat{Z}_1 \times B \hat{Z}_2}(S^{\kappa}, F(S^{V^+_r  + V_{12} +  W_{12} + V_q^-}, S^{V^-_r + \ell_{q} + TX + W_{12} + V^+_q })),
\end{tikzcd}
\end{equation}
and the fourth map is given by
\begin{equation}
  \begin{tikzcd}[row sep = small, column sep=tiny, transform shape, nodes={scale=.9}]
    F_{B\hat{Z}_1}(S^{V_1+ \ell_q}, (S^{I_1})^\mfib)  \wedge F_{B\hat{Z}_2}(S^{ V^+_r  + \ell_{p} + W_2 + V_2 + V_p^- }, (S^{V^+_r + I_2 + W_2 + V^-_p })^{\mfib})
  \ar[d] \\
 F_{B \hat{Z}_1 \times B \hat{Z}_2}(S^{ V^+_r + \ell_q + \ell_{p} + W_2 + V_{12} + V_p^- }, (S^{V^+_r + I_{12} + W_2 + V^-_p })^{\mfib})  \ar[d, "\wedge S^{W_1}"] \\
F_{B \hat{Z}_1 \times B \hat{Z}_2}(S^{\ell_p},  F(S^{ V^+_r + \ell_q + \ell_{p} + W_{12} + V_{12} + V_p^- }, (S^{V^+_r + I_{12} + W_{12} + V^-_p })^{\mfib})),
  \end{tikzcd}
\end{equation}
where we use the isomorphism $\kappa \cong \ell_p$ to get the desired statement.

The commutativity of Diagrams
\eqref{eq:left-module-complex-orientation} and
\eqref{eq:right-module-complex-orientation} are easy to check. The
commutativity of Diagram \eqref{eq:bimodule-complex-orientation}
follows from Equation \eqref{eq:compatible_complex-orientation}, which
directly implies the following diagram commutes: 
\begin{equation}
\begin{tikzcd}
\begin{aligned}
F_{B\hat{Z}_1}(S^{V_1 + \ell_q}, (S^{I_1})^\mfib) \\ \wedge  F_{B\hat{Z}_2}( S^{V_2},
(S^{TX_2})^{\mfib}) \\ \wedge S^{V_{r}^{-} +  V^+_p - (V^+_r + V_p^-) }
\end{aligned}
\ar[r] \ar[d] 
&  
\begin{aligned}
F_{B\hat{Z}_1}(S^{V_1 + \ell_q}, (S^{I_1})^\mfib )\\ \wedge
F_{B\hat{Z}_2}(S^{ V^+_r  + \ell_{p} + W_2 + V_2 + V_p^- },
\\ (S^{V^+_r + I_2 + W_2 + V^-_p})^{\mfib}) 
\end{aligned} \ar[d] \\
\begin{aligned}
F_{B\hat{Z}_1}(S^{V^+_p  + \ell_{q} + W_1 + V_1 + V_q^-},\\ (S^{V^-_p + \ell_{q} + T
  X_1 + W_1 + V^+_q})^{\mfib})  \\ \wedge  F_{B\hat{Z}_2}( S^{V_2},
(S^{TX_2})^{\mfib}) \\ \wedge S^{V_{r}^{-} +  V^+_p - (V^+_r + V_p^-) }
\end{aligned} \ar[r] 
&  
\begin{aligned}
F_{B \hat{Z}_1 \times B \hat{Z}_2}(S^{\kappa}, F(S^{V^+_r  + V_{12} +
  W_{12} + V_q^-}, \\ S^{V^-_r + \ell_{q} + TX + W_{12} + V^+_q })). 
\end{aligned}
\end{tikzcd}
\end{equation}
\end{proof}

Note that the three diagrams above are very similar to the diagrams
expressing the structure of $\scrI$ and $\scrT$ as spectral
categories, a $\scrI$-$\scrT$-bimodule structure on $\scrP$, and the maps $\scrI
\to \scrP$ and $\scrT \to \scrP$ as bimodule maps, except for the
fact that we see a desuspension by $\ell_q$.  Identifying $\ell_q$
with the collar direction (as in
Equation~\eqref{eqn:compactly-supported-boundary-map}), we immediately
obtain maps expressing these structures (except for associativity).

\begin{cor}
The data of an isomorphism of the product $\bX_1 \times \bX_2$ of a
pair of flag smooth Kuranishi charts with a boundary stratum of a
Kuranishi chart $\bX$ induces maps 
  \begin{align}
    C^{*,c}_{\rel \partial}( B \hat{Z}_1,   \scrT(\bX_1) )  \wedge C^{*,c}_{\rel \partial}( B \hat{Z}_2,  \scrT(\bX_2))  & \to C^{*,c}_{\rel \partial}( B \hat{Z},   \scrT(\bX)) \\
    C^{*,c}_{\rel \partial}( B \hat{Z}_1,   \scrI(\bX_1) )  \wedge C^{*,c}_{\rel \partial}( B \hat{Z}_2,  \scrI(\bX_2))  & \to C^{*,c}_{\rel \partial}( B \hat{Z},   \scrI(\bX)) \\    
C^{*,c}_{\rel \partial}( B \hat{Z}_1,   \scrP(\bX_1) )  \wedge C^{*,c}_{\rel \partial}( B \hat{Z}_2,  \scrT(\bX_2))  & \to C^{*,c}_{\rel \partial}( B \hat{Z},   \scrP(\bX)) \\
C^{*,c}_{\rel \partial}( B \hat{Z}_1,    \scrI(\bX_1)   ) \wedge C^{*,c}_{\rel \partial}( B \hat{Z}_2,  \scrP(\bX_2)) & \to C^{*,c}_{\rel \partial}( B \hat{Z},   \scrP(\bX) ),
\end{align}
such that the following three diagrams commute:
\begin{equation}
 \begin{tikzcd}
   \begin{gathered}
     C^{*,c}_{\rel \partial}( B \hat{Z}_1,  \scrI(\bX_1) ) \wedge \\
     C^{*,c}_{\rel \partial}( B \hat{Z}_2,  \scrI(\bX_2) )
   \end{gathered}
   \ar[d] \ar[r]  &
   \begin{gathered}
     C^{*,c}_{\rel \partial}( B \hat{Z}_1,  \scrI(\bX_1)) \wedge \\
     C^{*,c}_{\rel \partial}( B \hat{Z}_2,  \scrP(\bX_2))
   \end{gathered}
 \ar[d]  \\
C^{*,c}_{\rel \partial}( B \hat{Z}, \scrI(\bX)   ) \ar[r] &  C^{*,c}_{\rel \partial}( B \hat{Z}, \scrP(\bX))
  \end{tikzcd}
\end{equation}
\begin{equation}
 \begin{tikzcd}
   \begin{gathered}
     C^{*,c}_{\rel \partial}( B \hat{Z}_1,   \scrT(\bX_1)) \wedge \\
     C^{*,c}_{\rel \partial}( B \hat{Z}_2,  \scrT(\bX_2))
   \end{gathered}
\ar[d] \ar[r] &
\begin{gathered}
  C^{*,c}_{\rel \partial}( B \hat{Z}_1,  \scrP(\bX_1) ) \wedge \\ C^{*,c}_{\rel \partial}( B \hat{Z}_2,  \scrT(\bX_2)) 
\end{gathered}
 \ar[d]  \\
C^{*,c}_{\rel \partial}( B \hat{Z},  \scrT(\bX) )  \ar[r] & C^{*,c}_{\rel \partial}( B \hat{Z},  \scrP(\bX))
  \end{tikzcd}
\end{equation}
\begin{equation}
 \begin{tikzcd}
   \begin{gathered}
     C^{*,c}_{\rel \partial}( B \hat{Z}_1,     \scrI(\bX_1) )   \wedge \\
     C^{*,c}_{\rel \partial}( B \hat{Z}_2, \scrT(\bX_2))
   \end{gathered}
  \ar[d] \ar[r]  &
  \begin{gathered}
    C^{*,c}_{\rel \partial}( B \hat{Z}_1, \scrI(\bX_1))  \wedge \\ C^{*,c}_{\rel \partial}( B \hat{Z}_2, \scrP(\bX_2) )
  \end{gathered}
 \ar[d]  \\
    C^{*,c}_{\rel \partial}( B \hat{Z}_1,  \scrP(\bX_1))  \wedge  C^{*,c}_{\rel \partial}( B \hat{Z}_2, \scrT(\bX_2))   \ar[r] &  C^{*,c}_{\rel \partial}( B \hat{Z}, \scrP(\bX))
  \end{tikzcd}
\end{equation}
\qed
\end{cor}

\subsection{Signpost: Construction of the homotopy type for global charts}
\label{sec:signp-constr-homot}

We summarise the constructions of this section in the case of a global, closed flag smooth Kuranishi chart, equipped with a stable complex lift of its spherical tangent fibration. The outcome of the previous sections is that we have a zig-zag of equivalences
\begin{equation}
  X|Z^{-V-d} \wedge \bbk  \leftarrow \cdots  \to   C^*(BZ;   S^{TX -V-d} \wedge \bbk)  \leftarrow \cdots  \to  C^*(BZ; \Omega^\ell \bbk),
\end{equation}
where the first zig-zag arises from the Spanier-Whitehead duality map, and the equivalences of different models of the tangent spherical fibration, and the second arises from orientation data.

If we compose these equivalences on the right side with the map 
\begin{align}
  \Omega^\ell \bbk & \to   C^*(BZ; \Omega^\ell \bbk), \\
\end{align}
and compose with the map
\begin{equation}
  X|Z^{-V} \wedge \bbk  \to \Omega^d \bbk
\end{equation}
on the left, we obtain a zig-zag representing a map
\begin{equation}
\Omega^\ell \bbk  \to \Omega^d \bbk
\end{equation}
whose homotopy cofiber is the homotopy type associated to this chart.

\section{Coherent comparisons: support and duality}
\label{sec:local-glob-constr}

The goal of this section is to extend the constructions of the first
half of Section~\ref{sec:comp-coch-kuran} to the level of Kuranishi
presentations by describing their functorial properties.
Specifically, by checking that the behaviour of boundary and product
maps are functorial, we extract spectral categories which interpolate
between virtual cochains and Milnor-twisted cochains.

\subsection{Comparing compactly supported and ordinary cochains}
\label{sec:comp-comp-supp}

\subsubsection{Orbispace flow categories with unique factorisation}
\label{sec:relat-coch-kuran}

By passing to zero-loci, a Kuranishi flow category $\bX \co A \to
\Kur$ determines an orbispace flow category, and hence a pair of
topological flow categories $B \cZ$ and $ \cM$, which are respectively obtained by taking the
Borel construction and the quotient of charts. These categories are equipped with
a $\Pi$-equivariant functor $B \cZ \to \cM$, as explained in Section
\ref{sec:orbisp-flow-categ}.  Applying the construction of Section
\ref{sec:relat-coch-coll}, we obtain a spectrally enriched category of
relative cochains which we denote
\begin{equation}
  C^*_{\rel \partial}(\widehat{B \cZ}; \Omega \bbk).
\end{equation}
Recall that the morphisms in this category are defined in Equation \eqref{eq:definition_relative_cochains_Borel-naive} using the stratification of the morphism spaces $BZ(p,q)$ arising from the compositions 
\begin{equation}
  BZ(p,q) \times BZ(q,r) \to BZ(p,r).
\end{equation}

In the construction of virtual cochains, we instead used the
stratification induced by the inverse image of strata in the original
flow category $\cM$. This section explains a procedure for
constructing a category built from relative cochains with respect to
this geometric stratification.

The fundamental problem for defining a category whose morphisms are
\begin{equation}
    C^*( B \hat{Z}(p,q), B \partial \hat{Z}(p,q); \Omega \bbk)     
  \end{equation}
  is that the map
  \begin{equation}
         BZ(p,q) \times BZ(q,r) \to B \partial^q Z(p,r). 
  \end{equation}
  need not be a homeomorphism. One way to ensure that this condition holds is to consider the following variant of orbispace flow categories:
  \begin{defin} \label{def:collared_orbispace_flow_cat} An
    \emph{orbispace flow category with collars} consists of a
    $\Pi$-equivariant $2$-category $A$ over $\cP$ and a
    $\Pi$-equivariant strict $2$-functor $\hat{A} \to \Orb$, so that
    we have orbispace presentations 
\begin{equation}
(Z, G) \co \hat{A}(p,q) \to \Orb \langle \cP(p,q) \rangle
\end{equation}
such that the following properties hold:
\begin{enumerate}
\item The functors $\bX$ are $\Pi$-equivariant in the sense that the
  following diagram commutes:
      \begin{equation}
        \begin{tikzcd}
          \hat{A}(p,q) \ar[d] \ar[r] & \Orb \langle \cP(p,q) \rangle \ar[d] \\
          \hat{A}(\pi \cdot p, \pi \cdot q) \ar[r] &  \Orb \langle \cP(\pi \cdot p, \pi \cdot q) \rangle
        \end{tikzcd}
      \end{equation}
    \item For each triple $(p,q,r)$ the following diagram commutes:
       \begin{equation}
        \begin{tikzcd}
          \hat{A}(p,q) \times \hat{A}(q,r) \ar[d] \ar[r] & \Orb \langle \cP(p,q) \rangle \times \Orb \langle \cP(q,r) \rangle \ar[d] \\
          \hat{A}(p,r) \ar[r] &  \Orb \langle \cP(p,r)  \rangle,
        \end{tikzcd}
      \end{equation}
      where the right vertical map is given by taking the product with a collar labelled by $q$.
    \end{enumerate}
    We say that this flow category has \emph{unique factorisation} if the following additional condition holds:
    \begin{equation}
       \label{eq:boundatry-strata-are-uniquely-products}
\parbox{31em}{each component of $\partial^q Z_\alpha$ for $\alpha \in \hat{A}(p,r)$ is the homeomorphic image of a single component of $Z_{\alpha_-} \times  Z_{\alpha_+}$, for uniquely determined elements $\alpha_+  \in \hat{A}(p,q)$ and $\alpha_- \in \hat{A}(q,r)$.}
    \end{equation}
  \end{defin}
  \begin{rem}
We can define a notion of Kuranishi flow category with collars
replacing $\Orb$ by $\Kur$ everywhere in the above definition, but we
shall not use this notion. The main reason for avoiding it is that it
is awkward to formulate flag smoothness (see Section
\ref{sec:fiber-kuran-pres}) for Kuranishi flow categories with
collars.
\end{rem}

We can associate to an orbispace flow category with collars a
$\Pi$-equivariant topologically enriched category $BZ$ with morphism
spaces  
\begin{equation}
BZ(p,q) \equiv \hocolim_{\alpha \in \hat{A}(p,q)} BZ_\alpha
\end{equation}
as before.  The following result asserts that
Condition~\eqref{eq:boundatry-strata-are-uniquely-products} ensures
that the two stratifications of the Borel construction agree: 
\begin{lem}\label{lem:stratification_of_Borel_agree}
If an orbispace flow category with collars has unique factorisation,
then for each totally ordered subset $P$ with minimum $p$ and maximum
$q$, the inverse image of $\partial^P \cM(p,q)$ in $B\cZ(p,q)$ agrees
with $\partial^P B \cZ(p,q) $. 
\end{lem}

\begin{proof}
The fact that the image of $ \partial^P B \cZ(p,q)  $ lies in $
\partial^P \cM(p,q)$ follows from the construction. The reverse
inclusion is an immediate consequence of Condition
\eqref{eq:boundatry-strata-are-uniquely-products} that all strata are
products. Indeed, the inverse image of $\partial^P \cM(p,q)$ in
$B\cZ(p,q)$ is obtained as a homotopy colimit of all charts that
intersect $ \partial^P \cM(p,q)$, and these are all product charts by
our definition of an orbispace flow category with collars.
\end{proof}

It is clear that passing from a flow category to its collared
completion yields the data of an orbispace flow category with collars.  In
the remainder of this section, we explain how to associate to any
orbispace flow category $(Z, G) \co A \to \Orb$ an orbispace flow category
with collars, so that Condition
\eqref{eq:boundatry-strata-are-uniquely-products} holds as well. 

To begin, we define a new bicategory $\vec A$ with the $0$-cells again
given by the elements of $\cP$, and with $1$-cells given by categories
$\vec A(p,q)$ that we now define.  

\begin{defin}
For $p,q \in \cP$, the category $\vec A(p,q)$ has 
\begin{itemize}
\item objects given composable sequences $\vec \alpha$ with source $p$ and target $q$, i.e. a sequence
\begin{equation}
\vec \alpha = (\alpha_0, \ldots, \alpha_{k-1}),
\end{equation}
with $\alpha_i \in
A(p_i,p_{i+1})$ for some totally ordered subset $P \subseteq \cP(p,q)$, where we set $p_0 = p$ and $p_k = q$.
\item A morphism
\begin{equation}
\vec f \co \vec \alpha \to \vec \beta
\end{equation}
is specified as follows:  writing $P = (p_1, \ldots, p_k)$ and $Q = (q_1, \ldots, q_\ell)$ for the sequences of objects associated to $\vec \alpha$ and $\vec \beta$, we assume that we have an inclusion $Q \subset P$, which induces a decomposition of $\cP$ into subsets $\{P_j\}$
where $P_j = \{p \in \cP \, \mid \, q_j \leq p < q_{j+1}\}$.  We will
write $P_j = \{p_{j_1}, p_{j_2}, \ldots, p_{j_m}\}$ in what
follows.  The morphism is then specified by morphisms for each $j$
\begin{equation}
f_j \colon \mu(\alpha_{j_1} \times \alpha_{j_2} \ldots \times \alpha_{j_m}) \to \beta_j
\end{equation}
in $A(q_j,q_{j+1})$, where $\mu$ denotes the composition
\begin{equation}
A(q_j, p_{j_1}) \times A(p_{j_1}, p_{j_2}) \times \ldots \times 
A(p_{j_m},q_{j+1}) \to A(q_j, q_{j+1}).
\end{equation}
\end{itemize} 
\end{defin}
Note that there is a natural functor
\begin{equation}
  \vec A(p,q) \to A(p,q),
\end{equation}
which assigns to each sequence of objects $\vec \alpha$ their product
which we denote $\mu(\vec \alpha)$.  We now assemble these into a
$2$-category structure:  
\begin{lemma}
For $p,q,r \in \cP$, there are strictly associative functors
\begin{equation}
\vec A(p,q) \times   \vec A(q,r) \to \vec A(p,r),
\end{equation}
so that the following diagram commutes:
\begin{equation}
  \begin{tikzcd}
    \vec A(p,q) \times   \vec A(q,r) \ar[d] \ar[r] &  \vec A(p,r) \ar[d] \\
    A(p,q) \times    A(q,r)  \ar[r] &   A(p,r).
  \end{tikzcd}
\end{equation}
\end{lemma}
\begin{proof}
The functors in question are specified on objects by the assignment
 \begin{equation}
\vec \alpha \times \vec \beta \mapsto
  \vec \alpha \times \vec \beta. 
\end{equation}
On morphisms, given inclusions $Q_1 \subseteq P_1$ and $Q_2 \subseteq
P_2$, we have an induced inclusion $Q_1 \amalg \{ q \} \amalg Q_2 \subseteq P_1
\amalg \{q\} \amalg P_2$.  The induced partition of $P_1 \amalg \{q\}
\amalg P_2$ is the disjoint union of the partitions induced on $P_1$ and $P_2$,
and so we take the product of the collections of morphisms
$\{f^1_j\}$ and $\{f^2_j\}$. 
\end{proof}

To express the functoriality of the relative cochains with respect to
this construction, we need to use the twisted arrow category (see
Section~\ref{sec:twisted}).  Given an orbispace flow category $(Z,
G)$, we will replace our indexing $2$-category $A$ by applying the
twisted arrow category construction to each morphism category;
Lemma~\ref{lem:morphism-wise-twisting} shows that this produces a new
$2$-category.  Specifically, we will use the $\Pi$-equivariant
bicategory $\Tw \vec A$ with morphism categories given by $\Tw \vec A(p,q)$.
Note that this $2$-category is equipped with a natural $\Pi$-equivariant
$2$-functor
\begin{equation}
  \Tw \vec A \to A.
\end{equation}
We now define our replacement functor.

\begin{prop} \label{prop:collared-presentation-from-usual-presentation}
Given an orbispace flow category $(Z, G) \co A \to \Orb$, there is a
collared orbispace flow category with unique factorisation, given by functors 
\begin{equation}
(\vec{Z}, \vec G) \co  \Tw \vec A(p,q) \to \Orb \langle \cP(p,q) \rangle,   
\end{equation}
and a natural functor
\begin{equation} \label{eq:compare_cochains_to_collared_cochains}
C^*_{\rel \partial}(B \cZ; \Omega \bbk) \to C^*_{\rel \partial}(B \vec{\cZ}; \Omega \bbk)
\end{equation}
which is an equivalence of spectrally enriched, $\Pi$-equivariant categories.
\end{prop}
\begin{proof}
  We begin by constructing the functor $ (\vec{Z}, \vec G)$.  Using
  the composite functor 
\begin{equation}
\Tw \vec A(p,q) \to \vec A(p,q) \to A(p,q)
\end{equation}
we obtain a product chart
 $(\hat{Z}_{\mu \vec  \alpha}, \hat{G}_{\mu \vec \alpha})$, considered as an object of $\Orb \langle \cP(p,q) \rangle$ via the
external monoidal product. In particular, $\hat{Z}_{ \mu \vec  \alpha} $ is the collared completion of the product
\begin{equation}
Z_{\alpha_0} \times Z_{\alpha_1} \times \ldots \times Z_{\alpha_{k-1}}.
\end{equation}

Now assume that we are given an arrow $\vec f \co \vec \alpha \to \vec \beta$.  We define an object $\vec{Z}_{\vec{f}}$ of $\Orb \langle
\cP(p,q) \rangle$ by taking the union of the interiors of the strata
of $\hat{Z}_{\mu \vec \alpha} $ which are labeled by totally ordered
subsets of $Q$
\begin{equation}
\vec{Z}_{\vec{f}} \equiv \bigcup_{R \subset Q} \mathrm{Int}(\partial^R
\hat{Z}_{\mu \vec \alpha}).
\end{equation}
This is a $G_{\vec \alpha}$-equivariant stratified submanifold of
$\hat{Z}_{\mu \vec \alpha}$, and hence defines an object of $\Orb \langle
\cP(p,q) \rangle$.  Furthermore, it is straightforward to see that by
construction each component of $\partial^r \vec{Z}_{\vec{f}}$ is the
homeomorphic image of a single component of $\vec{Z}_{\vec{g}_-} \times
\vec{Z}_{\vec{g}_+}$, for uniquely determined elements $\vec{g}_+ \in \vec A(p,r)$
and $\vec{g}_- \in \vec A(r,q)$.

Next, we want to verify that this construction is functorial.
Explicitly, 
given a morphism in $\Tw \vec A(p,q)$ from $\vec{f} \co \vec \alpha \to \vec \beta$ to $\vec{g} \co \vec \alpha' \to \vec \beta'$, we need to construct a map $\vec{Z}_{\vec{f}} \to \vec{Z}_{\vec{g}}$.
The map is induced by the arrow from $\vec \alpha \to \vec
\alpha'$.

The fact that these functors are compatible with composition and
therefore that the above construction assembles into a $2$-functor again follows from the
fact that the twisted arrow category is monoidal and the
multiplicative structure on $\Orb \langle - \rangle$ ultimately comes
from the cartesian product.

Next, we want to map the spectral category of relative cochains
$C^*_{\rel \partial}(B \cZ; \Omega \bbk)$ to $C^*_{\rel   \partial}(B
\vec{\cZ}; \Omega \bbk)$.  Denoting the homotopy colimit of the spaces
$ B \vec{Z}_{\vec{f}}$ by 
\begin{equation} \label{eq:homotopy_colim_over_Tw-vec-A-spaces}
 B \vec{Z}(p,q) \equiv \hocolim_{\vec{f} \in \Tw \vec A(p,q)} B \vec{Z}_{\vec{f}},     
\end{equation}
we have a natural composition
\begin{equation}
 B \vec{Z}(p,q) \to  \hocolim_{\vec{f} \in \Tw \vec A(p,q)} B Z_{\mu \vec
  \alpha} \to \hocolim_{\alpha \in A(p,q)} B Z_{\alpha} = B Z(p,q),
\end{equation}
where the first arrow is induced by composing the inclusion $ 
  \vec{Z}_{\vec{f}} \subset \hat{Z}_{\vec \alpha}$ with the projection
  along the collar directions, and the second by the functor $\Tw
  \vec{A}(p,q) \to A(p,q)$. 

In particular, for each $p$ and $q$, pullback of relative
cochains yields a map of spectra 
\begin{equation}
C^*_{\rel \partial}\big(B Z; \Omega \bbk\big)(p,q) \to C^*_{\rel
  \partial}\big(B \vec{\cZ}; \Omega \bbk\big)(p,q).
\end{equation}
This map is compatible with the composition and clearly
$\Pi$-equivariant, so the pullback induces the functor in Equation
\eqref{eq:compare_cochains_to_collared_cochains}. 
Finally, it is straightforward to check that these morphisms are
pointwise equivalences of spectra.
 \end{proof}
\subsubsection{Compactly supported cochains}
\label{sec:comp-supp-coch-2}
\label{sec:twist-arrow-categ}

To study the functoriality of compactly supported cochains, note that
a map $f \co \alpha \to \beta$ of orbispace charts induces a map
$Z_\alpha \to Z_\beta$ which is the composition of a finite-to-one map
and an open inclusion. Since compactly supported functions are
contravariantly functorial with respect to proper maps, and
covariantly functorial with respect to open inclusions, this suggests
that we factor this map as 
\begin{equation}
  Z_\alpha \to Z_f \to Z_\beta,
\end{equation}
where $Z_f$ the quotient of $Z_\alpha$ by the kernel $G_f^\perp$ of
the map $G_\alpha \to G_\beta$. This quotient admits a residual action
of $G_\alpha / G_f^{\perp} \cong G_\beta$, so that we can form the
classifying space  
\begin{equation}
  BZ_f \equiv Z_f \times_{G_\beta} E G_\beta \cong B(Z_\alpha /
  G_f^{\perp}, G_\beta, \ast).
\end{equation}
The projection $Z_\alpha \to Z_f $ and the surjection $G_\alpha \to G_\beta$
induce a map
\begin{align}
BZ_{\alpha} = B(Z_\alpha, G_\alpha, \ast) \to B(Z_\alpha /
G_f^{\perp}, G_\beta, \ast) = BZ_f
\end{align}
with homotopy fibre $E G_f^\perp$, and which is therefore an equivalence.

We define
\begin{equation}
  C^{*,c}(BZ_f; \Omega \bbk)
\end{equation}
to be the cochains which are compactly supported over $Z_f$ (as in
Definition~\ref{defin:relative-compactly-supported}). The
functoriality of compact support yields a diagram 
\begin{equation}
  C^{*,c}(BZ_\alpha; \Omega \bbk ) \leftarrow  C^{*,c}(BZ_f; \Omega \bbk ) \to C^{*,c}(BZ_{\beta};\Omega \bbk)
\end{equation}
where the first map is pullback along the finite projection $Z_\alpha \to Z_f$, and the second is pushforward along the open inclusion $Z_f \to Z_\beta$. In order to formulate the functoriality of compactly supported
cochains, we are thus again led to pass from the category of orbispace charts to its twisted
arrow category.

For the statement of the next result, recall that $\Orb\langle \cS
\rangle$ is the category with objects orbispace charts stratified by
$\cS$, and morphisms maps of orbispace charts respecting the
stratification. In Equation
\eqref{eq:relative_cochains_stratified_orbichart}, we introduced the
relative cochains of such an orbispace chart.
  \begin{prop}
    The assignment $f \mapsto C^{*,c}_{\rel \partial}(B  Z_f; \Omega \bbk) $ extends to a functor
\begin{equation}
  \Tw   \Orb\langle \cS \rangle \to \Sp.
\end{equation}
  \end{prop}
  \begin{proof}
  Given arrows $f_i \co \alpha_i \to \beta_i$ in $\Orb\langle \cS \rangle$, and a factorisation of $f_0$ as $h \circ f_1 \circ g $, the inclusion $  Z_{f_0} \to Z_{\beta_0}$ factors through $Z_h$. The resulting map $Z_{f_0} \to Z_h$ is an open inclusion, and we obtain a corresponding map
  \begin{equation}
      C^{*,c}_{\rel \partial}(B  Z_{f_0}; \Omega \bbk)  \to C^{*,c}_{\rel \partial}(B  Z_h; \Omega \bbk).  
  \end{equation}
  Composing with the pullback map from $Z_h$ to $Z_{f_1}$, we obtain a map from the compactly supported (relative) cochains of $BZ_{f_0}$ to those of $BZ_{f_1}$:
  \begin{equation}
     C^{*,c}_{\rel \partial}(B Z_{f_0}; \Omega \bbk) \to    C^{*,c}_{\rel \partial}(B Z_{f_1}; \Omega \bbk).      
   \end{equation}
   In order to check the compatibility of this construction with composition, it will be important to have an alternative construction: the projection $Z_{\alpha_0} \to Z_{f_0}$ factors through $Z_g$, so we have a pullback map
   \begin{equation}
     C^{*,c}_{\rel \partial}(B  Z_{f_0}; \Omega \bbk)  \to C^{*,c}_{\rel \partial}(B  Z_h; \Omega \bbk).      
   \end{equation}
   Composing with the pushforward map associated to the open inclusion of $Z_h$ in $Z_{f_1}$, we also obtain a map of compactly supported (relative) cochains. The compatibility between proper pullback and open pushforward maps in compactly supported cochains is then encoded by the commutativity of the diagram
       \begin{equation}
        \begin{tikzcd}
          C^{*,c}_{\rel \partial}(B  Z_{f_0}; \Omega \bbk) \ar[r] \ar[d] & C^{*,c}_{\rel \partial}(B  Z_h; \Omega \bbk) \ar[d] \\
          C^{*,c}_{\rel \partial}(B  Z_g; \Omega \bbk) \ar[r] & C^{*,c}_{\rel \partial}(B  Z_{f_1}; \Omega \bbk).
    \end{tikzcd} 
  \end{equation}
  \end{proof}
  \begin{rem}
In Equation \eqref{eq:relative_cochains_stratified_orbichart}, we
defined the relative cochains by using the collared completion. If the
boundary is already equipped with a collar, or more generally if it is
a closed inclusion, we need not pass to the collared completion
first. This will be used below to avoid an unnecessary double collar. 
\end{rem}

This construction admits functorial boundary maps associated to
collars: given an element $q$ of $ \cS$, there is a natural functor 
\begin{equation}
\Orb \langle \partial^q \cS \rangle \to \Orb \langle \partial \cS
\rangle
\end{equation}
which maps a chart stratified by a subset of $\partial^q \cS$ to the
product with the collar $\kappa^q$ indexed by $q$. Letting $\ell^q$
denote a line associated to $q$, we have a natural map 
\begin{equation} \label{eq:boundary_map_compactly_supported}
   C^{*,c}_{\rel \partial}(B {Z}; \Omega^{\ell + \ell_q}  \bbk) \to    C^{*,c}_{\rel \partial}(B  {Z \times \kappa^q}; \Omega^{\ell}  \bbk),
\end{equation}
which defines a natural transformation between these two functors from
$  \Orb \langle \partial^q \cS \rangle$ to spectra. 

This construction is also functorial for products:  given a pair
$\cS_1$ and $\cS_2$ of partially ordered sets, there is a natural
functor 
\begin{equation}
  \Orb\langle \cS_1 \rangle \times   \Orb\langle \cS_2 \rangle \to \Orb\langle \cS_1 \cup \cS_2 \rangle
\end{equation}
given by taking the product of the underlying spaces and group
actions. This induces a functor on twisted arrow categories.  Given a
pair of real lines $\ell_1$ and $\ell_2$, the homeomorphism $B
Z_{f_1} \times B  Z_{f_2} \cong  BZ_{f_1 \times f_2}$ induces a
natural map 
\begin{equation} \label{eq:product_map_compactly_supported}
  C^{*,c}_{\rel \partial}(B  Z_{f_1}; \Omega^{\ell_1} \bbk ) \wedge
  C^{*,c}_{\rel \partial}(B  Z_{f_2};  \Omega^{\ell_2} \bbk)  \to
  C^{*,c}_{\rel \partial}( BZ_{f_1 \times f_2};  \Omega^{\ell_1 +
    \ell_2} \bbk) 
\end{equation}
which is also functorial in the inputs.

Given an orbispace presentation $\bX \co \hat{A} \to \Orb\langle \cS
\rangle $, and a ring spectrum $\bbk$, we define the compactly
supported cochains of the presentation as the homotopy colimit of the
compactly supported cochains for all arrows.  First, we assume that
the input is collared.

\begin{defin}\label{defin:compactly-supported-relative-cochains}
Given a collared orbispace flow category as in
Definition~\ref{def:collared_orbispace_flow_cat}, we define a
$\Pi$-equivariant spectral category $C^{*,c}_{\rel \partial}(B {\cZ};
\Omega \bbk) $ with objects elements of $\cP$, with morphism spectra
for a pair $(p,q)$ given by the compactly supported relative cochains
\begin{equation} \label{eq:definition_spectral_category_compact_collared}
 C^{*,c}_{\rel \partial}(B  {\cZ}; \Omega  \bbk)(p,q) \equiv    \hocolim_{f \in \Tw \hat{A}(p,q)}  C^{*,c}_{\rel \partial}(B {Z}_f; \Omega^{\ell_q}  \bbk)
\end{equation}
of the corresponding presentation, and composition given by the map
\begin{equation}
  C^{*,c}_{\rel \partial}(B  {\cZ}; \Omega  \bbk)(p,q) \wedge  C^{*,c}_{\rel \partial}(B  {\cZ}; \Omega  \bbk)(q,r)  \to    C^{*,c}_{\rel \partial}(B {\cZ} ; \Omega  \bbk)(p,r)
\end{equation}
induced by the product map from Equation
\eqref{eq:product_map_compactly_supported}, and the boundary map from
Equation \eqref{eq:boundary_map_compactly_supported} associated to the
stratum labelled by $q$. 
\end{defin}

Recall that the attachment of a collar associates to each  Kuranishi
flow category $\bX \co A \to \Kur$ an orbispace flow category with
collars 
\begin{equation}
(\hat{Z}, G) \co A(p,q) \to \Orb \langle \cP(p,q)
\rangle.
\end{equation}
Applying the above construction, we have: 

\begin{defin}
The \emph{category of compactly supported relative cochains}
associated to a Kuranishi flow category $\bX \co A \to \Kur$ is the
$\Pi$-equivariant spectral category  $ C^{*,c}(B {\cZ}; \Omega  \bbk)
$, with objects the elements of $\cP$ and with morphism spectra
\begin{equation}
 C^{*,c}_{\rel \partial}(B  {\cZ}; \Omega  \bbk)(p,q) \equiv
 \hocolim_{f \in \Tw A(p,q)}  C^{*,c}(B \hat{Z}_f, B \partial
 \hat{Z}_f; \Omega^{\ell_q}  \bbk). 
\end{equation}
\end{defin}

Starting with the orbispace flow category with unique factorisation
introduced in Proposition
\ref{prop:collared-presentation-from-usual-presentation}, and applying
Equation \eqref{eq:definition_spectral_category_compact_collared} to
$\hat{A} = \Tw \vec A$,  we can produce another spectrally enriched
category $C^{*,c}(B \vec{\cZ}; \Omega  \bbk) $, with morphism spectra 
\begin{equation}
C^{*,c}_{\rel \partial}(B  \vec{\cZ}; \Omega  \bbk)(p,q) \equiv    \hocolim_{\phi \in \Tw^2 \vec A(p,q)}  C^{*,c}(B \vec{Z}_{\phi}, B \partial \vec{Z}_{\phi}; \Omega^{\ell_q}  \bbk),
\end{equation}
where $\Tw^2 \vec A$ is the twisted arrow category of the category
$\Tw \vec A$ from Section \ref{sec:relat-coch-kuran}, and we write
$\phi$ for a morphism in this category (this consists of a commutative
diagram as in Equation \eqref{eq:morphism_in_twisted_arrow}). In
analogy with Proposition
\ref{prop:collared-presentation-from-usual-presentation}, we have the
following result: 
\begin{lem} \label{lem:equivalence_two_compactly_supported_cochains}
There is a $\Pi$-equivariant spectrally enriched equivalence
\begin{equation}
C^{*,c}_{\rel \partial}(B \vec{\cZ}; \Omega  \bbk) \to C^{*,c}_{\rel \partial}(B {\cZ}; \Omega  \bbk).
\end{equation}
\end{lem}

\begin{proof}[Sketch of proof]
Let
  \begin{equation}
    \mu \co   \Tw^2 \vec A(p,q) \to \Tw A(p,q), 
     \end{equation}
denote the (covariant) functor which is specified on objects by the
assignment to a commutative diagram representing an arrow in $\Tw \vec
A(p,q)$ 
  \begin{equation}
    \begin{tikzcd}
      \vec \alpha_0 \ar[r, "\vec f_0"] \ar[d,swap,"\vec g"] & \vec \beta_0 \\
       \vec \alpha_1 \ar[r,swap,"\vec f_1"] & \vec \beta_1. \ar[u,swap,"\vec h"]
    \end{tikzcd}
  \end{equation}
of the arrow in $A(p,q)$
  \begin{equation}
    \mu \vec g  \co \mu \vec \alpha_0 \to \mu \vec \alpha_1.     
  \end{equation}
It is clear that this assignment is functorial.

The desired comparison functor is then most easily described as the
composition of two maps: the first is the natural map 
  \begin{multline}
    \hocolim_{\phi \in   \Tw^2 \vec A(p,q) }   C^{*,c}(B \hat{Z}_{\mu
      \phi},  B \partial \hat{Z}_{\mu \phi} ; \Omega^{\ell_q}  \bbk)
    \to  \\
 \hocolim_{f \in \Tw A(p,q)}  C^{*,c}(B \hat{Z}_f, B \partial
 \hat{Z}_f; \Omega^{\ell_q}  \bbk).
  \end{multline}
associated to the pullback of the functor $ C^{*,c}(B \hat{Z}_f,
\partial \hat{Z}_f; \Omega^{\ell_q}  \bbk)$ to $ \Tw^2 \vec A(p,q) $
along $\mu$.  The second functor is the natural map
  \begin{multline}
\hocolim_{\phi \in \Tw^2 \vec A(p,q)}   C^{*,c}(B \vec{Z}_{\phi},  B \partial \vec{Z}_{\phi}; \Omega^{\ell_q}  \bbk) \to  \\ \hocolim_{\phi \in \Tw^2 \vec A(p,q)}  C^{*,c}(B \hat{Z}_{\mu \phi},   B \partial \hat{Z}_{\mu \phi}; \Omega^{\ell_q}  \bbk)
 \end{multline}
 induced by the open inclusion
 \begin{equation}
       \vec{Z}_{\phi} \subset \hat{Z}_{\mu \phi}.
     \end{equation}
We can see these functors are equivalences by applying Quillen's
theorem A; the fact that the categories $A_{[u]}$ have contractible
nerve implies that the categorical fibers are contractible.
\end{proof}

\subsubsection{A different model for compactly supported cochains}
\label{sec:diff-model-comp}

In the remainder of Section \ref{sec:comp-comp-supp}, we shall work with an orbispace flow category with unique factorisation as in Definition \ref{def:collared_orbispace_flow_cat}, with the goal of comparing its compactly supported and ordinary relative cochains.

\begin{rem}
  The essential difficulty in implementing the desired comparison is that the map
  \begin{equation}
    B Z_\alpha \to BZ \equiv \hocolim_{\alpha} BZ_\alpha    
  \end{equation}
  is not the inclusion of the inverse image of an open subset of $\colim Z_\alpha/G_\alpha$, hence does not induce a map of compactly supported cochains over this space.
\end{rem}

Let $\bX \co A \to \Orb \langle \cS \rangle $ be an orbispace presentation of $\cM$. For each object $\alpha \in A$, define 
\begin{equation}
 Z^{\alpha} \equiv \hocolim_{\beta \in A} Z_\beta \times_{\cM} (Z_{\alpha}/G_{\alpha}).
\end{equation}
In other words, for each chart $\beta$, we take the open subset of $Z_\beta$ consisting of points whose projection to the space $\cM$ lies in the footprint of the chart $\alpha$, and (homotopy) glue these spaces together over all objects $\beta \in A$. This construction is covariantly functorial in $\alpha$.
\begin{lem}
The map   $Z^{\alpha} \to Z^{\beta}$ induced by a morphism $f \co \alpha \to \beta$  is an open embedding. \qed
\end{lem}

Passing to Borel constructions, we define  
\begin{equation}
 BZ^{\alpha} \equiv \hocolim_{\beta \in A} BZ_\beta \times_{Z} Z_{\alpha}/G_{\alpha}.
\end{equation}
The corresponding functor $A \to \Top$ maps each arrow to an embedding lying over the open embedding of footprints, so we obtain an induced covariant functor
\begin{equation}
    \alpha \mapsto  C^{*,c}_{\rel \partial}(BZ^\alpha; \Omega \bbk)
  \end{equation}
  where the support condition is again that cochains are required to
  vanish away from the inverse image of a compact subset of
  $Z^\alpha$. 
  
Passing to the homotopy colimit over the indexing category, we have
the following definition.

\begin{defin}
Let $\bX \co A \to \Orb \langle \cS \rangle $ be an orbispace
presentation of $\cM$.  We define the compactly supported cochains as
 \begin{equation} \label{eq:compactly_supported_intermediate}
C^{*,c}_{\rel \partial}( BZ^{\bullet}; \Omega \bbk) \equiv   \hocolim_{\alpha \in A}    C^{*,c}_{\rel \partial}(BZ^\alpha; \Omega \bbk).
 \end{equation} 
\end{defin}

This construction is functorial with respect to collars, as we have a natural map
\begin{equation}
    C^{*,c}( B Z^{\bullet}; \Omega^{\ell +\ell_q} \bbk)  \to  C^{*,c}_{\rel \partial}( B(Z\times \kappa^q) ^{\bullet}; \Omega^\ell \bbk)
\end{equation}
which is analogous to Equation \eqref{eq:boundary_map_compactly_supported}. 

As to multiplicativity, we note that for a product presentation $A(1)
\times A(2) \to \Orb$, there is a natural homeomorphism 
\begin{equation}
 BZ^{\alpha_1} \times BZ^{\alpha_2} \to BZ^{\alpha_1 \times \alpha_2} 
\end{equation}
for each chart, which induces a map
\begin{equation}
 C^{*,c}_{\rel \partial}( BZ^{\bullet}(1); \Omega^{\ell_1} \bbk) \wedge   C^{*,c}_{\rel \partial}( BZ^{\bullet}(2); \Omega^{\ell_1} \bbk)  \to  C^{*,c}_{\rel \partial}( BZ^{\bullet}(12); \Omega^{\ell_1+\ell_2} \bbk) 
\end{equation}

\begin{rem} \label{rem:functoriality_modified_compactly_supported}
The functoriality of
Equation~\eqref{eq:compactly_supported_intermediate} with respect to
the indexing category is more delicate: given a functor $F \co A(1)
\to A(2)$, we have a natural map 
\begin{equation} \label{eq:inclusion_enlarged_chart_functor}
    BZ^{\alpha} \to BZ^{F(\alpha)} 
\end{equation}
which is compatible with the projection to $Z_\alpha$. We thus obtain a diagram
\begin{multline}
   C^{*,c}_{\rel \partial}( BZ^{\bullet}(1); \Omega \bbk)   \leftarrow
   \hocolim_{(f : \alpha \to \beta)\in A(1)}    C^{*,c}_{\rel \partial}(BZ^{F(\alpha)}; \Omega \bbk) \\ \to  C^{*,c}_{\rel \partial}( BZ^{\bullet}(2); \Omega \bbk).  
 \end{multline}
 We shall essentially avoid appealing to functoriality with respect to changing the index category by arranging for Equation \eqref{eq:inclusion_enlarged_chart_functor} to be an isomorphism.
\end{rem}

We now put together the above ingredients to construct a spectral
category $ C^{*,c}_{\rel \partial}( B{\cZ}^{\bullet}; \Omega \bbk) $
associated to each orbispace flow category with unique factorization. 
The factorization condition implies that the map
\begin{equation} \label{eq:isomorphism_homotopy_colim_zero_loci}
\hocolim_{\alpha_1 \times \alpha_2 \in \hat{A}(p,q) \times
\hat{A}(q,r)} Z_{\alpha_1} \times Z_{\alpha_2} \to  \hocolim_{\alpha
\in \hat{A}(p,r)} \partial^q Z_{\alpha}
\end{equation}
is a homeomorphism. We have a corresponding homeomorphism after
passing to Borel constructions, and taking fibre products over $Z(p,q)
\times Z(q,r) \cong \partial^q Z(p,r)$ yields a homeomorphism 
\begin{equation}\label{eq:Z-homeo}
Z^{\alpha_1} \times Z^{\alpha_2} \cong Z^{\alpha_1 \times \alpha_2}.
\end{equation}
This allows us to bypass the zig-zag discussed in Remark
\ref{rem:functoriality_modified_compactly_supported}.

\begin{defin}
Given an orbispace flow category with unique factorization, we define
\begin{equation}
 C^{*,c}_{\rel \partial}( B{\cZ}^{\bullet}; \Omega \bbk)(p,q) \equiv
 \hocolim_{(f \colon \alpha \to \beta) \in \hat{A}(p,q)}    C^{*,c}_{\rel \partial}(B{Z}^\alpha; \Omega^{\ell_q} \bbk). 
\end{equation}
The composition is induced by Equation~\eqref{eq:Z-homeo} and yields a
composition map 
\begin{equation}
  C^{*,c}_{\rel \partial}( B{\cZ}^{\bullet}; \Omega \bbk) (p,q) \wedge C^{*,c}_{\rel \partial}( B{\cZ}^{\bullet}; \Omega \bbk) (q,r) \to C^{*,c}_{\rel \partial}( B{\cZ}^{\bullet}; \Omega \bbk) (p,r)
\end{equation}
for triples in $p,q,r \in \cP$, which is associative and unital.
\end{defin}

\subsubsection{Comparison with ordinary cochains}
\label{sec:comp-with-ordin}

The starting point of the comparison with ordinary cochains is to
consider an orbispace presentation $A \to \Orb \langle \cS \rangle$:
there is a natural map 
\begin{equation}\label{eq:compact-comparison-sections}
C^{*,c}_{\rel \partial}(BZ^\alpha; \Omega \bbk) \to C^*_{\rel
  \partial}(BZ; \Omega \bbk) 
\end{equation}
induced by the inclusion $BZ^\alpha \to BZ$ and the inclusion of
the space of compactly-supported sections in the space of all
sections.  This map incudes a weak equivalence of spectral categories.

\begin{lem}
For each orbispace presentation, the map
 \begin{equation} \label{eq:compactly_supported_intermediate_to_cochains}
C^{*,c}_{\rel \partial}(BZ^\bullet; \Omega \bbk) \to C^*_{\rel
  \partial}(BZ; \Omega \bbk) 
\end{equation}
is an equivalence.
\end{lem}

\begin{proof}
Expanding, we need to show that the map
\begin{equation}
\hocolim_{\alpha \in A} C^{*,c}_{\rel \partial}(BZ^\alpha; \Omega \bbk)
 \to C^*_{\rel \partial}(\hocolim_{\alpha \in A} BZ_\alpha; \Omega \bbk)
\end{equation}
is an equivalence.  Since $\cM$ is compact, we can choose a finite
subcover of $\{Z^\alpha\}$ with indexing set $J$ and a partition of
unity subordinated to this subcover such that the each function $f_i$
has compact support.  We can now conclude that the evident restriction
map
\begin{equation}
C^{*,c}_{\rel \partial}(\hocolim_{\alpha \in A} BZ^\alpha; \Omega
\bbk) \to \hocolim_{\alpha \in A} C^{*,c}_{\rel \partial}(BZ^\alpha; \Omega
\bbk) 
\end{equation}
induced by the partition of unity is a weak equivalence.  Since the
comparison between globally supported compactly-supported sections and
all sections is a weak equivalence, the result follows.
\end{proof}

The comparison map from Equation~\eqref{eq:compact-comparison-sections} is
compatible both with collars and with the product of presentations. 
Therefore, applying it to the collared categories from Section
\ref{sec:relat-coch-kuran}, we conclude:
\begin{lem} \label{lem:compacre_strange_model_ordinary}
Given an orbispace flow category with unique factorization, the
comparison map induces a $\Pi$-equivariant equivalence 
\begin{equation}
  C^{*,c}_{\rel \partial}( B{\cZ}^{\bullet}; \Omega \bbk)    \to  C^*_{\rel \partial}( B{\cZ}; \Omega \bbk)
\end{equation}
of spectral categories.
\qed
\end{lem}

\subsubsection{Comparison with compactly supported cochains}
\label{sec:comp-with-comp}
Observe that, for each arrow $f \co \alpha \to \beta $, we have a homeomorphism
\begin{equation}
Z_f \equiv     Z_\beta \times_{X} (Z_{\alpha}/G_{\alpha}),
\end{equation}
which induces an inclusion
\begin{equation}
 Z_f \to Z^\alpha
\end{equation}
lying over the natural homeomorphism $Z_f/G_\beta \to Z_\alpha/G_\alpha$.

Passing to the compactly supported relative cochains of the Borel
constructions, we obtain a map 
\begin{equation} \label{eq:map_compactly_supported_two_models}
C^{*,c}_{\rel \partial}( BZ^\alpha, \Omega \bbk) \to   C^{*,c}_{\rel \partial}( BZ_f , \Omega \bbk).
\end{equation}

\begin{lem}
The map of Equation~\eqref{eq:map_compactly_supported_two_models}
defines a natural transformation between these functors from $\Tw
\hat{A}$ to $\Sp$ and hence on passage to homotopy colimits maps
\begin{equation} \label{eq:comparison_models_compactly_supported_cochains}
C^{*,c}_{\rel \partial}( BZ^{\bullet}; \Omega \bbk) \to C^{*,c}_{\rel
  \partial}( BZ; \Omega \bbk) 
\end{equation}
for each orbispace presentation. \qed
\end{lem}

These maps are again compatible with boundaries and products, and so
we conclude: 

\begin{lem} \label{lem:equivalence_strange_model_compact_cochains}
The map in Equation
\eqref{eq:comparison_models_compactly_supported_cochains} induces a
$\Pi$-equivariant equivalence of spectral categories
\begin{equation}
C^{*,c}_{\rel \partial}( B{\cZ}^{\bullet}; \Omega \bbk) \to C^{*,c}_{\rel \partial}( B{\cZ}; \Omega \bbk) .
\end{equation} \qed
\end{lem}

\subsubsection{Signpost: Compactly supported cochains}
\label{sec:signp-comp-supp}

We now continue the discussion of Section~\ref{sec:signpost}: given a
Kuranishi flow category, we have extended the maps from
Equation~\eqref{eq:sign_post_Kuranishi-spaces} in one direction, to
obtain a zig-zag of equivalences of $\Pi$-equivariant spectral
categories 
\begin{equation}
C^*(\cP ; \Omega \bbk) \to C^*_{\rel \partial}(B {\cZ} ; \Omega \bbk)
\leftarrow  \cdots  \to C^{*,c}_{\rel \partial}( B {\cZ}; \Omega
\bbk), 
\end{equation}
where the omitted arrows arise by applying
Proposition~\ref{prop:collared-presentation-from-usual-presentation},
and
Lematta~\ref{lem:equivalence_two_compactly_supported_cochains},~\ref{lem:compacre_strange_model_ordinary},
and~\ref{lem:equivalence_strange_model_compact_cochains}.

What remains to be done to construct a homotopy type is to produce a
further zig-zag  
\begin{equation}
C^{*,c}_{\rel \partial}( B {\cZ}; \Omega \bbk) \dashleftarrow \cdots \dashrightarrow  B \cX|\cZ^{-V-d} \wedge \bbk.
\end{equation}
In the absence of group actions, this is essentially asserting the
existence of a coherent (untwisted) Poincar\'e duality equivalence.

\subsection{The Milnor fibrations for flow categories}
\label{sec:span-whit-dual-3}

Our goal in this section is to prepare the ground to lift the
equivalence from Section~\ref{sec:span-whit-dual-5} to Kuranishi flow
categories, by constructing an associated spectral category built from
the spectra of sections of Milnor fibrations. We begin by addressing
the technical problem alluded to in Remark
\ref{rem:lack-of-functoriality-Milnor}.

\subsubsection{Inessential charts}
\label{sec:inessential-charts}

\begin{defin}
A map $f \co \bX \to \bX'$ of Kuranishi charts is \emph{inessential}
if the covering map $X \to X/G_f$ is trivial.  
\end{defin}
We extend the terminology to presentations and flow categories by
saying that they are inessential if all morphisms satisfy this
property.  

\begin{prop}
Let $\bX \co A \to \Kur$ be a $\Pi$-equivariant Kuranishi flow
category lifting a topological flow category $\cM$.  There exists an
inessential $\Pi$-equivariant flow category $\underline{\bX} \co
\underline{A} \to \Kur$ lifting the collared completion $\hat{\cM}$.
\end{prop}

\begin{proof}
For each pair $(p,q)$ of objects of $\cP$, we define a category
$\underline{A}(p,q)$ with 
\begin{enumerate}
\item objects the pairs $ \underline{\alpha} =
(\alpha, X_{\underline{\alpha}}) $, with  $X_{\underline{\alpha}}$ an
open subset of $\hat{X}_\alpha$ such that the inclusion map $
X_{\underline{\alpha}} \to \hat{X}_\alpha$  induces the trivial map on
fundamental groups, and
\item morphisms from $\underline{\alpha}$  to $\underline{\beta}$ given by
  a map $\alpha \to \beta$ such that the image of
  $X_{\underline{\alpha}}$ in $X_\beta$ is contained in
  $X_{\underline{\beta}}$. 
\end{enumerate}
The assignment 
\begin{equation}
\underline{\alpha} \mapsto  \underline{\bX}_{\underline{\alpha}}
\equiv  (G_\alpha, X_{\underline{\alpha}}, V_\alpha,  s_\alpha) 
\end{equation}
thus specifies a functor $\underline{\bX}$ to $\Kur$.

By construction, $\underline{\bX}$ is a functor over
$\hat{\cM}(p,q)$. The condition on fundamental groups implies that all
morphisms in the image of this functor are inessential.  To see that
$\underline{\bX}$ is a Kuranishi presentation, we consider 
the evident functor $\underline{A}(p,q) \to A(p,q)$ which assigns
$\alpha$ to $\underline{\alpha}$.  Since any sufficiently small
neighbourhood of a point in a manifold is inessential, Quillen's
Theorem A implies that the induced map on nerves
\begin{equation}
N_\bullet \underline{A}(p,q)[z] \to N_\bullet A(p,q)[z] 
\end{equation}
is a weak equivalence and hence $N_\bullet \underline{A}(p,q)[z]$ is
contractible for each $[z] \in \hat{\cM}(p,q)$.  That is,
$\underline{\bX}$ is a Kuranishi presentation. 

We define the product by assigning to a pair $\underline{\alpha}_1 \in  \underline{A}(p,q)$ and $\underline{\alpha}_2 \in  \underline{A}(q,r)$ the object
\begin{equation}
   \underline{\alpha}_1 \times \underline{\alpha}_2 \equiv (\alpha_1 \times \alpha_2 ,  X_{\underline{\alpha}_1} \times X_{\underline{\alpha}_2} \times \kappa^{q})
\end{equation}
of $A(q,r)$ consisting of the collar on the inclusion of $X_{\underline{\alpha}_1} \times X_{\underline{\alpha}_2}$ in the boundary stratum associated to $q$. It is straightforward to see that this yields a Kuranishi flow category $ \underline{\bX} \co \underline{A} \to \Kur$.
\end{proof}
\begin{rem} The only reason to pass to the collared category is that it provides us with a natural formula for the product of inessential charts. It is plausible that one can use an inductive scheme to construct a flow category lifting $\cM$ itself. One reason that it may not be worth it to pursue such a result is that it is not difficult to strengthen the above result by proving that all the cochains models for $\underline{\bX}$ are equivalent to the corresponding ones for $\bX$, thus allowing us to produce from the inessential presentation for $\hat{\cM}$, the desired data for $\cM$ itself. 
\end{rem}

We shall henceforth assume that every Kuranishi flow category we consider is inessential.

\subsubsection{One point compactifications and virtual cochains}
\label{sec:one-point-comp}

Recall that we introduced a (partial) compactification $X_\sigma
\subset \bar{X}_\sigma$ by adding a point at infinity for the vector
space direction of all facets of the form $X_\tau \times V_f^\perp$
(c.f. Lemma~\ref{lem:pushout_description_compactified_degeneration}). Since
our 
discussion of Spanier-Whitehead duality for manifolds with boundary
required passing to collared completions (see Section
\ref{sec:span-whit-dual-2}), we revise the construction of virtual
cochains to incorporate collars, and thus define
$\hat{\bar{X}}_{{\sigma}}$ to be the pushout 
  \begin{equation}
    \begin{tikzcd}
      \displaystyle{\coprod_{\rho \overset{g}{\to} \tau \overset{f}{\to} \sigma}  (\hat{X}_{{f}\circ {g}}  \times V_g^\perp)_+ \sma S^{V_f^\perp} }  \ar[r] \ar[d] &      \displaystyle{\coprod_{\rho \overset{f \circ g}{\to} \sigma}    \hat{X}_{{f}\circ {g}} \times S^{V_{f \circ g}^\perp} }\ar[d] \\
  \displaystyle{\coprod_{ \tau \overset{f}{\to} \sigma}  (\hat{X}_{{f}})_+ \sma S^{V_f^\perp} }\ar[r] & \hat{\bar{X}}_{{\sigma}}
    \end{tikzcd}
  \end{equation}
in the category of based spaces, where the coproducts are taken over
boundary strata (i.e., injective maps).  By construction, the
embedding $\hat{Z}_{\sigma} \subset \hat{\bar{X}}_{{\sigma}}$ avoids
the basepoint. 

This construction is compatible with the maps to $V_\tau$ and
$V_\sigma$, and thus yields a functor
\begin{align}
\square {\Kur}\langle \cS \rangle & \to \Sp \\
\nonumber {\sigma} & \mapsto \hat{\bar{X}}_{{\sigma}}| \hat{Z}_{\sigma}^{-V_\sigma}.
\end{align}
Returning to the construction from Section \ref{sec:virt-coch-homot},
we introduce a new model 
\begin{equation}
B\hat{\bar{X}}|\hat{Z}^{-V}({A})  \equiv   \hocolim_{ {\sigma} \in \square {A}}   B \hat{\bar{X}}_{{\sigma}}|\hat{Z}_{\sigma}^{-V_\sigma}   
\end{equation}
for the virtual cochains of a Kuranishi presentation.  This model of
the virtual cochains is multiplicative, in that there is an
associative product map 
\begin{multline}
B\hat{\bar{X}}|\hat{Z}^{-V}({A(1)}) \sma
B\hat{\bar{X}}|\hat{Z}^{-V}({A(2)}) \to \\ B\bar{X}|Z^{-V}(A(1))
\sma B\bar{X}|Z^{-V}(A(2)).
\end{multline}

For each $\sigma \in \square A$,  the map $ \hat{\bar{X}}_\sigma \to
\bar{X}_\sigma$ collapsing the collars induces an equivalence 
\begin{equation}\label{eq:collar-collapse}
\hat{\bar{X}}_{{\sigma}}|\hat{Z}_{\sigma} \to \bar{X}_\sigma|Z_\sigma,
\end{equation}
so we conclude:
\begin{lem}
For a Kuranishi presentation $A$, the collapse map induces an
equivalence of spectra
  \begin{equation} \label{eq:map_inessential_charts_virtual_cochains}
   B\hat{\bar{X}}|\hat{Z}^{-V}({A}) \to B\bar{X}|Z^{-V}(A).
  \end{equation} \qed
\end{lem}

The above equivalence is multiplicative.

\begin{lem}
The following diagram commutes
\begin{equation}
\begin{tikzcd}
B\hat{\bar{X}}|\hat{Z}^{-V}({A(1)}) \sma
B\hat{\bar{X}}|\hat{Z}^{-V}({A(2)}) \ar[r] \ar[d] & \ar[d] B\hat{\bar{X}}|\hat{Z}^{-V}({A(1) \times A(2)}) \\
B\bar{X}|Z^{-V}(A(1)) \sma B\bar{X}|Z^{-V}(A(2))  \ar[r] & B\bar{X}|Z^{-V}(A(1) \times A(2)).\\
\end{tikzcd}
\end{equation} \qed
\end{lem}

Equation~\eqref{eq:collar-collapse} is also compatible with the
inclusion of boundary strata, which justifies the following definition, where we use brackets to denote shifts using the multiplicative
spheres of Appendix~\ref{sec:action}. 

\begin{defin}
Given a Kuranishi flow category $\bX$ with objects $\cP$, we define
the $\Pi$-equivariant spectral category
$B\hat{\bar{\cX}}|\hat{Z}^{-V-d}$ to have morphism spectra
\begin{equation}
B\hat{\bar{\cX}}|\hat{Z}^{-V-d}(p,q)  \equiv
B\hat{\bar{X}}|\hat{Z}^{- V + V_{q} - V_{p} }(p,q)[ \deg p  - \deg q
], 
\end{equation}
as in Equation~\eqref{eq:morphisms_category_virtual_cochains}.   The composition is defined via
the evident extension of
Lemma~\ref{lem:category_virtual_cochains_composition} and the unit is
defined using the unit map $\bS \to \bS[0]$.
\end{defin}

We now have the following comparison.

\begin{lem}
The map in Equation~\eqref{eq:map_inessential_charts_virtual_cochains}
induces a $\Pi$-equivariant equivalence of spectral categories 
\begin{equation}
 B \hat{\bar{\cX}}|\hat{Z}^{-V-d} \to B\bar{\cX}|\cZ^{-V-d}.
\end{equation} \qed
\end{lem}

\subsubsection{Milnor fibration and cubical diagrams}
\label{sec:miln-fibr-cubic}

As in Section~\ref{sec:one-point-comp}, consider a cube ${\sigma}$ in
$\square {\Kur}\langle \cS \rangle$, with domain $\bbone^n$.
Let $ \widetilde{X}_{{\sigma}}$ denote the completion of
$\hat{X}_{{\sigma}}$ to a manifold without boundary, obtained by
attaching infinite ends, and let $\widetilde{\bar{X}}_{{\sigma}}$ be
the corresponding completion of $\hat{\bar{X}}_{{\sigma}}$.  Given a
map ${f} \co {\sigma} \to  {\tau}$ of cubes in $\square {\Kur}\langle \cS
\rangle$, define   $ \hat{\bar{X}}_{{f}}$ and
$\widetilde{\bar{X}}_{{f}}$ to be the quotients of $
\hat{\bar{X}}_{{\tau}}$ and  $\widetilde{\bar{X}}_{{\tau}}$ by
$G_f^\perp$. 

\begin{defin}\label{defn:milnor-spherical-fibration}
The \emph{Milnor spherical fibration} $\Milnor \hat{\bar{X}}_{{f}} \to
\hat{X}_{{f}}$ is the fibrewise cone of the complement of the
diagonal section of the projection $\widetilde{\bar{X}}_{{f}} \times
\hat{X}_{{f}} \to \hat{X}_{{f}}$. 
\end{defin}

We consider $\Milnor \hat{\bar{X}}_{{f}} $ as a fibration over a
zero-locus $\hat{s}^{-1}_f(0) \equiv \hat{Z}_{f}$. The group $G_f
\equiv G_\sigma$ naturally acts on $ \hat{\bar{X}}_{{f}}$ and
$\widetilde{\bar{X}}_{{f}}$, and the maps to $V_f \equiv V_\tau$ are
equivariant, so the pair $(\hat{Z}_f,G_f)$ is an orbispace chart. Note
that $\hat{Z}_{f}$ is the image of $\hat{Z}_\tau$ under the projection
map, which is obtained by taking the free quotient by $G_f^\perp$.  

Thus, we can associate to a map $f$ three models for the spectrum of
compactly supported relative sections of the Milnor fibration:
\begin{multline} \label{eq:map_models_compact_support_Milnor}
C^{*,c_0}_{\rel \partial}(B \hat{Z}_{{f}};\Milnor
\hat{\bar{X}}_{{f}}^{-V } \wedge \bbk) \leftarrow C^{*,c'_0}_{\rel
  \partial}(B \hat{Z}_{{f}};\Milnor \hat{\bar{X}}_{{f}}^{-V } \wedge
\bbk) \\ \to  C^{*,c}_{\rel \partial}(B \hat{Z}_{{f}};\Milnor
\hat{\bar{X}}_{{f}}^{-V } \wedge \bbk). 
\end{multline}
As in Section \ref{sec:span-whit-dual-4}, $c_0$ refers to sections
with value in the conical part outside a compact set, and which are
inward pointing at the boundary; $c$ refers to sections which strictly
vanish outside a compact set and on the boundary. The support
condition $c'_0$ is a subset of the support condition $c$, and agrees
with the support condition $c_0$ upon restricting to the union of the
interior with the union of the products of the boundary strata with
half-collars $[0,1/2]$. In order to minimise needless repetitions, we
write 
\begin{equation}
C^{*,\tilde{c}}_{\rel \partial}(B \hat{Z}_{{f}};\Milnor
\hat{\bar{X}}_{{f}}^{-V } \wedge \bbk)
\end{equation}
for any of these models.

The proof of the following result is elementary, but it is the key
reason for introducing inessential charts: 
\begin{lem} \label{lem:map_pullback_Milnor}
Given maps $f \colon \sigma \to \tau$ and $g \colon \rho \to \sigma$,
there is a natural equivalence of spherical fibrations over
$\hat{Z}_{{g}}$  
\begin{equation}
\Milnor \hat{\bar{X}}_{{f}\circ {g}}   \to \Milnor \hat{\bar{X}}_{{g}},
\end{equation}
where the first fibration is obtained by pullback under the natural
map $\hat{Z}_{{g}} \to \hat{Z}_{{f}\circ {g}}$. 
\end{lem}
\begin{proof}
By assumption, the covering space $X_{{\rho}} \to X_{{f}\circ {g} } $
is the product of $X_{{f}\circ {g} } $ with $G_{f \circ g}$, hence the
covering space $X_{{g}} \to X_{{f}\circ {g} }  $ is the product of
$X_{{f}\circ {g} }$ with $G_f^\perp$. The map of fibrations is then
obtained on each point $z$ in $\hat{Z}_{ {g}} $ by mapping
$\hat{X}_{{f}\circ {g} }  $ to the covering sheet of $\hat{X}_{ {g}} $
containing the image of $z$. 
\end{proof}

Thus, for each composition $f \circ g$, pullback defines maps of
twisted cochains 
\begin{equation}
\begin{tikzcd}
C^{*,\tilde{c}}_{\rel \partial}(B \hat{Z}_{{f}\circ {g}};\Milnor
\hat{\bar{X}}_{{f}\circ {g}}^{-V } \wedge \bbk)   \ar[r] &
C^{*,\tilde{c}}_{\rel \partial}(B \hat{Z}_{{g}};\Milnor
\hat{\bar{X}}_{{g}}^{-V } \wedge \bbk)  
\end{tikzcd}
\end{equation}
that fit into the commutative diagram
\begin{equation}\label{eq:pullback_compatibility}
\begin{tikzcd}
C^{*,c_0}_{\rel \partial}(B \hat{Z}_{{f}\circ {g}};\Milnor \hat{\bar{X}}_{{f}\circ {g}}^{-V } \wedge \bbk)   \ar[r] &    C^{*,c_0}_{\rel \partial}(B \hat{Z}_{{g}};\Milnor \hat{\bar{X}}_{{g}}^{-V } \wedge \bbk) \\
C^{*,c'_0}_{\rel \partial}(B \hat{Z}_{{f}\circ {g}};\Milnor \hat{\bar{X}}_{{f}\circ {g}}^{-V } \wedge \bbk)   \ar[r] \ar[d] \ar[u] &    C^{*,c'_0}_{\rel \partial}(B \hat{Z}_{{g}};\Milnor \hat{\bar{X}}_{{g}}^{-V } \wedge \bbk) \ar[d] \ar[u]  \\
C^{*,c}_{\rel \partial}(B \hat{Z}_{{f}\circ {g}};\Milnor \hat{\bar{X}}_{{f}\circ {g}}^{-V } \wedge \bbk)   \ar[r]  &    C^{*,c}_{\rel \partial}(B \hat{Z}_{{g}};\Milnor \hat{\bar{X}}_{{g}}^{-V } \wedge \bbk). 
  \end{tikzcd}
\end{equation}

On the other hand, the inclusion of the collar together with the
section of the fibration in the collar direction (see Equation
\eqref{eq:section_Milnor_collar_vanish_boundary}) with appropriate
values at the ends, induce pushforward maps of twisted cochains 
\begin{equation}
C^{*,\tilde{c}}_{\rel \partial}(B \hat{Z}_{{f}\circ {g}};\Milnor
\hat{\bar{X}}_{{f}\circ {g}}^{-V } \wedge \bbk)   \to
C^{*,\tilde{c}}_{\rel \partial}(B \hat{Z}_{{f}};\Milnor
\hat{\bar{X}}_{{f}}^{-V } \wedge \bbk),
\end{equation}
which fit into an analogous compatibility diagram to
Equation~\eqref{eq:pullback_compatibility}.

To state the compatibility between these two constructions, consider a
factorization ${f} = {{g}} \circ {f}' \circ {h}$, which induces a
commutative diagram of pairs  
\begin{equation} \label{eq:arrow_category_cubes_diagram_pairs}
\begin{tikzcd}
(\hat{X}_{{f}' \circ {h}} \times V^\perp_h \times V^\perp_f \times
  V_{g}^\perp, \hat{Z}_{{f}' \circ {h}}) \ar[twoheadrightarrow,r]
  \ar[hookrightarrow,d] &  (\hat{X}_{{f}} \times V^\perp_h \times
  V^\perp_{g \circ f'}, \hat{Z}_{{f}}) \ar[hookrightarrow,d] \\ 
(\hat{X}_{{f}'} \times V_f^\perp \times V_{g}^\perp, \hat{Z}_{{f}'})
  \ar[twoheadrightarrow,r]  & (\hat{X}_{{g} \circ {f}'} \times
  V^\perp_{g \circ f'}, \hat{Z}_{{g} \circ {f}'}) 
\end{tikzcd}
\end{equation}
where the horizontal maps are quotient maps by the action of a finite
group, and the vertical maps are inclusions of boundary strata. 

The outcome of the previous discussion is that we have a commutative diagram
\begin{equation} \label{eq:arrow_category_cubes_diagram_Milnor_twisted_cochains}
\begin{tikzcd}
    C^{*,\tilde{c}}(B \hat{Z}_{{f}' \circ {h}}; \Milnor \hat{\bar{X}}_{{f}' \circ {h}}^{-V } \wedge \bbk) \ar[d]  & \ar[l] C^{*,\tilde{c}}(B \hat{Z}_{{f}}; \Milnor \hat{\bar{X}}_{{f}}^{-V } \wedge \bbk) \ar[d]  \\
    C^{*,\tilde{c}}(B \hat{Z}_{{f}'}; \Milnor \hat{\bar{X}}_{{f}'}^{-V } \wedge \bbk) & \ar[l]  C^{*,\tilde{c}}(B \hat{Z}_{{g} \circ {f}'}; \Milnor \hat{\bar{X}}_{{g} \circ {f}'}^{-V } \wedge \bbk).
  \end{tikzcd}
\end{equation}
We can thus assign to each arrow from ${f}$ to ${f}'$ in $\Tw \square
{\Kur} \langle \cS \rangle$ the map 
\begin{equation}
C^{*,\tilde{c}}(B \hat{Z}_{{f}}; \Milnor \hat{\bar{X}}_{{f}}^{-V }
\wedge \bbk) \to C^{*,\tilde{c}}(B \hat{Z}_{{f}'}; \Milnor
\hat{\bar{X}}_{{f}'}^{-V } \wedge \bbk) 
\end{equation}
obtained by composition around this diagram.  To see that composition
is associative, we observe that for a composable pair of
factorizations ${f} = {{g}} \circ {f}' \circ {h}$ and $f' = g' \circ
f'' \circ h'$ we have the following commutative diagram arising from
an elaboration of
Equation~\eqref{eq:arrow_category_cubes_diagram_pairs}:
\begin{equation}\label{eq:composition_cube_diagram}
\begin{tikzpicture}
\node[scale=.75]{
\begin{tikzcd}
C^{*,\tilde{c}}(B \hat{Z}_{f'' \circ h' \circ {h}}; \Milnor
\hat{\bar{X}}_{{f}'' \circ h' \circ {h}}^{-V } \wedge \bbk) \ar[d]
& C^{*,\tilde{c}}(B \hat{Z}_{{f}' \circ {h}}; \Milnor
\hat{\bar{X}}_{{f}' \circ {h}}^{-V } \wedge \bbk) \ar[d]  \ar[l] & \ar[l]
C^{*,\tilde{c}}(B \hat{Z}_{{f}}; \Milnor \hat{\bar{X}}_{{f}}^{-V }
\wedge \bbk) \ar[d]  \\ 
C^{*,\tilde{c}}(B \hat{Z}_{{f}'' \circ {h'}}; \Milnor
\hat{\bar{X}}_{{f}'' \circ {h'}}^{-V } \wedge \bbk) \ar[d] 
& \ar[l] \ar[d] C^{*,\tilde{c}}(B \hat{Z}_{{f}'}; \Milnor \hat{\bar{X}}_{{f}'}^{-V }
\wedge \bbk) & \ar[l]  C^{*,\tilde{c}}(B \hat{Z}_{{g} \circ {f}'};
\Milnor \hat{\bar{X}}_{{g} \circ {f}'}^{-V } \wedge \bbk) \ar[d] 
\\
C^{*,\tilde{c}}(B \hat{Z}_{{f}''}; \Milnor \hat{\bar{X}}_{{f}''}^{-V }
\wedge \bbk) & \ar[l]  C^{*,\tilde{c}}(B \hat{Z}_{{g'} \circ {f}''};
\Milnor \hat{\bar{X}}_{{g'} \circ {f}''}^{-V } \wedge \bbk)
& \ar[l]
C^{*,\tilde{c}}(B \hat{Z}_{g \circ {g'} \circ {f}''};
\Milnor \hat{\bar{X}}_{g \circ {g'} \circ {f}''}^{-V } \wedge \bbk)
. 
\end{tikzcd}
};
\end{tikzpicture}
\end{equation}
We then have the following lemma, which expresses the functoriality of
pullback and pushforward in this context.
\begin{lemma}
Given composable maps $f_2 \circ f_1 \circ f_0$, the composite of
pushforward maps 
\begin{equation}
\begin{tikzpicture}
\node[scale=.8]{
\begin{tikzcd}
C^{*,\tilde{c}}(B \hat{Z}_{f_2}; \Milnor
\hat{\bar{X}}_{f_2}^{-V } \wedge \bbk) 
& C^{*,\tilde{c}}(B \hat{Z}_{f_2 \circ f_1}; \Milnor
\hat{\bar{X}}_{f_2 \circ f_1}^{-V } \wedge \bbk) \ar[l] & \ar[l]
C^{*,\tilde{c}}(B \hat{Z}_{f_2 \circ f_1 \circ f_0}; \Milnor
\hat{\bar{X}}_{{f_2 \circ f_1 \circ f_0}}^{-V } 
\wedge \bbk)
\end{tikzcd}
};
\end{tikzpicture}
\end{equation}
coincides with the map induced by pushforward of the composite
\begin{equation}
\begin{tikzcd}
C^{*,\tilde{c}}(B \hat{Z}_{f_2}; \Milnor
\hat{\bar{X}}_{f_2}^{-V } \wedge \bbk) & \ar[l] 
C^{*,\tilde{c}}(B \hat{Z}_{f_2 \circ f_1 \circ f_0}; \Milnor
\hat{\bar{X}}_{f_2 \circ f_1 \circ f_0}^{-V }
\wedge \bbk)
\end{tikzcd}
\end{equation}
and the composite of pullback maps
\begin{equation}
\begin{tikzpicture}
\node[scale=.8]{
\begin{tikzcd}
C^{*,\tilde{c}}(B \hat{Z}_{f_0}; \Milnor
\hat{\bar{X}}_{f_0}^{-V } \wedge \bbk) 
& C^{*,\tilde{c}}(B \hat{Z}_{f_1 \circ f_0}; \Milnor
\hat{\bar{X}}_{f_1 \circ f_0}^{-V } \wedge \bbk) \ar[l] & \ar[l]
C^{*,\tilde{c}}(B \hat{Z}_{f_2 \circ f_1 \circ f_0}; \Milnor
\hat{\bar{X}}_{{f_2 \circ f_1 \circ f_0}}^{-V } 
\wedge \bbk)
\end{tikzcd}
};
\end{tikzpicture}
\end{equation}
coincides with the pullback of the composite
\begin{equation}
\begin{tikzcd}
C^{*,\tilde{c}}(B \hat{Z}_{f_0}; \Milnor
\hat{\bar{X}}_{f_0}^{-V } \wedge \bbk) & \ar[l] 
C^{*,\tilde{c}}(B \hat{Z}_{f_2 \circ f_1 \circ f_0}; \Milnor
\hat{\bar{X}}_{f_2 \circ f_1 \circ f_0}^{-V }
\wedge \bbk)
\end{tikzcd}
\end{equation}
\end{lemma}
Thus, since the lemma implies that the compositions along the outside
of the square of Equation~\eqref{eq:composition_cube_diagram} coincide 
with the maps in the square
\begin{equation}
\begin{tikzcd}
C^{*,\tilde{c}}(B \hat{Z}_{{f}'' \circ h' \circ {h}}; \Milnor
\hat{\bar{X}}_{{f}'' \circ h' \circ {h}}^{-V } \wedge \bbk) \ar[d]  & \ar[l]
C^{*,\tilde{c}}(B \hat{Z}_{{f}}; \Milnor \hat{\bar{X}}_{{f}}^{-V }
\wedge \bbk) \ar[d]  \\ 
C^{*,\tilde{c}}(B \hat{Z}_{{f}''}; \Milnor \hat{\bar{X}}_{{f}''}^{-V }
\wedge \bbk) & \ar[l]  C^{*,\tilde{c}}(B \hat{Z}_{{g} \circ g' \circ
  {f}''}; 
\Milnor \hat{\bar{X}}_{{g} \circ g' \circ {f}''}^{-V } \wedge \bbk), 
  \end{tikzcd}
\end{equation}
we conclude that this assignment is compatible with the composition in
$\Tw \square {\Kur} \langle \cS \rangle$, which proves the following
lemma.

\begin{lem}\label{lem:Milnor_fibration_twisted_cochains_chart_functorial}
The assignment $f \mapsto   C^{*,\tilde{c}}(B \hat{Z}_{{f}}; \Milnor
\hat{\bar{X}}_{{f}}^{-V } \wedge \bbk)$ defines a functor $\Tw
\square {\Kur} \langle \cS \rangle \to \Sp$, and the maps in
Equation \eqref{eq:map_models_compact_support_Milnor} give rise to
natural transformations of the three models for compact
support. \qed 
\end{lem}

\subsubsection{Milnor-twisted cochains of flow categories}
\label{sec:miln-twist-coch}

As a first step to apply Lemma
\ref{lem:Milnor_fibration_twisted_cochains_chart_functorial} to flow
categories, we note that the inclusion of a stratum $Q \subset \cS$
gives rise to a map of cochains 
\begin{equation}
 C^{*,\tilde{c}}(B \partial^Q \hat{Z}_{{f}}; \Milnor\partial^Q  \hat{\bar{X}}_{{f}}^{-V } \wedge \bbk) \to  C^{*,\tilde{c}}(B \hat{Z}_{{f}}; \Milnor \hat{\bar{X}}_{{f}}^{-V } \wedge \bbk)
\end{equation}
induced by the inclusion of the collar $\kappa^Q \times \partial^Q \hat{Z}_{{f}}$ in $\hat{Z}_{{f}}$ and the aforementioned choice of section of the Milnor fibration of the collar.

The second step is to observe that the construction is multiplicative:
given a pair of arrows ${f} \in \square {\Kur} \langle \cS_1 \rangle $
and ${g} \in \square {\Kur} \langle \cS_2 \rangle$, the homeomorphism
$\hat{X}_{{f}} \times \hat{X}_{{g}} \cong \hat{X}_{{f} \times {g}}$
induces a product map 
\begin{multline} \label{eq:product_map_Milnor_cochains}
 C^{*,\tilde{c}}(B \hat{Z}_{{f}}; \Milnor \hat{\bar{X}}_{{f}}^{-V } \wedge \bbk) \wedge  C^{*,\tilde{c}}(B \hat{Z}_{{g}}; \Milnor \hat{\bar{X}}_{{g}}^{-V } \wedge \bbk) \\ \to  C^{*,\tilde{c}}(B \hat{Z}_{f \times g}; \Milnor \hat{\bar{X}}_{{f} \times {g}}^{-V } \wedge \bbk).
\end{multline}
Moreover, this product map is associative.

Putting these two ingredients together, we associate to each Kuranishi
flow category $\bX$ the following: 
\begin{defin}
The $\Pi$-equivariant
spectral category $C^{*,\tilde{c}}(B
\hat{\cZ}; \Milnor \hat{\bar{\cX}}^{-V -d} \wedge \bbk)$ has object set $\cP$ and morphism
spectra  
\begin{multline}
 C^{*,\tilde{c}}(B \hat{\cZ}; \Milnor \hat{\bar{\cX}}^{-V -d} \wedge
 \bbk)(p,q) \equiv \\ \hocolim_{f \in \square A(p,q) } C^{*,\tilde{c}}(B
 \hat{Z}_{{f}}; \Milnor \hat{\bar{X}}_{{f}}^{-V + V_{q} - V_{p}}
 \wedge \bbk)[ \deg p - \deg q].    
\end{multline}
The product map
\begin{equation}
\begin{tikzcd}
C^{*,\tilde{c}}(B \hat{\cZ}; \Milnor \hat{\bar{\cX}}^{-V -d} \wedge
\bbk)(p,q) \wedge  C^{*,\tilde{c}}(B \hat{\cZ}; \Milnor
\hat{\bar{\cX}}^{-V -d} \wedge \bbk)(q,r) \ar[d] \\
 C^{*,\tilde{c}}(B \hat{\cZ}; \Milnor \hat{\bar{\cX}}^{-V -d} \wedge
 \bbk)(p,r)
\end{tikzcd} 
\end{equation}
is defined by combining Equation
\eqref{eq:product_map_Milnor_cochains}, the boundary map associated to
the $q$-stratum of cubes of Kuranishi charts indexed by $A(p,r)$, and
the induced map of homotopy colimits. 
\end{defin}

Having constructed categories for each of the support conditions
$\tilde{c} \in \{c, c_0, c'_0\}$, the compatibility of the maps in
Equation~\eqref{eq:map_models_compact_support_Milnor} with the product
yields comparison functors:

\begin{lem}  
The comparison maps assemble to $\Pi$-equivariant spectral functors 
\begin{equation}
\begin{tikzcd}
C^{*,c_0}_{\rel \partial}(B \hat{\cZ};\Milnor \hat{\bar{\cX}}^{-V
   -d} \wedge \bbk) & \ar[l] \ar[d] C^{*,c'_0}_{\rel \partial}(B
 \hat{\cZ};\Milnor \hat{\bar{\cX}}^{-V -d} \wedge \bbk) \\
& C^{*,c}_{\rel \partial}(B \hat{\cZ};\Milnor \hat{\bar{\cX}}^{-V -d }
 \wedge \bbk) 
\end{tikzcd}
\end{equation}
which are DK-equivalences of spectral categories. \qed
\end{lem}

\subsection{Spanier-Whitehead duality}
\label{sec:span-whit-dual-1}

In Section~\ref{sec:span-whit-dual-5}, we constructed a comparison
between virtual cochains and cochains twisted by the Milnor model of
the tangent bundle.  The ingredients used were Spanier-Whitehead
duality, the Adams isomorphism, and the fact that the norm map with coefficients in Morava $K$-theory is an
equivalence. 

In this section, we realize the Spanier-Whitehead duality equivalence
as a functorial map, allowing us to prove an equivalence between
categories constructed from cochains twisted by the Milnor model and
the homotopy fixed points of the spaces underlying the construction of
virtual cochains. At the end, we explain how to apply the norm map and
the Adams isomorphism to complete the comparison, though the details
of most of the constructions (especially for the Adams isomorphism)
are postponed to Appendix~\ref{sec:adams-isom-ambid}.

\subsubsection{Functorial Spanier-Whitehead duality for manifolds}
\label{sec:funct-span-whit}

For expositional clarity, we begin by discussing the case of
manifolds. The essential problem is that, if $N$ is a finite group and
$X$ is a manifold which is a trivial cover of $X/N$ (with fiber $N$),
the diagram  
\begin{equation}
  \begin{tikzcd}
    X_+ \ar[r] \ar[d] & \Map_{X}(X, MX) \\
    X_+/N \ar[r] & \Map_{X/N}(X/N, MX/N) \ar[u]
  \end{tikzcd}
\end{equation}
does not commute. In fact, it fails to commute already for $X = N$. So
the comparison between the maps used to build the virtual cochains and
those used to build the twisted cochains requires a more elaborate
construction.

We introduce the based space   
\begin{equation}
 \cMap_{X/N}(X, S^0) \equiv X_+  \wedge_{N} \Map(N_+, S^0)
\end{equation}
whose points we think of as represented by a point in $X/N$, and a map
from the fiber to $S^0$. We have natural maps 
\begin{equation}
X_+ \to \cMap_{X/N}(X, S^0) \leftarrow X_+/N,
\end{equation}
where the rightward arrow is \emph{fibrewise Spanier-Whitehead duality}
\begin{equation}
X_+ \cong X_+ \wedge_{N} N \to X_+ \wedge_N \Map(N_+, S^0),
\end{equation}
induced by the map $\theta_n \colon N \to \Map(N_+, S^0)$ specified by
$\theta_n(x) = \ast$ unless $x = n$, and the leftward arrow is the
inclusion of constant maps with value on the distinguished orbit the
non-basepoint in $S^0$.

We also have a natural map
\begin{equation}
\cMap_{X/N}(X, S^0) \to \Map_X(X, MX)
\end{equation}
given by smashing the composition
\begin{equation}
\cMap_{X/N}(X, S^0) \to X/N  \to \Map_{X/N}(X/N, MX/N)  \to \Map_X(X, MX)
\end{equation}
with the map
\begin{equation}
\cMap_{X/N}(X, S^0) \to \Map(X, S^0)
\end{equation}
induced by the fact that $X \to X/N$ is assumed to be a trivial
cover.  The point of all of this is the following compatibility:

\begin{lem}
  The following diagram commutes
  \begin{equation}
    \begin{tikzcd}
    X_+ \ar[r] \ar[d]  & \Map_{X}(X, MX) \\
  \cMap_{X/N}(X, S^0) \ar[ur]  &  \\
  X_+/N \ar[r] \ar[u] & \Map_{X/N}(X/N, MX/N) \ar[uu]
    \end{tikzcd}
  \end{equation} \qed
\end{lem}

None of the map above are equivalences, but this can be addressed by
passing to the stable category. As an initial step, we consider the
spectrum 
\begin{equation}
  \ssF_{X/N}(X, \bS)  \equiv X_+  \wedge_{N} F(N_+, (\bS)^{\mfib}) .
\end{equation}
\begin{rem}
Note that the fibrant replacement functor used here involves only
trivial representations. This facilitates comparisons with the twisted
cochains, and computes the desired homotopy type because we shall
eventually pass to homotopy fixed points. In Appendix
\ref{sec:adams-isom-ambid}, we shall have occasion to use a fibrant
replacement functor involving all representations. 
\end{rem}

Once again, we have natural maps
\begin{equation}
\begin{tikzcd}
\Sigma^{\infty} X_+ \ar[r] & \ssF_{X/N}(X, \bS) & \ar[l]
\Sigma^{\infty} X_+/N
\end{tikzcd}
\end{equation}
induced by duality and the inclusion of the zero-sections; these are
underlying equivalences (see Lemma~\ref{lem:fiberwise-comparisons}).
There is also a natural map
\begin{equation}
\ssF_{X/N}(X, \bS) \to C^{*,c}(X; MX)
\end{equation}
defined as above; this is induced from the composite
\begin{equation}
X_+ \sma_{N} F(N_+, \bS^{\mfib}) \to X_+/N \to C^{*,c}(X/N; MX/N) \to
C^{*,c}(X; MX),
\end{equation}
where the first map is induced by the collapse map, the second by the
Spanier-Whitehead map (the inclusion of constant sections), and the
last by pullback of sections. 

The following lemma now expresses the compatibility of these
constructions; it is a specialization of the equivariant version we
actually use, which is stated below.

\begin{lem}\label{lem:Spanier-Whitehead-comparison}
The following diagram commutes, and the three maps in the top right
triangle are (underlying) equivalences of spectra: 
\begin{equation}
\begin{tikzcd}
\Sigma^{\infty} X_+ \ar[r] \ar[d] & C^{*,c}(X; MX) \\ \ssF_{X/N}(X,
\bS)  \ar[ur] &  \\ 
\Sigma^{\infty} X_+/N \ar[r] \ar[u] & C^{*,c}(X/N; MX/N) \ar[uu]
\end{tikzcd}
\end{equation}
\end{lem}

\begin{proof}
The bottom square commutes essentially by construction, and the top
triangle commutes from the definition of the Spanier-Whitehead duality
map, since the inclusion of constant sections factors through the
projection to $X_+/N$.  We have already explained why both the top
and downward arrow are underlying equivalences; the remaining arrow
immediately follows.
\end{proof}

The constructions above can be carried out assuming an ambient
$G$-action, extend to closed subsets as well as to manifolds with
boundary, and can be formulated with coefficients in any spectrum
$\bbk$.  Specifically, if we suppose that $X$ is a manifold with boundary
with an action of a finite group $G$, $N \subset G$ is a subgroup such
that $X \to X/N$ is a trivial cover, and $Z \subset X$ is a closed
subset, then we have a $G$-spectrum $\ssF_{X|Z/N}(X|Z,\bbk)$ defined
in analogy with non-equivariant version (see
Definition~\ref{defin:fiberwise-dual-total-space}).  The following
result follows via the same argument as
Lemma~\ref{lem:Spanier-Whitehead-comparison}.

\begin{lem} \label{lem:SW-functorial-closed-subset}
The following diagram of $G$-spectra commutes, and the three maps in
the top right triangle are equivalences of spectra:
  \begin{equation}
    \begin{tikzcd}
  \Sigma^{\infty} X|Z_+ \ar[r] \ar[d] & C^{*,c_0}_{\rel \partial}(Z;
  MX \wedge \bbk) \\ \ssF_{X|Z/N}(X|Z, \bbk) \ar[ur] &
  \\ \Sigma^{\infty} X|Z_+/N \ar[r] \ar[u] & C^{*,c_0}_{\rel
    \partial}(Z/N, MX/N \wedge \bbk) \ar[uu]
    \end{tikzcd}
  \end{equation} \qed
\end{lem}

We need to study the compatibility of this construction with
suspensions, products, and boundaries: first, if $S^{V^\perp}$ is a
$G$-sphere, the inclusion 
\begin{equation}
  S^{V^\perp} \wedge X|Z \to (S^{V^\perp} \sma X_+)|Z
\end{equation}
(see Equation \eqref{eq:inclusion-wedge-cone-to-cone-wedge}), induces a map
\begin{equation}
\begin{tikzcd}
S^{V^\perp} \wedge \ssF_{X|Z/N}(X|Z, \bbk) \ar[r] & \ssF_{(S^{V^\perp}
  \wedge X|Z)/N}(S^{V^\perp} \wedge X|Z, \bbk) \ar[d] \\
& \ssF_{(S^{V^\perp} \sma X_+)|Z/N}( (S^{V^\perp} \sma X_+)|Z, \bbk).
\end{tikzcd}
\end{equation}
On the other hand, since $Z$ lies away from the cone point of
$S^{V^\perp} \sma X_+$, we can make sense of the Milnor fibration $M
\left( S^{V^\perp} \sma X_+ \right)$ as a spherical fibration over
$Z$, and we have natural maps 
\begin{multline}
 \ssF_{(S^{V^\perp} \wedge X|Z)/N}(S^{V^\perp} \wedge X|Z, \bbk) \to
 C^{*,c_0}_{\rel \partial}(Z ; M (S^{V^\perp} \sma X_+) \wedge \bbk)
 \\ \leftarrow  \Sigma^{V^\perp} C^{*,c_0}_{\rel \partial}(Z ; M X \wedge
 \bbk). 
\end{multline}
Passing to the adjoint of these maps, and introducing a direct sum
decomposition $W = V \oplus V^\perp$, we have:
\begin{lem} \label{lem:SW-closed-subset-functorial-desuspend}
The following diagram commutes:
\begin{equation}
\begin{tikzcd}
  X|Z^{-V} \wedge \bbk \ar[r] \ar[d] & (S^{V^\perp} \sma X_+)|Z^{-W}  \wedge \bbk \ar[d] \\
  \ssF_{X|Z/N}(X|Z, \bbk)^{-V} \ar[r] \ar[d]  &  \ssF_{(S^{V^\perp} \sma X_+)|Z/N}( (S^{V^\perp}_+ \sma X)|Z,
 \bbk)^{-W} \ar[d] \\
 C^{*,c_0}_{\rel \partial}(Z; MX^{-V} \wedge \bbk)
 \ar[r] & C^{*,c_0}_{\rel \partial}(Z ; M (S^{V^\perp}_+
 \sma X)^{-W} \wedge \bbk). 
\end{tikzcd}
\end{equation} \qed
\end{lem}

Next, we assume that $X$ is stratified, and consider the inclusion of
a boundary stratum $\partial^Q X \subset X$. Passing to the collared
completion as in Section~\ref{sec:span-whit-dual-4}, and using the
modified map from $\hat{X}|\hat{Z}$ to $C^{*,c_0}(\hat{Z}, M \hat{X})$
from that section yields: 
\begin{lem} \label{lem:SW-closed-subset-functorial-boundary}
The following diagram of $G$-spectra commutes:
\begin{equation}
\begin{tikzcd}
\partial^Q \hat{X}| \partial^Q \hat{Z} \ar[r] \ar[d] &  \ssF_{\partial^Q \hat{X}| \partial^Q \hat{Z}/N}(\partial^Q \hat{X}| \partial^Q \hat{Z}, \bbk)  \ar[r] \ar[d]  & C^{*,c_0}_{\rel \partial}(\partial^Q  \hat{Z}; M\hat{X}) \ar[d] \\
 \hat{X}|\hat{Z} \ar[r]  &  \ssF_{\hat{X}|\hat{Z}/N}(\hat{X}|\hat{Z}, \bbk) \ar[r]  &  C^{*,c_0}_{\rel \partial}(\hat{Z}; M\hat{X})
\end{tikzcd}
\end{equation} \qed
\end{lem}

Finally, we note the compatibility of this construction with products:
for each pair of triples $(X_i, Z_i, N_i)$ for $i = 0,1$, we have a
natural map
\begin{equation}
\begin{tikzcd}
\ssF_{X_0| Z_0 /N_0}(X_0|Z_0 , \bbk) \sma \ssF_{X_1| Z_1
     /N_1}(X_1|Z_1 , \bbk) \ar[d] \\
\ssF_{X_0 \times X_1| Z_0 \times Z_1
     /(N_0 \times N_1)}( X_0 \times X_1| Z_0 \times Z_1, \bbk).
\end{tikzcd}
\end{equation}
given by the smash product in the base and the fiber (this uses the
fact that the fibrant replacement functor we use is multiplicative).
As discussed in Lemma~\ref{lem:fibmapmult}, these maps are associative
and unital.

\begin{lem} \label{lem:SW-duality-closed-subset-monoidal}
The diagrams
\begin{equation}
\begin{tikzcd}
X_0|Z_0 \wedge X_1|Z_1  \ar[r] \ar[d] &  \ssF_{X_0| Z_0 /N_0}(X_0|Z_0
, \bbk) \wedge \ssF_{X_1| Z_1 /N_1}(X_1|Z_1 , \bbk) \ar[d]  \\ 
X_0 \times X_1| Z_0 \times Z_1 \ar[r]  &  \ssF_{X_0 \times X_1| Z_0
  \times Z_1 /(N_0 \times N_1)}( X_0 \times X_1| Z_0 \times Z_1,
\bbk). 
\end{tikzcd}
\end{equation}
and
\begin{equation}
\begin{tikzcd}[row sep=small, column sep=tiny, transform shape, nodes={scale=.85}]
\ssF_{X_0| Z_0 /N_0}(X_0|Z_0 , \bbk) \wedge \ssF_{X_1| Z_1
  /N_1}(X_1|Z_1 , \bbk)  \ar[r] \ar[d]  & C^{*,c_0}_{\rel
  \partial}(Z_0; MX_0) \wedge C^{*,c_0}_{\rel \partial}(Z_1; MX_1)
\ar[d] \\ 
\ssF_{X_0 \times X_1| Z_0 \times Z_1 /(N_0 \times N_1)}( X_0 \times
X_1| Z_0 \times Z_1, \bbk) \ar[r]  &  C^{*,c_0}_{\rel \partial}( Z_0
\times Z_1; M X_0 \times X_1) 
\end{tikzcd}
\end{equation}
commute.
\end{lem}

\begin{proof}
The first square commutes by Lemma~\ref{lem:SW-map-multiplicative}.
That lemma also implies that the second square commutes, when combined
with the easy check that the natural inclusion of sections is
multiplicative.
\end{proof}

\subsubsection{Functorial Spanier-Whitehead duality for Kuranishi charts}
\label{sec:funct-span-whit-1}

In the previous section, we assembled all the properties of the
Spanier-Whitehead equivalence that are required to set up a functorial
natural transformation from the virtual cochains of Kuranishi
presentations to the cochains of the zero-loci twisted by the Milnor
model of the tangent fibration. In order to formulate the comparison,
we consider a refinement of the category of cubes: 

\begin{defin}
The category $\square_{\Sub} \Kur \langle \cS \rangle$ of
\emph{cubes of Kuranishi charts with a choice of freely acting
  subgroup} is the category with 
\begin{enumerate}
\item objects given by pairs $(\sigma, N)$,
where $\sigma$ is a cube of inessential Kuranishi charts 
and $N \subset G_\sigma$ is a subgroup acting freely on
$X_\sigma$, and
\item morphisms from $(\sigma,N)$ to $(\sigma',N')$ given by
a map $f \colon \sigma \to \sigma'$ such that $N'$ is contained in
$f(N)$.
\end{enumerate}
\end{defin}

Consider the category $\Sp_{\eq}$ of equivariant spectra, which can be thought of as a
stabilization of the category $\Orbo$; objects consist of pairs $(G,
X)$ where $G$ is a finite group and $X$ is a $G$-spectrum.  We now
construct a functor from $\Tw \square_{\Sub} \Kur\langle \cS \rangle$
to the category of equivariant spectra which is given on objects by
the assignment
\begin{equation} \label{eq:functor-Kuranishi-free-action-spectra-virtual-cochains-uncompleted}
f \mapsto (G_f, \ssF_{\frac{{\bar{X}}_f|{Z}_f}{N'}}({\bar{X}}_f|{Z}_f, \bbk)^{-V_f} )
\end{equation}
where $f$ is an arrow from $(\sigma,N)$ to $(\sigma',N')$, and we
recall that $V_f \equiv V_\sigma$. Observe that $ \frac{X_f|Z_f}{N'}$
is also the quotient of $X|Z$ by the free action of $f^{-1}(N')$. On
morphisms, the functor is given by the following special case of Lemma
\ref{lem:fibrewise-mapping-spectra-functorial}: 

\begin{lem}\label{lem:fibrewise-mapping-spectra-functorial-special} 
Each factorization
\begin{equation}
\begin{tikzcd}
(\sigma_0, N_0) \ar[r, "f_0"] \ar[d,"g"] & (\sigma'_0, N'_0) \\
(\sigma_1, N_1)  \ar[r,"f_1"] & (\sigma'_1, N'_1) \ar[u,"h"] 
\end{tikzcd}
\end{equation}
induces a natural map
\begin{equation} \label{eq:map_internal-mapping-spaces-twisted-arrow-S^0}
\ssF_{{\bar{X}}_{f_0}|{Z}_{f_0}/N'_0}({\bar{X}}_{f_0}|{Z}_{f_0} ,
\bbk)^{-V_{f_0}} \to
\ssF_{{\bar{X}}_{f_1}|{Z}_{f_1}/N'_1}({\bar{X}}_{f_1}|{Z}_{f_1},
\bbk)^{-V_{f_1}}.    
\end{equation}
\end{lem}

\begin{proof}
The orthogonal complement of the inclusion $V_{f_0} \to V_{f_1}$ is
$V_{g}^{\perp}$. By adjunction, it suffices to construct a map 
\begin{equation}
S^{V_g^\perp} \wedge
\ssF_{{\bar{X}}_{f_0}|{Z}_{f_0}/N'_0}({\bar{X}}_{f_0}|{Z}_{f_0} ,
\bbk) \to
\ssF_{{\bar{X}}_{f_1}|{Z}_{f_1}/N'_1}({\bar{X}}_{f_1}|{Z}_{f_1},
\bbk).
\end{equation}
Recall that $\bX_{f_0}$ is the quotient of $\bX_{\sigma_0}$ by
$G_{f_0}^{\perp}$. Factoring $f_0 = h \circ f_1 \circ g$, we can
describe  $\bX_{f_0}$ instead as the quotient of $\bX_{g}$ by $G_{h
  \circ f_1}^{\perp}$. Recalling that the product of $X_g$ with
$V_g^{\perp}$ has a natural embedding in $X_{\sigma_1}$, we have a
natural map  
\begin{equation}\label{eq:map_desuspended_fibrewise-spectra-factor-thickening-0}
\begin{tikzcd}
S^{V_g^\perp} \wedge
\ssF_{{\bar{X}}_{f_0}|{Z}_{f_0}/N'_0}({\bar{X}}_{f_0}|{Z}_{f_0} ,
\bbk) \ar[d] \\
\ssF_{(S^{V_g^\perp} \sma  ({\bar{X}}_{g})_+)|{Z}_{g}/(h
  \circ f_1)^{-1}N'_0}((S^{V_g^\perp} \sma  ({\bar{X}}_{g})_+)|{Z}_{g}
/ G_{h \circ f_1}^{\perp} , \bbk). 
\end{tikzcd}
\end{equation}
Thus, we need to construct a further map from the target of the above
map to $\ssF_{{\bar{X}}_{f_1}|{Z}_{f_1}/N'_1}({\bar{X}}_{f_1}|{Z}_{f_1},
\bbk)$.

We first construct a map $V_g^{\perp} \times (X_{f_0}/N'_0) \to
X_{f_1} /N'_1 $: we have a natural map 
\begin{equation}
V_g^{\perp} \times X_{g} \to X_{\sigma_1} \to X_{f_1},
\end{equation}
so it suffices to show that this map descends to the quotients. This
follows from the inclusion $h^{-1}(N'_0) \subset N'_1$, and the
description of $  X_{f_1} /N'_1$ as the quotient of $X_{\sigma_1}$ by
the inverse image of $N'_1$ under $f_1$, and of $X_{f_0}/N'_0$ as the
quotient of $X_{g}$ by the inverse image of $N'_0 $ under $h \circ
f_1$. 

Passing to fibrewise compactifications and taking the cone on the
complement of zero-loci, we obtain a map 
\begin{equation}
(S^{V_g^\perp} \sma  ({\bar{X}}_{g})_+)|{Z}_{g}/(h \circ f_1)^{-1}N'_0
  \to  {\bar{X}}_{f_1}|{Z}_{f_1}/N'_1. 
\end{equation}
 
To complete the construction, we must relate the quotient of
$(S^{V_g^\perp} \sma  ({\bar{X}}_{g})_+)|{Z}_{g} $ by $ G_{h \circ
  f_1}^{\perp} $ (i.e., the total space of the fibrewise spectrum of
maps in the target of
Equation~\eqref{eq:map_desuspended_fibrewise-spectra-factor-thickening-0}) 
to the subset of ${\bar{X}}_{f_1}|{Z}_{f_1}$ lying over the image of
the above map. This subset is precisely the image of $(S^{V_g^\perp}
\sma ({\bar{X}}_{g})_+)|{Z}_{g}$ under $f_1$, i.e., 
\begin{multline}
(S^{V_g^\perp} \sma ({\bar{X}}_{g})_+)|{Z}_{g}/(h \circ f_1)^{-1}N'_0
  \times_{{\bar{X}}_{f_1}|{Z}_{f_1}/N'_1} {\bar{X}}_{f_1}|{Z}_{f_1}
 \\ \cong   (S^{V_g^\perp} \sma ({\bar{X}}_{g})_+)|{Z}_{g} /
  G_{f_1}^{\perp}. 
\end{multline}
The evident inclusion $G_{f_1}^{\perp} \to G_{h \circ f_1}^{\perp}$
yields a covering map 
\begin{equation}\label{eq:map-covering}
(S^{V_g^\perp} \sma  ({\bar{X}}_{g})_+)|{Z}_{g} /  G_{f_1}^{\perp}  \to
  (S^{V_g^\perp} \sma  ({\bar{X}}_{g})_+)|{Z}_{g} /  G_{h \circ
    f_1}^{\perp}, 
\end{equation}
and hence a pullback map
\begin{equation} \label{eq:pullback-map-fibrewise-functors}
\begin{tikzcd}
    \ssF_{(S^{V_g^\perp} \sma  ({\bar{X}}_{g})_+)|{Z}_{g}/(h \circ
      f_1)^{-1}N'_0}((S^{V_g^\perp} \sma  ({\bar{X}}_{g})_+)|{Z}_{g} /
    G_{h \circ f_1}^{\perp} , \bbk) \ar[d] \\ 
\ssF_{(S^{V_g^\perp} \sma ({\bar{X}}_{g})_+)|{Z}_{g}/(h \circ f_1)^{-1}N'_0}((S^{V_g^\perp}
    \sma  ({\bar{X}}_{g})_+)|{Z}_{g} / G_{f_1}^{\perp} , \bbk).
\end{tikzcd}
\end{equation}
On the other hand, the inclusion of $(S^{V_g^\perp} \sma
({\bar{X}}_{g})_+)|{Z}_{g} /  G_{f_1}^{\perp}$ in
${\bar{X}}_{f_1}|{Z}_{f_1}$ yields an inclusion 
\begin{equation}\label{eq:inclusion-map}
\begin{tikzcd}
\ssF_{(S^{V_g^\perp} \sma ({\bar{X}}_{g})_+)|{Z}_{g}/(h \circ
  f_1)^{-1}N'_0}((S^{V_g^\perp} \sma ({\bar{X}}_{g})_+)|{Z}_{g} /
G_{f_1}^{\perp} , \bbk) \ar[d] \\
\ssF_{{\bar{X}}_{f_1}|{Z}_{f_1}/N'_1}({\bar{X}}_{f_1}|{Z}_{f_1},
\bbk).
\end{tikzcd}
\end{equation}
Thus, composing Equation~\eqref{eq:inclusion-map} with
Equation~\eqref{eq:pullback-map-fibrewise-functors} and
Equation~\eqref{eq:map_desuspended_fibrewise-spectra-factor-thickening-0} 
yields the result. 
\end{proof}

This assignment on morphisms is compatible with the composition in the
twisted arrow category.

\begin{prop}
The assignment
\begin{equation} \label{eq:functor-Kuranishi-free-action-spectra-virtual-cochains-uncompleted-prop}
    f \mapsto (G_f, \ssF_{\frac{{\bar{X}}_f|{Z}_f}{N'}}({\bar{X}}_f|{Z}_f, \bbk)^{-V_f} )
\end{equation}
specifies a functor
\begin{equation}
\Tw \square_{\Sub} \Kur\langle \cS \rangle \to \Sp_{\eq}
\end{equation}
\end{prop}

\begin{proof}
Given factorizations $f_0 = \id \circ f_1 \circ g_0$ and $f_1 = \id
\circ f_2 \circ g_1$, we need to show that the diagram
\begin{equation}
\begin{tikzcd}[column sep=tiny]
  \begin{gathered}
S^{V_{g_1}^\perp} \sma S^{V_{g_0}^\perp} \sma  \\
\ssF_{{\bar{X}}_{f_0}|{Z}_{f_0}/N'_0}({\bar{X}}_{f_0}|{Z}_{f_0} ,
\bbk)    
  \end{gathered}
 \ar[d] \ar[r] & S^{V_{g_1}^\perp} \sma
\ssF_{{\bar{X}}_{f_1}|{Z}_{f_1}/N'_1}({\bar{X}}_{f_1}|{Z}_{f_1}, \bbk)
\ar[dd] \\
\begin{gathered}
S^{V_{g_1}^\perp \oplus V_{g_0}^\perp} \sma  \\
\ssF_{{\bar{X}}_{f_0}|{Z}_{f_0}/N'_0}({\bar{X}}_{f_0}|{Z}_{f_0} ,
\bbk)  
\end{gathered}
 \ar[d] & \\
S^{V_{g_1 \circ g_0}^\perp} \sma 
\ssF_{{\bar{X}}_{f_0}|{Z}_{f_0}/N'_0}({\bar{X}}_{f_0}|{Z}_{f_0} ,
\bbk) \ar[r] &
\ssF_{{\bar{X}}_{f_2}|{Z}_{f_2}/N'_2}({\bar{X}}_{f_2}|{Z}_{f_2}, \bbk)
\end{tikzcd}
\end{equation}
commutes, where here we are using the identification $V_{g_1}^\perp
\oplus V_{g_0}^\perp \cong V_{g_1 \circ g_0}^{\perp}$.

As discussed in Section~\ref{sec:twisted}, it suffices to check that
the assignment to arrows is compatible with composition in a handful
of special cases.  We begin by considering the composite expressed by
the factorizations $f_0 = \id \circ f_1 \circ g_0$ and $f_1 = \id
\circ f_2 \circ g_1$.  In this case, the natural maps of
Equation~\eqref{eq:map_desuspended_fibrewise-spectra-factor-thickening-0}
are induced in part by the identification of $X_{f_0} = X_{g_0} /
G_{f_1}^{\perp}$, $X_{f_1} = X_{g_1} / G_{f_2}^{\perp}$, and $X_{f_0}
= X_{g_0} / G_{f_1 \circ f_0}^{\perp}$.

Fixing $h = \id$ in the construction of
Lemma~\ref{lem:fibrewise-mapping-spectra-functorial-special}, observe
that we are considering maps induced by
\begin{equation}
V_{g_0}^{\perp} \times X_{g_0} \to X_{\sigma_1} \to X_{f_1}
\end{equation}
and
\begin{equation}
V_{g_1}^{\perp} \times X_{g_1} \to X_{\sigma_2} \to X_{f_2}
\end{equation}
with
\begin{equation}
V_{g_1 \circ g_0}^{\perp} \times X_{g_1 \circ g_0} \to X_{\sigma_2
  \circ \sigma_1} \to X_{f_2}.
\end{equation}
Here note that $V_{g_1 \circ g_0}^{\perp}$ is the orthogonal
complement of the image of $V_{g_1}^{\perp}$.

Because the covering maps of Equation~\eqref{eq:map-covering} is the
identity in this case, the pullback maps of
Equation~\eqref{eq:pullback-map-fibrewise-functors} are isomorphisms as well and so
the relevant composite is determined by the inclusions of
Equation~\eqref{eq:inclusion-map}.  Putting this all together, the two
maps coincide, using the functoriality of cone construction and the
fact that the inclusions do.

Next, we consider the composite expressed by the factorizations $f_0 =
h_0 \circ f_1 \circ \id$ and $f_1 = h_1 \circ f_2 \circ \id$.  In this
case, the maps of
Equation~\eqref{eq:map_desuspended_fibrewise-spectra-factor-thickening-0}
are determined by the identifications of $\bX_{f_0}$ as the quotient
of $\bX$ by $G_{h_0 \circ f_1}^{\perp}$ and as the quotient of $\bX$
by $G_{h_1 \circ h_0 \circ f_2}^{\perp}$, and $\bX_{f_1}$ as the
quotient of $\bX$ by $G_{h_1 \circ f_2}^{\perp}$.

Here, the maps of
Equation~\eqref{eq:map_desuspended_fibrewise-spectra-factor-thickening-0} 
are induced by the maps $X_{\id} \to X_{\sigma_0} \to X_{f_1}$,
$X_{\id} \to X_{\sigma_0} \to X_{f_2}$, and $X_{\id} \to X_{\sigma_0}
\to X_{f_2 \circ f_1}$.
In this case, when considering the covering maps of
Equation~\eqref{eq:map-covering}, we are considering the coverings
associated to the inclusions $G_{f_1}^{\perp} \to G_{h_0 \circ
  f_1}^{\perp}$, $G_{f_2}^{\perp} \to G_{h_1 \circ f_2}^{\perp}$, and
$G_{f_2}^{\perp} \to G_{h_0 \circ h_1 \circ f_2}$.  Since $f_1 = h_1
\circ f_2$, the composite of the first two maps
\begin{equation}
G_{f_2}^{\perp} \to G_{h_1 \circ f_2}^{\perp} = G_{f_1}^{\perp} \to
G_{h_0 \circ f_1}^{\perp} = G_{h_0 \circ h_1 \circ f_2}^{\perp}
\end{equation}
coincides with the last map, and so the composite of the associated
coverings coincides with the covering of the last map.

Finally, we consider a composite expressed by the factorizations $f_0 =
h_0 \circ f_1 \circ g_0$ and $f_1 = h_1 \circ f_2 \circ \id$.  In
this case, small variations on the arguments above suffice to check
that the assignment is compatible with the composition.
\end{proof}

Passing to collared completions, the same construction yields another
functor.

\begin{prop}
There is a functor from $\Tw \square_{\Sub} \Kur\langle \cS \rangle$
to $\Sp_{\eq}$ specified by the assignment on objects
\begin{equation}\label{eq:functor-Kuranishi-free-action-spectra-virtual-cochains}
f \mapsto (G_f,
\ssF_{\frac{\hat{\bar{X}}_f|\hat{Z}_f}{N'}}(\hat{\bar{X}}_f|\hat{Z}_f,
\bbk)^{-V_f} ),
\end{equation}
equipped with a natural transformation to the functor specified by
Equation~\eqref{eq:functor-Kuranishi-free-action-spectra-virtual-cochains-uncompleted}. \qed
\end{prop}

Forgetting the subgroup $N$ induces a natural functor
\begin{equation}
\Tw \square_{\Sub} \Kur \langle \cS \rangle\to \Tw \square \Kur\langle \cS \rangle.
\end{equation}
Thus precomposing pulls back the functor $\Tw \square \Kur\langle \cS \rangle
\to \Sp_{\eq}$ specified by the assignment 
\begin{equation}
f \mapsto (G_f, C^{*,c_0}_{\rel \partial}(\hat{Z}_f, M \hat{\bar{X}}^{-V}_f
\wedge \bbk))
\end{equation}
(described in Lemma~\ref{lem:Milnor_fibration_twisted_cochains_chart_functorial}) to
a composite functor 
\begin{equation}
\Tw \square_{\Sub} \Kur\langle \cS \rangle \to 
\Tw \square \Kur\langle \cS \rangle \to \Sp_{\eq}.
\end{equation}
For each object $f$ of $\Tw \square_{\Sub} \Kur\langle \cS \rangle$,
Lemma~\ref{lem:SW-functorial-closed-subset} yields a map 
\begin{equation}\label{eq:cochain-comparison}
\ssF_{\frac{\hat{\bar{X}}_f|\hat{Z}_f}{N'}}(\hat{\bar{X}}_f|\hat{Z}_f,
\bbk)^{-V_f}   \to C^{*,c_0}_{\rel \partial}(\hat{Z}_f, M
\hat{\bar{X}}^{-V}_f\wedge \bbk) 
\end{equation}
of $G_f$-spectra.

Lemmas~\ref{lem:SW-closed-subset-functorial-desuspend}
and~\ref{lem:SW-closed-subset-functorial-boundary} generalize to the
following result.

\begin{lem}
The map of Equation~\eqref{eq:cochain-comparison} induces a natural
transformation of functors $\Tw \square_{\Sub} \Kur\langle \cS \rangle
\to \Sp_{\eq}$, i.e., for each factorization $f_0 = h \circ f_1 \circ
g$, the following diagram commutes:
\begin{equation}
\begin{tikzcd}
\ssF_{{\bar{X}}_{f_0}|{Z}_{f_0}/N'_0}({\bar{X}}_{f_0}|{Z}_{f_0} ,
\bbk)^{-V_{f_0}} \ar[r] \ar[d] & C^{*,c_0}_{\rel
  \partial}(\hat{Z}_{f_0}, M \hat{\bar{X}}^{-V_0}_{f_1} \wedge \bbk)
\ar[d] \\ 
\ssF_{{\bar{X}}_{f_1}|{Z}_{f_1}/N'_1}({\bar{X}}_{f_1}|{Z}_{f_1},
\bbk)^{-V_{f_1}} \ar[r] & C^{*,c_0}_{\rel \partial}(\hat{Z}_{f_1}, M
\hat{\bar{X}}^{-V_1}_{f_1} \wedge \bbk). 
\end{tikzcd}
\end{equation}
\end{lem}

\begin{proof}
We follow the steps in the construction of the map from Lemma
\ref{lem:fibrewise-mapping-spectra-functorial-special}: we start with
the commutative diagram 
\begin{equation}
\begin{tikzcd}[column sep=small, transform shape, nodes={scale=.85}]
S^{V_g^\perp} \wedge    \ssF_{{\hat{\bar{X}}}_{f_0}|{\hat{Z}}_{f_0}/N'_0}({\hat{\bar{X}}}_{f_0}|{\hat{Z}}_{f_0} , \bbk)  \ar[r] \ar[d] & C^{*,c_0}_{\rel \partial}(\hat{Z}_{f_0}, S^{V_g^\perp} \wedge  M \hat{\bar{X}}_{f_0} \wedge \bbk) \ar[d] \\
   \ssF_{(S^{V_g^\perp} \sma ({\hat{\bar{X}}}_{g})_+)|{\hat{Z}}_{g}/(h \circ f_1)^{-1}N'_0}((S^{V_g^\perp} \sma  ({\hat{\bar{X}}}_{g})_+)|{\hat{Z}}_{g} / G_{h \circ f_1}^{\perp} , \bbk)   \ar[r] \ar[d] & C^{*,c_0}_{\rel \partial}(\hat{Z}_{g}, S^{V_g^\perp} \sma (M \hat{\bar{X}}_{g})_+ \wedge \bbk) \\
    \ssF_{(S^{V_g^\perp} \sma  ({\bar{X}}_{g})_+)|{Z}_{g}/(h \circ
      f_1)^{-1}N'_0}((S^{V_g^\perp} \sma  ({\bar{X}}_{g})_+)|{Z}_{g} /
    G_{f_1}^{\perp} , \bbk)  \ar[ur] &  
    \end{tikzcd}
  \end{equation}
for which the commutativity of the top square arises from the
compatibility of our Spanier-Whitehead map with smashing with the
sphere $S^{V_g^\perp}$, and the commutativity of the bottom triangle
from its compatibility with pullback.  Next, we have a commutative
diagram  
\begin{equation}
\begin{tikzcd}[column sep=small, transform shape, nodes={scale=.9}]
 \ssF_{(S^{V_g^\perp} \sma ({\bar{X}}_{g})_+)|{Z}_{g}/(h \circ f_1)^{-1}N'_0}((S^{V_g^\perp} \sma ({\bar{X}}_{g})_+)|{Z}_{g} / G_{f_1}^{\perp} , \bbk)    \ar[d]   \ar[r] & C^{*,c_0}_{\rel \partial}(\hat{Z}_{g}, S^{V_g^\perp} \sma (M \hat{\bar{X}}_{g})_+ \wedge \bbk) \ar[d] \\
 \ssF_{{\bar{X}}_{f_1}|{Z}_{f_1}/N'_1}({\bar{X}}_{f_1}|{Z}_{f_1}, \bbk) \ar[r] &  C^{*,c_0}_{\rel \partial}(\hat{Z}_{f_1}, M \hat{\bar{X}}_{f_1} \wedge \bbk),
\end{tikzcd}
\end{equation}
given by the compatibility with pushforward.
\end{proof}

On the other hand, Lemma~\ref{lem:SW-duality-closed-subset-monoidal}
implies that these maps are compatible with products. To state the
compatibility, we write $\square \Kur $ for the following category
which combines the categories $\square \Kur \langle \cS \rangle$ as
the poset $\cS$ varies.  Note that this is not the category of cubes
on $\Kur$; the definition below reflects the restrictions on the
interaction of the cubes with the posets in the combined category that
arise in our setting.

\begin{defin}\label{defin:augmented-Kur-cubes}
Let $\square \Kur $ denote the category with
\begin{enumerate}
\item objects cubes in $\Kur \langle \cS \rangle$ for some partially
  ordered set $\cS$, and
\item morphisms from $\square \tau_0 \in \Kur \langle \cS_0 \rangle$ to
$\square \sigma_1 \in \Kur \langle \cS_1 \rangle$ given by a map
  of cubes $\tau_1 \to \sigma_1$ in $\Kur \langle \cS_1 \rangle$, an
identification $\cS_0 \cong \cS_1 \setminus Q$ for a totally ordered
set $Q$, and an identification $\tau_0 \cong \partial^Q \tau_1$ of
cubes in $\langle \cS_0 \rangle $, where $\partial^Q \tau_1 $ is the
restriction to the stratum labelled by $Q$. 
\end{enumerate}
We define $\square_{\Sub} \Kur$ to be the category which has objects
those of $\square \Kur$ along with the additional choice for each
chart of a freely acting subgroup $N$ of $G$; for morphisms we impose
the condition that $N_1$ be contained in the image of $N_0$. 
\end{defin}

It is straightforward to check that $\square \Kur$ is a monoidal
category under the evident product induced by the product of cubes and
of charts in $\Kur \langle \cS \rangle$.  We then have the following
lemma:

\begin{lem}
Given an object $f$ of $\Tw \square_{\Sub} \Kur$, the assignments 
\begin{equation}
f \mapsto (G_f, \ssF_{\frac{\hat{\bar{X}}_f|\hat{Z}_f}{N'}}(\hat{\bar{X}}_f|\hat{Z}_f,
  \bbk)^{-V_f})
\end{equation}
and
\begin{equation}
f \mapsto (G_f, C^{*,c_0}_{\rel
    \partial}(\hat{Z}_f, M \hat{\bar{X}}^{-V}_f\wedge \bbk)) 
\end{equation}
specify lax monoidal functors $\Tw \square_{\Sub} \Kur \to
\Sp_{\eq}$. \qed
\end{lem}

Finally, we can now record the comparison induced by the
Spanier-Whitehead duality map; this is a direct consequence of
Lemma~\ref{lem:SW-duality-closed-subset-monoidal}.

\begin{lem} \label{SW-duality-lax-monoidal}
Spanier-Whitehead duality defines a lax monoidal natural
transformation 
\begin{equation}
(G_f,
  \ssF_{\frac{\hat{\bar{X}}_f|\hat{Z}_f}{N'}}(\hat{\bar{X}}_f|\hat{Z}_f,
  \bbk)^{-V_f})       \Rightarrow  (G_f, C^{*,c_0}_{\rel
    \partial}(\hat{Z}_f, M \hat{\bar{X}}^{-V}_f\wedge \bbk)) 
\end{equation}
of functors from $\Tw \square_{\Sub} \Kur$ to $\Sp_{\eq}$. \qed
\end{lem}

\subsubsection{Spanier-Whitehead duality for flow categories}
\label{sec:span-whit-dual-6}

We now consider a Kuranishi flow category $\bX \co A \to \Kur$ as
before, where $A$ is a $2$-category.  In
Section~\ref{sec:miln-twist-coch}, we used the $2$-category $\square
A$ to define twisted cochains.  Here, we begin by introducing the
$2$-category $\square_{\Sub} A$, with the categories of $1$-cells 
$\square_{\Sub}A(p,q)$ consisting of the category with objects a
cube $\sigma$ in $A(p,q)$ together with a subgroup $N$ of $G_\sigma$
acting freely on $X_\sigma$, and maps such that $N$ is contained in
$N'$; to be precise, this specifies a bicategory, and we are
implicitly rectifying as in the discussion surrounding
Theorem~\ref{thm:2catrect}.  It is straightforward to see that $\bX$
lifts to a natural $\Pi$-equivariant $2$-functor 
\begin{equation}
\bX \colon \square_{\Sub} A \to \square_{\Sub} \Kur
\end{equation}
that induces a $\Pi$-equivariant $2$-functor
\begin{equation}
\bX \colon \Tw \square_{\Sub} A \to \Tw \square_{\Sub} \Kur.
\end{equation}

Composing with the functor of (fibrewise) virtual cochains from
Equation~\eqref{eq:functor-Kuranishi-free-action-spectra-virtual-cochains},
and taking homotopy fixed points, for each pair of objects $p$ and $q$
we obtain a functor 
\begin{align}
\Tw \square_{\Sub}A(p,q) & \to \Sp \\
f & \mapsto C^{*}(B G_f,
\ssF_{\frac{\hat{\bar{X}}_f|\hat{Z}_f}{N'}}(\hat{\bar{X}}_f|\hat{Z}_f,
\bbk)^{-V_f} ),
\end{align}
where $f$ is an arrow with target $(\sigma', N')$.  Taking the
homotopy colimit over this functor, and shifting by the degrees of $p$
and $q$ yields a spectrum.  Since the work above shows that these
functors are part of the data of a $\Pi$-equivariant $2$-functor, the
homotopy colimits assemble to yield a spectral category.

\begin{defin}
For a Kuranishi flow category $\bX$, we define the $\Pi$-equivariant
spectral category 
$C^{*}(B G,
\ssF_{\frac{\hat{\bar{\cX}}|\hat{\cZ}}{N}}(\hat{\bar{\cX}}|\hat{\cZ},
\bbk)^{-V-d} )$ to have objects $\cP$ and morphism spectra 
\begin{multline}
C^{*}(B G,
\ssF_{\frac{\hat{\bar{\cX}}|\hat{\cZ}}{N}}(\hat{\bar{\cX}}|\hat{\cZ},
\bbk)^{-V-d} )  (p,q) \equiv  \\ 
\hocolim_{f \in \Tw   \square_{\Sub}A(p,q) } C^{*}(B G_f,
\ssF_{\frac{\hat{\bar{X}}_f|\hat{Z}_f}{N'}}(\hat{\bar{X}}_f|\hat{Z}_f,
\bbk)^{-V_f +  V_{q} - V_{p} } )[ \deg p  - \deg q ], 
\end{multline}
\end{defin}

We also have an uncollared version.

\begin{defin}
For a Kuranishi flow category $\bX$, we define the $\Pi$-equivariant
spectral category 
$C^{*}(B G,
\ssF_{\frac{{\bar{\cX}}|{\cZ}}{N}}({\bar{\cX}}|{\cZ}, \bbk)^{-V-d} )$
to have objects $\cP$ and morphism spectra
\begin{multline} \label{eq:homotopy-fixed-points-virtual-cochains-uncollared}
  C^{*}(B G, \ssF_{\frac{{\bar{\cX}}|{\cZ}}{N}}({\bar{\cX}}|{\cZ}, \bbk)^{-V-d} )  (p,q) \equiv  \\
  \hocolim_{f \in \Tw   \square_{\Sub}A(p,q) } C^{*}(B G_f, \ssF_{\frac{{\bar{X}}_f|{Z}_f}{N'}}({\bar{X}}_f|{Z}_f, \bbk)^{-V_f+  V_{q} - V_{p}} )[ \deg p  - \deg q ].
\end{multline}
\end{defin}

For each $p$ and $q$, there is an evident comparison between the
mapping spectra of the collared and uncollared versions, which is
given by the projection map 
\begin{multline} \label{eq:homotopy-fixed-points-virtual-cochains-uncollared-comparison}
C^{*}(B G,
\ssF_{\frac{\hat{\bar{\cX}}|\hat{\cZ}}{N}}(\hat{\bar{\cX}}|\hat{\cZ},
\bbk)^{-V-d} )  (p,q) 
\to
C^{*}(B G, \ssF_{\frac{{\bar{\cX}}|{\cZ}}{N}}({\bar{\cX}}|{\cZ},
\bbk)^{-V-d} )  (p,q) 
\end{multline}
Lemma~\ref{lem:SW-duality-closed-subset-monoidal} implies that the
maps of
Equation~\eqref{eq:homotopy-fixed-points-virtual-cochains-uncollared-comparison} 
assemble to a spectral functor
\begin{equation}
C^{*}(B G,
\ssF_{\frac{\hat{\bar{\cX}}|\hat{\cZ}}{N}}(\hat{\bar{\cX}}|\hat{\cZ},
\bbk)^{-V-d} ) \to C^{*}(B G,
\ssF_{\frac{{\bar{\cX}}|{\cZ}}{N}}({\bar{\cX}}|{\cZ}, \bbk)^{-V-d} ) .  
\end{equation}

On the other hand, applying Lemma \ref{SW-duality-lax-monoidal} yields:

\begin{lem}
There is a natural equivalence
\begin{equation}
C^{*}(B G,
\ssF_{\frac{\hat{\bar{\cX}}|\hat{\cZ}}{N}}(\hat{\bar{\cX}}|\hat{\cZ},
\bbk)^{-V-d} )  \to     C^{*,c_0}_{\rel \partial}(B \hat{\cZ};\Milnor
\hat{\bar{\cX}}^{-V  -d} \wedge \bbk)  
\end{equation}
of $\Pi$-equivariant spectral categories.
\end{lem}

\begin{proof}
The functor factors through the $\Pi$-equivariant spectral category where
the morphism spectra are defined as the homotopy colimits 
of the Milnor twisted cochains over $\Tw \square_{\Sub} A(p,q)$.  The
first map is then a level-wise equivalence, and a straightforward
application of Quillen's Theorem A implies that the second map is an
equivalence as well.
\end{proof}

\subsubsection{Using the Adams isomorphism and the norm map}
\label{sec:using-adams-isom}

We now explain how to apply the results of
Appendix~\ref{sec:adams-isom-ambid} to give an equivalence between the
category $C^{*}(B G,
\ssF_{\frac{{\bar{\cX}}|{\cZ}}{N}}({\bar{\cX}}|{\cZ}, \bbk)^{-V-d} )$ 
constructed above and the category $B\bar{\cX}|\cZ^{-V-d} \wedge \bbk$
of virtual cochains defined in Section~\ref{sec:virt-coch-homot} (see
Definition~\ref{def:virtual-cochains}).

We begin by casting virtual cochains in terms of the constructions of
Appendix \ref{sec:adams-isom-ambid}: let $ \Sp^{-\Vect}_{\eq}$ denote
the category introduced in Appendix~\ref{sec:funct-fibr-mapp} (see
Definition~\ref{defn:all-spectra-shifted}), with: 
\begin{enumerate}
\item Objects specified by triples $(G,
Y, V)$, with $Y$ a (cofibrant) $G$-spectrum and $V$ a
finite-dimensional $G$-representation.
\item Morphisms $f \colon (G,Y, V)
\to (G', Y', V')$ given by  
\begin{enumerate} 
\item a surjection $G \to G'$ with kernel $G_f^\perp$ acting freely on
$Y$, 
\item a $G$-embedding $V \to V'$ which is an isometric
embedding, with cokernel $V^\perp_f$, 
\item a $G$-equivariant map 
\begin{equation}
    S^{V^\perp_f} \sma Y \to Y'.
\end{equation}
\end{enumerate}
\end{enumerate}

We can interpret the virtual cochains in terms of the category
$\Sp^{-\Vect}_{\eq}$.

\begin{lem}
There is a lax monoidal functor
\begin{equation}\label{eq:kur-all-spectra-compare}
\begin{aligned}
  \square \Kur & \to \Tw \Sp^{- \Vect}_{\eq} \\
\sigma & \mapsto (G_\sigma, \bar{X}|Z_{\sigma}, V_\sigma).
\end{aligned}
\end{equation} \qed
\end{lem}

The functor $B\bar{X}_\sigma|Z^{-V_\sigma}_{\sigma}$ used to
construct the spectral category of virtual cochains factors through
the functor $B Y^{-V}$ considered in
Appendix~\ref{sec:funct-fibr-mapp}. 

To express the construction of the spectrum of sections in similar
terms, consider the category $\Sp^{- \Vect}_{\eq, \Sub}$ that refines
$\Sp^{-\Vect}_{\eq}$ (see
Definition~\ref{defn:all-spectra-shifted-subgroup}), with: 
\begin{enumerate}
\item Objects specified by tuples $(G, Y, V,N)$, consisting of a
  $G$-equivariant cofibrant spectrum $Y$, a finite-dimensional
  $G$-representation $V$, and a subgroup $N \subseteq G$.
\item Morphisms $f \colon (G_0, Y_0, V_0, N_0) \to (G_1, Y_1, V_1,
  N_1)$ specified by a morphism in $\Sp^{-\Vect}_{\eq}$ such that
  $N_1 \subseteq f(N_0)$.
\end{enumerate}

We have the following functor out of this category.

\begin{lem}
The assignment
\begin{equation}
  (\sigma,N)  \mapsto  (G_\sigma, {\bar{X}}_\sigma|{Z}_\sigma, V, N)
\end{equation}
defines a lax monoidal functor from $ \square_{\Sub} \Kur$ to $\Sp^{-
  \Vect}_{\eq, \Sub} $, which gives rise to a commutative diagram:
\begin{equation}
  \begin{tikzcd}
    \Tw \square_{\Sub} \Kur \ar[r] \ar[d] & \Tw \Sp^{- \Vect}_{\eq, \Sub} \ar[d] \\
    \square \Kur \ar[r] & \Sp^{- \Vect}_{\eq}.
  \end{tikzcd}
\end{equation} \qed
\end{lem}

The key remaining point is that given a Kuranishi flow category $\bX$,
the functor 
\begin{equation}
\begin{aligned}
\Tw \square_{\Sub}A(p,q) &\to \Sp \\
f &\mapsto C^{*}(B
G_f, \ssF_{\frac{{\bar{X}}_f|{Z}_f}{N'}}({\bar{X}}_f|{Z}_f,
\bbk)^{-V_f+  V_{q} - V_{p}} )[ \deg p  - \deg q ]
\end{aligned}
\end{equation}
applied in Equation~\eqref{eq:homotopy-fixed-points-virtual-cochains-uncollared})
factors through the functor 
\begin{equation}
\begin{aligned}
\Tw \Sp^{- \Vect}_{\eq, \Sub} &\to \Sp \\
(G, Y, V, N) &\mapsto C^{*}(B G_f, \ssF_{Y_f/{N'}}(Y_f, \bbk)^{-V_f} ) 
\end{aligned}
\end{equation}
considered in Appendix~\ref{sec:funct-fibr-mapp}.  At this point the
following proposition now completes the comparison.

\begin{prop}\label{prop:adams-comparison}
There is a $\Pi$-equivariant zig-zag of equivalences of spectral
categories 
\begin{equation}
B\bar{\cX}|\cZ^{-V-d} \wedge \bbk \htp C^{*}(B G,
\ssF_{\frac{{\bar{\cX}}|{\cZ}}{N}}({\bar{\cX}}|{\cZ}, \bbk)^{-V-d}).
\end{equation} 
\end{prop}

\begin{proof}
This comparison is deduced from a combination of several lax monoidal
equivalences of functors with domain $\Tw \Sp^{- \Vect}_{\eq, \Sub}$,
which are constructed in Appendix~\ref{sec:adams-isom-ambid} and
summarized in Section~\ref{sec:complete-comparison}.
\end{proof}

\subsubsection{Signpost: Spanier-Whitehead duality for Kuranishi flow categories}
\label{sec:signp-span-whit-1}

We have constructed a $\Pi$-equivariant zig-zag of equivalences of
spectral categories
\begin{equation}
C^{*,c}_{\rel \partial}(B \hat{\cZ}; \Milnor \hat{\bar{\cX}}^{-V-d})
\to \cdots \leftarrow B \hat{\bar{\cX}}|\hat{Z}^{-V -d},
\end{equation}
relating the virtual cochains to the Milnor-twisted cochains of any $\Pi$-equivariant Kuranishi flow category. In the next section, we explain how, by lifting to a complex oriented Kuranishi flow category, and studying the cochains with coefficients in Morava $K$-theory, we obtain an isomorphism between the twisted and untwisted cochains.

\section{Coherent comparisons: tangent spaces and orientations}
\label{sec:coher-comp-tang}

The purpose of this section is to implement, in a functorial context,
the comparison of models for the tangent space considered in Section
\ref{sec:smoothly-fiber-kuran-1} and the construction of orientations
from Section \ref{sec:compl-lifts-tang}.  This requires us to work
with Kuranishi charts and presentations equipped with tangent
bundles and complex oriented Kuranishi charts and presentations, 
as introduced in Section \ref{sec:compl-orient-kuran-1}.  Throughout
the section, we will work with the internal category $\Tw \Box
\Kur^{\scrT}$; recall from Appendix~\ref{sec:internal} that this
category inherits its internal structure from $\Kur^{\scrT}$.  We will
also work with twisted arrow categories on cubes on enriched
categories, which inherit enrichments in the usual way.

\subsection{Tangent spaces and Whitney spherical fibrations}
\label{sec:tang-spac-whitn}

\subsubsection{Whitney spherical fibration of cubical diagrams}
\label{sec:whitn-spher-fibr}

Let $\sigma$ and $\tau$ be cubes in $\Kur^{fs}\langle \cS \rangle$,
and $f \co \tau \to \sigma$ a map of cubes.  We will write
$\square_\sigma$ for the geometric realization of the domain of a cube
$\sigma$, and $\square_f$ for the image of $\square_\tau$ in
$\square_\sigma$ under a map $f \colon \tau \to \sigma$ of cubes. We
also write $\kappa_f$ for the collar of $\square_\tau$ in
$\hat{\square}_\sigma$; this is a cube of dimension $\dim \sigma -
\dim \tau$.

Forgetting down from $\Box\Kur^{fs}\langle \cS \rangle$ to $\Box\Kur
\langle \cS \rangle$, recall from Section~\ref{sec:cubic-degen-norm} 
that we can associate to $f \co \tau \to \sigma$ a cubical
degeneration to the normal cone $X_f$ and a map 
$Z_f \to X_f$.  We have a natural map $X_\sigma \to \square_\sigma
\times X_{\sigma(1^n)}$, and we let $B_\sigma$ denote
$B_{\sigma(1^n)}$: the fibre of $X_\sigma \to B_\sigma$ has a natural
fibrewise smooth structure, so we define
\begin{equation}
T X_\sigma \equiv T^\sigma X_\sigma \oplus TB_\sigma.
\end{equation}

We now explain the identification that results from lifting $f$ to the
category $\Kur^{\scrT} \langle \cS \rangle$ of Kuranishi charts with
tangent bundles.  We begin by studying the case of a cube $\sigma$.

\begin{lem}\label{lem:tangent-cube}
Let $\sigma$ be an object of $\Box \Kur^{fs}\langle \cS \rangle$.  A
lift of $\sigma$ to $\Box \Kur^{\scrT} \langle \cS \rangle$ determines
an equivariant isomorphism
\begin{equation}
T X_\sigma \cong T X_{\sigma(1^n)} \oplus T \square_\sigma.
\end{equation}
Given a map $f \co \tau \to \sigma$ in $\square \Kur^{\scrT} \langle
\cS \rangle$, there is a natural equivariant isomorphism along
$X_\tau$  
\begin{equation}
f^* T X_\sigma \cong TX_\tau \oplus T \kappa_f \oplus V_f^\perp,
\end{equation}
so that the following diagram commutes:
\begin{equation}\label{eq:boundary-compatibility}
\begin{tikzcd}
f^* T X_\sigma \ar[d] \ar[r] & f^*\left(T^\sigma X_{\sigma(1^n)} \oplus T
\square_\sigma\right) \ar[d] \\ 
T X_\tau \oplus T \kappa_f \oplus V_f^\perp \ar[r] & T X_{\tau(1^m)}
\oplus T \square_\tau \oplus T \kappa_f \oplus V_f^\perp. 
\end{tikzcd}
\end{equation}
Given maps $f \co \tau \to \sigma$ and $g \co \nu \to \tau$, the
natural equivariant isomorphism
\begin{equation}
(g \circ f)^* T X_\sigma \cong TX_\nu \oplus T \kappa_{g \circ f}
  \oplus V_{g \circ f}^\perp
\end{equation}
coincides with the composite isomorphism
\begin{equation}\label{eq:composition-isomorphism}
\begin{aligned}
g^* (f^* T X_\sigma) &\cong g^* (TX_\tau \oplus T \kappa_f \oplus
V_f^\perp) \\
&\cong (TX_\nu \oplus T\kappa_{g} \oplus V_{g}^\perp) \oplus g^* T
\kappa_f \oplus g^* V_f^{\perp} \\
&\cong TX_\nu \oplus T \kappa_{g \circ f}
  \oplus V_{g \circ f}^\perp.
\end{aligned}
\end{equation}
and the analogue of Equation~\eqref{eq:boundary-compatibility}
commutes.
\end{lem}

\begin{proof}
By definition, the fibrewise tangent space sits in a short exact
sequence of vector bundles 
\begin{equation}
 T^\sigma X_\sigma \to  T \square_\sigma \oplus T^\sigma X_{\sigma(1^n)}  \oplus V_\sigma \to V_\sigma;
\end{equation}
the key point is that the map in the right depends on the chosen point in $\square_\sigma$. Taking the quotient by $V_\tau$ in the middle and right terms yields the sequence
\begin{equation}
 T^\sigma X_\sigma \to  T \square_\sigma \oplus T^\sigma X_{\sigma(1^n)}  \oplus V^\perp_f  \to V^\perp_f.
\end{equation}
The inner product on $T^\sigma X_{\sigma(1^n)}$ decomposes the above
as the direct sum of the kernel of the projection from $T^\sigma
X_{\sigma(1^n)} $ to $V^\perp_f $ with the kernel of 
\begin{equation}
  T \square_\sigma \oplus V^\perp_f  \oplus V^\perp_f  \to V^\perp_f.
\end{equation}
At this stage, we observe that this kernel can be described as the
graph of a map $T \square_\sigma \oplus V^\perp_f  \to V^\perp_f$,
whose target is the third factor in the source of the above map. This
yields the desired trivialisation. The compatibility with boundary
strata then follows from the compatibility of morphisms in
$\Kur^{\scrT}\langle \cS \rangle$ with strata.  The identification of
the isomorphism in Equation~\eqref{eq:composition-isomorphism} with
the isomorphism induced by the composite follows from the fact that
the decompositions involved are constructed in terms of the inner
products and directions in the cubes; despite the ordering suggested
by the notation, the terms of the decompositions coincide.
\end{proof}

We now extend Lemma~\ref{lem:tangent-cube} to maps of cubes as follows.
We define $TX_f$ to be the vector bundle over $X_f$ obtained from
$TX_\tau$ by taking the quotient by $G_f^\perp$.  We have the
following straightforward generalization:

\begin{cor} \label{cor:comparison-tangent-space-maps-of-cubes}
Each arrow $f\colon \sigma \to \tau$ in $\Box \Kur^{\scrT} \langle \cS
\rangle$ gives rise to a natural isomorphism 
  \begin{equation}
     T X_f \cong T X_{f(1^m)} \oplus T \square_f.
   \end{equation}
A factorization $f_0 = h \circ f_1 \circ g$ in $\Box \Kur^{\scrT}
\langle \cS \rangle$ determines a natural isomorphism 
\begin{equation}
T X_{f_0} \cong T X_{f_1} \oplus T \kappa_g \oplus V_g^\perp.
\end{equation} 
Given another factorization $f_1 = h' \circ f_2 \circ g'$ in $\Box
\Kur^{\scrT}$, the isomorphism associated to the factorization $f_2 =
(h' \circ h) \circ f_0 \circ (g \circ g')$ 
\begin{equation}
T X_{f_0} \cong T X_{f_2} \oplus T \kappa_{g \circ g'} \oplus V_{g
  \circ g'}^\perp.
\end{equation}
is equal to the composite isomorphism
\begin{equation}
\begin{aligned}
T X_{f_0} &\cong T X_{f_1} \oplus T \kappa_{g'} \oplus V_{g'}^\perp \\
          &\cong \left(T X_{f_2} \oplus T \kappa_{g} \oplus
V_{g}^{\perp}\right) \oplus T \kappa_{g'} \oplus V_{g'}^\perp \\
          &\cong TX_{f_2} \oplus (T \kappa_g \oplus T \kappa_{g'})
\oplus (V_g^{\perp} \oplus T_{g'}^{\perp})
\end{aligned}
\end{equation}
\qed
\end{cor}

\begin{rem}
Since the action of  $G_f^\perp$ on $B_\tau$ may not be free, directly
making sense of $T X_f$ requires introducing the tangent space of the
orbifold $B_f \equiv B_\tau/G^\perp_f$. We return to this point at the
beginning of Section~\ref{sec:nash-spher-fibr-1}.
\end{rem}

We next shall pass to collared cubical degenerations. In order to
study their tangent space, we need to equip the corresponding collared
base with a smooth structure along the region where the collars are
attached. To this end, we pick a compatible choice of smooth collar on
$B_f$ (in the sense of Definition \ref{def:smooth_collars_compatible})
for each arrow $f$ in $\Box \Kur^{\scrT} \langle \cS \rangle$.  
By induction on the dimension of the cube, we can arrange that these
collars are preserved by morphisms in $\Tw \Box \Kur^{\scrT} \langle
\cS \rangle $.

Having made this choice, we associate to each such arrow the vector bundle
\begin{equation}
  T\hat{X}_f \cong T^f \hat{X}_f \oplus T \hat{B}_f.
\end{equation}
over $\hat{X}_f$.

We will write $S^{T\hat{X}_f}$ for the spherical fibration over
$\hat{X}_f$ (and by restriction over $\hat{Z}_f$) obtained from the
fibrewise $1$-point compactification of $T \hat{X}_f$. Since the study
of vector bundles and their associated sphere bundles goes back to
Whitney \cite{Whitney1935}, we shall use the following terminology: 

\begin{defin}
Given an object $f \in \Tw \Box \Kur^{\scrT} \langle \cS \rangle$, the
\emph{Whitney spherical fibration} $S^{T\hat{X}_f}|0$ on $\hat{X}_f$
is the fibrewise cone of the complement of the zero section in
$S^{T\hat{X}_f}$. 
\end{defin}

Although $\Box \Kur^{\scrT} \langle \cS \rangle$ has a topologized space
of objects, the Whitney spherical fibration only depends on the
component of $f$ in the mapping space of $\Box \Kur^{\scrT} \langle
\cS \rangle$.

\begin{lem}\label{lem:Whitney-locally-constant}
Suppose that $f_0$ and $f_1$ belong to the same component of the space
of maps in $\Box \Kur^{\scrT} \langle \cS \rangle$.  Then the
spherical fibrations $S^{T\hat{X}_{f_0}}|0$ and $S^{T\hat{X}_{f_1}}|0$
coincide.
\end{lem}

\begin{proof}
Clearly $\hat{X}_{f_0}$ and $\hat{X}_{f_1}$ are equal, as the
quotient in question does not depend on the choice of inner product,
and similarly $T\hat{X}_{f_0}$ and $T\hat{X}_{f_1}$ coincide.  Since
the zero section is also independent of the inner product, the result
follows.
\end{proof}

Given a spectrum $\bbk$, we can stabilize the Whitney spherical
fibration.  See Section~\ref{sec:parametrized-spectra} for a concise
review of the aspects of the theory of parametrized spectra we require.

\begin{defin}
For each arrow $f\colon \sigma \to \tau$ in $\Box \Kur^{\scrT} \langle
\cS \rangle$, we have parametrized Whitney spectra 
\begin{equation}
S^{T\hat{X}_f - V} \sma \bbk \quad\textrm{and}\quad S^{T \hat{X}_f}|0^{-V} \wedge
\bbk
\end{equation}
over $\hat{X}_f$.  Pulling back along the map $\hat{Z}_f \to
\hat{X}_f$, we obtain parametrized spectra over $\hat{Z}_f$.
\end{defin}

By Lemma~\ref{lem:Whitney-locally-constant}, these spectra only depend
on the component of $f$ in the space of maps in $\Box \Kur^{\scrT}
\langle \cS \rangle$.

We consider as before the spectra of Borel equivariant compactly
supported sections,  which vanish at the
boundary,  of these parametrized spectra.  

\begin{defin}\label{defn:sections-Whitney-stable}
Let $f \colon \sigma \to \tau$ be a map in $\Box \Kur^{\scrT}\langle
\cS \rangle$.  Then there are assignments
\begin{equation}
\begin{aligned}
f &\mapsto C^{*,c}_{\rel \partial} \left( B\hat{Z}_f; S^{T \hat{X}_f
  -V} \wedge \bbk\right) \qquad \textrm{and} \\ f &\mapsto
C^{*,c}_{\rel \partial} \left( B\hat{Z}_f; S^{T \hat{X}_f}|0^{-V} \wedge
\bbk\right).
\end{aligned}
\end{equation}
Introducing the notation $\dhatS{X}_{f \circ g}$ for the double
collar (which we will require for technical reasons below), we also
have an assignment of global sections
\begin{equation}
f \mapsto C^{*,c}_{\rel \partial} \left( B\dhatS{Z}_{f};
S^{T\dhatS{X}_{f}}|0^{-V} \wedge \bbk\right).
\end{equation}
\end{defin}

We now explain the functoriality and compatibility of these
constructions.  We start by explaining how to assemble them into
spectral categories.  As usual, we begin with the pointwise
functoriality necessary to build the mapping spectra.  In the
following result, we use the twisted arrow category and categories of
cubes constructions in the context of internal categories; see
Definition~\ref{defn:internal-twisted-arrow} and the surrounding
discussion for a review of these constructions.  Throughout, we will
tacitly regard $\Sp$ as an internal category in topological spaces
with the discrete topology on objects.  

\begin{prop}
The  three assignments
of Definition~\ref{defn:sections-Whitney-stable} assemble to 
topologically enriched functors 
\begin{equation}
\Tw  \square \Kur^{\scrT}\langle \cS \rangle \to \Sp.
\end{equation}
\end{prop}

\begin{proof}
We give the argument for the functor specified by $C^{*,c}_{\rel
  \partial} \left(  B\hat{Z}_{{h}}; S^{T \hat{X}_{{h}}}|0^{-V} \right)$;
the case of the other two spectra is entirely analogous.  Suppose that
we have a factorization $f = g \circ f' \circ h$.  Then 
Corollary~\ref{cor:comparison-tangent-space-maps-of-cubes} shows that 
the diagram 
\begin{equation}\label{eq:maps_cochains_Whitney_compositions}
\begin{tikzpicture}
\node[scale=.7]{
\begin{tikzcd}
 C^{*,c}_{\rel \partial} \left(  B\hat{Z}_{{h}}; S^{T
   \hat{X}_{{h}}}|0^{-V} \right)     &  C^{*,c}_{\rel \partial} \left(
 B\hat{Z}_{{f}' \circ h}; S^{T \hat{X}_{{f}' \circ h}}|0^{-V}  \right)
 \ar["{f}' \circ",l]  \ar["\circ h",d] & C^{*,c}_{\rel \partial}
 \left(  B\hat{Z}_{f}; S^{T \hat{X}_{f}}|0^{-V} \right)   \ar["g
   \circ",l] \ar["\circ h",d]  \ar[ "\circ {f}' \circ h", bend
   left=80,dd] \ar["g \circ {f}' \circ", bend right=10,ll,swap]  \\ 
 & C^{*,c}_{\rel \partial} \left(  B\hat{Z}_{{f}'}; S^{T
   \hat{X}_{{f}'}}|0^{-V}  \right)     & C^{*,c}_{\rel \partial}
 \left(  B\hat{Z}_{g \circ {f}'}; S^{T \hat{X}_{g \circ {f}'}}|0^{-V}
 \right) \ar[ "\circ {f}'", d]  \ar["g \circ",l] \\ 
 & & C^{*,c}_{\rel \partial} \left(  B\hat{Z}_{g}; S^{T \hat{X}_g}|0^{-V} \right) , 
\end{tikzcd}
};
\end{tikzpicture}
\end{equation} 
with horizontal arrows given by pullback along compositions and
vertical arrows by pushforward along collars, commutes.  
(Similarly, the analogous versions of
Equation~\eqref{eq:maps_cochains_Whitney_compositions} 
 for $C^{*,c}_{\rel \partial} \left(
B\hat{Z}_f; S^{T \hat{X}_f -V} \wedge \bbk\right)$ 
and $C^{*,c}_{\rel \partial} \left( B\dhatS{Z}_{f};
S^{T\dhatS{X}_f}|0^{-V} \wedge \bbk\right)$
also commute.)

The functor is then given on the morphism in $\Tw \square
\Kur^{\scrT}\langle \cS \rangle \to \Sp$ from $f$ to $f'$ represented
by $f = g \circ f' \circ h$ by the composition 
\begin{equation}
C^{*,c}_{\rel \partial}  \left(  B\hat{Z}_{f}; S^{T
  \hat{X}_{f}}|0^{-V} \right) \to 
C^{*,c}_{\rel \partial} \left(  B\hat{Z}_{{f}'}; S^{T
   \hat{X}_{{f}'}}|0^{-V}  \right)
\end{equation}
around the square in
Equation~\eqref{eq:maps_cochains_Whitney_compositions}.  To see that
this is a functor, assume we have the composition of factorizations $f
= g \circ f' \circ h$ and $f' = g' \circ f'' \circ h'$.  Then the
compatibility with composition comes from the fact that the diagram 
\begin{equation}
\begin{tikzpicture}
\node[scale=.7]{
\begin{tikzcd}
 &  C^{*,c}_{\rel \partial} \left(
 B\hat{Z}_{{f}' \circ h}; S^{T \hat{X}_{{f}' \circ h}}|0^{-V}  \right)
 \ar["\circ h",d] & C^{*,c}_{\rel \partial}
 \left(  B\hat{Z}_{f}; S^{T \hat{X}_{f}}|0^{-V} \right)  \ar["g
   \circ",l] \ar["\circ h",d]  \\ 
\ar[d,"\circ h'"] C^{*,c}_{\rel \partial} \left(
 B\hat{Z}_{f'' \circ h'}; S^{T \hat{X}_{f'' \circ h'}}|0^{-V}  \right)
 & C^{*,c}_{\rel \partial} \left(  B\hat{Z}_{{f}'}; S^{T
   \hat{X}_{{f}'}}|0^{-V}  \right) \ar[l,"g' \circ"] \ar[d,"\circ h'"] & C^{*,c}_{\rel \partial}
 \left(  B\hat{Z}_{g \circ {f}'}; S^{T \hat{X}_{g \circ {f}'}}|0^{-V}
 \right) \ar["g \circ",l] \\ 
C^{*,c}_{\rel \partial} \left(
 B\hat{Z}_{f''}; S^{T \hat{X}_{f''}}|0^{-V}  \right) & C^{*,c}_{\rel \partial} \left(
 B\hat{Z}_{g' \circ f''}; S^{T \hat{X}_{g' \circ f''}}|0^{-V}  \right)
 \ar[l,"g' \circ"] & .
\end{tikzcd}
};
\end{tikzpicture}
\end{equation} 
and the diagram from the composite factorization $f = g \circ g' \circ
f'' \circ h \circ h'$
\begin{equation}\label{eq:maps_cochains_Whitney_compositions_functoriality}
\begin{tikzpicture}
\node[scale=.7]{
\begin{tikzcd}
C^{*,c}_{\rel \partial} \left(
 B\hat{Z}_{f'' \circ h \circ h'}; S^{T \hat{X}_{f'' \circ h \circ
     h'}}|0^{-V}  \right) 
\ar["\circ h \circ h'",d] & C^{*,c}_{\rel \partial}
 \left(  B\hat{Z}_{f}; S^{T \hat{X}_{f}}|0^{-V} \right)   \ar["g \circ
   g' \circ",l] \ar["\circ h \circ h'",d] \\
C^{*,c}_{\rel \partial} \left(  B\hat{Z}_{f''}; S^{T
   \hat{X}_{f''}}|0^{-V}  \right)     & C^{*,c}_{\rel \partial}
 \left(  B\hat{Z}_{g \circ {f}'}; S^{T \hat{X}_{g \circ {f}'}}|0^{-V}
 \right) \ar["g \circ g' \circ",l] 
\end{tikzcd}
};
\end{tikzpicture}
\end{equation} 
have the same outer composite.  This follows from the
composition statement in 
Corollary~\ref{cor:comparison-tangent-space-maps-of-cubes}.

To see that these functors are topologically enriched, we need to
check that functors are locally constant and induce continuous maps of
mapping spaces.  Lemma~\ref{lem:Whitney-locally-constant} implies the
first part of this, and the second condition follows from the
definition of morphisms in terms of the pushforward along collars and
the pullback maps.
\end{proof}

In light of this, we can make the following definition.  The homotopy
colimits here are topologized, as we are working with internal
categories.  However, since the range category has the discrete
topology on objects, this construction simplifies; see
Definition~\ref{defn:internal_hocolim} for a review of the
definition.

\begin{defin}\label{defin:whitney-mapping-spectra}
Given a Kuranishi presentation with tangent bundle $\bX \co A^{\scrI} \to
\Kur^{\scrT}\langle \cS \rangle$, we construct the following homotopy
colimits of spectra  
\begin{equation}
\begin{aligned}
 C^{*,c}_{\rel \partial} \left(  B\hat{Z}; S^{T \hat{X}}|0^{-V} \right) (A^{\scrI}) & \equiv \hocolim_{f \in \Tw \square A^{\scrI}}  C^{*,c}_{\rel \partial} \left( B\hat{Z}_f; S^{T\hat{X}_f}|0^{-V} \wedge \bbk\right)  \\
 C^{*,c}_{\rel \partial} \left(  B\hat{Z}; S^{T \hat{X} -V} \right) (A^{\scrI}) & \equiv \hocolim_{f \in \Tw \square A^{\scrI}}  C^{*,c}_{\rel \partial} \left( B\hat{Z}_f; S^{T \hat{X}_f -V} \wedge \bbk\right)  \\
 C^{*,c}_{\rel \partial} \left(  B\dhatS{Z}; S^{T \dhatS{X}}|0^{-V} \right) (A^{\scrI}) & \equiv \hocolim_{f \in \Tw \square A^{\scrI}}  C^{*,c}_{\rel \partial} \left( B\dhatS{Z}_f; S^{T\dhatS{X}_f}|0^{-V} \wedge \bbk\right) ,
\end{aligned}
\end{equation}
where there are two different models of the first spectrum,
corresponding to the choice of collar section which is hidden from the
notation. 
\end{defin}

We now assemble the pointwise constructions of
Definition~\ref{defin:whitney-mapping-spectra} to produce
$\Pi$-equivariant spectral categories; this amounts to showing that
there are natural associative composition maps which are
$\Pi$-equivariant.  We can choose the smooth collars to be
strictly compatible with products.  Then the multiplicative diagrams
established in 
Lemma~\ref{lem:multiplicativity-comparison-Whit-tangent} provide
associative composition maps.  By construction, these mapping spectra are
strictly compatible with the action of $\Pi$ and the composition of
the lemma clearly are as well, so we conclude the following.

\begin{prop}
Given a Kuranishi flow category with tangent bundle $\bX$, each of the
four models of twisted cochains in
Definition~\ref{defin:whitney-mapping-spectra} gives rise to a
$\Pi$-equivariant spectral category with objects those of
$\bX$ and morphism spectra 
\begin{equation}
\begin{aligned}
C^{*,c}_{\rel \partial} \left(  B\hat{\cZ}; S^{T \hat{\cX}}|0^{-V-d} \right)(p,q) & \equiv  C^{*,c}_{\rel \partial} \left(  B\hat{Z}; S^{T \hat{X}}|0^{-V + V_{q} - V_{p}} \right)(p,q)[ \deg p  - \deg q ], \\
C^{*,c}_{\rel \partial} \left(  B\dhatS{\cZ}; S^{T
  \dhatS{\cX}}|0^{-V-d} \right)(p,q) & \equiv  C^{*,c}_{\rel \partial}
\left(  B\dhatS{Z}; S^{T \dhatS{X}}|0^{-V + V_{q} - V_{p}}
\right)(p,q)[ \deg p  - \deg q ] \\ 
C^{*,c}_{\rel \partial} \left(  B\hat{\cZ};  S^{T \hat{\cX} -V -d}
\right)(p,q) & \equiv  C^{*,c}_{\rel \partial} \left(  B\hat{Z}; S^{T
  \hat{X} -V + V_{q} - V_{p}} \right)(p,q)  [ \deg p  - \deg q ]. 
\end{aligned}
\end{equation}
\qed
\end{prop}

Next, we will compare these $\Pi$-equivariant spectral categories.
The inclusion induces a pointwise comparison:

\begin{prop}
For each map $f\colon \sigma \to \tau$ in $\Box \Kur^{\scrT} \langle \cS
\rangle$, there is a natural comparison map  
\begin{equation} \label{eq:Whitney-cochains-definition}
C^{*,c}_{\rel \partial} \left( B\hat{Z}_f; S^{T \hat{X}_f -V}
\wedge \bbk\right) \to C^{*,c}_{\rel \partial} \left( B\hat{Z}_f;
S^{T \hat{X}_f}|0^{-V} \wedge \bbk\right) 
\end{equation}
of spectra of compactly supported sections which vanish at the
boundary. \qed
\end{prop}

The analogue of the argument for Lemma~\ref{lem:map_pullback_Milnor}
now establishes the following compatibility result.

\begin{lem}
For composable maps $f$ and $g$ in $\Box \Kur^{\scrT} \langle \cS
\rangle$, there is a commutative diagram 
\begin{equation}
  \begin{tikzcd}
    C^{*,c}_{\rel \partial} \left( B\hat{Z}_{f \circ g};   S^{T
      \hat{X}_{f \circ g} -V} \wedge \bbk\right) \ar[r] \ar[d] &
    C^{*,c}_{\rel \partial} \left( B\hat{Z}_{f \circ g}; S^{T
      \hat{X}_{f \circ g}}|0^{-V} \wedge \bbk\right) \ar[d] \\ 
C^{*,c}_{\rel \partial} \left( B\hat{Z}_g;   S^{T \hat{X}_g -V} \wedge
\bbk\right) \ar[r] & C^{*,c}_{\rel \partial} \left( B\hat{Z}_g;   S^{T
  \hat{X}_g}|0^{-V} \wedge \bbk\right) 
  \end{tikzcd}
\end{equation}
where the vertical maps are induced by pullback. \qed
\end{lem}

The comparison between the spectra associated to $f \circ g$ and $f$
is not as straightforward. There are, as discussed in
Section~\ref{sec:tang-twist-coch}, two natural sections that one can
use along the collar, which are respectively adapted to the comparison
with the Milnor model and to the standard spherical fibration.  Using
the double collar $\dhatS{X}_{f \circ g}$, we can use the
concatenation of these two sections to provide the desired comparison,
as in Lemma~\ref{lem:multiplicativity-comparison-Whit-tangent}.

\begin{prop}\label{prop:comparison-3-models}
For composable maps $f$ and $g$ in $\Box \Kur^{\scrT} \langle \cS
\rangle$, there is a commutative diagram
\begin{equation} \label{eq:comparison-3-models-Whit-pushforward-arrow}
\begin{tikzcd}
C^{*,c}_{\rel \partial} \left( B\hat{Z}_{f \circ g};   S^{T\hat{X}_{f \circ g}}|0^{-V} \wedge \bbk\right) \ar[r] & C^{*,c}_{\rel \partial} \left( B\hat{Z}_f;   S^{T\hat{X}_f}|0^{-V} \wedge \bbk\right)  \\ 
    C^{*,c}_{\rel \partial} \left( B\dhatS{Z}_{f \circ g};  S^{T
      \dhatS{X}_{f \circ g}}|0^{-V}  \wedge \bbk\right) \ar[r] \ar[d]
    \ar[u]  & C^{*,c}_{\rel \partial} \left( B\dhatS{Z}_{f};
    S^{T\dhatS{X}_{f}}|0^{-V} \wedge \bbk\right) \ar[d] \ar[u]  \\ 
C^{*,c}_{\rel \partial} \left( B\hat{Z}_{f \circ g}; S^{T \hat{X}_{f
    \circ g}}|0^{-V} \wedge \bbk\right) \ar[r] &   C^{*,c}_{\rel
  \partial} \left( B\hat{Z}_f;   S^{T \hat{X}_f}|0^{-V} \wedge
\bbk\right)  \\  
C^{*,c}_{\rel \partial} \left( B\hat{Z}_{f \circ g};  S^{T \hat{X}_{f
    \circ g} -V} \wedge \bbk\right) \ar[r] \ar[u] & C^{*,c}_{\rel
  \partial} \left( B\hat{Z}_f; S^{T \hat{X}_g -V} \wedge \bbk\right)
\ar[u] . 
\end{tikzcd}  
\end{equation} \qed
\end{prop}

To compare the spectral categories, we need to show that the
comparison maps above pass to homotopy colimits.

\begin{lem}\label{lem:pointwise-comparison-Whit-STX}
Given a Kuranishi presentation with tangent bundle $\bX \co A^{\scrI} \to
\Kur^{\scrT}\langle \cS \rangle$, there are natural equivalences of spectra
\begin{equation} \label{eq:Kuranishi-presentation-comparison-Whit-STX-pointwise}
\begin{tikzcd}
C^{*,c}_{\rel \partial} \left(  B\hat{Z}; S^{T \hat{X}-V}|0^{-V}
\right) (A^{\scrI}) & \ar[l] \ar[d]      C^{*,c}_{\rel \partial}
\left(  B\hat{Z}; S^{T \dhatS{X}}|0^{-V} \right) (A^{\scrI}) \\ 
C^{*,c}_{\rel \partial} \left(  B\hat{Z}; S^{T \hat{X} -V} \right)
(A^{\scrI}) \ar[r] &   C^{*,c}_{\rel \partial} \left(  B\hat{Z}; S^{T
  \hat{X}}|0^{-V} \right) (A^{\scrI})  
\end{tikzcd}
\end{equation}
\end{lem}

\begin{proof}
To see that the comparisons of
Proposition~\ref{prop:comparison-3-models} are functorial, we consider
composable factorizations $f = g \circ f' \circ h$ and $f' = g' \circ
f'' \circ h'$ in $\Box \Kur^{\scrT} \langle \cS \rangle$.  Putting these
together yields the composition $f = g \circ g' \circ {f}'' \circ h'
\circ h$ and inspection of the induced diagram of cochains associated
to the commutative diagram
\begin{equation}
\begin{tikzcd}
X_{{f}'' \circ h' \circ h}  \ar[twoheadrightarrow,"g'  \circ ",r]
\ar[twoheadrightarrow,"g \circ g'  \circ ",bend left=25,rr]
\ar[hookrightarrow, "\circ h",d] \ar[hookrightarrow, "\circ h' \circ
  h"',bend right=50,dd] &  X_{{f}' \circ h} \ar[twoheadrightarrow,"g
  \circ",r] \ar[hookrightarrow,"\circ h",d]  &   X_f  \\  
X_{{f}'' \circ h'}  \ar[twoheadrightarrow,"g'  \circ ",r]
\ar[hookrightarrow, "\circ h'",d]  & X_{{f}'}   &   \\ 
X_{{f}''}. & &
\end{tikzcd}
\end{equation}
yields the result; the key thing to check is that the concatenation of
sections that give rise to Proposition~\ref{prop:comparison-3-models}
is compatible with this diagram.  This follows because the comparison
is ultimately induced by collapsing along collar directions and
restricting sections (see Equation~\eqref{eq:collar-collapse-Milnor}).
Specifically, since the collapsing occurs along specified directions
for each map, to see that it is compatible with the composition, it
suffices to observe that by construction the labeling of these
directions is compatible with the composition. 
\end{proof}

The evident extension of the construction of
Section~\ref{sec:tang-twist-coch} and
Lemma~\ref{lem:multiplicativity-comparison-Whit-tangent} 
now shows that these comparison maps are multiplicative.  Moreover, it
is straightforward to check that the multiplicative structure is
compatible with the functoriality in the twisted arrow category since
the collapsing of collars is evidently multiplicative.  Putting this
all together, we have a comparison of spectral categories.

\begin{lem}
There are $\Pi$-equivariant DK-equivalences of spectral categories
\begin{equation} \label{eq:Kuranishi-presentation-comparison-Whit-STX}
\begin{tikzcd}
C^{*,c}_{\rel \partial} \left(  B\hat{\cZ}; S^{T \hat{\cX}}|0^{-V-d}
\right)    & \ar[l] \ar[d]    C^{*,c}_{\rel \partial} \left(
B\dhatS{\cZ}; S^{T \dhatS{\cX}}|0^{-V-d} \right)    \\ 
C^{*,c}_{\rel \partial} \left(  B\hat{\cZ};  S^{T \hat{\cX} -V -d}
\right) \ar[r] &  C^{*,c}_{\rel \partial} \left(  B\hat{\cZ}; S^{T
  \hat{\cX}}|0^{-V-d} \right). 
\end{tikzcd}
\end{equation} \qed
\end{lem}

\subsection{The Nash spherical fibration}
\label{sec:nash-spher-fibr-1}

Recall from Section~\ref{sec:nash-spher-fibr} (notably
Definition~\ref{defin:nash-fibration}) that we defined the Nash
spherical fibration of a smoothly fibred Kuranishi chart $\alpha$ in
terms of a subspace of the space $N X_\alpha$ of maps from a quadrant
to $X_\alpha$, whose projection to $B_\alpha$ depends only on the
first coordinate, which are appropriately differentiable near the
origin, and whose restriction to the diagonal gives a Nash path. In
this section, we formulate the functoriality of this construction. The
key idea that allows us to achieve suitable functoriality is to define
the Nash fibration on collared cubes in a piecewise way.

\subsubsection{The Nash tangent space}
\label{sec:spaces-nash-cubes}
Given a map $f \co \tau \to \sigma$ in $\Box \Kur^{fs}\langle  \cS
\rangle$, consider the projection map 
\begin{equation} \label{eq:map_X_f-to-quotient-B_f}
  X_f \to B_\tau/G^\perp_f \equiv B_f.
\end{equation}
Since the action of $G_f^\perp$ on $B_\tau$ may not be free, the space
$B_\tau/G^\perp_f$ is not necessarily a manifold.  However, it is an
orbifold, and so we can still make sense of the tangent space of $B_f$
(as an orbibundle), and of the notion of derivatives of paths.

\begin{rem}
Since we are only interested in the corresponding notions for paths
which factor through $X_f$, we shall lay things out explicitly: the
tangent space $TB_f$ makes sense as a vector bundle over $X_f$, since
it can be defined as the quotient of the pullback of $TB_\tau$ as a
vector bundle on $X_\tau$ under the action of $G^\perp_f$, as we did
above in
Corollary~\ref{cor:comparison-tangent-space-maps-of-cubes}. In
particular, we say that a path in $X_f$ projects to a smooth path in
$B_f$ if it lifts to a path in $X_\tau$ whose projection to $B_\tau$
is smooth. In such a situation, the derivative of the path is an
element of $TB_f$.
\end{rem}

The Nash tangent space of $\hat{X}_f$, denoted $N \hat{X}_f$, consists
as in Section~\ref{sec:nash-spher-fibr} of maps 
\begin{equation}
 [0,\infty)^{2} \to \widetilde{X}_f,
\end{equation}
where $\widetilde{X}_f$ is the space obtained by attaching infinite
collars to $X_f$, satisfying the following conditions:
\begin{enumerate}
\item  the restriction to the diagonal is a Nash path,
  \item the map takes the origin to $\hat{X}_f$,
\item the projection to the
infinite completion $\tilde{B}_{f}$ of the base is a smooth path that
depends only on the first variable, and 
\item the family of paths in the fibres of the projection map to
  $\tilde{B}_f$, parametrized by the first coordinate, is continuously
  differentiable.
\end{enumerate}

This construction is unfortunately not sufficiently functorial; the
inclusion $\hat{X}_f \to \hat{X}_\sigma$ does not induce a map of Nash
tangent spaces because we do not have a splitting of the projection
$B_f \to B_\sigma$.  Our strategy to handle functoriality will be to
construct a Nash tangent space over $\dhatS{X}_f$ in a piecewise way. 

\begin{figure}[h]
  \centering
\begin{tikzpicture}
  \draw (-3,-3) -- (-3,3) -- (3,3) -- (3,-3) -- cycle;
  \draw (-1.5,-3) -- (-1.5,3);
  \draw (1.5,-3) -- (1.5,3);
  \draw (-3,-1.5) -- (3,-1.5);
  \draw (-3,1.5) -- (3,1.5);
  \draw[] (-2.25,-2.25) node[] {$\dhatC{\square}^{g \circ h}_{f}$};
  \draw[] (0,-2.25) node[] {$\dhatC{\square}^{g}_{f}$};
   \draw[] (0,0) node[]  {$ \dhatC{\square}_{f}^{\id} = \hat{\square}_{f}$};
\end{tikzpicture}  
\caption{Decomposition of square}
\label{fig:decomposition-square}
\end{figure}
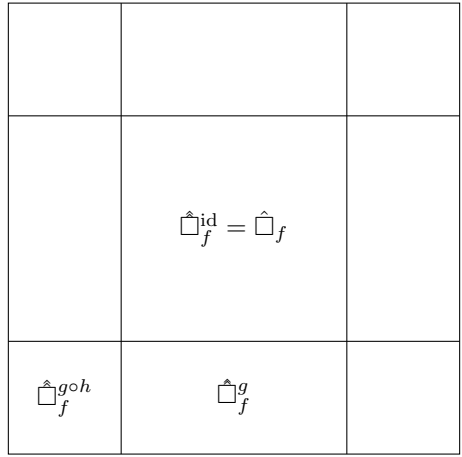

Recall that the strata of $X_f$ are indexed by compositions $f \circ
g$. We denote by $\dhatS{X}_{f}^{g}$ the subset of $\dhatS{X}_f $
given as  
\begin{equation}
  \dhatS{X}_{f}^{g} \equiv \hat{X}_{f \circ g} \times V_g^\perp \times \kappa_g.
\end{equation}
We have a natural projection map to a cover of $\dhatC{\square}_f$
\begin{equation}
  \dhatS{X}_{f}^{g} \to \dhatC{\square}^{g}_{f} \cong \hat{\square}_g \times  \kappa_g,
\end{equation}
as shown in Figure \ref{fig:decomposition-square}.
For the statement of the next result, we recall that we have a submersion of smooth orbifolds with corners $  B_{f \circ g} \to B_{f}$.
\begin{defin}
  \label{def:Nash-paths-orbifold}
The \emph{Nash tangent space} $N \dhatS{X}_{f}^{g}$ over $ \dhatS{X}_{f}^{g}$ consists of maps
\begin{equation}
 [0,\infty)^{2} \to \widetilde{X}_f
\end{equation}
such that the restriction to the diagonal inclusion is a Nash path, mapping the origin to $\dhatS{X}_{f}^{g}$, and such that the following properties hold near the origin:
\begin{align} \label{eq:nash-cube-projection-diff}
 & \parbox{31em}{The projection to $\tilde{B}_{f \circ g} $ is a continuously differentiable family of paths.}  \\ \label{eq:nash-cube-fibre-diff}
  & \parbox{31em}{The family of paths parametrised by the first coordinate is a continuously differentiable family of paths mapping to the fibres of the projection to  $\tilde{B}_{f}$.} 
\end{align}
\end{defin}
\begin{rem}
Note that we do not have a projection map $\dhatS{X}_{f} \to
\tilde{B}_{f \circ g}$, but the first condition above nonetheless
makes sense because the restriction to a neighbourhood of
$\dhatS{X}_{f}^{g}$ admits such a map. This is in fact the reason we
are working with the double collared space, as it would not be
sufficient to have the projection map defined on $\dhatS{X}_{f}^{g}$
if we did not know that it extends to a neighbouhood.
\end{rem}

Condition~\eqref{eq:nash-cube-fibre-diff} implies that for a
composition $f \circ g \circ h$ we have a map 
\begin{equation} \label{eq:inclusion-Nash-fibration-cover}
N \dhatS{X}_{f}^{g \circ h} \to  N \dhatS{X}_{f}^{g} 
\end{equation}
along the intersection of $\dhatS{X}_{f}^{g \circ h}$ and
$\dhatS{X}_{f}^{g}$, given by the projection from $\tilde{B}_{f \circ
  g \circ h}  $ to $\tilde{B}_{f \circ g} $. 
This map is natural in the sense that, for a quadruple composition $f
\circ g \circ h \circ k$, we have a commutative diagram 
\begin{equation}
  \begin{tikzcd}
  N \dhatS{X}_{f}^{g \circ h \circ k}   \ar[r] \ar[dr] &  N \dhatS{X}_{f}^{g \circ h}    \ar[d] \\
  &    N \dhatS{X}_{f}^{g}.
  \end{tikzcd}
\end{equation}
This leads us to the following definition:
\begin{defin}\label{defn:Nash-tangent-space} 
The \emph{Nash tangent space} $N \dhatS{X}_{f}$ of $ \dhatS{X}_{f}$ is
the union of the spaces $ N \dhatS{X}_{f}^{g}  $ over all compositions
of cubes $f \circ g$, glued along the maps in Equation
\eqref{eq:inclusion-Nash-fibration-cover}.
\end{defin}

Note that the space $N\dhatS{X}_{f}$ is homeomorphic to the homotopy
colimit of the diagram of spaces $N \dhatS{X}_f^{g}$ over the maps of
Equation~\eqref{eq:inclusion-Nash-fibration-cover}, as these are
closed inclusions.

\begin{lem}\label{lem:nash-inclusion}
Each composition $f \circ g$ in $\Box \Kur^{fs}\langle  \cS
\rangle$ induces a natural inclusion
\begin{equation}
 N \dhatS{X}_{f \circ g} \times V_g^\perp \times \tilde{\kappa}_g \to  N \dhatS{X}_{f}
\end{equation}
of spaces over $\dhatS{X}_{f \circ g}$. 
\end{lem}

\begin{proof}
The map is induced by the inclusion into the colimit coming from the
natural map $\dhatS{B}_{f \circ g} \to \dhatS{B}_{f}$; it is an
inclusion since the fibers over $\dhatS{B}_{f \circ g}$ are contained
in the fibers over $\dhatS{B}_{f}$.
\end{proof}

On the other hand, associated to a composition $f \circ g$ in $\Box
\Kur^{fs}\langle \cS \rangle$ we have a map
\begin{equation}\label{eq:nash-pullback}
 N \dhatS{X}_{f \circ g} \to N \dhatS{X}_{g}
\end{equation}
of Nash tangent spaces over $\dhatS{X}_{g}$, induced by pulling back
paths from $B_{f \circ g}$ to $B_g$.  These constructions are
compatible in the following sense:

\begin{lem}
For each factorization $f_0 = h \circ f_1 \circ g$ in $\Box
\Kur^{fs}\langle \cS \rangle$, there is a natural map
\begin{equation}
N \dhatS{X}_{f_0} \times V_g^\perp \times \tilde{\kappa}_g \to  N  \dhatS{X}_{f_1}
\end{equation}
of spaces over $\dhatS{X}_{f_1}$, which is functorial with respect to
composition in $\Tw \Box \Kur^{fs} \langle \cS \rangle$.
\end{lem}

\begin{proof}
The map in question is constructed as the composite
\begin{equation}
N \dhatS{X}_{f_0} \times V_g^\perp \times \tilde{\kappa}_g \to
N\dhatS{X}_{h \circ f_1} \to N\dhatS{X}_{f_1},
\end{equation}
where the first map is from Lemma~\ref{lem:nash-inclusion} and the
second map is from Equation~\eqref{eq:nash-pullback}.  Checking that
this map is compatible with the composition in the twisted arrow
category amounts to a straightforward verification of the
compatibility of pullback of paths and inclusion of fibers.
\end{proof}

We next formulate the multiplicativity of this construction: as in the
discussion preceding Lemma~\ref{SW-duality-lax-monoidal} (see
Definition~\ref{defin:augmented-Kur-cubes}), we start by defining a
variant of the category of cubes (and abuse notation by referring to
it with the symbol $\square$).

\begin{defin} \label{def:square_of_Kuranishi_tangent}
The category $\square \Kur^{\scrT}$ has
\begin{enumerate}
\item objects given by cubes in $\square \Kur^{\scrT}\langle \cS
\rangle$ for some partially ordered set $\cS$, and
\item morphisms given by a map of cubes, inclusions of boundary
  strata associated to totally ordered subsets of $\cS$, and
  inclusions of strata of cubes.
\end{enumerate}
\end{defin}

We write $\Tw \square \Kur^{\scrT}$ for the twisted arrow category as
always.  The following lemma describing the interaction of the Nash
tangent space with products is now essentially an immediate consequence of
the definitions.

\begin{lem} \label{lem:Nash-tangent-space-monoidal}
Associated to each pair $(f_0, f_1)$ of objects of $\square
\Kur^{\scrT} $ is a homeomrphism
\begin{equation}
N \dhatS{X}_{f_0} \times N \dhatS{X}_{f_1} \to   N \dhatS{X}_{f_0 \times f_1}
\end{equation}
of spaces over $\dhatS{X}_{f_0 \times f_1}$, which is functorial in
each variable in with respect to morphisms in $\Tw \square
\Kur^{\scrT}$ and is associative for triple products in the sense that
the diagram
\begin{equation}
\begin{tikzcd}
N \dhatS{X}_{f_0} \times N \dhatS{X}_{f_1} \times N \dhatS{X}_{f_2}
\ar[r] \ar[d] & \ar[d] N \dhatS{X}_{f_0 \times f_1} \times N
\dhatS{X}_{f_2} \\
N \dhatS{X}_{f_0} \times N \dhatS{X}_{f_1 \times f_2} \ar[r] & N
\dhatS{X}_{f_0 \times f_1 \times f_2} 
\end{tikzcd}
\end{equation}
commutes. 
\qed
\end{lem}

\subsubsection{Construction of the Nash spherical fibration}
\label{sec:constr-spher-fibr}

We begin with the observation that there is a natural map
\begin{equation}
  N \dhatS{X}_{f}^{g} \to   T^{f} X_f \oplus TB_{f \circ g}
\end{equation}
given by the derivatives at the origin, and where we identify the
tangent spaces of collared completions with those of the underlying
manifolds (we recall that $TB_{f \circ g}$ is not globally defined on
$\dhatS{X}_f$, but makes sense on $ \dhatS{X}_{f}^{g} $ where the left
hand side is given). Our proximate goal is to map the right hand side
to $TX_{f}$; we can then define the Nash spherical fibration as the cone
of the complement of the inverse image of $0$.

We start with the restriction of the vector bundle $T^{f} X_f \oplus
TB_{f \circ g}$ to $\dhatS{Z}_{f \circ g} \subset \dhatS{Z}_{f}$.  We
have the natural projection 
\begin{equation}
TB_{f \circ g} \to T B_{f},
\end{equation}
which provides the composite
\begin{equation}\label{eq:comp-proj}
N \dhatS{X}_{f}^{g} \to T^{f} X_f \oplus TB_{f \circ g} \to T^{f} X_f
\oplus T B_{f} \cong T X_f 
\end{equation}
that we use along this subset of $\dhatS{Z}_{f}$.  The compatibility
of these splittings of $T X_f$ with the maps induced on the Nash
tangent space by morphisms in $\Kur^{\scrT}$ implies the following:
 
\begin{lem}
Given a composition $f \circ g \circ h$, the following diagram
\begin{equation}
\begin{tikzcd}
N \dhatS{X}_{f}^{g \circ h}   \ar[dr] \ar[r] &     N \dhatS{X}_{f}^{g}  \ar[d] \\
&  T X_f
\end{tikzcd}
\end{equation}
commutes along $\dhatS{Z}_{f \circ g \circ h} \subset
\dhatS{Z}_{f}$, where the top map is the map of
Equation~\eqref{eq:nash-pullback} and the vertical and diagonal maps
are the projections of Equation~\eqref{eq:comp-proj}. \qed
\end{lem}

Recalling that $\dhatS{X}_{f}^{g}$ is the product of $\hat{X}_{f \circ
  g} \times V_{g}^{\perp}$ with a collar cube $\kappa_g$,
Equation~\eqref{eq:comp-proj} specifies the map from the Nash tangent
space to the ordinary tangent space along the corner stratum of the
cube in which all coordinates vanish (after restricting to the
vanishing locus). We now consider the corner stratum along which all
coordinates are equal to $1$.  We define a vector bundle
\begin{equation}
 T^{f \circ g} X_f \equiv T^{f \circ g}  X_{f \circ g} \oplus
 V_{g}^{\perp} \oplus T \kappa_g
\end{equation}
along this stratum, and observe that the data of morphisms in
$\Kur^{\scrT} \langle \cS \rangle$ gives a surjection  
\begin{equation}
T^{f} X_f  \to  T^{f \circ g} X_f.
\end{equation}
This yields a map
\begin{equation}\label{eq:nash-projection}
\begin{aligned}
N \dhatS{X}_{f}^{g} & \to  T^{f \circ g} X_f  \oplus T B_{f \circ g} \\
& \cong  T X_{f \circ g}  \oplus V_{g}^{\perp} \oplus T \kappa_g   \cong T X_f, 
\end{aligned}
\end{equation}
where we use
Lemma~\ref{lem:tangent-space-split-morphisms-well-defined} to identify
the tangent spaces.

In order to compare the map of Equation~\eqref{eq:nash-projection} to
previous constructions, we consider the map 
$\dhatS{Z}_{f \circ g} \to \dhatS{Z}_f$ associated to setting the
collar coordinates equal to $1$, and the projection $\dhatS{Z}_{g} \to
\dhatS{Z}_{f \circ g}$. 
 
\begin{lem}
The following diagrams of spaces over $\dhatS{Z}_{f \circ g} $ and
over $\dhatS{Z}_g$ commute: 
  \begin{equation}
    \begin{tikzcd}
      N \dhatS{X}^{\id}_{f \circ g} \times V_{g}^{\perp} \times \tilde{\kappa}_g \ar[r] \ar[d] & N \dhatS{X}_{f}^{g} \ar[d] \\
      T X_{f \circ g} \oplus V_{g}^{\perp} \oplus T \kappa_g \ar[r] & T X_f 
    \end{tikzcd}
  \end{equation}
   \begin{equation}
    \begin{tikzcd}
      N \dhatS{X}^{\id}_{f \circ g}  \ar[r] \ar[d] & N \dhatS{X}^{\id}_{g} \ar[d] \\
      T X_{f \circ g} \ar[r] & T X_g.
    \end{tikzcd}
  \end{equation} \qed
\end{lem}

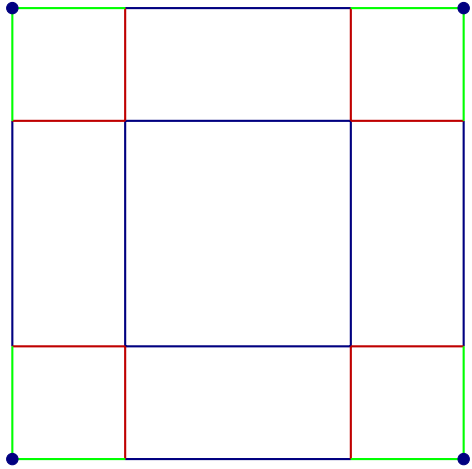
\begin{figure}[h]
  \centering
\begin{tikzpicture}
  \draw [thick, blue!50!black] (-1.5,-1.5) -- (1.5, -1.5) -- (1.5,1.5) -- (-1.5, 1.5) -- cycle;
  \draw [thick, blue!50!black] (-1.5,-3) -- (1.5,-3);
  \draw [thick, blue!50!black] (-1.5,3) -- (1.5,3);
  \draw [thick, blue!50!black] (-3,-1.5) -- (-3,1.5);
  \draw [thick, blue!50!black] (3,-1.5) -- (3,1.5);

  \draw [thick, red!75!black] (-1.5,-3) -- (-1.5,-1.5);
  \draw [thick, red!75!black] (-1.5,1.5) -- (-1.5,3);
  \draw [thick, red!75!black] (1.5,-3) -- (1.5,-1.5);
  \draw [thick, red!75!black] (1.5,1.5) -- (1.5,3);

   \draw [thick, red!75!black] (-3,-1.5) -- (-1.5,-1.5);
  \draw [thick, red!75!black] (1.5,-1.5) -- (3,-1.5);
  \draw [thick, red!75!black] (-3,1.5) -- (-1.5,1.5);
  \draw [thick, red!75!black] (1.5,1.5) -- (3,1.5);

    \draw [thick, green] (-3,-3) --(-3,-1.5);
  \draw [thick, green] (-3,-3) -- (-1.5,-3);
  \draw [thick, green] (-3,1.5) -- (-3,3);
  \draw [thick, green] (-3,3) -- (-1.5,3);

   \draw [thick, green] (3,-3) --(3,-1.5);
  \draw [thick, green] (3,-3) -- (1.5,-3);
  \draw [thick, green] (3,1.5) -- (3,3);
  \draw [thick, green] (3,3) -- (1.5,3);

    \draw[blue!50!black,fill=blue!50!black] (-3,-3) circle (.5ex);
  \draw[blue!50!black,fill=blue!50!black] (3,-3) circle (.5ex);
  \draw[blue!50!black,fill=blue!50!black] (-3,3) circle (.5ex);
  \draw[blue!50!black,fill=blue!50!black] (3,3) circle (.5ex);
  
  \end{tikzpicture}  
\caption{Boundary conditions along the square.\label{fig:colored_boundary-conditions}}

\end{figure}

The above procedure fixes the map of tangent spaces along the (dark)
blue region in Figure \ref{fig:colored_boundary-conditions} (as well
as in the interior of the middle square). We shall now extend to the
remainder by interpolation.  This will involve choices; we need to
ensure that the space of choices is nonempty and contractible.

\begin{defin}
We define $\scrS_{f}^g$ to be the space of splittings of the diagonal
inclusion 
\begin{equation}
  T^{f} B_{f \circ g} \to   T^{f} B_{f \circ g} \oplus   T^{f} B_{f \circ g},
\end{equation}
topologized as a subspace of the space of maps.
\end{defin}

The space $\scrS_{f}^g$ will control the choices of the map on
$\dhatS{X}_{f}^{g}$.  It is clear that $\scrS_{f}^g$ is contractible, 
and that the two choices used so far correspond to the projection to
the first and the second factors, via the decomposition 
\begin{equation}
T^f X_f \oplus T B_{f \circ g} \cong  T^{f \circ g} X_f  \oplus  T^{f}
B_{f \circ g} \oplus T^{f} B_{f \circ g} \oplus T B_{f}, 
\end{equation}
which makes sense in $\dhatS{X}_{f}^{g}$.

Given a composition $g \circ h$, we have a splitting 
\begin{equation}\label{eq:tangent-base-splitting}
T^{f} B_{f \circ g \circ h} \cong  T^{f} B_{f \circ g} \oplus T^{f
  \circ g} B_{f \circ g \circ h}. 
\end{equation}
Along the intersection of $\dhatS{X}_{f}^{g}$ with $\dhatS{X}_{f}^{g
  \circ h}$, recall that we have a natural map 
\begin{equation}
 T^{f \circ g} X_f  \to T^{f \circ g \circ h} X_f .
\end{equation}
We fix the map
\begin{equation} \label{eq:map_splittings-interior-compatibility}
  \scrS_{f}^g \to \scrS_{f}^{g \circ h}
\end{equation}
given by assigning to a splitting $T^{f} B_{f \circ g} \oplus   T^{f}
B_{f \circ g} \to T^{f} B_{f \circ g}$ the direct sum with the map 
\begin{equation}
\pi_2 \colon  T^{f \circ g} B_{f \circ g  \circ h} \oplus T^{f \circ g} B_{f \circ
    g  \circ h} \to T^{f \circ g} B_{f \circ g  \circ h} 
\end{equation}
given by projection to the second factor and using
Equation~\eqref{eq:tangent-base-splitting}.  Here, we recall that the
two factors are ordered in such a way that the first corresponds to
the splitting of the fibrewise tangent spaces of $X_f$, and the second
to the splitting of the (fibrewise) tangent spaces of $B_f$.  This
specifies the compatibility conditions along the red (medium dark)
labelled arcs in Figure~\ref{fig:colored_boundary-conditions}.

Next, consider the image of $\dhatS{X}^{g \circ h}_{g} \subset
\dhatS{X}_g$ in $ \dhatS{X}^{f \circ g \circ h}_{f}$. We have an
associated natural map 
\begin{equation}
T^{g \circ h} X_{g} \oplus   V_g^\perp \oplus T \kappa_g \cong   T^{f
  \circ g \circ h} X_{f}, 
\end{equation}
as well as a direct sum decomposition
\begin{equation}\label{eq:tangent-base-splitting2}
  T^{f} B_{f \circ g \circ h} \cong T^g B_{g \circ h} \oplus T^f B_{f \circ g}.
\end{equation}
We fix the map
\begin{equation} \label{eq:map-splitting-boundary-compatibility}
  \scrS_{g}^h \to \scrS_{f}^{g \circ h}
\end{equation}
given by assigning to a splitting $T^{g} B_{g \circ h} \oplus T^{g}
B_{g \circ h} \to T^{g} B_{g \circ h}$ its direct sum with the map
\begin{equation}
\pi_1 \colon T^f B_{f \circ g}  \oplus T^f B_{f \circ g}  \to  T^f B_{f \circ g}. 
\end{equation}
given by projection to the first factor and using
Equation~\eqref{eq:tangent-base-splitting2}.  This specifies the
compatibility conditions along the (light) green labelled arcs in
Figure~\ref{fig:colored_boundary-conditions}.

We now turn to the question of how to make a global choice of elements
of $\scrS_f^g$ for each composable pair $(f,g)$ that are suitably
compatible.  We will write
\begin{equation}
\Ar(\square \Kur^{\scrT})_{2} \equiv \Ar(\square \Kur^{\scrT})
\times_{\ob(\square \Kur^{\scrT})} \Ar(\square \Kur^{\scrT})
\end{equation}
to denote the space of pairs of composable arrows in $\square
\Kur^{\scrT}$ and more generally $\Ar(\square \Kur^{\scrT})_{k}$ for
the space of $k$-tuples of composable arrows.  We write $\scrS$ for
the total space of the parametrized space over the indexing space
$\Ar(\square \Kur^{\scrT})_{2}$ with fiber $\scrS_f^g$.  A section of
this parametrized space is precisely a continuous choice of section in
$\scrS_f^g$ for each composable pair $(f,g)$.  There are analogous
parametrized spaces $\scrS_k$ over $\Ar(\square \Kur^{\scrT})_{k}$,
and there is a natural continuous map from sections of $\scrS$ to
sections of $\scrS_k$ for each $k$.

Composition in $\Ar(\square \Kur^{\scrT})$ specifies two maps
\begin{equation}
\Ar(\square \Kur^{\scrT})_{3} \to \Ar(\square \Kur^{\scrT})_{2},
\end{equation} and together with Equations~\eqref{eq:map_splittings-interior-compatibility}
and~\eqref{eq:map-splitting-boundary-compatibility} we obtain induced
maps $\circ^{\scrS}_1$ and $\circ^{\scrS}_2$ from sections of the
parametrized space $\scrS_3$ to sections of $\scrS$.

\begin{defin}
We say that a section of $\scrS$ (i.e., choice of sections in 
$\scrS_f^g$ for all composable morphisms $f \circ g$) is
\emph{compatible} if, for each triple composition $f \circ g \circ h$, 
the values of the induced section of $\scrS_3$ coincide under the maps
$\circ^{\scrS}_1$ and $\circ^{\scrS}_2$.
\end{defin}

A priori, it is not clear that the space of compatible sections is
nonempty.  The folllowing result that it is both nonempty and
contractible, so in particular we can choose one to work with in our
constructions.

\begin{lem}\label{lem:contractible-space-of-sections}
The space of compatible sections of $\scrS$ for all composable
morphisms $f \circ g$ is nonempty and contractible. 
\end{lem}

\begin{proof}
We begin by showing that the space is nonempty.  The proof proceeds by
(double) induction: we simultaneously induct on the dimension of the
collar cubes and the collection of composable morphisms.
For cubes of dimension $0$, the space of choices of compatible
sections is a point.  Now, assuming that we have chosen compatible
sections of all cubes of dimension strictly less than $n$, we extend
the choices as follows.  We pick an arbitrary order on the cubes $f$
of dimension $n$.  For each such cube, we pick an order on the
composable morphisms $f \circ g$, and proceed to choose compatible
extensions by decreasing induction on the dimension of the domain of
$g$.  We can do this by interpolating arbitrarily between the
boundaries.  Since the choices for each cube are independent of one
another, there is no obstruction to completing the selection of
compatible choices in dimension $n$.  This establishes that the space
of compatible sections is nonempty.

To show that the space of compatible sections is contractible, we
choose a compatible section $s$.  We construct a homotopy from the
identity map to the constant map at $s$ using the same inductive
approach as above.  Specifically, we proceed by a double induction,
choosing at each stage compatible homotopies to the constant map on
$s$; we can choose any homotopy on each cube that is compatible on the
boundaries, and it is clear that this can always be done.  The
independence of the choices again mean that the choice of such lifts
is unobstructed.
\end{proof}

Now assuming that we have fixed a choice of compatible section of
$\scrS_f^g$, we see that for each composition $f \circ g$, we have a
natural commutative diagram  
\begin{equation} \label{eq:product-pullback-map-Nash-paths}
\begin{tikzcd}
N \dhatS{X}_{f \circ g} \times V_g^\perp \times \tilde{\kappa}_{g}
\ar[r] \ar[d] &  N \dhatS{X}_f \ar[d] \\ 
T X_{f \circ g}  \times V_g^\perp \times T \kappa_{g} \ar[r] & T X_f
\end{tikzcd}
\end{equation}
of spaces over $\dhatS{X}_{f \circ g}$, and a commutative diagram
\begin{equation}
  \begin{tikzcd}
    N \dhatS{X}_{f \circ g} \ar[r] \ar[d] & N \dhatS{X}_g \ar[d] \\
    T X_{f \circ g} \ar[r] & T X_g
  \end{tikzcd}
\end{equation}
of spaces over $\dhatS{Z}_g$.

\begin{defin} \label{def:Nash-path-cubical-degeneration}
The \emph{Nash spherical fibration} $\Nash \dhatS{X}_f$ over
$\dhatS{X}_f$ is the spherical fibration obtained by taking the cone
of the complement of the inverse image of $T X_f \setminus 0$.
\end{defin}

A key observation about this construction is that it only depends on
the component of $f$.

\begin{lem}
Let $f_1$ and $f_2$ be morphisms in $\square \Kur^{\scrT}$ such that
$f_1$ and $f_2$ are in the same component of the mapping space.  Then
$\Nash \dhatS{X}_{f_1}$ and $\Nash \dhatS{X}_{f_2}$ coincide.
\end{lem}

\begin{proof}
It is clear that $\dhatS{X}_{f_0}$ and $\dhatS{X}_{f_1}$ coincide; the
relevant quotients are independent of the choice of inner product on
the tangent space.  Next, observe that the total space $N
\dhatS{X}_{f}$ of the Nash tangent space as constructed in
Definition~\ref{defn:Nash-tangent-space} is also manifestly
independent of the choice of inner product.  The map to $\dhatS{X}_f$
constructed in this section does depend on the use of the inner
product data in the construction of the map in
Equation~\eqref{eq:nash-projection} (via
Lemma~\ref{lem:tangent-space-split-morphisms-well-defined}).  However,
the definition of the Nash spherical fibration does not depend on the
inner product, because the constructed map is independent of this data
(and the splitting identification) near $0$.
\end{proof}

The continuity condition for derivatives at the origin implies that
cubes whose directional derivatives do not vanish in $T X_f $ restrict
to a non-constant path along the diagonal. Evaluation of these
derivatives and of the value at the point $(1,1)$ thus induces natural
maps
\begin{equation} \label{eq:evaluation_maps_Nash-before-completion}
  \Milnor \dhatS{X}_f \leftarrow \Nash \dhatS{X}_f \to S^{T X_f}|0
\end{equation}
of parametrized spaces over $\dhatS{Z}_f$ (here recall the definition
of the Milnor spherical fibration from
Definition~\ref{defn:milnor-spherical-fibration}).  We argue below
that the induced maps are fiberwise equivalences.  We would like to
immediately conclude that therefore we have induced equivalences on
spaces of sections.  However, the projection $\Nash \dhatS{X}_f \to
\dhatS{Z}_f$ is a quasifibration but not a fibration.  For any
quasifibration, the path space construction gives rise to an
associated fibration, and it is standard to compute spaces of sections
in terms of this fibration; in general, spaces of sections computed
directly using the quasifibration might not have the correct homotopy
type.  Nonetheless, in the current setting, we do not have to perform
the replacement in order to work with sections.

\begin{lem}
The projection map $\Nash \dhatS{X}_f  \to \dhatS{Z}_f$ is a
quasifibration and the canonical map from the space of sections of the
projection to the space of sections of the associated fibration is a
weak equivalence.  The induced maps to the Milnor and Whitney
spherical fibrations are equivalences.
\end{lem}

\begin{proof}
We begin by arguing that $\Nash \dhatS{X}_f \to \dhatS{Z}_f$ is a
quasifibration such that the space of sections has the correct
homotopy type.  We rely on
Proposition~\ref{prop:gluing-good-quasifibrations}.  Considering the
construction of the Nash spherical fibration, observe that it is built
by gluing together fibrations over the cubical decomposition of the
double collar $\dhatS{X}_f$; in fact, by
Proposition~\ref{prop:Nash-equivalence-square}, restricted to the
pieces of the cubical cover we have fiber bundles.  Since subcubes in
the decomposition include by cofibrations essentially by construction,
the argument of Proposition~\ref{prop:gluing-good-quasifibrations}
applies to show that the comparison map of sections spaces is a weak
equivalence.  (The argument also shows that the induced parametrized
space is an ex-quasifibration.)

Next, we show that the comparison maps are equivalences, again
extending the argument of
Proposition~\ref{prop:Nash-equivalence-square}.
The projections maps
$\dhatS{X}_f^g \to \dhatS{P}_{g} \to \dhatS{P}_{f}$ are submersions,
hence are locally modelled after
  \begin{equation}
    \bR^{n+k+\ell} \to \bR^{n+k} \to \bR^{k}. 
  \end{equation}
At each point, we can choose such a local model so that the choices of
splitting in $\scrS^g_f$ identify the map 
  \begin{equation}
 T^{f} X_f \oplus T B_{f \circ g} \to   T^{f} X_f \oplus T B_f  
  \end{equation}
  with the projection
  \begin{equation}
    \bR^{\ell + n} \times  \bR^{n + k} \to  \bR^{\ell} \times \bR^{n} \times \bR^{k},
  \end{equation}
where the map on the factors $\ell$ and $k$ are the identity, and on
the $\bR^n \times \bR^n$ factors is a linear splitting of the identity
map.   Given a pair $(v_1,v_2)$ of vectors in the left hand side, we
define $\gamma_{v_1,v_2}$ to be the map from $[0,\infty)^2$ to
  $\bR^{n+k+\ell}$ given by radially extending 
  \begin{equation}
    (t_1, t_2) \mapsto t_1 v_2 + t_2 v_2    
  \end{equation}
  from $[0,1]^2$. This defines a splitting of the map 
  \begin{equation}
    N \dhatS{X}^g_f \to T X_f,
  \end{equation}
 thus inducing a map $S^{T X_f}|0 \to \Nash \hat{X}_f$.  Nash's
 argument then shows that this is a fibrewise homotopy
 equivalence.  We can now evaluate further to $T X_f | 0$ and glue.  
As in the argument above, we glue these fiberwise homotopy
equivalences over the cubical decomposition of $\dhatS{X}_f$; by
Proposition~\ref{prop:gluing-good-quasifibrations}, the assembled maps
are fiberwise equivalences.
\end{proof}

\subsubsection{Functoriality of the Nash fibrations}
\label{sec:funct-nash-fibr}

It is straightforward to see that the maps constructed in the previous
section induce a pullback map 
\begin{equation}
 \Nash \dhatS{X}_{f \circ g} \to \Nash \dhatS{X}_g.  
\end{equation}
As with the Milnor model of the tangent fibration, the key problem is
to make these spaces compatible with desuspension to make sense of
pushforward.

We start by adding a third collar, and constructing a Nash fibration
\begin{equation}
 \Nash \thatS{X}_f \to \thatS{X}_f  
\end{equation} 
by extending the definition of $N \thatS{X}_f$ in an invariant way across the third set of collars, as well as the map to $T X_f$. In this way, we obtain a canonical map
\begin{equation}
   \Nash \thatS{X}_{f \circ g}    \barwedge V_g^\perp|0 \barwedge   M \kappa_g \to \Nash \thatS{X}_f,
 \end{equation}
 of spaces over $\thatS{Z}_{f \circ g} \times \kappa_g$, extending the
 map from Equation \eqref{eq:product-pullback-map-Nash-paths} to a
 collar, and taking the appropriate cones. While the spaces $
 V_g^\perp|0 $ and $S^{V_g^\perp}$ are homotopy equivalent, the lack
 of a canonical equivalence between them makes it difficult to use
 this.

Lemma~\eqref{eq:arrow_category_cubes_diagram_pairs} leads us instead
to consider the following pushout diagram 
\begin{equation}
\begin{tikzpicture}
\node[scale=.9]{
\begin{tikzcd}
\displaystyle{\coprod_{\pi \overset{h}{\to} \rho \overset{g}{\to} \tau
    \overset{f}{\to} \sigma} \Nash  \thatS{X}_{f \circ g \circ
    h}\barwedge M \kappa_{g} \barwedge M \kappa_{h} \barwedge
  (V_h^\perp|0 \wedge S^{V_g^\perp}|0 ) } \ar[r] \ar[d] &
\displaystyle{\coprod_{\rho \overset{f \circ g}{\to} \sigma}     \Nash
  \thatS{X}_{f \circ g} \barwedge M \kappa_{g} \barwedge
  S^{V_g^\perp}|0   }\ar[d] \\ 
\displaystyle{\coprod_{\rho \overset{f \circ g \circ h}{\to} \sigma}
  \Nash  \thatS{X}_{g \circ h} \barwedge M \kappa_{g \circ h}
  \barwedge S^{V_{g \circ h}^\perp}|0   }\ar[r] & \Nash
\thatB{\bar{X}}_f.
\end{tikzcd}
};
\end{tikzpicture}
\end{equation}

\begin{lem}
The inclusion $ \Nash \thatS{X}_f \to  \Nash \thatB{\bar{X}}_f $ is an
equivalence of spherical fibrations over $\thatS{X}_f$. 
\end{lem}
\begin{proof}
The space $\Nash \thatS{X}_f$ is the pushout of the diagram 
\begin{equation}
\begin{tikzpicture}
\node[scale=.9]{
\begin{tikzcd}
      \displaystyle{\coprod_{\pi \overset{h}{\to} \rho \overset{g}{\to} \tau \overset{f}{\to} \sigma} \Nash  \thatS{X}_{f \circ g \circ h}\barwedge M \kappa_{g} \barwedge M \kappa_{h} \barwedge (V_h^\perp|0 \wedge V_g^\perp|0 ) } \ar[r] \ar[d] &      \displaystyle{\coprod_{\rho \overset{f \circ g}{\to} \sigma}     \Nash  \thatS{X}_{f \circ g} \barwedge M \kappa_{g} \barwedge V_g^\perp|0   }\ar[d] \\
  \displaystyle{\coprod_{\rho \overset{f \circ g \circ h}{\to} \sigma} \Nash  \thatS{X}_{g \circ h} \barwedge M \kappa_{g \circ h} \barwedge V_{g \circ h}^\perp|0   }\ar[r] & \Nash \thatS{X}_f.
\end{tikzcd}
};
\end{tikzpicture}
\end{equation}
The natural map $\Nash \thatS{X}_f \to \Nash \thatB{\bar{X}}_f $ is
induced by the inclusion of the first factor in the bottom left corner
of the above diagram; specifically the map $\Nash \thatS{X}_\sigma \to
\Nash \hat{\bar{X}}_\sigma $ is induced by the map $V|0 \to S^V|0$.
Since this latter map is a bijection, so is the comparison map,
although it is not usually a homeomorphism.  To see that it is a
homotopy equivalence, we use excision; as observed in
Remark~\ref{rem:change-topology}, the map can be regarded as changing
the topology near the basepoint in each fiber.  Since we can choose a
neighborhood $U$ of the basepoint section such that the comparison map
is a homeomorphism on the complement, the inclusion of $U$ into the
fiber is an NDR-pair, and $U$ is contractible, the map is a homotopy
equivalence.
\end{proof}

The naturality of this construction implies that the restriction of
$\Nash \thatB{\bar{X}}_f $ to the collar labelled by $f \circ g$ is
equipped with a natural map 
\begin{equation} \label{eq:map_Nash-spaces-ordinary-suspension}
\Nash \thatB{\bar{X}}_{f \circ g}   \barwedge  S^{V_g^\perp} \barwedge
M \kappa_g  \to \Nash \thatB{\bar{X}}_f. 
\end{equation}

In fact, it is straightforward to verify that this construction is
functorial in the twisted arrow category.

\begin{lem}
Each factorization $f_0 = h \circ f_1 \circ g$ induces a natural
equivalence 
\begin{equation}
\Nash \thatB{\bar{X}}_{f_1}\wedge S^{V_g^\perp} \wedge M \kappa_g
\to   \Nash \thatB{\bar{X}}_{f_0} 
\end{equation}
of spherical fibrations over $ \thatS{Z}_{f_1} \times \kappa_g $, and
these equivalences are compatible with the composition of
factorization.  That is, given composable factorizations $f_0 = h
\circ f_1 \circ g$ and $f_1 = h' \circ f_2 \circ g'$, the equivalence 
\begin{equation}
\Nash \thatB{\bar{X}}_{f_2}\wedge S^{V_{g' \circ g}^\perp} \wedge M
\kappa_{g' \circ g} \to \Nash \thatB{\bar{X}}_{f_0} 
\end{equation}
coincides with the composite 
\begin{equation}
\begin{aligned}
\Nash \thatB{\bar{X}}_{f_2}\wedge S^{V_{g' \circ g}^\perp} \wedge M
\kappa_{g' \circ g} &\to
\Nash \thatB{\bar{X}}_{f_2}\wedge S^{V_{g'}^\perp} \wedge M
\kappa_{g'} \sma S^{V_{g}^\perp} \sma M\kappa_{g} \\
&\to \Nash
\thatB{\bar{X}}_{f_1} \sma S^{V_{g}^\perp} \sma M\kappa_{g} \to  
\Nash \thatB{\bar{X}}_{f_0}.
\end{aligned} 
\end{equation}
\qed
\end{lem}

\subsubsection{Multiplicativity of Nash fibrations}
\label{sec:mult-nash-fibr}

We now consider the Nash fibration for products of charts.  Building
on Lemma~\ref{lem:Nash-tangent-space-monoidal}, we have the following
product map. 

\begin{cor}\label{cor:nash-product-choices}
For each pair $(f_0 \circ g_0, f_1 \circ g_1)$ of composable morphisms
in the category $\square \Kur^{\scrT}\langle  \cS_i \rangle$, there is a natural map
\begin{equation}
\scrS_{f_0}^{g_0} \times  \scrS_{f_1}^{g_1} \to \scrS_{f_0 \times
  f_1}^{g_0 \times g_1}
\end{equation}
which is associative for composable triples in the sense that the following
diagram commutes:
\begin{equation}
\begin{tikzcd}
\scrS_{f_0}^{g_0} \times  \scrS_{f_1}^{g_1} \times \scrS_{f_2}^{g_2}
\ar[r] \ar[d] & \scrS_{f_0 \times f_1}^{g_0 \times g_1} \times
\scrS_{f_2}^{g_2} \ar[d] \\
\scrS_{f_0}^{g_0} \times \scrS_{f_1 \times f_2}^{g_1 \times g_2}
\ar[r] & \scrS_{f_0 \times f_1 \times f_2}^{g_0 \times g_1 \times g_2} .
\end{tikzcd}
\end{equation} \qed
\end{cor}

In addition, for $Q$ a totally ordered subset of $\cS$, we have a
restriction map  
\begin{equation}\label{eq:boundary-section-restriction}
\scrS_{f}^{g} \to \scrS_{\partial^Q f}^{\partial^Q g}, 
\end{equation}
and these restriction maps are functorial under inclusions of totally
ordered subsets $Q' \subset Q \subset \cS$.

We can use the maps of Corollary~\ref{cor:nash-product-choices} and
Equation~\eqref{eq:boundary-section-restriction} to inductively choose
sections of $\scrS_{f}^g$ which are compatible with products.  The
multiplicative compatibility condition amounts to the following: the
product of the sections of $\scrS_{f_1}^{g_1}$ and $\scrS_{f_2}^{g_2}$
on $\dhatS{Z}_{f_1}^{g_1} \times \dhatS{Z}_{f_2}^{g_2} $ agrees with
the restriction of the section of $\scrS_{f_1 \times f_2}^{g_1 \times
  g_2} $ on $\dhatS{Z}_{f_1 \times f_2}^{g_1 \times g_2}$.  A small
modification of the double induction in
Lemma~\ref{lem:contractible-space-of-sections} then lets us find
sections which are compatible with the product structure, since the
data of the morphisms in $\square \Kur^{\scrT}$ ensures that the
required inductive choices are unobstructed.  Finally, suppose that
$\cP$ is equipped with a free action of $\Pi$.  Given any finite
totally ordered subset $Q \subset \cP$, we can make a choice of
inductive extension of the section for the orbit of each poset in
$\cP(Q)$ under $\Pi$ which is compatible with products, as these are
independent of one another.  

\begin{lem}\label{lem:fix_contractible_choice-sections-category}
The space of compatible sections of $\scrS_f^g$ for all partially
ordered sets $\cS$ and all composable pairs $f \circ g$ in $\square
\Kur^{\scrT}\langle \cS \rangle$ which are compatible with restriction
maps to boundary strata and with products is nonempty and
contractible.  Given a partially ordered set $\cP$ equipped with a
free action by a group $\Pi$, the choices for the partially ordered
sets $\cP(Q)$ obtained from finite totally ordered subsets $Q$ of
$\cP$ can be made equivariant under the action of $\Pi$.
\qed
\end{lem}

\subsubsection{The Nash-twisted cochains}
\label{sec:nash-twist-coch}

We now compare the constructions of Section \ref{sec:whitn-spher-fibr}
to the analogous cochains built from the Nash spherical fibration.
Given a spectrum $\bbk$, 

\begin{defin}
For each arrow in $f \in \square \Kur^{\scrT}\langle \cS \rangle $, we
have the spectrum
\begin{equation}
C^{*,c}_{\rel \partial} \left( B\thatS{Z}_f;
\Nash\thatB{\bar{X}}^{-V}_f \wedge \bbk\right)
\end{equation}
of Borel equivariant compactly supported sections which
vanish at the boundary of the parametrized spectrum
$\Nash\thatB{\bar{X}}^{-V}_f \wedge \bbk$ over $\thatS{Z}_f$.
\end{defin}

There is a natural pullback map
\begin{equation}
C^{*,c}_{\rel \partial} \left( B\thatS{Z}_{f \circ g};
\Nash\thatB{\bar{X}}^{-V}_{f \circ g} \wedge \bbk\right) \to
C^{*,c}_{\rel \partial} \left( B\that{Z}_g;
\Nash\thatB{\bar{X}}^{-V}_g \wedge \bbk\right)    
\end{equation}
for each composition. Fixing the section of $\Milnor \kappa$ chosen in
Equation \eqref{eq:section_Milnor_collar_vanish_boundary}, and using
the map from Equation \eqref{eq:map_Nash-spaces-ordinary-suspension},
we also obtain a map 
\begin{equation}
C^{*,c}_{\rel \partial} \left( B\thatS{Z}_{f \circ g};
\Nash\thatB{\bar{X}}^{-V}_{f \circ g} \wedge \bbk\right) \to
C^{*,c}_{\rel \partial} \left( B\thatS{Z}_f;
\Nash\thatB{\bar{X}}^{-V}_f \wedge \bbk\right)
\end{equation}
for each composition.

A diagram chase analogous to the one for
Diagram~\eqref{eq:maps_cochains_Whitney_compositions} implies:
\begin{lem}
The assignment
\begin{equation}\label{eq:internal-functor}
\begin{aligned}
 \Tw \square \Kur^{\scrT}\langle \cS \rangle  & \to \Sp \\
  f & \mapsto C^{*,c}_{\rel \partial} \left( B\thatS{Z}_f; \Nash\thatB{\bar{X}}^{-V}_f \wedge \bbk\right) 
\end{aligned}
\end{equation}
specifies a topologically enriched functor of internal categories
equipped with natural equivalences to the functors associated to
Milnor and Whitney twisted cochains.
\end{lem}

\begin{proof}
The only part of the proof that requires additional comment is the
verification that we have a topological functor.  Recall that we
are regarding $\Sp$ as an internal category in topological spaces with
the discrete topology on objects and the usual enrichment on mapping
spaces.  Therefore, the assignment must be locally constant in order
to specify a topological functor.  This follows from the fact that the
assignment of Equation~\eqref{eq:internal-functor} is constant on the
components of $\ob(\square \Kur^{\scrT}\langle \cS \rangle)$, as any
choice of inner product leads to the same spectrum.  From this, we can
deduce that the functor is constant on the components of $\ob(\Tw
\square \Kur^{\scrT}\langle \cS \rangle)$.
\end{proof}

We now pass to (topologized) homotopy colimits.  We review the
required notion of homotopy colimit over a topological category in
Section~\ref{sec:internal}; note that since the functor we are taking
the homotopy colimit of depends only on the component of the argument,
the situation is significantly simpler than the case of a general
functor.

\begin{defin}\label{defin:nash-twisted-cochains}
Given a topologically enriched functor $A^{\scrT} \to \square
\Kur^{\scrT}\langle \cS \rangle$, we define a spectrum of cochains in
terms of the homotopy colimit 
\begin{equation}
C^{*,c}_{\rel \partial} \left( B\thatS{Z}; \Nash\thatB{\bar{X}}^{-V}
\wedge \bbk\right)(A^{\scrT}) \equiv \hocolim_{f \in \Tw \square
  A^{\scrT}}   C^{*,c}_{\rel \partial} \left( B\thatS{Z}_f;
\Nash\thatB{\bar{X}}^{-V}_f \wedge \bbk\right). 
\end{equation}
\end{defin}
Applying the previous discussion, we conclude:
\begin{lem}
There is a natural equivalence from $C^{*,c}_{\rel \partial} \left(
B\thatS{Z}; \Nash\thatB{\bar{X}}^{-V} \wedge \bbk\right)(A^{\scrT}) $
to the corresponding construction involving the Whitney model. \qed
\end{lem}

Fixing a choice as in
Lemma~\ref{lem:fix_contractible_choice-sections-category} ensures that the 
assignment of Definition~\ref{defin:nash-twisted-cochains} is
multiplicative, using the natural product map on topological homotopy
colimits.  Therefore, we can conclude the following result. 

\begin{prop} 
Given a Kuranishi flow category with tangent bundles $\bX$, there is a $\Pi$-equivariant spectral category
$C^{*,c}_{\rel \partial} \left( B\thatS{\cZ};
\Nash\thatB{\bar{\cX}}^{-V-d} \sma \bbk \right)$, with objects those of $\bX$ and with
morphism spectra 
\begin{equation}
\begin{aligned}
C^{*,c}_{\rel \partial} \left(  B\thatS{\cZ};
\Nash\thatB{\bar{\cX}}^{-V-d} \sma \bbk \right)(p,q) &\equiv \\
& C^{*,c}_{\rel \partial} \left(  B\thatS{Z}; \Nash
  \thatB{\bar{X}}^{-V - V_{p} + V_{q}} \sma \bbk \right)(p,q)[ \deg p
    - \deg q ].
\end{aligned}
\end{equation} \qed
\end{prop}
       
Returning to Lemma \eqref{eq:arrow_category_cubes_diagram_pairs}, we
conclude that the evaluation maps in Equation
\eqref{eq:evaluation_maps_Nash-before-completion} extend to maps of
these completions 
  \begin{equation} \label{eq:evaluation_map_Nash-Milnor_Whitney}
\Milnor  \thatB{\bar{X}}_f     \leftarrow  \Nash \thatB{\bar{X}}_f \to S^{T X_f}|0
  \end{equation}
over $\thatS{Z}$.  We use these completions to compare $C^{*,c}_{\rel
  \partial} \left( B\thatS{\cZ}; \Nash\thatB{\bar{\cX}}^{-V-d} \sma
\bbk\right)$ to the other spectral categories we have constructed.

\begin{prop}\label{prop:nash-whitney-comparison-zigzag}
The maps in Equation~\eqref{eq:evaluation_map_Nash-Milnor_Whitney} are
functorial and yield $\Pi$-equivariant equivalences of spectral
categories 
\begin{equation}
\begin{tikzcd} 
 C^{*,c}_{\rel \partial} \left( B\thatS{\cZ}_f;   S^{T\cX -V-d} \wedge
 \bbk\right) & \ar[l] \ar[d] C^{*,c}_{\rel \partial} \left( B\thatS{\cZ};
 \Nash\thatB{\bar{\cX}}^{-V-d} \wedge \bbk\right)   \\
& C^{*,c}_{\rel
   \partial} \left( B\thatS{\cZ};   \Milnor\thatB{\bar{\cX}}^{-V-d}
 \wedge \bbk\right),
\end{tikzcd}
\end{equation}
where the first and last categories are obtained from the constructions of Sections \ref{sec:tang-spac-whitn} and \ref{sec:span-whit-dual-3} by using three collars. 
\end{prop}

\begin{proof}
The comparisons in
Equation~\eqref{eq:evaluation_map_Nash-Milnor_Whitney} are
multiplicative and independent of the choice of compatible section in
Lemma~\ref{lem:fix_contractible_choice-sections-category}.  Moreover,
they are evidently $\Pi$-equivariant, since the action of $\Pi$ simply
permutes the morphism spectra.
\end{proof}

Collapsing the inner two collars is multiplicative and
$\Pi$-equivariant, and provides a natural comparison with
the previously constructed categories:

\begin{prop}
Given a Kuranishi flow category with tangent bundle $\bX$, there is a
zig-zag of $\Pi$-equivariant equivalences of spectral categories
induced by collapsing collars:
\begin{equation}
\begin{tikzcd}
C^{*,c}_{\rel \partial} \left(  B\hat{\cZ}; S^{T \hat{\cX} - V -d}
\sma \bbk \right) & \ar[l] C^{*,c}_{\rel \partial} \left( B\thatS{\cZ};
 \Nash\thatB{\bar{\cX}}^{-V-d} \wedge \bbk\right) \ar[d] \\
&  C^{*,c}_{\rel \partial}
\left( B\hat{\cZ}; \Milnor\hat{\bar{\cX}}^{-V-d} \wedge \bbk\right)
\end{tikzcd}
\end{equation}
\end{prop}

\begin{proof}
The zig-zag is the composition of the comparison of
Proposition~\ref{prop:nash-whitney-comparison-zigzag} with the
comparison maps
\begin{equation}
C^{*,c}_{\rel \partial} \left( B\thatS{\cZ}_f;   S^{T\cX -V-d} \wedge
 \bbk\right)
\to
C^{*,c}_{\rel \partial} \left(  B\hat{\cZ}; S^{T \hat{\cX} - V -d}
\sma \bbk \right)
\end{equation}
and
\begin{equation}
C^{*,c}_{\rel
   \partial} \left( B\thatS{\cZ};   \Milnor\thatB{\bar{\cX}}^{-V-d}
 \wedge \bbk\right)
\to
C^{*,c}_{\rel \partial}
\left( B\hat{\cZ}; \Milnor\hat{\bar{\cX}}^{-V-d} \wedge \bbk\right)
\end{equation}
induced by the collapse of the inner collars and the induced maps on
homotopy colimits coming from the multiplicative functors $A^{\scrT}
\to A$.  Since by hypothesis (recall
Definition~\ref{def:flag-smooth-split-presentation}) the fibers of
$A^{\scrT} $ over $A$ are contractible, the source and target maps in
$A^{\scrT}$ are contractible, and the object space in $A$ is discrete,
the hypotheses of Proposition~\ref{prop:internal-Quillen-A} (the
internal version of Quillen's theorem A) are satisfied and we can
conclude that these maps are weak equivalences.
\end{proof}

\subsubsection{Signpost: Comparison of tangentially twisted spherical fibrations}
\label{sec:signp-comp-tang}

At this point, we have constructed a $\Pi$-equivariant zig-zag of
DK-equivalences of spectral categories
\begin{equation}
C^{*,c}_{\rel \partial}(B \hat{\cZ}; S^{T \hat{\cX}-V-d}) \to \cdots
\leftarrow B \bar{\cX}|\cZ^{-V -d}.
\end{equation}
What remains is the comparison between tangentially twisted and
ordinary cochains.

\subsection{Complex-oriented flow categories}
\label{sec:compl-orient-flow-1}

The purpose of this section is to compare the tangentially twisted
cochains of a flag smooth Kuranishi presentation equipped with a
relative complex structure in the sense of
Section~\ref{sec:compl-orient-kuran-1} with the ordinary cochains.  The 
comparison proceeds in three steps: 

\begin{enumerate}

\item We construct a category consisting of cochains twisted by the
  (stable) complex vector bundle appearing in the definition of
  relative complex orientations.

\item We compare this to the tangentially twisted cochains via a
  bimodule representing an equivalence.

\item We complete the argument by using complex-orientability
  to trivialize the cochains twisted by a complex vector bundle. 

\end{enumerate}

\subsubsection{Bimodule comparison between twisted cochains}
\label{sec:compl-stable-struct}

Let $V_0 = (V_0^+, V_0^-)$ and $V_1 = (V_1^+, V_1^-)$ be a pair of
stable vector spaces, with complex structures on $V^-_0$ and $V^-_1$.
Associated to $V_0$ and $V_1$ is the internal category
$\Kur^{\ori}(V_0,V_1)$ of Kuranishi charts equipped with stable
complex structures relative $V_0$ and $V_1$ from
Definition~\ref{def:Complex-oriented-Kuranishi-chart}.  We will work
with a distinguished subcategory of the category of cubes on
$\Kur^{\ori}(V_0,V_1)$.

\begin{defin}
In abuse of notation, we will write $\square \Kur^{\ori}(V_0,V_1)$ to
denote the subcategory of cubes on $\Kur^{\ori}(V_0,V_1)$ with:
\begin{enumerate}
\item Objects with fixed stratifying
set (i.e., each cube lands in a given subcategory
$\Kur^{\ori}(V_0,V_1) \langle \cS \rangle$). 
\item Morphisms defined as
in Definition~\ref{def:square_of_Kuranishi_tangent} by a morphism of
cubes with a fixed stratification along with a choice of boundary
stratum.
\end{enumerate}
The topology on the objects and morphisms is the subspace topology
induced from $\Kur^{\ori}(V_0, V_1)$.
\end{defin}

Given an object $\sigma$ of $\square \Kur^{\ori}(V_0,V_1)$ we set
\begin{equation}
TX_{\sigma} = TX_{\sigma(1^n)} \quad\textrm{and}\quad V_{\sigma} =
V_{\sigma(1^n)}
\end{equation}
as before, and define 
\begin{equation}
W_{\sigma} = W_{\sigma(1^n)}, \, I^{\bC}_{\sigma} =
I^{\bC}_{\sigma(1^n)}, \quad\textrm{and}\quad O_\sigma =
O_{\sigma(1^n)}
\end{equation}
(In fact, the set $O$ is independent of the vertex of $\sigma$, since
we have restricted the domain category to a fixed stratification.)
Given a morphism $f \co \tau \to \sigma$, we define 
\begin{equation}
W_{f} = W_{\tau}, \, I^{\bC}_f = I^{\bC}_{\tau}, \quad\textrm{and}\quad
O_f = O_\tau.
\end{equation}

Combining Corollary~\ref{cor:comparison-tangent-space-maps-of-cubes}
with the definition of complex-oriented charts, we find that there is
a natural isomorphism  

\begin{equation}
V_0^- \oplus \bR^{O_f}  \oplus W \oplus TX_f \oplus V_1^+ \cong  V_0^+
\oplus I^{\bC}_f \oplus T \square_f \oplus W \oplus V_1^-.
\end{equation}

We adopt the notation from Equation
\eqref{eq:morphisms-categories-TX-I-and-bimodule}, and define 
for each morphism $f \in \square \Kur^{\ori}(V_0,V_1)$ the spectra
{\footnotesize
\begin{align}
C^{*,c}_{\rel \partial}(B \hat{Z}_f;  \scrT(\bX)) \equiv  & C^{*,c}_{\rel \partial}(B \hat{Z}_f;  F( S^{V}, (S^{TX})^{\mfib})
\wedge S^{V_{0}^{-} +  V^+_1 - (  V^+_0 + V_1^-) }  \wedge \bbk) \\
C^{*,c}_{\rel \partial}(B \hat{Z}_f;  \scrP(\bX))  \equiv  & C^{*,c}_{\rel \partial}(B \hat{Z}_f; F(S^{V^+_0  + O_f + W + V + V_1^-}, (S^{V^-_0 + O_f + T X + W + V^+_1}\wedge \bbk)^{\mfib}) ) \\
\cong & \nonumber C^{*,c}_{\rel \partial}(B \hat{Z}_f;  F(S^{ V^+_0  + O_f + W
  + V + V_1^- }, (S^{V^+_0 + I + T \square_f + W + V^-_1 } \wedge
\bbk)^{\mfib})) \\\label{eq:twisted_cochains-complex-bundle-cube}
C^{*,c}_{\rel \partial}(B \hat{Z}_f;    \scrI(\bX))  \equiv  & C^{*,c}_{\rel \partial}(B \hat{Z}_f; F(S^{V}, (S^{I} \wedge
\bbk)^{\mfib})  \wedge S^{ T \square_f - O_f }  ). 
\end{align}} 
Note, that in contrast to Section~\ref{sec:compl-lifts-tang}, the ring
spectrum $\bbk$ enters in the construction of these spectra.  It is
straightforward to check that these constructions define topological
functors with domain $\Tw  \square \Kur^{\ori}(V_0,V_1)$; the key
point is that they are all locally constant on the object space.
Summarizing, we have the following proposition:

\begin{prop} 
The assignments 
\begin{equation}
\begin{aligned}
f  \mapsto & C^{*,c}_{\rel \partial}(B \hat{Z}_f;  \scrT(\bX)) \\
f \mapsto & C^{*,c}_{\rel \partial}(B \hat{Z}_f; \scrP(\bX)) \\
f \mapsto & C^{*,c}_{\rel \partial}(B \hat{Z}_f; \scrI(\bX)) \\ 
\end{aligned}
\end{equation}
specify topological functors
\begin{equation}
 \Tw  \square \Kur^{\ori}(V_0,V_1) \to \Sp.
\end{equation} \qed
\end{prop}

The smash product induces comparisons between these constructions.

\begin{lem}
For each morphism $f \in \square \Kur^{\ori}(V_0,V_1)$, the smash
product with the identity on $S^{W + O_f} $ and on $S^{V^+_0 + W +
  V_1^- }$ defines natural equivalences of spectral functors:
\begin{equation}
C^{*,c}_{\rel \partial}(B \hat{Z}_f;  \scrT(\bX)) \to C^{*,c}_{\rel
  \partial}(B \hat{Z}_f;  \scrP(\bX)) \leftarrow C^{*,c}_{\rel
  \partial}(B \hat{Z}_f;    \scrI(\bX)). 
\end{equation} \qed
\end{lem}

We now consider the multiplicativity of this construction: recall that 
$\Kur^{\ori}$ is a topological bicategory. 

\begin{prop}\label{thm:bimodule-comparison-tangent}
There are natural maps 
\begin{align}
C^{*,c}_{\rel \partial}(B \hat{Z}_{f_1};    \scrT(\bX)) \wedge
C^{*,c}_{\rel \partial}(B \hat{Z}_{f_2};    \scrT(\bX)) & \to
C^{*,c}_{\rel \partial}(B \hat{Z}_{f_1 \times f_2};    \scrT(\bX)) \\ 
C^{*,c}_{\rel \partial}(B \hat{Z}_{f_1};    \scrI(\bX)) \wedge
C^{*,c}_{\rel \partial}(B \hat{Z}_{f_2};    \scrI(\bX)) & \to
C^{*,c}_{\rel \partial}(B \hat{Z}_{f_1 \times f_2};    \scrI(\bX)), 
\end{align}
which define strict functors
\begin{equation}
\Tw  \square \Kur^{\ori} \to \Sp.
\end{equation} \qed
\end{prop}

We can now state the functorial analogue of Lemma
\ref{lem:toy-case-bimodule-comparison-tangent}: 
\begin{lem}\label{lem:bimodule-comparison-tangent}
There are natural maps
\begin{align}  
C^{*,c}_{\rel \partial}( B \hat{Z}_{f_1};   \scrP(\bX) )  \wedge
C^{*,c}_{\rel \partial}( B \hat{Z}_{f_2};  \scrT(\bX))  & \to
C^{*,c}_{\rel \partial}( B \hat{Z}_{f_1 \times f_2};   \scrP(\bX)) \\ 
C^{*,c}_{\rel \partial}( B \hat{Z}_{f_1};    \scrI(\bX)   ) \wedge
C^{*,c}_{\rel \partial}( B \hat{Z}_{f_2};  \scrP(\bX)) & \to
C^{*,c}_{\rel \partial}( B \hat{Z}_{f_1 \times f_2}; \scrP(\bX) ), 
\end{align}
such that the following three diagrams commute:
\begin{equation}
\begin{tikzpicture}
\node[scale=.65]{
\begin{tikzcd}
C^{*,c}_{\rel \partial}( B \hat{Z}_{f_1};  \scrI(\bX) ) \wedge
C^{*,c}_{\rel \partial}( B \hat{Z}_{f_2}; \scrP(\bX)) \wedge
C^{*,c}_{\rel \partial}( B \hat{Z}_{f_3};  \scrT(\bX) )   \ar[d]
\ar[r]  & C^{*,c}_{\rel \partial}( B \hat{Z}_{f_1};  \scrI(\bX) )
\wedge  C^{*,c}_{\rel \partial}( B \hat{Z}_{f_2 \times f_3};
\scrP(\bX)) \ar[d]  \\ 
 C^{*,c}_{\rel \partial}( B \hat{Z}_{f_1 \times f_2}; \scrP(\bX))
 \wedge  C^{*,c}_{\rel \partial}( B \hat{Z}_{f_3};  \scrT(\bX) )
 \ar[r] &  C^{*,c}_{\rel \partial}( B \hat{Z}_{f_1 \times f_2 \times
   f_3}; \scrP(\bX)) 
  \end{tikzcd}
};
\end{tikzpicture}
\end{equation}
\begin{equation}
\begin{tikzpicture}
\node[scale=.65]{
\begin{tikzcd}
C^{*,c}_{\rel \partial}( B \hat{Z}_{f_1};  \scrP(\bX) ) \wedge
C^{*,c}_{\rel \partial}( B \hat{Z}_{f_2}; \scrT(\bX)) \wedge
C^{*,c}_{\rel \partial}( B \hat{Z}_{f_3};  \scrT(\bX) )   \ar[d]
\ar[r]  & C^{*,c}_{\rel \partial}( B \hat{Z}_{f_1};  \scrP(\bX) )
\wedge  C^{*,c}_{\rel \partial}( B \hat{Z}_{f_2 \times f_3};
\scrT(\bX)) \ar[d]  \\ 
 C^{*,c}_{\rel \partial}( B \hat{Z}_{f_1 \times f_2}; \scrP(\bX))
 \wedge  C^{*,c}_{\rel \partial}( B \hat{Z}_{f_3};  \scrT(\bX) )
 \ar[r] &  C^{*,c}_{\rel \partial}( B \hat{Z}_{f_1 \times f_2 \times
   f_3}; \scrP(\bX)) 
\end{tikzcd}
};
\end{tikzpicture}
\end{equation}
\begin{equation}
\begin{tikzpicture}
\node[scale=.65]{
\begin{tikzcd}
C^{*,c}_{\rel \partial}( B \hat{Z}_{f_1};  \scrI(\bX) ) \wedge
C^{*,c}_{\rel \partial}( B \hat{Z}_{f_2}; \scrI(\bX)) \wedge
C^{*,c}_{\rel \partial}( B \hat{Z}_{f_3};  \scrP(\bX) )   \ar[d]
\ar[r]  & C^{*,c}_{\rel \partial}( B \hat{Z}_{f_1};  \scrI(\bX) )
\wedge  C^{*,c}_{\rel \partial}( B \hat{Z}_{f_2 \times f_3};
\scrP(\bX)) \ar[d]  \\ 
C^{*,c}_{\rel \partial}( B \hat{Z}_{f_1 \times f_2}; \scrI(\bX))
\wedge  C^{*,c}_{\rel \partial}( B \hat{Z}_{f_3};  \scrP(\bX) )
\ar[r] &  C^{*,c}_{\rel \partial}( B \hat{Z}_{f_1 \times f_2 \times
  f_3}; \scrP(\bX)) 
\end{tikzcd}
};
\end{tikzpicture}
\end{equation}
\qed
\end{lem}

\begin{rem}
We note the slight differences in formulation from
Section~\ref{sec:mult-relat-orient}, which is due to the fact that we
incorporated in $\Kur^{\ori}$ a choice of finite set for each object,
with respect to which we destabilise. In the case of Kuranishi flow
categories, this set will be the singleton $q$ for a chart of the
space of morphisms from $p$ to $q$. When taking products in
$\Kur^{\ori} $ we take disjoint unions of finite sets, and if we start
with singletons, one of the two elements ends up corresponding to the
collar direction.
\end{rem}

\subsubsection{The category of stable complex vector bundles}
\label{sec:categ-stable-compl}

For clarity of exposition, we now forget most of the data of a complex
oriented Kuranishi chart by passing to the underlying zero-locus.  To
formalize this, we make the following definition.

\begin{defin}
Let $S\VB^{\bC}$ denote the topological category of \emph{charts of
  stratified stable complex orbibundles}.  
\begin{enumerate}
\item Objects consist of a stratified orbispace chart $(\cS, G, Z)$
  (recall Definition~\ref{defin:stratorbicharts}), a finite set $Q$, a
  $G$-equivariant complex bundle $I$ over $Z$, and a $G$-equivariant
  complex inner product space $V$.

\item A morphism $f$ is given by a map $(G',\cS', Z') \to (G,\cS, Z)$
  of the underlying orbispace charts, a choice of finite set $Q_f$, a  
bijection $Q' \cong Q \amalg Q_f$, an inclusion $V' \to V$ of $G'$
representations, equivariant isomorphisms $ I/I' \cong V/V' $ of
vector bundles over $Z'$, and complex-linear equivariant splittings of
the exact sequences 
\begin{equation}
\begin{aligned} 
  I' & \to I  \to I/I' \\
  V' & \to V  \to V/V'.
\end{aligned}
\end{equation}
The topology on the morphism spaces is induced by the topology on the
space of such splittings.  Composition is defined by composition of
splittings.
\end{enumerate}
\end{defin}

There is a natural monoidal structure on $S\VB^{\bC}$ given by the product of
underlying charts, the external direct sum of the corresponding vector
bundles, and the direct sum of splittings.

The category $S\VB^{\bC}$ is a coarsening of $\Kur^{\ori}$, in the
following sense:
\begin{lem}\label{lem:relative-oriented-chart-forget-to-complex-orbibundle}
The forgetful map 
\begin{equation}
\alpha \mapsto  (G_\alpha, \cS_\alpha, Q_\alpha, Z_\alpha,
I^{\bC}_\alpha , V_\alpha)
\end{equation}
defines a strict functor
\begin{equation}
\Kur^{\ori} \to  S\VB^{\bC}.
\end{equation}
\qed
\end{lem}

We now define the subcategory of cubes on $S\VB^{\bC}$ we work with.

\begin{defin}
\hspace{5 pt}

\begin{enumerate}
\item Let $\square S\VB^{\bC}$ denote the subcategory of cubes in
$S\VB^{\bC}$ with fixed stratifying partially ordered set $\cS$ and
labelling finite set $Q$.
\item Let $\Tw \square S\VB^{\bC}$ denote the subcategory of the
  twisted arrow category whose objects are pairs consisting of an
  arrow in $\square S\VB^{\bC}$ and a choice of totally ordered subset
  of the corresponding partially ordered set $\cS$.
\end{enumerate}
\end{defin}

As before, we associate to each morphism $f \in \square S\VB^{\bC}$ a
cubical degeneration $Z_f$ which is an open subset of $Z_{f(1^n)}
\times \square_f$.  We now consider the lax monoidal functor
\begin{equation}
\Tw \square S\VB^{\bC} \to \Sp
\end{equation}
defined on objects by the assignment
\begin{equation} \label{eq:complex-twisted-cochains-functor}
f \mapsto C^{*,c}_{\rel \partial}(B \hat{Z}_f; F(S^{V}, (S^{I} \wedge
\bbk)^{\mfib})  \wedge S^{ T \square_f - Q_f }  ).
\end{equation}

Comparing with Equation
\eqref{eq:twisted_cochains-complex-bundle-cube}, we conclude: 
\begin{lem}\label{lem:reduce-to-SVB}
The lax functor $C^{*,c}_{\rel \partial}(B \hat{Z}_f; \scrI(\bX))$
from $\Tw \square \Kur^{\ori} \to \Sp$ factors through $\Tw \square
S\VB^{\bC}$.  
\qed 
\end{lem}

\begin{rem}
It is straightforward to construct a lax functor $\Tw \Kur^{\ori} \to
\Sp$ that factors through $\Tw S\VB^{\bC}$ without passing to the
category of cubes, but we shall require cubes to be able to formulate
the functoriality of choices of classifying maps for the given complex
vector bundles, which we do in the next section.
\end{rem}

\subsubsection{Functorial and multiplicative classifying maps}
\label{sec:funct-mult-triv}

We now impose the assumption that $\bbk$ is a multiplicative complex
oriented spectrum.  Given an object $f$ of $\Tw \square S\VB^{\bC}
$, we consider the vector bundle $I^{\bC}_f $ and the vector space $
V_f $ as complex vector bundles over $B Z_f$. We shall apply the
constructions of Section~\ref{sec:mult-relat-orient} to trivialize the
corresponding spherical fibrations, but to do this coherently requires
suitable choices of classifying maps.  

\begin{defin}
Classifying data for an object $f$ of $\Tw \square S\VB^{\bC}$ is a
collection of classifying maps for the complex bundles $I^{\bC}_f$,
$V_f$, and $V_f^\perp$, by which we mean maps to $BU(W)$ (for
appropriate $W$) along with specified isomorphisms from the pullback
of the universal bundle. 
\end{defin}

The key fact that allows us to choose the classifying maps we need is
the following standard lemma.

\begin{lem}\label{lem:contractible-classifying-for-bundle}
Let $f \colon E \to B$ be a complex vector bundle with fiber $V$.  Let
$\xi_V$ denote the universal bundle over $BU(V)$.  The space of pairs
$(\tilde{f}, \gamma)$, where $\tilde{f} \colon B \to BU(V)$ is a
classifying map for $f$ and $\gamma$ is a choice of isomorphism $f
\cong \tilde{f}^* \xi_V$ of complex bundles, is nonempty and
contractible.
\end{lem}

\begin{example}
When $B = \ast$, there is a natural map from the space of pairs
$(\tilde{f}, \gamma)$ to the space $\Map(\ast, BU(V))$; this is a
model of the universal bundle $EU(V) \to BU(V)$ with fiber $U(V)$.
\end{example}

In order to state the desired properties that such classifying maps
should satisfy, we consider a composition $f \circ g \co \rho \to \tau
\to \sigma$ of cubes, and note that the inclusion $V_\tau \to V_\sigma
$ implies that $V_\tau$ has a natural action of $G_\sigma$, hence that
the action of $G_g$ on $V_g^\perp$ is pulled back from an action of
$G_{f \circ g}$. In particular, thinking of $V_g^\perp $ as a vector
bundle over $B Z_g$, we find that it is naturally isomorphic to the
pullback of a bundle on $B Z_{f \circ g} $.

\begin{defin}
Let $f$ and $g$ be composable maps in $\square S\VB^{\bC}$.  A triple
of classifying maps for $(I^{\bC}_f, V_f, V_f^\perp)$, $(I^{\bC}_{f
  \circ g}, V_{f \circ g}, V_{f \circ g}^\perp)$, and $(I^{\bC}_g,
V_g, V_g^\perp)$ are \emph{compatible} if
\begin{enumerate}
\item The classifying map for $V_g^\perp$ factors through $B Z_{f \circ g}$.
\item The classifying maps for $ I^{\bC}_{f \circ g}$ and $ V_{f \circ
  g} $ agree with those for $I^{\bC}_g$ and $V_g$ under pullback with
  respect to the projection $Z_g \to Z_{f \circ g}$. 
\item The restriction of the classifying maps for $I^{\bC}_f$ and
  $V_f$ to $Z_{f \circ g}$ are given by the direct sum of the
  classifying maps for $ I^{\bC}_{f \circ g}$ and $ V_{f \circ g} $
  with the classifying map for $V_g^\perp $. 
\item The classifying map for $V_{f \circ g}^\perp$ is the direct sum
  of the classifying maps for $V_f^\perp$ and $V_g^\perp$. 
\end{enumerate}
\end{defin}

While the space of choice of classifying maps associated to each arrow
is contractible, it is not clear how to make sure that the choices are
compatible: the natural way to construct such choices is by induction
on the dimension of the target of each arrow, but the notion of
compatibility makes such an inductive scheme impossible because the
choice for $g$ is constrained by the one for $f \circ g$. This leads
us to consider the following category.

\begin{defin}
Let  $S\VB^{\bC}_{sub}\langle \cS \rangle$ be the category with:
\begin{enumerate}
\item Objects the pairs $(\alpha, \overline{G}_\alpha)$ with 
$\alpha$ an object of $ S\VB^{\bC}$ stratified by $\cS$,  and
$\overline{G}_\alpha$ a subgroup of $G_\alpha$ acting freely on
$Z_\alpha$ and trivially on $V_\alpha$. 
\item A morphism consists of a map
$f \co \alpha \to \beta$ in $ S\VB^{\bC}$, such that  
\begin{equation} \label{eq:map_subgroups_inclusion}
 \overline{G}_\beta \subset f( \overline{G}_\alpha).  
\end{equation}
(This construction is evidently closed under compositions.)
\end{enumerate}
\end{defin}

We assemble these categories into the category $\square
S\VB^{\bC}_{sub}$, with objects the union of the objects of
$S\VB^{\bC}_{sub}\langle \cS \rangle $ for all ordered sets $\cS$, and
morphisms given as in
Definition~\ref{def:stratified-Kuranishi-charts}.  The fiber category
of the forgetful functor $\square S\VB^{\bC}_{sub} \to \square
S\VB^{\bC}$ at each object is contractible, as it has a final
object given by taking the distinguished subgroup to be $\{e\}$.
Thus, Quillen's theorem A implies we can work with $\square
S\VB^{\bC}_{sub}$ in place of $\square S\VB^{\bC}$:

\begin{lem}
The forgetful functor
\begin{equation}
\square S\VB^{\bC}_{sub} \to \square S\VB^{\bC}
\end{equation}
is homotopy cofinal. \qed
\end{lem}

For each cube $\sigma$ in $S\VB^{\bC}_{sub}\langle \cS \rangle$, we
define $\overline{G}_{\sigma}$ to be the group associated to
$\sigma(1^n)$. Given a map $f \co \tau \to \sigma$ of cubes, we set
$\overline{G}_{f} \equiv \overline{G}_{\sigma}$. We define
$\underline{G}_f$ to be the quotient $ G_f/ \overline{G}_f  $. Since $
\overline{G}_f $ acts freely on $Z_f$, and $I^{\bC}_f$ is a
$G_f$-equivariant bundle, we conclude that it is pulled back from a $
\overline{G}_f$ equivariant bundle on the quotient $
Z_f/\overline{G}_f  $, and hence that the corresponding bundle over
$BZ_f$ is obtained by pull back under the projection map 
\begin{equation}
B Z_f \to E  \underline{G}_f \times_{\underline{G}_f } Z_f/\overline{G}_f .
\end{equation}
At the same time, the vector spaces $V_f$ and $V_f^\perp$ are $G_f$
representations which are pulled back under the surjection to
$\underline{G}_f$, and hence the corresponding bundle over $BZ_f$ is
pulled back from $B \underline{G}_f $. Given the commutative diagram 
\begin{equation}
\begin{tikzcd}
B Z_f \ar[r] \ar[d] &  E  \underline{G}_f \times_{\underline{G}_f }
Z_f/\overline{G}_f \ar[d] \\ 
B G_f \ar[r] & B \underline{G}_f
\end{tikzcd}
\end{equation}
we conclude:
\begin{lem}
The complex vector bundles $I^{\bC}_f$, $V_f$ and $V_f^\perp$ over $B
Z_f $ are naturally isomorphic to the pullback of vector bundles over
$E  \underline{G}_f \times_{\underline{G}_f }
Z_f/\overline{G}_f$. \qed 
\end{lem} 
Note that a composition $f \circ g$ induces a commutative diagram 
\begin{equation}
  \begin{tikzcd}[column sep=10]
B Z_{g} \ar[r] \ar[d] & B Z_{f \circ g} \ar[d] \ar[r] & B Z_{f} \ar[d]  \\
E  \underline{G}_g \times_{\underline{G}_g } Z_g/\overline{G}_g &
\ar[l] E  \underline{G}_{f \circ g} \times_{\underline{G}_{f \circ g}
} Z_{f \circ g}/\overline{G}_{f \circ g} \ar[r] & E  \underline{G}_f
\times_{\underline{G}_f } Z_f/\overline{G}_f, 
\end{tikzcd}
\end{equation}
where the right arrow on the bottom uses the fact that $
\overline{G}_{f \circ g} = \overline{G}_f $, and the left arrow uses
the surjection 
\begin{equation}
\underline{G}_{f \circ g}  \to   \underline{G}_g 
\end{equation}
induced by Equation \eqref{eq:map_subgroups_inclusion}.  
This lets us inductively choose compatible classifying maps, as
follows.

\begin{lem}\label{lem:classifying-compatible-composition}
The space of compatible choices of classifying maps for all arrows $f$
in $\square S\VB^{\bC}_{sub}\langle \cS \rangle$, with the property
that they are pulled back from $E \underline{G}_f
\times_{\underline{G}_f } Z_f/\overline{G}_f$, is nonempty and
contractible.
\end{lem}

\begin{proof}
We proceed by induction on the dimension of the target of an arrow.  For
$0$-dimensional targets, the choice is clearly unconstrained.  Now,
given a choice of compatible classifying maps for all arrows whose
target has dimension strictly less than $n$, the compatibility
conditions dictate the choices of classifying maps for $V_f$ and
$I^{\bC}_f$ whenever the domain of $f$ has dimension strictly less
than $n$ and the target has dimension $n$.  Thus, choosing classifying
maps for $V_f^\perp$ by induction on the codimension of $f$ determines
the choices of classifying maps for all maps with $n$-dimensional
target and domain of dimension strictly less than $n$.  Finally, we
can extend this choice to $n$-dimensional cubes, completing the
inductive step.  

This induction shows the space of compatible classifying maps is
nonempty.  To see that it is contractible, observe that the same
inductive procedure allows us to construct a homotopy to any
particular point in this space.
\end{proof}

Since the bundle $I^{\bC}_{\partial^Q f}$ and the representation
$V_{\partial^Q f}$ associated to a stratum labelled by a totally
ordered subset $Q$ of $\cS$ are defined to agree with the pair
$(I^{\bC}_f,V_f)$, it is straightforward to extend this construction
to the category $\square S\VB^{\bC}_{sub}$, so that classifying maps
are compatible with restriction to boundary strata. 

\begin{prop}\label{prop:classifying-compatible}
We can choose compatible classifying maps for the category $\square
S\VB^{\bC}_{sub}$. \qed
\end{prop}

We can also extend the choice of classifying data to $\Tw \square
S\VB^{\bC}_{sub}$, and for the remainder of the section we will work
with the twisted arrow category.  We now describe how to arrange for
the choice of classifying maps to be compatible with the external
product.

\begin{defin}
A choice of classifying maps is multiplicative if for each product
$f_1 \times f_2$ of charts, the classifying maps for $(I^{\bC}_{f_1
  \times f_2}, V_{f_1 \times f_2})$ coincide with the products of the
classifying maps for $(I^{\bC}_{f_1}, V_{f_1})$ and $(I^{\bC}_{f_2},
V_{f_2})$.  
\end{defin}

Choosing multiplicative classifying maps would involve keeping track 
of product decompositions of charts, so we instead consider the
following categories which encode the space of choices.  First, we
have the following pointwise definition, which is justified by
Proposition~\ref{prop:classifying-compatible}.

\begin{defin}\label{defin:add-classifying-data}
Let $\Tw \square S\VB^{\bC}_{sub, class}$ be the internal category in
spaces where the objects are given by the objects of $\Tw \square
S\VB^{\bC}_{sub}$ along with compatible classifying data.  The
topology on the objects is given by the topology on the space of
classifying data.  Morphisms are specified by maps in $\Tw \square
S\VB^{\bC}_{sub}$ that are compatible with the classifying data.
\end{defin}

By Lemma~\ref{lem:contractible-classifying-for-bundle}, the projection 
\begin{equation}\label{eq:classifying-data}
\Tw \square S\VB^{\bC}_{sub,class} \to \Tw \square
S\VB^{\bC}_{sub}
\end{equation}
induces a DK-equivalence of internal categories.  Moreover, as we explain
in more detail below, the topologized classifying map information can
be pulled back to $\Tw \square \Kur^{\ori}$.

The point of introducing $\Tw \square S\VB^{\bC}_{sub, class}$ is that
taking the product of classifying maps induces a continuous map from
the product of spaces of classifying data to classifying data on the
product.

\begin{prop}
The symmetric monoidal structure on $\Tw \square S\VB^{\bC}_{sub}$
coupled with the product of classifying maps induces a symmetric
monoidal structure on the topological category $\Tw \square
S\VB^{\bC}_{sub, class}$. \qed
\end{prop}

\subsubsection{Trivializing twisted cochains}
\label{sec:triv-twist-coch}

Now that we have produced compatible multiplicative product maps, we
can apply the results of Section~\ref{sec:multiplicative-orientations}
to compatibly trivialize the bundles.

\begin{defin}
Let $\wtOrb$ denote the topological category with  
\begin{enumerate}
\item objects consisting of a $\cS$-stratified orbispace chart $(G,Z,
  \cS)$, a finite set 
$Q$, a pair of complex inner-product spaces $W$ and $V$. 
\item morphisms $f$ from $(G',\cS', Z') \to (G,\cS, Z)$ given by a map
  of the underlying orbispace charts, a choice of finite set $Q_f$, a
  bijection $Q' \cong Q \amalg Q_f$, isometries $V' \to V$ and $W' \to
  W$, and a complex isomorphism $(V')^{\perp} \cong (W')^{\perp}$.
\end{enumerate}
\end{defin}

This is a monoidal category with monoidal structure induced by the
product of Kuranishi charts and direct sum of vector spaces.

\begin{lem}
We have a lax monoidal topological functor $\Tw \square
\wtOrb \to \Sp$ given on objects by 
\begin{equation} \label{eq:functor_untwisted_cohains-pair-spheres}
f \mapsto C^{*,c}_{\rel \partial}(B\hat{Z}_f; F(S^V,
S^W \sma \bbk)^{\mfib})  \wedge S^{T\square_f -Q_f} ).   
\end{equation} \qed
\end{lem}

There is a natural topological functor
\begin{equation}
\Tw \square S \VB^{\bC}_{sub, class} \to \Tw \square \wtOrb 
\end{equation}
induced from the projection $S \VB^{\bC}_{sub,class} \to
\wtOrb$ that passes to the fiber of the bundle $I$, and so
we obtain a topological functor   
\begin{equation}\label{eq:untwisted-cochains-pullback}
\begin{aligned}
\Tw \square S \VB^{\bC}_{sub,class} \to \Sp
\end{aligned}
\end{equation}
by pulling back
Equation~\eqref{eq:functor_untwisted_cohains-pair-spheres}.  We begin
by  comparing this to the complex twisted cochains functor of
Equation~\eqref{eq:complex-twisted-cochains-functor}, using the
complex Thom isomorphism (see
Section~\ref{sec:multiplicative-orientations}) and the fact that
$\bbk$ is complex-oriented.

\begin{prop}\label{prop:trivialize}
There is a natural zig-zag of lax monoidal equivalences connecting the
complex twisted cochains functor 
\begin{equation}
C^{*,c}_{\rel \partial}(B \hat{Z}_f; F(S^{V},
(S^{I} \wedge \bbk)^{\mfib}) \sma S^{ T \square_f - Q_f })   
\end{equation}
and the pullback functor in
Equation~\eqref{eq:untwisted-cochains-pullback}.
\end{prop}

\begin{proof}
Given an object of $\Tw \square S \VB^{\bC}_{sub}$, consider the
spectrum 
\begin{equation} \label{eq:intermediate-functor-complex-and-trivial}
f \mapsto C^{*,c}_{\rel \partial}(B \hat{Z}_f; F(S^{V}, (S^{\tilde{I}}
\wedge \bbk)^{\mfib}) \sma S^{T \square_f - Q_f} 
), 
\end{equation}
where $\tilde{I}$ denotes the complex inner-product space that is the
fiber of the vector bundle $I$.  This assignment specifies a lax
monoidal topological functor from $\Tw \square S \VB^{\bC}_{sub}$ to
spectra, and it is straightforward to check that the functor described
in Equation~\eqref{eq:intermediate-functor-complex-and-trivial} is
naturally equivalent to the pullback functor described in
Equation~\eqref{eq:untwisted-cochains-pullback}.

Therefore, it suffices to compare
Equation~\eqref{eq:intermediate-functor-complex-and-trivial} to the
complex twisted cochains.  This follows using the equivariant Thom
isomorphism of Theorem~\ref{thm:thom-iso}.  In particular,
Proposition~\ref{prop:thom-direct} implies that for a bundle $I$ with
base $Z$ and fiber $\tilde{I}$, there is a natural zig-zag
\begin{equation}
TI \sma \bbk \to \Sigma^{\tilde{I}}_+ Z \sma \bbk,
\end{equation}
where $TI$ denotes the Thom space of $I$, which is compatible with the
external multiplication.  Interpreted in terms of the corresponding
parametrized spectra, the Thom zig-zag yields a natural multiplicative
comparison
\begin{equation}
S^I \sma \bbk \htp S^{\tilde{I}}_{Z} \sma \bbk.
\end{equation}
Moreover, the Thom isomorphism is realized as a composite consisting
of enriched functors of classifying maps, and is in particular
continuous for the topology on the classifying data.
\end{proof}

\subsubsection{Comparison with untwisted cochains}
\label{sec:comp-with-untw}

The last step in the comparison between tangentially twisted and
ordinary cochains is to consolidate the sphere coordinates.

\begin{defin}
Let $\wtOrb^{disc}$ denote the category with 
\begin{enumerate}
\item objects $(\cS,
Z, G, Q, d_I, d_V)$ consisting of an $\cS$-stratified orbispace chart
$(Z,G)$, a finite set $Q$, and a pair $(d_I,d_V)$ of integers, and
\item morphisms specified by a morphism of stratified orbispace charts
  together with a non-negative integer $k$ such that the pairs of
  integers are related by
\begin{equation}\label{eq:dimension-shift}
  (d_I, d_V) \to (d_I +k, d_V +k).
\end{equation}
\end{enumerate}
This is a monoidal category with monoidal structure induced by the
product of Kuranishi charts and addition of integers. 
\end{defin}

The category $\wtOrb^{disc}$ discretizes the maps between
spheres in $\wtOrb$.  There is an evident lax monoidal
topological functor 
\begin{equation}\label{eq:orb-discretization}
\wtOrb \to \wtOrb^{disc}
\end{equation}
that takes $W$ and $V$ to their dimension and any morphism to the
unique corresponding map in $\wtOrb^{disc}$, using the fact
that a morphism exists in $\wtOrb$ only when
Equation~\eqref{eq:dimension-shift} holds.  Note that this projection
functor is not a DK-equivalence because the spaces of isometries
between inner product spaces are not contractible unless they are
$0$-dimensional.

\begin{lem}
The assignment for an arrow $f$ in $\square \wtOrb$ 
\begin{equation} \label{eq:last-untwisting}
f \mapsto C^{*,c}_{\rel \partial}(B\hat{Z}_f; S^{T\square_f -Q_f}
\wedge \bbk )[d_I - d_V]
\end{equation}
specifies a lax monoidal functor
\begin{equation}
\Tw \square \wtOrb^{disc} \to \Sp
\end{equation}
\qed
\end{lem}

We would now like to compare the pullback of
Equation~\eqref{eq:last-untwisting} along the projection
\begin{equation}
\Tw \square \wtOrb \to \Tw \square \wtOrb^{\disc},
\end{equation}
to the functor
\begin{equation}
f \mapsto C^{*,c}_{\rel \partial}(B\hat{Z}_f; F(S^V, S^W \sma
\bbk)^{\mfib})  \wedge S^{T\square_f -Q_f} ).    
\end{equation}
of Equation~\eqref{eq:functor_untwisted_cohains-pair-spheres}.

For $\bbk = \bS$, there would be no reason to expect a comparison of
this form, since the objects of $\wtOrb$ have nontrivial
automorphisms coming from unitary groups.  The key point, as discussed
in Lemma~\ref{lem:complex-stability-equivariant}, is that the zig-zag
representing the complex Thom isomorphism 
\begin{equation}
\begin{tikzcd}
MUP_G \ar[r] & \Sh_V MUP_G & \ar[l] S^V \sma MUP_G
\end{tikzcd}
\end{equation}
is $U(V)$-equivariant where we give $MUP_G$ the trivial $U(V)$-action
(and here $\Sh_V$ is the $V$-shift functor from
Definition~\ref{defin:V-shift}).

Given the trivialization of the unitary group actions, pointwise we
are looking at the equivalence
\begin{equation}
F(S^V, S^W) \htp \bS[|W| - |V|].
\end{equation}
To produce a model of this equivalence which is compatible with the
functoriality and monoidal structure on $\wtOrb^{disc}$, we
proceed as in Section~\ref{sec:compl-stable-struct} by constructing a
bimodule representing an equivalence.  Specifically, we consider the
functor $\Tw \square \wtOrb \to \Sp$ specified by the
assignment
\begin{equation}
f \mapsto C^{*,c}_{\rel \partial}(B\hat{Z}_f, F(S^{U},\bS[|V|])) \sma
S^{T\square_f -Q_f}.
\end{equation}

The argument in Section~\ref{sec:rigidifying-spheres} then proves the
following proposition, which establish there is a strictly
$\Pi$-equivariant bimodule structure on passage to homotopy colimits.

\begin{prop}\label{prop:bimodules-sphere-comparison}
There are natural associative and unital maps
\begin{equation}
\begin{tikzpicture}
\node[scale=.75]{
\begin{tikzcd}
C^{*,c}_{\rel \partial}(B\hat{Z}_f; F_{\bbk}(S^V \sma \bbk,
S^W \sma \bbk)^{\mfib})  \wedge S^{T\square_f -Q_f} ) \sma
C^{*,c}_{\rel \partial}(B\hat{Z}_f; F(S^{U},\bS[|V|])) \sma  
S^{T\square_f -Q_f} \ar[d] \\
 C^{*,c}_{\rel \partial}(B\hat{Z}_f; F(S^{U},\bS[|V|])) \sma
 S^{T\square_f -Q_f} 
\end{tikzcd}
};
\end{tikzpicture}
\end{equation}
\begin{equation}
  \begin{tikzpicture}
\node[scale=.75]{
\begin{tikzcd}
C^{*,c}_{\rel \partial}(B\hat{Z}_f; F(S^{U},\bS[|V|])) \sma
S^{T\square_f -Q_f} \sma C^{*,c}_{\rel \partial}(B\hat{Z}_f; S^{T\square_f -Q_f}
\wedge \bbk )[d_I - d_V] 
\ar[d] \\
C^{*,c}_{\rel \partial}(B\hat{Z}_f; F(S^{U},\bS[|V|])) \sma S^{T\square_f -Q_f}
\end{tikzcd}
};
\end{tikzpicture}
\end{equation} \qed
\end{prop}

Specifically, we get $\Pi$-equivariant comparisons of the associated
spectral categories with morphism spectra computed via $\hocolim_{\Tw
  \square A^{\scrT}}$.

\subsubsection{Complex-oriented Kuranishi presentations}
\label{sec:compl-orient-kuran}

Finally, we consider a $\Pi$-equivariant $2$-category $A^\ori$
equipped with a strictly $\Pi$-equivariant $2$-functor 
\begin{equation}
A^\ori \to  \Kur^{\ori}
\end{equation}
which lifts a Kuranishi presentation of a flow category over a
partially ordered set $\cP$, forming a \emph{complex-oriented
  Kuranishi presentation} as in
Definition~\ref{def:complex-oriented-flow-category}.  First, observe
that we can lift the classifying data attached to each object of $S
\VB^{\bC}_{sub,class}$ to $\Kur^{\ori}$ so that it is preserved under
the projection functor; this simply involves augmenting the
objects and morphisms of $\Kur^{\ori}$ in analogy with
Definition~\ref{defin:add-classifying-data}.  We denote the resulting
topological bicategory by $\Kur^{\ori, class}$.  We then have a 
$\Pi$-equivariant lift of the given complex-oriented Kuranishi
presentation to a strictly $\Pi$-equivariant topological $2$-functor 
\begin{equation}\label{eq:class-lift-flow-category}
A^{\ori,class} \to \Kur^{\ori, class} 
\end{equation} 
that sits in the commutative diagram
\begin{equation}
\begin{tikzcd}
A^{\ori,class} \ar[r] \ar[d] & \ar[d] \Kur^{\ori, class} \\
A^\ori \ar[r] & \Kur^{\ori}.
\end{tikzcd}
\end{equation}

Since the classifying data is contractible, passing to this lift does
not change the realization:

\begin{lem}
The projection $2$-functor $A^{\ori,class} \to A^{\ori}$ induces
homotopy cofinal functors $A^{\ori,class}(p,q) \to A^{\ori}(p,q)$ for
each pair $p,q \in \cP$. \qed
\end{lem}

We obtain from the complex-oriented Kuranishi presentation a
$\Pi$-equivariant assignment 
\begin{equation}
  p \to V_{p} = (V^+_{p}, V^-_{p})
\end{equation}
of a stable vector space to every element of $\cP$. The relative
isomorphism from
Equation~\eqref{eq:stable-isomorphism-tangent-with-complex} 
between $TX_\alpha$ and $I^\bC_\alpha$ implies that 
\begin{multline}
1 + \dim X_\alpha - \dim V_\alpha -  \dim V^+_{p} + \dim V^+_{q} = \\ 2
(\dim_{\bC} I^\bC_\alpha - \dim_{\bC} V_\alpha -  \dim_{\bC} V^-_{p} +
\dim_{\bC} V^-_{q}). 
\end{multline}
Applying Equation \eqref{eq:dimension-moduli-space-degree-difference},
we conclude that for each pair $(p,q)$ of objects of $\cP$, we have 
\begin{equation}\label{eq:shift-comparison}
\dim_{\bC} I^\bC_\alpha - \dim_{\bC} V_\alpha= \deg q - \deg p. 
\end{equation}

The following proposition is the main result of this section, which
completes the comparison zig-zag connecting $C^{*,c}_{\rel
  \partial}(B \hat{\cZ}; S^{T \hat{\cX}-V-d})$ to $C^{*,c}_{\rel
  \partial} \left( B\hat{\cZ}, \Omega \bbk \right)$, as defined in
Definition~\ref{defin:compactly-supported-relative-cochains}.

\begin{prop}
Given a complex oriented Kuranishi presentation, there is a
$\Pi$-equivariant zig-zag of DK-equivalences of spectrally enriched
categories 
\begin{equation}
\begin{tikzcd}[column sep = 10]
C^{*,c}_{\rel \partial} \left(  B\hat{\cZ};  S^{T \hat{\cX} -V
       -d} \wedge \bbk  \right)  \ar[r,dashed] &  C^{*,c}_{\rel
       \partial} \left(  B\hat{\cZ};  S^{\ell + \cI -V -d} \wedge \bbk
     \right)  & \ar[l,dashed] C^{*,c}_{\rel \partial} \left(B\hat{\cZ};  \Omega \bbk  \right).
\end{tikzcd}
\end{equation}
\end{prop}

\begin{proof}
First, Proposition~\ref{thm:bimodule-comparison-tangent} and
Lemma~\ref{lem:bimodule-comparison-tangent} establish an equivalence
of $\Pi$-equivariant spectral categories given by an invertible
bimodule 
\begin{multline}
C^{*,c}_{\rel \partial} \left( B\hat{\cZ}; S^{T \hat{\cX} -V -d}
\wedge \bbk \right) \htp \\ \hocolim_{\Tw \square A^{\scrT}}
C^{*,c}_{\rel \partial}(B \hat{Z}_f; F(S^{V}, (S^{I} \wedge
\bbk)^{\mfib}) \wedge S^{ T \square_f - O_f }),
\end{multline}
where here the notation on the righthand side denotes the spectral
category with morphisms from $p$ to $q$ given by the homotopy colimits
over $\Tw \Box A^{\scrT}(p,q)$.  Lemma~\ref{lem:reduce-to-SVB} now
shows that the natural map
\begin{equation}\label{eq:restrict-hocolim}
\begin{tikzcd}
\displaystyle{\hocolim_{\Tw \square A^{\scrT}}}
C^{*,c}_{\rel \partial}(B \hat{Z}_f; F(S^{V}, (S^{I} \wedge
\bbk)^{\mfib})  \wedge S^{ T \square_f - O_f })
\ar[d] \\
\displaystyle{\hocolim_{\Tw \square S\VB^{\bC}_{sub}}}
C^{*,c}_{\rel \partial}(B \hat{Z}_f; F(S^{V}, (S^{I} \wedge
\bbk)^{\mfib})  \wedge S^{ T \square_f - O_f })
\end{tikzcd}
\end{equation}
is a DK-equivalence of $\Pi$-equivariant spectral categories, where
the notation on the bottom abusively denotes the spectral category
with morphisms given by the homotopy colimit over the pullback of
diagrams along the projection $\Tw \square A^{\scrT} \to \Tw \square
S\VB^{\bC}_{sub}$.  Using the DK-equivalence of
Equation~\eqref{eq:classifying-data}, we can introduce classifying
data by considering $\Tw \square S\VB^{\bC}_{sub,class}$ and the lift
of the flow category from Equation~\eqref{eq:class-lift-flow-category}
to work with the composite
\begin{equation}
A^{\ori,class} \to \Kur^{\ori, class} \to S\VB^{\bC}_{sub,class}.
\end{equation}
Specifically, we have an equivalence
\begin{equation}
\begin{tikzcd}
\displaystyle{\hocolim_{\Tw \square A^{\scrT,class}}} C^{*,c}_{\rel
  \partial}(B \hat{Z}_f; F(S^{V}, (S^{I} \wedge \bbk)^{\mfib}) \wedge
S^{ T \square_f - O_f }) \ar[d] \\ \displaystyle{\hocolim_{\Tw \square
    S\VB^{\bC}_{sub,class}}} C^{*,c}_{\rel \partial}(B \hat{Z}_f;
F(S^{V}, (S^{I} \wedge \bbk)^{\mfib}) \wedge S^{ T \square_f - O_f })
\end{tikzcd}
\end{equation}
that lifts the map of Equation~\eqref{eq:restrict-hocolim} in the
sense that the square
\begin{equation}
\begin{tikzcd}
\displaystyle{\hocolim_{\Tw \square A^{\scrT,class}}} C^{*,c}_{\rel
  \partial}(-) \ar[d] \ar[r] & \displaystyle{ \hocolim_{\Tw \square
    S\VB^{\bC}_{sub, class}}} C^{*,c}_{\rel \partial}(-) \ar[d]
\\ \hocolim_{\Tw \square A^{\scrT}} C^{*,c}_{\rel \partial}(-) \ar[r]
& \hocolim_{\Tw \square S\VB^{\bC}_{sub}} C^{*,c}_{\rel \partial}(-)
\end{tikzcd}
\end{equation}
commutes, is $\Pi$-equivariant, and the vertical maps are
DK-equivalences.

Next, Proposition~\ref{prop:trivialize} trivializes the complex
bundles to produce a $\Pi$-equivariant equivalence
\begin{equation}
\begin{tikzcd}
\displaystyle{\hocolim_{\Tw \square S\VB^{\bC}_{sub,class}}}
C^{*,c}_{\rel \partial}(B \hat{Z}_f; F(S^{V}, (S^{I} \wedge
\bbk)^{\mfib})  \wedge S^{ T \square_f - O_f })
\ar[d] \\
\displaystyle{\hocolim_{\Tw \square S\VB^{\bC}_{sub,class}}}
C^{*,c}_{\rel \partial}(B\hat{Z}_f; F_{\bbk}(S^V \sma \bbk, (S^W \sma
\bbk)^{\mfib})  \wedge S^{T\square_f  -Q_f}).
\end{tikzcd}
\end{equation}
Finally, Proposition~\ref{prop:bimodules-sphere-comparison} assembles
into a $\Pi$-equivariant bimodule equivalence 
\begin{equation}
\begin{tikzcd}
\displaystyle{\hocolim_{\Tw \square S\VB^{\bC}_{sub,class}}}
C^{*,c}_{\rel \partial}(B\hat{Z}_f, F_{\bbk}(S^V \sma \bbk; S^W \sma
\bbk)^{\mfib})  \wedge S^{T\square_f  -Q_f}
\ar[dash,"\htp"]{d} \\
\displaystyle{\hocolim_{\Tw \square S\VB^{\bC}_{sub,class}}} C^{*,c}_{\rel
    \partial}(B\hat{Z}_f; S^{T\square_f -Q_f} \wedge \bbk )[d_I - d_V]
\end{tikzcd}
\end{equation}
combining the sphere coordinates, and applying
Equation~\eqref{eq:shift-comparison} and using the multiplicative
properties of the system of spheres $\{\bS[-]\}$ now completes the
argument.
\end{proof}

\part{Floer theoretic constructions}
\label{part:floer-theor-constr}

\section{Hamiltonian Floer theory and Morse theory}
\label{sec:hamilt-floer-theory-1}

Let $M$ be a closed symplectic manifold of real dimension $2n$. In this section, we construct the moduli spaces of Floer trajectories of a non-degenerate Hamiltonian function $H$ on $M$, and prepare the ingredients needed to construct a Kuranishi presentation of such spaces. As noted earlier, our approach is closest to Pardon's construction from \cite[Section 10]{Pardon2016}, and can be compared with the work of McDuff-Wehrheim \cite{McDuff2017}. All these points of view are variants of the Kuranishi construction of Fukaya-Ono \cite{FukayaOno1999} and Fukaya-Oh-Ohta-Ono \cite{FukayaOhOhtaOno2009}.

Since our applications rely on being able to (partially) compute the homotopy type associated to $H$, we shall simultaneously build a Kuranishi presentation of a larger category with four types of objects, three of which correspond to critical points of a fixed Morse functions, and one to Hamiltonian orbits of $H$. We then implement a variant of the comparison between Floer and Morse theory, via fibre products of moduli spaces with marked points (see  \cite{Fukaya1997b,PiunikhinSalamonSchwarz1996}).

\subsection{Stable maps and buildings}
\label{sec:stable-maps-build}

\begin{defin} \label{def:pre-stable}
A \emph{pre-stable cylinder} consists of the following data
  \begin{enumerate}
  \item A finite tree $T$ equipped with a distinguished root vertex, 
  \item For each vertex $v \in T$ other than the root, a genus $0$ closed Riemann surface $\Sigma_v$. We define $\Sigma_v = \bR \times S^1$ if $v$ is the root, which we can identify with the complement of the points $z_- = 0$ and $z_+ = \infty$ in $\bP^1$, 
  \item (Nodal points) For each flag $(v,e)$ (i.e. a pair consisting of an endpoint $v$ of an edge $e$ in $T$), a marked point $z_{v,e} \in \Sigma_v$.
  \end{enumerate}
  A \emph{pre-stable building} is a finite ordered collection of pre-stable cylinders.
\end{defin}
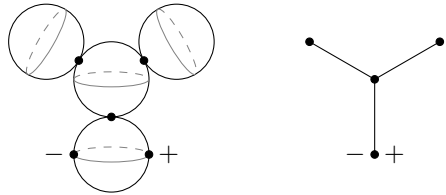
\begin{figure}[h]
  \centering
  \begin{tikzpicture}

    \begin{scope}
   \coordinate[label= left:$-$] (-) at (0,-0);
   \coordinate[label= right:$+$] (+) at (1,0);
   \draw (0.5,0) circle (.5);
   \draw[ultra thin, color=gray] (0,0) arc (-180:0:.5 and .1);
   \draw[ultra thin, color=gray,dashed] (0,0) arc (-180:-360:.5 and .1);
      \draw [fill,color=black] (+) circle (.05);
   \draw [fill,color=black] (0,0) circle (.05);
     \draw [fill,color=black] (.5,.5) circle (.05);
   \draw (0.5,1) circle (.5);
   \draw[ultra thin, color=gray] (0,1) arc (-180:0:.5 and .1);
   \draw[ultra thin, color=gray,dashed] (0,1) arc (-180:-360:.5 and .1);

   \draw (.5,1) [fill,color=black] +(30:.5) circle (.05);
   \draw (.5,1) [fill,color=black] +(150:.5) circle (.05);
   \draw ($(.5,1)+(30:1)$) circle (.5);
   \begin{scope}[shift={($(.5,1)+(30:1)+(120:.5)$)}, rotate=-60]
   \draw[ultra thin, color=gray] (0,0) arc (-180:0:.5 and .1);
   \draw[ultra thin, color=gray,dashed] (0,0) arc (-180:-360:.5 and .1); 
   \end{scope}
  \draw ($(.5,1)+(150:1)$) circle (.5);
   \begin{scope}[shift={($(.5,1)+(150:1)-(60:.5)$)}, rotate=60]
   \draw[ultra thin, color=gray] (0,0) arc (-180:0:.5 and .1);
   \draw[ultra thin, color=gray,dashed] (0,0) arc (-180:-360:.5 and .1); 
   \end{scope}
   \end{scope}
     \begin{scope}[shift={(4,0)}]
   \coordinate[label= left:$-$] (-) at (0,-0);
   \coordinate[label= right:$+$] (+) at (0,0);
   \draw [fill,color=black] (+) circle (.05);
   \draw (0,0) -- (0,1);
   \draw [fill,color=black] (0,1) circle (.05);

   \draw (0,1) -- ($(0,1)+(30:1)$);
\draw (0,1) [fill,color=black] +(30:1) circle (.05);
   \draw (0,1) -- ($(0,1)+(150:1)$);
   \draw (0,1) [fill,color=black] +(150:1) circle (.05);
   \end{scope}
  \end{tikzpicture}
  \caption{A representation of a pre-stable cylinder with automorphism group $\bZ/2$ and the corresponding tree; the great equators are supposed to help visualise that these are sphere components. The root carries the marked points $z_\pm$ which correspond to the ends $t = \pm \infty $ of $\Sigma_v \equiv \bR \times S^1$.}
  \label{fig:pre-stable-cylinder}
\end{figure}

We shall often write $\Sigma$ for a pre-stable building. The automorphism group $\Aut(\Sigma)$ is the direct product of the automorphism groups of the underlying pre-stable cylinders: these consist of a (rooted) automorphism $f$ of the corresponding tree and a biholomorphism $\phi_v$ of the Riemann surfaces $\Sigma_{v} \cong \Sigma_{f(v)}$ mapping the nodal point labelled by $e$ to the nodal point labelled by $f(e)$.

Consider a pair $x_\pm$ of free loops in $M$ (i.e., maps $S^1 \to M$). A map $u \co \bR \times S^1 \to M$ is \emph{asymptotic} to $x_\pm$ if $u(s,t)$ exponentially converges (in the $C^k$ norm for every $k$) to $x_\pm(t)$ in the limit $s \mapsto \pm \infty$.
\begin{defin} \label{def:stable_map_smooth}
  A \emph{stable map} $u$ with asymptotic conditions $x_\pm$ consists of a pre-stable cylinder and a map $u_v \co \Sigma_v \to M$ from each component such that:
  \begin{enumerate}
  \item (Asymptotic conditions) $u_v$ is asymptotic to $x_\pm$ if $v$ is the root,
  \item (Well-defined at the nodes) $u_v(z_{v,e}) = u_{v'}(z_{v',e})$ whenever $v$ and $v'$ are the two endpoints of an edge $e$, and
    \item (Stability) If $v$ is the root, and $u_v$ is independent of $\bR$, then the valence of $v$ is strictly positive. If $v$ is not the root, and $u_v$ is constant, then the valence of $v$ is strictly larger than two. 
    \end{enumerate}
\end{defin}
We refer to the map corresponding to the root as a cylinder, and the maps corresponding to other vertices as sphere bubbles. There is an equivalence relation on the set of stable maps defined as follows: an equivalence between two maps is a (rooted) isomorphism $f$ of the underlying tree and a biholomorphism $\phi_v$ of the Riemann surfaces $\Sigma_{v} \cong \Sigma_{f(v)}$ intertwining the marked points labelled by $e$ and $f(e)$ and the maps to $M$, such that the restriction to the root is a translation. In particular, the asymptotic conditions are locally constant in the space of stable maps. Note that, as soon as there is more than one vertex, the biholomorphism of the root must be the identity because no finite set is invariant under a non-trivial translation. An \emph{automorphism} of a stable map is a self-equivalence in the above sense.
\begin{rem}
In this generality, the group of automorphisms is not necessarily finite, as a composition of the height function $S^2 \to [0,1]$ with a path in $M$ will have $S^1$ as its group of automorphisms.
\end{rem}

Stable maps model the open subspace of the Gromov-Floer compactification of the moduli space of Floer cylinders obtained by allowing sphere bubbles. To obtain all strata of the compactification, we have to consider breaking of cylinders:
\begin{defin} 
A \emph{stable building} with asymptotic conditions given by loops $x_\pm\co S^1 \to M$ consists of a collection of loops $\{x_i\}_{i=1}^k$ and stable maps $\{u_i\}_{i=0}^k$ with asymptotic conditions $(x_{i}, x_{i+1})$, where we set $x_0 = x_-$ and $x_{k+1} = x_+$.

A pair of stable buildings are \emph{equivalent} if each of the corresponding stable maps are equivalent.
\end{defin}

\begin{rem}
The terminology of buildings goes back to the literature on symplectic field theory \cite{EliashbergGiventalHofer2000}. It is convenient in this setting to distinguish the bubbling phenomenon of pseudo-holomorphic curves, from the breaking phenomenon of solutions to Floer's equation.
\end{rem}

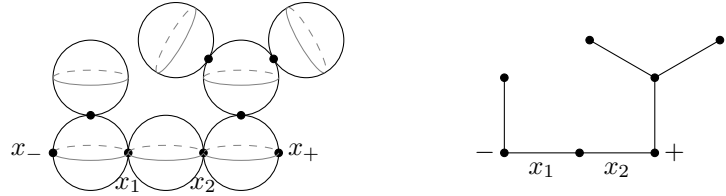
\begin{figure}[h]
  \centering
  \begin{tikzpicture}

    \begin{scope}
   \coordinate[label= below:$x_2$] (r) at (0,-0.2);
   \coordinate[label= right:$x_+$] (s) at (1,0);
   \draw [fill,color=black] (s) circle (.05);
   \draw [fill,color=black] (0,0) circle (.05);
   \draw (0.5,0) circle (.5);
   \draw [fill,color=black] (.5,.5) circle (.05);
   \draw (0.5,1) circle (.5);
 \draw[ultra thin, color=gray] (0,0) arc (-180:0:.5 and .1);
 \draw[ultra thin, color=gray,dashed] (0,0) arc (-180:-360:.5 and .1);
   \draw[ultra thin, color=gray] (0,1) arc (-180:0:.5 and .1);
   \draw[ultra thin, color=gray,dashed] (0,1) arc (-180:-360:.5 and .1);
   
   \draw (.5,1) [fill,color=black] +(30:.5) circle (.05);
   \draw (.5,1) [fill,color=black] +(150:.5) circle (.05);
   \draw ($(.5,1)+(30:1)$) circle (.5);
    \begin{scope}[shift={($(.5,1)+(30:1)+(120:.5)$)}, rotate=-60]
   \draw[ultra thin, color=gray] (0,0) arc (-180:0:.5 and .1);
   \draw[ultra thin, color=gray,dashed] (0,0) arc (-180:-360:.5 and .1); 
   \end{scope}
   \draw ($(.5,1)+(150:1)$) circle (.5);
 \begin{scope}[shift={($(.5,1)+(150:1)-(60:.5)$)}, rotate=60]
   \draw[ultra thin, color=gray] (0,0) arc (-180:0:.5 and .1);
   \draw[ultra thin, color=gray,dashed] (0,0) arc (-180:-360:.5 and .1); 
 \end{scope}
 
   \draw (-0.5,0) circle (.5);
   \coordinate[label= below:$x_1$] (q) at (-1,-0.2);
 \draw [fill,color=black] (-1,0) circle (.05);
 \draw (-1.5,0) circle (.5);
  \draw[ultra thin, color=gray] (-1,0) arc (-180:0:.5 and .1);
  \draw[ultra thin, color=gray,dashed] (-1,0) arc (-180:-360:.5 and .1);
 \draw[ultra thin, color=gray] (-2,0) arc (-180:0:.5 and .1);
 \draw[ultra thin, color=gray,dashed] (-2,0) arc (-180:-360:.5 and .1);

   \draw [fill,color=black] (-1.5,.5) circle (.05);
   \draw (-1.5,1) circle (.5);
      \draw[ultra thin, color=gray] (-2,1) arc (-180:0:.5 and .1);
   \draw[ultra thin, color=gray,dashed] (-2,1) arc (-180:-360:.5 and .1);
   \coordinate[label= left:$x_-$] (p) at (-2,-0);
   \draw [fill,color=black] (p) circle (.05);

   \end{scope}
   \begin{scope}[shift={(6,0)}]
   \coordinate[label= right:$+$] (+) at (0,0);
   \draw [fill,color=black] (+) circle (.05);
   \draw (0,0) -- (0,1);
   \draw [fill,color=black] (0,1) circle (.05);

   \draw (0,1) -- ($(0,1)+(30:1)$);
\draw (0,1) [fill,color=black] +(30:1) circle (.05);
   \draw (0,1) -- ($(0,1)+(150:1)$);
   \draw (0,1) [fill,color=black] +(150:1) circle (.05);

 \draw (0,0) -- (-1,0);
 \coordinate[label= below:$x_2$] (r) at (-.5,-0);
 \draw [fill,color=black] (-1,0) circle (.05);
 \draw (-1,0) -- (-2,0);
  \coordinate[label= below:$x_1$] (q) at (-1.5,-0);
  \draw (-2,0) -- (-2,1);
  \draw [fill,color=black] (-2,0) circle (.05);
   \draw [fill,color=black] (-2,1) circle (.05);
  \coordinate[label= left:$-$] (-) at (-2,-0);
   \end{scope}
  \end{tikzpicture}
  \caption{A representation of a stable building with automorphism group $\bZ/2$ and the corresponding tree.}
  \label{fig:stable-building}
\end{figure}

It shall be useful to encode the combinatorics of stable buildings by trees as follows: consider a tree equipped with a pair $v_\pm$ of distinguished vertices. A stable building thus consists of the following data: (i) for each vertex along the minimal path from $v_+$ to $v_-$, a cylinder mapping to $M$ with marked points labelled by the adjacent edges which do not lie on this path, and (ii) for each other vertex, a closed Riemann surface of genus $0$ with  marked points labelled by all adjacent edges, together with a map to $M$. Removing the edges along the path from $v_+$ to $v_-$, we obtain a collection of rooted trees; we require that the data associated to each such tree define a stable map as in Definition \ref{def:stable_map_smooth}.

\subsection{Lifts of Hamiltonian orbits and moduli spaces of stable Floer cylinders}
\label{sec:moduli-spaces-stable}

Let $H \co S^1 \times M \to \bR$ be a Hamiltonian function all of whose time-$1$ periodic orbits are non-degenerate. Let $\cP(H)$ denote the set of lifts of the contractible Hamiltonian orbits of $H$ to an intermediate regular cover $\tilde{\cL}M$ of the free loop space of $M$ on which the action and index of loops is well-defined (i.e., so that the torus in $M$ corresponding to any loop in this cover has trivial area and Chern number). Let $\Pi$ denote the group of deck transformations (i.e. the image of the surjective homomorphism from the fundamental group of $\cL M $ which is associated to this chosen cover). 

An element $p \in \cP(H)$ consists of an orbit of $H$ together with a homotopy class of bounding discs passing through a fixed basepoint in $M$ (capping disc). The action (the difference between the integral of $H$ over the orbit and the integral of the symplectic form $\omega$ over the bounding disc) and the (normalised) Conley-Zehnder index define maps
\begin{equation}
    \cA \co \cP(H) \to \bR \textrm{ and } \deg \co \cP(H) \to \bZ.
\end{equation}
We normalise the degree so that it is given by the dimension of the positive-definite subspace of the Hessian when computed for the constant capping discs of the orbits corresponding to the critical points of a Morse function (in particular, any minimum has degree equal to the dimension of $M$, and any maximum has degree $0$). There is a natural map $\Pi \to H_2(M, \bZ)$ which associates to each homotopy class of free loops the corresponding homotopy class of tori mapping to $M$, and we define 
\begin{equation}
    \cA \co \Pi \to \bR \textrm{ and } \deg \co \Pi \to \bZ.
\end{equation}
to be the composition with the maps on $H_2(M, \bZ) $ defined by $[\omega] \in H^2(M, \bR)$ and $c_1(M) \in  H^2(M, \bZ) $.

We write $[p]$ for the class of an element of $\cP(H)$ in the quotient $\cP(H)/\Pi$, and identify it with the corresponding orbit.  The action of $\Pi$ on $\cP(H)$ is free and the quotient is finite because a non-degenerate Hamiltonian function on a closed manifold admits only finitely many time-$1$ orbits. Given an $\omega$-tame almost complex structure $J$ on $M$, the following definition is standard, and is a variant of Floer's definition \cite{Floer1989b}, in the sense that it takes sphere bubbling into account:
\begin{defin}
A \emph{pseudo-holomorphic stable building}  is a stable building such that each sphere bubble is $J$-holomorphic (i.e. satisfies $du \circ j = J \circ du$), the Floer operator 
  \begin{equation} \label{eq:Floer-equation}
   \dbar_H \equiv  \left(du - X_H \otimes dt \right)^{0,1}   = \frac{1}{2} \left( ( du - X_H \otimes dt ) + J ( du - X_H \otimes dt) \circ j \right)
 \end{equation}
 vanishes on each cylinder, and the energy
   \begin{equation} \label{eq:Floer-equation2}
   E(u) = \int \|  du - X_H \otimes dt  \|^2
 \end{equation}
is finite, where the integral is taken over all components of the domain, and the inhomogeneous term $X_H \otimes dt$ vanishes on the sphere bubbles.

The \emph{moduli space  $\fcMR([p],[q])$ of Floer cylinders} is the space of equivalence classes of pseudo-holomorphic stable buildings, which are asymptotic to $[p]$ at $-\infty$, and $[q]$ at $+\infty$.
\end{defin}
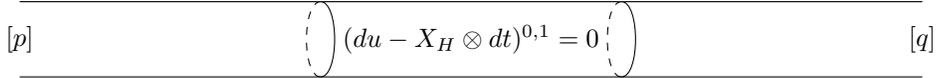
\begin{figure}[h]
  \centering
  \begin{tikzpicture}
    \coordinate[label= center:{$[p]$}] (p) at (-6,0);
    \coordinate[label= center:{$[q]$}] (q) at (6,0);
    \coordinate[label= center:{$(du - X_H \otimes dt)^{0,1}=0$}] (i) at (0,0);
    \begin{scope}[shift={(4,0)}]
     \draw [dashed,thin] (-2,-.5) .. controls (-2.25,-.5) and (-2.25,.5) .. (-2,.5);
    \draw [thin] (-2,-.5) .. controls (-1.75,-.5) and (-1.75,.5) .. (-2,.5); 
    \end{scope}
    \draw [dashed,thin] (-2,-.5) .. controls (-2.25,-.5) and (-2.25,.5) .. (-2,.5);
    \draw [thin] (-2,-.5) .. controls (-1.75,-.5) and (-1.75,.5) .. (-2,.5);
    \draw (-6,.5) -- (6,.5);
    \draw (-6,-.5) -- (6,-.5);
    
  \end{tikzpicture}
  \caption{The asymptotic conditions on elements of $\Mbar([p],[q])$.}
  \label{fig:Floer-cylinder}
\end{figure}
\begin{rem}
It is traditional to write $\Mbar([p],[q])$ for this moduli space, but we find it convenient to indicate the fact that the symmetry of the Floer components is broken from $\bC^*$ to $\bR$.
\end{rem}
The (Hausdorff) topology on $\fcMR([p],[q]) $ is obtained by implementing ideas of Gromov \cite{Gromov1985} and Kontsevich \cite{Kontsevich1995a} in this specific context, and is described in detail in \cite[Chapter 3]{FukayaOno1999}. We summarise the key ideas: the topology incorporates the possibility that a sequence of maps with domain a cylinder converge to a stable map with domain containing spheres, or a multi-level map with multiple cylinders. The key fact we shall use is the existence of a canonical decomposition
\begin{equation}
 \fcMR([p],[q])  \equiv \coprod_{b \in \pi_1(\cL M,[p], [q])}  \fcMR([p],[q]; b) 
\end{equation}
indexed by the  set $\pi_1(\cL M,[p], [q])$ of homotopy classes of paths in the free loop space from $[p]$ to $[q]$. There is a natural map $\pi_1(\cL M,[p], [q]) \to \bR$ given by the energy of any representing cylinder in $M$, which is the difference between the integrals of $\omega$ and of $d H \wedge ds$. The fundamental result we shall use is Gromov's compactness theorem, which appears in this context as \cite[Theorem 11.1]{FukayaOno1999}:
\begin{thm}
The map $ \fcMR([p],[q]) \to [0, \infty)$ induced by $\omega$ is proper, and the origin has a neighbourhood whose inverse image is trivial. \qed 
\end{thm}
\begin{rem}
The statement above is compatible with the existence of the translation-invariant cylinder associated to $[p]=[q]$ because the stability condition we imposed excludes it.
\end{rem}

We now return to consider lifted orbits $p,q \in \cP(H)$. As such lifts are points in $\tilde{\cL} M$, they determine a collection $\pi_1(\tilde{\cL} M,[p], [q])$ of homotopy classes of paths in $\cL M$ from $[p]$ to $[q]$ which lift to this space with the given endpoints. We define
\begin{equation} \label{eq:moduli_space_Floer-cylinders}
 \fcMR(p,q)  \equiv \coprod_{b \in \pi_1(\tilde{\cL} M,[p], [q])} \fcMR([p],[q]; b), 
\end{equation}
and note that this is a compact Hausdorff topological space by Gromov compactness. It is evident that each element $\pi \in \Pi$ induces an identification
\begin{equation}\label{eq:Mequivariant}
 \fcMR(p,q)  =  \fcMR(\pi \cdot p,\pi \cdot q) .  
\end{equation}

By construction, there is a natural associative map
\begin{equation}
 \fcMR(p,q) \times \fcMR(q,r) \to \fcMR(p,r)
\end{equation}
which concatenates stable buildings, and which strictly commutes with the $\Pi$-action.

In order for this construction to fit with our framework for constructing homotopy types from topologically enriched categories as described in Section \ref{sec:topol-enrich-categ}, it remains to equip $\cP(H)$ with the desired partial order, which we define by 
\begin{equation}
  \parbox{30em}{$p<q$ whenever $\fcMR(p,q)$ is non-empty.}
\end{equation}
Gromov's compactness theorem readily implies that, for each element $p$ of $\cP(H)$, and constant $E \in \bR$, there are finitely many elements $q$ such that $p < q$ and $\cA(q) \leq E$. Since the ordering is clearly preserved by $\Pi$, we conclude:
\begin{lem}
The moduli spaces $\fcMR(p,q)$ are the morphism spaces of a $\Pi$-equivariant topological flow category $\fcMR(H)$ with object set $\cP(H)$. \qed
\end{lem}

\subsection{Abstract moduli spaces of cylinders}
\label{sec:abstr-moduli-spac}

Consider a genus-$0$ Riemann surface $\Sigma$ with points marked by $\{+,-\}$. An \emph{angular lift} at $z_{\pm} $ is a choice of oriented real line in $T_{z_\pm} \Sigma$. Note that such a choice determines a biholomorphism to $\bP^1$ mapping the marked points to $0$ and $\infty$, uniquely up to positive real dilation, and in particular at choice at one end determines a choice at the other end.

We shall construct a moduli space of stable Riemann surfaces with marked points labelled by the union of $\pm$ with an ordered collection $S$ of finite sets, equipped with compatible angular lift at the points $z_\pm$. In order to construct a smooth structure on this moduli space (and later a stable almost complex structure), we start with the smooth structure on the Deligne-Mumford space. 

\begin{rem}
A delicate point in Floer theory is the construction of smooth structures on moduli spaces of abstract curves in such a way that moduli spaces of pseudo-holomorphic maps acquire a smooth structure for which the forgetful map is smooth. We shall not construct a smooth structure on the moduli space of  pseudo-holomorphic maps, so these considerations are irrelevant for our purpose.
\end{rem}

Given a finite set $r$, we shall consider the Deligne-Mumford space $\Mbar_{r \amalg \pm}$ of stable genus $0$ Riemann surfaces with points marked by the set $r \amalg \{\pm\}$. It is a standard fact going back to Knudsen \cite[Theorem 2.7]{Knudsen1983a} (see e.g. \cite[Appendix D.4]{McDuffSalamon2012} for a symplectic topology reference) that this space is a smooth complex projective variety stratified by the topological type of the underlying curve, which is determined by a tree $T$, and a partition of $r \amalg \pm$ among the vertices of $T$, so that every vertex is stable (i.e. the sum of the valence and the number of marked points associated to the vertex is not less than $3$).  In particular, there are vertices $v_\pm$ carrying the marked points $\{z_\pm\}$. We can think of each element of $\Mbar_{r \amalg \pm} $ as a building with levels given by a component corresponding to a vertex lying on the path from $+$ to $-$, together with all the components whose paths to $v_\pm$ pass through this vertex.  We write $\partial \Mbar_{r \amalg \pm} $ for the local normal crossings divisor corresponding to strata with at least two levels. This is further stratified by the number of edges between $v_-$ and $v_+$, and we write $\partial^k \Mbar_{r \amalg \pm} $ for the stratum with $k$ such edges. We say that curves on this stratum have \emph{$(k+1)$-levels,} which we order starting at the level containing the marked point $-$.

\begin{rem}
Our definition of the boundary divisor $\partial \Mbar_{r \amalg \pm} $ is not standard: one usually defines it to consist of the locus where the underlying curve has more than one component, so that our boundary divisor is a subset (in fact, an irreducible component) of the usual one. The reason for this choice is that, for the purpose of studying Floer theory, we need to treat edges separating $v_\pm$ differently from other edges.
\end{rem}

Note that the symmetry group on $r$ letters, which we denote $G_r$, acts on $\Mbar_{r \amalg \pm}$ preserving the boundary. We shall (partially) break this symmetry as follows (we invite the reader to look at Remark \ref{rem:why_the_moduli_space_is_so_complicated} below for some explanation for why all this data is required):
\begin{defin} \label{def:index-set-moduli-space-cylinders-marked-points}
  Given a pair $p < q$ of elements of $\cP(H)$,  we define $\cD(p,q)$ to be the set of data consisting of
  \begin{enumerate}
  \item a pair of positive integers $S$ and $S'$, and an injective order preserving map $\{1, \ldots, S'\} \to  \{1, \ldots, S \}$,
  \item a sequence $\{r_i\}_{i=1}^{S}$ of positive integers indexed by $S$,
  \item  an inclusion $P \subset P'$ of elements of $2^{\cP(H)}(p,q)$ (i.e. totally ordered subsets of $\cP(H)$ strictly between $p$ and $q$), and
    \item order preserving
  surjective maps from $\{1, \ldots, S\}$ to the successive elements of $\{p\} \amalg P \amalg \{q\}$, and from $\{1, \ldots, S'\}$ to the successive elements of $\{p\} \amalg P' \amalg \{q\}$.
\end{enumerate}
We require that these order preserving maps satisfy the following condition:  if $p_i < q_i$ and $p'_i < q'_i$ are the successive elements of $P$ and $P'$ associated to $i \in S'$ and its image in $S$, we require that
\begin{equation}
  p_i \leq p'_i < q'_i \leq q_i.
\end{equation}
\end{defin}
By abuse of notation, we write $S$ for the set $ \{1, \ldots, S\}$, and similarly for $S'$. In other words, we consider each natural number $S$ as an object of the category whose maps are order preserving injections. 

The above data thus consists of an assignment $r_i$ of a natural number for each $i \in S$; we write $r_i$ as well for the set of numbers $\{1, \ldots, r_i\}$. We write $r_{S}$ for the (disjoint) union of the sets $r_i$ indexed by the members of the sequence $S$, and similarly for $r_{S'}$.

For later purposes, we note the following immediate consequence of our construction:
\begin{lem} \label{lem:Pi-action-D}
  Each element $\pi \in \Pi$ induces a canonical bijection
  \begin{equation}
      \pi \cdot \_ \co   \cD(p,q) \to \cD(\pi \cdot p, \pi \cdot q),
      \end{equation}
      such that the following diagram
      \begin{equation}
        \begin{tikzcd}
          \cD(p,q) \ar[r," \pi \cdot \_"] \ar[dr," (\pi' \cdot \pi) \cdot \_"] &  \cD(\pi \cdot p, \pi \cdot q) \ar[d," \pi' \cdot \_"] \\
          & \cD(\pi' \cdot \pi \cdot p,\pi' \cdot  \pi \cdot q) 
        \end{tikzcd}
      \end{equation}
      commutes. 
      This assignment is strictly associative as is unital in that the identity of $\Pi$ acts by the identity map.
\qed
\end{lem}

\begin{defin} \label{def:moduli-cylinder-marked-points-quotient} 
  For each ${\underline{\alpha}} \in \cD(p,q)$, the moduli space $\Mbar_{\underline{\alpha}}$ (respectively $\Mbar'_{\underline{\alpha}}$) is the subset of  $\Mbar_{r_{S} \amalg \pm}$ (resp. $\Mbar_{r_{S'} \amalg \pm}$)  given by stable curves satisfying the following constraints:
  \begin{itemize}
  \item for each $i \in S$ (resp. $S'$), any two points marked by elements of $r_i$ lie in the same level
  \item the induced map from $S$ (resp. $S'$) to levels factors through an order preserving map from  the set of successive elements of $P$ (resp. $P'$) to the set of levels (here, $P \subset P'$ are the sets of orbits that are part of the data of ${\underline{\alpha}}$).

    \item the fibre of the universal curve is obtained by pullback from $\Mbar_{r_{S'} \amalg \pm}$.
    \end{itemize}
  \end{defin}

  As a consequence of the last condition, note that the projection map $\Mbar_{\underline{\alpha}} \to \Mbar'_{\underline{\alpha}} $  is a submersion, and that each curve in these moduli spaces is equipped with a labelling of the Floer edges of the corresponding tree by elements of $P$, given by the minimum of the labels appearing in the adjacent level closest to the output (this agrees with the maximum of the labels appearing in the adjacent level closest to the input, because we have assumed that the map from $S$ to successive elements of $P$ is surjective). 
  \begin{rem} \label{rem:why_the_moduli_space_is_so_complicated}
   To justify the above definition, the expert reader should anticipate that the basic charts on the interior of the moduli space $\fcMR(p,q)$ are associated to a choice $r$ of a number of marked points to stabilise the domains of trajectories, and the group acting on such a chart is the symmetric group on $r$ letters. From this point of view, the difficulty is to construct maps of charts associated to changing the number of marked points; the natural construction yields a correspondence for each inclusion $f \co \{1, \ldots, r'\} \to  \{1, \ldots, r\}$, i.e. the partition $r = \mathrm{Im} f \amalg (r \setminus \mathrm{Im} f)$, with automorphism group given by the products of the symmetric groups of order $r'$ and $r- r'$. Since our point of view is to consider maps of charts rather than correspondences, we are thus led naturally to consider charts labelled by a finite collection $S$ of integers, in which case the automorphism group will be the product of symmetric groups indexed by $S$. It is easy in this case to associate forgetful maps of Kuranishi charts to inclusions $S' \subset S$.

  The reason for labelling our basic charts by inclusions $S' \subset S$, rather simply by such a set, ultimately is related to the fact that we will need to lift our Kuranishi charts to flag smooth charts (see Section \ref{sec:fiber-kuran-pres}); in order to avoid discussions of smoothness of gluing maps, it is convenient to ensure that all maps of abstract moduli spaces that we consider are smooth fibrations, and forgetful maps which change the domain of the underlying map may not satisfy this property.

  In the above discussion, we considered only the interior of the moduli space. The description of a general Kuranishi chart of the moduli space requires a choice of a totally ordered set $P$ of $\cP(p,q)$ of intermediate orbits. Passing to pairs $P \subset P'$ will allow us to more clearly formulate the flag smooth structure discussed above, but will also be crucial in constructing functorial and multiplicative (relative) stable complex orientations in Section \ref{sec:stable-almost-compl-1} (see specifically Subsection \ref{sec:oriented-locus-space}). 

  Finally, the fact that the data of intermediate orbits and choices of marked points are not chosen independently is related to the fact that we only consider strict $2$-functors in our definition of a Kuranishi presentation. We shall later observe the existence of a natural map $\cD(p,q) \times \cD(q,r) \to \cD(p,r)$, and use this correspondence to associate to a product of charts for the pairs $(p,q)$ and $(q,r)$, a chart for the pair $(p,r)$. The map from $S$ to the successive elements of $P$ will ensure that this is in fact the product of charts (rather than a chart containing the product chart as an open subset). We could drop this choice from our data at the cost of re-defining Kuranishi presentations using lax $2$-functors.
\end{rem}
\begin{example}
  Consider a pair $p < r$ of elements of $\cP(H)$. Let $P = \emptyset$, and assume that we have an element $q$ of $\cP(H)$ such that $p < q <r$; let $P' = \{q\}$. Let $S = \{1 < 2 < 3\}$, and $S' = \{1 < 3\}$, so that we have a triple of totally ordered sets $r_1$, $r_2$, and $r_3$, which we assume all consist of a pair of elements, which we denote $z_{i}^{j}$ for $1 \leq i \leq 3$, and $1 \leq j \leq 2$. Since $P$ is empty, we assign to each $i \in S$ the pair of elements $p < r$. We refine this for elements of $S'$ by assigning $p < q$ to the element $1$ and $q < r$ to the element $2$.

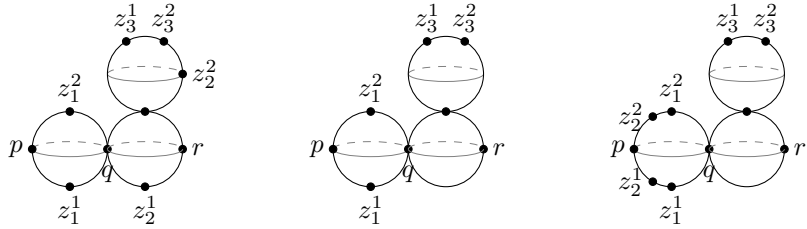
\begin{figure}[h]
  \centering
  \begin{tikzpicture}
    \begin{scope}
   \coordinate[label= left:$p$] (p) at (-1,0);
   \coordinate[label= below:$q$] (q) at (0,-0.1);
   \coordinate[label= right:$r$] (r) at (1,0);
   \draw (-0.5,0) circle (.5);
 \draw[ultra thin, color=gray] (-1,0) arc (-180:0:.5 and .1);
 \draw[ultra thin, color=gray,dashed] (-1,0) arc (-180:-360:.5 and .1);
   \draw [fill,color=black] (p) circle (.05);
   \draw [fill,color=black] (r) circle (.05);
   \draw [fill,color=black] (0,0) circle (.05);
   \draw (0.5,0) circle (.5);
    \draw[ultra thin, color=gray] (0,0) arc (-180:0:.5 and .1);
    \draw[ultra thin, color=gray,dashed] (0,0) arc (-180:-360:.5 and .1);
     \draw[ultra thin, color=gray] (0,1) arc (-180:0:.5 and .1);
    \draw[ultra thin, color=gray,dashed] (0,1) arc (-180:-360:.5 and .1);
   \draw [fill,color=black] (.5,.5) circle (.05);
   \draw (0.5,1) circle (.5);

   \draw (.5,1) [fill,color=black] +(60:.5) circle (.05);
   \draw (.5,1) [fill,color=black] +(120:.5) circle (.05);
   \coordinate[label= above:$z_3^1$] (z31) at ($(.5,1)+(120:.5)$) ;
   \coordinate[label= above:$z_3^2$] (z32) at ($(.5,1)+(60:.5)$);

    \draw  [fill,color=black] (-.5,-.5) circle (.05);
    \draw [fill,color=black] (-.5,.5) circle (.05);
    \coordinate[label=  below:$z_1^1$] (z11) at (-.5,-.5) ;
   \coordinate[label= above:$z_1^2$] (z12) at (-.5,.5);   
 \end{scope}
   \begin{scope}[shift={(4,0)}]
   \coordinate[label= left:$p$] (p) at (-1,0);
   \coordinate[label= below:$q$] (q) at (0,-0.1);
   \coordinate[label= right:$r$] (r) at (1,0);
   \draw (-0.5,0) circle (.5);
    \draw[ultra thin, color=gray] (-1,0) arc (-180:0:.5 and .1);
 \draw[ultra thin, color=gray,dashed] (-1,0) arc (-180:-360:.5 and .1);
   \draw [fill,color=black] (p) circle (.05);
   \draw [fill,color=black] (r) circle (.05);
   \draw [fill,color=black] (0,0) circle (.05);
   \draw (0.5,0) circle (.5);
       \draw[ultra thin, color=gray] (0,0) arc (-180:0:.5 and .1);
    \draw[ultra thin, color=gray,dashed] (0,0) arc (-180:-360:.5 and .1);
     \draw[ultra thin, color=gray] (0,1) arc (-180:0:.5 and .1);
    \draw[ultra thin, color=gray,dashed] (0,1) arc (-180:-360:.5 and .1);
   \draw [fill,color=black] (.5,.5) circle (.05);
   \draw (0.5,1) circle (.5);

   \draw (.5,1) [fill,color=black] +(60:.5) circle (.05);
   \draw (.5,1) [fill,color=black] +(120:.5) circle (.05);
   \coordinate[label= above:$z_3^1$] (z31) at ($(.5,1)+(120:.5)$) ;
   \coordinate[label= above:$z_3^2$] (z32) at ($(.5,1)+(60:.5)$);

    \draw  [fill,color=black] (-.5,-.5) circle (.05);
    \draw [fill,color=black] (-.5,.5) circle (.05);
    \coordinate[label=  below:$z_1^1$] (z11) at (-.5,-.5) ;
   \coordinate[label= above:$z_1^2$] (z12) at (-.5,.5);   
  \draw (-.5,0) [fill,color=black] +(120:.5) circle (.05);
   \draw (-.5,0) [fill,color=black] +(-120:.5) circle (.05);
   \coordinate[label= left:$z_2^1$] (z31) at ($(-.5,0)+(-120:.5)$) ;
   \coordinate[label= left:$z_2^2$] (z32) at ($(-.5,0)+(120:.5)$);
 \end{scope}
    \begin{scope}[shift={(-4,0)}]
   \coordinate[label= left:$p$] (p) at (-1,0);
   \coordinate[label= below:$q$] (q) at (0,-0.1);
   \coordinate[label= right:$r$] (r) at (1,0);
   \draw (-0.5,0) circle (.5);
    \draw[ultra thin, color=gray] (-1,0) arc (-180:0:.5 and .1);
 \draw[ultra thin, color=gray,dashed] (-1,0) arc (-180:-360:.5 and .1);
   \draw [fill,color=black] (p) circle (.05);
   \draw [fill,color=black] (r) circle (.05);
   \draw [fill,color=black] (0,0) circle (.05);
   \draw (0.5,0) circle (.5);
     \draw[ultra thin, color=gray] (0,0) arc (-180:0:.5 and .1);
    \draw[ultra thin, color=gray,dashed] (0,0) arc (-180:-360:.5 and .1);
     \draw[ultra thin, color=gray] (0,1) arc (-180:0:.5 and .1);
     \draw[ultra thin, color=gray,dashed] (0,1) arc (-180:-360:.5 and .1);
     \draw [fill,color=black] (.5,.5) circle (.05);
   \draw (0.5,1) circle (.5);

   \draw (.5,1) [fill,color=black] +(60:.5) circle (.05);
   \draw (.5,1) [fill,color=black] +(120:.5) circle (.05);
   \coordinate[label= above:$z_3^1$] (z31) at ($(.5,1)+(120:.5)$) ;
   \coordinate[label= above:$z_3^2$] (z32) at ($(.5,1)+(60:.5)$);

    \draw  [fill,color=black] (-.5,-.5) circle (.05);
    \draw [fill,color=black] (-.5,.5) circle (.05);
    \coordinate[label=  below:$z_1^1$] (z11) at (-.5,-.5) ;
    \coordinate[label= above:$z_1^2$] (z12) at (-.5,.5);

  \draw  [fill,color=black] (.5,-.5) circle (.05);
   \draw  [fill,color=black] (1,1)  circle (.05);
   \coordinate[label= below:$z_2^1$] (z31) at (.5,-.5 ) ;
   \coordinate[label= right:$z_2^2$] (z32) at (1,1);
    
 \end{scope}
  
     \end{tikzpicture}
  \caption{A graphical representation of an element of $ \Mbar'_{\underline{\alpha}}$ (center), and two elements of its inverse image in $ \Mbar_{\underline{\alpha}}$ (left and right). The marked points can have arbitrary position within the given component.}
  \label{fig:DM-moduli-ordered-collection-marked}
\end{figure}

The moduli space $ \Mbar'_{\underline{\alpha}}$ is thus a subset of the moduli space of stable spheres with four marked points (corresponding to the union of $\{p,r\}$ with  $r_1 \amalg r_3$). The key condition is that, if the point $p$ and $r$ do not lie in the same component, then they are separated by exactly one nodal point (which we label $q$), and all points labelled by $r_1$ and $r_3$ respectively lie on the same side of this node as $p$ and $r$ (see Figure \ref{fig:DM-moduli-ordered-collection-marked}).
\end{example}

We define the \emph{moduli space of Floer cylinders} $\fcMR_{r_{S}} \to \Mbar_{r_{S} \amalg \pm} $ to be the moduli space of Riemann surfaces with points marked by $r_{S} \amalg \pm$, and an asymptotic marker at the positive end of each Floer node. This is a smooth manifold with boundary, which can be defined as a circle bundle over an oriented boundary blowup of Deligne-Mumford space along the boundary divisor $\partial \Mbar_{r_{S} \amalg \pm}$, but we shall give an alternative and explicit construction in Section \ref{sec:constr-moduli-space} below (the construction as a blowup would proceed along the lines of the construction of \cite{KimuraStasheffVoronov1995}).  For each ${\underline{\alpha}} \in \cD(p,q)$, we denote the inverse image of $\Mbar_{{\underline{\alpha}}}$ in this moduli space by $\fcMR_{\underline{\alpha}}$. We similarly obtain a moduli space $\Mbar^{\prime \bR}_{\underline{\alpha}} $ over $\Mbar'_{\underline{\alpha}}$, and a submersion $  \fcMR_{\underline{\alpha}} \to \Mbar^{\prime \bR}_{\underline{\alpha}} $ over the projection $\Mbar'_{\underline{\alpha}} \to \Mbar_{\underline{\alpha}}$.

\begin{rem}
The appearance of two moduli spaces $\fcMR_{\underline{\alpha}}$ and  $\Mbar^{\prime \bR}_{\underline{\alpha}}$ is particularly relevant for the functoriality of flag smooth structures as in Section \ref{sec:fiber-kuran-pres}.
\end{rem}

\subsection{Continuation maps}
\label{sec:cont-maps-morse}

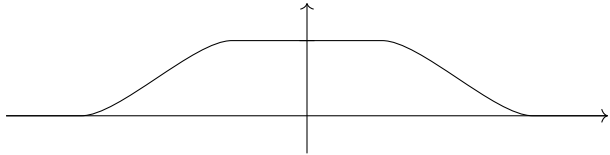
\begin{figure}[h]
  \centering
  \begin{tikzpicture}
    \draw [thin] (-.1,1) -- (.1,1);
    \draw  (-4,0) -- (-3,0) .. controls (-2.5,0) and (-1.5,1) .. (-1,1) -- (1,1) .. controls (1.5,1) and (2.5,0) .. (3,0) --   (4,0);
    \draw [->,thin] (-4,0) -- (4,0);
    \draw [->,thin] (0,-.5) -- (0,1.5);
  \end{tikzpicture}
  \caption{The graph of the cutoff functions $\chi_S$ for $S$ sufficiently large.}
  \label{fig:cutoff-function-graph-S}
\end{figure}

Floer's construction of homology groups associated to Hamiltonians extends to map associated to $1$-parametric families of Hamiltonians, i.e. maps from $\bR \times S^1$ to the space of functions on $M$. We only need the following special situation:
\begin{defin}
  An \emph{admissible continuation equation} is a map  $\bR \times S^1 \times M \to \bR$ whose restriction to each end of the cylinder either vanishes or agrees with $H$. 
\end{defin}
\begin{figure}[h]
  \centering
  \begin{tikzpicture}
    \coordinate[label=right:{$1$}] (1) at (0,1);
    \draw [thin] (-.1,1) -- (.1,1);
    \draw  (-4,0) -- (-1,0) .. controls (-.5,0) and (.5,1) .. (1,1) -- (4,1);
    \draw [->,thin] (-4,0) -- (4,0);
    \draw [->,thin] (0,-.5) -- (0,1.5);
  \end{tikzpicture}
  \caption{The graph of the cutoff function $\chi$.}
  \label{fig:cutoff-function-graph}
\end{figure}
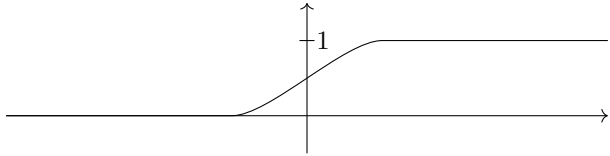
Using a cut-off function $\chi \co \bR \to [0,1]$ which vanishes at the negative end and equals $1$ at the positive end, we associated to $H$ a pair 
\begin{align}
H_+(s,t) & \equiv \chi(s)  H(t) \\
H_-(s,t) & \equiv \chi(-s)  H(t)
\end{align}
of admissible continuation equations. We also fix a family $\chi_S$ of functions from  $\bR$ to $[0,1]$, parametrised by $S \in [0,\infty)$ such that $\chi_0$ vanishes identically, and 
\begin{equation}
  \chi_{S}(s) = \begin{cases} \chi_+(s+S) & s \leq 0 \\
\chi_-(s-S) & 0 \leq s 
\end{cases}
\end{equation}
whenever $S \in [0,\infty)$ is sufficiently large (see Figure \ref{fig:cutoff-function-graph-S}).
We then define a family $H_{S}$ of continuation equations parametrised by $S \in [0,\infty)$, given by
\begin{equation}
  H_{S}(s,t) = \chi_{S}(s) H(t).
\end{equation}
Note that the definition makes sense because $H_{S}(0,t) = H(t) $ at the boundary of the domains of definition. In the limit $S \to +\infty$, the domain splits into two copies of the real line respectively carrying the Hamiltonians $H_+$ and $H_-$. We write $\{ H_S \}_{S \in [0,\infty]}$ when we allow this possibility.

We now consider pseudo-holomorphic maps associated to such equations. The key difference is that stability is defined relative the map $ \bR \times S^1 \times M \to \bP^1 \times M$, which is given by the natural inclusion on the first factor (see Definition \ref{def:moduli_extra_marked_point_continuation}); note that this stability condition is different even when the function $H$ identically vanishes.
\begin{defin}
A \emph{pseudo-holomorphic stable map} with respect to an admissible equation $H_{\bullet}$ is a stable map such that each sphere bubble is $J$-holomorphic, and the cylindrical component satisfies the equation
  \begin{equation} \label{eq:continuation-equation}
   \dbar_{H_\bullet} \equiv  \left(du - X_{H_\bullet} \otimes dt \right)^{0,1} = 0.    
 \end{equation}
\end{defin}
The moduli space of stable maps as defined above may not be compact because we have not yet incorporated breaking or bubbling along the ends.

The compactification is obtained by considering pseudo-holomorphic stable buildings as before: we consider a chain of pseudo-holomorphic stable maps with compatible asymptotic conditions (see Section \ref{sec:abstr-moduli-spac-1} below for the construction of the corresponding abstract moduli spaces).

Given a Hamiltonian orbit $[p]$ of $H$, we extract from the space of pseudo-holomorphic buildings for $H_+$ the moduli space $\fcMR(M,[p])$ of equivalence classes of maps which converge to $[p]$. These consist of a stable map whose cylindrical component solves Equation \eqref{eq:continuation-equation} for $H_\bullet = H_+$, and a collection of stable maps whose cylindrical component solves Floer's equation (for $H$) with matching asymptotic conditions, such that the last asymptotic condition is given by $[p]$ (along the positive end). Note that, since the inhomogeneous term vanishes at the end labelled $M$, the marked point $z_-$ may lie on a sphere bubble. The same construction for $H_-$ yields a moduli space $ \fcMR([p],M)$ with asymptotic condition $[p]$ along the negative end.

The identically vanishing Hamiltonian defines a moduli space $\fcMR(M,M;0)$ consisting of equivalence classes of maps with constant asymptotic conditions at both ends. Finally, the $1$-parametric family of continuation equations $\{H_{S}\}_{S \in  [0,\infty]}$ yields a moduli space $\fcMR(M,M)$ of equivalence classes of maps with constant asymptotic conditions at both ends. In the case $S = \infty$, an element of this moduli space consists of a chain of stable maps with stable asymptotic conditions, the first and last of which respectively correspond to $H_\bullet = H_\pm$, while the intermediate ones are solutions to Floer's equation.

For each orbit $[p]$, we have tautological inclusions
\begin{equation} \label{eq:glue-two-half-planes-to-sphere}
   \fcMR(M,M;0) \to  \fcMR(M,M) \leftarrow \fcMR(M,[p]) \times \fcMR([p],M),
 \end{equation}
and for each orbit $[q]$ we have as well natural inclusions of products
 \begin{align}
   \fcMR(M,[p]) \times \fcMR([p],[q]) & \to     \fcMR(M,[q]) \\ \label{eq:glue-positive-half-plane-cylinder}
\fcMR([p],[q])  \times   \fcMR([q],M) & \to \fcMR([p],M).
 \end{align}
In this context, Gromov's compactness theorem implies:
\begin{lem}
The subsets of $\fcMR(M,[p])$,  $ \fcMR([p],M)$, and $\fcMR(M,M)$ consisting of elements of bounded energy are compact. Moreover, there is a constant $\epsilon$ so that the only elements of $ \fcMR(M,M;0)$ of energy smaller than $\epsilon$ are constant. \qed
\end{lem}
\subsection{Combining continuation maps and Morse trajectories}
\label{sec:morse-theory}
Let $f$ be a Morse-Smale function as in Appendix \ref{sec:morse-theory-and-homotopy-types}; to harmonise notation, we shall write $[x]$ for a critical point of $f$. We shall presently use the moduli spaces $  \barT([x],M) $, $ \barT([x],[y])$, and $ \barT([y],M)$, associated to a pair $[x]$ and $[y]$ of critical points of $f$, and consisting of flow lines with convergence conditions at one or both ends given by fixed critical points. These are defined more precisely in the Appendix.

We introduce an additional moduli space $\barT([x],M,M)$ whose interior consists of a choice of a point on a (half-infinite) gradient flow line starting at $[x]$. There are again natural maps:
\begin{equation}
  \barT([x],M) \to    \barT([x],M,M) \leftarrow  \barT([x],[y]) \times  \barT([y],M)
\end{equation}
corresponding respectively to the locus where the additional point is the finite endpoint, and where it converges to the infinite end. We also introduce the moduli space
\begin{equation}
  \barT([x],M,[y])  \equiv    \barT([x],M)  \times_M  \barT(M,[y])
\end{equation}
consisting of a pair of gradient flow lines meeting at a point in $M$.

\begin{figure}[h]
  \centering
  \begin{tikzpicture}
   \coordinate[label= center:$x_-$] (x) at (-4,0);
   \coordinate[label= center:$z_+$] (z) at (4,0);
   \coordinate[label= center:$p$] (p) at (0,4);
   \coordinate[label= center:$y_0$] (y) at (0,-4);

   \draw [thick] (-4,0.25) .. controls (-4,4) .. (-2,4);
   \draw (-0.5,3.5) -- (-1.5,3.5) arc (-90:-270:.5) --  (-0.5,4.5);
 \coordinate[label= center:{$\dbar_{H_+}=0$}] (+) at (-1.25,4);

    \draw [thick] (4,0.25) .. controls (4,4) .. (2,4);
    \draw (0.5,3.5) -- (1.5,3.5) arc (-90:90:.5) --  (0.5,4.5);
    \coordinate[label= center:{$\dbar_{H_-}=0$}] (-) at (1.25,4);

\draw [thick] (3.75,0.25) .. controls (2,2) .. (1.5,2);
\draw  (1,1.5) arc (-90:90:.5) -- (-1,2.5) arc (90:270:.5) -- cycle;
\draw [thick] (-3.75,0.25) .. controls (-2,2) .. (-1.5,2);
\coordinate[label= center:{$\dbar_{H_S}=0$}] (-) at (0,2);

    \draw [thick] (-4,-0.25) .. controls (-4,-2) .. (-3,-3)  .. controls (-2,-4) .. (-0.25,-4);
    \draw [fill,color=black] (-3,-3) circle (.05);

     \draw [thick] (4,-0.25) .. controls (4,-2) .. (3+.5*sqrt{.5},-3+.5*sqrt{.5}); 
     \draw (3,-3) circle (.5);
 \coordinate[label= center:{$\dbar=0$}] (-) at (3,-3);    
     \draw [thick] (0.25,-4) .. controls (2,-4) .. (3-.5*sqrt{.5},-3-.5*sqrt{.5});

 \draw [thick] (-3.75,0) -- (-0.5,0); 
 \draw (0,0) circle (.5);
 \coordinate[label= center:{$\dbar=0$}] (0) at (0,0);
     \draw [thick] (3.75,0) -- (.5,0); 

\draw [thick] (3.75,-0.25) -- (2+.5*sqrt{.5},-2+.5*sqrt{.5}); 
\draw (2,-2) circle (.5);
 \coordinate[label= center:{$\dbar=0$}] (-) at (2,-2);  
     \draw [thick] (-3.75,-.25) -- (-2+.5*sqrt{.5},-2-.5*sqrt{.5}) .. controls (-1,-3) .. (0,-3) .. controls (1,-3) .. (2-.5*sqrt{.5},-2-.5*sqrt{.5});
     \draw [fill,color=black] (-2,-2) circle (.05);

  \end{tikzpicture}
  \caption{A graphical representation of the flow category $\fcMR(f,H,f)$.}
  \label{fig:PSS-flow-category}
\end{figure}
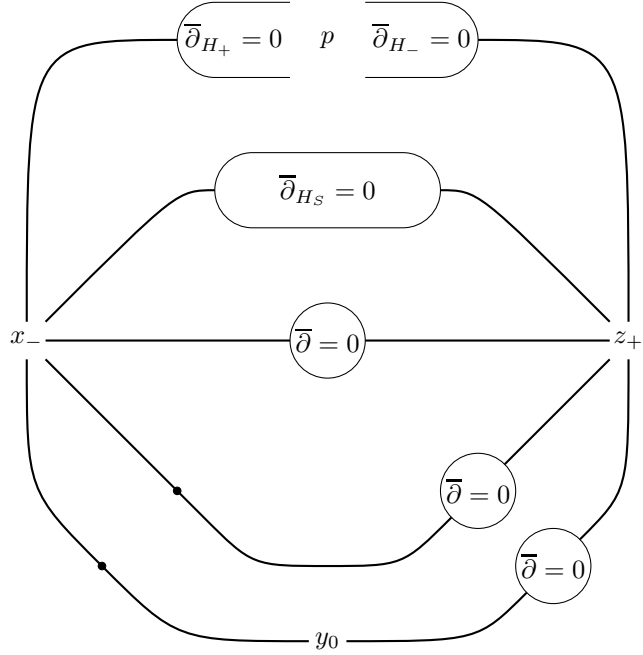

Let $\cP(f)$ denote the set of lifts of the critical points of $f$ (considered as constant maps with domain $S^1$) to $\tilde{\cL}M$. An element of $\cP(f)$ is thus represented as before by a homotopy class of bounding discs. We shall construct a flow category with objects three copies of $\cP(f)$, which we denote $\cP(f)_-$, $ \cP(f)_0$, and $\cP(f)_+$, and one copy of $\cP(H)$. To define the morphism spaces in this category, we associate to elements $(x,y,z)$ of $\cP(f)$ and  $p$ of $\cP(H)$ the  moduli spaces
\begin{align} \label{eq:moduli-space-Morse-to-Floer}
\fcMR(x_-,p)  & \subset \barT([x],M) \times_M  \fcMR(M,[p]) \\
\fcMR(p,z_+) & \subset \fcMR([p],M)  \times_M  \barT(M,[z]) \\
\fcMR(y_0,z_+) & \subset \barT([y],M) \times_M \fcMR(M,M;0) \times_M  \barT(M,[z]) \\ 
\fcMR(x_-,y_0) & \subset   \barT([x],M,[y])   \label{eq:moduli_space_Morse_0_+}
\end{align}
consisting of elements which represent paths in $\tilde{\cL}M$ with the prescribed asymptotic conditions. These moduli spaces are represented on the outer edges of Figure \ref{fig:PSS-flow-category}.

To construct the moduli space associated to a pair $(x_-,z_+)$, we introduce a space whose natural boundary consists of the products of spaces labelling the top and bottom of Figure \ref{fig:PSS-flow-category}. We proceed in two step: first, we introduce the moduli spaces
\begin{align} 
\fcMR(x_-,z_+)_{\mathrm{Floer}} & \subset \barT([x],M) \times_M \fcMR(M,M) \times_M  \barT(M,[z]) \\
\fcMR(x_-,z_+)_{\mathrm{Morse}} & \subset  \barT([x],M,M)   \times_M \fcMR(M,M;0) \times_M  \barT(M,[z])
\end{align}
where we again prescribe the homotopy class imposed by the input and output. Note that both of these spaces include, as a boundary stratum, the subset of the fibre product
\begin{equation}
   \barT([x],M) \times_M \fcMR(M,M;0) \times_M  \barT(M,[z]),
\end{equation}
consisting of configurations in the prescribed homotopy class (this is represented by the central diagram in Figure \ref{fig:PSS-flow-category}). Taking the union along this stratum, we obtain
\begin{equation} \label{eq:piecewise_definition_moduli_-_+}
\fcMR(x_-,z_+) \equiv \fcMR(x_-,z_+)_{\mathrm{Morse}}  \cup \fcMR(x_-,z_+)_{\mathrm{Floer}} .  
\end{equation}
For notational consistency, for each pair  $x$ and $y$ of lifts to  $\tilde{\cL}M$  of  critical points of $f$,  we define the moduli spaces $\barT(x,y) \subset \barT([x],[y])$ of gradient flow lines that lift as paths in $\tilde{\cL}M $ with the prescribed endpoints, and we then set
\begin{equation}
  \begin{rcases}
    &\fcMR(x_-,y_-)  \\
    & \fcMR(x_0,y_0) \\
    &  \fcMR(x_+, y_+)
  \end{rcases} \equiv \barT(x_0,y_0).
\end{equation}
\begin{rem}
It would be more natural to construct a category involving only the moduli spaces appearing in the upper half of Figure \ref{fig:PSS-flow-category} (i.e. use only the moduli space $\fcMR(x_-,z_+)_{\mathrm{Floer}}$, and not $ \fcMR(x_-,z_+)_{\mathrm{Morse}}$). However, the corresponding flow category does not admit a fundamental chain in the sense of this paper. We choose not to develop the formalism of \emph{relative fundamental chains of flow categories} that would be required to obtained the desired results in this context.
\end{rem}

The moduli spaces described above are compact because prescribing the lift determines the action. For the statement of a more refined consequence of Gromov's compactness theorem, we introduce the set
\begin{equation}
  \cP(f,H,f) \equiv \cP(f)_- \amalg \cP(H) \amalg \cP(f)_0 \amalg \cP(f)_+,
\end{equation}
where $\cP(f)_{\bullet}$, for $\bullet \in \{+,0,-\}$ is a copy of $\cP(f)$ with elements denoted $x_\bullet$. We equip this set with an ordering given by
\begin{equation}
  \parbox{30em}{$a<b$ whenever $\fcMR(a,b)$ is non-empty.}
\end{equation}
Gromov compactness again immediately implies that this ordering satisfies the property that, given any $a$, the set of elements $b$ which receive morphisms from $a$, and have bounded action, is finite.

The $\Pi$ action on $\cP(f,H,f)$ preserves this partial order. The moduli spaces $\fcMR(a,b)$ are the morphisms of a category whose composition maps are induced by Equation
\eqref{eq:glue-two-half-planes-to-sphere}-\eqref{eq:glue-positive-half-plane-cylinder},
as well as Equations
\eqref{eq:breaking-gradient-trajectory}-\eqref{eq:breaking-gradient-negative-flow-line}.
We conclude (recalling Equation~\eqref{eq:Mequivariant} for the $\Pi$-action):

\begin{lem}
The morphisms above define a strictly $\Pi$-equivariant topological
category which we denote $\fcMR(f,H,f)$ with objects $\cP(f,H,f)$.\qed
\end{lem}

There is one minor inconvenience with the above construction: while there is a natural map $\cA$ from $\cP(f,H,f)$ to $\bR$ given by the action of each lift to an orbit to  $\tilde{\cL}M$, this map does not respect the order, because, unlike the moduli space of solutions to Floer's equation, the moduli space of continuation maps need not be empty when its topological energy is negative. We shall remedy this by recalling the following standard result:
\begin{lem} \label{lem:control-energy-mixed-moduli-spaces}
  There is a constant $C_H$ with the property that,  for any pair of elements $a$ and $b$ of $\cP(f,H,f)$, the moduli space $\fcMR(a,b)$ is empty whenever
  \begin{equation}
  \cA(b) - \cA(a)   \leq - C_H.
  \end{equation}
  Moreover, there is a positive constant $\epsilon$ so that if $x_0 \in \cP(f)_0$ and $y_+ \in \cP(f)_+$ and $\cA(y) - \cA(x) < \epsilon$, then  the only moduli spaces $\fcMR(x_0,y_+)$ which are not empty are those consisting of a pair of gradient flow lines and a constant sphere. \qed 
\end{lem}

\subsection{Abstract moduli spaces of half-planes and spheres}
\label{sec:abstr-moduli-spac-1}

The moduli spaces of Section \ref{sec:abstr-moduli-spac} are sufficient to construct the Kuranishi flow category associated to a Hamiltonian. In order to compare the corresponding homotopy type to the one obtained from Morse theory, we consider the following construction:

We begin by fixing the categories of data associated to each case, extending Definition \ref{def:index-set-moduli-space-cylinders-marked-points}. To this end, we introduce the set
\begin{equation}
  \underline{\cP}(f,H,f) \equiv \{ M_-, M_0, M_+ \} \amalg \cP(H),
\end{equation}
which is partially ordered by $M_- < M_0 < M_+$, and $M_- < p < M_+ $ for each element $p \in \cP(H)$, and the previously given ordering of elements $\cP(H)$. We have a natural order preserving map
\begin{align}
  \cP(f,H,f) & \to  \underline{\cP}(f,H,f) \\
  a & \mapsto \underline{a},
\end{align}
which is given by the identity on $\cP(H)$, and the evident projections
\begin{equation}
 \cP(f)_- \amalg  \cP(f)_0 \amalg \cP(f)_+ \to \{M_-,  M_0, M_+\}.
\end{equation}

\begin{defin}
  Given an ordered pair $a < b$ of elements of $\cP(f,H,f)$, we define a category $\cD(\underline{a},\underline{b})$ with objects consisting of an inclusion of totally ordered subsets $P \subset P'$ of $ \underline{\cP}(\underline{a},\underline{b}) \cap \cP(H)$,
  
  and the following additional data if the pair $(\underline{a},\underline{b})$ is not contained in $\{ M_-,  M_0 \}$, or in $\{M_+\}$:
  \begin{enumerate}
  \item an order preserving injective map $\{1, \ldots, S'\} \to \{1, \ldots, S\}$ of sets indexing a collection $\{r_i\}_{i=1}^{S}$  of natural numbers, and
    \item a surjective map from $S$ to the successive elements of $\{\underline{a}\} \amalg P \amalg \{\underline{b}\}$, and a refinement of the restriction to $S'$, valued in the set of successive elements of $\{\underline{a}\} \amalg P' \amalg \{\underline{b}\} $.
  \end{enumerate}

An arrow $\underline{\alpha}_0 \to \underline{\alpha}_1$ is given by the following data:
  \begin{enumerate}
  \item  a factorisation of the map $ S^\prime_0 \to S_0 $ into a composition of injective maps
\begin{equation} \label{eq:ordering-nested-inclusion}
     S^\prime_0 \to S^\prime_1 \to  S_1 \to S_0
\end{equation}
with the property that the sequence $\{r_i\}_{i=1}^{S_1}$ is obtained by pulling back  $\{r_j\}_{j=1}^{S_0}$, and
\item inclusions $P_0 \subset P_1 \subset P'_1 \subset P'_0$ which are compatible with the assignments of successive elements of these sets to each element of $S_i$ and $S'_i$.

\end{enumerate}
\end{defin}
Note that, if $\underline{a}$ and $\underline{b}$ both lie in $\{M_-, M_0, M_+\}$, the datum $P$ is uniquely determined to be the empty set. Moreover, if they both lie in $\{M_-, M_0\}$, or both agree with $ M_+$, the set $\cD(\underline{a},\underline{b})$ is a singleton.
\begin{rem}
The fact that we treat elements of $\cP(H)$ and $\cP(f)$ differently is ultimately a consequence of the fact that we appeal to the smoothness of moduli spaces of gradient trajectories in Morse theory, but not in Floer theory. It is possible to have an entirely parallel discussion for the two cases by studying Morse trajectories with \emph{hypersurface constraints}.
\end{rem}

For the construction of Kuranishi flow categories we need to define
composition functors on the categories
$\cD(\underline{a},\underline{b})$; for this, we recall that the
addition of natural number lifts to a monoidal structure on the
category with morphisms given by order preserving injections. This
category is \emph{strictly monoidal} in the sense that the comparison
map $S_1 \otimes (S_2 \otimes S_3) \cong (S_1 \otimes S_2) \otimes
S_3$ is the identity map. This notion is related to the notion of a
strict $2$-category, discussed in Appendix~\ref{sec:2categories}. 

\begin{lem}
  There is a natural map
  \begin{equation}
    \cD(\underline{a},\underline{b}) \times \cD(\underline{b},\underline{c}) \to \cD(\underline{a},\underline{c}),    
  \end{equation}
 yielding a strict $\Pi$-equivariant $2$-category $\cD(f,H,f)$ with objects $\underline{\cP}(f,H,f)$, and $1$-morphisms the categories $\cD(\underline{a},\underline{b})$.
\end{lem}
\begin{proof}
 Given  $\underline{\alpha}_1= (S'_1 \subset S_1, P_1 \subset P'_1) \in   \cD(\underline{a},\underline{b})$ and $ \underline{\alpha}_2= (S'_2 \subset S_2, P_2\subset P'_2) \in  \cD(\underline{b},\underline{c})$, we define $\underline{\alpha}_1 \times \underline{\alpha}_2$ to consist of (i) the inclusion $S'_1 \amalg S'_2 \to S_1 \amalg S_2$, and (ii) the inclusion $P_1 \amalg \{ \underline{b}\} \amalg P_2 \subset P'_1 \amalg \{ \underline{b}\} \amalg P'_2 $ if $b$ lies in $\cP(H)$, and inclusion of disjoint unions of sets otherwise. The triple products $(\underline{\alpha}_1 \times \underline{\alpha}_2)\times \underline{\alpha}_3 $ and $\underline{\alpha}_1 \times (\underline{\alpha}_2 \times \underline{\alpha}_3)$ are evidently equal. Equivariance follows immediately from Lemma \ref{lem:Pi-action-D}, and the compatibility of these compositions with the $\Pi$ action.
\end{proof}

We now extend the construction of the previous section to a $2$-functor defined on $\cD(f,H,f) $, valued in the monoidal category of  equivariant submersions of stratified manifolds (see Definition \ref{def:equivariant_submersion_monoidal_category}). In order to describe the groups appearing in this functor, we associate to each natural number $k$ the symmetric group on $k$ letters. Given a pair $a < b$ in $\cP(f,H,f)$, and $\underline{\alpha} \in \cD(\underline{a},\underline{b})$, we denote by $G_{\underline{\alpha}}$ the product of symmetric groups indexed by the sequence $S_{\underline{\alpha}}$.

To define this $2$-functor,  we proceed by considering all possibilities for a pair of elements $\underline{a} < \underline{b}$ in $\underline{\cP}(f,H,f) $. Our discussion below will explain how to assign a smooth $G_{\underline{\alpha}}$-equivariant submersion of $\langle \underline{\cP}(\underline{a},\underline{b}) \rangle $-manifolds to each object $\underline{\alpha}$ of $\cD(\underline{a},\underline{b})$. The functoriality and multiplicativity of the construction is left to the reader, and is a straightforward exercise in applying the explicit constructions  of Section \ref{sec:constr-moduli-space} below:

We start by considering the cases in which the set $S$ is assumed to be empty: 

  if $\underline{a}$ and $\underline{b}$ agree, or $\underline{a} = M_-$  and $\underline{b} = M_0$,  we define   $\fcMR_{\underline{\alpha}}$ and $\Mbar^{\prime \bR}_{\underline{\alpha}}$ to both be a point.

  Next, if $\underline{a}$ and $\underline{b}$ both lie in $\cP(H)$, we use the definition from Section \ref{sec:abstr-moduli-spac}. There are four remaining cases to consider.

We leave the case $\underline{a} = M_-$ and $\underline{b} = M_+$ for the end; for the other cases, we need the following:
\begin{defin} \label{def:moduli_extra_marked_point_continuation}
For each finite set $r$, the moduli space
\begin{equation}
  \Mbar^{\chi}_{r \amalg \pm}  \subset  \Mbar_{\pm \amalg r}(\bP^1,\{0,\infty\}) 
\end{equation}
consists of the degree $1$ component of the moduli space of stable maps into $\bP^1$, with $z_-$ mapping to $0$ and $z_+$ to $\infty$.  
\end{defin}
Having restricted to degree $1$ maps, the map is a biholomorphism on some component, on which there is a distinguished point, which we denote $z_\chi$, which maps to $1 \in \bP^{1}$; we refer to this as a marked point even though it is allowed to agree with a node. Note that there are natural maps
\begin{equation}
   \Mbar^{\chi}_{r \amalg \pm} \to  \Mbar^{\chi}_{r' \amalg \pm}
\end{equation}
associated to forgetting marked points.

We define the underlying universal curve
\begin{equation}
  \Cbar^{\chi}_{r \amalg \pm}  \to   \Mbar^{\chi}_{r \amalg \pm} 
\end{equation}
by restricting the universal curve over $ \Mbar_{\pm \amalg r}(\bP^1) $.  
We shall use this moduli space to model moduli spaces with marked points needed to relate Morse and Floer theory. 

\begin{defin}[{\bf Case 1}] \label{def:Case1-definition-abstract-moduli}
  Given
  $\underline{\alpha} \in \cD(M_0,M_+)$, define
\begin{equation}
  \Mbar^{\prime \bR} _{\underline{\alpha}} \equiv \Mbar^{\chi}_{r_{S'} \amalg \pm}
\end{equation}
and set
\begin{equation}
  \fcMR_{\underline{\alpha}} \subset \Mbar^{\chi}_{r_{S} \amalg \pm}
\end{equation}
to consist of those curves whose domain is pulled back from $  \Mbar^{\chi}_{r_{S'} \amalg \pm}$ under the forgetful map.
\end{defin}
Note that, in this case, the sets $P$ and $P'$ are empty, so that there is only one level.

\begin{figure}[h]
  \centering
  \begin{subfigure}[t]{.25\linewidth}
    \centering
      \begin{tikzpicture}
   \coordinate[label= left:$-$] (p) at (-1,0);
   \coordinate[label= right:$+$] (r) at (1,0);
   \draw (-0.5,0) circle (.5);
 \draw[ultra thin, color=gray] (-1,0) arc (-180:0:.5 and .1);
 \draw[ultra thin, color=gray,dashed] (-1,0) arc (-180:-360:.5 and .1);
   \draw [fill,color=black] (p) circle (.05);
   \draw [fill,color=black] (r) circle (.05);
   \draw [fill,color=black] (0,0) circle (.05);
   \draw (0.5,0) circle (.5);
    \draw[ultra thin, color=gray] (0,0) arc (-180:0:.5 and .1);
    \draw[ultra thin, color=gray,dashed] (0,0) arc (-180:-360:.5 and .1);
     \draw[ultra thin, color=gray] (0,1) arc (-180:0:.5 and .1);
    \draw[ultra thin, color=gray,dashed] (0,1) arc (-180:-360:.5 and .1);
   \draw [fill,color=black] (.5,.5) circle (.05);
   \draw (0.5,1) circle (.5);

   \draw (.5,1) [fill,color=black] +(60:.5) circle (.05);
   \draw (.5,1) [fill,color=black] +(120:.5) circle (.05);
   \coordinate[label= above:$z^1$] (z31) at ($(.5,1)+(120:.5)$) ;
   \coordinate[label= above:$z^3$] (z32) at ($(.5,1)+(60:.5)$);

    \draw  [fill,color=black] (-.5,-.5) circle (.05);
    \draw [fill,color=black] (-.5,.5) circle (.05);
    \coordinate[label=  below:$z^4$] (z11) at (-.5,-.5) ;
   \coordinate[label= above:$z^2$] (z12) at (-.5,.5);   
 \end{tikzpicture}
\caption{Case 1: $|S|=1$}\label{fig:Continuation-constant}
\end{subfigure}
\begin{subfigure}[t]{.30\linewidth}
\centering   \begin{tikzpicture}

    \begin{scope}
   \coordinate[label= below:$p$] (p) at (0,-0.1);
   \coordinate[label= right:$+$] (q) at (1,0);
   \draw [fill,color=black] (q) circle (.05);
   \draw [fill,color=black] (0,0) circle (.05);
   \draw (0.5,0) circle (.5);
   \draw [fill,color=black] (.5,.5) circle (.05);
   \draw (0.5,1) circle (.5);
 \draw[ultra thin, color=gray] (0,0) arc (-180:0:.5 and .1);
 \draw[ultra thin, color=gray,dashed] (0,0) arc (-180:-360:.5 and .1);
   \draw[ultra thin, color=gray] (0,1) arc (-180:0:.5 and .1);
   \draw[ultra thin, color=gray,dashed] (0,1) arc (-180:-360:.5 and .1);
   
   \draw (.5,1) [fill,color=black] +(30:.5) circle (.05);
   \draw (.5,1) [fill,color=black] +(150:.5) circle (.05);
    \coordinate[label= right:$z_1^3$] (z13) at ($(.5,1)+(30:.5)$);
   \draw ($(.5,1)+(150:1)$) circle (.5);
 \begin{scope}[shift={($(.5,1)+(150:1)-(60:.5)$)}, rotate=60]
   \draw[ultra thin, color=gray] (0,0) arc (-180:0:.5 and .1);
   \draw[ultra thin, color=gray,dashed] (0,0) arc (-180:-360:.5 and .1);
   \draw (.5,0) [fill,color=black] +(30:.5) circle (.05);
   \draw (.5,0) [fill,color=black] +(150:.5) circle (.05);
     \coordinate[label= above:$z_1^1$] (z11) at ($(.5,0)+(30:.5)$);
   \coordinate[label= left:$z_1^2$] (z12) at  ($(.5,0)+(150:.5)$);
 \end{scope}
 
 \draw (-0.5,0) circle (.5);
\draw (-0.5,-.1) [thick] +(-.05,-.05)--+(.05,.05);
\draw (-0.5,-.1) [thick] +(.05,-.05)--+(-.05,.05);
 \coordinate[label= below:$z_\chi$] (zchi) at  (-0.5,-.1);
\draw [fill,color=black] (-1,0) circle (.05);
 \draw (-1.5,0) circle (.5);
 \coordinate[label= above:$z_2^1$] (z21) at  (-1.5,.5);
 
  \draw[ultra thin, color=gray] (-1,0) arc (-180:0:.5 and .1);
  \draw[ultra thin, color=gray,dashed] (-1,0) arc (-180:-360:.5 and .1);
 \draw[ultra thin, color=gray] (-2,0) arc (-180:0:.5 and .1);
 \draw[ultra thin, color=gray,dashed] (-2,0) arc (-180:-360:.5 and .1);

   \draw [fill,color=black] (-1.5,.5) circle (.05);
   \coordinate[label= left:$-$] (p) at (-2,-0);
   \draw [fill,color=black] (p) circle (.05);

   \end{scope}
  \end{tikzpicture}
\caption{ Case 2: $|S|=2$. There are two levels. The inverse image of $1 \in \bP^1$ is labelled by a cross.}\label{fig:continuation-to-HF}
\end{subfigure}
\begin{subfigure}[t]{.30\linewidth}
  \centering
   \begin{tikzpicture}

    \begin{scope}
   \coordinate[label= below:$q$] (q) at (0,-0.2);
   \coordinate[label= right:$+$] (+) at (1,0);
   \draw [fill,color=black] (+) circle (.05);
   \draw [fill,color=black] (0,0) circle (.05);
   \draw (0.5,0) circle (.5);
   \draw [fill,color=black] (.5,.5) circle (.05);
  \coordinate[label= above:$z_1^1$] (z11) at  (.5,.5);

   \draw[ultra thin, color=gray] (0,0) arc (-180:0:.5 and .1);
 \draw[ultra thin, color=gray,dashed] (0,0) arc (-180:-360:.5 and .1);
\draw (0.5,-.1) [thick] +(-.05,-.05)--+(.05,.05);
\draw (0.5,-.1) [thick] +(.05,-.05)--+(-.05,.05);
 \coordinate[label= below:$z_\chi$] (zchi) at  (0.5,-.1); 
\draw [fill,color=black] (-1,0) circle (.05);

   \coordinate[label= below:$p$] (p) at (-1,-0.2);
\draw (-.5,0) circle (.5);
\draw [fill,color=black] (-.5,.5) circle (.05);
 \coordinate[label= above:$z_2^1$] (z21) at  (-.5,.5);
  \draw[ultra thin, color=gray] (-1,0) arc (-180:0:.5 and .1);
  \draw[ultra thin, color=gray,dashed] (-1,0) arc (-180:-360:.5 and .1);

   \draw [fill,color=black] (-1.5,.5) circle (.05);
 \draw (-1.5,0) circle (.5);
 \draw [fill,color=black] (-1.5,.5) circle (.05);
  \coordinate[label= above:$z_3^1$] (z31) at  (-1.5,.5);
 \draw[ultra thin, color=gray] (-2,0) arc (-180:0:.5 and .1);
 \draw[ultra thin, color=gray,dashed] (-2,0) arc (-180:-360:.5 and .1);

   \coordinate[label= left:$-$] (-) at (-2,-0);
   \draw [fill,color=black] (-) circle (.05);

   \end{scope}
  \end{tikzpicture}
\caption{Case 3: $|S|=3$. There are three levels, which are separated by the nodes labelled $p$ and $q$.}\label{fig:continuation-from-HF}
\end{subfigure}
  \caption{A graphical representation of elements of $ \Mbar_{\underline{\alpha}}$ (left and right) for Cases 1, 2, and 3. The point labelled with $z_\chi$ may agree with any node or marked point.} \label{fig:Moduli-space-continutation-1}
\end{figure}
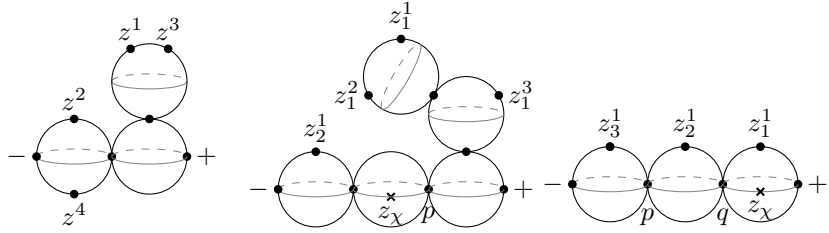

We now consider the next two cases. If $\underline{a} = M_-$ and $ \underline{b} = p \in \cP(H)$, then for each tree labelling a curve in $  \Mbar^{\chi}_{r_{S'} \amalg \pm}$, we define a Floer edge to be an edge along the path between the vertex carrying the marked point $\chi$ and the positive end. If $a = p  \in\cP(H)$ and $b = M_+$, we define a Floer edge to be an edge along the path between the vertex carrying the marked point $\chi$ and the negative end. In either case, we define a level to be a connected component of the complement of the Floer edges. The levels are again ordered.
\begin{defin}[\bf{Cases 2 and 3}]
  \label{def:case2-3-abstract-moduli}
If $ p$ lies $\cP(H)$, then for each element $\underline{\alpha}$ in $\cD(M_-,p)$ or in $\cD(p, M_+)$ we define
\begin{equation}
  \Mbar'_{\underline{\alpha}} \subset  \Mbar^{\chi}_{r_{S'} \amalg \pm} 
\end{equation}
to consist of those curves such that:
\begin{enumerate}
\item all points labelled by the same element of $S'$ lie in the same level, and 
\item the induced map from $S'$ to levels preserves order, and factors through an order preserving surjective map from the set of successive elements of $\{ - \} \amalg  P'  \amalg \{+\}$ to the set of levels, with $\{ - , + \}$ corresponding to $(M_-, p)$ in the first case, and to $(p, M_+)$ in the second.
\end{enumerate}
 We define
\begin{equation}
  \Mbar_{\underline{\alpha}} \subset  \Mbar^{\chi}_{r_{S} \amalg \pm}
 \end{equation}
to consist of curves whose domain is pulled back from $ \Mbar'_{\underline{\alpha}}$ under the forgetful map, such that all points labelled by the same element of $S$ lie in the same level, and such that the induced map from $S$ to the set of levels factors through successive elements of $\{ - \} \amalg   P  \amalg \{+\}$.  
\end{defin}
Note that $\Mbar_{\underline{\alpha}}$ and $ \Mbar'_{\underline{\alpha}}$ are respectively stratified by $P$ and $P'$.

We define a Floer vertex to be a vertex with two adjacent Floer edges, and define
\begin{equation}
   \fcMR_{\underline{\alpha}} \to   \Mbar_{\underline{\alpha}}  \textrm{ and }   \Mbar^{\prime \bR}_{\underline{\alpha}} \to   \Mbar'_{\underline{\alpha}}
\end{equation}
to be the spaces of curves equipped with compatible asymptotic markers at the two ends of each component corresponding to a Floer vertex of the underlying tree (see Section \ref{sec:constr-moduli-space}). The spaces $ \fcMR_{\underline{\alpha}}  $ and $  \Mbar^{\prime \bR}_{\underline{\alpha}} $ are  respectively  smooth $\langle P \rangle$ and  $\langle P' \rangle$ manifolds, and in particular
both are $\langle \underline{\cP}(\underline{a},\underline{b}) \rangle$-manifolds (the boundary stratum associated to each element of $\underline{\cP}(\underline{a},\underline{b})$ that does not lie in $P$ is empty). There is a natural forgetful map
\begin{equation}
 \fcMR_{\underline{\alpha}}  \to \Mbar^{\prime \bR}_{\underline{\alpha}}  
\end{equation}
which is a stratified submersion (see Lemma \ref{lem:submersion-floer-cylinder-forgetful-map} below).

For the final case, we must introduce new moduli spaces of Deligne-Mumford type: given a set $r$, we consider a moduli space
\begin{equation}
 \Mbar^{\chi_\pm}_{r \amalg \pm}  \subset  \Mbar_{\pm \amalg r}( (\bP^1 , \{0, \infty\})^{2} ),
\end{equation}
consisting of stable curves in $\bP^1 \times \bP^1 $  of degree $(1,1)$, mapping $z_-$ to $(0,0)$ and $z_+$ to $(\infty,\infty)$, and which is the closure of the set of maps such that
\begin{equation}
  \parbox{31em}{the second factor is given by composing the first factor by multiplication by a complex number $\lambda$ with $1 \leq |\lambda|$.}
\end{equation}
When $r = 0$, we may identify $\Mbar^{\chi_\pm}_{\pm}$ with a closed unit disc $D^2$.

There is a natural forgetful map
\begin{equation}
   \Mbar^{\chi_\pm}_{r \amalg \pm}  \to  \Mbar^{\chi_\pm}_{r' \amalg \pm} 
\end{equation}
as before, and in particular a map to $D^2$ obtained by forgetting all marked points except $z_{\pm}$. The fibres of this map to $D^2$ over the point $0$ and $1$ 
\begin{align}
  \partial^0   \Mbar^{\chi_\pm}_{r \amalg \pm} \subset   \Mbar^{\chi_\pm}_{r \amalg \pm} \supset \partial^1   \Mbar^{\chi_\pm}_{r \amalg \pm}
\end{align}
respectively consist of maps which factor into the union of two lines, and those which pass through the point $(1,1)$.

\begin{lem}
  The  moduli space $\Mbar^{\chi_\pm}_{r \amalg \pm} $ is a smooth manifold with boundary $ \partial^1   \Mbar^{\chi_\pm}_{r \amalg \pm}$ which is canonically identified with the product:
  \begin{equation}
  \Mbar^{\chi}_{r \amalg \pm} \times S^1.
  \end{equation}

\end{lem}
\begin{proof}
  We note that multiplication by $S^1$ acts freely on the set of pairs of maps which differ by multiplication by $\lambda$ with modulus $1$. Since every map from $\bP^1$ to $\bP^1 \times \bP^1$ in degree $(1,1)$, passing through the points $(0,0)$, $(1,1)$ and $(\infty,\infty)$ factors through the diagonal, we conclude that those for which $\lambda = 1$ exactly correspond to curves in $ \Mbar^{\chi}_{r \amalg \pm}$. 
\end{proof}
We define the underlying universal curve
\begin{equation}
  \Cbar^{\chi_\pm }_{r \amalg \pm}  \to   \Mbar^{\chi_\pm}_{r \amalg \pm} 
\end{equation}
as before by pulling back the universal curve over the moduli spaces of maps with target $\bP^1 \times \bP^1$.

\begin{figure}[h]
  \centering
\begin{subfigure}[b]{.45\linewidth}
\centering   \begin{tikzpicture}
      \coordinate[label= left:$+$] (-) at (-1,0);
      \coordinate[label= right:$+$] (+) at (1,0);
   \draw [fill,color=black] (+) circle (.05);
   \draw [fill,color=black] (-) circle (.05);
   \draw (0,0) circle (1);
 \draw[ultra thin, color=gray] (-1,0) arc (-180:0:1 and .2);
 \draw[ultra thin, color=gray,dashed] (-1,0) arc (-180:-360:1 and .2);

   \draw [fill,color=black] (30:1) circle (.05);
   \draw [fill,color=black] (150:1) circle (.05);
    \coordinate[label= right:$z_1^2$] (z12) at (30:1);
    \coordinate[label= above:$z_1^1$] (z11) at (150:1);

     \coordinate[label= below:$z_{\chi_-}$] (zchi-) at  (-120: 1 and .2);
\draw (zchi-) [thick] +(-.05,-.05)--+(.05,.05);
\draw (zchi-) [thick] +(.05,-.05)--+(-.05,.05);

\coordinate[label= below:$z_{\chi_+}$] (zchi+) at  (-60: 1 and .2);
\draw (zchi+) [thick] +(-.05,-.05)--+(.05,.05);
\draw (zchi+) [thick] +(.05,-.05)--+(-.05,.05);
  \end{tikzpicture}
\caption{A one level curve}\label{fig:homotopy_composition_continuation-1}
\end{subfigure}
\begin{subfigure}[b]{.45\linewidth}
  \centering
   \begin{tikzpicture}

    \begin{scope}
   \coordinate[label= below:$q$] (q) at (0,-0.2);
   \coordinate[label= right:$+$] (+) at (1,0);
   \draw [fill,color=black] (+) circle (.05);
   \draw [fill,color=black] (0,0) circle (.05);
   \draw (0.5,0) circle (.5);
   \draw [fill,color=black] (.5,.5) circle (.05);
  \coordinate[label= above:$z_1^1$] (z11) at  (.5,.5);

   \draw[ultra thin, color=gray] (0,0) arc (-180:0:.5 and .1);
 \draw[ultra thin, color=gray,dashed] (0,0) arc (-180:-360:.5 and .1);
\draw (0.5,-.1) [thick] +(-.05,-.05)--+(.05,.05);
\draw (0.5,-.1) [thick] +(.05,-.05)--+(-.05,.05);
 \coordinate[label= below:$z_{\chi_+}$] (zchi+) at  (0.5,-.1); 
\draw [fill,color=black] (-1,0) circle (.05);

   \coordinate[label= below:$p$] (p) at (-1,-0.2);
\draw (-.5,0) circle (.5);
\draw [fill,color=black] (-.5,.5) circle (.05);
 \coordinate[label= above:$z_2^1$] (z21) at  (-.5,.5);
  \draw[ultra thin, color=gray] (-1,0) arc (-180:0:.5 and .1);
  \draw[ultra thin, color=gray,dashed] (-1,0) arc (-180:-360:.5 and .1);

   \draw [fill,color=black] (-1.5,.5) circle (.05);
 \draw (-1.5,0) circle (.5);
 \draw [fill,color=black] (-1.5,.5) circle (.05);
  \coordinate[label= above:$z_3^1$] (z31) at  (-1.5,.5);
 \draw[ultra thin, color=gray] (-2,0) arc (-180:0:.5 and .1);
 \draw[ultra thin, color=gray,dashed] (-2,0) arc (-180:-360:.5 and .1);
  
  \draw (-1.5,-.1) [thick] +(-.05,-.05)--+(.05,.05);
\draw (-1.5,-.1) [thick] +(.05,-.05)--+(-.05,.05);
 \coordinate[label= below:$z_{\chi_-}$] (zchi-) at  (-1.5,-.1); 
   \coordinate[label= left:$-$] (-) at (-2,-0);
   \draw [fill,color=black] (-) circle (.05);

   \end{scope}
  \end{tikzpicture}
\caption{A 3 level curve with $|S|=3$}\label{fig:continuation-from-HF-2}
\end{subfigure}
  \caption{A graphical representation of two elements of $ \Mbar_{\underline{\alpha}}$ (left and right) for Case 4.} \label{fig:Moduli-space-continutation-2}
\end{figure}
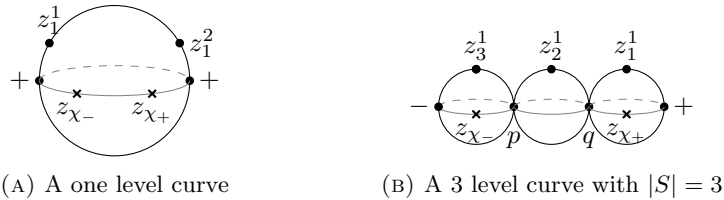

We finally return to the case $a = M_-$ and $b= M_+$. Given a tree labelling a curve in $  \Mbar^{\chi_\pm}_{r \amalg \pm}$, we define a Floer edge to be an edge along the path between the vertex carrying the marked point $\chi_-$ and the vertex carrying $\chi_+$, yielding a decomposition into levels as before. 
 \begin{defin}[{\bf Case 4}]
  \label{def:Case4_abstract}
For each element $\underline{\alpha} \in \cD(M_-,M_+)$ we define
\begin{equation} \label{eq:inclusion_mbar_alpha-to-moduli-with-two-marked-points}
  \Mbar'_{\underline{\alpha}} \subset    \Mbar^{\chi_\pm}_{r_{S'} \amalg \pm} 
\end{equation}
such that (i) all points labelled by the same element of $S'$ lie in the same level, and (ii) the induced map from $S'$ to levels factors through an order preserving surjective map from  the set of successive elements of $\{-\} \amalg P' \amalg \{+\} $ to the set of levels.

We let
\begin{equation}\label{eq:inclusion_mbar_alpha-to-moduli-with-two-marked-points-prime}
   \Mbar_{\underline{\alpha}} \subset  \Mbar^{\chi}_{r_{S} \amalg \pm}
\end{equation}
denote the set of curves whose domain is pulled back from $\Mbar'_{\underline{\alpha}} $ under the forgetful map, imposing again the condition above on the domain, and in particular that the map from $S$ to the set of levels factors through  successive elements of $\{-\} \amalg P \amalg \{+\} $.
\end{defin}

We define a Floer vertex to be a vertex with two adjacent Floer edges, and in the next section, we shall define
\begin{equation}
  \fcMR_{\underline{\alpha}} \to   \Mbar_{\underline{\alpha}}  \textrm{ and }   \Mbar^{\prime \bR}_{\underline{\alpha}} \to   \Mbar'_{\underline{\alpha}}
\end{equation}
to be the spaces of curves in $\Mbar_{\underline{\alpha}}$ equipped with compatible asymptotic markers at the two ends of each component corresponding to a Floer vertex of the underlying tree. These are  again $\langle \underline{\cP}(\underline{a},\underline{b}) \rangle$-manifolds, equipped with an equivariant smooth submersion $  \fcMR_{\underline{\alpha}} \to \Mbar^{\prime \bR}_{\underline{\alpha}} $.

\subsection{Construction of abstract moduli spaces}
\label{sec:constr-moduli-space}
\label{sec:funct-abstr-moduli}

In this section, we give a precise construction of the abstract moduli spaces of holomorphic curves needed for Hamiltonian Floer theory. The basic idea is that the moduli spaces we need have the additional datum of a real line connecting the marked points $\pm$, and that we shall obtain them from moduli spaces of maps to products of $\bP^1$ on which we impose a reality condition. Before proceeding further, we remind the reader that we have defined a notion of Floer component for all moduli spaces introduced in Sections \ref{sec:abstr-moduli-spac} (Floer trajectories) and \ref{sec:abstr-moduli-spac-1} (Continuation maps). Given a pair of elements $\underline{a}$ and $\underline{b}$ of $\underline{\cP}(f,H,f)$, and $\underline{\alpha} \in \cD(\underline{a},\underline{b})$, we shall define several auxiliary moduli spaces before constructing $\fcMR_{\underline{\alpha}}$ and $\Mbar^{\prime \bR}_{\underline{\alpha}}$.

To begin, in the case (i) $\underline{a} = \underline{b}$, (ii) $\underline{a} = M_-$ and $\underline{b} = M_0$, or (iii) $\underline{a} = M_0$ and $\underline{b} = M_+$, we have already defined these moduli spaces (see in particular Definition \ref{def:Case1-definition-abstract-moduli}).

We must therefore discuss the cases of Floer trajectories, and the remaining cases discussed in Section  \ref{sec:abstr-moduli-spac-1}. To unify the discussion, let $\{\chi_a^b\}$ denote (i) the empty set in the setting of Floer trajectories, (ii) the singleton $\chi_-$ or $\chi_+$ in Cases 2 and 3 of Section  \ref{sec:abstr-moduli-spac-1}, or (ii) the pair $(\chi_-,\chi_+)$ in Case 4.

Before proceeding further, given $\underline{\alpha} \in \cD(\underline{a},\underline{b})$, we define $S^H$ and $S^{\prime H}$ to be the subsets of $S$ and $S'$ mapping to  successive elements of $P \cap \cP(H)$ and $P' \cap \cP(H)$. We have a natural ordering on the union of these sets with $\{\chi_a^b\}$, given by the ordering on $S^H$ induced by its inclusion in $S$, and the requirement that
\begin{equation} \label{eq:ordering_continuation_and_successive_elements}
\parbox{30em}{  $\chi_-$ is the initial, and $\chi_+$ the terminal element of $S^H \amalg \{\chi_a^b\}$,} 
\end{equation}
and similarly for $S^{\prime H}$.
\begin{rem}
  For many purposes, it is sufficient to replace $S^H$ by its image in the set of successive elements of $P$. However, note that we have a forgetful map $S^H \to S^{\prime H}$, but no such map when considering successive elements of $P$ and $P'$. This makes our choice easier to use.
\end{rem}

Consider the moduli space 
\begin{equation}
  \Mbar_{r_{S} \amalg \pm}((\bP^{1}, \{ 0,\infty \})^{\{\chi_a^b\} \amalg S^H})    
\end{equation}
of stable maps from a genus $0$ Riemann surface with points marked by $r_{S} \amalg \pm $ to $(\bP^{1})^{\{\chi_a^b\} \amalg S^H}$, mapping $z_-$ to $(0, \ldots, 0)$ and $z_+$ to $(\infty, \ldots, \infty) $. We have a natural forgetful map
\begin{equation}
   \Mbar_{r_{S} \amalg \pm}((\bP^{1}, \{0,\infty\})^{\{\chi_a^b\} \amalg S^H})     \to  \Mbar_{r_{S'} \amalg \pm} ((\bP^{1}, \{0,\infty\})^{\{\chi_a^b\} \amalg S^H}).
  \end{equation}

\begin{defin} \label{def:moduli_extra_maps_S}
  For each pair $\underline{a}$ and $\underline{b}$ of elements of $\underline{\cP}(f,H,f)$, we define the moduli spaces
  \begin{align}
    \Mbar^{S^H}_{\underline{\alpha}} & \subset    \Mbar_{r_{S} \amalg \pm}((\bP^{1}, \{0,\infty\})^{\{\chi_a^b\} \amalg S^H}) 
    \\
\Mbar^{S^{\prime H}}_{\underline{\alpha}} & \subset  \Mbar_{r_{S'} \amalg \pm}((\bP^{1}, \{0,\infty\})^{\{\chi_a^b\} \amalg S^{\prime H}})  
\end{align}
to consist of maps satisfying the following properties:
\begin{align} \label{eq:degree_1_on_each_factor_S}
  & \parbox{30em}{ The projection to each component of the target has degree $1$.} 
    \\ \label{eq:ordering_components_agrees_labels}
& \parbox{30em}{The ordering along the path from $z_-$ to $z_+$ of the components which map non-trivially to each factor of the target respects the ordering in Equation \eqref{eq:ordering_continuation_and_successive_elements}. Moreover, any two components associated to factors of $S^H$ mapping to the same pair of sucessive elements of $P$ agree.}  \\\label{eq:tie_projection_to_S_factor_to_marked_points}
& \parbox{30em}{For each $i \in S^H$ (or $S^{\prime H}$), the image of the points marked by elements of $r_i$ in the factor of $(\bP^{1})^{S^H}$ corresponding to $i$ lies in $\bP^{1} \setminus \{0,\infty\}$ .} \\\label{eq:tie_projection_to_S_factor_to_marked_points-minus}
& \parbox{30em}{ For each $i \in S$ (or $S'$) mapping to a pair $(M_-,p)$ of successive elements of $P$ (or $P'$), the image of the points marked by elements of $r_i$ in the factor  $(\bP^{1})^{\{\chi_-\}}$ lies in $\bP^{1} \setminus \{\infty\}$ .} \\\label{eq:tie_projection_to_S_factor_to_marked_points-plus}
& \parbox{30em}{ For each $i \in S$ (or $S'$) mapping to a pair $(p,M_+)$ of successive elements of $P$ (or $P'$), the image of the points marked by elements of $r_i$ in the factor  $(\bP^{1})^{\{\chi_+\}}$ lies in $\bP^{1} \setminus \{0\}$ .}
\\\label{eq:forgetful_map_does_not_collapse_with_marked_s}
& \parbox{30em}{The domain of each curve in $ \Mbar^{S^H}_{\underline{\alpha}}$ is pulled back from its image in $\Mbar^{S^{\prime H}}_{\underline{\alpha}} $.} 
\end{align}
\end{defin}

To understand the above moduli spaces, we consider the following commutative diagram of forgetful map
\begin{equation} \label{eq:diagram_forgetful_maps_with_r_and_S}
  \begin{tikzcd}
    \Mbar_{r_{S} \amalg \pm}((\bP^{1}, \{0,\infty\})^{\{\chi_a^b\} \amalg S^H})    \ar[d] \ar[r] &   \Mbar_{r_{S'} \amalg \pm} ((\bP^{1}, \{0,\infty\})^{\{\chi_a^b\} \amalg S^H}) \ar[d] \\
    \Mbar_{r_S \amalg \pm} ((\bP^{1}, \{0,\infty\})^{\{\chi_a^b\}}) \ar[r] &   \Mbar_{r_{S'} \amalg \pm} ((\bP^{1}, \{0,\infty\})^{\{\chi_a^b\}})
  \end{tikzcd}
  \end{equation}
Recalling that we have natural inclusions
  \begin{align}
    \Mbar_{\underline{\alpha}}  & \to \Mbar_{r_S \amalg \pm} ((\bP^{1}, \{0,\infty\})^{\{\chi_a^b\}}) \\
     \Mbar'_{\underline{\alpha}}  & \to \Mbar_{r_{S'} \amalg \pm} ((\bP^{1}, \{0,\infty\})^{\{\chi_a^b\}}) 
  \end{align}
we can relate the moduli spaces from Definition \ref{def:moduli_extra_maps_S} with those from previous sections:

\begin{lem}
  The moduli spaces $\Mbar^{S^H}_{\underline{\alpha}}$ and $\Mbar^{S^{\prime H}}_{\underline{\alpha}}$ are smooth manifolds with corners, equipped with smooth $G_{\underline{\alpha}}$ actions, and Diagram \eqref{eq:diagram_forgetful_maps_with_r_and_S} induces a commutative diagram
  \begin{equation}
    \begin{tikzcd}
      \Mbar^{S^H}_{\underline{\alpha}} \ar[r] \ar[d] & \Mbar^{S^{\prime H}}_{\underline{\alpha}} \ar[d] \\
     \Mbar_{\underline{\alpha}} \ar[r] & \Mbar'_{\underline{\alpha}} 
    \end{tikzcd}
  \end{equation}
in which  all arrows are equivariant submersions.
\end{lem}
\begin{proof}
  Condition \eqref{eq:tie_projection_to_S_factor_to_marked_points} gives a decomposition of the domains of curves in $   \Mbar^{S^H}_{\underline{\alpha}}$ and  $\Mbar^{S^{\prime H}}_{\underline{\alpha}}$ into levels, which by Conditions \eqref{eq:tie_projection_to_S_factor_to_marked_points}--\eqref{eq:tie_projection_to_S_factor_to_marked_points-plus}  agrees with the decomposition into levels of their images under the map which forgets the factors labelled by elements of $S^H$ or $S^{\prime H}$.  This shows that the image of $\Mbar^{S^{\prime H}}_{\underline{\alpha}} $ under the forgetful map lies in $    \Mbar'_{\underline{\alpha}}$, and Condition \eqref{eq:forgetful_map_does_not_collapse_with_marked_s} implies the same for $\Mbar^{S^H}_{\underline{\alpha}} $.  Moreover, since no component is collapsed by any projection map, they are all submersions.
\end{proof}
If we instead forget the marked points $r_S$, we obtain a map to the moduli space $ \Mbar_{ \pm}((\bP^{1}, \{0,\infty\})^{S^H \amalg \{\chi_a^b\}}) $. We introduce the subset
\begin{equation}
  \Rbar^{\bC}_{\alpha}  \subset \Mbar_{ \pm}((\bP^{1}, \{0, \infty\})^{S^H \amalg \{\chi_a^b\}}) 
\end{equation}

consisting of stable maps of degree $1$ in each component, which map $z_\pm$ to $(0, \ldots, 0)$ and $ (\infty, \ldots, \infty)$, and such that
\begin{equation}
  \parbox{30em}{the ordering of components satisfies Condition \eqref{eq:ordering_components_agrees_labels}. Moreover, if $\{\chi_a^b\} = \{\chi_-, \chi_+\}$, then the projection to $ \Mbar_{\pm}\left((\bP^1, \{0,\infty\})^{\{\chi_-, \chi_+\}} \right)$ maps to $ \cM^{\chi_\pm}_{\pm} $.}
\end{equation}

We note that these moduli spaces have straightforward explicit descriptions:
\begin{lem} \label{lem:moduli_space_affine_space}
  If $\{\chi_a^b\}$ is empty, or consists of a singleton, the moduli space $\Rbar^\bC_\alpha$ is biholomorphic to a product $(\bC^*)^{|S^H| - |P \cap \cP(H)|} \times \bC^{|P \cap \cP(H)| -1}$. If $\{\chi_a^b\} = \{ \chi_-, \chi_+\}$, it is biholomophic to the product of $(\bC^*)^{|S^H| - |P \cap \cP(H)|} $ with the subset $| \prod_i z_i| \leq 1$ in a complex affine space of dimension $ |P| -1$. Moreover, the forgetful map
  \begin{equation}
    \Rbar^\bC_\alpha \to \Rbar^{\prime \bC}_\alpha
  \end{equation}
 is a submersion, and there is a canonical short exact sequence
  \begin{equation}
  0 \to  \bC^{S^H \setminus S^{\prime H}} \to T   \Rbar^\bC_\alpha \to T \Rbar^{\prime \bC}_\alpha \to 0
\end{equation}
of tangent spaces.
  \qed
\end{lem}

For the statement of the next result, recall that $\Mbar^{S^H}_{\underline{\alpha}}  $ and $\Mbar^{S^H}_{\underline{\alpha}}$ have natural actions by the product $G_{\underline{\alpha}}$ of permutation groups. We equip $\Rbar^\bC_\alpha $ and $\Rbar^{\prime \bC}_\alpha$ with the trivial action of this group.
        \begin{lem}
          The following diagram of forgetful maps commutes, and all arrows are $G_{\underline{\alpha}}$ equivariant smooth submersions:
          \begin{equation}
            \begin{tikzcd}
              \Mbar^{S^H}_{\underline{\alpha}} \ar[d] \ar[r] & \Mbar^{S^{\prime H}}_{\underline{\alpha}} \ar[d] \\
              \Rbar^\bC_\alpha \ar[r] & \Rbar^{\prime \bC}_\alpha  .
            \end{tikzcd}
          \end{equation}
        \end{lem}
        \begin{proof} The key point is that, for any target $X$, the critical points of a forgetful map $\Mbar_{r}(X) \to \Mbar_{r'}(X) $ consist of points lying on a stratum where there is a node in the target whose inverse image in the source is non-trivial (i.e. there is at least one component that is collapsed to this node). Condition \ref{eq:tie_projection_to_S_factor_to_marked_points} precludes this.
        \end{proof}

        At this stage, we 
        impose the desired reality condition
                by considering the subspaces
        \begin{align}
        \Rbar_{\underline{\alpha}}   & \subset  \Rbar^\bC_\alpha \\
           \Rbar'_{\underline{\alpha}} & \subset  \Rbar^{\prime \bC}_\alpha
        \end{align}
        consisting of maps from a pre-stable curve $\Sigma$ to $(\bP^{1})^{S^H \amalg \{\chi_a^b\}}$ or $(\bP^1)^{ S^{\prime H} \amalg \{\chi_a^b\}}$ satisfying the following property:
        \begin{align}
          \label{eq:maps_differ_positive_rescaling}
          & \parbox{30em}{on any component of $\Sigma$, all non-constant factors of the map to $(\bP^{1})^{S^H \amalg \{\chi_a^b\}}$ differ by multiplication by a positive real number.}  
        \end{align}
        
For the proof of the next result, which is the real analogue of Lemma \ref{lem:moduli_space_affine_space}, it is convenient to note that  every domain is a chain of rational curves with one end carrying the marked point $z_-$ and the other the marked point $z_+$. 
        \begin{lem} \label{lem:real_moduli_spaces_no_marked_points}
            The moduli space $ \Rbar_\alpha$ is a smooth manifold with corners, whose interior is diffeomorphic to a Euclidean space. The forgetful map to $\Rbar^{\prime \bC}_\alpha$ defines a submersion
            \begin{equation}
                    \Rbar_\alpha \to \Rbar'_{\underline{\alpha}},
                  \end{equation}
                  with fibre the product $\bR^{S^H \setminus S^{\prime H}}$ (canonically up to translation).
          \end{lem}
          \begin{proof}
            Condition \eqref{eq:maps_differ_positive_rescaling} implies that each such stratum of $\Rbar_\alpha$ is a product of the quotient by the diagonal action of $\bR$ on a product of euclidean spaces, and in particular, that the interior stratum is $ \bR^{S+ |\{\chi_a^b\}|  -1}$. The smooth structure is inherited from the ambient space. 

To prove that the forgetful map has the desired target, recall from Definition \ref{def:index-set-moduli-space-cylinders-marked-points} that we have required for the map from $S'$ to the  successive elements of $\{a\} \amalg P' \amalg \{b\}$ to be surjective, which implies that the map from $S'$ to the successive elements of $\{a\} \amalg P' \amalg \{b\}$ is also surjective, hence that the forgetful map from  $ \Rbar_\alpha$ to $\Rbar^{\prime \bC}_\alpha$ does not collapse any component. From this, it follows that the image lies in $\Rbar'_{\underline{\alpha}}$, with fibre a product of euclidean spaces indexed by $S^H \setminus S^{\prime H}$.
          \end{proof}

Consider  the fibre products
        \begin{align}
          \Mbar^{S^H,\bR}_{\underline{\alpha}} & \equiv \Mbar^{S^H}_{\underline{\alpha}} \times_{ \Rbar^\bC_\alpha }  \Rbar_\alpha\\
          \Mbar^{S^{\prime H},\bR}_{\underline{\alpha}} & \equiv \Mbar^{S^{\prime H}}_{\underline{\alpha}} \times_{ \Rbar^{\prime \bC}_\alpha  }  \Rbar'_{\underline{\alpha}},
      \end{align}
       which have real dimension $2r_{S} + S^H + |\{\chi_a^b\} | - 1 $ and $2 r_{S'} + S^{\prime H} + |\{\chi_a^b\} | - 1 $.
      The notation indicates that we have, in addition to the marked points labelled by the union of the sets labelled by elements of $S $ or $S'$, an additional $S^H\amalg \{\chi_a^b\}$ or $S^{\prime H}\amalg \{\chi_a^b\}$ special points which are \emph{real}, obtained as the inverse image of $1 \in \bP^1$ in the corresponding factor.

      \begin{example}
If $r_S =  r_{S'} = 1$ is a point, and in the absence of additional marked points corresponding to continuation, the moduli spaces $  \Rbar^\bC_\alpha = \Rbar_\alpha$ are given by a point, so that $\Mbar^{S_H,\bR}_{\underline{\alpha}}$ is a copy of $\bC^*$.

If $S_H = S'_H$ has two elements, and $|r_{S}| = 2$, then $\Mbar^{S_H}_{\underline{\alpha}}$ is a $6$ dimensional manifold. The moduli space $\Rbar_\alpha$ is an interval embedded in $   \Rbar^\bC_\alpha   \cong \bP^1$. Thus, $\Mbar^{S_H,\bR}_{\underline{\alpha}}$ has codimension $1$ in $\Mbar^{S_H,\bR}_{\underline{\alpha}}$, so that its real dimension is $5$.
      \end{example}

      We now consider the forgetful map in the other direction:  there is a natural action of $\bR^{S^H}$ and $\bR^{S^{\prime H}}$ on $ \Mbar^{S^H,\bR}_{\underline{\alpha}}$ and $\Mbar^{S^{\prime H},\bR}_{\underline{\alpha}} $ by composing the given map to $(\bP^1)^{S^H} $ and  $(\bP^1)^{S^{\prime H}} $ with the action of $\bR$ on each factor of the target, given by
      \begin{equation}
        (t,z) \mapsto e^t z.
      \end{equation}
\begin{defin}
The  \emph{moduli spaces of abstract Floer cylinders,}  $\fcMR_{\underline{\alpha}}$ and $ \Mbar^{\prime \bR}_{\underline{\alpha}} $ are the quotients of $\Mbar^{S^H,\bR}_{\underline{\alpha}} $ and $\Mbar^{S^{\prime H},\bR}_{\underline{\alpha}} $ by the action of $\bR^{S^H}$ and $\bR^{S^{\prime H}}$.
\end{defin}

      \begin{prop} \label{prop:moduli-space-compact-smooth-Haussdorff}
        The spaces  $\fcMR_{\underline{\alpha}}$ and $ \Mbar^{\prime \bR}_{\underline{\alpha}}$ are smooth (Hausdorff) manifolds with corners, equipped with a $G_{\underline{\alpha}}$ action. There is an equivariant commutative diagram
        \begin{equation}
          \begin{tikzcd}
            \Mbar^{S^H,\bR}_{\underline{\alpha}} \ar[r] \ar[d] & \Mbar^{S^{\prime H},\bR}_{\underline{\alpha}} \ar[d] \\
            \fcMR_{\underline{\alpha}} \ar[r] & \Mbar^{\prime \bR}_{\underline{\alpha}}
          \end{tikzcd}
        \end{equation}
in which all arrows are submersions, and the vertical maps have a contractible space of sections.

        The kernel $T^{\fcMR}  \Mbar^{S^H,\bR}_{\underline{\alpha}} $ of the projection map $T \Mbar^{S^H,\bR}_{\underline{\alpha}} \to T \fcMR_{\underline{\alpha}}$ is naturally isomorphic to $\bR^{S^H}$, the kernel $T^{\Rbar}  \Mbar^{S^H,\bR}_{\underline{\alpha}}$ of the projection map $ T \Mbar^{S^H,\bR}_{\underline{\alpha}} \to T \Rbar_\alpha $ is equipped with a natural complex structure, and there is a canonical short exact sequence
        \begin{equation} \label{eq:moduli-space-with-real-points-two-kernels-span}
          \bR \to   T^{\fcMR}  \Mbar^{S^H,\bR}_{\underline{\alpha}} \oplus   T^{\Rbar}  \Mbar^{S^H,\bR}_{\underline{\alpha}} \to  T \Mbar^{S^H,\bR}_{\underline{\alpha}}  ,
        \end{equation}
        whose kernel corresponds to translation. \end{prop}
      \begin{proof}
        The properties of the two moduli spaces are entirely analogous, so we focus on $\fcMR_{\underline{\alpha}}$.       The action of $\bR^{S^H}$ is free because we have assumed that, whenever a component of the domain is equipped with a non-trivial map to one of the factors of $(\bP^1)^{S^H}$, it must carry at least one marked point.  Properness of the action can be seen from analysing the extension to the full moduli space of maps (fixed points are given by configurations in which the components carrying degree $1$ maps are distinct). Each orbit of this action lies in a fibre of the projection map $ \Mbar^{S^H,\bR}_{\underline{\alpha}} \to \Mbar_{\underline{\alpha}}$. In the interior of $\Mbar_{\underline{\alpha}}$, a local identification of the complement of the marked points $\pm$ with $\bR \times S^1$ fixes a local diffeomorphism between fibres and a product of an open subset of $\bR^{S^H}$ (corresponding to the $\bR$ coordinates of the real marked points), with a circle (corresponding to their common angular parameter).  

        The boundary of $\Mbar^{S^H,\bR}_{\underline{\alpha}} $ is mapped to the boundary divisor of $\Mbar_{\underline{\alpha}}$. The boundary strata are indexed by a decomposition $\underline{\alpha} = \underline{\alpha}_1 \times \cdots \times \underline{\alpha}_k$, which gives rise to a commutative diagram
        \begin{equation}
          \begin{tikzcd}
            \Mbar^{S_1^H, \bR}_{\underline{\alpha}_1} \times \cdots \times  \Mbar^{S_k^H, \bR}_{\underline{\alpha}_k} \ar[d] \ar[r] & \Mbar^{S^H,\bR}_{\underline{\alpha}} \ar[d] \\ 
          \Mbar_{\underline{\alpha}_1} \times \cdots \times \Mbar_{\underline{\alpha}_k} \ar[r] &   \Mbar_{\underline{\alpha}}.
          \end{tikzcd}
        \end{equation}
        This diagram is a fibre product, and the image of the top map is a codimension-$k$ boundary stratum. This locally identifies the fibre over codimension-$k$ boundary strata as a product of $(S^1)^k$ with $\bR^{S^H}$.  The above analysis implies that the map $\fcMR_{\underline{\alpha}} \to \Mbar_{\underline{\alpha}} $  is surjective, and that its restriction to the interior of each codimension-$k$ stratum is an $(S^1)^k$ bundle.

        Finally, the contractibility of the space of sections follows from the fact that, over a fixed point of $G_{\underline{\alpha}}$ in $\fcMR_{\underline{\alpha}}$, the action on the fibres is trivial.

        The descriptions of the kernels of tangent spaces are straightforward exercises, which follow from the statement that the projection $\Mbar^{S^H,\bR}_{\underline{\alpha}} \to \Mbar^{\bR}_{\underline{\alpha}} $ does not collapse any component, while the projection $\Mbar^{S^H,\bR}_{\underline{\alpha}} \to  \Rbar_\alpha$ does not collapse any component to a node.
      \end{proof}
     
        \begin{cor}
          For each decomposition $\underline{\alpha} = \underline{\alpha}_1 \times \underline{\alpha}_2$, there is a commutative diagram of $G_{\underline{\alpha}_1} \times G_{\underline{\alpha}_2}$ manifolds
        \begin{equation} 
          \begin{tikzcd}[column sep=small]
    \Mbar^{\prime \bR}_{\underline{\alpha}_1} \times  \Mbar^{\prime \bR}_{\underline{\alpha}_2}   \ar[d]     &  \fcMR_{\underline{\alpha}_1} \times  \fcMR_{\underline{\alpha}_2} \ar[d] \ar[l]  &             \Mbar^{S_1^H, \bR}_{\underline{\alpha}_1} \times   \Mbar^{S_2^H, \bR}_{\underline{\alpha}_2} \ar[d] \ar[r] \ar[l] &\Rbar_{\underline{\alpha}_1}\times \Rbar_{\underline{\alpha}_2}   \ar[d] \\
     \Mbar^{\prime \bR}_{\underline{\alpha}}     &   \fcMR_{\underline{\alpha}} \ar[l]  & \Mbar^{S^H,\bR}_{\underline{\alpha}}  \ar[l]  \ar[r] &    \Rbar_\alpha 
          \end{tikzcd} 
        \end{equation}
        such that the vertical maps enumerate the codimension $1$-boundary strata of the moduli spaces in the bottom row. The construction is associative  in the sense that the following diagram commutes:
  \begin{equation}
    \begin{tikzcd}
      \fcMR_{\underline{\alpha}_1} \times \fcMR_{\underline{\alpha}_2} \times \fcMR_{\underline{\alpha}_3} \ar[r] \ar[d] &  \fcMR_{\underline{\alpha}_1 \times \underline{\alpha}_2} \times \fcMR_{\underline{\alpha}_3} \ar[d]  \\
      \fcMR_{\underline{\alpha}_1} \times \fcMR_{\underline{\alpha}_2 \times \underline{\alpha}_3} \ar[r] &  \fcMR_{\underline{\alpha}_1 \times \underline{\alpha}_2 \times \underline{\alpha}_3}.
    \end{tikzcd}
  \end{equation}
   \qed
 \end{cor}

 Given a pair of elements $\underline{a} < \underline{c}$  in $\underline{\cP}(f,H,f)$, and an element $\underline{b} \in   \cP(\underline{a},\underline{c})$, we define $\partial^{\underline{b}} \fcMR_{\underline{\alpha}}$ for each $\underline{\alpha} \in \cD(\underline{a},\underline{c})$ to be the union of the image of all maps $\fcMR_{\underline{\alpha}_1} \times  \fcMR_{\underline{\alpha}_2} $ such that  $\underline{\alpha}_1 \in \cD(\underline{a},\underline{b})$ and $\underline{\alpha}_2 \in \cD(\underline{b},\underline{c})$.
        \begin{lem} \label{lem:moduli_spaces_abstract_holo_stratified}
          The moduli space $\fcMR_{\underline{\alpha}} $ is a $G_{\underline{\alpha}}$ equivariant smooth $\langle \cP(\underline{a},\underline{c}) \rangle$ manifold with corners. 
        \end{lem}
        \begin{proof}
          The only case which in not clear from the definition is the description of the stratum $\partial^{M_0}  \fcMR_{\underline{\alpha}}$, whenever $\underline{a} = M_-$ and $\underline{c} = M_+$. We define this stratum to be the locus where the two maps to $\bP^1$, labelled by $\chi_-$ and $\chi_+$, agree. 
        \end{proof}

        We now consider an arrow $\underline{\alpha} \to {\underline{\beta}}$ in $\cD(\underline{a},\underline{b})$.
\begin{lem} \label{lem:submersion-floer-cylinder-forgetful-map}
  The following diagram of forgetful maps is commutative, and all arrows except the leftmost vertical map are $G_{\underline{\alpha}}$-equivariant submersions:
\begin{equation} \label{eq:diagram_forgetful_maps_real_moduli_with_marked_points_forget}
  \begin{tikzcd}
 \Mbar^{\prime \bR}_{\underline{\beta}} \ar[d] & \ar[l]    \fcMR_{{\underline{\beta}}}    & \ar[r] \ar[l] \Mbar^{S_{\underline{\beta}}^H,\bR}_{{\underline{\beta}}}  &  \Rbar_{\underline{\beta}} \\
 \Mbar^{\prime \bR}_{\underline{\alpha}}  & \ar[l]   \fcMR_{\underline{\alpha}}  \ar[u] & \ar[r]\ar[u] \ar[l]  \Mbar^{S_{\underline{\alpha}}^H,\bR}_{\underline{\alpha}} & \Rbar_{\underline{\alpha}} \ar[u]
\end{tikzcd}
\end{equation}
The induced map of kernels  $T^{{\underline{\beta}}} \Mbar^{S_{\underline{\alpha}}^H,\bR}_{\underline{\alpha}} \to T^{\underline{\beta}}  \fcMR_{\underline{\alpha}} \oplus T^{\underline{\beta}} \Rbar_{\underline{\alpha}}$ is an isomorphism. Writing $ T^{\underline{\beta}, \Rbar}  \Mbar^{S^H,\bR}_{\underline{\alpha}}$ for the kernel of the induced map 
\begin{equation}
  T^{\Rbar}  \Mbar^{S^H,\bR}_{\underline{\alpha}} \to   T^{\Rbar}  \Mbar^{S^H,\bR}_{\underline{\beta}},
\end{equation}
there are canonical identifications
\begin{align}
  T^{\underline{\beta}}  \fcMR_{\underline{\alpha}}  \cong T^{\underline{\beta}, \Rbar}  \Mbar^{S^H,\bR}_{\underline{\alpha}} & \cong  \bigoplus_{z \in r_{S_{\underline{\alpha}}} \setminus r_{S_{\underline{\beta}}} } T_z \Sigma \\ 
T^{\underline{\beta}} \Mbar^{\bR}_{\pm}(S^H_{\underline{\alpha}}, \{\chi_a^b\}) & \cong \bR^{S_{\underline{\alpha}}^H  \setminus S_{\underline{\beta}}^H }  ,
\end{align}
where $\Sigma$ is the domain of the given element of $ \fcMR_{\underline{\alpha}} $.
\end{lem}
\begin{proof}
The inclusion $ S_1 \subset S_0 $ induces a map
  \begin{equation}
     \Mbar_{r_{S_0} \amalg \pm} \to  \Mbar_{r_{S_1} \amalg \pm}  
   \end{equation}
   by forgetting marked points. The condition that $ S'_0 \subset S'_1 $ implies that points whose underlying curve is pulled back from $ \Mbar_{r_{S'_0} \amalg \pm}$ map to points whose underlying curve is pulled back from  $ \Mbar_{r_{S'_1} \amalg \pm}$. Moreover, given a curve in the domain which is in the image of $\fcMR_{\underline{\alpha}_0}$ the condition that the map from $S'_0$ to the set of levels factors through $\{-\} \amalg P_0 \amalg \{+\}$ implies in particular that the map from $S_0$ to the set of levels is surjective, hence that the forgetful map on this locus preserves levels. This immediately implies that the set of points labelled by elements of $S_1$ lie in the same level, and that the map from $S'_1$ to the set of levels is surjective.  Finally, the condition that the inclusion $P_0 \subset P_1$ is compatible with indexing gives the desired factorisation of the map from the set of successive elements of $\{-\} \amalg P_1 \amalg \{+\}$ to the set of levels.  This gives a map $\Mbar_{\underline{\alpha}} \to \Mbar_{\underline{\beta}}$, which fits in a commutative diagram
   \begin{equation}
     \begin{tikzcd}
       \Mbar_{\underline{\alpha}} \ar[r] \ar[d] & \Mbar_{\underline{\beta}} \ar[d] \\
       \Mbar'_{\underline{\alpha}} & \ar[l]  \Mbar'_{\underline{\beta}}.
     \end{tikzcd}
   \end{equation}
   In order to lift this diagram to the moduli spaces we study, observe that the forgetful maps $\Mbar^{S^H_{\underline{\alpha}},\bR}_{\underline{\alpha}} \to  \Mbar^{S^H_{\underline{\beta}},\bR}_{{\underline{\beta}}} $ and $ \Rbar_\alpha \to \Rbar_\alpha $ are submersions because each curve represented by a point in the domain is biholomorphic to the curve represented by its image. The tangent space of the fibres of the map $  \Rbar_\alpha \to \Rbar_\alpha $ is canonically isomorphic to $ \bR^{S_{\underline{\alpha}}^H \setminus S_{\underline{\beta}}^H}$, while the relative tangent space of the third column of Diagram \eqref{eq:diagram_forgetful_maps_real_moduli_with_marked_points_forget} is the direct sum of $\bigoplus_{z \in r_{\underline{\alpha}} \setminus r_{\underline{\beta}} } T_z \Sigma $ with $ \bR^{S_{\underline{\alpha}}^H \setminus S_{\underline{\beta}}^H}$. The first factor is transverse to the orbits of the action by $\bR^{S_{\underline{\alpha}}^H}$, and the second is contained in it.  The compatibility of this submersion with the foliation thus induces a submersion $ \fcMR_{\underline{\alpha}}  \to  \fcMR_{{\underline{\beta}}} $, with the stated relative tangent space. 
\end{proof}
Note that the map $ \Mbar^{\prime \bR}_{\underline{\beta}} \to \Mbar^{\prime \bR}_{\underline{\alpha}}$ may not be a submersion because we have not imposed enough conditions on the fibre of the universal curve of elements of $\Mbar^{\prime \bR}_{\underline{\beta}}$. Comparing with Definition \ref{def:equivariant-submersion-category}, we have:
\begin{cor}
  The assignment $\underline{\alpha} \mapsto (G_{\underline{\alpha}},\Mbar^{\prime \bR}_{\underline{\alpha}} \to \fcMR_{\underline{\alpha}})$ defines a functor from $\cD(\underline{a},\underline{b})$ to the category of equivariant submersions of $\langle \underline{\cP}(\underline{a},\underline{b}) \rangle$-manifolds. \qed
\end{cor}

Next, we state the compatibility between the maps associated to arrows and to products:
\begin{lem}
  Given maps $\underline{\alpha}_1 \to {\underline{\beta}}_1$ and $\underline{\alpha}_2 \to {\underline{\beta}}_2$, the following diagram commutes:
  \begin{equation}
    \begin{tikzcd}[column sep=small]
      \fcMR_{\underline{\alpha}_1} \times \fcMR_{\underline{\alpha}_2} \ar[r] \ar[d] &  \fcMR_{{\underline{\beta}}_1} \times \fcMR_{{\underline{\beta}}_2} \ar[r] \ar[d] &  \Mbar^{\prime \bR}_{{\underline{\beta}}_1} \times \Mbar^{\prime \bR}_{{\underline{\beta}}_2} \ar[r] \ar[d] &  \Mbar^{\prime \bR}_{\underline{\alpha}_1} \times \Mbar^{\prime \bR}_{\underline{\alpha}_2} \ar[d]  \\
\fcMR_{\underline{\alpha}_1 \times \underline{\alpha}_2}  \ar[r]     &  \fcMR_{{\underline{\beta}}_1 \times {\underline{\beta}}_2}  \ar[r]  & \Mbar^{\prime \bR}_{{\underline{\beta}}_1 \times {\underline{\beta}}_2}  \ar[r] & \Mbar^{\prime \bR}_{\underline{\alpha}_1 \times \underline{\alpha}_2}
\end{tikzcd}
  \end{equation}
Moreover, the corresponding associativity diagrams commute.
\qed
\end{lem}

\begin{cor}
The assignment of the submersion $ \fcMR_{\underline{\alpha}} \to \Mbar^{\prime \bR}_{\underline{\alpha}}$ to each object $\underline{\alpha} \in \cD(\underline{a},\underline{b})$ extends to a strictly $\Pi$-equivariant lax $2$-functor from $\cD(f,H,f)$  to the category of stratified equivariant submersions.  \qed 
\end{cor}

\section{Kuranishi flow categories from Floer theory}
\label{sec:kuran-flow-categ-1}
\label{sec:hamilt-floer-theory}

\subsection{Categories of thickening data}
\label{sec:categ-thick-data}

Given a pair $a < b$ of objects of $\cP(f,H,f)$ such that $\cD(\underline{a}, \underline{b})$ is not a point we shall define a Floer thickening datum $\alpha$ for a stratum of the moduli space $\fcMR(a,b)$ to consist of:
\begin{enumerate}[(i)]
\item an object $\underline{\alpha} $ of $\cD(\underline{a}, \underline{b})$. We shall write $\underline{\alpha} = (S'_\alpha \subset S_\alpha, P_\alpha \subset P'_\alpha) $ for the constituent sets.
  \item (Choice of stratum) a subset $ Q_\alpha$ of  $\cP(a,b)$, such that $Q_\alpha \cap \cP(H)$ is contained in $P_\alpha$. We define $\underline{Q}_\alpha $ to be the projection to $\underline{\cP}(\underline{a}, \underline{b})$ of the intersection of $Q_\alpha$ with $\cP(H) \cap \cP(f)_0$,  
  \item a commutative diagram
\begin{equation} \label{eq:inclusion_open_subset_stratum-abstract-moduli-discs}
  \begin{tikzcd}
\fcMR_{\alpha}  \ar[r] \ar[d] & \partial^{\underline{Q}_\alpha} \fcMR_{\underline{\alpha}} \ar[d] \\
  \Mbar^{\prime \bR}_{\alpha} \ar[r] & \partial^{\underline{Q}_\alpha} \Mbar^{\prime \bR}_{\underline{\alpha}} 
  \end{tikzcd}
\end{equation}
in which the horizontal arrows are  $G_\alpha$ invariant open inclusions. We write $\fcCR_{\alpha}\to  \fcMR_\alpha$ for the pullback of the universal curve over $  \fcMR_{\underline{\alpha}}$ to $\fcMR_\alpha $. ,
  
\item a complex linear finite dimensional subspace
\begin{equation}
  V_\alpha \subset \bC^\infty \otimes_{\bC} \bigoplus_{i \in S} \bC [ G_{r_i}],
\end{equation}
which is invariant under the group $G_\alpha \equiv \bigoplus G_{r_i}$. 
\item a $G_\alpha$ equivariant complex-linear map
  \begin{equation}
   \lambda_\alpha \co V_\alpha \to C^{\infty}_c( \fcCR_{\alpha}  \times M, T^{0,1}_{\fib}\fcCR_{\alpha}   \otimes TM),
  \end{equation}
  whose image consists of sections which are supported away from all nodes, and away from the positive and negative ends $\epsilon_\alpha$. Here,
\begin{equation}
  T^{0,1}_{\fib}\fcCR_{\alpha}   \otimes TM
\end{equation}
is the vector bundle over $ \fcCR_{\alpha} \times M$ consisting of
complex anti-linear maps from the relative tangent bundle of the
projection $\fcCR_{\alpha}\to  \fcMR_\alpha$ to $TM$, and we write 
\begin{equation}
   C^{\infty}_c( \fcCR_{\alpha} \times M, T^{0,1}_{\fib}\fcCR_{\alpha}   \otimes TM)
\end{equation}
for the space of sections which are supported away from all nodes and all ends.
\item a collection $D_\alpha$ of compact codimension $2$ submanifolds with boundary $D_i \subset M$, indexed by $S_\alpha$, and equipped with a path of sub-bundles of $TM$ along $D_i$ from $TD_i$ to an almost complex sub-bundle.
\end{enumerate}
Whenever $\cD(\underline{a}, \underline{b})$ is a singleton, the only non-trivial choice above is that of the subset $Q_\alpha$; we set  $\fcMR_{\alpha} = \{ *\}$ and $V_\alpha = \{0\}$ in this case.

At this stage, we are ready to introduce a category of thickening data:

\begin{defin}\label{def:thickening-data} 
  For each pair $(a,b)$ of elements of $\cP(f,H,f)$, define a category $A(a,b)$ of  \emph{Floer thickening data} whose objects are given by the data
  \begin{equation}
    \alpha = (\underline{\alpha}, \fcMR_{\alpha}, \Mbar^{\prime \bR}_{\alpha}, Q_\alpha,  V_\alpha, \lambda_\alpha,D_\alpha ),
  \end{equation}
  as above.
  A morphism $\alpha \to \beta$ consists of
\begin{enumerate}
\item a map $\underline{\alpha} \to \underline{\beta}$ in  $\cD(a,b)$ (and in particular a sujection  $G_\alpha \to G_\beta$),
\item an inclusion $ Q_\beta \subset Q_\alpha$, such that the induced maps of moduli spaces restrict to  maps
  \begin{align}
    \fcMR_{\alpha} & \to \fcMR_{\beta} \\
    \Mbar^{\prime \bR}_\beta & \to  \Mbar^{\prime \bR}_\alpha,  \end{align}
\item a complex linear isometric embedding $V_\alpha \to V_\beta$ of $G_\alpha$-representations, so that the diagram
  \begin{equation}
    \begin{tikzcd}
    V_\alpha \ar[d] \ar[r] & V_\beta \ar[d] \\
    C^{\infty}_c( \fcCR_{\alpha} \times M, T^{0,1}_{\fib}\fcCR_{\alpha}   \otimes TM) & \ar[l] C^{\infty}_c( \fcCR_{\beta} \times M, T^{0,1}_{\fib}\fcCR_{\beta}   \otimes TM)  
    \end{tikzcd}
      \end{equation}
      commutes.
    \end{enumerate}
    We require that the submanifold $D_i$ associated to each element of $S_\beta$ agree with the submanifold associated to its image in $S_\alpha$, and similarly for the path of subbundles of $TM$ over $D_i$.
\end{defin}
The definition of composition in this category is straightforward once
we observe the following: if we are given a pair of arrows $f \co
\alpha \to \beta$ and $g \co \beta \to \gamma$, and subsets
$\fcMR_{{\alpha}} \subset \fcMR_{\underline{\alpha}}$,
$\fcMR_{{\beta}} \subset \fcMR_{\underline{\beta}}$, and
$\fcMR_{{\gamma}} \subset \fcMR_{\underline{\gamma}}$ such that
$\fcMR_{{\alpha}}$ maps to $\fcMR_{{\beta}}$ under $f$ and
$\fcMR_{{\beta}}$  maps to $\fcMR_{{\gamma}} $ under $g$, we can
conclude that $\fcMR_{{\alpha}}$ maps to $\fcMR_{{\gamma}} $ under the
composite $g \circ f$.

\begin{rem}
The requirement that $V_\alpha$ be contained in the $G_{\alpha}$-universe $\bC^\infty \otimes_{\bC} \bigoplus_{i \in S} \bC [ G_{r_i}]$ ensures that $A(a,b)$ is a small category, and will imply that the compositions we shall presently define are strictly associative. The reason for passing to open subsets is related to the need to later construct stable complex structure in Section \ref{sec:oriented-locus-space}
\end{rem}
In analogy with Lemma \ref{lem:Pi-action-D}, we have:
\begin{lem} 
  Each element $\pi \in \Pi$ induces a canonical isomorphism of categories
  \begin{equation}
      \pi \cdot \_ \co   A(a,b) \to A(\pi \cdot a, \pi \cdot b),
      \end{equation}
      such that the following diagram
      \begin{equation}
        \begin{tikzcd}
          A(a,b) \ar[r," \pi \cdot \_"] \ar[dr," (\pi' \cdot \pi) \cdot \_"] &  A(\pi \cdot a, \pi \cdot b) \ar[d," \pi' \cdot \_"] \\
          & A(\pi \cdot a, \pi \cdot b) 
        \end{tikzcd}
      \end{equation}
      commutes.  This assignment is unital in the sense that $e
      \in \Pi$ is taken to the identity map and it is strictly
      associative.
\qed
\end{lem}

We now consider the multiplicativity of this construction: given a pair $\alpha_1$ and $\alpha_2$ of Floer thickening data for strata of the moduli spaces $\fcMR(a,b)$ and $\fcMR(b,c)$, we define the product Floer thickening datum $\alpha_1 \times \alpha_2$ for a stratum of the moduli space $\fcMR(a,c)$ to consist of:
\begin{enumerate}[(i)]
\item the composition $\underline{\alpha}_1 \times \underline{\alpha}_2 \in \cD(\underline{a}, \underline{c})$,
\item the  subset of $\cP(a,c)$ defined as $ Q_{\alpha_1 \times \alpha_2} \equiv  Q_{\alpha_1} \amalg \{ b \} \amalg Q_{\alpha_2}$, 
\item the commutative diagram
\begin{equation}
  \begin{tikzcd}
 \fcMR_{\alpha_1 \times \alpha_2}  \equiv \fcMR_{\alpha_1} \times \fcMR_{\alpha_2}  \ar[r] \ar[d] & \partial^{\underline{Q}_{\alpha_1 \times \alpha_2}} \fcMR_{\underline{\alpha}_1 \times \underline{\alpha}_2} \ar[d] \\
 \Mbar^{\prime \bR}_{\alpha_1 \times \alpha_2}  \equiv  \Mbar^{\prime \bR}_{\alpha_1} \times \Mbar^{\prime \bR}_{\alpha_2} \ar[r] & \partial^{\underline{Q}_{\alpha_1 \times \alpha_2}} \Mbar^{\prime \bR}_{\underline{\alpha}_1 \times \underline{\alpha}_2},
  \end{tikzcd}
\end{equation}
\item the finite dimensional subspace
  \begin{equation}
  V_{\alpha_1 \times \alpha_2} \cong V_{\alpha_1} \oplus V_{\alpha_2} \subset  \bC^\infty \otimes_{\bC} \bC [ G_{\alpha_1 \times \alpha_2}],
  \end{equation}
\item the map 
  \begin{equation}
   \lambda_{\alpha_1 \times \alpha_2} \co V_{\alpha_1 \times \alpha_2} \to C^{\infty}_c( \fcCR_{\alpha_1 \times \alpha_2} \times M, T^{0,1}_{\fib}\fcCR_{\alpha_1 \times \alpha_2}  \otimes TM) , 
 \end{equation}
given by the direct sums of $\lambda_1$ and $\lambda_2$, and
\item the collection $D_{\alpha_1 \times \alpha_2}$ which is the union of $D_{\alpha_1}$ with $D_{\alpha_2}$.
\end{enumerate}

Given an ordered triple $(a,b,c)$ of elements of $\cP(f,H,f)$, the product of Floer thickening data defines a functor
\begin{equation} \label{eq:product_of_charts}
  A(a,b) \times A(b,c) \to A(a,c)
\end{equation}
such that the associativity diagram
\begin{equation}
  \begin{tikzcd}
    A(a,b) \times A(b,c) \times A(c,d) \ar[r] \ar[d] & A(a,b) \times A(b,d)  \ar[d] \\
     A(a,c) \times A(c,d) \ar[r] & A(a,d)
  \end{tikzcd}
\end{equation}
strictly commutes and the functor is unital.  Moreover, the action of
$\Pi$ is strict.  Summarizing, we have the following result.

\begin{lem}
  The functor in Equation \eqref{eq:product_of_charts} equips the collection $A(f,H,f)$ of categories $A(a,b)$ with the structure of a strict $\Pi$-equivariant $2$-category. \qed 
\end{lem}

\subsection{Regular thickenings of moduli spaces}
\label{sec:regul-thick-moduli}

Next, we consider a fibre $\Sigma$ of the universal curve $\fcCR_{\alpha} \to \fcMR_{\alpha} $ for $\alpha \in A(a,b)$, and a pseudo-holomorphic stable building $u$, with domain $\Sigma$. Denoting the direct sum of the spaces of smooth anti-holomorphic $1$-forms on each component of $\Sigma$, with values in the pullback of $TM$ by  $ \Omega^{0,1}(\Sigma, u^*TM)$, we have an element
\begin{equation}
  \dbar_{(a,b)} u \in  \Omega^{0,1}(\Sigma, u^*TM)
\end{equation}
given on each Floer cylinder by Equation~\eqref{eq:Floer-equation}, on components which carry the marked point $\chi$ by Equation~\eqref{eq:continuation-equation}, and on all other components by the homogeneous $\dbar$ operator. For the next definition, recall that a choice of element in $\alpha \in \cD(a,b)$ entails a choice of totally ordered subset $P_\alpha$ of $\cP(a,b)$, and that the Floer edges of trees labelling strata of $\fcMR_\alpha$ are decorated by elements of $P_\alpha$. In addition, for $a \in \cP(f,H,f)$,  we introduce the notation
\begin{equation}
  [\underline{a}] =
  \begin{cases}
    [a] & \textrm{ if } a \in \cP(H) \\
    M & \textrm{ otherwise. } 
  \end{cases}
\end{equation}
\begin{defin} \label{def:thickened_moduli_space}
For each $\alpha \in A(a,b)$, the \emph{thickened moduli space} $\fcMR_{\alpha}([\underline{a}],[\underline{b}])$ is the space whose elements $(\Sigma, u, v)$ consist of
\begin{enumerate}
\item a fibre $\Sigma$ of the universal curve  $\fcCR_{\alpha} \to \fcMR_{\alpha}$,
\item a stable map $u \co \Sigma \to M$ with asymptotic conditions  $[\underline{a}]$ along the negative end, and $[\underline{b}]$ along the positive end, and given along each Floer node by the corresponding element of $P_\alpha$. For each $i \in S$, we require that the set of points labelled by $r_i$ be invariant under the action of the automorphism group of the underlying map without marked points, and that their images under $u$ lie in  $D_i$.
\item an element $v$ of $V_\alpha$ 
such that $u$ satisfies the inhomogeneous Cauchy-Riemann equation
\begin{equation} \label{eq:inhomegeneous-CR}
 \dbar_{(a,b)} u +  \lambda (v)  = 0.
\end{equation}
\end{enumerate}
\end{defin}
To clarify the definition, whenever $[\underline{a}]$ or $[\underline{b}] $  is given by $M$, the asymptotic condition in the definition is an unspecified point in $M$. In the case $a$ lies in $\cP(f)_-$ and $b$ in $\cP(f)_0$, we interpret the above definition to say that $\fcMR_{\alpha}([\underline{a}],[\underline{b}]) = M$, corresponding to constant maps.

By construction, we have a continuous projection map
\begin{equation}
  \fcMR_{\alpha}([\underline{a}],[\underline{b}]) \to \fcMR_{\alpha}.
\end{equation}
Since the target of this projection is stratified by $\underline{\cP}(\underline{a}, \underline{b})$, the inverse images induce a stratification of $\fcMR_{\alpha}([\underline{a}],[\underline{b}])$.

Note that we have a natural map
\begin{equation}
s \co      \fcMR_{\alpha}([\underline{a}],[\underline{b}]) \to V_\alpha,
\end{equation}
which is $G_\alpha$ equivariant, providing us with all the ingredients for a Kuranishi chart. However, the total space is not necessarily a manifold. In order to ensure this we introduce some notation: given an element of $  \fcMR_{\alpha}([\underline{a}],[\underline{b}]) $, we write
\begin{equation}
   \cF(u^* TM, u^* T D_\alpha) 
 \end{equation}
 for the space of sections of $u^* TM$ on each component of the domain, whose values at the marked points labelled by $r_i$ lie in $u^* T D_i$. We require these sections to have appropriate regularity (e.g., of $W^{k,2}$ class for $2 \leq k$), and to decay exponentially along the ends. The linearisation of the Floer equation (and the holomorphic curve equation on sphere bubbles), defines a map
 \begin{equation}
      \cF(u^* TM, u^* T D_\alpha ) \to     \cE^{0,1}( u^* TM)  
    \end{equation}
    where the right hand side is the direct sum of a copy of the fibre of $TM$ at each node with the spaces of complex antilinear $1$-forms, valued in the pullback of $TM$, on each component of the domain  (c.f. Section \ref{sec:based-floer-traj} below). The regularity of the target is assumed to be one lower than that of the domain.
\begin{defin} \label{def:regular_moduli_space}
The regular part $\fcMRreg_{\alpha}([\underline{a}],[\underline{b}]) \subset \fcMR_{\alpha}([\underline{a}],[\underline{b}])$ is the locus of points where (i)  $u$ is transverse to $D_i$ at each point in $r_i$, and (ii) the linearised $\dbar$ operator
  \begin{equation} \label{eq:dbar_operator}
    V_\alpha \oplus   \cF(u^* TM, u^* T D_\alpha ) \to     \cE^{0,1}( u^* TM) 
  \end{equation}
is surjective.
\end{defin}
\begin{rem}
Our notion of regularity amounts to the condition that the linearisation of the $\dbar$ operator \emph{with constraints on the modulus of the corresponding points in $\fcMR_{\alpha}$} is surjective. There is a weaker notion of regularity where one would require only surjectivity of the operator for the parametrised problem over this moduli space. Moreover, the fact that we assumed surjectivity of the operator with exponentially decaying conditions implies that the evaluation map to $M$ is a submersion at each end with vanishing Hamiltonian conditions. In fact, if there are two such ends, then the evaluation map to $M^2$ at both ends is a submersion.
\end{rem}

Recall that our goal is to construct Kuranishi charts for the moduli spaces in Equations \eqref{eq:moduli_space_Floer-cylinders}, \eqref{eq:moduli-space-Morse-to-Floer}--\eqref{eq:moduli_space_Morse_0_+} and \eqref{eq:piecewise_definition_moduli_-_+}. To simplify the notation, we shall write $Q_{\alpha}(f)_0$ and $Q_{\alpha}(f)_\pm$ for the intersections of $Q_{\alpha}$ with $\cP(f)_0$ and $\cP(f)_\pm$. Given elements $(x,y,z)$ of $\cP(f)$ and  $(p,q)$ of $\cP(H)$, we then define
\begin{align}
  \begin{rcases}
    &\fcMRreg_{\alpha}(x_-,y_-)  \\
    &\fcMRreg_{\alpha}(x_0,y_0)  \\
    &\fcMRreg_{\alpha}(x_+,y_+)
  \end{rcases}
    & \equiv \partial^{Q_\alpha} \barT(x,y) \\
  \fcMRreg_{\alpha}(p,q)  & \subset \fcMRreg_{\alpha}([p],[q]) \\
\fcMRreg_{\alpha}(x_-,p)  & \subset \partial^{Q_\alpha (f)_-} \barT([x],M) \times_M  \fcMRreg_{\alpha}(M,[p]) \\
\fcMRreg_{\alpha}(p,z_+) & \subset \fcMRreg_{\alpha}([p],M)  \times_M  \partial^{Q_\alpha (f)_+} \barT(M,[z]) \\
\fcMRreg_{\alpha}(y_0,z_+) & \subset  \\
\notag & \hspace{-.2in}  \partial^{Q_\alpha (f)_0} \barT([y],M) \times_M \fcMRreg_{\alpha}(M,M) \times_M \partial^{Q_\alpha (f)_+} \barT(M,[z])  \\
\fcMRreg(x_-,y_0) & \subset \partial^{Q_\alpha}  \barT([x],M,[y])  ,
\end{align}
to consist of those curves which lift to the cover $\tilde{\cL}M$ of the loop space as paths from the first to the second element of $\cP(f,H,f)$. It remains to define the moduli space for a pair $x_- \in \cP(f)_-$ and $z_+ \in \cP(f)_+$, which we break up into several cases (in all cases, the given subset consists of all elements of the right hand side in the given homotopy class):
\begin{enumerate}[(i)]
\item If  $Q_{\alpha}(f)_0 \neq \emptyset$, we have
  \begin{multline}
 \fcMRreg_{\alpha}(x_-,z_+)  \subset \\ \partial^{Q_\alpha(f)_- \cup \Q_{\alpha}(f)_0} \barT([x],M,M) \times_M  \partial^{M_0} \fcMRreg_{\alpha}(M,M) \times_M \partial^{Q_\alpha(f)_+}  \barT(M,[z]) ,   
\end{multline}
where $ \partial^{M_0} \fcMRreg_{\alpha}(M,M)$ is given in Lemma \ref{lem:moduli_spaces_abstract_holo_stratified}.
\item  If  $Q_{\alpha} \cap \cP(H) \neq \emptyset$, we have
  \begin{multline}
 \fcMRreg_{\alpha}(x_-,z_+)  \subset \\ \partial^{Q_\alpha (f)_-}  \barT([x],M) \times_M \fcMRreg_{\alpha}(M,M) \times_M \partial^{Q_\alpha (f)_+}  \barT(M,[z]),
\end{multline}
where we note that $\underline{Q}_\alpha =  Q_{\alpha} \cap \cP(H)$ in this case, so that this intersection enters in the definition of the moduli space $\fcMRreg_{\alpha}(M,M) $ via the top line in Diagram \eqref{eq:inclusion_open_subset_stratum-abstract-moduli-discs}.
\item If $Q_{\alpha}(f)_0  = \emptyset =  Q_{\alpha} \cap \cP(H)$, we have
  \begin{multline}
    \fcMRreg_{\alpha}(x_-,z_+)  \subset \\ \partial^{Q_\alpha (f)_-}  \barT([x],M) \times_M \fcMRreg_{\alpha}(M,M) \times_M \partial^{Q_\alpha (f)_+}  \barT(M,[z]) \cup  \\ 
 \partial^{Q_\alpha(f)_- } \barT([x],M,M) \times_M  \partial^{M_0} \fcMRreg_{\alpha}(M,M) \times_M \partial^{Q_\alpha(f)_+}  \barT(M,[z]),
\end{multline}
where the union is taken along the common stratum
\begin{equation}
 \partial^{Q_\alpha (f)_-}  \barT([x],M) \times_M  \partial^{M_0}  \fcMRreg_{\alpha}(M,M) \times_M \partial^{Q_\alpha (f)_+} \barT(M,[z]).
\end{equation}
\end{enumerate}
In the above definition, we have implicitly used the fact that the three cases are mutually exclusive: the key point is that elements of $\cP(H)$ and $\cP(f)_0$ are not comparable, so that the first two cases have no overlap. Keeping this in mind, we can in fact give a uniform definition as
\begin{multline} \label{eq:Kuranishi-chart-minus-to-+}
  \fcMRreg_{\alpha}(x_-,z_+)  \subset \\
  \partial^{Q_\alpha (f)_-}  \barT([x],M) \times_M \fcMRreg_{\alpha}(M,M) \times_M \partial^{Q_\alpha (f)_+}  \barT(M,[z])  \cup \\ \partial^{Q_\alpha(f)_- \cup \Q_{\alpha}(f)_0} \barT([x],M,M) \times_M  \partial^{M_0} \fcMRreg_{\alpha}(M,M) \times_M \partial^{Q_\alpha(f)_+}  \barT(M,[z]).
\end{multline}

We now state the main consequence of regularity and gluing theory for holomorphic curves, which follows from the standard methods as can be found in \cite{McDuffSalamon2012}, and:
\begin{prop}[c.f. \cite{Pardon2016}]
 \label{prop:Ham_Kuranishi_chart}
  The  quadruple $(\fcMRreg_{\alpha}(a,b), V, s, G_\alpha)$ is a $ \partial^{Q_\alpha} \cP(a,b)$-stratified Kuranishi chart.
\end{prop}
\begin{proof}
  We discuss the case in which $a$ and $b$ are both orbits; the general case follows by the same analysis, keeping in mind that the evaluation map at the ends labelled by $M$ is submersive, so that the fibre product is necessarily transverse.

  It is convenient to use the forgetful map from $\fcMRreg_{\alpha}(a,b) $ to $\fcMR_{\alpha}  $, which by definition respects stratification. The inverse image of the interior of a stratum in $\fcMR_{\alpha} $ is a smooth manifold, as it can be locally described as the zero-locus of a Fredholm section of a Banach bundle over a base which is the product of a smooth chart for the given stratum of $\fcMR_{\alpha}  $ (given as a prestable Riemann surface with marked points, equipped with varying conformal structure) with the vector space $V_\alpha$, and with the space of maps of appropriate regularity.

  The gluing theorem for moduli spaces of pseudo-holomorphic curves, which is implemented in detail in \cite[Appendix C]{Pardon2016} shows that the total space $ \fcMRreg_{\alpha}(a,b)$ is a topological manifold with the stated stratification.
\end{proof}

The regular locus is compatible with the $G_{\alpha}$ action, so that the data
\begin{equation}
\bM_{\alpha}(a,b) \equiv \left( \fcMRreg_{\alpha}(a,b),   V_{\alpha},  s_\alpha, G_{\alpha} \right)
\end{equation}
defines an object of $\Kur^{\cP(a,b)}$; the stratification is given by orbits and critical points that appear as limits of each level.

\begin{rem}
We shall often abuse notation and write $u$ for an element of $  \fcMR_{\alpha}(a,b)$, rather than specifying that it (possibly) consists (of a pair) of (half)-gradient flow lines as well.
\end{rem}

\subsection{Kuranishi presentation of moduli spaces}
\label{sec:kuran-pres-moduli}

We begin by considering the functoriality of the construction in the previous section:
\begin{lem} The assignment $\alpha \mapsto \bM_{\alpha}(a,b)  $ extends to a functor
  \begin{equation}
    \bM(a,b) \co   A(a,b) \to \Kur .
  \end{equation}
\end{lem}
\begin{proof}
The key point is to assign a natural map of Kuranishi charts
\begin{equation}
 \bM_{\alpha_0}(a,b) \to  \bM_{\alpha_1}(a,b)
\end{equation}
to each map from $\alpha_0$ to $\alpha_1$. The surjection $G_{\alpha_0} \to G_{\alpha_1}$ and the inclusion $V_{\alpha_0} \to V_{\alpha_1}$ are given by the definition of a morphism, and we also have a forgetful map $\fcMR_{\alpha_0} \to \fcMR_{\alpha_1}$. It remains to lift this to a map
\begin{equation}
    \fcMR_{\alpha_0}(a,b) \to  \fcMR_{\alpha_1}(a,b)
  \end{equation}
  and then check compatibility with regular loci.

  For the first part, fix a point in $ \fcMR_{\alpha_1} $. This determines the domain for curves lying in its inverse image in  $ \fcMR_{\alpha_1}(a,b)$, and similarly for its inverse image in $ \fcMR_{\alpha_0}(a,b)$ via the projection to $\fcMR_{\alpha_0}$. Since the forgetful map of abstract moduli spaces induces an isomorphism on the pullback of the universal curve, it suffices to show that we impose less conditions on a map to obtain an element of  $ \fcMR_{\alpha_1}(a,b)$ than we need to obtain an element of $ \fcMR_{\alpha_0}(a,b) $: (i) the asymptotic conditions at each Floer node are the same by the inclusion $P_0 \subset P_1$, (ii) the inclusion $S_1 \subset S_0$ and the compatibility on the choice of  divisors implies that fewer marked points conditions are imposed in the target, and (iii) the  inclusion $V_{\alpha_0} \to V_{\alpha_1}$ implies that the allowed space of inhomogeneous terms is  larger. We conclude that we have a commutative diagram
  \begin{equation}
    \begin{tikzcd}
      \fcMR_{\alpha_0}(a,b ) \ar[r] \ar[d]&   \fcMR_{\alpha_1}(a,b) \ar[d] \\
      \fcMR_{\alpha_0} \ar[r] &  \fcMR_{\alpha_1}.
    \end{tikzcd}
  \end{equation}

  To check that the map $\fcMR_{\alpha_0}(a,b) \to \fcMR_{\alpha_1}(a,b) $ preserves regular loci, we observe that the condition for regularity in the domain are stronger than in the range, because of the inclusion $V_{\alpha_0} \to V_{\alpha_1}$, and the fact that we impose divisorial conditions at marked points by the set $S_0$ which contains $S_1$.

 Checking that the conditions for a map of Kuranishi charts hold is now straightforward: the kernel of the map $G_{\alpha_0} \to  G_{\alpha_1}$ acts freely on $ \fcMR_{\alpha_0}(a,b)$ because the forgetful map of domains does not collapse any component, and the implicit function theorem together with the gluing theorem imply that the space $\fcMR_{\alpha_1}(a,b)$ is locally the product of $V_{\alpha_1}/V_{\alpha_0} $ with the image of $\fcMR_{\alpha_0}(a,b) $. 
\end{proof}

For each $\alpha \in A(a,b)$, we have a natural equivariant projection
\begin{equation}
     \fcMR_{\alpha}(a,b)  \times_{ V_{\alpha}} \{0\} \to \fcMR(a,b),
\end{equation}
which is compatible with morphisms. Thus, the functor $\bM(a,b)$ is naturally a functor over the moduli space $\fcMR(a,b)$. 

\begin{prop} \label{prop:Ham-Kuranishi-presentation}
  The functor $\bM(a,b) $ is a Kuranishi presentation of $\fcMR(a,b)$.
\end{prop}

\begin{proof}
Given a point $[u] \in \fcMR(a,b)$, let $A(a,b)[u]$ denote the
subcategory of $A(a,b)$ consisting of elements $\alpha$ for which
$[u]$ lifts to $ \fcMRreg_{\alpha}(a,b)$. We shall prove that the
nerve of $A(a,b)[u]$ is contractible.  To see this, it will suffice to
show that for any functor $F \co \cC \to A(a,b)[u]$ from a category
with finitely many objects, the induced map of nerves is
null-homotopic.  By taking such a category be to the category of simplices
associated to subdivisions of the spheres, we conclude that the
homotopy groups of the nerve of $A(a,b)[u]$ vanish, and by the
Whitehead theorem, we have that $A(a,b)[u]$ is contractible.
(c.f. the proof of~\cite[1.6.7]{Waldhausen1985}).

To show that the induced map of nerves $N_\bullet F$ is
null-homotopic, we recall that a natural transformation $F \to G$
induces a homotopy of maps of nerves $N_\bullet F \htp N_\bullet G$.
Thus, our null-homotopy will be constructed as a zig-zag of natural
transformations.  More precisely, we pick a convenient chart $\omega
\in A(a,b)$ so that $u$ lifts to the regular locus of
$\fcMR_{\omega}(a,b)$, and construct a zig-zag of natural
transformations between $F$ and the constant functor at $\omega$. The
most important condition we impose is that $D_\omega$ be disjoint from
all divisors $D_{F(c)}$ appearing in the image of $F$. For simplicity,
we also assume that the element of $\cD(a,b)$ corresponding to
$\omega$ is of the form $(S'_{\omega} = S_{\omega}, Q_\omega =
P'_\omega = P_\omega = P_u)$, where $P_u$ is the set of asymptotic
conditions appearing as intermediate orbits in $u$. We also require
that $V_\omega$ be sufficiently large that it surjects onto the
linearised $\dbar$ operator at $u$, with constraints along the
divisors $D_\omega \amalg D_{F(c)}$ for each $c \in \cC$.

The first step is to observe that $P_u$ is a subset of $P_{F(c)}$ for
each $c \in \cC$. The triple $(S'_c \subset S_c, P_u \subset
P_{F(c)})$ defines an object of $\cD(a,b)$, and there is a
corresponding inclusion of moduli spaces which is locally a
homeomorphism near the image of $u$ in $\fcMR_{F(c)}$. We thus obtain
a functor from $\cC$ to $A(a,b)$, equipped with a natural
transformation to $F$, with the property that the set of asymptotic
orbits appearing for each object of $\cC$ is $P_u\subset P_{F(c)}$. A
similar argument reduces the problem to the situation in which $
Q_{F(c)} = P'_{F(c)} = P_{F(c)} = P_u$ for each $c \in \cC$. 

Next, for each object $\alpha$ of $A(a,b)$ which in the image of $F$,
consider a map  $\alpha \to \alpha'$ in $A(a,b)$ where $\alpha'$ given
by the same data as $\alpha$ except that the representation
$V_{\alpha'}$ properly includes $V_\alpha$. This construction yields
an inclusion 
  \begin{equation}
    \fcMR_{\alpha}(a,b)\to \fcMR_{\alpha'}(a,b).  
  \end{equation}
We may choose $V_{\alpha'}$ large enough so that any lift of $u$
remains regular if we impose divisorial conditions along the
hypersurfaces $D_\omega$ associated to the object $\omega$ fixed at
the beginning of the proof. This in particular means that these points
lie in the image of the Kuranishi chart associated to the object of
$A(a,b)$ given by the same data as $\alpha'$, except that we take the
disjoint union of $D_\alpha$ with $D_\omega$, and add the sequence
$S_\omega$ to $S_\alpha$. 

Since $A(a,b) $ is filtered by the number of elements of the set
$S_\alpha$, we can use an inductive argument over the image of $\cC$
under $F$ to choose $\alpha'$ together with a prescribed commutative
diagram 
  \begin{equation}
    \begin{tikzcd}
      V_\alpha \ar[r] \ar[d] & V_{\beta} \ar[d] \\
      V_{\alpha'} \ar[r] & V_{\beta'}
    \end{tikzcd}
  \end{equation}
whenever $\alpha \to \beta$ is the image of an arrow in $\cC$. This
yields a zig-zag of natural transformations from $F$ to a functor $F'$
which now has the property that $S_{F'(c)}$ contains $S_\omega$, and
$D_{F'(c)} $ contains $D_\omega$ (for all $c $ in the domain). By
further considering the natural transformations  associated to the
inclusion $S'_{F(c)} \subset S'_{F(c)} \amalg S_{\omega} \supset
S_{\omega}$, we have a further zig-zag of natural transformations to a
functor $F''$ with the property that $S_{F''(c)} = S_\omega$.  

For the next step, we pick an embedding of $V_{F(c)}$ and $V_\omega$
in a representation $V_{\omega(c)}$ which contains their direct sum,
and keep the remaining data unchanged. This yields a zig-zag to a
chart with marked points given by $S_{F(c)} \amalg S_{\omega}$, and
obstruction space given by $V_\omega$. Forgetting the points marked by
$S_{F(c)}$ yields the final natural transformation.  This completes
the proof that the map induced by $F$ on nerves is null-homotopic.
\end{proof}

\subsection{Kuranishi presentations of product moduli spaces}
\label{sec:kuran-pres-prod}

Given an ordered triple $(a,b,c)$ of elements of $\cP(f,H,f)$, recall that the product of Floer thickening data associates to a thickening $\alpha$ of $ \partial^{Q_\alpha} \fcMR(a,b)$ and a thickening $\beta$ of $\partial^{Q_\beta} \fcMR(b,c)$ a thickening $\alpha \times \beta$ of $\partial^{Q_{\alpha \times \beta}} \fcMR(a,c)$, where $Q_{\alpha \times \beta} \equiv  Q_{\alpha} \amalg \{ b \} \amalg Q_{\beta}$. The following result is immediate from the definition of the product thickening:

\begin{lem} \label{lem:Kuranishi_functor_Ham_multiplicative}
  There is a canonical isomorphism of Kuranishi charts 
  \begin{equation} \label{eq:product_map_Kuranishi_charts_Hamiltonian}
    \bM_{\alpha}(a,b) \times \bM_{\beta}(b,c)    \to   \bM_{\alpha \times \beta}(a,c),
  \end{equation}
  which defines a natural isomorphism of functors 
  in the following diagram
  \begin{equation}
    \begin{tikzcd}
      A(a,b) \times A(b,c) \ar[r] \ar[d] & A(a,b) \ar[d] \\
      \Kur \times \Kur \ar[r] & \Kur.
    \end{tikzcd}
  \end{equation}
  \qed 
\end{lem}
Lemma \ref{lem:Kuranishi_functor_Ham_multiplicative} provides the
remaining data for a Kuranishi flow category indexed by
$A(f,H,f)$.  The $2$-associativity property again follows from the
corresponding property for disjoint unions of sets, the condition that
strata are products is implied by  Lemma
\ref{lem:Kuranishi_functor_Ham_multiplicative}, and equivariance
follows from the naturality of our construction.  We summarize this
discussion (and recall Definition~\ref{def:Kuranishi-flow-category}): 

\begin{lem}
The topological flow category $\fcMR(f,H,f)$ associated to a non-degenerate Hamiltonian, an almost complex structure, and a Morse function lifts to a $\Pi$-equivariant Kuranishi flow category $\bM$, with domain the $2$-category $A(f,H,f)$. \qed
\end{lem}

\subsection{Flag smooth presentation}
\label{sec:fibr-smooth-pres}
We now lift the Kuranishi presentation above to a flag smooth presentation (see Section \ref{sec:fibr-smooth-struct}, and in particular Definition \ref{def:smoothly_fibred_complete}).

We start by defining the smooth parameter space which serves as the base for all charts that we consider. Thus, given elements $(x,y,z)$ of $\cP(f)$ and  $(p,q)$ of $\cP(H)$, and for an object $\alpha$ of the appropriate category of Floer thickening data, we define
\begin{align} \label{eq:base_moduli_space_Morse}
  \begin{rcases}
    &   \cB_{\alpha}(x_-,y_-) \\
    & \cB_{\alpha}(x_0,y_0) \\
  & \cB_{\alpha}(x_+,y_+)   
  \end{rcases}
  & \equiv  \partial^{Q_\alpha}\barT(x,y) \\
  \cB_{\alpha}(p,q)  & \equiv \fcMR_{\alpha} \\
\cB_{\alpha}(x_-,p)  & \equiv  \partial^{Q_\alpha (f)_-} \barT([x],M) \times  \fcMR_{\alpha} \\
\cB_{\alpha}(p,z_+) & \equiv \fcMR_{\alpha}  \times  \partial^{Q_\alpha (f)_+}  \barT(M,[z]) \\
\cB_{\alpha}(y_0,z_+) & \equiv \partial^{Q_\alpha (f)_0} \barT([y],M) \times \fcMR_{\alpha} \times \partial^{Q_\alpha (f)_+}  \barT(M,[z])  \\ \label{eq:base_moduli_space_minus-to-0}
\cB(x_-,y_0) & \equiv   \barT([x],M,[y]).
\end{align}
Finally, for a pair of elements of $\cP(f)_-$ and $\cP(f_0)$, we define
\begin{multline}
  \cB_{\alpha}(x_-,z_+) \equiv  \partial^{Q_\alpha (f)_-} \barT([x],M) \times \fcMR_{\alpha} \times_M \partial^{Q_\alpha (f)_+}  \barT(M,[z]) \label{eq:base-chart-minus-to-+}\\ 
 \hspace{-.1in}\cup \partial^{Q_\alpha(f)_- \cup \Q_{\alpha}(f)_0} \barT([x],M,M) \times  \partial^{M_0} \fcMR_{\alpha} \times  \partial^{Q_\alpha (f)_+}  \barT(M,[z]),
\end{multline}
where the union 
is taken along the common boundary stratum
\begin{equation}
  \partial^{Q_\alpha (f)_-}    \barT([x],M) \times  \partial^{M_0}  \fcMR_{\alpha} \times \partial^{Q_\alpha (f)_+} \barT(M,[z]). 
 \end{equation}
 \begin{lem}
For each pair $(a,b)$ in $\cP(f,H,f)$, the space $\cB_\alpha(a,b) $ is a $G_\alpha$-equivariant $\langle \partial^{Q_\alpha} \cP(a,b) \rangle$-smooth manifold.
 \end{lem}
 \begin{proof}
   This follows immediately form the smoothness of the moduli spaces $ \fcMR_{\alpha}$, and the moduli spaces of gradient flow lines. The only exceptional situation is that of $(a,b) = (x_-, z_+)$, where our definition is piecewise, and one would in principle only obtain a smooth manifold by choosing a smoothing along the common boundary stratum. Nonetheless, since $ \fcMR_{\alpha} $ and $\barT([x],M,M) $ are respectively canonically locally diffeomorphic to the products of $\partial^{M_0}\fcMR_{\alpha}$ and $\barT([x],M) $ with an interval near these boundary strata, the smooth structure on their union is canonical. 
 \end{proof}
 We have a projection map
\begin{equation}
  \fcMR_{\alpha}(a,b) \to \cB_{\alpha} 
\end{equation}
which is $G_{\alpha}$-equivariant and forgets the map $u$. Moreover, replacing the moduli spaces $\fcMR_{\alpha}$ in Equations \eqref{eq:base_moduli_space_Morse}--\eqref{eq:base-chart-minus-to-+} by $\Mbar^{\prime \bR}_{\alpha}$, we obtain a moduli space  $\cB'_\alpha$, with a smooth stratified $G_\alpha$ equivariant submersion
\begin{equation}
  \cB_{\alpha} \to \cB'_\alpha.
\end{equation}

\begin{lem} \label{lem:fibrewise_smooth_Hamiltonian_chart}
  For each $\alpha \in A(a,b)$, the forgetful maps
  \begin{equation}
    \fcMRreg_{\alpha}(a,b) \to \cB_\alpha \to \cB'_\alpha 
  \end{equation}
 are equipped with a natural relative smooth structure, lifting $ \bM_{\alpha}(a,b)$ to a flag smooth Kuranishi chart. Each map $ f \co \underline{\alpha} \to \underline{\beta}$ in $\cD(a,b)$ lifts to the category of flag smooth charts.
\end{lem}
\begin{proof}
  By definition, $\fcMRreg_{\alpha}(a,b) $,  $\cB_{\alpha}$ and $\cB'_{\alpha}  $ are stratified by $\partial^{Q_\alpha} \cP(a,b)$. The forgetful map from  $\fcMRreg_{\alpha}(a,b) $  to the abstract moduli spaces is $G_\alpha$-equivariant, and is compatible with the stratification. The moduli spaces $\cB_{\alpha} $ and $\cB'_{\alpha}  $ are smooth, and the forgetful map  $\cB_{\alpha}  \to \cB'_{\alpha}  $ is a smooth submersion by Proposition \ref{prop:moduli-space-compact-smooth-Haussdorff}.  Proposition \ref{prop:Ham_Kuranishi_chart} immediately implies that the map to $\cB_{\alpha}$ is a  $G_{\alpha}$-equivariant topological submersion.

  Since the domain of an element of $\fcMRreg_{\alpha}(a,b)$ is identified with its image in $\cB'_{\alpha}$, we can repeat the proof of Proposition \ref{prop:Ham_Kuranishi_chart} using the smaller moduli space $\cB'_{\alpha}$ as our base: this implies that the restriction of this projection to each stratum is equipped with a relative smooth structure, obtained as the zero-locus of a Fredholm bundle over a Banach manifold, whose base contains $V_\alpha$ as a factor. The gluing theorem for moduli spaces of pseudo-holomorphic curves is smooth for fixed parameter (see, e.g. \cite[Theorem 10.1.2]{McDuffSalamon2012}, or the universal characterisation in \cite{Swaminathan2019} for the analogous problem in Gromov-Witten theory), which gives rise to a global relative smooth structure over $\cB'_{\alpha}$, with the property that the map to $V_\alpha$ is smooth.

  On a fibre of the projection map $ \fcMRreg_{\alpha}(a,b) \to \cB'_{\alpha}$, the map to $\cB_{\alpha}$ is obtained by recording the position in the domain of inverse images of real codimension-$2$ submanifolds in $M$, hence is smooth by the implicit function theorem. It is a submersion onto the fibres of $\cB_{\alpha} \to \cB'_{\alpha}$ because of the assumption that the corresponding curves with constraints in $ \cB_{\alpha}$ are regular.

  Given a morphism  $f \co \alpha \to \beta$, the assumption that $V_\alpha $ surjects onto the cokernel of the linearised Cauchy-Riemann operator implies that the restriction of the map $ \fcMRreg_{\beta}(a,b) \to V_\beta$ to each fibre over $ \cB_{\beta}$ is (smoothly) transverse to $V_\alpha$ along the image of $\fcMRreg_{\alpha}(a,b) $. Finally, the two smooth structures on the fibres of the map $ \fcMRreg_{\alpha}(a,b) \to \cB'_{\beta}$ agree because they can be described in terms of a zero-locus of a Fredholm operator. This completes the construction of a lift to the category of flag-smooth charts.

\end{proof}

Combining the above result with Lemma \ref{lem:submersion-floer-cylinder-forgetful-map}, we conclude:
\begin{cor}
The projection maps $ \fcMRreg_{\alpha}(a,b) \to \cB_{\alpha} \to \cB'_{\alpha}$ lift $\bM(a,b)$ to a flag smooth Kuranishi presentation of $ \fcMR(a,b)$. \qed 
\end{cor}

This construction is compatible with products of charts: given elements $ \alpha \in A(a,b)$ and $\beta \in A(b,c) $, we have a 
commutative diagram
\begin{equation}
   \begin{tikzcd}
   \fcMR_{\alpha}(a,b)  \times \fcMR_{\beta}(b,c) \ar[r]  \ar[d] & \fcMR_{\alpha \times \beta}(a,c) \ar[d] \\
   \cB_{\alpha}  \times \cB_{\beta}  \ar[r]  & \cB_{\alpha \times \beta} 
  \end{tikzcd}
\end{equation}
in which the horizontal arrows are homeomorphisms. 
Since the
above construction is functorial in $\alpha$ and $\beta$, the
associated strict unitality and associativity diagrams commute, and
$\Pi$ strictly acts, we conclude (recall
Definition~\ref{def:flag-smooth-presentation}):

\begin{lem}
The Kuranishi flow category $\bM$ lifts to a flag smooth Kuranishi flow category. \qed
\end{lem}

\section{Stable complex structures}
\label{sec:stable-almost-compl-1}
In this section, we lift the flag-smooth Kuranishi flow category constructed in the previous section to a complex-oriented flow category. The construction proceeds in three steps: we first construct a stable almost complex structure on the abstract moduli spaces of Floer cylinders $\fcMR_\alpha$ and a stable framing of the moduli spaces of Morse trajectories, then construct a stable almost complex structure on the fibrewise tangent space of the projection from the moduli spaces of Floer trajectories $\fcMR_\alpha(p,q)$ to $\cB_\alpha$, and finally combine the two to obtain the desired stable almost complex structure on the tangent space of each chart, which we define to be the direct sum of this fibrewise tangent space with the pullback of $T\cB_\alpha$.

\subsection{Stable complex structures on abstract moduli spaces}
\label{sec:stable-compl-struct}

We assume that we are in the non-trivial situations studied in  Sections \ref{sec:constr-moduli-space}, namely, we consider a pair $(a,b) \in \cP(f,H,f)$ such that  $\underline{a}$  and $ \underline{b}$ are distinct elements of $\underline{\cP}(f,H,f)$, and we are not in the situation $\underline{a} = M_-$ and $\underline{b} = M_0$, or $\underline{a} = M_0$ and $\underline{b} = M_+$. This ensures that the moduli spaces $\fcMR_{\underline{\alpha}}$ are not defined to be a point for trivial reasons.  We shall in particular use the description of tangent spaces in Proposition \ref{prop:moduli-space-compact-smooth-Haussdorff} and Lemma \ref{lem:submersion-floer-cylinder-forgetful-map}. We shall apply these results to the subsets $\fcMR_\alpha \subset \fcMR_{\underline{\alpha}} $ considered in Diagram \eqref{eq:inclusion_open_subset_stratum-abstract-moduli-discs}. For this purpose, we write
\begin{align}
  \Mbar^{S,\bR}_\alpha & \subset  \Mbar^{S,\bR}_{\underline{\alpha}} 
\end{align}
for the inverse image of $\fcMR_\alpha $ in $ \Mbar^{S,\bR}_{\underline{\alpha}}  $. This is a smooth manifold with corners, which submerses over $\fcMR_\alpha$. We write
\begin{equation}
  \Rbar_\alpha  \subset \Rbar_{\underline{\alpha}}
\end{equation}
for the image of this manifold under the submersion $ \Mbar^{S,\bR}_{\underline{\alpha}} \to  \Rbar_{\underline{\alpha}}$. This is a stratum of $\Rbar_{\underline{\alpha}} $, hence again a smooth manifold with corners. 

\begin{defin}
The space $\scrJ_{\base}(\alpha)$ consists of the following data:
\begin{enumerate}
\item a $G_\alpha$-invariant inner product on $T \Mbar^{S,\bR}_\alpha$ which is   \begin{equation}
    \label{eq:conditions-metric-on-moduli-with-real-marked-points}
    \parbox{30em}{induced from  its restriction to $ T^{\fcMR}\Mbar^{S,\bR}_\alpha \cong \bR^{S_\alpha^H}$ and $T^{\Rbar}\Mbar^{S,\bR}_\alpha$, via Equation \eqref{eq:moduli-space-with-real-points-two-kernels-span}. We assume that the first restriction splits orthogonally as direct sum of the factors of $\bR^{S_\alpha^H}$, and the second restriction is (the real part of) a hermitian inner product on $T^{\Rbar}\Mbar^{S,\bR}_\alpha $.} 
  \end{equation}
\item an isomorphism $ \ell_{\underline{b} } \oplus \bR^{\underline{Q}_\alpha} \oplus T \Rbar_\alpha \cong \bR^{S_\alpha^H \amalg \{\chi_a^b\}} $ of vector bundles over $\Rbar_\alpha$, where $\ell_{\underline{b}}$ is a trivialised real line associated to $\underline{b}$, and
\item a $G_\alpha$-equivariant section of the projection map $\Mbar^{S^H,\bR}_\alpha \to \fcMR_{\alpha}$.

\end{enumerate}
\end{defin}

We begin by noting that a stable complex structure on $\fcMR_\alpha$ may be obtained from the space $\scrJ_{\base}(\alpha)$.  Indeed, the inner product determines an isomorphism
\begin{align}
    T \fcMR_\alpha \oplus  \bR^{S_\alpha^H}  \cong  T  \Mbar^{S^H,\bR}_\alpha \cong  T \Rbar_\alpha \oplus T^{\Rbar}  \Mbar^{S^H,\bR}_\alpha 
  \end{align}
  of vector bundles over $\Mbar^{S^H,\bR}_\alpha$. Taking the direct sum with $\ell_{\underline{b}} \oplus \bR^{\underline{Q}_\alpha}$, and using the trivialisation of  $\bR \oplus T \Rbar_\alpha $ yields an isomorphism
\begin{equation} \label{eq:basic_stable_isomorphism_abstract}
  \ell_{\underline{b}} \oplus \bR^{\underline{Q}_\alpha} \oplus T \fcMR_\alpha \oplus  \bR^{S_\alpha^H} \cong \bR^{S_\alpha^H \amalg \{\chi_a^b\}} \oplus T^{\Rbar} \Mbar^{S^H,\bR}_\alpha
\end{equation}
over $\Mbar^{S^H,\bR}_\alpha$, and the choice of section yields an isomorphism of vector bundles over $\fcMR_\alpha$.

Before proceeding further, we note the following:
\begin{lem}
 $\scrJ_{\base}(\alpha)  $ is contractible. 
\end{lem}
\begin{proof}
  This is straightforward to see for the choice of framing of $\Rbar_\alpha$ because this manifold is contractible, and for the choice of  section of the projection from $\Mbar^{S,\bR}_{\alpha} $ because the fibre is contractible and the isotropy group acts trivially on the fibre over any point with non-trivial stabiliser. Regarding the inner product, recall that the kernels $T^{\Rbar} \Mbar^{S,\bR}_\alpha $ and $T^{\fcMR}\Mbar^{S,\bR}_\alpha $ span the tangent space $T \Mbar^{S,\bR}_\alpha$, and that their intersection is canonically isomorphic to $\bR$, corresponding to the diagonal inclusion in $\bR^{S_\alpha^H} $ and to the (real) translation of all components in $T^{\Rbar} \Mbar^{S,\bR}_\alpha$.

To describe the space of inner products, start with the real part of a hermitian inner product on $T^{\Rbar} \Mbar^{S,\bR}_\alpha$; this determines a metric on the subspace spanned by translation of all components. A choice of extension to $\bR^{S_\alpha^H} $ compatible with the inner product is a fibre of the map $\bR_+^{S_\alpha^H} \to \bR_+$ given by adding all components, which is evidently contractible. This completes the construction.
\end{proof}

We conclude:
\begin{lem} 
Every element of $\scrJ_{\base}(\alpha)$ determines an isomorphism
\begin{equation} \label{eq:stable-almost-complex-structure-base}
  \ell_{\underline{b}} 
  \oplus T \fcMR_\alpha \oplus  \bR^{\underline{Q}_\alpha \amalg S_\alpha^H} \cong \bR^{\{\chi_a^b\} \amalg S_\alpha^H} 
  \oplus T^{\Rbar} \Mbar^{S,\bR}_\alpha
\end{equation}
of vector bundles over $\fcMR_{\alpha}$. \qed  
\end{lem}

Recall that, for each arrow ${f} \co {\alpha} \to {\beta}$ in $A(a,b)$, the map of moduli spaces takes $\fcMR_{\alpha}$ to $\fcMR_{{\beta}}$. It similarly maps $\Rbar_\alpha$ to $\Rbar_{\beta}$, and $\Mbar^{S^H,\bR}_\alpha$ to $\Mbar^{S^H,\bR}_\beta$. Since these are all open subsets of strata of the moduli spaces studied in Section \ref{sec:constr-moduli-space}, Lemma \ref{lem:real_moduli_spaces_no_marked_points} and Proposition \ref{prop:moduli-space-compact-smooth-Haussdorff} describe the fibres of the induced maps of tangent space. Recalling that the strata are prescribed by the choice of subsets $Q_\beta \subset Q_\alpha$, we have:
\begin{lem}
  The projections to the cokernels of the maps of tangent spaces
  \begin{align}
    T \fcMR_{\alpha} & \to T \fcMR_{{\beta}} \\ \label{eq:map_tangent_spaces_Rbar}
    T\Rbar_\alpha & \to T \Rbar_{\beta} \\
    T \Mbar^{S^H,\bR}_\alpha &\to  T\Mbar^{S^H,\bR}_\beta
  \end{align}
  are canonically split, and admit natural isomorphisms to $\bR^{\underline{Q}_\alpha
  \setminus \underline{Q}_\beta}$ which are uniquely determined up to positive rescaling of each factor. \qed
\end{lem}

\begin{defin} \label{def:compatible_ac_base}
For each arrow ${f} \co {\alpha} \to {\beta}$, the
space  $\scrJ_{\base}(f)$ consists of a pair of elements of $\scrJ_{\base}(\alpha)$ and
$\scrJ_{\base}(\beta)$ such that: 
\begin{enumerate}
\item the projection map $ T \Mbar^{S^H,\bR}_{\alpha} \to  T \Mbar^{S^H,\bR}_{\beta}  $  is orthogonal, and 
\begin{equation}
  \parbox{30em}{the induced metric on the kernel $T^\beta \Mbar^{S^H,\bR}_{\alpha} $  splits orthogonally as a direct sum of the factors of $\bR^{S_\alpha^H \setminus S_\beta^H}$  and  $\bigoplus T_{z} \Sigma$,
and the induced metric on the cokernel respects the decomposition of $\bR^{\underline{Q}_\alpha
  \setminus \underline{Q}_\beta} $ into summands.}
\end{equation}
\item there is an identification of the cokernel in Equation \eqref{eq:map_tangent_spaces_Rbar} with  $\bR^{\underline{Q}_\alpha
  \setminus \underline{Q}_\beta}$, in the prescribed class, such that the induced map to $T \Rbar_\beta $ yields a commutative diagram
\begin{equation}\label{eq:splitting-map-of-tangent-spaces-Stasheff-compatible}
  \begin{tikzcd}
     \bR^{S_\alpha^H \setminus S_\beta^H} \ar[r] \ar[d] &  \bR^{S_\alpha^H} \ar[r] \ar[d] &   \bR^{S_\beta^H } \ar[d]  \\ 
 \bR^{S_\alpha^H \setminus S_\beta^H} \ar[r]  &   \ell_{\underline{b}} \oplus \bR^{\underline{Q}_\alpha}   \oplus T \Rbar_\alpha \ar[r] &  \ell_{\underline{b}} \oplus \bR^{\underline{Q}_\beta}   \oplus T \Rbar_\beta,
  \end{tikzcd}
\end{equation}
\item the sections gives rise to a commutative diagram
\begin{equation}
  \begin{tikzcd}
    \fcMR_{\alpha} \ar[d] \ar[r] &   \fcMR_{{\beta}} \ar[d] \\
    \Mbar^{S^H,\bR}_{\alpha} \ar[r] & \Mbar^{S^H,\bR}_{\beta}.
  \end{tikzcd}
\end{equation}
\end{enumerate}
We equip $\scrJ_{\base}(f) $ with the subspace topology.
\end{defin}

To state the associativity of this construction, we introduce a
category internal to the category of topological
spaces (see Appendix~\ref{sec:internal} for a brief review of the
definitions of internal categories and internal functors that we use):

\begin{defin}
The category $A^{\base}(a,b)$ has:
  \begin{itemize}
  \item objects given by the disjoint union, over $\alpha \in A(a,b)$, of the spaces $\scrJ_{\base}(\alpha)$.
  \item morphisms given by the disjoint union, over morphisms $f \in A(a,b)$, 
    of the spaces $  \scrJ_{\base}(f)$.
    \item composition given by the natural continuous map
  \begin{equation}
     \scrJ_{\base}(f) \times_{\scrJ_{\base}(\beta)} \scrJ_{\base}(g) \to \scrJ_{\base}(g \circ f),
   \end{equation}
   obtained from the fact that composition is compatible with commutative squares and orthogonal projections.
  \end{itemize}
\end{defin}

Having established earlier that the object spaces of this category are contractible, we now consider morphisms. For the next result, we use the fact that the kernel $ T^\beta \fcMR_\alpha  $ of the projection map from $ T \fcMR_\alpha $ to $T \fcMR_\beta $, is naturally isomorphic to the kernel of the projection map from $ T^{\Rbar} \Mbar^{S^H,\bR}_\alpha $ to $ T^{\Rbar} \Mbar^{S^H,\bR}_\beta  $, as can be seen by noting that both are given by the direct sum over all points in $S_\alpha^H \setminus S_\beta^H$ of the tangent space at the underlying Riemann surface (see Lemma \ref{lem:submersion-floer-cylinder-forgetful-map}).
\begin{lem}
  The space $\scrJ_{\base}(f) $ is contractible, and
  each element induces a commutative diagram
  \begin{equation} 
    \begin{tikzcd}
      T^\beta \fcMR_\alpha \ar[r] \ar[d, dashed, bend right] &      T \fcMR_\alpha \oplus \bR^{\underline{Q}_\alpha \setminus \underline{Q}_\beta}  \ar[r] \ar[d, dashed, bend right]  &  T \fcMR_\beta  \ar[l,dashed,bend right]  \ar[d, dashed, bend right] \\
 T^\beta \Mbar^{S^H,\bR}_\alpha  \ar[r] \ar[u] \ar[d] &  T \Mbar^{S^H,\bR}_\alpha  \oplus \bR^{\underline{Q}_\alpha \setminus \underline{Q}_\beta} \ar[r] \ar[u] \ar[d] & T\Mbar^{S^H,\bR}_\beta  \ar[l,dashed,bend right]   \ar[u] \ar[d] \\
T^\beta \Rbar_\alpha \ar[r]  &  T \Rbar_\alpha  \ar[r]   & T \Rbar_\beta  \ar[l,dashed,bend right]         \end{tikzcd}
  \end{equation}
  of maps vector bundles over $\fcMR_{\alpha}$, in which the rows are split exact sequences. Moreover, the induced splitting of the sequence
  \begin{equation} \label{eq:exact_sequence_abstract-moduli-stable-complex} 
    \begin{tikzcd}
T^\beta \fcMR_\alpha \ar[r]  &  T^{\Rbar} \Mbar^{S^H,\bR}_\alpha  \ar[r]  & T^{\Rbar} \Mbar^{S^H,\bR}_\beta  \ar[l,dashed,bend right]   
    \end{tikzcd}
  \end{equation}
  respects the complex structure. \qed
\end{lem}
  \begin{cor}
     The following diagram, in which the horizontal maps are induced by the stable complex structures on $\fcMR_{\alpha}$ and on $\fcMR_{{\beta}}$, is commutative:
\begin{equation} 
  \begin{tikzcd}[row sep=small, column sep=tiny, transform shape, nodes={scale=.75}]
  \ell_{\underline{b}} \oplus 
  T \fcMR_\beta \oplus  \bR^{\underline{Q}_\beta \amalg S_\beta^H } \oplus \bR^{S_\alpha^H \setminus S_\beta^H} \oplus T^\beta \fcMR_\alpha    \ar[r] \ar[d] & T^{\Rbar} \Mbar^{S^H,\bR}_\beta  \oplus \bR^{ \{\chi_a^b \} \amalg S_\beta^H }
  \oplus \bR^{S_\alpha^H \setminus S_\beta^H} \oplus  T^\beta \fcMR_\alpha   \ar[d] \\
   \ell_{\underline{b}} 
   \oplus T \fcMR_\alpha \oplus  \bR^{\underline{Q}_\alpha \amalg  S_\alpha^H} \ar[r]  & T^{\Rbar}  \Mbar^{S^H,\bR}_\alpha   \oplus \bR^{S_\alpha^H  \amalg \{\chi_a^b \}}. 
  \end{tikzcd}
\end{equation} \qed
  \end{cor}

For the next result, we note that a splitting of the surjections $T \fcMR_\alpha \to T \fcMR_\beta$ and $T \fcMR_\beta \to     T \fcMR_\gamma $ induces a splitting of the kernel of the composite map 
   \begin{equation} \label{eq:splitting-kernels-moduli-spaces-with-extra-real-points}
     T^\gamma \fcMR_\alpha \cong     T^\beta \fcMR_\alpha  \oplus T^\gamma \fcMR_\beta     .
   \end{equation}

\begin{lem}
The splitting in Equation \eqref{eq:splitting-kernels-moduli-spaces-with-extra-real-points} associated to each element of the fibre product $ \scrJ_{\base}(f) \times_{\scrJ_{\base}(\beta)} \scrJ_{\base}(g) $ gives rise to a commutative diagram:
  \begin{equation}
     \begin{tikzcd}
       T^\beta \fcMR_\alpha \oplus  T^\gamma \fcMR_\beta \oplus   T \fcMR_\gamma \ar[r] \ar[d] &  T^\beta \fcMR_\alpha  \oplus T \fcMR_\beta  \oplus \bR^{\underline{Q}_\beta \setminus \underline{Q}_\gamma} \ar[d] \\
  T^\gamma \fcMR_\alpha  \oplus T \fcMR_\gamma \ar[r] & T \fcMR_\alpha \oplus \bR^{\underline{Q}_\alpha \setminus \underline{Q}_\gamma},
\end{tikzcd}
\end{equation}
and similarly for the tangent spaces of $\Rbar_\alpha$ and $ \Mbar^{S^H,\bR}_\alpha$. \qed
 \end{lem}
\begin{cor}
The stable isomorphism between the stable complex structures on $\fcMR_{\alpha}$ and $\fcMR_{{\gamma}} $ associated to the composition $g \circ f$ coincides with the composition of the stable isomorphisms associated to $g$ and $f$.  \qed
\end{cor}

We now consider the multiplicativity of this construction: recall that, given a pair of inclusions $\{\fcMR_{{\alpha}_i} \subset  \fcMR_{\underline{\alpha}_i}\}_{i=1}^{2} $, with $\alpha_1 \in A(a,b)$ and $\alpha_2 \in A(b,c)$,  we have defined $  \fcMR_{{\alpha}_1 \times {\alpha}_2} $ to be the product of these spaces. We similarly define $\Rbar_{{\alpha}_1 \times {\alpha}_2}$ and $ \Mbar^{S,\bR}_{{\alpha}_1 \times {\alpha}_2}$ to be the product of the corresponding moduli spaces.

There is a natural map 
\begin{equation}
\scrJ_{\base}( {{\alpha}_1}) \times \scrJ_{\base}( {{\alpha}_2}) \to
\scrJ_{\base}( {{\alpha}_1 \times {\alpha}_2})
\end{equation}
defined as the direct sums of: 
\begin{enumerate}
\item The inner products on  $ T
\Mbar^{S,\bR}_{{\alpha}_i}$, 
\item the stable framings of  $\Rbar_{
  {\alpha}_i}$, and 
\item the sections of the projections $\Mbar^{S,\bR}_{{\alpha}_i} \to
  \fcMR_{{\alpha}_i} $. 
\end{enumerate}
In particular, for the framings, we note that 
\begin{equation}
  Q_{\alpha_1} \amalg \{ b \} \amalg Q_{\alpha_2} \amalg \{c\} \equiv
  Q_{\alpha_1 \times \alpha_2} \amalg \{c\},
\end{equation}
so we obtain the desired stable framing of the product.

\begin{lem} \label{lem:multiplicativity_base-stable-complex}
The stable almost complex structures induced by elements of $\scrJ_{\base}( {{\alpha}_1})$, $\scrJ_{\base}( {{\alpha}_2})$, and their image in $\scrJ_{\base}( {{\alpha}_1 \times {\alpha}_1})$ fit in a commutative diagram
\begin{equation}
  \begin{tikzcd}
 \ell_{\underline{c}} \oplus   \bR^{\underline{Q}_{\alpha_1 \times \alpha_2} 
   \amalg S^H_{\alpha_1\times \alpha_2} } \oplus T  \fcMR_{{\alpha_1 \times \alpha_2}}   \ar[r] \ar[d] & T^{\Rbar}  \fcMR_{\alpha_1 \times \alpha_2}  \oplus  \bR^{\{\chi_a^c\} \amalg S^H_{\alpha_1\times \alpha_2} }
 \ar[d] \\
 \begin{gathered}
   \ell_{\underline{b}} \oplus    \bR^{\underline{Q}_{\alpha_1} 
     \amalg S^H_{\alpha_1} }  \oplus   T \fcMR_{\alpha_1} \oplus  \\ \ell_{\underline{c}} \oplus \bR^{\underline{Q}_{\alpha_2} 
     \amalg S^H_{\alpha_2} }  \oplus     T \fcMR_{\alpha_2} 
 \end{gathered}
       \ar[r] &
       \begin{gathered}
        T^{\Rbar}  \fcMR_{\alpha_1} \oplus  \bR^{\{\chi_a^b\} \amalg S^H_{\alpha_1} } \\
   \oplus  T^{\Rbar}  \fcMR_{\alpha_2}    \oplus  \bR^{\{\chi_b^c\} \amalg S^H_{\alpha_2}}.  
       \end{gathered}
  \end{tikzcd}
\end{equation}
\qed
\end{lem}

Next, we state the associativity of the product construction: to start, we observe that, given a triple $ \{\fcMR_{{\alpha}_i} \subset   \fcMR_{\underline{\alpha}_i}\}_{i=1}^{3}$,  the products
\begin{align}
  \fcMR_{{\alpha}_1 \times {\alpha}_2 \times {\alpha}_3}  & \subset \fcMR_{\underline{\alpha}_1 \times \underline{\alpha}_2 \times \underline{\alpha}_3} \\
  \Mbar^{S^H, \bR}_{{\alpha}_1 \times {\alpha}_2 \times {\alpha}_3}  & \subset \Mbar^{S^H, \bR}_{\underline{\alpha}_1 \times \underline{\alpha}_2 \times \underline{\alpha}_3} \\
  \Rbar_{{\alpha}_1 \times {\alpha}_2 \times {\alpha}_3}  & \subset \Rbar_{\underline{\alpha}_1 \times \underline{\alpha}_2 \times \underline{\alpha}_3}
\end{align}
are independent of parenthesisation. The fact that the comparison of orientations on moduli spaces and their products is associative then amounts to:
\begin{lem}
The following diagram commutes:
  \begin{equation}
    \begin{tikzcd}
   \scrJ_{\base}( {{\alpha}_1}) \times \scrJ_{\base}( {{\alpha}_2})\times \scrJ_{\base}( {{\alpha}_3})  \ar[r] \ar[d] &  \scrJ_{\base}( {{\alpha}_1}) \times \scrJ_{\base}( {{\alpha}_2 \times{\alpha}_3})    \ar[d] \\
    \scrJ_{\base}( {{\alpha}_1 \times {\alpha}_2}) \times \scrJ_{\base}( {{\alpha}_3})     \ar[r] & \scrJ_{\base}({{\alpha}_1 \times {\alpha}_2 \times {\alpha}_3} ).
    \end{tikzcd}
  \end{equation} \qed
\end{lem}

Finally, we consider the functoriality of the product maps:   given maps $\{\alpha_i \to \beta_i\}_{i=1}^2$,   and open subsets $\{ \fcMR_{{\alpha}_i} \subset \fcMR_{\underline{\alpha}_i} \}_{i=1}^{2}$ and $\{\fcMR_{{\beta}_i} \subset \fcMR_{\underline{\beta}_i}\}_{i=1}^{2}$ included in the strata labelled by $\underline{Q}_{\alpha_i}$ and $\underline{Q}_{\beta_i}$, such that $\fcMR_{{\alpha}_i}$ maps to $\fcMR_{{\beta}_i}$, we observe that the following diagram
\begin{equation}
  \begin{tikzcd}
       \fcMR_{\alpha_1} \oplus      \fcMR_{\alpha_2}    \ar[r] \ar[d] &  \fcMR_{\beta_1} \oplus  \fcMR_{\beta_2}    \ar[d] \\
\fcMR_{{\alpha_1 \times \alpha_2}}   \ar[r] &   \fcMR_{{\beta_1 \times \beta_2}} 
  \end{tikzcd}  
\end{equation}
commutes, as do the corresponding diagrams
\begin{equation}
  \begin{tikzcd}
       \Mbar^{S^H, \bR}_{\alpha_1} \oplus      \Mbar^{S^H, \bR}_{\alpha_2}    \ar[r] \ar[d] &  \Mbar^{S^H, \bR}_{\beta_1} \oplus  \Mbar^{S^H, \bR}_{\beta_2}    \ar[d] \\
\Mbar^{S^H, \bR}_{{\alpha_1 \times \alpha_2}}   \ar[r] &   \Mbar^{S^H, \bR}_{{\beta_1 \times \beta_2}} 
  \end{tikzcd}  
\end{equation}
and
\begin{equation}
  \begin{tikzcd}
       \Rbar_{\alpha_1} \oplus      \Rbar_{\alpha_2}    \ar[r] \ar[d] &  \Rbar_{\beta_1} \oplus  \Rbar_{\beta_2}    \ar[d] \\
\Rbar_{{\alpha_1 \times \alpha_2}}   \ar[r] &   \Rbar_{{\beta_1 \times \beta_2}} .
  \end{tikzcd}  
\end{equation}
Moreover, the identifications of the cokernel of each of these maps with
\begin{equation}\label{eq:cokernel_agree}
  \bR^{\underline{Q}_{\alpha_1 \times \alpha_2} \setminus \underline{Q}_{\beta_1 \times \beta_2}}
\end{equation}
agree (up to the same real positive dilation ambiguity in each factor).

The functoriality of the product maps allows us to pass to
arrows: 
\begin{lem}
   Given a pair of arrows ${f}_1 \co {\alpha}_1 \to {\beta}_1$ and  ${f}_2 \co {\alpha}_2 \to {\beta}_2$, there is a natural map
  \begin{equation}
     \scrJ_{\base}({{f}_1}) \times    \scrJ_{\base}({{f}_2}) \to \scrJ_{\base}({{f}_1 \times {f}_2}),
  \end{equation}
  which is functorial in the sense that, for each pair of compositions $ {g}_1 \circ  {f}_1$ and $ {g}_2 \circ  {f}_2$, the following diagram commutes:
  \begin{equation}
    \begin{tikzcd}[row sep=small, column sep=tiny, transform shape, nodes={scale=.9}]
      \begin{gathered}
         \scrJ_{\base}({{f}_1}) \times_{\scrJ_{\base}({{\beta}_1})} \scrJ_{\base}({{g}_1}) \\ \times    \scrJ_{\base}({{f}_2}) \times_{\scrJ_{\base}({{\beta}_2})} \scrJ_{\base}({{g}_2}) 
      \end{gathered}
      \ar[d] \ar[r] &        \scrJ_{\base}({{g}_1 \circ {f}_1}) \times  \scrJ_{\base}({{g}_2 \circ {f}_2}) \ar[d] \\
 \scrJ_{\base}({{f}_1 \times {f}_2})  \times_{\scrJ_{\base}({{\beta}_1 \times {\beta}_2})}  \scrJ_{\base}({{g}_1 \times {g}_2})  \ar[r]  &    \scrJ_{\base}({{g}_1 \circ {f}_1 \times {g}_2 \circ {f}_2}).
    \end{tikzcd}
  \end{equation}
Moreover, this construction is multiplicative in the sense that the following diagram commutes:
  \begin{equation}
    \begin{tikzcd}
   \scrJ_{\base}( {{f}_1}) \times \scrJ_{\base}( {{f}_2})\times \scrJ_{\base}( {{f}_3})  \ar[r] \ar[d] &  \scrJ_{\base}( {{f}_1}) \times \scrJ_{\base}( {{f}_2 \times{f}_3})    \ar[d] \\
    \scrJ_{\base}( {{f}_1 \times {f}_2}) \times \scrJ_{\base}( {{f}_3})     \ar[r] & \scrJ_{\base}({{f}_1 \times {f}_2 \times {f}_3} ).
    \end{tikzcd}
  \end{equation}
  \qed
\end{lem}
Reformulating the above result in more abstract terms, we have:

\begin{cor}
The categories $A^{\base}(a,b)$ are the $1$-cells of a bicategory (in
fact, a $2$-category) internal to topological spaces, with $0$-cells
the elements of $\cP(f,H,f)$.  \qed
\end{cor}

\subsection{Stable complex structures and Morse stable framings}
\label{sec:stable-compl-struct-2}

In this subsection, we combine the constructions of the preceding subsection, with those of Appendix~\ref{sec:stable-fram-moduli}, to construct a stable complex structure on the moduli spaces $\cB_\alpha$, which are the bases of the flag smooth Kuranishi presentations introduced in Section \ref{sec:fibr-smooth-pres}.

If $x$ is a critical point of $f$, let $V_x^+$ denote the positive definite subspace of the Hessian of $f$ at $x$. The starting point is to fix the choice of trivialisation of the ascending manifold at $x$
\begin{equation}
T \Tbar([x],M) \cong V_x^+  
\end{equation}
in Equation \eqref{eq:fix-trivialisation}. As explained in Appendix~\ref{sec:morse-theory-and-homotopy-types}, this induces a stable framing
\begin{equation} \label{eq:stable_framing_moduli_trajectories}
  T \barT(x,y) \oplus \ell_y \oplus V^+_y \cong V_x^+.
\end{equation}

Next, we consider the moduli spaces $\Tbar(M,[x])$. Proceeding by induction on the Morse index of $[y]$, and letting $\ell_M = \bR^{\{M\}}$, we pick stable isomorphisms
\begin{equation} \label{eq:stable_framing_ascending_manifold}
  T \Tbar(M,[x]) \oplus \ell_y \oplus V_x^+ \cong \ell_M \oplus TM
\end{equation}
subject to the following constraint: along the codimension $1$
boundary stratum $\Tbar(M,[x]) \times \Tbar([x],[y])$ of
$\Tbar(M,[y])$, we require that the following diagrams, where the
arrows are given by Equation
\eqref{eq:stable_framing_moduli_trajectories} and
\eqref{eq:stable_framing_ascending_manifold}, are commutative:
\begin{equation}
  \begin{tikzcd}[column sep = small]
    T \Tbar(M,[x]) \oplus \ell_x \oplus T \Tbar([x],[y]) \oplus \ell_y \oplus V^+_y \ar[r] \ar[d] &  T \Tbar(M,[y]) \oplus \ell_y \oplus V^+_y \ar[d] \\
     T \Tbar(M,[x]) \oplus \ell_x \oplus V_{x}^{+} \ar[r] & \ell_M \oplus TM.
  \end{tikzcd}
\end{equation}
\begin{rem}
  Observe that Equation \eqref{eq:stable_framing_ascending_manifold} provides a stable almost complex structure on the moduli spaces $\Tbar(M,[x])$. In fact, these moduli spaces have natural stable framings by the negative-definite subspace of the Hessian. We do not use these stable framings because we have already framed $\Tbar([x],M)$, and these two families of moduli spaces cannot be simultaneous framed, compatibly with Poincar\'e duality, unless $M$ is itself a stably framed manifold.
\end{rem}

Since the moduli spaces $\Tbar([x],M,[y])$ are defined as the fibre products of $\Tbar([x],M)$ and $\Tbar(M,[y])$ over the evaluation map to $M$, the framing of $\Tbar([x],M)$ and the stable isomorphism of $T\Tbar(M,[y])$ with $TM$ yields a stable framing:
\begin{equation}
  T \Tbar([x],M,[y]) \oplus V_y^+ \oplus \ell_y \cong V_x^+ \oplus \ell_M,
\end{equation}
for which it is straightforward to check that we have a commutative diagram \begin{equation}
  \begin{tikzcd}
    \begin{gathered}
       T \Tbar([x],M,[y]) \oplus \ell_y  \oplus \\ T \Tbar([y], [z]) \oplus V_z^+ \oplus \ell_z
    \end{gathered}
    \ar[r] \ar[d] &    T \Tbar([x],M,[z]) \oplus V_z^+ \oplus \ell_z \ar[d] \\
  T \Tbar([x],M,[y]) \oplus \ell_y  \oplus  V_y^+  \ar[r]  &    V_x^+ \oplus \ell_M
  \end{tikzcd}
\end{equation}
associated to the boundary stratum $\Tbar([x],M,[y]) \oplus T \Tbar([y], [z]) $, and a commutative diagram
\begin{equation}
  \begin{tikzcd}
    \begin{gathered}
      T \Tbar([x],[y]) \oplus \ell_y  \oplus \\
      T \Tbar([y],M,[z]) \oplus V_z^+ \oplus \ell_z
    \end{gathered}
    \ar[r] \ar[d] &    T \Tbar([x],M,[z]) \oplus V_z^+ \oplus \ell_z \ar[d] \\
  T \Tbar([x],[y]) \oplus \ell_y  \oplus  V_y^+ \oplus \ell_M  \ar[r]  &    V_x^+ \oplus \ell_M
  \end{tikzcd}
\end{equation}
associated to the boundary stratum $\Tbar([x],[y]) \oplus T \Tbar([y],M,[z]) $.

Finally, we consider the moduli space $\Tbar([x],M,M)$, for which we state the following result, whose proof we leave to the reader:
\begin{lem}
There is a choice of stable framings
\begin{equation} \label{eq:framing_moduli_space_Morse-two_marked_points}
  T \Tbar([x],M,M) \cong V_x^+ \oplus \ell_M
 \end{equation}
 subject to the following constraints:
 \begin{itemize}
 \item Along the boundary stratum $T \Tbar([x],M)$, Equation \eqref{eq:framing_moduli_space_Morse-two_marked_points} is obtained from the framing $T \Tbar([x],M) \cong V_x^+$, and the product decomposition
   \begin{equation}
     T \Tbar([x],M,M) \cong T \Tbar([x],M) \oplus \ell_M.
   \end{equation}
 \item Along each boundary stratum $ \Tbar([x],M,[y]) \times  \Tbar([y],M) $, we have a commutative diagram
   \begin{equation}
     \begin{tikzcd}
     T  \Tbar([x],M,[y]) \oplus \ell_y \oplus  \Tbar([y],M) \ar[r] \ar[d] &   T  \Tbar([x],M,M) \ar[d] \\
    T  \Tbar([x],M,[y]) \oplus \ell_y  \oplus V_y^+ \ar[r]   &  V_x^+ \oplus \ell_M.  
     \end{tikzcd}
        \end{equation}
   \item Along each boundary stratum $ \Tbar([x],[y]) \times  \Tbar([y],M,M) $, we have a commutative diagram
   \begin{equation}
     \begin{tikzcd}
     T  \Tbar([x],[y]) \oplus \ell_y \oplus  \Tbar([y],M,M) \ar[r] \ar[d] &   T  \Tbar([x],M,M) \ar[d] \\
    T  \Tbar([x],[y]) \oplus \ell_y  \oplus V_y^+ \oplus \ell_M \ar[r]   &  V_x^+ \oplus \ell_M.  
     \end{tikzcd}
        \end{equation}
 \end{itemize} \qed  
\end{lem}

Finally, we state the desired result concerning the stable complex
structure of the moduli spaces which serve as the base of the flag
smooth Kuranishi structures in Hamiltonian Floer theory. To this end,
we write 

\begin{equation} \label{eq:keep_Morse_space_set_Floer_0}
  V_{\overline{a}}^+ \equiv
  \begin{cases}
    0 & \textrm{if } a \in \cP(H) \\
    V_a^+ & \textrm{otherwise.}
  \end{cases}
\end{equation}

Given $\alpha \in A(a,b)$, we also define

\begin{align}
  O_\alpha & \equiv  Q_\alpha \amalg \{\chi_a^b\} \amalg \{b\} \\
  I^{\base,\bC}_{\alpha} & \equiv   \begin{cases}
    T^{\Rbar}  \Mbar^{S^H,\bR}_\alpha  \oplus \bC^{\{\chi_a^b\}} \oplus TM  & \textrm{ if } b \in \cP(f)_+ \\
    T^{\Rbar}  \Mbar^{S^H,\bR}_\alpha  \oplus \bC^{\{\chi_a^b\}}  &  \textrm{otherwise.}
  \end{cases} \\
  W^\base_\alpha & \equiv \bR^{S_\alpha^H}.
\end{align}

Given an arrow $f$ from $\alpha$ to $\beta$ in $A(a,b)$, we define

\begin{align}
  I^{\base,\bC}_f & \equiv T^\beta \fcMR_\alpha \\
    W^\base_f & \equiv \bR^{S_\alpha^H \setminus S_\beta^H},
\end{align}

and note that we have a natural isomorphism

\begin{equation}
   W^\base_\alpha \cong  W^\base_\beta \oplus  W^\base_f.
\end{equation}
\begin{prop}
  The space $\scrJ_{\base}(\alpha)$ parameterises a continuous family of isomorphisms
  \begin{equation} \label{eq:stable-almost-complex-structure-base-1}
    V_{\overline{b}}^+  
  \oplus T \cB_\alpha \oplus  \bR^{O_\alpha}  \oplus  W^\base_\alpha  \cong  W^\base_\alpha  \oplus I^{\base,\bC}_{\alpha} \oplus V_{\overline{a}}^+.
\end{equation}
Given an arrow $f$ in $A(a,b)$ from $\alpha$ to $\beta$, the space   $\scrJ_{\base}(f)$ parametrises a continuous family of isomorphisms of vector bundles
\begin{align}
  T \cB_\alpha  \oplus \bR^{Q_\alpha \setminus Q_\beta} &  \cong    T \cB_\beta \oplus I^{\base,\bC}_f \\
  I^{\base,\bC}_{\alpha} &  \cong   I^{\base,\bC}_{\beta} \oplus I^{\base,\bC}_f,
\end{align}
which are compatible with composition in the category $A^{\base}(a,b)$. In addition, the following diagram commutes (we use the evident isomorphism $ O_\alpha = O_\beta \amalg (Q_\alpha \setminus Q_\beta)$ in the left vertical arrow):
\begin{equation}\label{eq:base_compat}
  \begin{tikzcd}
    \begin{gathered}
       V_{\overline{b}}^+ 
    \oplus T \cB_\beta \oplus  \bR^{O_\beta}  \oplus  W^\base_\beta \\  \oplus W^\base_f   \oplus I^{\base,\bC}_f
    \end{gathered}
    \ar[r] \ar[d]&
    \begin{gathered}
      W^\base_\beta  \oplus I^{\base,\bC}_{\beta} \oplus V_{\overline{a}}^+ \oplus  W^\base_\beta  \\ \oplus W^\base_f   \oplus I^{\base,\bC}_f 
    \end{gathered}
 \ar[d] \\
    V_{\overline{b}}^+  
  \oplus T \cB_\alpha \oplus  \bR^{O_\alpha}  \oplus  W^\base_\alpha  \ar[r] &  W^\base_\alpha  \oplus I^{\base,\bC}_{\alpha} \oplus V_{\overline{a}}^+.
  \end{tikzcd}
\end{equation}
This data is associative and multiplicative in the sense that it
specifies an enriched $2$-functor to the topological bicategory of
Definition~\ref{defin:bicat_vectorspaces} below. 
\end{prop}

\begin{proof}
We use the isomorphism $\ell_b \cong \bR^{\{b\}}$ to incorporate this factor into $\bR^{O_\alpha}$.    The only non-trivial case to discuss occurs whenever $a \in \cP(f)_-$ and $b \in \cP(f)_+$: we identify the copy of $\ell_M$ appearing in Equation \eqref{eq:framing_moduli_space_Morse-two_marked_points} with $\bR^{\{\chi_-\}}$. 
\end{proof}

We define an auxiliary topological bicategory to articulate the
associativity of the preceding definition; this is a simplified
version of the construction of
Definition~\ref{def:Complex-oriented-Kuranishi-chart}. 

\begin{defin}\label{defin:bicat_vectorspaces}
We define a topologically enriched bicategory $\aC$ with:
\begin{enumerate}
\item objects given by finite-dimensional complex inner-product spaces $V$,
\item morphism categories define as the internal categories $\aC(V_1,
  V_2)$ with the space of objects specified by:
\begin{enumerate}
\item A partially ordered set $\cS$ and a totally-ordered subset $Q_{\alpha} \subset \cS$,  
\item A smooth $\langle \partial^{Q_{\alpha}} \cS \rangle$-manifold
  $B_\alpha$ with action by $G_\alpha$, 
\item A complex vector bundle $I_\alpha$ on $B_\alpha$,
\item A set $O_{\alpha}$ and a complex vector space $W_\alpha$,
\item and an isomorphism
\begin{equation}
V_0 \oplus TB_\alpha \oplus \bR^{O_{\alpha}} \oplus W_\alpha \cong
I_\alpha \oplus W_\alpha \oplus V_1.
\end{equation}
\end{enumerate}
\end{enumerate}
The topology is induced from the topology on the space of
isomorphisms.

The space of morphisms $f \colon \alpha \to \beta$ is specified by the
following data: 
\begin{enumerate}
\item An inclusion $Q_\beta \subset Q_\alpha$ and a surjection $G_\beta \to G_\alpha$, 
\item An equivariant submersion $B_\alpha \to \partial^{Q_\alpha}
  B_\beta$ (and a choice of cokernel identification as in
  Equation~\eqref{eq:cokernel_agree}),
\item A complex vector bundle $I_f$ on $B_\alpha$, and
\item isomorphisms
\begin{equation}
I_\alpha \cong I_\beta \oplus I_f
\end{equation}
and
\begin{equation}
TB_\alpha \oplus \bR^{Q_\alpha \setminus Q_\beta} \cong TB_\beta
\oplus I_f
\end{equation}
such that the analogous diagram to Equation~\eqref{eq:base_compat} commutes.
\end{enumerate}
Composition is specified in the evident way, following the work of
Section~\ref{sec:bicat-kuran-charts}.
\end{defin}

\subsection{Stable fibrewise complex structure}
\label{sec:stable-almost-compl}

In this section, we shall construct stable almost complex structures on moduli spaces of pseudo-holomorphic cylinders with marked points, and with fixed conformal structure. These are the fibres of the projection map from the Kuranishi charts we consider to the abstract moduli space of cylinders with marked points, so we refer to this construction as a \emph{fibrewise (stable almost) complex structure}.

\subsubsection{Fredholm operators associated to orbits}
\label{sec:constr-stable-compl-1}

\begin{notn}
  Given a Riemann surface $\Sigma$, obtained from a closed connected Riemann surface by removing finitely many points, and a map $u \co \Sigma \to M$, we write $\cF(u^* TM)$ for sections of the pullback of $TM$, and $\cE^{0,1}(u^* TM)$ for complex anti-linear $1$-forms on $\Sigma$ valued in $TM$. For concreteness, we require the sections of $\cF$ to be in $W^{2,2}$ and $\cE^{0,1}$ in $W^{1,2}$,
  where the norms are taken with respect to cylindrical metrics near the punctures. We write $\Omega^{0,1}_c(u^* TM)$ for the subspace of $\cE^{0,1}(TM)$ consisting of smooth $1$-forms of compact support.
\end{notn}

Given a connection $\nabla$ on $u^* TM$, which is complex-linear in the sense that $\nabla_v J \xi = J \nabla_v \xi$, we begin by considering the complex-linear Fredholm operator
\begin{align} \label{eq:linear_operator-curve}
  \nabla^{0,1} \co \cF(u^* TM) & \to \cE^{0,1} (u^* TM) \\
  \xi & \mapsto \left(\nabla \xi \right)^{0,1} \equiv \left( v \mapsto \nabla_{jv} \xi - \nabla_v J \xi \right)
\end{align}
associated to each map $u \co \Sigma \to TM$. If the domain is not compact, we assume that $u$ converges along each puncture to a map $x \co S^1 \to M$, that we have a chosen unitary trivialisation of $x^* TM$, and that the chosen connection on $u^* TM$ extends the connection induced by this trivalisation. Since we are only interested in asymptotic conditions given by time-$1$ Hamiltonian orbits, we shall fix:
\begin{equation}
  \label{eq:fix_complex_linear_connection}
  \parbox{30em}{a complex linear connection on the pullback of $TM$ to $S^1 \times TM$, whose restriction to the graph of every Hamiltonian orbit is induced by a trivialisation,}
\end{equation}
and assume that $\nabla$ converges to this connection along cylindrical ends associated to the punctures of $\Sigma$.

Let us now consider a (not-necessarily complex) linear map
\begin{equation}
  Y \co  \cF(u^* TM) \to \cE^{0,1} (u^* TM),
\end{equation}
of order $0$ (in the sense that $Y(\xi)$, evaluated at a point $z \in \Sigma$ depends only on the value $\xi(z)$ rather than on the derivatives at this point), and consider the operator 
\begin{align} \label{eq:linear_operator-curve2}
  \nabla^{0,1}_{Y} \co \cF(u^* TM) & \to \cE^{0,1} (u^* TM) \\
  \xi & \mapsto \left(\nabla \xi   \right)^{0,1} - Y(\xi).
\end{align}
In order to control the behaviour of this operator, we shall assume that 
\begin{equation}
  \label{eq:fix_asymptotic_conditions}
  \parbox{30em}{near each puncture, $Y = (A \otimes dt)^{0,1}$ for an endomorphism $A$ of $u^* TM$, and the connection $\nabla_A$ induced by $A$ and the chosen connection preserves the symplectic form. Moreover, either $A$ vanishes or $\nabla_A$  exponentially converges in all derivatives to a flat connection whose monodromy is a symplectic matrix that does not admit $1$ as an eigenvalue.}
\end{equation}
Explicitly, the above convergence conditions say that there is a loop  $S(t)$ of symmetric matrices such that $A$ converges to $S(t)$ with respect to a choice of cylindrical end. A loop of such matrices generates a path of symplectomorphisms, whose endpoint has constrained eigenvalues.

We shall make use of the following standard fact, which amounts to the statement that the operator considered above is elliptic:
\begin{lem} \label{lem:complement_support_arbitrary}
  Let $G$ be a finite group of automorphisms of $u$ preserving $u^*Y$. There is a finite dimensional complex $G$-representation $V$, equipped with an equivariant complex-linear map
\begin{equation} 
   \lambda \co V \to  \Omega^{0,1}_c(\Sigma, u^*TM),
 \end{equation}
 whose image may be assumed to be supported in any arbitrary $G$-invariant non-empty open subset of the domain, such that the direct sum
 \begin{equation}
 \nabla_{Y}^{0,1}\oplus \lambda \co  \cF(u^* TM)  \oplus V\to \cE^{0,1}(u^* TM)
\end{equation}
is surjective.  Moreover, if this operator is surjective at $u$, then any continuous extension to a sufficiently small neighbourhood  is surjective (using the $C^\infty$-topology with exponential decay along the ends).  \qed 
\end{lem}

We shall use this construction to associate to each Hamiltonian orbit $[p]$ (with chosen trivialisation of $[p]^* TM$) a stable vector space as follows: abusing notation, we write $ \cF([p]^* TM)$ and $ \cE^{0,1} ([p]^* TM) $ for the Banach spaces obtained from $[p]$ by considering the induced translation-independent map $\bR \times S^1 \to M$.  Let
\begin{equation} \label{eq:cutoff_equation}
  \chi \co \bP^1 \to [0,1]
\end{equation}
be a cutoff function which vanishes outside the circle of radius $\exp(1)$, and is identically equal to $1$ in the circle of radius $\exp(-1)$, and consider the $1$-form on $\bR \times S^1 = \bP^1 \setminus \{ 0, \infty \}$
\begin{equation}
Y_{[p]}(s,t)  \equiv  \left(\chi(e^{s+it}) \cdot B_H \otimes dt \right)^{0,1} \in  \cE^{0,1} ([p]^* TM),
\end{equation}
where the matrix $B_H$ is obtained by differentiating the Hamiltonian flow with respect to the chosen connection.  By construction, this $1$-form vanishes whenever $t$ is sufficiently large, hence satisfies  Condition \eqref{eq:fix_asymptotic_conditions} along the positive end. Along the negative end, the assumption that each orbit be non-degenerate is equivalent to the requirement that Condition \eqref{eq:fix_asymptotic_conditions} hold, so that we obtain a Fredholm operator $\nabla_{[p]}^{0,1} $ from Equation \eqref{eq:linear_operator-curve2}.

\begin{figure}[h]
  \centering
  \begin{tikzpicture}
    \coordinate[label= right:$\nabla^{0,1} \xi$] (h) at (4,0);
    \coordinate[label= right:$(\nabla \xi - B_H \otimes dt)^{0,1}$] (i) at (-6,0);
    \begin{scope}[shift={(4,0)}]
     \draw [dashed,thin] (-2,-.5) .. controls (-2.25,-.5) and (-2.25,.5) .. (-2,.5);
    \draw [thin] (-2,-.5) .. controls (-1.75,-.5) and (-1.75,.5) .. (-2,.5); 
    \end{scope}
    \draw [dashed,thin] (-2,-.5) .. controls (-2.25,-.5) and (-2.25,.5) .. (-2,.5);
    \draw [thin] (-2,-.5) .. controls (-1.75,-.5) and (-1.75,.5) .. (-2,.5);
    \draw (-6,.5) -- (6,.5);
    \draw (-6,-.5) -- (6,-.5);
    
  \end{tikzpicture}
  \caption{The operator $\nabla_{[p]}^{0,1} $ associated to each Hamiltonian orbit.}
  \label{fig:orbit_operator}
\end{figure}
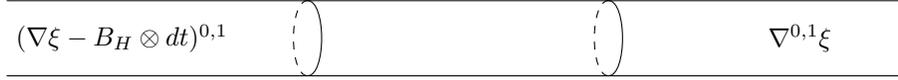

Applying Lemma \ref{lem:complement_support_arbitrary}, we find a finite dimensional complex vector space $V^-_{[p]}$ equipped with a map
\begin{equation} \label{eq:obstruction_map_orbit}
  \lambda_{[p]} \co V^-_{[p]} \to \Omega^{0,1}_c([p]^*  TM),
\end{equation}
so that the direct sum $\nabla_{[p]}^{0,1}  \oplus \lambda_x$ is surjective. For concreteness, we also assume that the support, expressed in cylindrical coordinates,  is contained in $[-1,1] \times S^1$. Letting $V^+_{[p]}$ denote the kernel of this map, we obtain a stable vector space $(V^+_{[p]}, V^-_{[p]})$ associated to each orbit.
\begin{rem}
In the abstract theory discussed in Sections \ref{sec:compl-lifts-tang} and \ref{sec:compl-orient-kuran-1}, the complex structure on $V^-_{[p]}$ is used in formulating the appropriate notion of multiplicativity of orientations, but we shall need it to before that stage, to define the stable complex structure.
\end{rem}

The above discussion will be used to define stable complex orientations for moduli spaces of pseudo-holomorphic cylinders.
If $[x]$ is a critical point of the Morse function $f$, we have already defined $V^+_{[x]}$ in Section \ref{sec:stable-compl-struct-2} to be the tangent space of the ascending manifold of $f$ at $[x]$; we set $V^-_{[x]} = 0$. To have a consistent notation, given $a \in \cP(f,H,f)$, we write $V^\pm_{a}$ for the vector spaces associated to the critical point or orbit corresponding to $a$, and then define
\begin{equation}
  V^\pm_{\underline{a}}  \equiv
  \begin{cases}
     V_a^\pm  & \textrm{if } a \in \cP(H) \\
    0 & \textrm{otherwise,}
  \end{cases}
\end{equation}
noting that comparing with Equation \eqref{eq:keep_Morse_space_set_Floer_0} gives a canonical isomorphism
\begin{equation}
   V^\pm_{\underline{a}} \oplus  V_{\overline{a}}^\pm \cong V_a^\pm.
\end{equation}

\subsubsection{Moduli spaces of maps and Fredholm operators}
\label{sec:based-floer-traj}

We shall presently construct a $1$-parameter family of stable vector bundles which interpolate between the (fibrewise) tangent space of the moduli space of Floer trajectories and a stable complex vector bundle. The construction will proceed via an interpolating family of Fredholm operators, but the usual difficulties of bubbling and breaking in Floer theory will make the construction clearer if we introduce an auxilliary moduli space. For the next definition, we let $(a,b)$ be elements of $\cP(f,H,f)$, and recall that the Floer components are the irreducible components of the domain of a pseudo-holomorphic curve on which the Cauchy-Riemann operator that we consider has a non-vanishing inhomogeneous term. In the case of Floer trajectories, this consists of all components lying between the input and the output: 

\begin{defin}
  The moduli space $\tilde{\cM}^{\bR}(a,b)$ of \emph{based Floer trajectories} consists of an element of $\fcMR(a,b)$, together with a lift of the underlying map $u \co \Sigma \to M$ to the space of equivalence classes of stable maps
  \begin{equation}
     \tilde{u} \co \tilde{\Sigma} \to M \times \bP^1, 
   \end{equation}
   whose domain is a pre-stable curve $\tilde{\Sigma}$,  such that  (i) the stabilisation of the projection to $M$ has domain $\Sigma$, and agrees with $u$, and (ii) the projection to $\bP^1$ has degree $1$, maps the marked points $z_\pm$ to $0$ and $\infty$, and respects the angular structure on the non-constant component. 
 \end{defin}
 To clarify the last condition, observe that since holomorphic maps have non-negative degree, and are constant when the degree vanishes, there is a unique component of $\tilde{\Sigma}$ mapping non-trivially to $\bP^1$. The assumptions imply that, if this component is not collapsed by the projection to $M$, then it has an identification with $\bR \times S^1$ that is canonical up to translation; the additional datum in a based trajectory thus fixes the ambiguitity in this identification, and in this case $\tilde{\Sigma} = \Sigma$. The only other possibility is the existence of a single unstable component on which the map is non-constant, which is labelled by an element of $\cP$ appearing as an asymptotic condition of one of the cylindrical components of $\Sigma$. We conclude:
\begin{lem}
The forgetful map $\tilde{\cM}^{\bR}(a,b) \to  \fcMR(a,b)$ is a fibre bundle with fibre homeomorphic to a closed interval. \qed
\end{lem}
We orient this interval according to the position of the distinguished component in the arc connecting $-$ to $+$, i.e. so that there is an orientation preserving identification with $[0,1]$ mapping $0$ to the point in the fibre where this component is leftmost in the chain (if both $a$ and $b$ lie in $\cP(H)$, then this is necessarily a Floer component). If $a = p$ is an orbit of $H$, this corresponds to the component mapping to $M$ via the projection to $S^1$ and the map $p \co S^1 \to M$. Otherwise, this corresponds to the case where the map from this component to $M$ is constant.

We now place ourselves in the context of Section \ref{sec:regul-thick-moduli},  i.e. we consider the zero locus $Z_\alpha$ of the moduli space $\fcMR_\alpha(a,b)$, whose points are represented by a Riemann surface $\Sigma$ with marked points, together with a map $u \co \Sigma \to M$. The moduli space is stratified by trees, and if $\Sigma$ lies in the stratum indexed by a tree $T$, we write $\{\Sigma_v\}_{v \in V(T)}$ for the underlying curves labelled by the vertices of a tree $T$.

Pulling back $\tilde{\cM}^{\bR}(a,b) $ under the map $Z_\alpha \to \fcMR(a,b)$, we obtain a space which we denote $\tilde{Z}_\alpha $, which is stratified by trees $\tilde{T}$, with vertices labelled $v_+$ and $v_-$ as before, and a distinguished vertex $v_\bullet$ along the path from $v_+$ to $v_-$, corresponding to the component carrying a non-trivial map to $\bP^1$. We may canonically assign to $\tilde{T}$ a tree $T$ which labels a stratum of $ \fcMR_\alpha(a,b) $, by forgeting the vertex $v_\bullet$ (and collapsing it if the resulting map is unstable). Note that, if this vertex is a stable component, then $T$ is isomorphic to $\tilde{T}$, and otherwise is has one fewer vertex.

There are two distinguished sections $\{\tilde{u}_\bullet = \underline{a}\}$ and $\{\tilde{u}_\bullet = \underline{b} \}$ of the map from unbased to based Floer trajectories, given by adding an unstable component mapping either to $\underline{a}$ or to $\underline{b}$, and corresponding to the stratum for which $v_\bullet$ labels an unstable component which is either leftmost or rightmost along the arc of Floer components. Our goal is to construct a stable vector bundle on $\tilde{Z}_{\alpha} $ whose restriction to the section labelled by $a$ is the direct sum of the stable tangent space of $\fcMR_{\alpha}(a,b)$ with $(V^+_{\underline{a}},V^-_{\underline{a}})$, and whose restriction to the other section is the direct sum of a stable complex vector bundle with $(V^+_{\underline{b}}, V^-_{\underline{b}})$. The contractibility of the fibres will then induce an isomorphism between these vector bundles on $Z_\alpha$ which is canonical up to contractible choice. As discussed earlier, the construction will proceed via families of Fredholm operators whose domains and range we begin by introducing.

For each element $i \in S_\alpha$ of the set of Floer thickening data, recall that we fixed codimension $2$ submanifolds $D_i$ of $M$, and a path of subbundles of the restriction of $T M$ to $D_i$. Given the identification of the fibre of $\tilde{\cM}^{\bR}(a,b) \to \fcMR(a,b)$ with an interval, we thus obtain an assignment $\Delta_i(\tilde{u})$ of a subbundle of the pullback of $TM$ along this path for each marked point labelled by $i$ of the domain of a curve $\tilde{u} \in \tilde{Z}_{\alpha}$, which agrees with the tangent space of $D_i$ along the section labelled by $a$, and with the complex subbundle along the section labelled by $b$.

We obtain a Banach space 
\begin{equation}
  \cF\left(\tilde{u}^*_v TM, \{ \Delta_i \}_{i \in S_\alpha} \right)   
\end{equation}
associated to each vertex $v \in {T}$, consisting of sections whose values at a point marked by $r_i$ lies in $\Delta_i(\tilde{u})$ for each $i \in S_\alpha$. We specify that we require exponential decay at the ends $z_\pm$, even if these are asymptotic to constants (rather than Floer trajectories), but that we do not impose any condition at nodes.  If $v_\bullet$ labels an unstable component corresponding to an orbit, then this is the Banach space associated in Section \ref{sec:constr-stable-compl-1} to an orbit. By taking the direct sum over all vertices of ${T}$, we define
\begin{equation}
  \cF(\tilde{u}^* TM,  \Delta_\alpha ) \equiv \bigoplus_{v}  \cF(\tilde{u}^*_v TM, \{ \Delta_i \}_{i \in S_\alpha}).
\end{equation}

This construction is functorial in $\alpha$, in the sense that, for each map $\alpha \to \beta$,  the forgetful map  induces an inclusion
\begin{equation} \label{eq:inclusion_adding_marked_points_space_of_maps}
    \cF(\tilde{u}^* TM,  \Delta_\alpha )  \subset     \cF(\tilde{u}^* TM, \Delta_\beta ).
  \end{equation}
  of Banach bundles over $\tilde{Z}_{\alpha} $.
\begin{lem} \label{lem:quotient_inclusion_tangent_spaces_marked}
If  $\Delta_i(\tilde{u})$ is transverse to the image of the tangent space of the domain under the differential of $u$, the cokernel of the inclusion in Equation \eqref{eq:inclusion_adding_marked_points_space_of_maps} is naturally isomorphic, as a vector bundle over $\tilde{Z}_{\alpha} $, to the pullback of the fibrewise tangent space $ T^{\beta} \fcMR_{\alpha} $ of the projection from $\fcMR_{\alpha}$ to $\fcMR_{\beta}$.
\end{lem}
\begin{proof}
  The fibrewise tangent space $T^{\beta} \fcMR_\alpha$ is naturally identified with a direct sum, over all marked points forgotten by the map $\fcMR_\alpha \to \fcMR_{\beta}$, of the tangent space of the underlying curve. Since the domain in Equation \eqref{eq:inclusion_adding_marked_points_space_of_maps} is obtained from the range by imposing divisorial conditions at these marked points, the transversality condition implies that the cokernel is isomorphic to this direct sum as well. 
\end{proof}

We also define
\begin{equation} 
  \cE^{0,1}(\tilde{u}^* TM) \equiv \bigoplus_{v}  \cE^{0,1}(\tilde{u}^*_v TM) \oplus \bigoplus_{e}   T_{\tilde{u}(e)} M 
\end{equation}
where the second direct sum is taken over all edges of $T$ which are not Floer edges (i.e. the nodes of the Riemann surface) and, as before, $\tilde{u}(e)$ is the common value of the two sides of the node associated to the edge $e$.

Recall that we fixed a cutoff function $\chi \co \bR \to [0,1]$ in Equation \eqref{eq:cutoff_equation}. This function induces a map $\chi \circ \tilde{u}_v$ on each Floer component of the domain of $\tilde{u}$ by composing $\tilde{u}$ with the projection to $\bP^1$, and pulling back under the projection $\bR \times S^1 \to \bR$, and the map $\tilde{u} \to \bP^1$. By construction, $\chi \circ \tilde{u}_v$ is a non-trivial cutoff function if $v = v_\bullet$, vanishes on all component between $v_+$ and $v_\bullet$, and is identically $1$ on the components between $v_\bullet$ and $v_-$.

We can write the linearisation of the pseudo-holomorphic curve equation along each component $u_v$ of $u$ in the form of Equation~\eqref{eq:linear_operator-curve2}. On a lift $\tilde{u}_v$,  we consider the  $1$-form
\begin{equation}
  Y_{\tilde{u}_v} \equiv \left( \chi \circ\tilde{u}_v \right)  Y_v \in C^{\infty}( \Sigma_v \times M, T^{0,1} \Sigma_v  \otimes TM).
\end{equation}

We obtain an operator
\begin{equation} \label{eq:dbar_operator_stable_map}
 \nabla_{\tilde{u}}^{0,1} \co    \cF(\tilde{u}^* TM,  \Delta_\alpha)  \to    \cE^{0,1}(\tilde{u}^* TM) 
\end{equation}
given on each component by the restriction of
\begin{align} \label{eq:linear_operator-lift}
  \xi & \mapsto \left(\nabla \xi   \right)^{0,1} - Y_{\tilde{u}_v}.
\end{align}
to the space of sections with values  constrained to lie in $\Delta_i(\tilde{u})$ at the points marked by $r_i$, by the evaluation map
\begin{equation}
 \cF(\tilde{u}^*_v TM, \Delta_\alpha)  \to  T_{\tilde{u}(e)} M
\end{equation}
at each flag $(v,e)$ pointing towards the root (with $e$ corresponding to a node), and its negative for flags pointing away from the root (recall that a flag is a pair consisting of an edge and one of its endpoints, so that it can be thought of as a direction along the edge). Since the operator $ \nabla^{0,1}_{Y_{\tilde{u}_v}}$ 
 satisfies Condition \eqref{eq:fix_asymptotic_conditions} on the cylindrical components, the operator $ \nabla_{\tilde{u}}^{0,1}$ is Fredholm.

 \begin{lem}
   The pullback of the operator \eqref{eq:dbar_operator_stable_map} to the section $\{\tilde{u}_\bullet = \underline{b}\}$ is canonically isomorphic to the direct sum of the linearisation of the $\dbar$ operator on $u$ with the operator $\nabla_{b}^{0,1}$ (we interpret this operator to vanish if $b$ is not a Hamiltonian orbit). The pullback of the operator \eqref{eq:dbar_operator_stable_map} to the section $\{\tilde{u}_\bullet = \underline{a}\}$ is the direct sum of a complex linear operator with the operator $\nabla_{a}^{0,1}$.  \qed
 \end{lem}

Note that $\tilde{Z}_\alpha$ is naturally equipped with a $G_\alpha$ action with the property that the projection to the Floer moduli space is $G_\alpha$ equivariant.  The operators $ \nabla_{\tilde{u}}^{0,1}$ are compatible  with this action in the sense that each element $g$ induces a commutative diagram
\begin{equation}
  \begin{tikzcd}
    \cF(\tilde{u}^* TM, \Delta_\alpha)  \ar[r] \ar[d] &    \cE^{0,1}(\tilde{u}^* TM) \ar[d] \\
     \cF((g \cdot \tilde{u}) ^* TM,  \Delta_\alpha)  \ar[r] &    \cE^{0,1}((g \cdot \tilde{u})^* TM),
  \end{tikzcd}
\end{equation}
and the vertical maps associated to a product $h \cdot g$ agree with the composition.

\subsubsection{The fibrewise oriented locus of the space of stable maps}
\label{sec:oriented-locus-space}

The purpose of this section is to incorporate the obstruction spaces $V_\alpha$ that enter in the construction of the moduli space $\Mbar^{\bR,\reg}_\alpha(a,b)$, and those which we associate to each orbit, into a single obstruction space mapping to the restriction of $\cE^{0,1}(\tilde{u}^* TM)$  to a subset of the space $\tilde{Z}_\alpha$  constructed in the previous section. The essential missing ingredient in the construction of this larger obstruction space is a choice of cylindrical ends for Floer trajectories, i.e. an identification of a neighbourhood of each Floer end with either $[0,+\infty) \times S^1$ or $(-\infty, 0] \times S^1$. Since the component carrying each end is canonically identified with $\bR \times S^1$ up to translation, this choice can locally be identified with a choice of a positive real number.

The key fact we need about cylindrical ends is that,  if curves $\Sigma_1$ and $\Sigma_2$ are equipped with cylindrical ends along the positive and negative end, then the result
\begin{equation}
  \Sigma_1 \#_R \Sigma_2  
\end{equation}
of gluing the positive end of $\Sigma_1$ to the negative end of $\Sigma_2$ with a finite gluing parameter $R \in [0,\infty)$ is again canonically equipped with cylindrical ends. We say that such a cylindrical end is \emph{obtained by gluing.} If $\Sigma_i$ are equipped with marked points, so that the pair $(\Sigma_1, \Sigma_2)$ define an element of a codimension $1$ boundary stratum of some moduli space $\fcMR_\alpha$, then the construction of the glued surface above gives a neighbourhood of this boundary stratum in $\fcMR_\alpha$, when $\Sigma_i$ are allowed to vary in modulus.

To systematically use the above observation, recall that a totally ordered subset $P_\alpha \subset \underline{\cP}(\underline{a},\underline{b})$ is part of the datum of an element of $A(a,b)$. We are particularly interested in the case $P_\alpha$ is a subset of $\cP(H)$, so we restrict our attention to that situation. 

\begin{defin}
  The space $\Ends(\alpha)$ of \emph{parametrised ends} on $\fcMR_\alpha$, is the space of smoothly varying choices $\epsilon$ of cylindrical ends for each Floer flag $(v,e)$ with $e$ labelled by an element of $P_\alpha \setminus Q_\alpha$, on the stable curves parametrised by the stratum $\partial^{P_\alpha} \fcMR_{\alpha}$ which lie in the canonical $\bR$-family of ends induced by the cylindrical structure, are
  \begin{equation}
    \label{eq:ends_disjoint_special}
    \parbox{30em}{separated from all marked points and special points, by an annulus of modulus at least $1$,}
  \end{equation}
 and such that
  \begin{equation}
    \label{eq:condition-U-obtained-by-gluing} \parbox{30em}{every point in $\fcMR_{\alpha}$ is represented as a glued surface in the image of the embedding of $\partial^{P_\alpha} \fcMR_{\alpha} \times [0,\infty)^{P_\alpha \setminus \underline{Q}_\alpha}$ in $\partial^{\underline{Q}_\alpha}\fcMR_{\underline{\alpha}}$.}
  \end{equation}
  \end{defin}

  As discussed above, every surface in $  \fcMR_{\alpha}$ inherits a choice of parametrised ends induced by gluing.  Because of Condition \eqref{eq:condition-U-obtained-by-gluing}, the space $\Ends(\alpha)$ may be empty. However, the following observation will greatly simplify the discussion in Section \ref{sec:categ-orient-charts} below:
  \begin{lem} \label{lem:ends_alpha_contractible}
    If $\Ends(\alpha)$ is non-empty, then it is contractible, and the projection map
    \begin{equation}
      \fcMR_{\alpha} \to \partial^{P_\alpha} \fcMR_{\alpha} 
    \end{equation}
    induced by forgetting the gluing parameters is independent of this choice of element of $\Ends(\alpha)$. \qed
  \end{lem}

In Section \ref{sec:constr-stable-compl-1}, we fixed a complex vector space $V^-_{[p]}$ for each orbit $[p]$,  and a map $\lambda_{[p]}$ in Equation \eqref{eq:obstruction_map_orbit}. Letting $p$ denote a lift to $\cP(H)$, we write $V^-_p$ for $ V^-_{[p]}$. Given an element $u \in Z_\alpha$ we obtain a map
\begin{equation} \label{eq:obstruction_map_orbit_unstable_component}
  \lambda_p \co V^-_p \to  \cE^{0,1}(\tilde{u}^* TM)
\end{equation}
whenever $\tilde{u}_{v_\bullet}$ corresponds to the orbit $p$, given by the composition
\begin{equation}
 V^-_p \to \Omega^{0,1}_c(p^*  TM) \to  \cE^{0,1}(\tilde{u}^*_{v_\bullet} TM) \to \cE^{0,1}(\tilde{u}^* TM).
\end{equation}
Given an element $\epsilon \in \Ends({\alpha})$, we shall extend this map to an open subset of $\tilde{Z}_\alpha$. To begin, note that the choice of ends determines, for each curve $\Sigma$ represented by a point in $\fcMR_{\alpha}$ and each orbit $p \in P_\alpha$,  a \emph{thin part labelled by $p$} which we denote $\Sigma_p \subset \Sigma$ which is either (i) a half-infinite cylinder, (ii) the union of two half-infinite cylinders, or (iii) an annulus $[0,R] \times S^1$ which is the region along which the two surfaces are glued. By definition, $\Sigma_p$ admits a canonical projection to $S^1$. For the next definition, we fix, for each orbit $[p]$, a smooth map from $S^1 \times D^{2n}$ to $M$, extending the map $[p]$ for $S^1 \times \{0\}$, and with the property that the image of $\{t\} \times D^{2n}$ is a geometrically convex neighbourhood of $[p](t)$ for each $t \in S^1$.

\begin{defin} \label{def:space-cylinders-oriented-not-regular}
  The space of Floer trajectories $Z^\epsilon_{\alpha}$, is the space of of maps $u \in Z_\alpha$ such that each point in $\Sigma_p$ lying over a point $t \in S^1$ maps to the fixed geodesically convex neighbourhood of $[p](t)$.
\end{defin}

We now define a map
\begin{equation}
 V^-_p \to  \cE^{0,1}(\tilde{u}^* TM)
\end{equation}
for each lift $\tilde{u} \in  \tilde{Z}^\epsilon_\alpha $  of such a map as follows: the inverse image of $[-1,1] \times S^1  \subset \bP^1$ under the map $\tilde{\Sigma} \to \bP^1$ which determines the lift yields a map $[-1,1] \times S^1 \to  \Sigma$ which is an embedding unless the lift has an additional component, corresponding to the inverse image of the corresponding subset of $\bP^1$. If the image of $[-1,1] \times S^1$ is not contained in $\Sigma_p$, we set the extension of $\lambda_p$ to vanish, and we define it by parallel transport from the image of $p$ to the thin region whenever the distance from $[-1,1] \times S^1$ to every point in  $\partial \Sigma_p$ is greater than $1$. We interpolate between these choices using a fixed cutoff function depending on the distance to $\partial \Sigma_p$.

The construction of the thickened moduli space $\fcMR_{\alpha}(a,b)$ also provides us with a map $\lambda_{\alpha}$ from $V_{\alpha}$ to $\cE^{0,1}(\tilde{u}^* TM)$. For the next definition, we set
\begin{align} \label{eq:W_S}
  V'_{\alpha} & \equiv   V_{\alpha} \oplus \bigoplus_{p \in P_\alpha} V^-_p
\end{align}
\begin{defin} \label{def:regular_locus_based_curve}
  The \emph{regular locus} $\tilde{Z}^{\epsilon,\reg}_{\alpha}$ is the set of based Floer trajectories $\tilde{u}$ in $\tilde{Z}^\epsilon_\alpha  $ such that (i) the image under $u$ of the tangent space of the domain at a point labelled by $r_i$ is transverse to $\Delta(\tilde{u})$ for each $i \in S$, and (ii) the map
  \begin{equation} \label{eq:extended-operator-based-curve}
 \cF(\tilde{u}^* TM, \{ \Delta_i \}_{i \in S_\alpha}) \oplus V'_{\alpha} \to    \cE^{0,1}(\tilde{u}^* TM)  
\end{equation}
is surjective.
\end{defin}

 For the section associated to $b$, Equation \eqref{eq:extended-operator-based-curve} is, by construction,  the direct sum of three factors: (i) the operator $ \nabla^{0,1}_b \oplus \lambda_b $ constructed in the previous section, (ii)  the trivial map on  the direct sum of $V^-_{\underline{a}}$ with the vector space
\begin{equation}
  W^{\fib}_{\alpha} \equiv  \bigoplus_{p \in P_\alpha} V^-_p,
\end{equation}
and (iii) the inhomogeneous operator 
\begin{equation} \label{eq:linearisation-operator-curve}
 \cF(u^* TM,  \Delta_\alpha) \oplus V_{\alpha} \to    \cE^{0,1}(u^* TM) 
  \end{equation}
obtained by linearising the Floer equation. 

For the section associated to $a$, we also have a direct sum of three factors: (i) the operator $ \nabla^{0,1}_a \oplus \lambda_a $, (ii) the trivial map on the vector space $ V^-_{\underline{b}} \oplus W^{\fib}_{\alpha}$, and (iii) a complex-linear map
\begin{equation} \label{eq:complex-linear-operator-curve}
 \cF(u^* TM, \Delta_\alpha) \oplus V_{\alpha} \to    \cE^{0,1}(u^* TM) 
  \end{equation}
whose restriction to the first factor is homogeneous. 

We define the  \emph{index bundle} $ I^{\fib,\bC}_{\alpha}$  to be the complex vector bundle defined as the kernel of Equation \eqref{eq:complex-linear-operator-curve} over the locus where this complex linear operator is surjective.  The above discussion implies: 
\begin{prop} \label{prop:equivariant_stabilisation_vector_bundle}
There is a  $G_\alpha$-equivariant vector bundle over $ \tilde{Z}^{\epsilon, \reg}_{\alpha}$ such that (i) the pullback under section associated to $b$ is naturally isomorphic to the direct sum
  \begin{equation}  \label{eq:t=1-vector-bundle-tangent-space}
  V^+_{\underline{b}} \oplus W^{\fib}_\alpha \oplus  T^{\alpha} \fcMR_\alpha(a,b) \oplus  V^-_{\underline{a}} ,     
  \end{equation}
  while (ii) the pullback under the section associated to $a$ is naturally isomorphic to a direct sum
  \begin{equation} \label{eq:t=0-vector-bundle-complex}
  V^-_{\underline{b}} \oplus W^{\fib}_\alpha \oplus I^{\fib,\bC}_{\alpha} \oplus V^+_{\underline{a}} .  
  \end{equation}
  \qed
\end{prop}

This result leads us to the following:
\begin{defin}
The \emph{stably fibrewise oriented locus} 
\begin{equation}
  Z_{\alpha}^{\epsilon,\ori} \subset   Z_{\alpha}^{\epsilon, \reg}
\end{equation}
is the set of points whose inverse image $\tilde{Z}_{\alpha}^{\epsilon,\ori} $ in $\tilde{Z}^{\epsilon}_{\alpha} $ lies in the regular locus $\tilde{Z}^{\epsilon,\reg}_{\alpha}$.
\end{defin}

As an immediate consequence of Proposition \ref{prop:equivariant_stabilisation_vector_bundle}, we conclude that the restriction of $ T^{\alpha} \fcMR_\alpha(a,b)$ to the stably oriented zero-locus is equipped with a stable isomorphism to the complex vector bundle $ I^{\fib,\bC}_{\alpha}$, which is defined up to contractible choice.

To be more precise, it is convenient to introduce a parametrisation of the fibres of the maps $\tilde{\cM}^{\bR}(a,b) \to \Mbar^{\bR}(a,b) $, which is compatible with the structure maps. Note that the map $\fcMR(a,b) \times \fcMR(b,c) \to \fcMR(a,c)$ for a triple $(a,b,c)$ lifts to an identification of each fibre of $ \tilde{\cM}^{\bR}(a,c)$ over this boundary stratum with the union of the fibres of $ \tilde{\cM}^{\bR}(a,b)$ and  $\tilde{\cM}^{\bR}(b,c)$, glued at the common endpoint. This leads to the following definition:

\begin{defin} \label{def:Moore-parametrisation-multiplicative}
A \emph{multiplicative Moore parametrisation} of the fibres of
$\tilde{\cM}^{\bR}(a,b) \to  \fcMR(a,b)$ is specified by: 
\begin{enumerate}
\item A function $m_{ab} \co \fcMR(a,b) \to (0,\infty)$ for each pair
  $(a,b) \in \cP(f,H,f)$, and 
\item an identification of each fibre with the interval $[0,m_{ab}]$
  such that, for each triple $(a,b,c)$, 
\begin{enumerate}
\item the map $m_{ac}$ restricts
  on the stratum $\fcMR(a,b) \times \fcMR(b,c)$  to $m_{ab} + m_{bc}$,
  and 
\item the identification of fibres is compatible with the natural
  map 
\begin{equation} \label{eq:concatenation_is_addition}
[0,m_{ab}] \cup [0,m_{bc}] \to [0,m_{ab}] \cup [m_{ab}, m_{bc}
    + m_{ab}] \to [0, m_{ab} + m_{bc}].
\end{equation}
\end{enumerate}
\end{enumerate} 
\end{defin}

From the above construction, we obtain a family of
$G_\alpha$-equivariant vector bundles on $Z_{\alpha}^{\epsilon, \ori}$, 
\begin{equation}
\{I^{t,\epsilon}_{\alpha}, t \in [0,m_{ab}]\},
\end{equation}
from the identification of the fibre of $\tilde{\cM}^{\bR}(a,b) \to \fcMR(a,b)$ with intervals $[0,m_{ab}]$, and such that $I^{0,\epsilon}_{\alpha} $ is given by Equation \eqref{eq:t=0-vector-bundle-complex}, and $I^{m_{ab},\epsilon}_{\alpha} $ by Equation \eqref{eq:t=1-vector-bundle-tangent-space}.

\begin{defin}
  For each choice $\epsilon \in \Ends({\alpha})$ of strip-like ends,
  we define $\scrJ_{\fib}(\epsilon) $ to be the space consisting of
  continuous families of 
\begin{enumerate}
\item inner products on the vector bundles
  $I^{t,\epsilon}_{\alpha}$, and 
\item continuous families of
  inner-product preserving isomorphisms 
  \begin{equation} \label{eq:trivialisation_family_of_bundles}
    I^{t,\epsilon}_{\alpha} \cong I^{0,\epsilon}_{\alpha },
  \end{equation}
\end{enumerate}
parametrised by $t \in [0,m_{ab}]$ which are the identity for $t=0$.
\end{defin}

\begin{rem}
  It is important to note that the space $\scrJ_{\fib}(\epsilon) $ does not depend on the choice of parametrisation in the sense that there is a canonical isomorphism associated to a changing the parametrisation. In that sense, the only use of Definition \ref{def:Moore-parametrisation-multiplicative} is that it will allows us to use the notation $t_1 + t_2$ for concatenation of the points in the fibres, via Equation \eqref{eq:concatenation_is_addition}.
\end{rem}

\begin{cor}
  Each element of $\scrJ_{\fib}(\epsilon)$ determines an isomorphism
  \begin{equation} \label{eq:stable_fibrewise-complex-structure}
    V^+_{\underline{b}} \oplus W^{\fib}_\alpha \oplus  T^{\alpha} \fcMR_\alpha(a,b)  \oplus V^-_{\underline{a}}  \cong V^-_{\underline{b}} \oplus W^{\fib}_\alpha  \oplus I^{\fib,\bC}_{\alpha}  \oplus  V^+_{\underline{a}} ,
  \end{equation}
 of vector bundles over $  \tilde{Z}_{\alpha}^{\epsilon,\ori}$. \qed
\end{cor}

\subsubsection{Functoriality of stable complex fibrewise orientations}
\label{sec:funct-stable-compl}

In order to formulate the functoriality of the construction of the previous section, observe that a map $f \co \alpha \to \beta$ associates to each choice of ends $\epsilon$ on $\fcMR_{\beta}$, a choice of ends on $\fcMR_{\alpha}$. However, this choice does not necessarily satisfy Conditions \eqref{eq:ends_disjoint_special} and  Condition \eqref{eq:condition-U-obtained-by-gluing} because the additional marked points of curves in $\fcMR_\alpha $ may lie in the cylindrical ends specified over $\fcMR_{{\beta}} $.
\begin{defin}
  The space $\Ends({f})$ is the subset of element $\epsilon \in \Ends({{\beta}}) $ for which the pullback $f^* \epsilon$ satisfies Condition \eqref{eq:condition-U-obtained-by-gluing} when restricted to $\fcMR_{\alpha} $. 
\end{defin}

For the next result, we observe that an arrow $f \co \alpha \to \beta$ induces a natural direct sum decomposition
\begin{equation}
  V'_\beta \cong   V'_\alpha \oplus V'_\beta/V'_\alpha,
\end{equation}
where the quotient is the direct sum of $V_\beta/V_\alpha$ with
\begin{equation}
  W^{\fib}_\beta/W^{\fib}_\alpha =  \bigoplus_{r \in P_\beta \setminus P_\alpha} V^+_r.
\end{equation}

\begin{lem}
  The restriction of the map $ \tilde{Z}_{\alpha} \to \tilde{Z}_{\beta}$ maps the oriented locus $ \tilde{Z}^{\epsilon,\ori}_{\alpha} $ to  $ \tilde{Z}^{\epsilon,\ori}_{\beta} $, and induces a short-exact sequence
  \begin{equation} \label{eq:short-exact-sequence-family-vector-bundles-change-chart}
0 \to     I^{t,\epsilon}_{\alpha} \to   I^{t,\epsilon}_{\beta} \to  V'_\beta/V'_\alpha \oplus T^\beta \fcMR_\alpha  \to 0
\end{equation}
for each $t \in [0,m_{ab}] $.  Given a composition $\alpha \to \beta \to \gamma $, the induced exact sequences fit in a commutative diagram:
  \begin{equation}
    \begin{tikzcd}
      I^{t,\epsilon}_{\alpha} \ar[r] \ar[d,"="] & I^{t,\epsilon}_{\beta} \ar[r] \ar[d] & V'_\beta/V'_\alpha \oplus T^\beta \fcMR_\alpha   \ar[d]  \\
      I^{t,\epsilon}_{\alpha} \ar[r] &  I^{t,\epsilon}_{{\gamma}}  \ar[r] \ar[d] &  V'_\gamma/V'_\alpha \oplus T^\gamma \fcMR_\alpha \ar[d] \ar[d] \\
      & V'_\gamma/V'_\beta \oplus T^\gamma \fcMR_\beta  \ar[r,"="]  &  V'_\gamma/V'_\beta \oplus T^\gamma \fcMR_\beta.
    \end{tikzcd}
  \end{equation}

\end{lem}
\begin{proof}
Since $P_\alpha \subset P_\beta$, Definition \ref{def:regular_locus_based_curve} implies that $ Z^\epsilon_{\alpha} $ is mapped to  $ Z^\epsilon_{\beta} $.

  Next, we note that the operators in Equation \eqref{eq:extended-operator-based-curve} fit in a commutative diagram
  \begin{equation}
    \begin{tikzcd}
      \cF(\tilde{u}^* TM,  \Delta_\alpha) \oplus V'_{\alpha} \ar[d] \ar[r] &   \cE^{0,1}(\tilde{u}^* TM) \\
\cF(\tilde{u}^* TM, \Delta_\beta) \oplus V'_{\beta} \ar[ur] &  
    \end{tikzcd}
  \end{equation}
  where the vertical arrow is an inclusion because $S_\beta \subset S_\alpha$ (i.e. there are more constraints in the source), and $P_\alpha \subset P_\beta$ (there are more obstruction bundles associated to orbits in the target). The cokernel of the vertical map is $V'_\beta/V'_\alpha$, from which the desired result follows.
\end{proof}

We now consider the restriction of the above exact sequence to the endpoints of the interval $[0,m_{ab}]$. For $t=m_{ab}$, it is identified  under the isomorphisms in Equation \eqref{eq:t=1-vector-bundle-tangent-space} with the direct sum of
\begin{equation}
  T^\alpha \fcMR_\alpha(a,b)  \to T^{\beta}\fcMR_\beta(a,b)  \to V_\beta/V_\alpha \oplus T^\beta \fcMR_\alpha   
\end{equation}
with the identity on $V^-_{\underline{a}}$, and the direct sum decomposition $W^{\fib}_\beta \cong  W^{\fib}_\alpha \oplus W^{\fib}_\beta  / W^{\fib}_\alpha$. 
On the other hand, at $t=0$, this short exact sequence is identified by Equation \eqref{eq:t=0-vector-bundle-complex} with the direct sum of a short exact sequence of complex bundles
\begin{equation}
  I^{\fib,\bC}_\alpha  \to I^{\fib,\bC}_\beta \to V_\beta/V_\alpha \oplus T^\beta \fcMR_\alpha   
\end{equation}
with the identity on $V^-_{\underline{b}}$. 
In order to use this observation, recall that we have a fixed inner product on $V'_\beta$ and $V'_\alpha$ and hence on the quotient. We consider the space
\begin{equation}
 \scrJ({f},  \epsilon)  \subset \scrJ_{\fib}(\epsilon) \times \scrJ_{\base}({\alpha})
\end{equation}
consisting of an element of $\scrJ_{\fib}(\epsilon) $  such that the metric on $V'_\beta/V'_\alpha \oplus T^\beta \fcMR_\alpha$ induced by Equation \eqref{eq:short-exact-sequence-family-vector-bundles-change-chart} splits as an orthogonal direct sum so that $ V'_\beta/V'_\alpha$ is equipped with its natural inner product, and an element of $\scrJ_{\base}({\alpha}) $ yielding the same inner product on $T^\beta \fcMR_\alpha $ (c.f. Section \ref{sec:splitt-exact-sequ}, where this condition of agreement of inner products appears in an abstract context). In addition, we require  the trivialisation of the family of bundles $I^{t,\epsilon}_{\beta} $, which we fixed in Equation \eqref{eq:stable_fibrewise-complex-structure}, to  preserve the family of subspaces $ I^{t,\epsilon}_{\alpha} $.

The above discussion implies:
\begin{lem} \label{lem:splitting_family_index_bundles_functorial}
There is a natural map $ \scrJ({f},  \epsilon) \to \scrJ_{\fib}(\epsilon)$. Moreover, each element of $\scrJ({f},  \epsilon) $  induces a family of splittings
  \begin{equation} \label{eq:splitting-family-index-bundles-map}
        I^{t,\epsilon}_{\beta} \cong  I^{t,\epsilon}_{\alpha} \oplus  V'_\beta/V'_\alpha \oplus T^\beta \fcMR_\alpha   
  \end{equation}
  parametrised by $t \in [0,m_{ab}]$, yielding a commutative diagram
  \begin{equation} \label{eq:diagram_induced_by_fibrewise_stable_complex_data_for_map}
    \begin{tikzcd}
      \begin{gathered}
          V^+_{\underline{b}} \oplus W^{\fib}_\alpha  \oplus    T^{\alpha} \fcMR_\alpha(a,b) \\ \oplus  V'_\beta/V'_\alpha \oplus T^\beta \fcMR_\alpha \oplus V^-_{\underline{a}}
      \end{gathered}
   \ar[r] \ar[d] &
   \begin{gathered}
     V^-_{\underline{b}}  \oplus W^{\fib}_\alpha  \oplus  I^{\fib,\bC}_{\alpha} \\ \oplus  V'_\beta/V'_\alpha \oplus  T^\beta \fcMR_\alpha  \oplus V^+_{\underline{a}} 
   \end{gathered}
 \ar[d] \\ 
 V^+_{\underline{b}} \oplus W^{\fib}_\beta  \oplus  T^{\beta} \fcMR_\beta(a,b) \oplus V^-_{\underline{a}}  \ar[r] &  V^-_{\underline{b}} \oplus W^{\fib}_\beta \oplus I^{\fib,\bC}_{\beta} \oplus  V^+_{\underline{a}} 
    \end{tikzcd}
  \end{equation}
  of vector bundles over  $ \tilde{Z}^{\epsilon,\ori}_{\alpha} $. \qed
\end{lem}

We now discuss the functoriality of these constructions for composable morphisms  ${f} \co {\alpha} \to {\beta}$ and ${g} \co {\beta} \to {\gamma}$. First, we observe that there is an associative composition map
\begin{equation}
  \Ends({f}) \times_{\End(\beta)} \Ends({g}) \to \Ends(g \circ f),
\end{equation}
given by the fact that pullback of ends under forgetful maps is associative. For each $\epsilon \in \Ends({g}) \subset \Ends({\gamma}) $ such that $g^*(\epsilon) \in \Ends(\beta)$ lies in $\Ends({f}) $, we have natural maps
\begin{equation}
  \scrJ({g},  \epsilon) \to    \scrJ_{\fib}({g}^*(\epsilon)) \times \scrJ_{\base}(\beta)  \leftarrow    \scrJ({f},  {g}^*(\epsilon)).
\end{equation}
\begin{lem} \label{lem:composition-morphisms-large-space-complex-structure}
  There is a natural map
  \begin{equation}
      \scrJ({g},  \epsilon) \times_{   \scrJ_{\fib}({g}^*(\epsilon)) \times \scrJ_{\base}(\beta)} \scrJ({f},  {g}^*(\epsilon)) \to \scrJ({g \circ f},  \epsilon) 
    \end{equation}
  which is associative up to natural isomorphism.\qed
\end{lem}

\subsubsection{Multiplicativity of  stable complex fibrewise orientations}
\label{sec:mult-stable-compl}

We now consider the multiplicativity of the construction from the previous section. Consider a triple of orbits $(a,b,c)$, and charts  $\alpha_1$ and $\alpha_2$ in $A(a,b)$ and $A(b,c)$.

Given choices of strip-like ends $\epsilon_i \in \Ends({{\alpha}_i})$ satisfying Condition \eqref{eq:condition-U-obtained-by-gluing}, the product family $\epsilon_1 \times \epsilon_2$ on $ \fcMR_{{\alpha}_1 \times {\alpha}_1}$ also satisfies this condition, so that we have a natural map
\begin{equation}
  \Ends(\alpha_1) \times \Ends(\alpha_2) \to \Ends(\alpha_1 \times \alpha_2)
\end{equation}

By construction, the product map of Kuranishi charts induces a map of oriented loci,
\begin{equation}
  Z^{\epsilon_1,\ori}_{ \alpha_1} \times Z^{\epsilon_2, \ori}_{\alpha_2}\to Z^{\epsilon_1 \times \epsilon_2,\ori}_{\alpha_1 \times \alpha_2}.
\end{equation}
Moreover, comparing the definition of the index bundles on based Floer trajectories of each factor with those in the total space yields a natural isomorphism
\begin{equation} \label{eq:family-vector-bundles-multiplicative}
      I^{t,\epsilon_1 \times \epsilon_2}_{\alpha_1 \times \alpha_2}  \cong
      \begin{cases}
        I^{t,\epsilon_1}_{\alpha_1} \oplus W^{\fib}_{\alpha_2}  \oplus I^{\fib,\bC}_{\alpha_2}  \oplus V_{\underline{c}}^{-} & 0 \leq t \leq m_{ab} \\
        
        \begin{gathered}
       V_{\underline{a}}^- \oplus T^{\alpha_1} \fcMR_{\alpha_1}(a,b) \oplus W^{\fib}_{\alpha_1} \\ \oplus         I^{t- m_{ab}, \epsilon_2}_{\alpha_2}  
        \end{gathered}   & m_{ab} \leq t \leq m_{ab} + m_{bc},
      \end{cases}
    \end{equation}
where we use the fact that there is a copy of $V^-_{\underline{b}}$ in $I^{t,\epsilon_1}_{\alpha_1}$ and one in $I^{t- m_{ab}, \epsilon_2}_{\alpha_2} $.    For the statement of the next result, we note the isomorphisms
    \begin{align}
      I^{\fib,\bC}_{\alpha_1 \times \alpha_2} & \cong I^{\fib,\bC}_{\alpha_1}  \oplus  I^{\fib,\bC}_{\alpha_2}   \\
W^{\fib}_{\alpha_1 \times \alpha_2} &     \cong W^{\fib}_{\alpha_1}  \oplus 
V^-_{\underline{b}} \oplus  W^{\fib}_{\alpha_2} .
    \end{align}
\begin{lem} 
There is a natural map  $\scrJ_{\fib}(\epsilon_1) \times \scrJ_{\fib}(\epsilon_2) \to \scrJ_{\fib}(\epsilon_1\times \epsilon_2)$, such that the product family of complex structures on the fibrewise tangent space fit in a commutative diagram: 
{\begin{equation} \label{eq:fibre-stable-complex-structure-product-chart}
  \begin{tikzcd}
    \begin{gathered}
   W^{\fib}_{\alpha_1} \oplus  T^{\alpha_1} \fcMR_{\alpha_1}(a,b) \oplus V^-_{\underline{a}}  \oplus \\    V^+_{\underline{c}}   \oplus W^{\fib}_{\alpha_2} \oplus   T^{\alpha_2} \fcMR_{\alpha_2}(b,c)  \oplus  V^-_{\underline{b}}   
    \end{gathered}
    \ar[r] \ar[d]  &  \ar[d]
    \begin{gathered}
      V^+_{\underline{c}} \oplus W^{\fib}_{\alpha_1 \times \alpha_2} \oplus \\  T^{\alpha_1 \times \alpha_2} \fcMR_{\alpha_1 \times \alpha_2}(a,c) \oplus V^-_{\underline{a}} 
    \end{gathered}\\
    \begin{gathered}
      W^{\fib}_{\alpha_1} \oplus  T^{\alpha_1} \fcMR_{\alpha_1}(a,b) \oplus V^-_{\underline{a}}   \oplus \\  V^-_{\underline{c}}   \oplus W^{\fib}_{\alpha_2} \oplus   I^{\fib,\bC}_{\alpha_2} \oplus V^+_{\underline{b}}
    \end{gathered}
      \ar[d] &
      \begin{gathered}
V^-_{\underline{c}}  \oplus W^{\fib}_{\alpha_1 \times \alpha_2} \oplus  \\ I^{\fib,\bC}_{\alpha_1 \times \alpha_2}   \oplus V^+_{\underline{a}}.
      \end{gathered} \\
      \begin{gathered}
        V^-_{\underline{b}}   \oplus W^{\fib}_{\alpha_1} \oplus I^{\fib,\bC}_{\alpha_1} \oplus V^+_{\underline{a}} \oplus \\   V^-_{\underline{c}}  \oplus W^{\fib}_{\alpha_2}  \oplus     I^{\fib,\bC}_{\alpha_2} 
      \end{gathered}
   \ar[ur] &
  \end{tikzcd}
\end{equation}} \qed
\end{lem}

This isomorphism is functorial in the sense that a pair of maps $f_i \co \alpha_i \to \beta_i$, induce a commutative diagram \begin{equation} \label{eq:family-vector-bundles-concatenation-is-functorial}
  \begin{tikzcd}[column sep=tiny]
    \begin{gathered}
       I^{t,\epsilon_1}_{\alpha_1}\oplus W^{\fib}_{\alpha_2}  \\ \oplus I^{\fib,\bC}_{\alpha_2}   \oplus V_{\underline{c}}^{-}
    \end{gathered}
 \ar[r] \ar[d]  &
 \begin{gathered}
I^{t,\epsilon_1}_{\beta_1}  \oplus W^{\fib}_{\beta_2}  \\ \oplus I^{\fib,\bC}_{\beta_2}  \oplus V_{\underline{c}}^{-}    
 \end{gathered}
 \ar[r] \ar[d] &
 \begin{gathered}
   V'_{\beta_1}/ V'_{\alpha_1} \oplus T^{\beta_1}  \fcMR_{\alpha_1 } \\ \oplus  V'_{\beta_2}/V'_{\alpha_2} \oplus T^{ \beta_2} \fcMR_{\alpha_2} 
 \end{gathered}
\ar[d] \\
 I^{t,\epsilon_1 \times \epsilon_2}_{\alpha_1 \times \alpha_2}  \ar[r]  & I^{t,\epsilon_1 \times \epsilon_2}_{\beta_1 \times \beta_2} \ar[r] &
 \begin{gathered}
V'_{\beta_1 \times \beta_2}/V'_{\alpha_1 \times \alpha_2} \oplus \\ T^{\beta_1 \times \beta_2}  \fcMR_{\alpha_1 \times \alpha_2} 
 \end{gathered} \\
 \begin{gathered}
  V_{\underline{a}}^- \oplus  T^{\alpha_1} \fcMR_{\alpha_1}(a,b) \\ \oplus W^{\fib}_{\alpha_1} \oplus         I^{t- m_{ab}, \epsilon_2}_{\alpha_2}   
 \end{gathered}
  \ar[u] \ar[r] &
  \begin{gathered}
    V_{\underline{a}}^- \oplus T^{\beta_1} \fcMR_{\beta_1}(a,b) \\ \oplus W^{\fib}_{\beta_1} \oplus         I^{t- m_{ab}, \epsilon_2}_{\beta_2}
  \end{gathered}
 \ar[u]\ar[r] &
 \begin{gathered}
V'_{\beta_1}/ V'_{\alpha_1} \oplus T^{\beta_1}  \fcMR_{\alpha_1 }\\  \oplus  V'_{\beta_2}/V'_{\alpha_2} \oplus T^{ \beta_2} \fcMR_{\alpha_2}
 \end{gathered} \ar[u]
  \end{tikzcd}
\end{equation} where  $ 0 \leq t \leq m_{ab}$ for the two squares at the top, and  $m_{ab} \leq t \leq m_{ab} + m_{bc} $ for the two squares at the bottom, and all vertical arrows are isomorphisms.   This implies that a pair of splittings for the short exact sequences in Equation \eqref{eq:splitting-family-index-bundles-map} for the arrows $f_i \co \alpha_i \to \beta_i$ determine a splitting for the product arrow $f_1 \times f_2$, characterised by the property that the diagram with backward pointing arrows commutes.

Building upon Lemma \ref{lem:splitting_family_index_bundles_functorial}, we have:
\begin{lem}
For each  pair of maps $f_i \co \alpha_i \to \beta_i$, the following diagram commutes
  \begin{equation}
    \begin{tikzcd}
      \scrJ_{\fib}({f}_1, \epsilon_1) \times \scrJ_{\fib}({f}_2, \epsilon_2)  \ar[r] \ar[d] &  \scrJ_{\fib}(\epsilon_1) \times  \scrJ_{\fib}(\epsilon_2) \ar[d] \\
      \scrJ_{\fib}({f}_1  \times {f}_2, \epsilon_1 \times \epsilon_2) \ar[r] & \scrJ_{\fib}(\epsilon_1 \times \epsilon_2) .
    \end{tikzcd}
  \end{equation}\qed
\end{lem}
From the commutativity of splittings of Diagram \eqref{eq:family-vector-bundles-concatenation-is-functorial}, we conclude that the comparison maps of the product of stable complex structures are functorial in the sense that they fit in a commutative cube, with one face given by the direct sum of the maps in Diagram \eqref{eq:diagram_induced_by_fibrewise_stable_complex_data_for_map} with
\begin{multline}
 V'_{\beta_1 \times \beta_2} /V'_{\alpha_1 \times \alpha_2} \oplus T^{\beta_1 \times \beta_2}  \fcMR_{\alpha_1 \times \alpha_2} \cong \\ V'_{\beta_1}/V'_{\alpha_1} \oplus T^{\beta_1} \fcMR_{\alpha_1}  \oplus V'_{\beta_2}/V'_{\alpha_2} \oplus T^{\beta_2} \fcMR_{\alpha_2}
\end{multline}
and the other the analogue of Equation \eqref{eq:diagram_induced_by_fibrewise_stable_complex_data_for_map} for the morphism $f_1 \times f_2$. In particular, we obtain an isomorphism
\begin{equation}
  \begin{tikzcd}
    \begin{gathered}
W^{\fib}_{\alpha_1}  \oplus    T^{\alpha_1} \fcMR_{\alpha_1}(a,b) \oplus V^-_{\underline{a}} \oplus  V^+_{\underline{c}} \oplus\\ W^{\fib}_{\alpha_2}  \oplus   T^{\alpha_2} \fcMR_{\alpha_2}(b,c)   \oplus  V^-_{\underline{b}} \oplus V'_{\beta_1 \times \beta_2} /V'_{\alpha_1 \times \alpha_2} \oplus T^{\beta_1 \times \beta_2}  \fcMR_{\alpha_1 \times \alpha_2}
    \end{gathered} \ar[d] \\
 V^-_{\underline{c}} \oplus W^{\fib}_{\alpha_1 \times \alpha_2} \oplus I^{\fib,\bC}_{\beta_1 \times \beta_2} \oplus V^+_{\underline{a}}  ,
  \end{tikzcd}     
\end{equation}
over $  Z^{\epsilon_1,\ori}_{ \alpha_1} \times Z^{\epsilon_2, \ori}_{\alpha_2} $, depending canonically on elements of $\scrJ_{\fib}({f}_1, \epsilon_1) \times  \scrJ_{\fib}({f}_2, \epsilon_2) $.

The isomorphism of Equation \eqref{eq:family-vector-bundles-multiplicative} is also multiplicatively associative in the sense that, if we are given four orbits $a$, $b$, $c$, and $d$, and a triple of charts with ends $( \fcMR_{{\alpha}_i}, \epsilon_i)_{i=1}^{3}$, with $\alpha_1 \in A(a,b)$, $\alpha_2 \in A(b,c)$, and $\alpha_3 \in A(c,d)$, then  we have commutative diagrams
\begin{equation} \label{eq:commutative-diagram-triple-concatenation-bundles-first-interval}
    \begin{tikzcd}
     I^{t,\epsilon_1 \times \epsilon_2 \times \epsilon_3}_{\alpha_1 \times \alpha_2 \times \alpha_3}  \ar[r] \ar[d] &  I^{t,\epsilon_1}_{\alpha_1}    \oplus W^{\fib}_{\alpha_2 \times \alpha_3} \oplus I^{\fib,\bC}_{\alpha_2 \times \alpha_3} \oplus V_{\underline{c}}^{-} \ar[d] \\
     \begin{gathered}
 I^{t,\epsilon_1 \times \epsilon_2}_{\alpha_1 \times \alpha_2}   \oplus W^{\fib}_{\alpha_3}  \\ \oplus I^{\fib,\bC}_{\alpha_3} \oplus V_{\underline{c}}^{-}       
     \end{gathered}
 \ar[r] &
 \begin{gathered}
   I^{t,\epsilon_1}_{\alpha_1}   \oplus W^{\fib}_{\alpha_2} \oplus  I^{\fib,\bC}_{\alpha_2} \\
   \oplus V_{\underline{b}}^{-}  \oplus I^{\fib,\bC}_{\alpha_3}    \oplus V_{\underline{c}}^{-}
 \end{gathered}
\end{tikzcd}
\end{equation}
for $0 \leq t \leq m_{ab}$, where we note that $V_{\underline{c}}^-$ appears because it is an intermediate orbit. For $m_{ab} \leq t \leq m_{ab} + m_{bc}$, we instead have a commutative diagram
 \begin{equation}    \begin{tikzcd}
     I^{t,\epsilon_1 \times \epsilon_2 \times \epsilon_3}_{\alpha_1 \times \alpha_2 \times \alpha_3}  \ar[r] \ar[d] & V_{\underline{a}}^- \oplus  T^{\alpha_1} \fcMR_{\alpha_1}(a,b) \oplus W^{\fib}_{\alpha_1} \oplus  I^{t- m_{ab}, \epsilon_2 \times \epsilon_3}_{\alpha_2 \times \alpha_3}   \ar[d] \\
     \begin{gathered}
I^{t,\epsilon_1 \times \epsilon_2}_{\alpha_1 \times \alpha_2}  \oplus W^{\fib}_{\alpha_3}\\ \oplus I^{\fib,\bC}_{\alpha_3}   \oplus V^-_{\underline{c}}
     \end{gathered} \ar[r] &
     \begin{gathered}
V_{\underline{a}}^- \oplus  T^{\alpha_1} \fcMR_{\alpha_1}(a,b) \oplus W^{\fib}_{\alpha_1} \oplus    \\  I^{t- m_{ab},\epsilon_2}_{\alpha_2}   \oplus W^{\fib}_{\alpha_3} \oplus I^{\fib,\bC}_{\alpha_3}   \oplus V^-_{\underline{c}}.
     \end{gathered}
\end{tikzcd}
\end{equation}
Finally, for $m_{ab} + m_{bc} \leq t \leq m_{ab} + m_{bc} + m_{cd}$, the following diagram commutes:
\begin{equation}\label{eq:commutative-diagram-triple-concatenation-bundles-third-interval}
  \begin{tikzcd}
     I^{t,\epsilon_1 \times \epsilon_2 \times \epsilon_3}_{\alpha_1 \times \alpha_2 \times \alpha_3}  \ar[r] \ar[d] &
     \begin{gathered}
V_{\underline{a}}^- \oplus  T^{\alpha_1} \fcMR_{\alpha_1}(a,b) \oplus \\ W^{\fib}_{\alpha_1} \oplus     I^{t- m_{ab}, \epsilon_2 \times \epsilon_3}_{\alpha_2 \times \alpha_3} 
     \end{gathered}  \ar[d] \\
     \begin{gathered}
V_{\underline{a}}^- \oplus  T^{\alpha_1 \times \alpha_2} \fcMR_{\alpha_1 \times \alpha_2}(a,c)\\ \oplus W^{\fib}_{\alpha_1 \times \alpha_2}    \oplus  I^{t -m_{ab} - m_{bc} ,\epsilon_3}_{\alpha_3}  
     \end{gathered} \ar[r] &
\begin{gathered}
V_{\underline{a}}^- \oplus  T^{\alpha_1} \fcMR_{\alpha_1}(a,b) \oplus W^{\fib}_{\alpha_1} \\ \oplus   V_{\underline{b}}^-   \oplus   T^{\alpha_2} \fcMR_{\alpha_2}(b,c) \oplus W^{\fib}_{\alpha_2} \\  \oplus         I^{t- m_{ab}- m_{bc}, \epsilon_3}_{\alpha_3}  .
\end{gathered}
\end{tikzcd}
\end{equation}

Returning to Lemma \ref{lem:splitting_family_index_bundles_functorial}, we have:
\begin{lem}
Given four orbits $a$, $b$, $c$, and $d$, a triple  $\alpha_1 \in A(a,b)$, $\alpha_2 \in A(b,c)$, and $\alpha_3 \in A(c,d)$, and a triple  $\{ \epsilon_i \}_{i=1}^{3}$ of strip-like ends $\epsilon_i \in \Ends(\alpha_i)$, the following diagram commutes:
  \begin{equation}
    \begin{tikzcd}
      \scrJ_{\fib}(\epsilon_1) \times  \scrJ_{\fib}(\epsilon_2)\times  \scrJ_{\fib}(\epsilon_3) \ar[r] \ar[d] & \scrJ_{\fib}(\epsilon_1 \times \epsilon_2)\times  \scrJ_{\fib}(\epsilon_3)  \ar[d] \\
      \scrJ_{\fib}(\epsilon_1)\times  \scrJ_{\fib}( \epsilon_2 \times \epsilon_3)  \ar[r] & \scrJ_{\fib}(\epsilon_1 \times \epsilon_2 \times \epsilon_3 )
    \end{tikzcd}
  \end{equation}
 \qed
\end{lem}

From the commutativity of Diagrams \eqref{eq:commutative-diagram-triple-concatenation-bundles-first-interval}--\eqref{eq:commutative-diagram-triple-concatenation-bundles-third-interval}, we conclude that the comparison maps for product orientations are associative, i.e. that they give rise to the same map
  \begin{equation} \label{eq:fibre-stable-complex-structure-product-chart2}
  \begin{tikzcd}
    \begin{gathered}
       W^{\fib}_{\alpha_1}\oplus  T^{\alpha_1} \fcMR_{\alpha_1}(a,b) \oplus V^-_{\underline{a}} \oplus  W^{\fib}_{\alpha_2}\oplus   T^{\alpha_2} \fcMR_{\alpha_2}(b,c)  \oplus  V^-_{\underline{b}}\\ \oplus  V^+_{\underline{d}} \oplus W^{\fib}_{\alpha_3}  \oplus  T^{\alpha_3} \fcMR_{\alpha_3}(c,d)   \oplus  V^-_{\underline{c}}
    \end{gathered}
    \ar[d]  \\
      V^-_{\underline{d}} \oplus W^{\fib}_{\alpha_1 \times \alpha_2 \times \alpha_3} \oplus  I^{\fib,\bC}_{\alpha_1 \times \alpha_2 \times \alpha_3} \oplus  V^+_{\underline{a}} 
  \end{tikzcd}
\end{equation}
of vector bundles over $ \fcMR_{{\alpha}_1 \times {\alpha}_2 \times {\alpha}_3}$.

\subsection{The category of oriented charts}
\label{sec:categ-orient-charts}

\subsubsection{Oriented charts}
\label{sec:oriented-charts}
We begin this section by noting that the constructions of Section
\ref{sec:oriented-locus-space} may be performed in families: we
associate to each object ${\alpha}$ of $A(a,b)$ the contractible space
$\scrJ_{\fib}({\alpha})$, mapping to $\Ends({\alpha})$, whose fibre at
$\epsilon$ is $\scrJ_{\fib}(\epsilon)$.  We define $\scrJ({\alpha})$
to be the product $\scrJ_{\fib}({\alpha}) \times
\scrJ_{\base}(\alpha)$. 

We have natural evaluation maps
\begin{equation}
\scrJ({\alpha}) \leftarrow  \scrJ({f}) \to \scrJ({\beta}),
\end{equation}
so we can consider the disjoint union of the spaces $\scrJ({f})$ as the
morphism spaces of a category.
Lemma~\ref{lem:composition-morphisms-large-space-complex-structure}
asserts that composition is well-defined and associative.
Summarizing, we can make the following definition.

\begin{defin}
The topological category $\scrJ(a,b)$ is specified by
\begin{equation}
\ob \scrJ(a,b) = \coprod_{\alpha} \scrJ({\alpha}) \qquad
\mor \scrJ(a,b) = \coprod_{\epsilon \in \Ends({f})} \scrJ({f},  \epsilon).
\end{equation}
\end{defin}

We shall presently introduce a category of subsets of $
\scrJ(A)(a,b)$, whose purpose is to solve a technical problem
explained in Remark~\ref{rem:why-we-pass-to-subsets-of-oriented-charts} below. 
In order to set up the problem, recall that $\Ends({\alpha})$ is a
subspace of the space of choices of cylindrical ends on $
\partial^{P_\alpha} \fcMR_{\alpha} $. Given a subset $\cE \subset
\Ends({\alpha})$, we define the space of elements of $\scrJ(\alpha) $
lying over $\cE$ by the pullback
\begin{equation}
  \scrJ_{\cE}(\alpha) \equiv  \scrJ(\alpha) \times_{\Ends({\alpha}) }  \cE .
\end{equation}
Given a morphism $f \co \alpha \to \beta$, and subsets $\cE \subset \Ends({\alpha})$ and $\cF \subset \Ends({\beta})$, we define
\begin{equation}\label{eq:morphisms_category_oriented_charts}
  \scrJ^{\cF}_{\cE}(f) \equiv
  \begin{cases}
   \cE \times_{\Ends(\alpha)} \scrJ_{\cF}(\beta)  & \cE \subset f^* \cF \\
    \emptyset & \textrm{otherwise.}
  \end{cases}
\end{equation}
Given a pair of composable morphisms $f \co \alpha \to \beta$ and $g \co \beta \to \gamma$, and subsets  $\cE \subset \Ends({\alpha})$, $\cF \subset \Ends({\beta})$, and $\cG \subset \Ends({\gamma})$, note that the map $\scrJ(f) \times_{\cF}   \scrJ(g) \to  \scrJ(g \circ f)$ restricts to a map
\begin{equation} \label{eq:composition_category_oriented_charts}
  \scrJ_{\cE}^{\cF}(f) \times_{\cF}   \scrJ_{\cF}^{\cG}(g) \to  \scrJ_{\cE}^{\cG}(g \circ f). 
\end{equation}

Putting this together, we make the following definition.

\begin{defin}
The category $A^\ori(a,b)$ is the category whose space of objects is
the disjoint union over all objects $\alpha \in A(a,b)$ and over
subset $\cE$ of $\Ends(\alpha)$ of the spaces $
\scrJ_{\cE}({\alpha})$: 
  \begin{equation}
    \ob A^\ori(a,b) \equiv \coprod_{\alpha \in A(a,b)} \coprod_{\cE \subset \Ends({\alpha}) } \scrJ_{\scrE}({\alpha}).
  \end{equation}
The space of morphisms is defined to be the disjoint union over all morphisms $f$ in $A(a,b)$, and over all pairs $\cE \subset \Ends({\alpha})$ and $  \cF \subset \Ends({\beta}) $ of the spaces $  \scrJ_{\cE}^{\cF}(f)$:
    \begin{equation}
    \mor A^\ori(a,b) \equiv \coprod_{f \in A(a,b)} \coprod_{ \cE, \cF } \scrJ_{\cE}^{\cF}(f).
  \end{equation}
The source and target maps from $\mor A^\ori(a,b)$ to $\ob
A^\ori(a,b)$ are induced by the evident projection
$\scrJ_{\cE}^{\cF}(f) \to \scrJ_{\cE}(\beta)$ for the target and the
pullback along $f$ and projection on the source.  Composition in
$A^\ori(a,b)$ is induced from $\scrJ(a,b)$ via
Equation~\eqref{eq:composition_category_oriented_charts}.
\end{defin}
 
It will be convenient to analyse the topological category
$A^\ori(a,b)$ via the discrete category $A^{\Ends}(a,b)$ whose objects 
consist of an object $\alpha$ of $A(a,b)$ and a subset $\cE$ of
$\Ends(\alpha)$, and whose morphisms are given by a morphism $f$ in
$A(a,b)$ and the datum of the inclusion of $\cF$ in $f^* \cE$ if $\cE$
lifts to $\Ends(f)$. Before proceeding further, we note a significant
technical advantage of $A^\ori(a,b)$, which essentially follows by
construction.  
\begin{lem}
The source and target maps $s,t\colon \mor A^\ori(a,b) \to  \ob
A^\ori(a,b)$ are fibrations. \qed
\end{lem}
Note that the previous result fails for $\scrJ(a,b)$ because
the map $\Ends({f})\to \Ends({\beta})$ is not necessarily surjective.

We associate to each object $(\alpha,\cE)$ in $A^{\Ends}(a,b)$, the subset 
\begin{equation}
  Z^{\cE}_{\alpha}(a,b) \subset  Z_{\alpha}(a,b),
\end{equation}
consisting of points which lie in $Z^{\epsilon',\ori}_{\alpha}(a,b) $ for each choice of cylindrical end $\epsilon'$ contained in $\cE$.

\begin{lem}\label{lem:pres}
The assignment $ (\alpha,\scrE) \mapsto (\partial^{Q_\alpha}\cP(a,b), Z^{\cE}_{\alpha}(a,b),G_\alpha)$ extends to a functor $A^{\Ends}(a,b) \to \Orb$.
\end{lem}
\begin{proof}
A morphism $f \co \alpha \to \beta$ maps $Z_{\alpha}(a,b)$ to $ Z_{\beta}(a,b)$, and this maps takes the subset $ Z^{f^*\epsilon,\ori}_{\alpha}(a,b)$ to  $ Z^{\epsilon,\ori}_{\beta}(a,b)$ whenever $\epsilon$ lies in $\Ends(f)$. The result follows. 
\end{proof}

\begin{rem} \label{rem:why-we-pass-to-subsets-of-oriented-charts}
It is tempting to streamline the above construction and consider only the subset $\cE = \Ends(\alpha)$ as our object. Unfortunately, this assignment is not functorial in $A(a,b)$, which is the reason for passing to $A^{\Ends}(a,b)$ and for introducing the piecewise definition in Equation \eqref{eq:composition_category_oriented_charts}.
\end{rem}

The functor given in Lemma~\ref{lem:pres} is an orbispace presentation.

\begin{lem}
The functor specified in Lemma~\ref{lem:pres} by the assignment
$(\alpha, \cE) \mapsto (\partial^{Q_\alpha}\cP(a,b),
Z^{\cE}_{\alpha}(a,b),G_\alpha)$ is an orbispace presentation of
$\fcMR(a,b)$. 
\end{lem}

\begin{proof}
The key point is that, for each $[u] \in  \fcMR(a,b)$ which lifts to a
space $\Mbar^{\bR}_\alpha(a,b)$ and for each $\epsilon \in
\Ends(\alpha)$, there is some chart $\alpha'$, which differs from
$\alpha$ only in having a larger obstruction space $V_{\alpha'}$
containing $V_{\alpha}$, so that  $[u]$ lies in the image of the
projection map from $ Z^{\epsilon }_{\alpha}$. Indeed, the condition
of lying in $ Z_{\alpha}^{\epsilon,\ori}$ is a surjectivity statement
for a family of Fredholm operators parametrised by a closed interval,
hence can be achieved by enlarging the obstruction. This implies that
each point in $ \fcMR(a,b)$ lies in the image of a chart with
associated to $(\alpha, \cE)$ for some $\cE$. To prove that the nerve
of the category of charts covering a given point is contractible, we
use the same argument as Proposition
\ref{prop:Ham-Kuranishi-presentation}: the key point is that, if the
data in charts $\alpha$ and $\alpha'$ agree except for the choice of
obstruction bundle $V_\alpha$, and $[u]$ lies in the image of the
charts associated to $(\alpha,\cE)$ and $(\alpha', \cE')$, then $[u]$
will lie in the image of a chart associated to $(\alpha'', \cE' \amalg
\cE) $, where the data for $\alpha''$ agrees with that of $\alpha$,
except that the obstruction bundle contains both $V_\alpha$ and
$V_{\alpha'}$. 
\end{proof}
For the next statement, we use the notion of a topologically enriched orbifold presentation, which is a mild extension of
Definition~\ref{defn:orbipres}; the indexing category is an internal
category, and the required homeomorphism compares $\fcMR(a,b)$ to the
colimit over $\pi_0 A^{\ori}(a,b)$
(c.f. Proposition~\ref{prop:internal_discretize}).
\begin{cor}
The composition $A^{\ori}(a,b) \to A^{\Ends}(a,b) \to \Orb$ is a
topologically enriched orbispace presentation of $\fcMR(a,b)$.
\end{cor}

\begin{proof}
It suffices to verify that the functor $A^{\ori}(a,b) \to
A^{\Ends}(a,b)$ induces an acyclic fibration of nerves.  But this
follows from Quillen's theorem B since each relevant overcategory is
contractible and the induced comparison maps are evidently
equivalences.
\end{proof}

Recall that we defined the tangent bundle of $\fcMRreg_{\alpha}(a,b)$
to be the vector bundle 
\begin{equation}
  T \fcMR_{\alpha}(a,b) \equiv T^\alpha \fcMR_{\alpha}(a,b) \oplus T \cB_{\alpha}.
\end{equation}
In a mild abuse of notation, we denote its restriction to $  Z^{\scrE}_{\alpha}(a,b)$ by $T \fcMR_{\alpha}(a,b) $. We analogously introduce the complex vector bundle
\begin{equation}
  I^\bC_{\alpha} \equiv I^{\fib,\bC}_{\alpha} \oplus  I^{\base,\bC}_{\alpha}
\end{equation}
and the vector space
\begin{align}
 W_{\alpha} & \equiv W^{\fib}_{\alpha} \oplus W^{\base}_{\alpha} . \end{align}
Taking the direct sum of the isomorphisms in Equations \eqref{eq:stable-almost-complex-structure-base-1} and \eqref{eq:stable_fibrewise-complex-structure}, we have:
\begin{lem} \label{lem:oriented-charts-family-iso}
  The space  $\scrJ_{\scrE}(\alpha)$ parametrises a continuous family of isomorphisms
  \begin{equation}\label{eq:oriented-charts-family-iso}
    V^+_b
    \oplus \bR^{O_\alpha}  \oplus  T \fcMR_{\alpha}(a,b)  \oplus    W_{\alpha}  \oplus    V^-_a \cong  V^-_b  \oplus I^\bC_{\alpha} \oplus    W_{\alpha} \oplus  V^+_a
    \end{equation}
of vector bundles over $ Z^{\scrE}_{\alpha}(a,b)$.   \qed
\end{lem}

We associate to a morphism in $A^\ori(a,b) $ the  vector space
\begin{align}
   W_{{f}} & \equiv  W^{\fib}_{{f}} \oplus W^{\base}_{{f}} .
 \end{align}
Each  element of $\scrJ_{\scrE}^{\scrF}(f)$ induces isomorphisms
  \begin{align}
  \bR^{Q_\beta \setminus Q_\alpha} \oplus  T \fcMR_{\alpha}(a,b) & \cong  T \fcMR_{{\beta}}(a,b) \oplus V_\alpha/V_\beta \\
    I^\bC_{\alpha} & \cong  I^\bC_{{\beta}} \oplus V_\alpha/V_\beta  \\
     W_{\alpha} & \cong W_{{\beta}} \oplus W_{{f}}
  \end{align}
  where the first two arise from the splitting of Diagram \eqref{eq:short-exact-sequence-family-vector-bundles-change-chart} that is part of the definition of an element of $\scrJ(f)$, and the second two are straightforward direct sum decompositions.
  
The next result asserts that the stable complex structure on $ T \fcMR_{\alpha} $ is identified with the direct sum of the stable complex structure on $ T \fcMR_{{\beta}} $ with $V_\alpha/V_\beta \oplus  W_{{f}}   $. \begin{lem}
For each morphism  ${f}$ from ${\alpha}$ to ${\beta}$, and for each element of  $\scrJ_{\scrE}^{\scrF}(f)$, the following diagram of $G_\alpha$ equivariant vector bundles over $ Z^{\cE}_{\alpha}(a,b) $ commutes:
  \begin{equation}
    \begin{tikzcd}
      \begin{gathered}
      V^+_b 
      \oplus \bR^{O_\beta} \oplus  T \fcMR_{{\beta}}(a,b)   \oplus    W_{{\beta}}  \oplus  V^-_a \\ \oplus V_\alpha/V_\beta \oplus  W_{{f}}  
      \end{gathered}
       \ar[r] \ar[d] &
       \begin{gathered}
V^+_b  
      \oplus \bR^{O_\alpha} \oplus  T \fcMR_{\alpha}(a,b) \\  \oplus    W_{\alpha}  \oplus  V^-_a          
       \end{gathered}
 \ar[d] \\
 \begin{gathered}
V^-_b   \oplus  I^\bC_{{\beta}}   \oplus     W_{{\beta}}   \oplus V^+_a  \\ \oplus   V_\alpha/V_\beta \oplus  W_{{f}}    
 \end{gathered}
 \ar[r] &
 \begin{gathered}
V^-_b  \oplus  I^\bC_{\alpha} \oplus    W_{{\beta}} \oplus V^+_a.   
 \end{gathered}
    \end{tikzcd}
  \end{equation}
  \qed
\end{lem}

For the statement of the next result, we specify that we interpret Lemma  \ref{lem:oriented-charts-family-iso} as the assignment of a stable isomorphism between $T \fcMR_{\alpha}(a,b)  \oplus    W_{\alpha}  $ and $I^\bC_{\alpha} \oplus    W_{\alpha}  $, relative the stable vector spaces $V_a$ and $V_b$, and the set $O_\alpha$. Returning to Definition \ref{def:Complex-oriented-Kuranishi-chart}, 
we have:

\begin{cor}
The family of isomorphisms in Lemma  \ref{lem:oriented-charts-family-iso} defines a (topologically enriched) lift of the functor $A^\ori(a,b) \to \Kur^{fs}(a,b) $ to $ \Kur^{\ori}(a,b)$. \qed
\end{cor}

We now consider the multiplicativity of this construction:  the results of  Sections \ref{sec:stable-compl-struct} and Section \ref{sec:mult-stable-compl} imply that the collection of categories $\scrJ(a,b)$ are the $1$-cells of a bicategory over $A(a,b)$.

Taking the product of subsets of $\Ends(\alpha_1)$ and $\Ends(\alpha_2)$, we obtain a (topologically enriched) functor
\begin{equation} \label{eq:product_oriented-charts}
    A^\ori(a,b) \times A^\ori(b,c) \to A^\ori(a,c),
  \end{equation}
which we think of as the product of oriented charts. There is a natural isomorphism between the two functors
\begin{equation}
    A^\ori(a,b) \times A^\ori(b,c)  \times A^\ori(c,d) \to A^\ori(a,d),
\end{equation}
arising from the associativity of products of sets. Since the construction of $A^\ori(a,b)$ is invariant under the action of $\Pi$, we obtain a $\Pi$-equivariant bicategory $A^\ori$, with $0$-cells $a \in \cP(f,H,f)$, $1$-cells $A^\ori(a,b)$ enriched in $\Top$, and $2$-compositions given by Equation \eqref{eq:product_oriented-charts}.

We now observe that we have natural isomorphisms 
\begin{align}
 T \fcMR_{{\alpha}_1\times {\alpha}_2}(a,c) & \cong  T \fcMR_{{\alpha}_1}(a,b)  \oplus  T \fcMR_{{\alpha}_2}(b,c) \\ 
      W_{{\alpha}_1 \times {\alpha}_2} & \cong W_{{\alpha}_1} \oplus V^-_b \oplus      W_{{\alpha}_2}  \\
      I^\bC_{{\alpha}_1 \times {\alpha}_2} & \cong  I^\bC_{{\alpha}_1} \oplus I^\bC_{{\alpha}_2} \\
      \bR^{O_{\alpha_1 \times \alpha_2}}   & \cong  \bR^{O_{\alpha_1}}  
      \oplus \bR^{O_{\alpha_2}},
     \end{align}

We finally conclude that we have constructed a complex oriented Kuranishi flow category (in the sense of Definition~\ref{def:complex-oriented-flow-category}):

\begin{lem}
The family of isomorphisms in Lemma~\ref{lem:oriented-charts-family-iso} defines a (topologically enriched) $\Pi$-equivariant lift of the Kuranishi presentation $A^\ori \to A \to \Kur$ to the bi-category  $ \Kur^{\ori}$ (equipped with the trivial $\Pi$-action). \qed
\end{lem}

\appendix
\part{Appendices}

\section{Groups, categories, and spectra}
\label{sec:groups-categ-spectra}

The purpose of this section is to provide necessary technical
background for our work.

\subsection{Background on the category of $G$-spectra}
\label{sec:backgr-categ-g}

In this subsection, we give a concise review of the point-set model of
the equivariant stable category we use, the category of orthogonal
$G$-spectra.  This is a symmetric monoidal category equipped with a
notion of stable equivalence detected by the equivariant stable
homotopy groups.  The canonical reference for this category
is~\cite{MandellMay}; see also the appendix
to~\cite{HillHopkinsRavenel} for another treatment with a slightly
different emphasis.

\subsubsection{The point-set category of orthogonal $G$-spectra}\label{sec:ortho}

We will assume throughout this section a choice of finite group $G$.
A universe $U$ for $G$ is a countably infinite-dimensional real inner
product space containing $\bR^{\infty}$, equipped with a linear
$G$-action that preserves the inner product, which is the direct sum
of finite-dimensional $G$-inner product spaces, and such that any
finite-dimensional $G$-inner product in $U$ occurs infinitely often.
Two extreme examples of universes are given by the trivial universe
$\bR^{\infty}$ and the complete universes which contain all
irreducible finite-dimensional $G$-representations.  We will often
work with the specific model of the complete universe given by the
countable sum of the regular representation which we denote $\rho$.
We write $\cV(U)$ for the set of finite dimensional $G$-inner product
spaces that are isomorphic to a $G$-stable subspace of $U$.

We write $\sI_G^U$ for the category with object set $\cV(U)$ and with
morphisms $ \sI_G^U(V,W)$ for $V, W \in \cV(U)$ given by the $G$-space
of isometric isomorphisms $V \to W$, where $G$ acts by conjugation.
Let $\cT^G$ denote the category of $G$-spaces and $G$-equivariant
maps, regarded as enriched over itself where the mapping $G$-space is
the space of non-equivariant maps equipped with the conjugation
action.

\begin{defin}[Definition II.2.6 of~\cite{MandellMay}]
An orthogonal $G$-spectrum is an enriched functor $X \colon \sI_G^U
\to \cT^G$ equipped with structure maps
\begin{equation}
\sigma_{V,W} \colon X(V) \sma S^W \to X(V \oplus W)
\end{equation}
that comprise a natural transformation of functors and are associative
and unital.  A map of orthogonal $G$-spectra is a natural
transformation that commutes with the structure map.  We denote the
category of orthogonal $G$-spectra by $\Sp_G$; we will write $\Sp_G^U$
when it is necessary to emphasize the universe $U$.
\end{defin}

\begin{rem}
The category $\Sp^U_G$ only depends on the isomorphism classes of finite dimensional representations appearing as invariant subspaces of $U$; it would thus be cleaner to label the category by the set of such isomorphism classes, but this conflicts with the established literature. The choice of universe will enter, later, in the construction of fibrant replacement functors. 
\end{rem}

The following proposition records the existence of enriched limits and
colimits in $\Sp_G^U$; this can be deduced as a consequence
of~\cite[II.4.3]{MandellMay}, as the following paragraph there
explains.

\begin{prop}\label{prop:cocomplete}
The category $\Sp_G^U$ is complete and cocomplete.  Moreover, it is
tensored and cotensored over based $G$-spaces; for a based $G$-space
$A$ and an orthogonal $G$-spectrum $X$, the tensor $A \sma X$ is given
by the levelwise smash product which has $V$th space $A \sma X(V)$ and
the cotensor $F(A, X)$ is given by the levelwise function space which
has $V$th space $\Map(A,X(V))$. 
\end{prop}

An important advantage of $\Sp_G^U$ over earlier models of the
equivariant stable category is that it is symmetric monoidal.  We
write $\bS$ for the orthogonal $G$-spectrum with $V$th space the
representation sphere $S^V$ (i.e., the one-point compactification of
the $G$-representation $V$).  To see this, we give another description
of orthogonal $G$-spectra as diagrams.

Let $\sJ_G^{\sV}$ denote the category with objects the elements of
$\sV(U)$ and morphism spaces defined as the following Thom spaces: let
$E(V,V')$ be the total space of the sub-bundle of the product bundle
$\sI(V,V') \times V'$ specified by taking pairs $(f,x)$ such that $x
\in V' - f(V)$, where $V' - f(V)$ denotes the orthogonal complement.
Then $\sJ_G^{\sV}(V,V')$ is the associated equivariant Thom space.
This is a symmetric monoidal category under the operation $\oplus$ given by
\begin{equation}
(f,x) \oplus (g,x') = (f \oplus g, x + x').
\end{equation}
The category of orthogonal $G$-spectra is precisely the symmetric
monoidal category of continuous (enriched) functors from $\sJ_G^{\sV}$
to spaces.  In light of this, the following result is just an
application of the Day convolution.

\begin{prop}[Theorem II.3.1 of \cite{MandellMay}]
The category $\Sp_G^U$ is a symmetric monoidal category with respect
to the smash product $\sma$ of orthogonal $G$-spectra, where the unit
is $\bS$.
\end{prop}

Next, we recall that although the universe $U$ will play a central role
in determining the homotopy theory of $\Sp_G^U$ which we discuss in
Appendix~\ref{sec:homot-theory-orth} below, its role in the 
point-set theory can be elided as follows:
\begin{prop}[Theorem V.1.5 of \cite{MandellMay}]\label{prop:change-of-universe}
For universes $U$ and $U'$, there is a change of universe functor
\begin{equation}
\sI_{U}^{U'} \colon \Sp_G^U \to \Sp_G^{U'}.
\end{equation}
that is strong symmetric monoidal.  When $U = U'$, this functor is the
identity, and given $U$, $U'$, and $U''$, there is a natural
isomorphism $\sI_{U'}^{U''} \circ \sI_U^{U'} \cong \sI_{U}^{U''}$.  In
particular, $\sI_{U}^{U'}$ is always an equivalence of categories.
\end{prop}

We refer to the monoid objects in $\Sp_G^U$ as associative ring
orthogonal $G$-spectra and the commutative monoid objects in $\Sp_G^U$
as commutative ring orthogonal $G$-spectra.  The following is an easy
consequence of Proposition~\ref{prop:cocomplete} and the analysis
of~\cite[\S II.7]{ElmendorfKrizMandellMay1997}.

\begin{prop}
The categories of monoids and commutative monoids in $\Sp_G^U$, which
we respectively denote $\Ass_G$ and $\Comm_G$, are complete and
cocomplete.  The categories $\Ass_G$ and $\Comm_G$ are tensored and
cotensored over unbased $G$-spaces; for an unbased $G$-space $A$, the
cotensor is created in $\Sp_G$.
\end{prop}

We make frequent use of an ``external'' smash product.  Given an
orthogonal $G_0$-spectrum $X_0$ and an orthogonal $G_1$-spectrum
$X_1$, we want to regard the smash product $X_0 \sma X_1$ as a $G_0
\times G_1$ spectrum.  One way to do this is simply to regard $X_0$ as
a $G_0 \times G_1$-spectrum with trivial $G_1$-action and $X_1$ as a
$G_0 \times G_1$-spectrum with trivial $G_0$-action; then $X_0 \sma
X_1$ is an orthogonal $G_0 \times G_1$-spectrum.  However, it is
technically convenient to have an explicit model of this process.

The external smash product~\cite[II.2.4]{MandellMay} produces an
orthogonal $G_0 \times G_1$-spectrum $X_1 \bar{\sma} X_2$ indexed on
the product category $\sJ_{G_0}^{\sV_0} \times \sJ_{G_1}^{\sV_1}$, so that
the value at $(V_1, V_1')$ is $X_1(V_1) \sma
X_2(V_1')$.  To internalize this, we use the direct sum map
\begin{equation} \label{eq:direct_sum_operation_external-smash}
\oplus \colon \sJ_{G_0}^{\sV_0} \times \sJ_{G_1}^{\sV_1} \to \sJ_{G_0
  \times G_1}^{\sV_0 \oplus \sV_1}.
\end{equation}
On morphisms, this functor is specified for $(f,x) \in
\sJ_{G_0}^{\sV_0}(V,V')$ and $(g,y) \in \sJ_{G_1}^{\sV_1}(W,W')$ by
the element of $\sJ_{G_0 \times G_1}^{\sV_0 \oplus \sV_1}(V \oplus W,
V' \oplus W')$
\begin{equation}
(f,x) \oplus (g,y) \mapsto (f \oplus g, x + y).
\end{equation}
Left Kan extension along the direct sum in Equation \eqref{eq:direct_sum_operation_external-smash} now produces the external
smash functor:

\begin{lem}\label{lem:external-smash-product}
For finite groups $G_1$ and $G_2$, there is an external smash product
functor 
\begin{equation}
\Sp^{U_1}_{G_1} \times \Sp^{U_2}_{G_2} \to \Sp^{U_1 \oplus U_2}_{G_1
  \times G_2}
\end{equation}
which is associative in the sense that for groups $G_1$, $G_2$, and
$G_3$, the diagram
\begin{equation}
\begin{tikzcd}
\Sp^{U_1}_{G_1} \times \Sp^{U_2}_{G_2} \times \Sp^{U_3}_{G_3} \ar[r]
\ar[d] & \Sp^{U_1 \oplus U_2}_{G_1 \times G_2} \times \Sp^{U_3}_{G_3}
\ar[d] \\
\Sp^{U_1}_{G_1} \times \Sp^{U_2 \oplus U_3}_{G_2 \times G_3} \ar[r] & 
\Sp^{U_1 \oplus U_2 \oplus U_3}_{G_1 \times G_2 \times G_3} 
\end{tikzcd}
\end{equation}
commutes.
\end{lem}

\begin{proof}
The associativity of the external smash product follows from the fact
that the diagram 
\begin{equation}
\begin{tikzcd}
\sJ_{G_0}^{\sV_0} \times \sJ_{G_1}^{\sV_1} \times \sJ_{G_2}^{\sV_2}
\ar[r,"\id \times \oplus"] \ar[d,swap,"\oplus \times \id"] &
\sJ_{G_0}^{\sV_0} \times \sJ_{G_1 \times G_2}^{\sV_1 \oplus \sV_2}
\ar[d,"\oplus"] \\ \sJ_{G_0 \times G_1}^{\sV_0 \oplus \sV_1} \times
\sJ_{G_2}^{\sV_2} \ar[r,swap,"\oplus"] & \sJ_{G_0 \times G_1 \times
  G_2}^{\sV_0 \oplus \sV_1 \oplus \sV_2}
\end{tikzcd}
\end{equation}
evidently commutes.
\end{proof}

We now turn to a discussion of the internal function objects in
$\Sp_G$; it is a closed symmetric monoidal category.  The category
$\Sp_G$ can be viewed as a category of presheaves of topological
spaces and therefore comes with a natural enrichment; we denote the
mapping space as $\Map_{\Sp_G}(X,Y)$.

\begin{defin}
For $X,Y$ orthogonal $G$-spectra, the mapping orthogonal $G$-spectrum
$F(X,Y)$ is specified by the formula
\begin{equation}
F(X,Y)(W) = \Map_{\Sp_G}(X, Y[W]),
\end{equation}
where $Y[W](V) = Y(W \oplus V)$.
\end{defin}

For spectra $X$ and $Y$, there is a natural evaluation map
\begin{equation}
X \sma F(X,Y) \to Y
\end{equation}
which is the adjoint of the identity map $F(X,Y) \to F(X,Y)$.  Given
spectra $X$, $Y$, $W$, and $Z$, the evaluation maps induce a map
\begin{equation}
X \sma W \sma F(X,Y) \sma F(W,Z) \to X \sma F(X,Y) \sma W \sma F(W,Z)
\to Y \sma Z.
\end{equation}
The adjoint of this map is the smash product map on function spectra
\begin{equation}
F(X,Y) \sma F(W,Z) \to F(X \sma W, Y \sma Z).
\end{equation}

The associativity of the smash product implies the following
proposition.

\begin{prop}\label{prop:associative-smash-function-spectra}
The smash product map of function spectra is associative and unital in
the sense that the diagrams
\begin{equation}
\begin{tikzcd}
F(X_1, Y_1) \sma F(X_2, Y_2) \sma F(X_3, Y_3) \ar[r] \ar[d] & F(X_1
\sma X_2, Y_1 \sma Y_2) \sma F(X_3, Y_3) \ar[d] \\ F(X_1, Y_1) \sma
F(X_2 \sma X_3, Y_2 \sma Y_3) \ar[r] & F(X_1 \sma X_2 \sma X_3, Y_1
\sma Y_2 \sma Y_3)
\end{tikzcd}
\end{equation}
and
\begin{equation}
\begin{tikzcd}
\bS \sma F(X,Y) \ar[r] \ar[dr] & F(\bS,\bS) \sma F(X,Y) \ar[d] \\ &
F(X,Y)
\end{tikzcd}
\end{equation}
commute.
\end{prop}

\subsubsection{Shift desuspension}\label{sec:negative-MM-spheres} 

We now introduce the functors relating $G$-spaces to $G$-spectra, and which are analogous in the non-equivariant setting to the functors which assign to each space its suspension spectrum and to each spectrum the corresponding space assigned to the spheres $S^0$. For convenience, we consider a slightly more general situation, wherein we choose a finite-dimensional real representation $V$ that is isomorphic to one in the universe $U$, to which we shall associate an adjoint pair 
\begin{equation}
\begin{tikzcd}
\Ev_V  \colon \Sp_G \ar[r, shift left=1ex] & \ar[l, shift left=1ex] \cT^G \colon F_V. 
\end{tikzcd}
\end{equation}

\begin{defin}
The functor $\Ev_V$ evaluates an orthogonal $G$-spectrum at $V$, i.e.,
$\Ev_V(X) = X(V)$ and the functor $F_V$ is the left adjoint to
$\Ev_V$.  For $V \subset W$, $F_V$ is given by the formula
\begin{equation}
F_V X (W) = O(W)_+ \sma_{O(W - V)} \Sigma^{W-V} X.
\end{equation}
\end{defin}
The left adjoint $F_V$ is known as the ``shift desuspension functor''.
When $V = 0$, $F_0 X$ is a model of the suspension spectrum functor
typically denoted by $\Sigma^{\infty}$ and $\Ev_V$ is the functor
typically denoted by $\Omega^\infty$.  Notice that there is an
isomorphism 
\begin{equation}
\Sigma^\infty X = F_0 X \cong F_0 S^0 X \cong S \sma X.
\end{equation}

Mapping out of suspension spectra is straightforward: when $X =
\Sigma^{\infty} A = F_0 A$, we have
\begin{align}
F(\Sigma^{\infty} A, Y)(W) &= F(F_0 A, Y)(W) = \Map_{\Sp_G}(F_0 A,
Y[W]) \\ &\cong \Map_{\cT^G}(A, Y[W](0)) \cong \Map_{\cT^G}(A, Y(W)).
\end{align}

When $V$ is not $0$, $F_V S^0$ gives a specific model of the
``negative $V$-sphere''.  The shift desuspension functors are
multiplicative in the sense that there is a natural isomorphism
\begin{equation}
F_V A \sma F_W B \cong F_{V \oplus W} A \sma B.
\end{equation}
(See Lemma II.4.8 of~\cite{MandellMay}.)

The shift desuspension and evaluation functors are compatible with the
external smash product.

\begin{lemma}
Let $A_0$ be a $G_0$-space and $A_1$ be a $G_1$-space.  For any
$G_0$-representation $U_0$ and $G_1$-representation $U_1$, there is a
natural homeomorphism
\begin{equation}
F_{U_0} A_0 \sma F_{U_1} A_1 \cong F_{U_0 \oplus U_1} A_0 \sma A_1
\end{equation}
\end{lemma}

\begin{lemma}
Let $A_0$ be a $G_0$-space and $A_1$ be a $G_1$-space.  For any
$G_0$-representation $U_0$ and $G_1$-representation $U_1$, there is a
natural map
\begin{equation}
\Ev_{U_0} A_0 \sma \Ev_{U_1} A_1 \to \Ev_{U_0 \oplus U_1} A_0 \sma A_1
\end{equation}
\end{lemma}

\subsubsection{Suspension, loops, and shift}\label{sec:shift}

For any finite-dimensional
$G$-representation $V$, there is an adjoint pair $(\Sigma^V,
\Omega^V)$ of suspension and loop endofunctors.  These functors are
defined in terms of the tensor and cotensor with the based space
$S^V$, i.e., $\Sigma^V X = S^{V} \sma X$ and $\Omega^V X = F(S^V, X)$.

\begin{lemma}
Let $X_0$ be an orthogonal $G_0$-spectrum and $X_1$ be an orthogonal
$G_1$-spectrum.  For any $G_0$-representation $U_0$ and
$G_1$-representation $U_1$, there is a homeomorphism
\begin{equation}
S^{U_0} X_0 \sma S^{U_1} X_1 \to S^{U_0 \oplus U_1} X_0 \sma X_1.
\end{equation}
and a natural map 
\begin{equation}
\Omega^{U_0} X_0 \sma \Omega^{U_1} X_1 \to \Omega^{U_0 \oplus U_1} X_0
\sma X_1.
\end{equation}
\end{lemma}

There is another model of the suspension functor which is sometimes
useful: 

\begin{defin}\label{defin:V-shift}
Let $V$ be a finite-dimensional $G$-representation and $X$ be an
orthogonal $G$-spectrum.  Then we define the $V$-shift functor $\Sh_V$
applied to $X$ via the formula 
\begin{equation}
(\Sh_V X)(W) = X(V \oplus W),
\end{equation}
where the structure maps are induced from those of $X$ and the
orthogonal action via direct sum with the identity on $V$.
\end{defin}

For each $V$, there is a natural transformation
\begin{equation}
S^V \sma (-) \to \Sh_V (-),
\end{equation}
induced by the structure map of $X$, and an easy computation with
stable homotopy groups shows that this is always a stable equivalence.

The functoriality of the shift functor in $V$ can be summarized as
follows.

\begin{lem}\label{lem:shift-functoriality}
Let $X$ be an orthogonal $G$-spectrum.  The construction
$\Sh_{(-)}(X)$ specifies a functor from $\sJ_G^{\sV}$ to orthogonal
$G$-spectra.
\end{lem}

In particular, iterating the shift functor can be identified with the
direct sum of representations: for representations $V$ and $W$, there
is a natural homeomorphism
\begin{equation}
\Sh_V \Sh_W X \cong \Sh_{V \oplus W} X.
\end{equation}

We can also obtain ``translation'' morphisms of the following form:

\begin{lem}\label{lem:shift-translation}
For $U \subset V$, we have a natural map
\begin{equation}
\Sh_U X \to \Omega^{V - U} \Sh_V X
\end{equation}
induced by adjunction from the map
\begin{equation}
S^{V-U} \sma \Sh_U X \to \Sh_V X.
\end{equation}
For any $U$, $V$, and $W$, this induces a natural map
\begin{equation}\label{eq:loop-shift}
\Omega^{U} \Sh_W X \to \Omega^{U \oplus V} \Sh_{W \oplus V} X
\end{equation}
constructed as the composite
\begin{equation}
\Omega^{U} \Sh_W X \to \Omega^{U} \Omega^{V} \Sh_{W \oplus V} X \cong
\Omega^{U \oplus V} \Sh_{W \oplus V} X,
\end{equation}
using Equation~\eqref{eq:loop-shift} applied to the inclusion $W \to W
\oplus V$ to produce the map $\Sh_W X \to \Omega^V \Sh_{W \oplus V}
X$.
\end{lem}

The shift functor is compatible with the external product.

\begin{lem}
There are natural transformations
\begin{equation}
\Sh_{W_1} \sma \Sh_{W_2} \to \Sh_{W_1 \oplus W_2}
\end{equation}
induced by the maps
\begin{equation}
\begin{aligned}
X(V_1 \oplus W_1) \sma X(V_2 \oplus W_2) &\longrightarrow X(V_1 \oplus W_1 \oplus
V_2 \oplus W_2) \\
&\longrightarrow X(V_1 \oplus V_2 \oplus W_1 \oplus W_2).
\end{aligned}
\end{equation}
These transformations are associative and unital.
\end{lem}

It is straightforward to check that the natural transformation $S^V
\sma (-) \to \Sh_V (-)$ is externally monoidal:

\begin{lem}\label{lem:shift-suspend-monoidal-compare}
For finite groups $G_1$ and $G_2$ and representations $V_1$ and $V_2$
respectively, for any $G_1$-spectrum $X$ and $G_2$-spectrum $Y$ the
diagrams
\begin{equation}
\begin{tikzcd}
(S^{V_1} \sma X) \sma (S^{V_2} \sma Y) \ar[r] \ar[d] & S^{V_1 \oplus
    V_2} \sma (X \sma Y) \ar[d] \\ \Sh_{V_1} X \sma \Sh_{V_2} Y \ar[r]
  & \Sh_{V_1 \oplus V_2} (X \sma Y)
\end{tikzcd}
\end{equation}
commutes.  Moreover, the analogous associativity and unitality
diagrams commute.
\end{lem}

\subsubsection{The category of equivariant spectra}

In our work, we will need to consider the category that combines all
equivariant spectra as the group $G$ varies.

\begin{defin}
Let $\Sp_{\eq}$ denote the category with: 
\begin{enumerate}
\item Objects pairs $(G, X)$ where $G$ is a finite group.
\item Morphisms from $(G,X) \to (G',X')$ specified by a group
  homomorphism $f\colon G \to G'$ and a map $X \to f^* X'$.
\end{enumerate}  
Here we are considering the equivariant spectra indexed on the trivial
universe.
\end{defin}

Lemma~\ref{lem:external-smash-product} implies that $\Sp_{\eq}$ is a
symmetric monoidal category.

\begin{prop}
The category $\Sp_{\eq}$ is a symmetric monoidal category product given by
the external product of spectra:
\begin{equation}
(G_1, X_1) \sma (G_2, X_2) = (G_1 \times G_2, X_1 \sma X_2),
\end{equation}
and the unit is $(\{e\}, \bS)$.
\end{prop}

\begin{proof}
Functoriality follows from the observation that given pairs $(G_1, X_1)$
and $(G_2, X_2)$ and a map $f \colon (G_1, X_1) \to (G_1', X_1')$, there
is a natural map of $(G_1 \times G_2)$-spectra
\begin{equation}
X_1 \sma X_2 \to (f \times \id)^*(X_1' \sma X_2) 
\end{equation}
induced from the identification
\begin{equation}
f^* X_1' \sma X_2 \cong (f \times \id)^*(X_1' \sma X_2). 
\end{equation}
\end{proof}

We will in practice often work with the subcategory where the maps $f
\colon (G_1, X_1) \to (G_1', X_2')$ are given by surjections $G_1 \to
G_1'$.

\begin{defin}
Let $\Sp_{\eq, \surj}$ denote the symmetric monoidal subcategory of
$\Sp_{\eq}$ with:
\begin{enumerate}
\item Objects pairs $(G, X)$ where $G$ is a finite group.
\item Morphisms from $(G,X) \to (G',X')$ specified by a surjective group
  homomorphism $f\colon G \to G'$ and a map $X \to f^* X'$.
\end{enumerate}  
\end{defin}

In order to express the compatibly of the enrichments on the different
mapping spectra in this category, is it useful to take a more
sophisticated view of $\Sp_{\eq}$ and $\Sp_{\eq, \surj}$ and regard
these as symmetric monoidal index enriched categories; these are
equivalently categorical fibrations of symmetric monoidal categories.
Here the indexing category is the category of finite groups.  We do
not make serious use of this perspective, but refer the reader
to~\cite{Shulman2013, Bohmann2014} for a discussion of this point.

\subsubsection{Fixed points and orbits}

For a $G$-space $X$, there are two reasonable notions of fixed-point
space for a subgroup $H \subset G$.  We can form the fixed set
\begin{equation}
X^H = \{x \in X \, | \, hx = x, \forall h \in H\},
\end{equation}
or equivalently consider the space of equivariant maps $G/H \to X$.
The space $X^H$ has a natural action of the Weyl group $WH = N_G H/H$,
where here $N_G H$ denotes the normalizer of $H$ in $G$.

On the other hand, we can also form the homotopy fixed-points
\begin{equation}
X^{hH} = \Map(EG, X)^H,
\end{equation}
i.e., the $H$-equivariant maps from $EG$ to $X$, where $EG$ is a free
contractible $G$-space.  There is a natural map $X^H \to X^{hH}$
induced by the projection $EG \to \ast$, which is not usually an
equivalence.

Analogously, we can form the orbits $X_{H} = X / H$ as the quotient
and the homotopy orbits 
\begin{equation}
X_{hH} = (EG \times X) / H,
\end{equation}
which we will often write $EG \times_H X$, and is usually called the
\emph{Borel construction.}  This can also be described using the
two-sided bar construction as 
\begin{equation}
X_{hH} \cong B(G/H, G, X).
\end{equation}  
We record the following technical lemma about the behavior of the
Borel construction in the context of free actions.

\begin{lem} \label{lem:equivalence_borel_free_action} 
Let $\theta \colon H \to G$ be a homomorphism and $Z$ an $H$-space
such that $\ker \theta$ acts freely on $Z$.  Then there is a natural
equivalence $(G \sma_{\theta} Z)_{hG} \htp Z_{hH}$, i.e.,
\begin{equation}
(EG_+ \sma (G \sma_{\theta} Z)) / G \htp (EH_+ \sma Z) / H.
\end{equation}
\end{lem}

\begin{proof}
For expositional clarity, we start with the case where $\theta$ is the
inclusion of a subgroup.  Then $G \sma_{\theta} Z$ is the usual
induction functor $G \sma_H Z$, and so we have an equivalence of
$G$-spaces
\begin{equation}
(EG_+ \sma (G \sma_H Z)) \htp G \sma_H (EH_+ \sma Z).
\end{equation}
Collapsing $G$ on the right, we obtain
\begin{equation}
* \sma_H (EH_+ \sma Z) \cong (EH_+ \sma Z) / H.
\end{equation}
Now for arbitrary $\theta$, we proceed as follows.  We can rewrite the
left-hand side of the desired equivalence as
\begin{equation}
EG_+ \sma \left(G \sma_{H / \ker \theta} ((H / \ker \theta) \sma_H
X)\right),
\end{equation}
where $H/\ker \theta \to G$ is (isomorphic to) an inclusion.  This
expression is naturally equivalent to
\begin{equation}
G \sma_{H / \ker \theta} \left(E (H / \ker \theta)_+ \sma ((H / \ker
\theta) \sma_H X)\right),
\end{equation}
Finally, since $\ker \theta$ acts freely on $X$, when we collapse $G$
we can use the iterated homotopy orbit formula to obtain an expression
that is naturally weakly equivalent to $Z_{hH}$.
\end{proof}

In the context of orthogonal $G$-spectra, we can form analogues of
these constructions.  The situation with orbits is slightly simpler,
so we begin there: given a subgroup $H \subset G$, we say that a
$G$-universe $U$ is $H$-trivial if it is trivial as an
$H$-representation.
\begin{defin}
The orbit spectrum of an $X$ orthogonal $G$-spectrum indexed on an
$H$-trivial universe is defined via the levelwise formula
\begin{equation}
X/H(V) = X(V)/H.
\end{equation}
\end{defin}

When applying the orbits to a spectrum indexed on a universe $U$ which
is not $H$-trivial, we tacitly precompose with the change-of-universe
functor.  Analogously, we can consider the homotopy orbits to be
defined as the orthogonal spectrum $X_{hH} = EG_+ \sma_H X$.

For fixed-points, there are analogues of the definitions given for
spaces.  Specifically, we can define the categorical fixed-points
$X^H$ and the homotopy fixed-points $X^{hH}$.

\begin{defin}
Let $X$ be an orthogonal $G$-spectrum.  For $H \subset G$, the
categorical fixed-points $X^H$ are defined via the levelwise formula
\begin{equation}
X^H(V) = X(V)^H.
\end{equation}
\end{defin}

The categorical fixed-points can also be described as the cotensor
$F(G/H_+, X)$ or equivalently the mapping spectrum $F(\Sigma^{\infty}
G/H_+, X)$.

The categorical fixed points are lax monoidal with respect to the
external smash product of spectra.  

\begin{lem}
For an orthogonal $G_0$-spectrum $X_0$ and an orthogonal
$G_1$-spectrum $X_1$, there is a natural map of spectra
\begin{equation}\label{eq:externalsma}
  X_0^{G_0} \wedge X_1^{G_1} \to (X_0 \wedge X_1)^{G_0 \times G_1}.
\end{equation}
which is associative and unital.
\end{lem}

\begin{proof}
For a $G_0$-space $A_0$ and $G_1$-space $A_1$, there is an evident
homeomorphism of spaces 
\begin{equation}\label{eq:spacefix}
A_0^{G_0} \times A_1^{G_1} \to (A_1 \times A_2)^{G_0 \times G_1}.
\end{equation}
The required space-level map
\begin{equation}
X_0^{G_0}(V) \times X_0^{G_0}(W) \to (X_0 \sma X_1)^{G_0 \times G_1}(V
\oplus W).
\end{equation}
is specified by the homeomorphism of Equation~\eqref{eq:spacefix} and
the inclusion into the colimit computing $(X_0 \sma X_1)(V \oplus W)$.
\end{proof}

We now turn to the analog of the homotopy fixed point functor.

\begin{defin}
Let $X$ be an orthogonal $G$-spectrum.  For $H \subset G$, the
homotopy fixed points are defined as
\begin{equation}
X^{hH} = F(EG_+, X)^H.
\end{equation}
\end{defin}

One might hope that the categorical fixed-points commute with the 
suspension spectrum functor.  However, this is false, as we recall
below.  As a consequence, we define another kind of fixed-point
functor, the geometric fixed-points $\Phi^H$.

\begin{thm}[Definition V.4.3, Corollary V.4.6, Proposition V.4.7
    of~\cite{MandellMay}] For a subgroup $H \subset G$, there is a
  geometric fixed-point functor
\begin{equation}
\Phi^H \colon \Sp_G \to \Sp_{WH}
\end{equation}
with the properties that $\Phi^H$ is strong symmetric monoidal and
that $\Phi^H \Sigma^{\infty A} \cong \Sigma^{\infty} A^H$.
\end{thm}

There is a natural relationship between the point-set geometric fixed
point and categorical fixed point functors.

\begin{lem}
For any $H \subset G$, there is a natural transformation
\begin{equation}
\Phi^H \to (-)^H.
\end{equation}
\end{lem}

\subsubsection{The homotopy theory of orthogonal $G$-spectra}
\label{sec:homot-theory-orth}
In order to control the homotopy category and compute derived
functors, we shall work with the standard stable model structures on
$\Sp_G^U$.  This model structure is a stabilization of the model
structure on the category of $G$-spaces where the weak equivalences
and fibrations are detected by passage to $H$-fixed points for all $H
\subset G$.  Specifically, the equivariant stable equivalences are
detected by the equivariant stable homotopy groups, where we stabilize
along the poset of representations in $U$.

\begin{defin}
The equivariant stable homotopy groups of an orthogonal $G$-spectrum
$X$ with respect to the universe $U$ are defined for a subgroup $H
\subseteq G$ and an integer $q$ to be
\begin{equation}
\pi_q^H(X) =
\begin{cases}
\quad\displaystyle \colim_{V \subset U}
\pi_{q}((\Omega^{V}X(V))^{H})&q\geq 0\\ \quad\displaystyle
\colim_{\bR^{-q} \subset V \subset U}
\pi_{0}((\Omega^{V-\bR^{-q}}X(V))^{H})&q< 0,\\
\end{cases}
\end{equation}
\end{defin}

As usual, the stable homotopy groups can be computed in terms of the
spaces of $X$ when the adjoint structure maps are equivalences of
$G$-spaces.

\begin{defin}
An orthogonal $G$-spectrum $X$ is said to be an $\Omega$-spectrum if
the adjoint structure maps
\begin{equation}
X(V) \to \Omega^W X(V \oplus W)
\end{equation}
are weak equivalences of $G$-spaces for all $V$ and $W$.
\end{defin}

We define the equivariant stable equivalences in terms of the stable
homotopy groups.

\begin{defin}
An equivariant stable equivalence $X \to Y$ of orthogonal $G$-spectra
is a map that induces isomorphisms
\begin{equation}
\pi_q^H(X) \to \pi_q^H(Y)
\end{equation}
for every $q$ and $H$.
\end{defin}

We would now like to define the equivariant stable category as the
localization of the category $\Sp_G$ at the equivariant stable
equivalences.  In order to maintain control on this localization, we
work with a model structure on $\Sp_G$.  In fact, there are several
relevant model structures on $\Sp_G^U$, but all of the ones we work
with have the same weak equivalences; the differences amount to
different choices of resolutions for computing derived functors.

\begin{prop}[Section III.4 of \cite{MandellMay}] \label{prop:standard_model_structure_G-spectra}
There is a stable model structure on $\Sp^G_U$ in which the weak
equivalences are the equivariant stable equivalences and the fibrant
objects are the $\Omega$-spectra.
\end{prop}

This model structure is arranged to be compatible with the standard
model structure on the category $\cT_G$ where the weak equivalences
are detected on passage to fixed points for all closed subgroups.  The
following theorem is now implicit in the construction of the model
structure on $\Sp_G$; see the proof of Theorem III.4.2
of~\cite{MandellMay}.

\begin{lemma}
The adjunction $(\Sigma^{\infty}, \Ev_0)$ is a Quillen adjunction.
\end{lemma}

Furthermore, the natural isomorphism
\begin{equation}
F_0 X \sma F_0 Y \cong F_0 (X \sma Y)
\end{equation}
and the fact that the subcategory of cofibrant $G$-spaces is closed
under smash product implies that this adjunction is in fact a
symmetric monoidal Quillen adjunction.

There are associated model structures on $\Ass_G$ and $\Comm_G$; see
Sections III.7 and III.8 of~\cite{MandellMay}, substituting Lemma
B.132 of~\cite{HillHopkinsRavenel} for Lemma III.8.4
of~\cite{MandellMay}.

\begin{thm}
There is a model structure on $\Ass_G$ where the weak equivalences and
fibrations are detected by the forgetful functor to $\Sp_G$.  There is
a model structure on $\Comm_G$ where the weak equivalences are the
stable equivalences and the fibrant objects are objects which satisfy
the conditions of an $\Omega$-prespectrum except at level $0$.
\end{thm}

As a warning, note that one troublesome aspect of the stable model
structure on $\Comm_G$ described above is that the underlying
orthogonal $G$-spectra of cofibrant commutative ring spectra are
usually not cofibrant.

\subsubsection{The derived smash product and mapping spectrum in
  $\Sp_G$}

We can use the model structure on $\Sp_G$ to compute derived functors.
In general, given a functor $F \colon \Sp_G \to \aC$, it is a natural
question to ask if $F$ factors through the canonical functor $\Sp_G
\to \Ho(\Sp_G)$.  This requires that $F$ preserve stable equivalences,
which does not happen in general.  But if $F$ is a left or right
Quillen adjoint, it will preserve weak equivalences between cofibrant
or fibrant objects, respectively.  In these cases, we can derive $F$
by precomposing with the cofibrant or fibrant replacement functor.

Of particular interest are the derived smash product and derived
mapping spectrum.  If $X$ is cofibrant, then $X \sma (-)$ preserves
weak equivalences~\cite[III.7.3]{MandellMay}, and so for orthogonal
$G$-spectra $X$ and $Y$ the derived smash product $X \sma^L Y$ can be
computed by cofibrantly replacing $X$ or $Y$ and forming the point-set
smash product.  We can also derive the external smash product in the
analogous fashion.

\begin{lemma}
Let $X_0$ be a cofibrant orthogonal $G_0$-spectrum.  Then for any map
$f \colon X_1 \to X_1'$ of orthogonal $G_1$-spectra, the induced map
\begin{equation}
X_0 \sma X_1 \to X_0 \sma X_1'
\end{equation}
of orthogonal $G_0 \times G_1$-spectra is a weak equivalence.
Analogously, if $X_1$ is a cofibrant $G_1$-spectrum, then for any map
$g \colon X_0 \to X_0'$ of orthogonal $G_0$-spectra, the induced map
\begin{equation}
X_0 \sma X_1 \to X_0' \sma X_1
\end{equation}
of orthogonal $G_0 \times G_1$-spectra is a weak equivalence.
\end{lemma}

On the other hand, to derive the mapping spectrum, we use the fact
that the spectral version of Quillen's axiom SM7 for simplicial model
categories implies that for cofibrant $X$ the functor $F(X,-)$
preserves stable equivalences between fibrant objects.

\begin{defin}
For $X$ and $Y$ orthogonal $G$-spectra, the derived mapping spectrum
$RF(X,Y)$ is given by the formula
\begin{equation}
RF(X,Y) = F(X',Y'),
\end{equation}
where $X'$ is a cofibrant orthogonal $G$-spectrum stably equivalent to
$X$ and $Y'$ is a fibrant orthogonal $G$-spectrum stably equivalent to
$Y$.
\end{defin}

We will also make use to the following result that provides
homotopical control when mapping out of dualizable spectra; in our
example, we will consider the loop spectrum $F(S^V,-)$.

\begin{lem}
Let $A$ be a finite $G$-CW complex.  Then the mapping spectrum
$F(A,-)$ preserves weak equivalences of orthogonal $G$-spectra.
\end{lem}

\subsubsection{Derived fixed-point functors}

We now review how to derive the various fixed-point functors on
orthogonal $G$-spectra.  We begin by discussing the categorical and
geometric fixed points.

\begin{prop}[Proposition V.3.4 and Proposition V.4.17 of~\cite{MandellMay}]
The categorical fixed-point functor $(-)^H$ is a Quillen right
adjoint; we derive by applying the point-set functor to a fibrant
replacement of the input.  The geometric fixed-point functor $\Phi^H$
is a homotopical left adjoint can be derived by cofibrantly replacing
the input.
\end{prop}

We can detect weak equivalences with either the categorical or
geometric fixed-points.

\begin{prop}[Corollary V.4.13 of \cite{MandellMay}]
For a map $X \to Y$ of orthogonal $G$-spectra, the following
conditions are equivalent: 
\begin{enumerate}
\item the map $X \to Y$ is an equivalence.
\item the induced maps
\begin{equation}
X^H \to Y^H
\end{equation}
are (non-equivariant) weak equivalences for all $H \subset G$.
\item the induced maps
\begin{equation}
\Phi^H X \to \Phi^H Y
\end{equation}
are (non-equivariant) weak equivalences for all $H \subset G$.
\end{enumerate}
\end{prop}

As mentioned above, the interaction of the categorical fixed-point
functor with the suspension spectrum functor is complicated: a
fundamental structural result about the equivariant stable category is
the tom Dieck splitting, which characterizes this relationship.

\begin{thm}[Theorem XIX.1.3 in ~\cite{MayAlaska}]\label{thm:tD}
There is a natural isomorphism in the (non-equivariant) stable
category
\begin{equation}
(\Sigma^{\infty}_+ X)^H \cong \bigvee_{K} (E(WK)_+ \sma_{WK} X^K),
\end{equation}
where $K$ varies over the conjugacy classes of subgroups of $H$, and $WK$ is the Weyl group of $K$ in $H$ (the quotient of the normaliser by $K$). 
\end{thm}

In terms of the tom Dieck splitting (Theorem~\ref{thm:tD}), the
natural transformation $\Phi^H \to (-)^H$ is the inclusion of the
summand corresponding to $H$. 

Next, we explain how to compute the derived fixed-points of the
external smash product.  The homomorphisms
\begin{equation}
i_0 \colon G_0 \cong G_0 \times \{e\} \to G_0 \times G_1
\qquad\textrm{and}\qquad i_1 \colon G_1 \cong \{e\} \times G_1 \to G_0
\times G_1
\end{equation}
determined by the evident inclusions induce functors
\begin{equation}
i_0^* \colon \Sp^{G_0 \times G_1} \to \Sp^{G_0}
\qquad\textrm{and}\qquad i_1^* \colon \Sp^{G_0 \times G_1} \to
\Sp^{G_1}.
\end{equation}
We use the fact that orthogonal $G$-spectra are tensored over
orthogonal spectra, and let $i^*_e$ denote the forgetful functor
$\Sp^G \to \Sp$.

\begin{lemma}
Let $X_0$ be an orthogonal $G_0$-spectrum and $X_1$ an orthogonal
$G_1$-spectrum.  There are natural isomorphisms
\begin{equation}
i_0^* (X_0 \sma X_1) \cong X_0 \sma i_e^* X_1 \qquad\textrm{and}\qquad
i_1^* (X_0 \sma X_1) \cong i^*_e X_0 \sma X_1.
\end{equation} 
of orthogonal $G_0$-spectra and $G_1$-spectra respectively.
\end{lemma}

By definition, there is a natural identification
\begin{equation}
(X_0 \sma X_1)^{G_0} \cong (i_0^* (X_0 \sma X_1))^{G_0}
\end{equation}
of orthogonal spectra.  Keeping track of the $G_1$-action, we have the
following homotopical version of this isomorphism.

\begin{lemma}
Let $X_0$ be an orthogonal $G_0$-spectrum and $X_1$ an orthogonal
$G_1$-spectrum.  Then there is a natural isomorphism in the stable
category of orthogonal $G_1$-spectra
\begin{equation}
(X_0^{G_0}) \sma X_1 \to (X_0 \sma X_1)^{G_0}
\end{equation}
and a natural isomorphism in the stable category of orthogonal
$G_0$-spectra
\begin{equation}
X_0 \sma (X_1^{G_1}) \to (X_0 \sma X_1)^{G_1}.
\end{equation}
\end{lemma}

\begin{cor}
Let $X_0$ be an orthogonal $G_0$-spectrum and $X_1$ be an orthogonal
$G_1$-spectrum.  Then there are natural isomorphisms in the stable
category
\begin{equation}
((X_0 \sma X_1)^{G_0})^{G_1} \cong (X_0 \sma X_1)^{G_0 \times G_1}
\cong ((X_0 \sma X_1)^{G_1})^{G_0}.
\end{equation}
\end{cor}

\subsubsection{Borel equivariant homotopy theory}
\label{sec:borel}

The derived functors of the homotopy orbits and homotopy fixed points
of an orthogonal $G$-spectrum are defined on the ``Borel'' equivariant
stable category, which we now review.  (See~\cite[\S
  1]{BlumbergMandellKunneth} for a more expansive discussion of this
homotopy category in these terms.)

\begin{defin}
A map $X \to Y$ of orthogonal $G$-spectra (on any universe) is a Borel
equivalence if it is an underlying equivalence of spectra, i.e., if it
induces an isomorphism on the underlying non-equivariant stable
homotopy groups $\pi^{e}_k$ for all $k$.
\end{defin}

Computing the geometric fixed-points makes it clear that if a map $X
\to Y$ of orthogonal $G$-spectra is an underlying equivalence of
spectra, then $EG_+ \sma X \to EG_+ \sma Y$ is an equivariant stable
equivalence.

\begin{thm}[Theorem IV.6.3 of~\cite{MandellMay}]
Fix a complete universe $U$.  There is a model structure on $\Sp^U_G$
where the weak equivalences are the Borel equivalences and the fibrant
objects are the spectra $X$ such that the maps $X_n \to \Omega
X_{n+1}$ are underlying equivalences.
\end{thm}

The Borel equivariant homotopy theory on $\Sp_G$ can also be described
as both a localization and a colocalization of the stable model
structure (for the complete universe $U$) on $\Sp_G$ at the Borel
equivalences.  The local objects in $\Sp_G$ are those $G$-spectra $X$
such that the evident map $X \to F(EG, X^{\fib})$ is a weak
equivalence.  The colocal objects are those $G$-spectra $X$ such that
the map $EG_+ \sma X \to X$ is a weak equivalence.  Equivalently, we
can regard this as the full subcategory of orthogonal $G$-spectra on
$U$ built solely from free cells of the form $G_+ \sma \Sigma^{\infty}
D^n$.  

These characterizations show that derived mapping spectra $RF^B(X,Y)$
in the Borel category can then be described in terms of derived
mapping spectra in $\Sp^U_G$: 
\begin{equation}
\begin{aligned}
RF^B(X,Y) \cong RF(EG_+ \sma X, EG_+ \sma Y) &\cong RF(EG_+ \sma X, Y) \\
&\cong RF(X, F(EG,Y)).
\end{aligned}
\end{equation}
Choosing a cofibrant model for $EG_+$, $EG_+ \sma X$ is cofibrant in
$\Sp_G$ if $X$ is and $F(EG_+,Y)$ is fibrant if $Y$ is, so we can
compute using a cofibrant replacement of $X$ and a fibrant replacement
of $Y$ in $\Sp_G^U$.

The following calculation characterizes the Borel homotopy type of the
derived mapping spectrum in the Borel category.

\begin{lem}
Let $X$ and $Y$ be orthogonal $G$-spectra.  Then there is an
equivalence of derived mapping spectra
\begin{equation}
EG_+ \sma RF(X,Y) \cong EG_+ \sma RF(EG_+ \sma X, Y) \cong EG_+ \sma
RF^B(X,Y).
\end{equation}
\end{lem}

\begin{proof}
Let $U$ denote the complete universe.  Suppose that $X$ is cofibrant
in $\Sp_G$ with respect to $U$.  Then the map in the equivariant
stable category
\begin{equation}
EG_+ \sma RF(X,Y) \to EG_+ \sma RF(EG_+ \sma X, Y)
\end{equation}
can be computed as the point-set map
\begin{equation}
EG_+ \sma F(X, (Y)^{\fib}_U) \to EG_+ \sma F(EG_+ \sma X,
(Y)^{\fib}_U).
\end{equation}
This map is an equivalence if the map $F(X, (Y)^{\fib}_U) \to F(EG_+
\sma X, (Y)^{\fib}_U)$ is an underlying equivalence, which follows
since $EG_+ \sma X \to X$ is an underlying equivalence and the
underlying spectrum of $(Y)^{\fib}_U$ is homeomorphic to
$(Y)^{\fib}_{\bR^{\infty}}$, i.e., to the non-equivariant
fibrant-replacement functor.
\end{proof}

\subsubsection{Fibrant replacement functors on
  $\Sp_G$}\label{sec:fibrant}
\label{sec:genu-fibr-repl}

There are two different families of fibrant replacement functors for
the stable model structure on $\Sp_G$ that are relevant to our work.
Both of these functors are equipped with a natural weak equivalence
from the identity functor on $\Sp_G$.

First, for any universe $U$, the construction of the model structure
produces a functorial fibrant replacement functor $(-)^{\fib}$; this
functor arises via the functorial factorization of the terminal map $X
\to \ast$ into an acyclic cofibration followed by a fibration.  When
it is necessary to emphasize the universe, we will write
$(-)^{\fib}_U$.  Although the construction of this functor is not very
explicit, it has the attractive property (evident from the description
in terms of factorization) that if $X$ is cofibrant then $X^{\fib}$ is
cofibrant.  However, there is no reason to expect that this functor is
strictly compatible with the symmetric monoidal structure.

Work of Kro~\cite[3.2]{Kro} (generalized to the equivariant context
in~\cite[\S 19]{BlumbergMandellKunneth}) constructs for any universe
$U$ a lax monoidal fibrant replacement functor, $(-)^{\mfib}$.  This
functor is extremely useful when constructing products on mapping
spectra.

\begin{defin}\label{def:lax-fibrant}
For a universe $U$ and an orthogonal $G$-spectrum $X$, define the
orthogonal spectrum $\cQ_U X$ via the assignment
\begin{equation}
X^{\mfib}(V) = \cQ_U X(V) = \hocolim_{W \in U} \Omega^{W \otimes V}
X((W \oplus \bR) \otimes V).
\end{equation}
Here the homotopy colimit is indexed over the partially ordered set of
finite dimensional subspaces of $U$ (ordered by inclusion).  There is a natural
transformation $\id \to (-)^{\mfib}$ induced by the inclusion of $X(V)
= \Omega^{0 \otimes V} X(\bR \otimes V)$.
\end{defin}

When used without specification of the universe, we will understand
$(-)^{\mfib}$ to be constructed with respect to a complete universe
$U$.  An important caveat is that $(-)^{\mfib}$, although given by an
explicit formula, does not preserve cofibrant objects in general.
However, the natural transformation $\id \to (-)^{\mfib}$ is a
Hurewicz cofibration, as we now explain.

\begin{lemma}
Let $X$ be an orthogonal $G$-spectrum.  The natural map $X \to
X^{\mfib}$ is a Hurewicz cofibration.
\end{lemma}

\begin{proof}
Standard arguments (e.g.,~\cite[X.3.5]{ElmendorfKrizMandellMay1997})
imply that the defining inclusion is an $h$-cofibration of $G$-spaces,
and therefore $\id \to (-)^{\mfib}$ is an $h$-cofibration of
orthogonal $G$-spectra.
\end{proof}

For a $G_1$-spectrum $X$ and a $G_2$-spectrum $Y$ indexed on universes
$U_1$ and $U_2$ respectively, we have the following composite
\begin{equation}
\begin{tikzcd}
\Omega^{W_1 \otimes V_1} X((W_1 \oplus \bR) \otimes V_2) \sma
\Omega^{W_2 \otimes V_2} Y((W_2 \oplus \bR) \otimes V_2) \ar[d]
\\ \Omega^{(W_1+W_2) \otimes V_1} X(((W_1+W_2) \oplus \bR) \otimes
V_2) \sma \Omega^{(W_1+W_2) \otimes V_2} Y(((W_1+W_2) \oplus \bR)
\otimes V_2) \ar[d] \\ \Omega^{(W_1+W_2) \otimes (V_1 \oplus V_2)} (X
\sma Y)((W_1+W_2) \oplus \bR) \otimes (V_1 \oplus V_2))
\end{tikzcd}
\end{equation}

On passage to homotopy colimits, this yields the following.

\begin{prop}\label{prop:Q-lax-mon}
There is a map of orthogonal $G_1 \times G_2$-spectra
\begin{equation}
\cQ_{U_1} X \sma \cQ_{U_2} Y \to \cQ_{U_1 \oplus U_2} (X \sma Y)
\end{equation}
which is associative and unital.
\end{prop}

When $G = G_1 = G_2$ and $U = U_1 = U_2$, internalizing via Kan
extension along the direct sum yields the following result for the
smash product on $\Sp_G^U$.

\begin{prop}
For any universe $U$, the fibrant replacement functor $\cQ_{U}$ is a
lax symmetric monoidal functor on $\Sp_G^U$.
\end{prop}

We record a simple consequence of this fact.  

\begin{cor}
For any pair $A$ and $B$ of orthogonal $G$-spectra, there is a natural
map 
\begin{equation}\label{eq:Qassem}
\cQ_{U} F(A,B) \to F(A, \cQ_{U} B)
\end{equation}
that arises as the adjoint of the composite
\begin{equation}
A \sma \cQ_{U} F(A,B) \to \cQ_{U} A \sma \cQ_{U} F(A,B) \to \cQ_{U}(A
\sma F(A,B)) \to \cQ_{U} B,
\end{equation}
where the first map is induced by the unit of $\cQ_{U} A$, the second
by the lax monoidal structure map, and the third by the evaluation map
on $F(-,-)$.
\end{cor}

We now discuss functoriality and multiplicativity in the context of
surjections $p\colon G \to H$.  We begin with a general result about
the interaction of the fibrant replacement functor with arbitrary
group homomorphisms.

\begin{lem}\label{lem:restriction}
Let $p \colon G_1 \to G_2$ be a group homomorphism, $X$ a orthogonal
$G_2$-spectrum, and $U$ a $G_2$-universe.  Then there is a natural map
of orthogonal $G_1$-spectra
\begin{equation}
p^* \cQ_U X \to \cQ_{p^* U} p^* X
\end{equation}
that is a homeomorphism.  Moreover, given $p_1 \colon G_1 \to G_2$ and
$p_2 \colon G_2 \to G_3$, the diagram
\begin{equation}\label{eq:restrictcom}
\begin{tikzcd}
(p_2 \circ p_1)^* \cQ_U X = p_2^* (p_1^* \circ \cQ_U X) \ar[r] \ar[d]
  & p_2^* \cQ_{p_1^* U} p_1^* X \ar[d] \\ \cQ_{(p_2 \circ p_1)^*} (p_2
  \circ p_1)^* X \ar["=",r] & \cQ_{p_2^* p_1^* U} p_2^* p_1^* X.
\end{tikzcd}
\end{equation}
commutes.
\end{lem}

\begin{proof}
For the first part, there are homeomorphisms
\begin{equation}
p^* \hocolim_{W \in U} \Omega^{W \otimes V} X((W \oplus \bR) \otimes
V) \cong \hocolim_{W \in U} \Omega^{p^*W \otimes p^* V} X((p^* W
\oplus \bR) \otimes V),
\end{equation}
where we are assuming $V$ has trivial action.  Since the poset $\{p^*
W\}$ is cofinal in the poset of finite-dimensional subspaces of
$p^*U$, this implies the result.  For the second part, this just
amounts to the fact that $(p_2 \circ p_1)^* = p_2^* \circ p_1^*$ as
functors.  Specifically, the inclusion of $\{(p_2 \circ p_1)^* W\}$ in
$(p_2 \circ p_1)^* U_3$ coincides with the composite of the inclusion
of $p_1^* U_1$ in $U_2$ and the pullback along $p_2$.
\end{proof}

It will turn out to be very useful to restrict to specific models of
$G$-universes defined in terms of the regular representation.  We now
review some facts about the regular representations $\rho_G$ and
$\rho_H$ and the consequences for the complete $G$ and $H$ universes
formed as countable direct sums of regular representations.  (See for
example~\cite{Wimmer2019} for similar constructions.)

There is a fixed point identification $(\rho_G)^G \cong \bR$, which we
can represent as the subspace spanned by any scalar multiple of the
vector $\sum_i g_i$; it is convenient to choose a representative that
is a unit vector.  The pullback $p^* \rho_H$ is a $G$-representation,
and there is a direct sum decomposition $\rho_G \cong p^* \rho_H \oplus
(p^* \rho_H)^{\perp}$.  We can define a specific isometric embedding
$\widehat{p} \colon p^* \rho_H \to \rho_G$ using the formula
\begin{equation}
h \mapsto \frac{1}{\sqrt{\# \ker(p)}} \sum_{g \in p^{-1}(h)} g.
\end{equation}
(Here we are again choosing the scaling to ensure we have a unit
vector.)  We have a similar embedding
\begin{equation}\label{eq:regembed}
f_{G_1, G_2} \colon \rho_{G_1} \oplus \rho_{G_2} \to \rho_{G_1}
\otimes \rho_{G_2} \cong \rho_{G_1 \times G_2}
\end{equation}
defined using the linear extension of the maps $\rho_{G_1} \to
\rho_{G_1} \otimes \rho_{G_2}$ and $\rho_{G_2} \to \rho_{G_1} \otimes
\rho_{G_2}$ specified as
\begin{equation}
g_1 \mapsto \frac{1}{\sqrt{\# G_2}} \sum_{g \in G_2} g_1 \otimes g
\end{equation}
and similarly for $G_2$.  It is convenient to write
\begin{equation}
\tilde{u}_{G_2} = \frac{1}{\sqrt{\# G_2}} \sum_{g \in G_2} g
\end{equation}
to simplify this formula to $g_1 \mapsto g_1 \otimes \tilde{u}_{G_2}$,
at which point we can write Equation~\eqref{eq:regembed} as
\begin{equation}
(g_1, g_2) \mapsto g_1 \otimes \tilde{u}_{G_2} + \tilde{u}_{G_1}
  \otimes g_2.
\end{equation} 

We now record the functorial and multiplicative properties of these
maps.  First, note that $\widehat{(-)}$ is a functor; the key point
here is that for $p \colon G \to H$, $p^* \tilde{u}_{H} =
\tilde{u}_{G}$.

\begin{lem}
For surjections $p_1 \colon G \to H$ and $p_2 \colon H \to K$,
\begin{equation}
\widehat{p_2} \circ \widehat{p_1} = \widehat{p_2 \circ p_1}.
\end{equation}
\end{lem}

Next, observe that $f_{(-,-)}$ is associative.  For this, we need the
fact that under the canonical isomorphism $\rho_{G_1} \otimes
\rho_{G_2} \cong \rho_{G_1 \times G_2}$, there is an induced
identification $\tilde{u}_{G_1} \otimes \tilde{u}_{G_2} \cong
\tilde{u}_{G_1 \times G_2}$.

\begin{lem}\label{lem:mult-universe-embedding}
For groups $G_1$, $G_2$, and $G_3$, the diagram
\begin{equation}
\begin{tikzcd}
\rho_{G_1} \oplus \rho_{G_2} \oplus \rho_{G_3} \ar[r] \ar[d] &
\rho_{G_1} \oplus \rho_{G_2 \times G_3} \ar[d] \\ \rho_{G_1 \times
  G_2} \oplus \rho_{G_3} \ar[r] & \rho_{G_1 \times G_2 \times G_3}
\end{tikzcd}
\end{equation}
commutes.
\end{lem}

\begin{proof}
Expanding, we have
\begin{multline}
(g_1 \otimes \tilde{u}_{G_2} + \tilde{u}_{G_1} \otimes g_2) \otimes
  \tilde{u}_{G_3} + \tilde{u}_{G_1 \times G_2} \otimes g_3 \\ = g_1
  \otimes \tilde{u}_{G_2} \otimes \tilde{u}_{G_3} + \tilde{u}_{G_1}
  \otimes g_2 \otimes \tilde{u}_{G_3} + \tilde{u}_{G_1 \times G_2}
  \otimes g_3 \\ \cong g_1 \otimes \tilde{u}_{G_2} \otimes
  \tilde{u}_{G_3} + \tilde{u}_{G_1} \otimes g_2 \otimes
  \tilde{u}_{G_3} + \tilde{u}_{G_1} \otimes \tilde{u}_{G_2} \otimes
  g_3.
\end{multline}
And similarly for the other way around the diagram.
\end{proof}

Finally, these embeddings are compatible.

\begin{lem}\label{lem:rho-surj-compatible}
For surjections $p_1 \colon G_1 \to H_1$ and $p_2 \colon G_2 \to H_2$,
the following diagram commutes
\begin{equation}
\begin{tikzcd}
p_1^* \rho_{H_1} \oplus p_2^* \rho_{H_2} \ar[r] \ar[d] & \rho_{G_1}
\oplus \rho_{G_2} \ar[d] \\ (p_1 \oplus p_2)_* (\rho_{H_1} \oplus
\rho_{H_2}) \ar[r] \ar[d] & \rho_{G_1} \oplus \rho_{G_2} \ar[d]
\\ (p_1 \times p_2)_* (\rho_{H_1} \otimes \rho_{H_2}) \ar[r] \ar[d] &
\rho_{G_1} \otimes \rho_{G_2} \ar[d] \\ (p_1 \times p_2)_* (\rho_{H_1
  \times H_2}) \ar[r] & \rho_{G_1 \times G_2}
\end{tikzcd}
\end{equation}
\end{lem}

\begin{proof}
The commutativity of the top square follows immediately from the fact
that the kernel and the inverse image commute with direct sums.  For
the middle square, going over and down yields
\begin{equation}
\left(\frac{1}{\# \ker(p_1)} \sum_{h \in p_1^{-1}(h_1)} h\right)
\otimes \tilde{u}_{H_2} + \tilde{u}_{H_1} \otimes \left(\frac{1}{\#
  \ker(p_2)} \sum_{h \in p_2^{-1}(h_2)} h\right).
\end{equation}
On the other hand, going down and over yields
\begin{equation}
\frac{1}{(\# \ker(p_1))(\# \ker(p_2))} \left(\sum_{h \in
  p_1^{-1}(h_1)} h \otimes \tilde{u}_{G_2} + \sum_{h \in
  p_2^{-1}(h_2)} \tilde{u}_{G_1} \otimes h\right).
\end{equation}
The bottom square commutes as a consequence of the comparison of
$\tilde{u}_{G_1} \otimes \tilde{u}_{G_2}$ and $\tilde{u}_{G_1 \times
  G_2}$.
\end{proof}

The point of these observations is that the monoidal fibrant
replacement functor (when constructed using the universes that are
countable sums of the regular representation) is monoidal with respect
to the external product and contravariantly functorial in surjections
of groups.  Specifically, for any group $G$ we write $\cU(G) = \rho_G
\otimes \bR^{\infty}$ for the $G$-universe which is the infinite
direct sum of copies of the regular representation equipped with a
natural inner product.  First, note that for an orthogonal
$G$-spectrum $Y$, the natural map $i_U \colon Y \to \cQ_{\cU(G)} Y$
arising from the inclusion of the trivial vector space in $\cU(G)$ has
a preferred model using the choice of fixed-point identification
$(\rho_G)^G \cong \bR$.  Next, observe that for a surjection $G \to
J$, the isometric embedding $\widehat{p} \colon p^* \rho_J \to \rho_G$
gives rise to an isometric embedding 
\begin{equation}
\widehat{p} \colon p^* \cU(J) \to \cU(G).
\end{equation}

\begin{prop}\label{prop:surj-funct}
Given a surjection $G \to J$, there is a natural transformation
\begin{equation}
\cQ_{p^* \cU(J)} \to \cQ_{\cU(G)}
\end{equation}
which is compatible with the external lax monoidal structure on
$\cQ_{(-)}$ in the sense that for surjections $G_1 \to J_1$ and $G_2
\to J_2$, there is a commutative diagram
\begin{equation}
\begin{tikzcd}
\cQ_{p^* \cU(J_1)} \sma \cQ_{p^* \cU(J_2)} \ar[r] \ar[d] &
\cQ_{\cU(G_1)} \sma \cQ_{\cU(G_2)} \ar[d] \\
\cQ_{p^* \cU(J_1 \times J_2)} \ar[r] & \cQ_{\cU(G_1 \times G_2)}.
\end{tikzcd}
\end{equation}
\end{prop}

\begin{proof}
The existence of the transformation follows from the isometric
embedding $\widehat{p}$, which induces maps
\begin{equation}
\hocolim_{W \in U} \Omega^{W \otimes V} X((W \oplus \bR) \otimes V)
\to 
\hocolim_{\widehat{p}W \in \widehat{p} U} \Omega^{\widehat{p}W \otimes
  V} X((\widehat{p} W \oplus \bR) \otimes V).
\end{equation}
Lemma~\ref{lem:rho-surj-compatible} then implies these transformations
are externally monoidal. 
\end{proof}

Let $f \co (G, Y) \to (G', Y')$ be an arrow in the category $\Sp_{\eq,
  \surj}$ of equivariant orthogonal spectra, and recall that a
morphism consists of a surjection $p \colon G \to G'$ together with a
$G$-equivariant map from $Y \to p^* Y'$.  We consider the assignment
\begin{equation} \label{eq:functor-equivariant-fibrant-replacement}
  f \mapsto (G, \cQ_{p^* \cU(G')} Y) \in \Sp_{\eq}
\end{equation}
of an equivariant spectrum to each such arrow.

\begin{lem}\label{lem:Q-arrow}
Each factorization $f_0 = h \circ f_1 \circ g$ in $\Sp_{\eq}$
\begin{equation}
    \begin{tikzcd}
       (G_0, Y_0) \ar[r, "f_0"] \ar[d,"g"] & (G'_0, Y'_0) \\ (G_1,
      Y_1) \ar[r,"f_1"] & (G'_1, Y'_1) \ar[u,"h"]
    \end{tikzcd}
  \end{equation}
  induces a natural map
\begin{equation}
  (G_0, \cQ_{f_0^* \cU(G'_0)} Y_0) \to (G_1, \cQ_{f_1^* \cU(G'_1)}
  Y_{1}).
\end{equation}
\end{lem}

\begin{proof}
To produce the map 
\begin{equation}
\cQ_{f_0^* \cU(G'_0)} Y_0 \to g^* \cQ_{f_1^* \cU(G'_1)} Y_{1}
\end{equation}
of orthogonal $G_0$-spectra, Lemma~\ref{lem:restriction} shows that is
suffices to obtain a map
\begin{equation}
\cQ_{f_0^* \cU(G'_0)} Y_0 \to \cQ_{g^* f_1^* \cU(G'_1)} g^* Y_{1} =
\cQ_{(g \circ f_1)^* \cU(G'_1)} g^* Y_{1}.
\end{equation}
This is given as the composite
\begin{equation}
\cQ_{f_0^* \cU(G'_0)} Y_0 \to \cQ_{f_0^* \cU(G'_0)} g^* Y_1 \to 
\cQ_{(g \circ f_1)^* \cU(G'_1)} g^* Y_{1}.
\end{equation}
where the first map is induced by the map $Y_0 \to g^* Y_1$ specified
by the arrow $g$ and the second by the factorization $f_0 = h \circ (g
\circ f_1)$ and Proposition~\ref{prop:surj-funct}.
\end{proof}

This assignment yields a functor from the twisted arrow category.

\begin{prop}\label{prop:fibrant-twisted-functoriality}
The assignment of
Equation~\eqref{eq:functor-equivariant-fibrant-replacement} yields a 
monoidal functor
\begin{equation}
\Tw \Sp_{\eq, \surj} \to \Sp_{\eq, \surj}.
\end{equation}
\end{prop}

\begin{proof}
To see that this assignment is functorial in the twisted arrow
category, we have to check that the description given by
Lemma~\ref{lem:Q-arrow} is compatible with composition in $\Tw
\Sp_{\eq, \surj}$.  This follows from the fact that the pullback
$(-)^*$ is functorial and the restriction map of
Lemma~\ref{lem:restriction} is functorial (as indicated in
Equation~\eqref{eq:restrictcom}).  The fact that this is a monoidal
transformation comes from Proposition~\ref{prop:surj-funct}.
\end{proof}

We will abusively denote by $(G,Y)$ the functor $\Tw
\Sp_{\eq} \to \Sp_{\eq}$ that projects onto the domain of $f$.  Then
the natural inclusions $i_U$ induce a comparison:

\begin{cor}\label{cor:unit-fibrant-funct}
There is a monoidal natural transformation
\begin{equation}
  (G, Y) \Rightarrow (G, \cQ_{\cU(G')} Y)
\end{equation}
of functors $\Tw \Sp_{\eq} \to \Sp_{\eq}$. \qed
\end{cor}

\subsubsection{An explicit model for the homotopy cofiber}

For a pair $(X,A)$ of unbased $G$-spaces, we will find it convenient
to work with an explicit model of the homotopy cofiber.  We record the
definition and some basic properties (c.f.~\cite[\S
  III.4]{LewisMaySteinberger1986}).

\begin{defin}\label{def:mapping_cone}
Let $C(X,A)$ denote the unreduced mapping cone specified as the union
$X \cup CA$ where the basepoint of $CA$ is given by the cone point
$1$.
\end{defin}

Given a pair of unbased $G$-spaces $(X,A)$ and a $G$-space $B$, there
is a natural $G$-homeomorphism of based spaces
\begin{equation}
B_+ \sma C(X,A) \to C(X \times B, A \times B).
\end{equation}

More generally, for pairs of unbased $G$-spaces $(X,A)$ and $(Y,B)$,
we have a natural weak equivalence
\begin{equation}\label{eqn:unbased-cone-prod}
C(X,A) \sma C(Y,B) \to C(X \times Y, (A \times Y) \cup (X \times B))
\end{equation}
specified by the formulas
\begin{equation}
\begin{cases}
(x,y) \qquad x \in X, y \in Y \\ ((x,b), t) \qquad x \in X, (b,t) \in
  CB \\ ((a,y), s) \qquad y \in Y, (a,s) \in CA \\ ((a,b), \max(s,t))
  \qquad (a,s) \in CA, (b,t) \in CB.\\
\end{cases}
\end{equation}

These maps are associative.

\begin{prop}
For pairs of unbased $G$-spaces $(X,A)$, $(Y,B)$, and $(Z,D)$, the
following diagram commutes:
\begin{equation}
\begin{tikzpicture}
\node[scale=.75]{
\begin{tikzcd}
C(X,A) \sma C(Y,B) \sma C(Z,D) \ar[r] \ar[d] & C(X \times Y, (A \times
Y) \cup (X \times B)) \sma C(Z,D) \ar[d] \\ C(X,A) \sma C(Y \times Z,
(B \times Z) \cup (Y \times A)) \ar[r] & C(X \times Y \times Z, (A
\times Y \times Z) \cup (X \times B \times Z) \cup (X \times Y \times
D)).\\
\end{tikzcd}
};
\end{tikzpicture}
\end{equation}
\end{prop}

\begin{proof}
This follows by tracing elements around the diagram using the explicit
formulas for the maps.  For example, given $x \in X$, $y \in Y$, and
$(d,t) \in CD$, the over and down composite goes to $((x,y), (d,t))$
and then to $((x,y,d),t)$.  The down and over composite goes to $(x,
((y,d), t))$ and then to $((x,y,d), t)$.  The other checks are
analogous, using the fact that $\max(s,\max(t,v)) = \max(\max(s,t), v)
= \max(s,t,v)$.
\end{proof}

When specializing Equation~\eqref{eqn:unbased-cone-prod} to the case
of $C(X,X-A)$ and $C(Y,Y-B)$ for $A \subset X$ and $B \subset Y$, we
obtain the following product map
\begin{equation}\label{eqn:cone-difference-prod}
\begin{aligned}
C(X,X-A) \sma C(Y,Y-B) &\to C(X \times Y, ((X-A) \times Y) \cup (X
\times (Y-A))) \\ 
&= C(X \times Y, (X \times Y) - (A \times B)),
\end{aligned}
\end{equation}
which we will make frequent use of.

We are particularly interested in the special case of $C(S^V,
S^V-{0})$ for a finite-dimensional real vector space $V$.  Then the
map in Equation~\eqref{eqn:unbased-cone-prod} produces a natural weak
equivalence
\begin{equation}\label{eqn:cone-sphere-smash}
C(S^V,S^V-\{0\}) \sma C(S^W,S^W-\{0\}) \to C(S^{V \oplus W}, S^{V
  \oplus W}-\{0\}),
\end{equation}
by composing the map
\begin{equation}
\begin{tikzcd}
C(S^V,S^V-\{0\}) \sma C(S^W,S^W-\{0\}) \ar[d] \\ C(S^V \times
S^W,(S^V-\{0\} \times S^W) \cup (S^V \times S^W-\{0\}))
\end{tikzcd}
\end{equation}
with the map induced by the basepoint collapse map $S^V \times S^W
\to S^{V \oplus W}$ and the map of cones
\begin{equation}
C\left( (S^V-\{0\} \times S^W) \cup (S^V \times S^W-\{0\})\right) \to
C\left(S^{V \oplus W} - \{0\}\right)
\end{equation}
induced by the inclusion and basepoint collapse.  It is
straightforward to check that these maps are also associative, since
the basepoint collapse is obviously associative and the cone maps
commute.  Summarizing, we have the following:
\begin{prop}
For vector spaces $U$, $V$, and $W$, the following diagram commutes:
\begin{equation}
\begin{tikzpicture}
\node[scale=.75]{
\begin{tikzcd}
C(S^U,S^U-\{0\}) \sma C(S^V,S^V-\{0\}) \sma C(S^W,W-\{0\}) \ar[r]
\ar[d] & C(S^{U\oplus V}, S^{U \oplus V} - \{0\}) \sma C(S^W,W-\{0\})
\ar[d] \\ C(S^U,U-\{0\}) \sma C(S^{V\oplus W}, S^{V \oplus W} - \{0\})
\ar[r] & C(S^{U \oplus V \oplus W}, S^{U \oplus V \oplus W} -
\{0\}).\\
\end{tikzcd}
};
\end{tikzpicture}
\end{equation}
\end{prop}

\subsection{Graded spectra, filtered spectra, and completions}\label{sec:graded}

In this section, we review the theory of filtered and graded spectra
and completions from this perspective.  For more discussion, we
recommend the treatment in~\cite{GwilliamPavlov2018} (and see
also~\cite[\S 2.4, \S 3]{Lurie2014a}, the quick review in ~\cite[\S
  5.1]{BhattMorrowScholze2019}, and~\cite[\S
  6]{BlumbergMandellKunneth}).

\subsubsection{Filtered and graded spectra}

Let $\bZ$ denote the category associated to the poset $(\bZ, \leq)$.

\begin{defin}
A filtered spectrum is a functor $\bZ^{\op} \to \Sp$.  Explicitly, we
can write this as a sequence
\begin{equation}
X_{\bullet} = \ldots \to X^n \to X^{n-1} \to \ldots
\end{equation}
of spectra.  We will denote the category of filtered spectra as $\Filt(\Sp)$.
\end{defin}

We regard $\Filt(\Sp)$ as giving models of decreasing filtrations, and
so we think of 
\begin{equation}
X_{-\infty} = \colim_n X_{\bullet}
\end{equation}
as the ``underlying object'' of the filtration.

Let $\bZ_{\disc}$ to denote the category with objects the elements of
$\bZ$ and only the identity morphisms.

\begin{defin}
A graded spectrum is a functor $\bZ_{\disc} \to \Sp$.  We will write
$\Gr(\Sp)$ to denote the category of graded spectra.
\end{defin}

The category of graded spectra is equipped with a functor to spectra
specified by the formula
\begin{equation}
X_\bullet \mapsto \bigvee_n X_n.
\end{equation}
We refer to this as the ``underlying object'' of the graded spectrum.

\begin{defin}
There is an associated graded functor
\begin{equation}
\Gr \colon \Filt(\Sp) \to \Gr(\Sp)
\end{equation}
specified on objects by the formula
\begin{equation}
\Gr(X_\bullet)_n = X_{n-1}/X_n.
\end{equation}
\end{defin}

Since $\bZ^{\op}$ and $\bZ_{\disc}$ are symmetric monoidal categories
with operation $+$ and unit $0$, the category $\Filt(\Sp)$ and
$\Gr(\Sp)$ are endowed with symmetric monoidal structures by the Day
convolution.

\begin{prop}
The categories $\Filt(\Sp)$ and $\Gr(\Sp)$ are symmetric monoidal under
the Day convolution.  The unit of $\Filt(\Sp)$ is given by the filtered
spectrum that is $\bS$ for $n \leq 0$ and $\ast$ otherwise (with all
maps the $\id$ or the trivial map).  The unit in $\Gr(\Sp)$ is the
graded spectrum with value $\bS$ at $0$ and $\ast$ everywhere else.
\end{prop}

The symmetric monoidal product on $\Filt(\Sp)$ can be written
\begin{equation}
(X_\bullet \sma Y_\bullet)_n = \colim_{p+q \geq n} X_p \sma Y_q
\end{equation}
and on $\Gr(\Sp)$ as
\begin{equation}
(X_\bullet \sma Y_\bullet)_n = \coprod_{p+q = n} X_p \sma Y_q.
\end{equation}
For example, a monoid in $\Sp^{\gr}$ is specified by associative
product maps 
\begin{equation}
\bigvee_{i + j = k} X_i \sma X_j \to X_k.
\end{equation}
In particular, for each $(i, j)$ such that $i + j = k$, we have maps 
\begin{equation}
X_i \sma X_j \to X_k
\end{equation}
which are associative and unital in the evident fashion.

The Day convolution product on $\Filt(\Sp)$ encodes a filtered
multiplication on $X_{-\infty}$.

\begin{lem}
The underlying spectrum functor
\begin{equation}
\colim_n \colon \Filt(\Sp) \to \Sp
\end{equation}
and the associated graded functor
\begin{equation}
\Gr \colon \Filt(\Sp) \to \Gr(\Sp)
\end{equation}
are strong symmetric monoidal.
\end{lem}

We now turn to the homotopy theory of filtered and graded spectra.
The category $\Filt(\Sp)$ can be equipped with the projective model
structure.

\begin{prop}
There is a model structure on $\Filt(\Sp)$ in which the fibrations and
weak equivalences are detected levelwise.
\end{prop}

Cofibrant objects in $\Filt(\Sp)$ in particular have the property that
the structure maps are cofibrations, and so we can derive the
underlying spectrum functor to obtain the homotopy colimit. 

The category $\Gr(\Sp)$ also admits a projective model structure,
which in this case is just the product model structure since there are
no nontrivial morphisms in $\bZ_{\disc}$.

\begin{prop}
There is a model structure on $\Gr(\Sp)$ in which the cofibrations,
fibrations, and weak equivalences are detected levelwise.
\end{prop}

The associated graded functor are evidently homotopical.

\begin{lem}
The functor $\Gr$ is the left adjoint of a Quillen adjunction.
\end{lem}

We use this to define a coarser notion of weak equivalence on
$\Filt(\Sp)$ that will be useful when we study completions.

\begin{defin}
A map $f \colon X_\bullet \to Y_\bullet$ of filtered spectra is a
graded (or filtered) equivalence if the induced map
\begin{equation}
\Gr(X_\bullet) \to \Gr(Y_\bullet)
\end{equation}
is an equivalence.
\end{defin}

Moreover, these products are compatible with the model structures, in
the following sense.

\begin{prop}[Theorem 3.50 of~\cite{GwilliamPavlov2018}]
The projective model structures on $\Filt(\Sp)$ and $\Gr(\Sp)$ are
monoidal model categories for the Day convolution product and $\Gr$ is
the left adjoint in a monoidal Quillen adjunction.
\end{prop}

\subsubsection{The completion of a filtered spectrum}\label{appendix:complete}

Let $X_\bullet$ be a filtered spectrum.  We can form the completion of
$X_\bullet$ as follows.  For each $n$, the commutative diagram 
\begin{equation}
\begin{tikzcd}
X_n \ar{r} \ar{d} & X_{-\infty} \\ X_{n-1} \ar{ur}
\end{tikzcd}
\end{equation}
implies that there are induced maps $X_{-\infty}/X_{n} \to
X_{-\infty}/X_{n-1}$.

\begin{defin}[see Section 3.10 of~\cite{GwilliamPavlov2018}]
The completion (as a spectrum) is defined to be
\begin{equation}
\widehat{X} = \holim_n \left( X_{-\infty} / X_n\right).
\end{equation}
The completion is itself equipped with a decreasing filtration; we
define a filtered spectrum
\begin{equation}
\widehat{X}_n = \hofib(\widehat{X} \to X_{-\infty} / X_n).
\end{equation}
\end{defin}

There is a natural map $X_n \to \widehat{X}_n$ for every $n$, and
these assemble into a natural transformation $\id \to \widehat{(-)}$.
Essentially by construction, this natural transformation induces an
equivalence on associated graded spectra.

\begin{lem}\label{lem:comparegraded}
For any $X_\bullet$, the natural map $X_\bullet \to
\widehat{X}_\bullet$ is a graded equivalence.
\end{lem}

As a consequence of Lemma~\ref{lem:comparegraded}, the completion can
be abstractly described by inverting the graded equivalences.

\begin{defin}[Definition 3.5 of~\cite{GwilliamPavlov2018}]\label{defn:completion}
Let $\Comp(\Filt(\Sp))$ denote the model category structure on
$\Filt(\Sp)$ obtained by the left Bousfield localization at the graded
equivalences.
\end{defin}

The formal model of completion given by localization coincides with
the explicit construction given above.

\begin{prop}[Proposition 3.31 of~\cite{GwilliamPavlov2018}]
The completion functor is naturally equivalent to the localization
functor $\Filt(\Sp) \to \Comp(\Filt(\Sp))$.
\end{prop}

The symmetric monoidal structure on $\Filt(\Sp)$ induces one on
$\Comp(\Filt(\Sp))$.  

\begin{thm}[Theorem 3.50 of~\cite{GwilliamPavlov2018}]\label{thm:compmonoid}
  The Day convolution product equips
$\Comp(\Filt(\Sp))$ with the structure of a symmetric monoidal model
category and the localization functor is a lax symmetric monoidal
Quillen functor.
\end{thm}

Explicitly, we can compute the completed smash product on $\Filt(\Sp)$
as 
\begin{equation}
X_\bullet \widehat{\sma} Y_\bullet = \widehat{(X_\bullet \sma
  Y_\bullet)}.
\end{equation}

As a consequence of Theorem~\ref{thm:compmonoid}, the completion
preserves multiplicative structures:
\begin{cor}
Let $R_\bullet$ be an associative or commutative ring object in
$\Filt(\Sp)$.  Then $\widehat{R}_\bullet$ is an associative or
commutative ring object, respectively, in $\Comp(\Filt(\Sp))$.
\end{cor}

In particular, we have the following useful corollary.

\begin{cor}
Let $R_\bullet$ be an associative ring object in $\Filt(\Sp)$.  Then
the underlying spectrum of the completion $\widehat{R} = \colim_n
\widehat{R}_\bullet$ is an associative ring orthogonal spectrum.
\end{cor}

In the previous discussion, we have focused on spectra with a
filtration indexed by $\bZ$.  However, it is often natural to consider
other directed systems, notably $\bR$.  In this context, a filtered
spectrum is a functor $\bR^{\op} \to \Sp$, which again we think of as
a decreasing filtration on $X_{-\infty} = \colim_n X_{\bullet}$.
However, the notion of associated graded is more complicated, as it
should be described in terms of the quotients
\begin{equation}
\Gr(X_\bullet)_n = X_n / \colim_{m > n} X_m.
\end{equation}
We can nonetheless define the completion of an $\bR$-filtered spectrum
as the filtered spectrum
\begin{equation}
\widehat{X}_n = \hofib(\widehat{X} \to X_{-\infty} / X_n).
\end{equation}
having underlying spectrum 
\begin{equation}
\widehat{X} = \holim_n \left( X_{-\infty} / X_n\right).
\end{equation}

However, in order to avoid dealing with some of the technical
complexities that arise in this setting, in the body of the paper we
will always work with filtered spectra induced by composites
\begin{equation}
\bZ^{\op} \to \bR^{\op} \to \Sp
\end{equation}
obtained by choosing cofinal indexing sets in $\bR^{\op}$.

\subsubsection{Graded spectra and the $2$-periodic sphere
  spectrum}\label{sec:action}

In this section, we explain how to obtain a coherent multiplicative
system of positive and negative spheres.  We do this by working with a
particular model of the $2$-periodic sphere spectrum, following Lurie~\cite[\S 3]{Lurie2014a}.
Conceptually, the best description is as the Thom spectrum of the
$E_2$ map
\begin{equation}
\bZ \htp \Omega^2 BU(1) \to \Omega^2 BU \to BU \times \bZ.
\end{equation}
Here the last map is specified by Bott periodicity; it takes some work
to show that this map is $E_2$.  For our purposes, we will in fact use
an intermediate form of the construction in the category of graded
spectra.  

Lurie constructs a graded $E_2$ ring spectrum (by which we mean an
$E_2$ algebra in the category of graded spectra) denoted
$\bS[\beta^{\pm}]$ with the property that $(\bS[\beta^{\pm}])_k \htp
\bS^{-2k}$~\cite[\S 3.4]{Lurie2014a}; the underlying $E_2$ ring
spectrum of $\bS[\beta^{\pm}]$ is $\bigvee_{n \in \bZ}
\bS^{-2n}$~\cite[\S 3.5.13]{Lurie2014a}.  We can use the constituent
spectra of $\bS[\beta^{\pm}]$ as our coherent family of (even)
spheres, as follows:
\begin{prop}\label{prop:luriesphere}
There is a collection of spectra $\{\bS[n]\}_{n \in \bZ}$ with the
following properties:
\begin{enumerate}
\item For $n \in \bZ$, $\bS[n]$ is cofibrant as a spectrum,
\item For $n \in \bZ$, $\bS[n] \htp S^{-2n}$,
\item For $n,m \in \bZ$, there are strictly associative and unital
  maps
\[
\bS[n] \sma \bS[m] \to \bS[n+m],
\]
which are models of the standard equivalences
\[
S^{-2n} \sma S^{-2m} \to S^{-2(n+m)}.
\]
\end{enumerate}
\end{prop}
\begin{proof}
Lurie constructs $\bS[\beta^{\pm}]$ as a graded $E_2$ ring spectrum in
the $\infty$-category of graded spectra.
Using~\cite[1.1]{NikolausSagave2017}, we rectify $\bS[\beta^{\pm}]$ to
a strictly associative graded ring spectrum.
\end{proof}

In fact, we can say something slightly stronger about the comparison
to the standard spheres; the following proposition is a consequence
of~\cite[3.4.5]{Lurie2014a}, and says that the multiplication maps on
the $\{\bS[n]\}$ are coherently compatible with the standard
multiplication of spheres.

\begin{prop}\label{prop:lurie-comparison}
There is a zig-zag of equivalences as graded associative ring spectra
between the free associative ring on $\bS^{-2}(-1)$ (where this
denotes the graded ring spectrum $X$ such that $X_{-1} = \bS^{-2}$)
and the negative truncation of $\bS[\beta^{\pm}]$.  Analogously, there
is a zig-zag of equivalences as graded associative ring spectra
between the free associative ring on $\bS^2(1)$ and the positive
truncation of $\bS[\beta^{\pm}]$. \qed
\end{prop}

\begin{rem}
The residual $E_2$ structure on the system of spheres $\{\bS[n]\}$
gives rise to a graded commutative structure when we pass to homotopy
groups.  Nonetheless, a natural question that arises (although the
answer is not germane to our enterprise) is whether we can do better
before passage to homotopy groups.  However, a calculation with power
operations shows that in fact there cannot be an $E_3$ structure.
(See~\cite[VII.6.1]{BrunerMayMcClureSteinberger1986} for an early
example of this kind of argument in the context of $H_\infty$ ring
spectra.)  In contrast, there are $E_\infty$ structures on the
analogous constructions of periodic cobordism and periodic $H\bZ$.
\end{rem}

We now turn to a generalization of the notion of graded spectrum
specified above that we use to describe our version of the spectral
twisted Novikov ring.  Let $\Sigma$ be a discrete monoid (which we
assume is countable), potentially non-unital.  In the following, we
will regard $\Sigma$ as a discrete monoidal category with object set
$\Sigma$, monoidal product specified by the product on $\Sigma$, and
the only morphisms the identity.  (This monoidal category is
non-unital in the case that $\Sigma$ is non-unital.)

\begin{defin}
The category $\Sp^{\Sigma,\gr}$ of $\Sigma$-graded spectra consists of
functors $\Sigma \to \Sp$; a graded spectrum $X_\bullet$ is specified
by a collection of spectra $\{X_\sigma\}$ for $\sigma \in \Sigma$, and
maps $X_\bullet \to Y_\bullet$ are given levelwise.
\end{defin}

Once again, the Day convolution endows $\Sp^{\Sigma,\gr}$ with a
monoidal structure.  If $\Sigma$ is unital, the Day convolution
monoidal structure has unit the $\Sigma$-graded spectrum that is $\bS$
at $0$ and $\ast$ elsewhere.  An explicit formula for the product is
given by the expression
\begin{equation}
(X_\bullet \sma Y_\bullet)_\sigma = \bigvee_{\sigma_1 \sigma_2 =
    \sigma} X_{\sigma_1} \sma Y_{\sigma_2}.
\end{equation}
Moreover, there is again a lax monoidal ``underlying spectrum''
functor $\Sp^{\Sigma,\gr} \to \Sp$ specified on objects by the
assignment
\begin{equation}
X_\bullet \to \bigvee_{\sigma \in \Sigma} X_\sigma.
\end{equation}

Now, assume that we have a degree homomorphism $\deg \colon \Sigma \to
\bZ$ and fix an associative ring orthogonal spectrum $\bbk$.

\begin{lemma}\label{lem:gradedtwisted}
Consider the $\Sigma$-graded spectrum specified by the formula
\begin{equation}
X_\sigma := \bS[-\deg(\sigma)] \sma \bbk.
\end{equation}
The product maps
\begin{multline*}
(\bS[-\deg(\sigma_1)] \sma \bbk) \sma (\bS[-\deg(\sigma_2)] \sma \bbk)
  \to \\ (\bS[-\deg(\sigma_1)] \sma \bS[-\deg(\sigma_2)]) \sma (\bbk
  \sma \bbk) \to \\ \bS[-\deg(\sigma_1)-\deg(\sigma_2)] \sma \bbk
  \cong \bS[-\deg(\sigma_1 \sigma_2)] \sma \bbk
\end{multline*}
make this into a monoid object in $\Sigma$-graded spectra.
\end{lemma}

\begin{defin}
We let $\Sigma^{\deg} \bbk[\Sigma]$ denote the underlying spectrum of
$X_\sigma$.  This is an associative ring orthogonal spectrum with
homotopy type
\[
\Sigma^{\deg} \bbk[\Sigma] \htp \bigvee_{\sigma \in \Sigma}
\bS[-\deg(\sigma)] \sma \bbk.
\]
(When $\Sigma$ is non-unital, so is the underlying ring spectrum.)
\end{defin}

\subsection{Homotopical algebra of enriched categories}
\label{sec:homotopical-algebra}

\subsubsection{Enriched categories}

Our definitions of flow category and virtual fundamental chain rely on
the notions of categories enriched in spaces and spectra, i.e.,
categories $\aC$ such that for every pair of objects $x,y \in
\ob(\aC)$ there is a mapping space or spectrum $\aC(x,y)$.
See~\cite{Kelly1982} for a comprehensive introduction to enriched
category theory.  In what follows, let $V$ be a symmetric monoidal
category with product $\boxtimes$ and unit $\mathbb{1}$; in our
applications $V$ will either be topological spaces under the cartesian
product, pointed topological spaces under the smash product, or
spectra under the smash product.

\begin{defin}
A $V$-enriched category $\aC$ is specified by a class of objects
$\ob(\aC)$ and for every pair of objects $x,y \in \ob(\aC)$, an object
$\aC(x,y)$ of $V$ satisfying the following conditions:
\begin{enumerate}
\item For every triple of objects $x,y,z \in \ob(\aC)$, there are
  composition maps
\begin{equation}
\aC(x,y) \boxtimes \aC(y,z) \to \aC(x,z).
\end{equation}
\item For every object $x \in \ob(\aC)$, there is a distinguished unit
  morphism $\mathbb{1} \to \aC(x,x)$.
\item The composition maps are associative and unital in the evident
  sense.
\end{enumerate}
We will say that an enriched category is {\em small} if it has a set
of objects.
\end{defin}

Associated to a $V$-enriched category $\aC$ we can extract an
underlying ordinary category.

\begin{defin}
Let $\aC$ be a $V$-enriched category.  The ordinary category
underlying $\aC$ has morphism sets specified as $\Map_V(\mathbb{1},
\aC(x,y))$.
\end{defin}

If $\aC$ and $\aD$ are $V$-enriched categories, an enriched functor $F
\colon \aC \to \aD$ is specified by: 
\begin{enumerate}
\item A function $\ob(\aC) \to \ob(\aD)$, and 
\item morphisms $\aC(x,y) \to \aD(Fx, Fy)$ in $V$ that are
  compatible with the unit and the composition.  
\end{enumerate}
Enriched natural transformations are defined analogously.

We will often use the following construction on enriched categories.

\begin{defin}\label{defn:prodenrich}
Let $\aC$ and $\aD$ be $V$-enriched categories.  The $V$-enriched
category $\aC \boxtimes \aD$ is defined to have objects $\ob(\aC)
\times \ob(\aD)$ and morphism objects given by the formula
\begin{equation}
(\aC \boxtimes \aD)((x,y), (x',y')) = \aC(x,x') \boxtimes \aD(y,y').
\end{equation}
\end{defin}

\subsubsection{Spectral categories}\label{sec:spectral}

We now specialize to a discussion of the theory of small categories
enriched in orthogonal spectra, which we refer to as spectral
categories.  Let $\Spcat$ denote the category of small spectral
categories and {\em spectral functors}.  

The category $\Spcat$ is symmetric monoidal, where the product is
given by the smash product of spectral categories (as in
Definition~\ref{defn:prodenrich}) and the unit is the spectral
category $\bS$ with a single object $x$ and morphism spectrum $\bS$.  

It is often useful to think of a spectral category as a ``ring with
many objects''; given a spectral category $\aC$ with a single object
$x$, the mapping spectrum $\aC(x,x)$ is an associative ring orthogonal
spectrum.  We define modules over a spectral category $\aC$ as
follows.

\begin{defin}
Let $\aC$ be a spectral category.  The category $\aC$-$\mod$ of left
$\aC$-module has:
\begin{enumerate}
\item objects the spectral functor $\aC \to \Sp$, and
\item morphisms the natural transformations.
\end{enumerate}
Analogously, the category $\mod$-$\aC$ of right $\aC$-module has
objects the spectral functors $\aC^{\op} \to \Sp$.  

Given spectral categories $\aC$ and $\aD$, a $(\aC,\aD)$-bimodule is a
left module over the spectral category $\aC \sma \aD^{\op}$.
\end{defin}

A spectral functor $F \colon \aC \to \aD$ induces an adjunction on
module categories.

\begin{prop}\label{prop:basechange}
Let $F \colon \aC \to \aD$ be a spectral functor.  There is an
adjoint pair 
\begin{equation}
\begin{tikzcd}
F_!  \colon \aC\textrm{-}\mod \ar[r, shift left=1ex] & \ar[l, shift
  left=1ex] \aD\textrm{-}\mod \colon F^*. 
\end{tikzcd}
\end{equation}
where $F^*$ is the pullback and $F_!$ is the enriched left Kan
extension.  (And analogously for the categories of right modules.)
\end{prop}

We now turn to discussion of the homotopy theory of spectral
categories and their modules.  Associated to any spectral category
$\aC$, we can form the homotopy category $\Ho(\aC)$.

\begin{defin}
For a spectral category $\aC$, let $\Ho(\aC)$ be the ordinary category 
with the same objects as $\aC$ and morphism sets specified as
\begin{equation}
\Ho(\aC)(x,y) = \pi_0 \aC(x,y).
\end{equation}
\end{defin}

This gives rise to the following homotopical notion of equivalence of
spectral categories:

\begin{defin}\label{defn:DK-equivalence}
A spectral functor $F \colon \aC \to \aD$ is a {\em DK-equivalence} of
spectral categories if
\begin{enumerate}
\item for each pair of objects $x,y \in \ob(\aC)$, the induced map of
  spectra $\aC(x,y) \to \aD(Fx,Fy)$ is a stable equivalence, and
\item the induced functor $\Ho(\aC) \to \Ho(\aD)$ is an equivalence of
  categories.
\end{enumerate}
\end{defin}

We work with spectral categories up to DK-equivalence; to express
this, it is useful to construct a model structure on $\Spcat$ in which
the DK-equivalences are the weak equivalences.  In the remainder of the section, we assume we
fixed a chosen model structure on orthogonal spectra (e.g., the stable
model structure described in Proposition
\ref{prop:standard_model_structure_G-spectra}).

\begin{thm}\label{thm:model}
There is a model structure on $\Spcat$ in which:
\begin{enumerate}
\item the weak equivalences are the DK-equivalences, and 
\item the fibrations are the functors such that the maps $\aC(x,y) \to
  \aD(Fx,Fy)$ are fibrations for every $x,y \in \ob(\aC)$.
\end{enumerate}
\end{thm}

\begin{proof}
When considering the category $\Spcat^{\Delta}$ of small spectral
categories enriched over symmetric spectra, this result is proved
in~\cite{Tabuada2009}.  There is an adjunction $(F,U)$ between
symmetric spectra and orthogonal spectra, where $F$ is given by left
Kan extension and $U$ is the restriction.  This adjunction is
monoidal, in the sense that $F$ is strong symmetric monoidal and $U$
is lax monoidal.  As a consequence, there is an induced adjunction
$(F,U)$ between $\Spcat^{\Delta}$ and $\Spcat$.  The model structure
on $\Spcat$ can be constructed as the transferred model structure
along the functor $U \colon \Spcat \to \Spcat^{\Delta}$; we check the
axioms for~\cite[1.10]{BergerMoerdijk2013}, all of which are clear
except for the existence of a generating set of intervals.  But that
follows by applying the functor $F$ to the corresponding set for
$\Spcat^{\Delta}$.
\end{proof}

Unfortunately, the model structure of the preceding theorem is not
compatible with the symmetric monoidal structure; it is not the case that if
$\aC$ is cofibrant then $\aC \sma (-)$ is a left Quillen functor.
Nonetheless, we can form the derived smash product of spectral
categories, using a theory of flat objects.

\begin{defin}
We say a small spectral category $\aC$ is {\em pointwise cofibrant} if
for all pairs of objects $x,y$ in $\aC$, the mapping spectrum
$\aC(x,y)$ is a cofibrant orthogonal spectrum.  We say $\aC$ is {\em
  pointwise fibrant} if each mapping spectrum $\aC(x,y)$ is a fibrant
orthogonal spectrum.
\end{defin}

We now have the following key result.

\begin{prop}\label{prop:derivedsmash}
Let $\aC$ be a pointwise cofibrant spectral category.  Then the
functor
\begin{equation}
\aC \sma (-) \colon \Spcat \to \Spcat
\end{equation}
preserves DK-equivalences.
\end{prop}

Since fibrant objects in Theorem~\ref{thm:model} are precisely the
pointwise fibrant objects and the cofibrant objects are in particular
pointwise cofibrant~\cite[3.3]{Tabuada2009} and this model structure
has functorial factorization, we can always arrange for these
properties to hold.  In fact, an even easier approach is to work with
the model structure on the category of small spectral categories with
fixed object set; see~\cite[\S 2]{BlumbergMandell2012} for discussion
of this point, following~\cite{SchwedeShipley2003}.

\begin{prop}\label{prop:specapprox}
There are functors
\begin{equation}
(-)^{f} \colon \Spcat \to \Spcat \qquad\textrm{and}\qquad (-)^{c}
  \colon \Spcat \to \Spcat
\end{equation}
such that $\aC^{f}$ is pointwise fibrant for all $\aC$, $\aC^c$ is
pointwise cofibrant for all $\aC$, and there are natural
transformations through $DK$-equivalences 
\begin{equation}
\id \to (-)^f \qquad\textrm{and}\qquad (-)^c \to \id
\end{equation}
that are the identity on object sets.
\end{prop}

As a consequence of Proposition~\ref{prop:derivedsmash}, given
spectral categories $\aC$ and $\aD$ we write $\aC \sma^{L} \aD$ to
denote $\aC' \sma \aD$, where $\aC' \to \aC$ is a pointwise cofibrant
replacement.  We refer to $\aC \sma^L \aD$ as the {\em derived smash
  product} of the spectral categories $\aC$ and $\aD$.

There are also (much simpler) model structures on the categories
$\aC$-$\mod$ and $\mod$-$\aC$; as with all presheaf categories, these
are lifted directly from the model structure on $\Sp$.

\begin{thm}[Theorem 6.1 of~\cite{SchwedeShipley2003b}]\label{thm:modulemodel}
There are model structures on $\aC$-$\mod$ and $\mod$-$\aC$ in which
the fibrations and weak equivalences are detected pointwise.
\end{thm}

Inspection of the cofibrant objects provides the following useful
corollary.

\begin{cor}
When $\aC$ is pointwise cofibrant, a cofibrant $\aC$-module is itself
pointwise cofibrant; in this case, there is a pointwise cofibrant
replacement functor on $\aC$-modules.
\end{cor}

The model structure of Theorem~\ref{thm:modulemodel} is compatible
with the adjunction of Proposition~\ref{prop:basechange}, since
fibrations and weak equivalences are clearly preserved by pullback
along a spectral functor $\aC \to \aD$.

\begin{prop}[Proposition 3.2 of~\cite{Toen2007}]
Let $F \colon \aC \to \aD$ be a spectral functor.  Then the adjunction
$(F_!, F^*)$ is a Quillen adjunction.  When $F$ is a DK-equivalence,
the adjunction is a Quillen equivalence.
\end{prop}

\subsubsection{Spectral presheaves}

Let $\aD$ be a small symmetric monoidal topologically enriched
category.  In many of our examples, $\aD$ will in fact have discrete
mapping spaces. Using the canonical topological enrichment on the
category of orthogonal spectra, we have:

\begin{defin}
The category $\Pre(\aD,\Sp)$ of spectral presheaves has objects the
functors $\aD \to \Sp$ and morphisms the natural transformations.
\end{defin}

Since $\aD$ is symmetric monoidal, Day convolution endows the category
$\Pre(\aD,\Sp)$ with the structure of a symmetric monoidal category,
where the unit is the presheaf represented by $\aD(\mathbb{1}, -)_+
\sma \bS$.  The category of (commutative) monoids in $\Pre(\aD,\Sp)$
is precisely the category of lax (symmetric) monoidal functors $\aD
\to \Sp$.

\begin{thm}\label{thm:proj-model}
The category $\Pre(\aD, \Sp)$ admits a projective model category
structures where the weak equivalences and fibrations are defined
pointwise, lifted from the stable and positive stable model structures
on $\Sp$.  When working with the positive stable model structure, the
projective model structure makes $\Pre(\aD,\Sp)$ into a symmetric
monoidal model category.
\end{thm}

The symmetric monoidal model structure on $\Pre(\aD, \Sp)$ can be
lifted to categories of monoids and commutative monoids.

\begin{thm}\label{thm:lax-monoidal-model}
There are cofibrantly generated model structures on the categories of
monoids and commutative monoids in $\Pre(\aD, \Sp)$ where the
fibrations and weak equivalences are determined by the forgetful
functor to $\Pre(\aD, \Sp)$ (with the positive stable model
structure).
\end{thm}

For us, the import of the preceding theorem is that it provides a
cofibrant replacement functor: given a lax symmetric monoidal functor
$\aD \to \Sp$, we can functorial replace it with a monoidal functor
that takes values in flat spectra, i.e., spectra for which the
point-set smash product computes the derived smash product.
Specifically, an orthogonal spectrum $X$ is flat if the functor $X
\sma (-)$ preserves weak equivalences.  Cofibrant spectra are flat,
but there are interesting examples of spectra that are flat but not
cofibrant.  Most notably, the underlying spectra of cofibrant
commutative ring orthogonal spectra are flat but not cofibrant.
More generally, the standard analysis of the underlying spectrum of
cofibrant associative and commutative ring spectra extends to
establish the following result for presheaves.

\begin{prop}
Cofibrant monoids and cofibrant commutative monoids in $\Pre(\aD, \Sp)$ are
pointwise flat.
\end{prop}

\subsubsection{The two-sided bar construction for spectral categories}\label{sec:spectralbar}

Let $\aC$ be a spectral category, $\aM$ a right $\aC$-module, and
$\aN$ a left $\aC$-module.  We can define the smash product of $\aM$
and $\aN$ over $\aC$ as the usual coequalizer
\begin{equation}
\begin{tikzcd}
\aM \sma \aC \sma \aN \ar[r, shift right] \ar[r, shift
  left] & \aM \sma \aN \ar[r] & \aM \sma_{\aC} \aN
\end{tikzcd}
\end{equation}
where the two parallel maps are the actions of $\aC$ on $\aM$ and
$\aN$ respectively.  For practical work, we use the resolution of $\aM
\sma_{\aC} \aN$ given by the two-sided bar construction.  In this
section we review the properties and definition of the bar
construction.

\begin{defin}[Definition 6.1 of~\cite{BlumbergMandell2012}]\label{defn:spectralbar}
Let $\aC$ be a spectral category, $\aM$ a right $\aC$-module and $\aN$
a left $\aC$-module.  Then the two-sided bar construction is the
geometric realization $B(\aM; \aC; \aN)$ of the simplicial spectrum
$B_\bullet(\aM; \aC; \aN)$ with simplices
\begin{equation}
[k] \to \bigvee_{c_1, c_2, \ldots, c_k} \aM(c_1) \sma \aC(c_1, c_2)
\sma \ldots \sma \aC(c_{k-1}, c_k) \sma \aN(c_k),
\end{equation}
degeneracy maps induced by the composition in $\aC$ and the module
structure maps, and face maps induced by the unit of $\aC$.
\end{defin}

The usual simplicial homotopy shows that the canonical map
\begin{equation}
B(\aC; \aC; \aN) \to \aN
\end{equation}
is always a homotopy equivalence of spectra (e.g.,
see~\cite[6.3]{BlumbergMandell2012}).  More generally, there is a
natural map 
\begin{equation}\label{eqn:barcomp}
B(\aM; \aC; \aN) \to \aM \sma_{\aC} \aN
\end{equation}
given by composition.  Equation~\eqref{eqn:barcomp} is an equivalence
under suitable cofibrancy hypotheses, which we now discuss.

First, recall that in order for the bar construction to be tractable,
it is essential to impose conditions to that it is a {\em proper}
simplicial object: levelwise equivalences of proper objects induce
weak equivalences on geometric realization.  See~\cite[\S
  X.2]{ElmendorfKrizMandellMay1997} for a careful discussion of the
geometric realization of proper simplicial spectra; adapting that
argument, in this context it suffices for $\aC$ to be pointwise
cofibrant.

\begin{lem}
Let $\aC$ be a pointwise cofibrant spectral category.  Then
$B_\bullet(\aM; \aC; \aN)$ is a proper simplicial spectrum.
\end{lem}

Given a spectral category $\aC$ and a $\aC$-module $M$, by the
discussion of the preceding section we can produce
\begin{enumerate}
\item a pointwise cofibrant spectral category $\aC'$ along with a
DK-equivalence $F \colon \aC' \to \aC$ and
\item a pointwise cofibrant $\aC'$-module $M'$ along with a
  DK-equivalence $M' \to F^*M$.
\end{enumerate}
Thus, we can conclude the following result about deriving the bar
construction.

\begin{prop}
Let $\aC$ be a spectral category, $\aM$ a right $\aC$-module, and
$\aN$ a left $\aC$-module.  If $\aC$ is pointwise cofibrant and $\aM$
is pointwise cofibrant, then $B(\aM; \aC; -)$ maps
DK-equivalences to weak equivalences.
\end{prop}

\begin{proof}
Let $\aN \to \aN'$ be a DK-equivalence of spectral categories.  Then
there is a levelwise equivalence
\begin{equation}
B_k(\aM'; \aC; \aN) \to B_k(\aM; \aC; \aN)
\end{equation}
and since both of these simplicial objects are proper, there is an
induced weak equivalence on geometric realizations.
\end{proof}

There is also a weaker condition than properness that is sometimes useful,
namely being a {\em split} simplicial object.  We follow the
discussion in~\cite[\S14.4]{Riehl2014} in our exposition of this
situation.

\begin{defin}[Definition 14.4.2 of~\cite{Riehl2014}]
A simplicial space $X_\bullet$ is split if there exist subspaces
$N_\bullet \subset X_\bullet$ such that the canonical map
\begin{equation}
\bigvee_{[n] \twoheadrightarrow [k]} N_k \to X_n
\end{equation}
is a homeomorphism for each $n$.
\end{defin}

The point of this condition is that the geometric realization can be
computed as the filtered colimit of skeleta which are constructed
iteratively via the pushouts
\begin{equation}\label{eq:skpush-appendix}
  \begin{tikzcd}
    \bar{N}_{n} \times \partial \Delta_n \ar[d] \ar[r] & \sk_{n-1}
    X_\bullet \ar[d] \\ N_n \times \Delta_n \ar[r] & \sk_n X_\bullet.
  \end{tikzcd}
\end{equation}

One consequence of this is that even when $X_\bullet$ and $Y_\bullet$
are not proper, a levelwise weak equivalence induces a weak
equivalence on geometric realizations~\cite[14.5.7]{Riehl2014}.

\subsubsection{Functors and bimodules}

We now turn to a discussion of the relationship between spectral
functors and bimodules.  Specifically, we will explain how to use the
representation of spectral functors as bimodules to rectify zigzags
for the purposes of computing the homotopy colimit.

\begin{defin}
Let $F \co \cC_0 \to \cC_1$ be a spectral functor.  Associated to $F$
are the two bimodules:
\begin{enumerate}
\item The $(\cC_0 \sma \cC_1^{\op})$-module ${}^F \cC_1$ specified on
  objects by the assignment
\begin{equation}
(x_0, x_1) \mapsto \cC_1(x_1, Fx_0)
\end{equation}
and 
\item the $(\cC_1 \sma \cC_0^{\op})$-module $\cC_1^{F}$ specified on
  objects by the assignment 
\begin{equation}
(x_1, x_0) \mapsto \cC_1(Fx_0, x_1).
\end{equation}
\end{enumerate}
\end{defin}

We think of the bimodule ${}^F \cC_1$ as encoding $F$ in the sense
that it specifies the data of an assignment of the representable
spectral presheaf $\cC_1(-, Fx)$ to an object $x \in \cC_0$.

These assignments are compatible with composition in the following
sense.

\begin{prop}
Let $F \co \cC_0 \to \cC_1$ and $G \co \cC_1 \to \cC_2$ be spectral
functors.  Then there are natural isomorphism in the homotopy category
\begin{equation}
{}^F \cC_1 \sma^{L}_{\cC_1} {}^G \cC_2 \cong {}^{FG} \cC_2
\end{equation}
and
\begin{equation}
\cC_2^{G} \sma^{L}_{\cC_1} \cC_1^{F} \cong \cC_2^{FG}
\end{equation}
\end{prop}

\begin{proof}
We explain the argument for the first comparison.  Fixing a pair of
objects $x \in \ob(\cC_0)$ and $y \in \ob(\cC_2)$, the derived smash
product on the lefthand side can be written as the bar construction
$B(\cC_1(-,Fx); \cC_1; \cC_2(y, G-))$.  Composition yields a natural
map
\begin{equation}
B(\cC_1(-,Fx); \cC_1; \cC_2(y, G-)) \to \cC_2(x_2, GFx)
\end{equation}
which is a weak equivalence by the usual simplicial
contraction~\cite[6.3]{BlumbergMandell2012}.  It is straightforward to
check that these pointwise equivalences assemble to a natural
transformation through weak equivalences of spectral functors.
\end{proof}

We now explain how to use this formalism to invert zig-zags for the
purposes of computing the bar construction.  Suppose that $F \colon
\aC_0 \to \aC_1$ is a DK-equivalence.  Then in particular, $F$ is
homotopically essentially surjective, and so for any $x_1 \in
\ob(C_1)$ we have $x_1 \htp Fz$ for some object $z$ in $\cC_0$.  Therefore,
\begin{equation}
\begin{tikzcd}
\cC_1(Fx_0, x_1) \htp \cC_1(Fx_0, Fz) & \ar[l,swap,"\htp"] \cC_0(x_0, z).
\end{tikzcd}
\end{equation}
That is, in this case the $(\cC_1 \sma \cC_0^{\op})$-module
$\cC_1^{F}$ can be regarded as specifying the homotopical inverse to
$F$.  Thus, we can compute the composite of $F$ and its inverse as the
bar construction $B({}^F \cC_1, \cC_0, \cC_1^{F})$; there is a natural
equivalence 
\begin{equation}
B({}^F \cC_1, \cC_0, \cC_1^{F}) \htp \cC_0,
\end{equation}
when $\cC_1$ and $\cC_0$ are pointwise cofibrant.
Note that we are describing a familiar phenomenon from Morita theory
here; $\cC_1^{F}$ and $^{F} \cC_1$ are invertible bimodules realizing
the equivalence between $\cC_0$ and $\cC_1$.  In some situations, it
is more convenient to construct an invertible bimodule than a functor
realizing an equivalence; we use this technique in a comparison
starting in Section~\ref{sec:mult-relat-orient}.

For a second example, suppose that we have a zig-zag of spectral
functors
\begin{equation}
\begin{tikzcd}
\cC_0 & \ar["F","\htp"']{l} \ar[swap, "G"]{r} \cC_1 & \cC_2 \\
\end{tikzcd}
\end{equation}
where $F$ is a DK-equivalence.  Then the bimodule representing the
homotopical composite functor $\cC_0 \to \cC_2$ represented by the
zig-zag can be computed as the $(\cC_0 \sma \cC_2^{\op})$-module given
by the derived smash product $\cC_0^{F} \sma^L_{\cC_1} {}^G \cC_2$. 
We can compute this bimodule using the bar construction:
For objects $x \in \ob(\cC_0)$ and $y \in \ob(\cC_2)$, we have 
\begin{equation}
(\cC_0^{F} \sma^L_{\cC_1} {}^G \cC_2)(x,y) \htp B(\cC_0(Fx,-); \cC_1(-,-); \cC_2(y,G-)),
\end{equation}
provided that $\cC_1$ and either $\cC_0$ or $\cC_2$ are pointwise cofibrant.

More generally, given a zig-zag
\begin{equation}
\begin{tikzcd}
\cC_0 & \ar["F_0","\htp"']{l} \ar[swap,"G_0"]{r} \cC_1 & \cC_2 &
\ar["F_1","\htp"']{l} \ar[swap,"G_1"]{r} \cC_3 & \ldots &
\ar["F_{k-1}","\htp"']{l} \cC_{k-1} \ar[swap,"G_{k-1}"]{r} & \cC_k,
\end{tikzcd}
\end{equation}
where the backward functors $\{F_i\}$ are DK-equivalences, we can
construct a model of the composite homotopical functor $\cC_0 \to
\cC_k$ as the $(\cC_0, \cC_k)$-bimodule given by the iterated derived
smash product
\begin{equation}
\left(\cC_0^{F_0} \sma^L_{\cC_1} {}^{G_0} \cC_2 \right) \sma^L_{\cC_2}
\left(\cC_2^{F_1} \sma^L_{\cC_3} {}^{G_1} \cC_4 \right) \sma^L_{\cC_4}
\ldots \sma^L_{\cC_{k-2}} \left(\cC_{k-2}^{F_{k-1}} \sma^L_{\cC_{k-1}}
       {}^{G_{k-1}} \cC_k \right).
\end{equation}
Explicitly, under pointwise cofibrancy hypotheses, we can compute this
as the geometric realization of the multisimplicial bar construction
\begin{align}\label{eq:barzig}
B_{\bullet}(\cC_0(F_0-,-); \cC_1(-,-); \cC_2(-,G_0-); \cC_2(-,-);
\\ \cC_2(F_1-, -); \cC_3(-,-); \cC_4(-,G_1-); \ldots ;
\cC_k(-,G_{k-1}-)) \nonumber
\end{align}
defined as follows.

Recall that a $m$-multisimplicial spectrum is a functor 
\begin{equation}
\underbrace{\Delta^{\op} \times \Delta^{\op} \times \ldots \times \Delta^{\op}}_{m} \to
\Sp.
\end{equation}
When $m=2$, such an object is more commonly referred to as a
bisimplicial spectrum.

\begin{defin}
Given spectral categories $\aC_0, \aC_1, \ldots, \aC_m$ and $\aC_i
\sma \aC_{i+1}^{\op}$-modules $\aM_{i, i+1}$ for $0 \leq i < m$, the
$m$-multisimplicial bar construction 
\begin{equation}
B_{\bullet,\bullet,\ldots,\bullet} (\cM_{0,1}; \cC_1; \cM_{1,2};
\cC_2; \cM_{2,3}; \cC_3; \cM_{3,4}; \ldots; \cC_{m-1}; \cM_{m-1,m}) 
\end{equation}
is the $m$-simplicial $\aC_0 \sma \aC_m^{\op}$-module with simplices
specified by the assignment 
\begin{equation}
\begin{aligned}
[k_0, &k_1, \ldots, k_m] \mapsto \\
\bigvee &\aM_{0,1}(x, y_{1,1}) \sma \aC_1(y_{1,1}, y_{1,2}) \sma \ldots \sma \aC_1(y_{1,k_0-1},
y_{1,k_0}) \\
&\sma \aM_{1,2}(y_{1,k_0},y_{2,1}) \sma \aC_2(y_{2,1},y_{2,2}) \sma
\ldots \sma \aC_2(y_{2,k_1-1},y_{2,k_1}) \\
&\sma \aM_{2,3}(y_{2,k_1},y_{3,1}) \sma \ldots \sma \ldots \sma
\aC_{m-2}(y_{m-2,k_{m-3}-1},y_{m-2,k_{m-3}}) \\ 
&\sma \aM_{m-2,m-1}(y_{m-2,k_{m-3}},y_{m-1,1}) \sma
\aC_{m-1}(y_{m-1,1},y_{m-2,2}) \sma \ldots \\
&\qquad\qquad \sma \aC_{m-1}(y_{m-1,k_{m-1}-1},y_{m-1,k_{m-1}}) \\
&\sma \aM_{m-1,m}(y_{m-1,k_{m-1}},z)
\end{aligned}
\end{equation} 
where the multisimplicial structure maps are induced by the bimodule
actions on $\cM_{i,i+1}$, the compositions in $\aC_i$, and the unit.
\end{defin}

We can form the geometric realization by passing to the diagonal
(yielding a simplicial spectrum) and then taking the usual geometric
realization, or by taking iterative geometric realizations in any
order.  We continue to denote the geometric realization of a
multisimplicial spectrum as $|-|$.

\begin{lem}
There is a natural map of 
\begin{equation}
\begin{tikzcd}
\relax|B_{\bullet,\bullet,\ldots,\bullet} (\cM_{0,1}; \cC_1; \cM_{1,2};
\cC_2; \cM_{2,3}; \cC_3; \cM_{3,4}; \ldots; \cC_{k-1}; \cM_{m-1,m})\relax|
\ar[d] \\
\cM_{0,1} \sma_{\cC_1} \cM_{1,2} \sma_{\cC_2} \sma \cM_{2,3}
\sma_{\cC_3} \sma_{\cM_{3,4}} \ldots \sma_{\cC_{i}} \ldots \sma_{\cC_{k-1}} \cM_{m-1,m} 
\end{tikzcd}
\end{equation}
that is a weak equivalence when all of $\{\cC_i\}$ and
$\{\cM_{j,j+1}\}$ are pointwise cofibrant.
\end{lem}

\subsubsection{Homotopy colimits}\label{sec:hoco}

We briefly review the practical theory of homotopy colimits and record
here some technical material required for manipulating them.  We will
work with a version of the Bousfield-Kan definition of the homotopy
colimit in terms of the bar construction, following~\cite[\S
  5]{Riehl2014}.  For a comprehensive discussion of the relationship
to the left derived functor of the colimit, also see the excellent
treatment in~\cite{Shulman2006}.

We begin by defining the relevant version of the bar construction.  We
will fix a symmetric monoidal category $\cV$ and consider
$\cV$-enriched categories; as usual, in our examples $\cV$ will be
either spaces or spectra.

\begin{defin}
Let $\aC$ be a $\cV$-enriched category that is tensored and $\cA$ be a
small category.  For functors $F \colon \cA^{\op} \to \cV$ and $G
\colon \cA \to \aC$, the simplicial bar construction
$B_\bullet(F;\aA;G)$ is the simplicial object of $\aC$ with
$k$-simplices
\begin{equation}
[k] \mapsto \coprod_{a_0 \to a_1 \to \ldots \to a_k} F(a_0) \otimes
G(a_k),
\end{equation}   
where $\otimes$ denotes the tensor of an object of $\cV$ and an object
of $\aC$.
\end{defin}

In order to form the geometric realization of $B_\bullet(F;\aC;G)$ we
need to assume that $\aC$ has more structure, specifically an
enrichment in either spaces or simplicial sets.  In this situation, we
denote by $B(F; \aC; G)$ the geometric realization
$|B_\bullet(F;\aC;G)|$.

\begin{rem}
The bar construction is a fattened-up version of the tensor product of
functors, which in the setting above is defined to be the coequalizer
\begin{equation}
\begin{tikzcd}
F \otimes_{\cA} G = \Coeq \Bigl( \coprod_{f \colon a \to a'} G(a')
\otimes F(a) \ar[shift left=2,"f^*"]{r} \ar[shift
  right=2,swap,"f_*"]{r} & \coprod_a G(a) \otimes F(a) \Bigr).  \\
\end{tikzcd}
\end{equation}
Specifically, we have natural isomorphisms
\begin{equation}
B(F; \aA; G) \cong B(F; \aA; \aA) \otimes_{\aA} G.
\end{equation}
and
\begin{equation}
B(F; \aA; G) \cong F \otimes_{\aA} B(\cA; \cA; G).
\end{equation}
Since the standard simplicial homotopies imply that the maps $B(\cA;
\cA; G) \to G$ and $B(F; \aA; \aA) \to F$ are equivalences, the bar
construction should be thought of as tensoring with a resolution of
$F$ or $G$, respectively. 
\end{rem}

We now assume that:
\begin{enumerate}
\item $\aC$ is a cofibrantly generated model category,
\item $\aC$ is enriched in simplicial sets or topological spaces and
  the enrichment is compatible with the model structure (i.e.,
  satisfies the analogue of Quillen's SM7), and
\item $\aC$ admits functorial cofibrant replacement.
\end{enumerate}
In this context, we can define the homotopy colimit in terms of the
bar construction.

\begin{defin}\label{defn:hoco}
Let $\cA$ be a small category, and $F \colon \cA \to \aC$ a functor
that takes values in cofibrant objects.  Then we define the homotopy
colimit via the formula
\begin{equation}
\hocolim_{\cA} F = B(F; \cA; \ast),
\end{equation}
where $\ast$ denotes the constant functor to simplicial sets at the
terminal object and the bar construction is computed as the geometric
realization of the simplicial object in $\aC$ with $k$-simplices
specified by the formula
\begin{equation}
[k] \mapsto \coprod_{a_0 \to a_1 \to \ldots \to a_k} F(a_0).
\end{equation}
\end{defin}
Notice that if $\aC$ is the category of simplicial sets and we take
$F$ to also be the functor $\ast$, then 
\begin{equation}
\hocolim_{\cA} \ast = |N_\bullet \aA|,
\end{equation} 
i.e., the nerve of $\aA$.

The hypothesis on $F$ guarantees that the construction of the homotopy
colimit is invariant under natural transformations of diagrams that
are pointwise weak equivalences.  In the event that we are considering
a functor $F$ that does not take values in cofibrant objects, we
precompose with the cofibrant replacement functor on $\aC$ to obtain a
homotopy-invariant construction from the formula of
Definition~\ref{defn:hoco}.

We will use the fact that when $\aC$ is a symmetric monoidal category
in which the tensor commutes with colimits in each variable, the
homotopy colimit inherits a natural external product structure.  We
specialize to the case of main interest.  

\begin{lem}
Suppose we have two functors $F \colon \cA \to \Sp$ and $G \colon \cA
\to \Sp$.  Then there is a natural map
\begin{equation}\label{eq:hocosmash}
\hocolim_{\cA} F \sma \hocolim_{\cA} G \to \hocolim_{\cA \times \cA}
(F \sma G)
\end{equation}
specified as
\begin{equation}
B(F; \cA; \ast) \sma B(G; \cA; \ast) \to B(F \sma G; \cA \times \cA;
\ast),
\end{equation}
defined on the $k$-simplices as
\begin{equation}
\coprod_{a_0 \to a_1 \to \ldots \to a_k} F(a_0) \sma \coprod_{a'_0 \to
  a'_1 \to \ldots \to a'_k} G(a'_0) \to \coprod_{\substack{(a_0 \to
    a_1 \to \ldots \to a_k), \\ (a'_0 \to a'_1 \to \ldots \to a'_k)}}
F(a_0) \sma G(a'_0).
\end{equation}
\end{lem}

Next, we quickly discuss some comparison results for homotopy
colimits.  Suppose that we have diagram categories $\cA$ and $\cB$, a
functor $G \colon \cA \to \cB$, and a functor $F \colon \cB \to \cC$
where $\aB$ and $\aC$ are model categories as above.  In this case,
there is a comparison map
\begin{equation}
\gamma \colon \hocolim_{\cA} F \circ G \to \hocolim_{\cB} F.
\end{equation}
The functor $G \colon \cA \to \cB$ is {\em homotopy final} if for each
$b \in B$ the comma category $b \downarrow F$ has contractible nerve.

\begin{lem}[Theorem 8.5.6 in~\cite{Riehl2014}]
If $G \colon \cA \to \cB$ is homotopy final, then comparison map
$\gamma$ is a weak equivalence for any functor $F \colon \aB \to \aC$.
\end{lem}

\begin{notn}
In some previous work, homotopy final functors had been referred to as
homotopy cofinal.  We adopt the recent consensus (e.g.,
see~\cite[8.3.3]{Riehl2014}) that the correct term is homotopy final;
under this convention, the dual condition (which we will not need)
refers to {\em homotopy initial} functors.
\end{notn}

Taking $F$ to be the functor $\cB \to \Set^{\Delta^{\op}}$ that is
constant at a point, we immediately obtain the following version of
Quillen's theorem A.

\begin{cor}
A homotopy final functor $G \colon \aA \to \aB$ induces a weak
equivalence of simplicial sets $N_{\bullet}G \colon N_{\bullet} \aA
\to N_{\bullet} \aB$.
\end{cor}

We will find the following condition for checking homotopy finality
very useful.

\begin{prop}[See Lemma 8.5.2 in~\cite{RadulescuBanu2006}]\label{prop:adjfinal} 
Let $G \colon \aA \to \aB$ be a right adjoint.  Then $G$ is homotopy
final.
\end{prop}

Interesting applications of Proposition~\ref{prop:adjfinal} arise when
$G$ is the inclusion of a full subcategory; in this case, the
conclusion is that the inclusion of a reflective subcategory is
homotopy final.

One of the advantages of Definition~\ref{defn:hoco} is that it can be
adapted to any suitable enriched category; we have already seen this
in the context of spectral categories.  More generally, we now assume
that $\aC$ is a $\cV$-enriched model category, where $\cV$ is itself a
model category equipped with a well-behaved notion of geometric
realization.  Specifically, we assume that we are given a cosimplicial
object $\Delta^{\bullet}$ in $\cV$.  Then given a simplicial object
$X_\bullet$ in $\cV$, the geometric realization of $X_\bullet$ can be
defined as
\begin{equation}
|X_\bullet| = \int_{\Delta} \Delta^n \otimes X_n.
\end{equation}
Moreover, the functor $|-|$ has a right adjoint (specified by the
mapping object out of $\Delta^n$).  If we assume that $|-|$ is strong
symmetric monoidal and a left Quillen functor, then $\aC$ and $\cV$
have the structure of simplicial model categories.  In this context,
we can define a weighted homotopy colimit; see~\cite[\S 9]{Riehl2014}
for further discussion of the properties of this definition.

\begin{defin}\label{defn:weightedhoco}
Let $\cA$ be a small $\cV$-category that is pointwise cofibrant, $F
\colon \cA \to \aC$ a functor that takes values in cofibrant objects,
and $G \colon \cA^{\op} \to \cV$ a functor that takes values in
cofibrant objects.  Then we define the weighted homotopy colimit of
$F$ with weights $G$ via the formula 
\begin{equation}
\hocolim_{\cA}^G F = B(F; \cA; G),
\end{equation}
by which we mean the geometric realization of the simplicial object
in $\cV$ with $k$-simplices
\begin{equation}
[k] \mapsto \coprod_{a_0, a_1, \ldots, a_k} F(a_0) \otimes \cA(a_0,
a_1) \otimes \ldots \otimes \cA(a_{k-1},a_k) \otimes G(a_k).
\end{equation}
\end{defin}

One problematic aspect of using Definition~\ref{defn:weightedhoco} is
the issue of ensuring that $F$ and $G$ take values in cofibrant
objects.  Cofibrant replacement functors are typically not enriched,
and so additional hypotheses are necessary in practice; in general, we
need to assume the existence of the projective model structure on
diagrams.  See~\cite[\S 23]{Shulman2006} and~\cite[9.2]{Riehl2014} for
a discussion of this point (and of the cofibrancy conditions required
in the definition more generally).

\subsection{$2$-categories and $2$-functors}\label{sec:2categories} 

In this section, we give a brief review of the definitions we need
from the theory $2$-categories and bicategories.  We refer the reader
to Lack's exposition in~\cite{Lack2010} for a more detailed treatment.

\begin{defin}\label{defin:2cat}
A $2$-category $\aC$ is a category enriched in categories:
\begin{enumerate}
\item A class $\ob(\aC)$ of objects (the $0$-cells).
\item For each pair $x,y \in \ob(\aC)$ a category $\aC(x,y)$; the
  objects of $\aC(x,y)$ are referred to as $1$-cells and the morphisms
  as $2$-cells.
\item For each triple $x,y,z$ of objects there is a strictly
  associative and unital composition functor.
\end{enumerate}
\end{defin}

For example, the category $\Cat$ of categories has an enrichment in
categories given by taking the functor category as the category of
morphisms.  Another natural class of examples comes from permutative
(strict) monoidal categories: there is a single object and the
category of morphisms is given by the objects of the monoid, with
composition the monoidal composition law.  However, many
natural examples are not $2$-categories because the composition isn't
strict; e.g., monoidal categories.

This leads to a weaker notion of a $2$-category, given by the theory
of bicategories.  Roughly speaking, the idea is that a bicategory
$\aC$ is a category enriched over categories in a weak sense; $\aC$ is
equipped with a mapping category $\aC(x,y)$ for each pair of objects
$x,y$ and composition and unit functors that satisfy associativity and
unitality conditions up to natural isomorphism.

\begin{defin}\label{defin:bicat}
A bicategory $\aC$ consists of the following data:
\begin{enumerate}
\item A class $\ob(\aC)$ of objects (the $0$-cells).
\item For each pair $x,y \in \ob(\aC)$ a category $\aC(x,y)$; the
  objects of $\aC(x,y)$ are referred to as $1$-cells and the morphisms
  as $2$-cells.
\item For each $x \in \ob(\aC)$, a distinguished $1$-cell $\id_x \in
  \ob(\aC(x,x))$.
\item For $x,y,z \in \ob(\aC)$, a composition functor
\begin{equation}
\aC(x,y) \times \aC(y,z) \to \aC(x,z).
\end{equation}
\item For $f \in \ob(\aC(w,x))$, $g \in \ob(\aC(x,y))$, and $h \in
  \ob(\aC(y,z))$, a natural isomorphism $(fg)h \to f(gh)$.

\item For $f \in \ob(\aC(x,y))$, natural isomorphisms $\id_x f \to f$
  and $f \id_y \to f$.

\item Associativity pentagons and unit diagrams that strictly commute;
  see e.g.,~\cite[1.0]{Leinster1998}.
\end{enumerate}
\end{defin}

When the associativity and unit isomorphisms are the identity, this
data just specifies an enrichment over $\Cat$, i.e., a $2$-category as
in Definition~\ref{defin:2cat}; to distinguish this situation, we will
refer to this as a {\em strict} $2$-category.  We shall presently see
that every bicategory is equivalent to a strict $2$-category, just as
every symmetric monoidal category is equivalent to a permutative
category.

\begin{defin}\label{defin:lax-functor}
Let $\aC$ and $\aD$ be bicategories.  A {\em lax functor} $F \colon
\aC \to \aD$ consists of the following data:
\begin{enumerate}
\item A function $F \colon \ob(\aC) \to \ob(\aD)$.
\item For every $x,y \in \ob(\aC)$, a functor $F_{xy} \colon \aC(x,y)
  \to \aD(Fx,Fy)$.
\item For $1$-cells $f \colon x \to y$ and $g \colon y \to z$, natural
  transformations (i.e., $2$-cells) $Fg \circ Ff \to F(g \circ f)$.
\item For $0$-cells $x \in \ob(\aC)$, natural transformations (i.e.,
  $2$-cells) $\id_{Fx} \to F(\id_x)$.
\item Associativity and unitality diagrams for the $2$-cells described
  in the preceding items; see e.g.,~\cite[1.1]{Leinster1998}.
\end{enumerate}
When the $2$-cells are natural isomorphisms, $F$ is a {\em
  pseudofunctor}.  When the $2$-cells are in fact identities, $F$ is a
{\em strict $2$-functor}.
\end{defin}

There is an evident analogue of a natural transformation.

\begin{defin}
Let $F,G \colon \aC \to \aD$ be lax functors between bicategories
$\aC$ and $\aD$.  A lax transformation $\tau \colon F \to G$ consists
of the following data:
\begin{enumerate}
\item For each $x \in \ob(\aC)$, a $1$-cell $\tau_x \colon F(x) \to
  G(x)$ in $\aD$.
\item For each $1$-cell $f \colon x \to y$ in $\aC$, a $2$-cell $Gf
  \circ \tau_x \to \tau_y \circ Ff$.
\item Associativity and unit diagrams for the $2$-cells described
  above.
\end{enumerate}
When the $2$-cells are natural isomorphisms, $F$ is a {\em strong}
transformation.  When the $2$-cells are in fact identities, $F$ is a
{\em strict} transformation.
\end{defin}

Using the preceding definition, we can now define an equivalence of
bicategories.

\begin{defin}
Let $\aC$ and $\aD$ be bicategories.  A {\em biequivalence} of bicategories
between $\aC$ and $\aD$ consists of pseudofunctors $F \colon \aC \to
\aD$ and $G \colon \aD \to \aC$ and strong transformations $G \circ F
\to \id_{\aC}$ and $\id_{\aD} \to F \circ G$.  If $\aC$ and $\aD$ are
strict $2$-categories, $F$ and $G$ are strict functors, and the
transformations are strict isomorphisms, we say this is a
{\em $2$-equivalence} of $2$-categories.
\end{defin}

Using a version of Isbell's construction~\cite{Isbell1969}, we have
the following rectification theorem.

\begin{thm}[e.g., see Section 2.3.3 of~\cite{Gurski2006}]\label{thm:2catrect}
Let $\aC$ be a bicategory.  There exists a strict $2$-category $\aC'$
and a biequivalence $\aC' \to \aC$.
\end{thm}

The rectification is functorial in the following sense.

\begin{thm}[e.g., see Section 2.4.3 of~\cite{Gurski2006}]\label{thm:funcrect}
Let $F \colon \aC \to \aD$ be a pseudofunctor between bicategories
$\aC$ and $\aD$.  Then there exists a strict functor $F' \colon \aC'
\to \aD'$, where $\aC'$ and $\aD'$ are the strictifications of $\aC$
and $\aD$ such that the square
\begin{equation}
\begin{tikzcd}
\aC \ar[d] \ar["F"]{r} & \aD \ar[d] \\ \aC' \ar[swap,"F'"]{r} & \aD'\\
\end{tikzcd}
\end{equation}
commutes.
\end{thm}

Finally, we note that there are two different possible generalizations
of the notion of the opposite of a category to the setting of
bicategories; these correspond to reversing the $1$-cells or the
$2$-cells, respectively.

\begin{defin}
Let $\aC$ be a bicategory.
\begin{enumerate}
\item The opposite bicategory $\aC^{\op}$ has
the same objects as $\aC$ and we define the category of morphisms
$\aC^{\op}(x,y)$ to be $\aC(y,x)$.  
\item The conjugate bicategory $\aC^{conj}$ has the
same objects as $\aC$ and we define the category of morphisms
$\aC^{c}(x,y)$ to be the opposite category $\aC(x,y)^{\op}$.
\end{enumerate}
In each case, the rest of the structure of the
bicategory is defined in the evident way.
\end{defin}

\subsection{Group actions on categories}\label{sec:groupact}

In this section, we discuss the foundations of the theory of groups
actions on categories.  In general, this is a technically demanding
subject; fortunately, in the work at hand, the actions we encounter
are very rigid.

Let $G$ be a discrete group and $BG$ the strict $2$-category with a
single object, $1$-cells the discrete category with objects the
elements of $G$ and identity morphisms, and $2$-cells given by the
composition in $G$ (i.e., the monoidal structure on the $1$-cells).
Let $\Cat$ denote the strict $2$-category of categories, functors, and
natural transformations.

\begin{defin}
An \emph{action} of $G$ on a category $\aC$ is a $2$-functor
\begin{equation}
\gamma \colon BG \to \Cat
\end{equation}
The action is {\em strict} when $\gamma$ is a strict functor, {\em
  pseudo} when $\gamma$ is a pseudofunctor, and {\em lax} when
$\gamma$ is a lax functor.
\end{defin}

Unpacking this data, the action of $G$ on a category $\aC$ is
specified by:
\begin{enumerate}
\item A collection of functors $F_g \colon \aC \to \aC$ indexed by $g
  \in G$ and
\item suitable associative and unital natural transformations
\begin{equation}
\alpha_{gh} \colon F_g \circ F_h \to F_{gh}
\end{equation}
for $g \in G$.
\end{enumerate}
When the action is strict, the natural transformations $\alpha_{gh}$
are the identity.  When the action is a pseudo-action, these natural
transformations are isomorphisms.  A lax action simply has natural
transformations.  Note that this definition requires that the
composition transformations are strictly associative and unital.

\begin{example}
Any category $\aC$ can be endowed with the trivial $G$-action by
considering the strict functor which takes the unique object of $BG$
to $\aC$, each $g \in G$ to the identity, and for which all the
$2$-cells are the identity.  We will denote this $2$-functor by
$\id_G$; the category $\aC$ will be clear from context.
\end{example}

\begin{rem}
Although we do not use this perspective in the paper, we note that
another way to encode such equivariant structures is in terms of
categorical fibrations over $BG$.
\end{rem}

In our context, we are interested in the generalization where $\Cat$
is replaced by a category of enriched categories.

\begin{defin}
Let $V$ be a symmetric monoidal category.  Denote by $\Cat_V$ denotes
the $2$-category of $V$-enriched categories, enriched functors, and
enriched natural transformations.  
\end{defin}

We now have an analogous definition of the action of a group $G$ on a
$V$-enriched category.  We give a general definition, although in our
work $V$ is typically either be the category of topological spaces or
the category of orthogonal spectra.

\begin{defin}
An \emph{action} of $G$ on a $V$-category $\aC$ is a $2$-functor
\begin{equation}
\gamma \colon BG \to \Cat_V.
\end{equation}
The action is {\em strict} when $\gamma$ is a strict functor, {\em
  pseudo} when $\gamma$ is a pseudofunctor, and {\em lax} when
$\gamma$ is a lax functor.
\end{defin}

Recall that for a pair of $V$-enriched categories $\aC$ and $\aD$, we
can form the tensor product $\aC \otimes \aD$, with objects
$\ob(\aC) \times \ob(\aD)$ and morphisms given by the tensor
$\aC(x,y) \otimes \aD(x',y')$.  The following lemma records the
compatibility of the tensor product of $V$-enriched categories with
group actions.

\begin{lem}
Let $\aC$ and $\aD$ be $V$-enriched categories with actions by $G$.
Then $\aC \otimes \aD$ is a $V$-enriched category with an action of $G$.
\end{lem}

\begin{proof}
Let $F_g^{\aC}$ and $F_g^{\aD}$ denote the functors encoding the
$G$-actions on $\aC$ and $\aD$ respectively.  The required functors
$\aC \otimes \aD \to \aC \otimes \aD$ are defined as
$F_g^{\aC} \otimes F_g^{\aD}$; the natural transformations expressing
composition are defined analogously.
\end{proof}

\subsubsection{Cofibrant replacement and group
  actions}\label{sec:equicof}

We now specialize to the case where $V = \Sp$, to discuss the
interaction of group actions with cofibrant replacement.
Specifically, in order to maintain homotopical control, it will be
necessary for us to cofibrantly replace spectral categories $\aC$ with
actions of $G$ such that the result inherits an equivariant structure.
This turns out to be straightforward for strict actions, since the
pointwise-cofibrant replacement is functorial.

Specifically, given spectral functors $\gamma_{g_1} \colon \aC \to
\aC$ and $\gamma_{g_2} \colon \aC \to \aC$, there are induced functors 
$\gamma_{g_1}^c \colon \aC^{c} \to \aC^{c}$ and $\gamma_{g_2} \colon
\aC^{c} \to \aC^{c}$ such that
\begin{equation}
\gamma^c_{g_1} \circ \gamma^c_{g_2} = (\gamma_{g_1} \circ
\gamma_{g_2})^c.
\end{equation}
(Here $\aC^{c}$ denotes the pointwise cofibrant replacement.)  When
these functors are part of a strict action, $\gamma_{g_1} \circ
\gamma_{g_2} = \gamma_{g_1 g_2}$ and so
\begin{equation}
\gamma^c_{g_1} \circ \gamma^c_{g_2} = \gamma_{g_1 g_2}^c.
\end{equation}
This suggests that the assignment $g \mapsto \gamma^c_g$ specifies the
strict action of $G$ on $\aC^c$.  To verify this, note that
analogously, $(-)^c$ carries the associativity diagrams for the action
of $G$ on $\aC$ to the associativity diagrams for an action of $G$ on
$\aC^c$, and similarly for the unitality diagrams.  This discussion
proves the following proposition.

\begin{prop}
Let $\aC$ be a spectral category with a strict action of $G$.  Then
$\aC^c$ is a pointwise-cofibrant spectral category with a strict
action of $G$.
\end{prop}

\subsubsection{Homotopical group action on categories}
\label{sec:group-acti-spectr}

The coherent system of spheres $\{\bS[n]\}$ constructed in
Proposition~\ref{prop:luriesphere} gives rise to a lax action of $\bZ$
on the category $\Sp$ of spectra or on the category of $R$-modules for
any ring spectrum $R$.

\begin{prop}\label{prop:actionspec}
The assignment $n \mapsto (-) \sma \bS[-n]$ specifies a lax functor
from $B\bZ$ to $\Cat$, where the unit and associativity
transformations are induced by the associative ring structure on the
system $\{\bS[-n]\}$.
\end{prop}

In fact, since the functor $(-) \sma \bS[-n] \colon \Sp \to \Sp$ is
evidently spectrally enriched, the proof of
Proposition~\ref{prop:actionspec} immediately extends to provide a lax
functor from $B\bZ$ to $\Spcat$.

Although the multiplication map $\bS[-n] \sma \bS[-m] \to \bS[-n-m]$
is not a homeomorphism, it is a weak equivalence.  We describe this
situation as follows.

\begin{defin}\label{defn:homoact}
A {\em homotopy action} of a group $G$ on a model category $\aC$ is
given by a lax action of $G$ on $\aC$ such that the unit and
associativity transformations are through weak equivalences in $\aC$.
\end{defin}

Specializing and rewriting, we have the following basic result.

\begin{thm}\label{thm:kmodact}
The assignment $n \mapsto (-) \sma \bS[-n]$ specifies a homotopy
action of $\bZ$ on $\bbk-\mod$.
\end{thm}

In the interests of concision, we do not give a complete treatment of
the theory of homotopical actions of discrete groups.  In fact, the
only thing we really need is a version of the following consistency
check, which provides a justification for
Definition~\ref{defn:homoact} (and our language in describing it).
(See~\cite[1.4.3]{Hovey1999} for discussion of the relevant
$2$-category of model categories and pseudofunctor induced by passage
to the homotopy category and derived functors.)

\begin{prop}
Let $F$ be a lax functor from $BG$ to the $2$-category of model
categories, Quillen adjunctions, and natural transformations such that
$F$ specifies a homotopy action of $G$ on $\aC$.  Then composition
with the pseudofunctor specified on objects by passage to the homotopy
category yields a pseudoaction of $G$ on $\Ho(\aC)$.
\end{prop}

\subsubsection{Groups actions, functors, and bimodules}\label{sec:equifunc}

We now turn to discussion of equivariant functors between categories
with $G$ actions.  In the following, we omit the modifiers on functors
and transformations except when necessary, as there are analogous
versions of the definitions and results for each degree of
strictness.

\begin{defin}\label{defn:equivariant-functors}
Suppose that $\aC$ and $\aD$ are categories with actions of $G$, i.e.,
there are $2$-functors $\gamma_{\aC} \colon BG \to \Cat$ and
$\gamma_{\aD} \colon BG \to \Cat$ that pick out $\aC$ and $\aD$.  Then
an equivariant functor $F \colon \aC \to \aD$ is specified by a
natural transformation $\gamma_{\aC} \to \gamma_{\aD}$.
\end{defin}

Spelling this out a little bit, an equivariant functor from $\aC$ to
$\aD$ is specified by:
\begin{enumerate}
\item A functor $F \colon \aC \to \aD$.
\item Coherent natural transformations $\gamma_g \circ F \to
F \circ \gamma_g$ for each $g \in G$.
\end{enumerate}

We are most interested in the setting of equivariant spectral
categories.  

\begin{defin}
Let $\aC$ be a spectral category with an action of $G$.  Then a
$G$-equivariant $\aC$-module is a $G$-equivariant functor from $\aC$
to $\Sp$, where $\Sp$ is given the trivial $G$-action.
\end{defin}

Writing out part of this data, a $G$-equivariant
$\aC$-module is specified by a $\aC$-module $\aM$ equipped with
natural transformations $\gamma_g \co \aM(-) \to \aM(g-)$ for each $g \in
G$ such that the diagrams
\begin{equation}
\begin{tikzcd}
\aC(x,y) \sma \aM(x) \ar[swap,"\gamma_g \sma \gamma_g"]{d} \ar{r} &
\aM(y) \ar["\gamma_g"]{d} \\ \aC(gx,gy) \sma \aM(gx) \ar{r} & \aM(gy)
\\
\end{tikzcd}
\end{equation} 
commute.

We now want to describe the interaction of equivariant spectral
functors $\aC \to \aD$ with $G$-equivariant bimodules.  Recall that
given a functor $F \colon \aC \to \aD$, there are associated
$\aC \sma \aD^{\op}$ and $\aD \sma \aC^{\op}$ bimodules $^F \aD$ and
$\aD^F$ specified on objects by the assignments $(c,d) \mapsto \aD(d,
Fc)$ and $(c,d) \mapsto \aD(Fc,d)$, respectively.

\begin{prop}
Let $\aC$ and $\aD$ be spectral categories with a $G$-action, and let
$F \colon \aC \to \aD$ be a $G$-equivariant spectral functor.  Then
the bimodules $^F \aD$ and $\aD^F$ are $G$-equivariant bimodules.
\end{prop}

\begin{proof}
We give the argument for $^F \aD$; the argument for $\aD^F$ is
analogous.  We are given a spectral functor $F \colon \aC \to \aD$ and
a lax transformation $\gamma_{\aC} \to \gamma_{\aD}$ that lifts $F$,
i.e., coherent natural transformations $\tau_g \colon \gamma_g \circ F
\to F \circ \gamma_g$ for each $g \in G$.  We need to construct a lax
transformation $\gamma_{\aC \sma \aD} \to \id_G$, which amounts to
providing natural transformations
\begin{equation}
\aD(-,F-) \to \aD(\gamma_g(-),F(\gamma_g(-)))
\end{equation}
that satisfy the required compatibilites.  These are constructed as
the composite
\begin{equation}
\begin{tikzcd}
\aD(-,F-) \ar["\gamma_g"]{r} & \aD(\gamma_g(-), \gamma_g(F(-)))
\ar["\tau_g"]{r} & \aD(\gamma_g(-), F(\gamma_g(-))).
\end{tikzcd}
\end{equation}
\end{proof}

We now turn to a discussion of the equivariant structure of the bar
construction in the setting of $G$-equivariant spectral categories.
In the following statement, recall that an orthogonal spectrum with
$G$-action is an object of the Borel category of equivariant spectra;
this is equivalent to a module spectrum over $\Sigma^{\infty}_+ G$.

\begin{thm}
Let $\aC$ be a spectral category with a strict $G$-action, and suppose
that we have a right $\aC$-module $M$ and a left $\aC$-module $N$
which are strictly $G$-equivariant.  Then the two-sided bar
construction $B(M; \aC; N)$ is an orthogonal spectrum with $G$-action.
\end{thm} 

\begin{proof}
Producing the structure of an orthogonal spectrum with $G$-action
given an orthogonal spectrum $X$ is equivalent to producing a map of
topological monoids $G_+ \to \Map(X,X)$.  For each $n$ and $g \in G$,
the $G$-actions on $M$, $N$, and $\aC$ give rise to maps of spectra
\begin{equation}
\begin{tikzcd}
\displaystyle \bigvee_{c_0, \ldots c_n} M(c_0) \sma \aC(c_0,c_1) \sma
\ldots \sma \aC(c_{n-1},c_{n}) \sma N(c_n) \ar[d,"\theta_g"] \\ \displaystyle
\bigvee_{gc_0, \ldots gc_n} M(gc_0) \sma \aC(gc_0,gc_1) \sma \ldots
\sma \aC(gc_{n-1},gc_{n}) \sma N(gc_n).
\end{tikzcd}
\end{equation}
Since the action maps are via functors, the maps $\theta_g$ assemble
to maps of simplicial spectra
\begin{equation}
\theta_g \colon B_\bullet(M;\aC; N) \to B_\bullet(M; \aC; N).
\end{equation}
Since the actions and functors are strict, we see that
$\theta_g \circ \theta_h = \theta_{gh}$, and so these maps specify the
structure of an orthogonal spectrum with $G$-action.
\end{proof}

More generally, the same argument implies we have the following basic
result that lets us handle the multisimplicial bar construction that
arises from zigzags of spectral functors.

\begin{prop}
Let $\cC_0, \cC_1, \ldots, \cC_k$ be strictly $G$-equivariant spectral
categories.  Given a zig-zag of strictly $G$-equivariant functors 
\begin{equation}
\begin{tikzcd}
\cC_0 & \ar[l,"F_0"] \ar[r, "G_0",swap] \cC_1 & \cC_2 & \ar[l,"F_1"]
\ar[r,"G_1",swap] \cC_3 & \ldots & \ar[l,"F_{k-1}"] \cC_{k-1}
\ar[r,"G_{k-1}",swap] & \cC_k,
\end{tikzcd}
\end{equation}
then the bar construction of Equation~\eqref{eq:barzig} is endowed
with an induced $G$-action.
\end{prop}

In our case of interest, the last category in the zig-zag has a
homotopical action rather than a strict action.  In this case, the
argument above shows that the multisimplicial bar construction
produces a homotopy coherent $G$-action on the spectrum.  As we review
in Section~\ref{sec:rect} below, a homotopy coherent action is simply
given by a homotopy coherent diagram in spectra indexed by the
category $BG$ with a single object and morphism set $G$.  This is a
simplicial diagram in spectra (regarded as a simplicial category via
the standard topological enrichment) indexed by $\fC BG$, the free
resolution of $BG$.  Homotopy coherent diagrams are equivalent to
strict diagrams in the sense that any orthogonal spectrum with a
homotopy coherent $G$-action can be rectified to an equivalent
spectrum with $G$-action.

\begin{prop}
Let $\cC_0, \cC_1, \ldots, \cC_{k-1}$ be strictly $G$-equivariant spectral
categories and let $\cC_k$ be a spectral category with a homotopical
$G$-action.  Given a zig-zag of strictly $G$-equivariant functors
\begin{equation}
\begin{tikzcd}
\cC_0 & \ar[l,"F_0"] \ar[r, "G_0",swap] \cC_1 & \cC_2 & \ar[l,"F_1"]
\ar[r,"G_1",swap] \cC_3 & \ldots & \ar[l,"F_{k-1}"] \cC_{k-1}
\ar[r,"G_{k-1}",swap] & \cC_k,
\end{tikzcd}
\end{equation}
then the bar construction of Equation~\eqref{eq:barzig} is endowed
with a homotopy coherent $G$-action.
\end{prop}

\begin{proof}
The key point is that the coherence data of the homotopical $G$-action
on $\aC$ shows that the maps $\theta_g$ specify a simplicial functor
from the free resolution of the category $BG$ specified by $G$.  This
is essentially immediate; $\fC(BG)$ has $k$-simplices the iterated
free group on the elements of $G$.
\end{proof}

\subsubsection{Group actions on bicategories}\label{sec:equibicat}

In addition to the notion of a group action on a category, we will
also need the notion of a $\Pi$-equivariant bicategory and of lax
functors between $\Pi$-equivariant bicategories, for a discrete group
$\Pi$.  Since we do not want to discuss tricategories in this paper,
we will write down explicitly the data we have in mind.  This is
manageable only because the version of the notion we are dealing with is
extremely strict.  (See~\cite{Bernaschini2019} for a more detailed
treatment of this notion.)

\begin{defin}
Let $\aC$ be a bicategory and $\Pi$ a discrete group.  A strict
$\Pi$-action on $\aC$ is determined by specifying for each $\pi \in
\Pi$, a strict $2$-functor $\gamma_{\pi} \colon \aC \to \aC$.  We
require that this assignment is
\begin{enumerate}
\item strictly unital, in the sense that the $2$-functors $\gamma_{e}
  \circ \gamma_{\pi_1}$, $\gamma_{\pi_1} \circ \gamma_e$ and
  $\gamma_{\pi_1}$ are equal.
\item strictly respects composition, in the sense that the
  $2$-functors $\gamma_{\pi_1} \circ \gamma_{\pi_2}$ and
  $\gamma_{\pi_1 \pi_2}$ are equal.
\item strictly associative, in the sense that the $2$-functors
  $\gamma_{\pi_1 \pi_2} \circ \gamma_{\pi_3}$ and $\gamma_{\pi_1}
  \circ \gamma_{\pi_2 \pi_3}$ are equal.
\end{enumerate}
\end{defin}

We record the following relationship between a strict action on a
bicategory and the associated topological category.

\begin{prop}
Let $\aC$ be a $2$-category (a strict bicategory) with a strict
$\Pi$-action.  Then the associated topological category formed by
applying the classifying space functor to each morphism category has a
strict $\Pi$-action.
\end{prop}

Let $\aC$ and $\aD$ be bicategories with strict $\Pi$-actions.  We
will work with two notions of equivariant $2$-functors between $\aC$
and $\aD$ in this paper.  First, we define a strictly
$\Pi$-equivariant lax $2$-functor $F \colon \aC \to \aD$ as follows.

\begin{defin}
A strictly $\Pi$-equivariant lax $2$-functor $F \colon \aC \to \aD$ is
a lax $2$-functor such that for each $\pi \in \Pi$, the $2$-functors
$F \circ \gamma_{\pi}$ and $\gamma_{\pi} \circ F$ are equal.  That is,
we require that the diagram
\begin{equation}
\begin{tikzcd}
\aC \ar["\gamma_{\pi}"]{r} \ar[swap,"F"]{d} & \aC \ar["F"]{d} \\ \aD
\ar[swap,"\gamma_{\pi}"]{r} & \aD \\
\end{tikzcd}
\end{equation}
commutes for each $\pi \in \Pi$.
\end{defin}

The strictness of the data we are requiring in this context makes it
possible to check that our examples satisfy these requirements by
checking a comparatively small number of explicit diagrams, even
though the underlying functor itself is lax.

We will also work with the variant where we relax the strictness of
the $\Pi$-action at the cost of making the underlying functor more
strict.  Specifically, we define a pseudo $\Pi$-equivariant strict
$2$-functor $F \colon \aC \to \aD$ as follows.

\begin{defin}
A pseudo $\Pi$-equivariant strict $2$-functor $F \colon \aC \to \aD$
is a strict $2$-functor such that for each $\pi \in \Pi$, there is a
pseudonatural equivalence connecting the $2$-functors $F \circ
\gamma_{\pi}$ and $\gamma_{\pi} \circ F$.  That is, we require that
the diagram
\begin{equation}
\begin{tikzcd}
\aC \ar["\gamma_{\pi}"]{r} \ar[swap,"F"]{d} & \aC \ar["F"]{d} \\ \aD
\ar[swap,"\gamma_{\pi}"]{r} & \aD \\
\end{tikzcd}
\end{equation}
commute up to natural isomorphism for each $\pi \in \Pi$.  Moreover,
we require a variety of coherence diagrams to commute,
see~\cite[2.3]{Bernaschini2019} for explicit details.
\end{defin}

Once again, it is possible to check that our examples satisfy these
requirements using a comparatively manageable amount of data.

\subsection{Internal categories and bicategories enriched in internal
  categories}\label{sec:internal}

In addition to topologically enriched categories, where the morphism
sets are topologized, we also work with categories where the set
of objects is also equipped with a topology.  Such objects are
referred to as internal categories in topological spaces (or sometimes
category objects in spaces).  In this section, we give a rapid review;
see~\cite[\S XII.1]{MacLane} and~\cite[\S A]{Lind2013} for more details.

\begin{defin}
An internal category $\aC$ in topological spaces consists of the data
of spaces $\aC_0 = \ob(\aC)$ and $\aC_1 = \mor(\aC)$ along with
\begin{enumerate}
\item Source and target maps $t,s \colon \aC_1 \to \aC_0$,
\item an identity map $\id \colon \aC_0 \to \aC_1$,
\item and composition maps $\aC_1 \times_{\aC_0} \aC_1 \to \aC_1$,
  where the pullback is over the source and target maps respectively, 
\end{enumerate}
such that the evident associativity and unitality diagrams commute.
\end{defin}

Notice that an internal category in topological spaces with a discrete
space of objects specifies the same data as a topologically enriched
category; given a topological category $\aC$, the associated internal
category has 
\begin{equation}
\aC_0 = \ob(\aC) \quad\textrm{and}\quad \aC_1 = \coprod_{\ob(\aC)
  \times \ob(\aC)} \aC(x,y).
\end{equation}  
The notion of an internal category can be defined in terms of any ambient
category that has enough pullbacks, although we will only require the
example of spaces.

Functors between internal categories are specified in terms of
suitable maps between object and morphism objects.

\begin{defin}\label{defn:internal-functor}
Let $\aC$ and $\aD$ be internal categories.  An internal functor $F
\colon \aC \to \aD$ is specified by the data of continuous maps $F_0
\colon \aC_0 \to \aD_0$ and $F_1 \colon \aC_1 \to \aD_1$ which are
suitably compatible with the source, target, composition, and identity
maps.  Specifically, we require that the diagrams
\begin{equation}
\begin{tikzcd}
\aC_1 \times_{\aC_0} \aC_1 \ar[r, "f_1 \times_{f_0} f_1"] \ar[d] &
\aD_1 \times_{\aD_0} \aD_1 \ar[d] \\ \aC_1 \ar[r,swap, "f_1"] & \aD_1
\end{tikzcd}
\end{equation}
and
\begin{equation}
\begin{tikzcd}
\aC_1 \ar[d, "f_1"] \ar[r, "s", shift left] \ar[r, swap, "t", shift
  right] & \aC_0 \ar[r, "\id"] \ar[d, "f_0"] & \aC_1 \ar[d, "f_1"]
\\ \aD_1 \ar[r, "s", shift left] \ar[r, swap, "t", shift right] &
\aD_0 \ar[r, swap, "\id"] & \aD_1
\end{tikzcd}
\end{equation}
commute.
\end{defin}

Similarly, there is a notion of internal natural transformations
between internal functors.

\begin{defin}
Let $\aC$ and $\aD$ be internal categories and $F$ and $G$ internal
functors $\aC \to \aD$.  An internal natural transformation $F \to G$
is specified by a map $\tau \colon \aC_0 \to \aD_1$ satisfying:
\begin{enumerate}
\item The composites $s \circ \tau = F$ and $t \circ \tau G$.
\item The diagram 
\begin{equation}
\begin{tikzcd}
\aC_1 \ar[r,"{(G_1,\tau \circ s)}"] \ar[d, swap, "{(\tau \circ t,F_1)}"] & \aD_1 \times_{\aD_0} \aD_1 \ar[d] \\ 
\aD_1 \times_{\aD_0} \aD_1 \ar[r] & \aD_1
\end{tikzcd}
\end{equation}
commutes.
\end{enumerate}
\end{defin}

Putting this together, we can organize the category of internal
categories as follows.

\begin{prop}
There is a $2$-category of internal categories in spaces, internal
functors, and internal natural transformations.  This is in fact a
symmetric monoidal $2$-category with respect to the cartesian product.
\end{prop}

We can define monoidal and symmetric monoidal internal categories in
terms of the symmetric monoidal structure on the $2$-category of
internal categories.

We also have an enrichment and internal hom on the category of
internal functors, which are constructed as follows.  First, given
internal categories $\aC$ and $\aD$, there is a space of internal
functors $\Fun(\aC, \aD)$ defined as the subspace of $\Map(\aC_0,
\aD_0) \times \Map(\aC_1, \aD_1)$ that satisfy the requirements of
Definition~\ref{defn:internal-functor}.

\begin{prop}
The mapping spaces $\Fun(\aC,\aD)$ specify an enrichment in spaces on
the category of internal categories and internal functors.
\end{prop}

The mapping spaces provide the object spaces of internal categories of
functors.  That is, we define $\Fun(\aC,\aD)_0 = \Fun(\aC, \aD)$ and
$\Fun(\aC,\aD)_1$ is defined to be the subspace of $\Fun(\aC,\aD)_0
\times \Fun(\aC,\aD)_0 \times \Map(\aC_0, \aD_1)$ satisfying the
requirement that for the tuple $(f_0, f_1, \gamma)$, $\gamma$
specifies the data of a natural transformation at each point between
$f_0$ and $f_1$.  (In the interest of concision, we do not write out
in detail the diagrams representing this compatibility.)

\begin{prop}
The category of internal categories and internal functors is Cartesian
closed.
\end{prop}

\begin{example}
We will be most interested in the category of internal diagrams.  That
is, for a discrete category $I$ and an internal category $\aC$, we
have an internal category of $I$-shaped diagrams in $\aC$.
\end{example}

Given an internal category $\aC$ in spaces, applying $\pi_0$ to both
the objects and the morphisms yields an ordinary category $\pi_0 \aC$;
this assignment is evidently functorial.

\begin{prop}\label{prop:internal_discretize}
There is a functor $\pi_0$ from the category of internal categories in
spaces to the category of categories.
\end{prop}

We will often deal with internal functors $F \colon \aC \to \aD$ where
$\aD$ has a discrete object set.  As a consequence, it is useful to be
explicit about what conditions the discrete object set of $\aD$
imposes on $F$.  
\begin{enumerate}
\item The requirement that $F \colon \aC_0 \to
\aD_0$ be continuous means that $F$ must factor through $\pi_0
\aC_0$.
\item The continuous maps $\aC(x,y) \to \aD(Fx,Fy)$ can only depend
on the image of $x$ and $y$ in $\aD_0$; e.g., when $F$ factors through
$\pi_0 \aC_0$, the maps on morphism spaces must depend only on the
component the objects belong to.
\end{enumerate}

We can define the classifying space of an internal category in spaces
as the evident generalization of the usual notion.

\begin{defin}
The nerve of an internal category in topological spaces $\aC$ is the
simplicial space $N_\bullet \aC$ with $n$-simplices
\begin{equation}
[n] \mapsto \underbrace{\aC_1 \times_{\aC_0} \aC_1 \times_{\aC_0}
  \times \ldots \times_{\aC_0} \aC_1}_{n},
\end{equation}
where we understand $N_0 \aC = \aC_0$ and $N_1 \aC = \aC_1$.  The
degeneracies are induced by the identity map and the degeneracies by
the composition maps.
\end{defin}

When either of the natural maps $s,t \colon \aC_1 \to \aC_0$ is a
fibration, the pullbacks that appear in $N_\bullet \aC$ are homotopy
pullbacks.  When the identity map $\aC_0 \to \aC_1$ is a cofibration,
the nerve is a proper simplicial space.  We then have the following
easy lemma.

\begin{lem}
Let $\aC$ and $\aD$ be internal categories in spaces such that at
least one of the maps $s$ and $t$ are fibrations and the identity maps
for $\aC$ and $\aD$ are cofibrations.  Then an internal functor $F
\colon \aC \to \aD$ such that $F_0$ and $F_1$ are weak equivalences
induces a weak equivalence $|N_\bullet \aC| \to |N_\bullet \aD|$.
\end{lem}

Next, we consider arrow categories in the context of internal
categories.  Let $\aC$ be an internal category.  The arrow category
$\Ar \aC$ of $\aC$ is the internal category with objects $(\Ar \aC)_0
= \aC_1$ and morphisms 
\begin{equation}
(\Ar \aC)_1 = (\aC_1 \times_{\aC_0} \aC_1) \times_{\aC_1} (\aC_1 \times_{\aC_0}
\aC_1),
\end{equation}
where the pullback is over the composition maps.  Composition is
induced by commuting pullbacks and the compositions maps in $\aC$.
Similarly, we define the twisted arrow category of $\aC$ as follows.

\begin{defin}\label{defn:internal-twisted-arrow}
Let $\aC$ be an internal category in spaces.  Then the internal
twisted arrow category $\Tw \aC$ has object space $(\Tw \aC)_0 =
\aC_1$ and morphisms specified by the pullback diagram
\begin{equation}
\begin{tikzcd}
(\Tw \aC)_1 \ar[r] \ar[d] & \ar[d] \aC_1 \times_{\aC_0} \aC_1
  \times_{\aC_0} \aC_1 \ar[d] \\ \ast \ar[r] & \aC_1, \\
\end{tikzcd}
\end{equation}
where the vertical map is the iterated composition and the bottom
horizontal map picks out a morphism in $\aC_1$.  Composition is
induced by the composition in $\aC$ and commuting pullbacks.
\end{defin}

We now quickly review the theory of homotopy colimits in a context
where the indexing category $\cX$ is an internal category in spaces.
We refer the reader to Appendix A of~\cite{Lind2013} for further
discussion. 

Without loss of generality we will assume that $i \colon \cX_0 \to
\cX_1$ is a Hurewicz cofibration; the usual whiskering
construction ensures that we can harmlessly replace $\cX$ if necessary
to ensure this condition holds.

\begin{defin}
Left and right $\cX$-modules are spaces $M$ and $N$ over $\cX_0$, respectively
equipped with action maps
\begin{equation}
\cX_1 \times_{\cX_0} M \to M \qquad N \times_{\cX_0}
\cX_1 \to N.
\end{equation}
\end{defin}

We can now give a bar construction version of the definition of a
homotopy colimit.

\begin{defin}
Given an internal category $\cX$ and right and left $\cX$-modules $M$
and $N$, the bar construction $B_\cdot(M, \cX, N)$ is the simplicial
space defined via the assignment
\begin{equation}
[k] \mapsto M \times_{\cX_0} \cX_1 \times_{\cX_0} \ldots
\times_{\cX_0} \cX_1 \times_{\cX_0} N
\end{equation}
and the obvious structure maps coming from composition, the unit, and
the module structures.
\end{defin}

In particular, we have the following definition.

\begin{defin}\label{defn:internal_hocolim}
Given a left $\cX$-module $F$, the (uncorrected) homotopy colimit of
$F$ over $\cX$ is specified to be the bar construction $B(\ast, \cX,
F)$, where $\ast$ denotes the constant module at the point.
\end{defin}

However, we often are interested in considering the homotopy colimit
of an enriched functor $F \colon \cX \to \Top$; recall that such an
enriched functor ignores the topology on the object set of $\cX$.
In particular, $F(x)$ depends only on $[x] \in \pi_0 \cX_0$.  In
this case, the homotopy colimit $\hocolim_{\cX} F$ is specified as the
homotopy colimit associated to the $\cX$-module $\hat{F}$ constructed
as the (enriched) coend: 
\begin{equation}\label{eq:simple-internal-hocolim}
\int^{x \in \cX_0} (\cX_1 \times_{\cX_0} \{x\}) \times
F(x),
\end{equation}
where here the pullback on the left indicates the subspace of
morphisms with source $x$.  The module structure is given by the
obvious composition and the target map.  The $\cX$-module $\hat{F}$ is
constructed to make the usual coend formulas (which ignore the
topology on the object set) correct.  That is, when it exists, we can
equivalently compute the homotopy colimit using the factorization
through the enriched category with objects $\pi_0 \cX_0$ and morphisms
the fibers of $\cX_1$ over the components.

In order to compare homotopy colimits over internal categories, we
will need a version of Quillen's theorem A that applies in this
context.  For our applications, it will suffice to consider the
situation where we have a functor $F \colon \cA \to \cB$ with $\cB$
having a discrete set of objects, but we give a general statement (see
e.g.,~\cite[4.7]{EbertRandalWilliams2019} for a discussion and proof).
In what follows, we write $(F/B)_{0,0}$ for the subspace of $\ob(\cA)
\times \mor(\cB)$ consisting of pairs $(a_0, F(a_0) \to b_0)$.

\begin{prop}\label{prop:internal-Quillen-A}
Let $\cA$ and $\cB$ be internal categories in spaces and $F$ a
continuous functor $F \colon \cA \to \cB$ such that
\begin{enumerate}
\item For each object $b \in \ob(\cB)$, the overcategory $F/b$ is
  contractible,
\item the target map $\hom(\aA) \to \ob(\aA)$ is a fibration,
\item the map $(F_D)_{0,0} \to \ob(B)$ specified by $(a_0, F(a_0) \to
  b_0) \mapsto b_0$ is a fibration,
\item either the source or the target map $\hom(\cB) \to \ob(\cB)$ is
  a fibration, and
\item the identity maps $\ob(\aA) \to \aA$ and $\ob(\aB) \to \cB$ are
  Hurewicz cofibrations.
\end{enumerate}
Then $F$ induces a weak equivalence on classifying spaces.
\end{prop}

Finally, we will need to work with $2$-categories and bicategories
where the morphism categories are given by internal categories.  The
former are simply categories enriched in internal categories in
spaces, and the latter are bicategories enriched in the monoidal
$2$-category of internal categories of
spaces~\cite{GarnerShulman2016}.  We refer to these as {\em
  topological $2$-categories} and {\em topological bicategories}
respectively.  We do not write out the definitions, since what we need
from the theory of such categories is quite limited.

\subsection{Homotopy limits over zig-zags and homotopy coherent diagrams}\label{sec:rect}

Let $\widetilde{I}$ denote a category with 
\begin{enumerate}
\item objects $\bN$ or $\bZ$, and 
\item non-identity morphisms generated by the requirement that
for each pair $(i,i+1)$ there is a unique morphism either from $i \to
i+1$ or from $i+1 \to i$.  
\end{enumerate}
If $\aC$ is a model category, a zig-zag diagram in $\aC$ over
$\widetilde{I}$ is a functor $\widetilde{I} \to \aC$ such that the
morphisms from $i$ to $i+1$ are taken to weak equivalences.  Such a
diagram are a homotopical model of a directed system; the zig-zag from
$i$ to $j$ represents a map in the homotopy category.  (Note that here
we are fixing a preferred direction for the diagram; the discussion is
entirely equivalent if we choose the other polarity.)

We have several options for computing the homotopy (co)limit of a
diagram $D \colon \widetilde{I} \to \aC$.  One possibility is simply
to directly compute the homotopy (co)limit of the zig-zag diagram, via
the Bousfield-Kan formula.  Another possibility is to rectify: to
replace the original diagram $D \colon \aC \to \widetilde{I}$ with a
weakly equivalent diagram (i.e., construct a natural transformation $D
\rightarrow D'$) such that each entry is a cofibrant-fibrant object of
$\aC$.  Then one can choose homotopy inverses for each backwards map,
i.e., convert $D'$ to a diagram indexed on $\bN$ or $\bZ$ regarded as
a poset.  The homotopy co(limit) of $D'$ can be computed directly, and
the work
of~\cite[\S3.3]{AngeltveitBlumbergGerhardtHillLawsonMandell2018} shows
that computing the homotopy limit in this fashion produces a weakly
equivalent result.  Another outcome of the arguments there is that the
existence of parallel zig-zags that are homotopic do not change the
homotopy limit.
 
Another way to study this situation is to use the formal theory of
rectification of homotopy coherent diagrams.  Specifically, let $I$ be
a small category and $\aC$ be an arbitrary category.  There is an
evident natural functor
\begin{equation}\label{eq:forgetcoherence}
\Ho(\aC^I) \to \Ho(\aC)^I.
\end{equation}
An object in $\Ho(\aC^I)$ can be thought of as represented by a
cofibrant-fibrant diagram in a suitable model structure on $\aC^I$,
when such a model structure is available.  Alternatively, we can
regard it as a homotopy coherent diagram on
$I$~\cite{Vogt1973,CordierPorter1997}.  
Given a category $I$, let $UI$
denote the underlying graph where the vertices are objects and the
edges specified by morphisms.  This forgetful functor has as left
adjoint the free category on a graph, and we write $\fC I$ to denote
the simplicial resolution associated to the monad of this adjunction.

\begin{defin}
When $\aC$ is a simplicial model category and $I$ is a small category,
a homotopy coherent diagram of shape $I$ is a simplicial functor $\fC
I \to \aC$, where $\fC I$ denotes the free simplicial resolution of
$I$.
\end{defin}

The main result of Vogt amounts to a proof that $\Ho(\aC^I)$ is
equivalent to the homotopy category of homotopy coherent diagrams.

\begin{thm}
There is an equivalence of categories 
\begin{equation}
\Ho(\aC^I) \cong \Coh(\aC,I),
\end{equation}
where the latter denotes the homotopy category with
objects the homotopy coherent $I$-shaped diagrams in $\aC$ and
morphisms the homotopy classes of natural transformations.
\end{thm}

A particular example of interest is when $I$ is the category $BG$
associated to a discrete group $G$, with a single object and morphisms 
the elements of $G$.  Then a homotopy-coherent action of $G$ on a
space or spectrum is given by a simplicial functor from $\fC BG$.

Since it is much easier to construct diagrams in the homotopy
category, i.e., objects of $\Ho(\aC)$, it is essential to understand
the obstruction to producing a section of the object function of
Equation~\eqref{eq:forgetcoherence}.  As one would hope, there is a
concrete obstruction theory for this problem~\cite{Cooke1978}.
However, it is clear from the definition that the obstruction theory
is trivial when $I$ is a free category.

\begin{prop}
Let $I$ be a free category on a graph $G$.  Then there is a section of
the comparison functor; any homotopy commutative diagram of shape $I$
can be rectified to a homotopy coherent diagram of shape $I$.\qed
\end{prop}

Recall that the free category on the graph $G$ has objects the
vertices of $G$ and non-identity arrows the strings of compatible
edges of $G$, i.e., strings $e_1, e_2, \ldots, e_n$ where $s(e_i) =
t(e_{i-1})$.  For our purposes, the most interesting free categories
are the posets $\bN$ and $\bZ$; these are the free categories on the
graphs with vertices in bijection with the objects and edges
corresponding to pairs $\{i,i+1\}$.

\begin{cor}\label{cor:Nrect}
Any homotopy commutative diagram over $\bN$ or $\bZ$ (or in fact any
countable totally ordered set) can be rectified to a homotopy coherent
diagram.
\end{cor}

Now, suppose we have a zig-zag diagram $\widetilde{I}$ as above.  This
discussion gives a formal rectification procedure: choosing homotopy
inverses, we obtain a diagram in the homotopy category, which we can
rectify using Corollary~\ref{cor:Nrect}.  The homotopy limit of the
rectified diagram coincides with the homotopy limit over
$\widetilde{I}$.

\subsection{The twisted arrow category}
\label{sec:twisted}

In this section, we record facts that we need about the twisted arrow
category construction (sometimes referred to as the subdivision of a
category).
See~\cite[7.2.9]{Riehl2014},~\cite[2.3.5]{Kerodon}, or~\cite[3.1]{Penney2017}
for more comprehensive treatments; the latter two in particular give
clear explanations of the relationship to the bicategory of spans.

\begin{defin}
The \emph{twisted arrow category} $\Tw \cC$ of a category $\cC$ is the
category whose objects are arrows $f \co \alpha \to \beta$ in $\cC$,
and whose morphisms from $f_0$ to $f_1$ are given by factorisations of
$f_0$ through $f_1$, i.e., by diagrams 
\begin{equation} \label{eq:morphism_in_twisted_arrow}
\begin{tikzcd}
\alpha_0 \ar[r, "f_0"] \ar[d,"g",swap] & \beta_0 \\
\alpha_1 \ar[r,"f_1"] & \beta_1. \ar[u,"h",swap]
\end{tikzcd}
\end{equation}
\end{defin}

\begin{example}
The twisted arrow category of $0 \leftarrow 01 \to 1$  is given by
  \begin{equation}
    \id_0 \leftarrow f_0 \to  \id_{01} \leftarrow f_1 \to \id_1.    
  \end{equation}  
\end{example}

\begin{lem}
There is a natural functor $\Tw \aC \to \aC$ which assigns
to each arrow its domain, and which yields a commutative diagram
\begin{equation}
  \begin{tikzcd}
    \Tw \aC \times \Tw \aD \ar[d] \ar[r] & \Tw (\aC \times \aD) \ar[dl]\\
     \aC \times \cD. & 
  \end{tikzcd}
\end{equation} \qed
\end{lem}

The construction of the twisted arrow category specifies a monoidal
endofunctor on the category of categories; this has the extremely
useful consequence that applying the twisted arrow category to the
morphism categories of a bicategory (or a $\Pi$-equivariant
bicategory) produces a new ($\Pi$-equivariant) bicategory.

\begin{lem}\label{lem:morphism-wise-twisting}
Given a $\Pi$-equivariant bicategory $\aC$, we can form a
$\Pi$-equivariant bicategory $\Tw \aC$ such that the objects of $\Tw
\aC$ are those of $\aC$ and the morphism categories are specified by
the assignment 
\begin{equation}
(\Tw \aC)(p,q) = \Tw (\aC(p,q)).
\end{equation} \qed
\end{lem}

As befits a subdivision, the geometric realization of the twisted
arrow category of $\aC$ is the same as that of $\aC$.

\begin{prop}
The functor $\Tw \aC \to \aC$ induces a weak equivalence of simplicial
sets 
\begin{equation}
N_\bullet \Tw \aC \to N_\bullet \aC.
\end{equation} \qed
\end{prop}

Pulling back along $\Tw \aC \to \aC$ induces an equivalence on
homotopy colimits.

\begin{lem}
Given a functor $F \colon \aC \to \aD$, the natural map
\begin{equation}
\hocolim_{\Tw \aC} F \to \hocolim_{\aC} F
\end{equation}
is a weak equivalence.\qed
\end{lem}

Next, we record here a standard strategy for checking that a function
$\ob(\Tw \aC) \to \ob(\aC)$ in fact is functorial in the twisted arrow
category.  Any arrow
\begin{equation}
\begin{tikzcd}
X_0 \ar[d,swap,"g"] \ar[r,"f_0"] & X_0' \\
X_1 \ar[r,swap,"f_1"] & X_1' \ar[u,swap,"h"]
\end{tikzcd}
\end{equation}
in $\Tw(\aC)$ can be factored as the composites
\begin{equation}\label{eq:factor-twisted-arrows}
\begin{tikzcd}
X_0 \ar[d,swap,"g"] \ar[r,"f_0"] & X_0' && X_0 \ar[d,swap,"\id"] \ar[r,"f_0"] & X_0' \\
X_1 \ar[d,swap, "\id"] \ar[r, "f_1 \circ h"] & X_0' \ar[u,swap,"\id"] && 
X_0 \ar[d,swap,"g"] \ar[r,"f_1 \circ g"] & X_1' \ar[u,swap,"h"] \\
X_1 \ar[r,swap, "f_1"] & X_1' \ar[u,swap,"h"] && X_1 \ar[r,swap,"f_1"] & X_1'
\ar[u,swap,"\id"].
\end{tikzcd}
\end{equation}

As a consequence, to check that an assignment is functorial, it
suffices to check compositions in $\Tw(\aC)$ of the form $(g,\id)
\circ (\id, h)$, $(g_1, \id) \circ (g_0, \id)$, and $(\id, h_0) \circ
(\id, h_1)$ and then a composition in which one map is generic and one
is either of the form $(g, \id)$ or $(\id, h)$.  This observation is
particularly useful when checking the existence of a natural
transformation $F \Rightarrow G$ for functors $F,G \colon \Tw \aC \to
\aD$.  By definition, such a transformation is specified by maps
$\nu_f \colon F(f) \to G(f)$ in $\aD$ that commute with maps $f \to g$
in $\Tw \aC$.

\begin{lem}\label{lem:twisted-arrow-natural-transformation}
Let $F$ and $G$ be functors $\Tw \aC \to \aD$.  A natural
transformation $F \Rightarrow G$ is specified by a collection of maps
$\nu_f \colon F(f) \to G(f)$ such that the naturality diagrams hold
for either of the classes of pairs in
Equation~\eqref{eq:factor-twisted-arrows}.
\end{lem}

To provide more context for the discussion above, we conclude our
discussion by explaining the connection between the twisted arrow
category and the $2$-category of spans.

\begin{defin}
For a category $\aC$, recall that the $2$-category of spans in $\aC$
has objects the zig-zags $x \leftarrow y \rightarrow z$, which
morphisms the category of commutative diagrams
\begin{equation}
\begin{tikzcd}
& y \ar[dl] \ar[dd] \ar[dr] & \\
x & & z \\
& y' \ar[ul] \ar[ur] 
\end{tikzcd}
\end{equation}
and the $2$-cells determined via pullback.  The associativity
isomorphisms are the usual comparisons of pullbacks.
\end{defin}

The connection between spans and the twisted arrow category is that
functors $\Tw(\aC)^{\op} \to \aD$ are in bijection with normal oplax
functors from $\aC$ to spans on $\aD$.  (Here recall that a normal
functor has unit transformation the identity.)  In fact, there is an
equivalence of suitable categories, as follows (e.g.,
see~\cite[2.3.5]{Kerodon} and~\cite[App.~A]{Errington1999}).

\begin{thm}\label{thm:twisted-to-span}
There is an equivalence between the category of functors $\Tw \aC \to
\aD$ with morphisms natural transformations and the category of normal
oplax functors $\aC \to \Span(\aD)$ with morphisms isomorphism classes
of oplax transformations such that the $1$-cells admit right adjoints.
\end{thm}

\begin{proof}[Sketch of proof]
We mainly explain how the bijection on objects is constructed.  Assume
we are given a functor $F \colon \Tw \aC \to \aD$.  We construct the
functor $\tilde{F} \colon \aC \to \Span(\aD)$ as follows.  On an
object $c \in \aC$, we set $\tilde{F}(c) = F(\id_c)$.  For a morphism
$f \colon c \to d$ in $\aC$, we have the span
\begin{equation}
\begin{tikzcd}
F(\id_c) & F(f) \ar[l] \ar[r] & F(\id_d) 
\end{tikzcd}
\end{equation}
formed from the arrows in $\Tw \aC$
\begin{equation}
\begin{tikzcd}
c \ar["f",r] \ar["\id",swap,d] & d && c \ar["f",swap,d] \ar["f",r] & d \\
c \ar["\id",swap,r] & c \ar["f",swap,u] && d \ar["\id",swap,r] & d \ar["\id",swap,u].
\end{tikzcd}
\end{equation}
Write this span as $\tilde{F}(f)$.

Given maps $f \colon c \to d$ and $g \colon d \to e$ in $\aC$, we
produce the morphism 
\begin{equation}
\tilde{F}(f \circ g) \to \tilde{F}(f) \circ \tilde{F}(g) 
\end{equation}
specified by the diagram
\begin{equation}
\begin{tikzcd}
& F(f \circ g) \ar[dr,"{F(f,\id_e)}"] \ar[dl,swap,"{F(\id_c, g)}"] & \\
F(f) \ar[dr,swap,"{F(f,\id_c)}"] & & F(g) \ar[dl,"{F(\id_d, g)}"] \\
& F(\id_d), &
\end{tikzcd}
\end{equation}
where here the notation $F(-,-)$ refers to the image of the map of
spans determined by the pair of arrows.
\end{proof}

\subsection{Categories of cubes}\label{sec:cubes} 

The purpose of this section is to record some basic results about the
categories of cubical diagrams.  Let $\b1$ denote the category
specified by the poset $\{0 \to 1\}$, and $\b1^n$ the product.  There
are face maps
\begin{equation}
d^{i,\epsilon} \colon \b1^n \to \b1^{n+1}
\end{equation}
which insert $\epsilon$ at position $i$, and degeneracy maps
\begin{equation}
s^i \colon \b1^n \to \b1^{n-1}
\end{equation}
which project away from the $i$th entry.  The category $\Box$ has
objects the categories $\b1^n$ and morphisms generated by the face and
degeneracy maps.  This is an elegant Reedy category with fibrant
constants.

\begin{defin}\label{defn:cubes}
For a (small) category $\aC$, let $\Box \aC$ denote the category
with objects the functors $\b1^n \to \aC$ and morphisms the commutative
diagrams 
\begin{equation}
\begin{tikzcd}
\b1^n \ar[d] \ar[r] & \b1^m \ar[dl] \\
\aC. &
\end{tikzcd}
\end{equation}
We refer to $\Box \aC$ as the category of cubes in $\aC$.
\end{defin}

Given a functor $F \colon \aC \to \aM$, where $\aM$ is a model
category, we can compute the homotopy colimit $\hocolim_{\aC} F$ by
precomposing with the functor $t$ that evaluates at the terminal object to
get a functor $\Box \aC \to \aM$, cofibrantly replacing in the Reedy
model structure, and computing the colimit $\hocolim_{\Box \aC} F$.

\begin{prop}
Let $\aM$ be a model category.  For any functor $F \colon \aC \to
\aM$, there is a natural equivalence
\begin{equation}
\hocolim_{\Box \aC} F \circ t \to \hocolim_{\aC} F.
\end{equation}
\end{prop}

In our work, constructing degeneracies introduces difficulties, so we
consider nondegenerate cubes.

\begin{defin}
For a small category $\aC$, let $\Box_{\nd} \aC$ denote the full
subcategory of $\Box \aC$ specified by the functors that do not factor
through any degeneracy map, i.e., for which no morphism is taken to
the identity of $\aC$.  Note that the morphisms here are specified by
the inclusions of cubes.  We refer to this category as the
nondegenerate cubes in $\aC$.
\end{defin}

There is a natural inclusion functor $\Box_{\nd} \aC \to \Box \aC$.
Corresponding to any element $\sigma \in \Box \aC$ is a unique
nondegenerate cube $\tilde{\sigma} \in \Box_{\nd} \aC$ that factors
maps out of $\sigma$ to elements of $\Box_{\nd} \aC$.  We now need to
impose an additional condition on $\aC$ in order to simplify the
homotopical properties of $\Box_{\nd} \aC$.

\begin{defin}
We will say that an element of $\Box_{\nd} \aC$ is totally
nondegenerate if all of its faces are also elements of $\Box_{\nd}
\aC$.  We will refer to $\Box_{\nd} \aC$ as nonsingular if all
elements are totally nondegenerate.
\end{defin}

Since the faces of $\Box_{\nd} \aC$ involves composition of maps, the
following proposition is clear.

\begin{lem}\label{lem:posnon}
Let $\aC$ be a (small) category.  Suppose that for any $f$ in $\aC$,
if $gf = \id$ or $fg = \id$ then $f = g = \id$; i.e., $\aC$ has no
nontrivial isomorphisms and no nontrivial retractions.  Then
$\Box_{\nd} \aC$ is nonsingular.
\end{lem}

The requirement that $\Box_{\nd} \aC$ be nonsingular makes it possible
to check that the assignment of the unique nondegenerate cube
corresponding to a degenerate cube is part of the data of an adjoint.

\begin{lem}
Suppose that $\Box_{\nd} \aC$ is nonsingular.  Then the
inclusion $\Box_{\nd} \aC \to \Box \aC$ has a left adjoint. 
\end{lem}

As a consequence of Proposition~\ref{prop:adjfinal}, we now conclude
the following comparison theorem.

\begin{thm}\label{thm:boxnd}
Assume that $\Box_{\nd} \aC$ is nonsingular.  Let $F \colon
\Box \aC \to \aM$ be a functor.  Then the inclusion functor $\iota
\colon \Box_{\nd} \aC \to \Box \aC$ induces a weak equivalence
\begin{equation}
\hocolim_{\Box_{\nd} \aC} F \circ \iota \to \hocolim_{\Box \aC} F
\end{equation}
\end{thm}

In particular, Lemma~\ref{lem:posnon} and Theorem~\ref{thm:boxnd}
implies that for a poset or more generally for a category without
non-identity retractions, we can compute the homotopy colimit of a
functor over $\aC$ using the pullback to $\Box_{\nd} \aC$.

\section{Parametrized spectra, orientation theory, and ambidexterity}
\label{sec:orientations}

\subsection{Parametrized spectra}\label{sec:parametrized-spectra}

In this section, we review the some basic definitions from the theory
of parametrized spaces and spectra over a topological space $B$ with a
$G$-action.  We will assume throughout that $B$ is locally compact,
Hausdorff, and second-countable.  However, we will not necessarily
assume that $B$ has the homotopy type of a $G$-$CW$ complex.  

We will denote by $\cT^G_{/ B}$ the category of ex-spaces over $B$,
i.e., $G$-spaces over $B$ equipped with a chosen section $s$ that
equips the fibers with basepoints.  For a finite-dimensional inner
product space $W$, we write $S_B^W$ for the parametrized $W$-sphere
that has total space $S^W \times B$ and fibers $S^W$.  The category
$\cT^G_{/ B}$ has a closed symmetric monoidal under the fiberwise
smash product; in particular, we will write $\Omega_B^V$ for the
fiberwise mapping spectrum from $S^V_B$.

Roughly speaking, recall that a parametrized orthogonal $G$-spectrum
over $B$ is defined to be a collection of parametrized spaces
$\{X(V)\}$ over $B$ equipped with $O(V)$-actions and $O(W) \times
O(V)$-equivariant structure maps
\begin{equation}
S^W_B \sma X(V) \to X(V \oplus W).
\end{equation}
Parametrized spectra can be given a diagrammatic interpretation akin
to the definition of orthogonal $G$-spectra.

\begin{defin}
A parametrized orthogonal $G$-spectrum over $B$ is an enriched functor
$X \colon \sI_G^U \to \cT^G$ equipped with structure maps
\begin{equation}
\sigma_{V,W} \colon X(V) \sma S_B^W \to X(V \oplus W)
\end{equation}
that comprise a natural transformation of functors and are associative
and unital.  A map of parametrized orthogonal $G$-spectra is a natural
transformation that commutes with the structure map.  
We will denote by $\Sp_{/ B, G}$ the category of parametrized
orthogonal $G$-spectra over $B$, writing $\Sp^U_{/ B, G}$ when
necessary to emphasize the universe $U$.
\end{defin}

We will often use an external version of the symmetric monoidal
structure on parametrized spaces and spectra.

\begin{defin} \label{def:external_smash_para}
The external smash product of parametrized spaces $f_1 \colon E_1 \to
B_1$ and $f \colon E_2 \to B_2$ is the parametrized space $E_1
\barwedge E_2$ over $B_1 \times B_2$, with fibers $F_1 \sma F_2$.
Specifically, the total space is defined as the pushout 
\begin{equation}
\begin{tikzcd}
(E_1 \times B_2) \cup_{B_1 \times B_2} (B_1 \times E_2) \ar[r] \ar[d]
  & E_1 \times E_2 \ar[d] \\
B_1 \times B_2 \ar[r] & E_1 \barwedge E_2,
\end{tikzcd}
\end{equation}
where the maps are induced by the sections and the structure maps
respectively.  If the $f_1$ is a $G_1$-ex-space and $f_2$ is a
$G_2$-ex-space, then $f_1 \barwedge f_2$ has a $(G_1 \times
G_2)$-action.

If $f_1$ and $f_2$ are parametrized orthogonal spectra, we form the
external smash product $f_1 \barwedge f_2$ using the left Kan
extension and the external smash product of spaces, just as we defined
the external smash product of orthogonal spectra.  If $f_1$ is an
orthogonal $G_1$-spectrum and $f_2$ is an orthogonal $G_2$-spectrum,
then $f_1 \barwedge f_2$ is a $(G_1 \times G_2)$-spectrum.
\end{defin}

We now turn to a concise statement of the results on the homotopy
theory of parametrized spaces and spectra that we need.  The
construction and use of model structures on parametrized spectra is
technical and significantly more difficult than in the absolute case.
We first record a very convenient condition on a parametrized space,
the notion of an ex-fibration (see~\cite[8.1]{MaySigurdsson2006}).

\begin{defin}
An ex-space $B \to X \to B$ is an ex-fibration if the structure map $X
\to B$ be a Hurewicz fibration and the section $B \to X$ is a
closed inclusion such that there exists a retraction $X \times I \to X
\coprod_B (B \times I)$ (over $B$).
\end{defin}

We now turn to the model structures on parametrized spectra.

\begin{defin}
The {\em level model structure} on $\Sp^U_{/ B, G}$ is specified by
taking the fibrations to be the fiberwise weak equivalences and the 
fibrations to be the fiberwise fibrations of orthogonal $G$-spectra.
\end{defin}

We define the stable homotopy groups of a parameterized $G$-spectrum
to be the collection of stable homotopy groups of the fibers.  A
fiberwise stable equivalence $f \colon X \to Y$ is a map of
parametrized $G$-spectra that induces an equivalence on stable
homotopy groups after fibrant replacement in the level model
structure.  

\begin{defin}
The {\em stable model structure} on $\Sp^U_{/ B, G}$ is specified by
taking the weak equivalences to be the fiberwise stable equivalences.
The fibrant objects are fibrant in the level model structure and are
$\Omega$-spectra, in that the adjoint structure maps are fiberwise
stable equivalences.
\end{defin}

For computing right-derived functors, one key construction we need is
a fiberwise lax monoidal fibrant replacement functor for parametrized
spectra.  We do this using a fiberwise analogue of the construction of
the non-parametrized lax monoidal fibrant replacement functor given in
Definition~\ref{def:lax-fibrant}.

\begin{defin}\label{def:lax-fibrant-param}
For a universe $U$ and an orthogonal parametrized $G$-spectrum $X$
over a base $B$, define the orthogonal spectrum $\cQ_U X$ via the
assignment
\begin{equation}
X^{\mfib}_B(V) = \cQ_U X(V) = \hocolim_{W \in U} \Omega_B^{W \otimes
  V} X_B((W \oplus \bR) \otimes V).
\end{equation}
Here the homotopy colimit is indexed over the poset under inclusion of
finite dimensional subspaces of $U$.  There is a natural
transformation $\id \to (-)^{\mfib}$ induced by the inclusion of $X_B(V)
= \Omega_B^{0 \otimes V} X_B(\bR \otimes V)$.
\end{defin}

This functor is lax monoidal with respect to the external product.

\begin{prop}\label{prop:Q-lax-mon-param}
There is a map of orthogonal $G_1 \times G_2$-spectra over $B_1 \times
B_2$
\begin{equation}
\cQ_{U_1} X \sma \cQ_{U_2} Y \to \cQ_{U_1 \oplus U_2} (X \sma Y)
\end{equation}
which is associative and unital.
\end{prop}

In contrast with the non-parametrized setting, the construction of
Definition~\ref{def:lax-fibrant-param} does not necessarily produce a
fibrant parametrized orthogonal spectrum.  This is true however when
the input is already fibrant in the level model structure; moreover,
if the input is ex-fibrant, so is the output of this process.
Furthermore, given a parametrized spectrum $X$, it follows
from~\cite[5.2.5]{Malkiewich2019} that one can construct a suitable
replacement of $X$ with the property that all of the constituent
spaces are ex-fibrations.

\subsubsection{Spaces and spectra of sections}

Starting in Section~\ref{sec:comp-coch-kuran}, we defined the twisted
cochains in terms of a spectrum of sections with compact support.  In
this section, we review the theory of spaces and spectra of sections.

\begin{defin}
The (based) space of sections of an ex-space
$X \to B$ with basepoint section $s$ is defined as 
\begin{equation}
\Gamma_B(X) = \Map_B(B,X),
\end{equation}
i.e., maps over $B$ from $B$ to $X$, with basepoint $s$.  This can be
described as $p_* X$, where $p \colon X \to *$ is the canonical map
and $p_*$ participates in the adjunction $(p^*, p_*)$.  If $X$ is a
parametrized $G$-space, $\Gamma_B(X)$ is a $G$-space via the
conjugation action.
\end{defin}

The construction of the space of sections is compatible with the
external smash product.

\begin{prop}
For ex-spaces $X_1 \to B_1$ and $X_2 \to B_2$, there is a natural map 
\begin{equation}
\Gamma_{B_1}(X_1) \sma \Gamma_{B_2}(X_2) \to \Gamma_{B_1 \times
  B_2}(X_1 \barwedge X_2).
\end{equation}
These maps are associative and unital.  When $X_1$ is a parametrized
$G_1$-space and $X_2$ is a parametrized $G_2$-space, the natural map
is a map of $(G_1 \times G_2)$-spaces.
\end{prop}

We will be interested in sections with compact support.  As is often
the case with mapping constructions in parametrized spaces, the
point-set topology requires a bit of care.  The required discussion is
folklore, but is not in the classical literature in a convenient way;
see~\cite[\S A.2]{Malkiewich2014} for a recent treatment.

\begin{defin}\label{defin:space-of-sections}
For an ex-space $X \to B$ with basepoint section $s$, we
define the following spaces of sections:
\begin{itemize}
\item The space of sections with compact support 
\begin{equation}
\Gamma^c_B(X) \subseteq \Map_B(B,X)
\end{equation}
consisting of those sections which vanish outside a compact subset of
$B$.

\item For a closed subset $A \subset B$, we define the subspace
\begin{equation}
\Gamma^{c}_{B,A}(X) \subseteq \Map_B(B,X)
\end{equation}
to consist of those sections that vanish off $A$.  When $s \in
\Gamma^{c}_{B,A}(X)$, this is a based space.

\item For $Y \subset X$, we define the compactly supported sections
relative to $Y$, denoted 
\begin{equation}
\Gamma^{c}_B(X,Y) \subseteq \Map_B(B,X),
\end{equation}
to be the space of sections that land in $Y$ off a compact subset of
$B$.
\end{itemize}
The topology on these spaces is described below.
\end{defin}

When $B$ is compact, we can use the subspace topology induced from the
topology on $\Map_B(B,X)$.  If $B$ is not compact, we consider the
natural inclusion $B \to B^+$, where $B^+$ denotes the one-point
compactification.  We can form the pushout $X \cup_{B} B^+$ using the
section $s \colon B \to X$.  Now we give $\Gamma^c_B(X)$ the subspace
topology as a subspace of $\Map_{B^+}(B^+, X \cup_{B} B^+)$; sections
with compact support over $B$ extend to sections over $B^+$.  The
topologies on $\Gamma^{c,A}_B(X)$ and $\Gamma^c_B(X,Y)$ are defined
analogously.

\begin{prop}\label{prop:multiplicative-compact-sections-spaces}
For parametrized spaces $X_1$ and $X_2$ over base spaces $B_1$ and
$B_2$ respectively, there is a natural map
\begin{equation}
\Gamma^{c}_{B_1}(X_1) \sma \Gamma^{c}_{B_2}(X_2) \to \Gamma^{c}_{B_1 \times
  B_2}(X_1 \barwedge X_2).
\end{equation}
For $Y_1 \subset X_1$ and $Y_2 \subset X_2$, there are natural maps 
\begin{equation}
\Gamma^{c}_{B_1}(X_1, Y_1) \sma \Gamma^{c}_{B_2}(X_2, Y_2) \to
\Gamma^{c}_{B_1 \times B_2}(X_1 \barwedge X_2, Y_1 \times Y_2).
\end{equation}
These maps are associative and unital.  When $X_1$ is a parametrized
$G_1$-space and $X_2$ is a parametrized $G_2$-space, these are maps of
$(G_1 \times G_2)$-spaces.
\end{prop}

We now turn to the analogous definitions for parametrized
spectra.  

\begin{defin}\label{defin:spectrum-of-sections}
For a parametrized spectrum $X$, the parametrized spectrum of sections
$\Gamma_B(X)$ is defined by the spacewise application of the space of
sections:
\begin{equation}
(\Gamma_B(X))_V = \Gamma_B(X_V).
\end{equation}
Equivalently, this can again be described as $p_* X$ where $p \colon B
\to \ast$ is the canonical map.  The parametrized spectrum of
compactly supported sections and compactly supported sections relative
to a closed subspace $A \in B$ are defined analogously.
\end{defin}

The fact that spacewise application of these functors on ex-spaces
yields functors on orthogonal spectra requires verification of
compatibility with suspension; see for
example~\cite[11.4.1]{MaySigurdsson2006} for this kind of argument.

Proposition~\ref{prop:multiplicative-compact-sections-spaces} has the
following analogue in the stable setting.

\begin{prop}\label{prop:multiplicative-compact-sections-spectra}
For parametrized spectra $X_1$ and $X_2$ over base spaces $B_1$ and
$B_2$ respectively, there is a natural map
\begin{equation}
\Gamma^{c}_{B_1}(X_1) \sma \Gamma^{c}_{B_2}(X_2) \to \Gamma^{c}_{B_1 \times
  B_2}(X_1 \barwedge X_2).
\end{equation}
For $Y_1 \subset X_1$ and $Y_2 \subset X_2$, there are natural maps 
\begin{equation}
\Gamma^{c}_{B_1}(X_1, Y_1) \sma \Gamma^{c}_{B_2}(X_2, Y_2) \to
\Gamma^{c}_{B_1 \times B_2}(X_1 \barwedge X_2, Y_1 \sma Y_2).
\end{equation}
These maps are associative and unital.  
\end{prop}

We record the interaction of the cotensor with the spectrum of
sections.

\begin{lemma}\label{lem:cotensor-inside}
Let $A$ be a $G$-space and $X$ a parametrized orthogonal $G$-spectrum
over $B$.  Then there is an adjunction homeomorphism
\begin{equation}
F(A, \Gamma^c_{B}(X)) \cong \Gamma^c_{B}(F_B(A,X)).
\end{equation} \qed
\end{lemma}

Next, we recall an equivariant interpretation of the spectrum of
sections.  Given an orthogonal $G$-spectrum $X$, the Borel
construction provides an associated parametrized spectrum $EG_+ \sma_G
X \to BG$, with fibers equivalent to $X$.

\begin{prop}\label{prop:sections-BG}
Let $X$ be an orthogonal $G$-spectrum.  Then there is a natural
equivalence
\begin{equation}
F(EG_+, X)^G \htp \Gamma_{BG}(EG_+ \sma_G X).
\end{equation} \qed
\end{prop}

We now recall how to compute derived functors of spaces and spectra of
sections using the point-set models.  In our context, the basic slogan
is that when working with ex-fibrations, the various spaces of
sections are homotopically well-behaved.  For compactly supported
sections, correct statements require some additional hypotheses to
control the point-set topology of the base spaces, mostly surrounding
normality hypotheses to permit arguments involving the construction of
homotopies.  For simplicity, we state the results for the absolute
cases of the spaces and spectra of sections and compactly supported
sections; we leave the relative cases to the reader.

\begin{lem}
Let $X \to Y$ be a fiberwise weak equivalence of ex-fibrations over a
locally compact Hausdorff space $B$.  Then the natural maps  
\begin{equation}
\Gamma_B(X) \to \Gamma_B(Y) \quad\textrm{and}\quad \Gamma_B^{c}(X) \to \Gamma_B^{c}(Y)
\end{equation}
are weak equivalences of spaces.

Let $X \to Y$ be a stable equivalence of parametrized spectra that are
spacewise ex-fibrations over a locally compact Hausdorff space $B$.  Then the natural maps
\begin{equation}
\Gamma_B(X) \to \Gamma_B(Y) \quad\textrm{and}\quad \Gamma_B^{c}(X) \to \Gamma_B^{c}(Y)
\end{equation}
are stable equivalences of spectra.
\end{lem}

Next, we turn to the question of the interaction of spaces of
sections with homotopy colimits in the category of parametrized
spectra.  We have the following excision result.

\begin{thm}
Suppose that $U$ and $V$ are normal subsets of a locally compact
Hausdorff space $B$.  Then there are homotopy pullback squares of
spaces 
\begin{equation}
\begin{tikzcd}
\Gamma_{U} (X) \ar[d] & \Gamma_{U \cup V}(X) \ar[l] \ar[d] \\
\Gamma_{U \cap V}(X) & \Gamma_{V} (X) \ar[l]
\end{tikzcd}
\end{equation}
and when $X$ is an ex-fibration there are homotopy pullback squares of spaces
\begin{equation}
\begin{tikzcd}
\Gamma_{U}^c (X) \ar[r] & \Gamma_{U \cup V}^c(X) \\
\Gamma_{U \cap V}^c(X) \ar[r] \ar[u] & \Gamma_{V}^c (X) \ar[u].
\end{tikzcd}
\end{equation}
The analogous results for parametrized spectra hold.
\end{thm}

We have the following result about filtered colimits.

\begin{prop}
Suppose we have a filtered diagram of parametrized spaces $\{f_i\}$
indexed by a filtered category $D$.  If each $f_i$ is an ex-fibration
and each space $B_i$ is normal, then the natural map
\begin{equation}
\Gamma^c_{\hocolim_D B_i}(\hocolim_D f_i) \to \hocolim_D \Gamma^c_{B_i}(f_i) 
\end{equation}
and the natural map
\begin{equation}
\Gamma^c_{\hocolim_D B_i}(\hocolim_D f_i) \to \Gamma^c_{\hocolim_D B_i}(\tilde{f})
\end{equation}
are weak equivalences, where $\tilde{f}$ is an ex-fibrant replacement
of $\hocolim_{D_i} f_i$.  The analogous result for parametrized
spectra holds.\qed
\end{prop}

In our applications, we will sometimes need to compute spaces of
sections for parametrized spaces that are quasifibrations but not
fibrations.  Recall that a map $p \colon E \to B$ is a quasifibration
if the canonical map $p^{-1}(b) \to Fp$ from the fiber to the homotopy
fiber is a weak equivalence for each $b \in B$.  Although not
explicitly phrased in terms of a lifting condition, this amounts to a
lifting criterion up to a homotopy.  The natural map from a 
quasifibration to the associated path fibration (or more generally an
ex-fibrant replacement in the model structure on parametrized spaces)
is a fiberwise weak equivalence, essentially by definition.

Quasifibrations are sufficient for computing the fiberwise derived
homotopy groups of a parametrized space or spectrum.  However, in
general, quasifibrations can not be used to compute the space of
sections; rather, it is standard to define the space of sections of a
quasifibration in terms of the associated fibration.  Nonetheless, it
is possible to give conditions under which the space of sections of a
quasifibration have the correct homotopy type.

\begin{prop}
Suppose that we have a pushout diagram of parametrized spaces
\begin{equation}
\begin{tikzcd}
(f_{01} \colon E_{01} \to B_{01}) \ar[r] \ar[d] & (f_0 \colon E_0 \to
  B_0) \ar[d] \\
(f_1 \colon E_1 \to B_1) \ar[r] & (f \colon E \to B)
\end{tikzcd}
\end{equation}
such that one of the maps $B_{01} \to B_0$ or $B_{01} \to B_1$ is a
cofibration and the maps $f_{01}$, $f_0$, and $f_{10}$ are
ex-fibrations.

Then the map $f$ is an ex-quasifibration and the natural map
\begin{equation}
\Gamma_B(E) \to \Gamma_B(\tilde{E}) 
\end{equation}
is a weak equivalences, where $\tilde{f} \colon \tilde{E} \to B$ is a
fibration equivalent to $f$.  If $B_{01}$, $B_0$, and $B_1$ are normal
spaces then the natural map
\begin{equation}
\Gamma^c_B(E) \to \Gamma^c_B(\tilde{E}) 
\end{equation}
is a weak equivalence.
\end{prop}

\begin{proof}
The fact that $f$ is an ex-quasifibration is essentially classical;
see e.g.,~\cite{Puppe1974}.  To see that $f$ has homotopically
meaningful sections, we argue as follows.  A section of $f$ is
precisely a pair $(s_0, s_1)$ of sections of $f_0$ and $f_1$ that
coincide on $B_{01}$; i.e., the space of sections $\Gamma_B(f)$ is the
pullback 
\begin{equation}\label{eq:section-pullback}
\Gamma_B(f) \cong \Gamma_{B_0}(f_0) \times_{\Gamma_{B_{01}}(f_{12})}
\Gamma_{B_1}(f_1).
\end{equation}
Without loss of generality, we assume that $B_{01} \to B_0$ is a
cofibration.  Then the map of section spaces $\Gamma_{B_0}(f_0) \to
\Gamma_{B_{01}}(f_{01})$ is a fibration, and so the pullback in
Equation~\eqref{eq:section-pullback} is a homotopy pullback.  By the
gluing lemma, this implies the natural map $\Gamma_B(f) \to
\Gamma_B(\tilde{f})$ is a weak equivalence.  For the case of compactly
supported sections, the normality assumption allows us to run an
essentially analogous argument.
\end{proof}

More generally, we have the following statement, which is a
consequence of the previous results by the standard induction over
homotopy pushouts along latching maps.

\begin{prop}\label{prop:gluing-good-quasifibrations}
Suppose that we have a homotopy colimit diagram of parametrized spaces
$\{f_i \colon E_i \to B_i\}$ indexed by a discrete Reedy category $D$.
If each $f_i$ is an ex-fibration and each space $B_i$ is normal, then
the natural map
\begin{equation}
\Gamma_{\hocolim_D B_i}^c(\hocolim_D f_i) \to \holim_D
\Gamma^c_{B_i}(f_i)
\end{equation}
is a weak equivalence and the natural map
\begin{equation}
\Gamma_{\hocolim_D B_i}^c(\hocolim_D f_i) \to \holim_D
\Gamma^c_{\hocolim_D B_i}(\tilde{f})
\end{equation}
is a weak equivalence, where $\tilde{f}$ is an ex-fibrant replacement
of $\hocolim_{D_i} f_i$.  The analogous result for parametrized
spectra holds.\qed
\end{prop}

\subsection{Multiplicative orientations and the equivariant complex
  cobordism spectrum}
\label{sec:multiplicative-orientations}
\label{sec:eqMU}

The purpose of this section is to establish the facts we need from the
theory of equivariant multiplicative complex orientations of vector
bundles.  Although the general theory of equivariant orientations
is very complicated, the situation for equivariant complex bundles
(i.e., $(G, U(n))$-bundles, where $G$ acts by complex bundle
maps) is much simpler.  We begin by describing a specific model of the
classifying space for equivariant complex bundles that encodes the
multiplicative structure, which arises from the theory of diagrammatic
spaces.  See for example~\cite{BlumbergCohenSchlichtkrull,
  Schlichtkrull2009, SagaveSchlichtkrull2012, SagaveSchlichtkrull2019,
Schwede2018book} for more detailed treatments of this theory.

\subsubsection{Grassmanians and diagram spaces}

Consider a complete complex $G$-universe $U$.  Often we will take 
\begin{equation}
U = \tilde{U} \otimes_{\bR} \bC
\end{equation} 
for a complete real universe $\tilde{U}$.  

\begin{defin}
Let $\sI_G^U$ denote the category of finite-dimensional complex vector
spaces in $U$ and complex isometries.  This category is symmetric
monoidal under direct sum, with unit the vector space $\{0\}$.  
\end{defin}

Moreover, there is the following external product:

\begin{lem}\label{lem:vector-space-associativity}
For groups $G_1$ and $G_2$ and complex universes $U_1$ and $U_2$,
respectively, there are external product maps 
\begin{equation}
\sI_{G_1}^{U_1} \times \sI_{G_2}^{U_2} \to \sI_{G_1 \times G_2}^{U_1
  \oplus U_2}.
\end{equation}
These maps are associative in the sense that the diagrams
\begin{equation}
\begin{tikzcd}
\sI_{G_1}^{U_1} \times \sI_{G_2}^{U_2} \times \sI_{G_3}^{U_3} \ar[r]
\ar[d] & \sI_{G_1 \times G_2}^{U_1 \oplus U_2} \times \sI_{G_3}^{U_3}
\ar[d] \\ \sI_{G_1}^{U_1} \times \sI_{G_2 \times G_3}^{U_2 \oplus U_3}
\ar[r] & \sI_{G_1 \times G_2 \times G_3}^{U_1 \oplus U_2 \oplus U_3}
\end{tikzcd}
\end{equation}
commute. \qed
\end{lem}

Note that $U_1 \oplus U_2$ is not typically a complete $(G_1 \times
G_2)$-universe even if we assume that $U_1$ and $U_2$ are complete
universes.  However, for a complete universe $U_{12}$, there will 
be a contractible space of isometries $U_1 \oplus U_2 \to U_{12}$.  An
arbitrary isometry will not necessarily compatible with the symmetric
monoidal structure, but sometimes the system of isometries can be made
compatible with direct sum.  For example, as discussed in
Lemma~\ref{lem:mult-universe-embedding}, we can construct such a
monoidal system when working with complete universes obtained as
direct sums of the regular representation.  Recall that for any group
$G$ we write $\cU_{\bC}(G) = \rho_G \otimes \bC^{\infty}$ for the
$G$-universe which is the infinite direct sum of copies of the regular
representation equipped with a natural inner product.

\begin{cor}\label{cor:vector-space-associativity}
For groups $G_1$ and $G_2$ and complex universes $U_1$ and $U_2$,
respectively, there are external product maps 
\begin{equation}
\sI_{G_1}^{U_1} \times \sI_{G_2}^{U_2} \to \sI_{G_1 \times G_2}^{U_1
  \oplus U_2}.
\end{equation}
These maps are associative in the sense that the diagrams
\begin{equation}
\begin{tikzcd}
\sI_{G_1}^{\cU_{\bC}(G_1)} \times \sI_{G_2}^{\cU_{\bC}(G_2)} \times
\sI_{G_3}^{\cU_{\bC}(G_3)} \ar[r] \ar[d] & \sI_{G_1 \times
  G_2}^{\cU_{\bC}(G_1 \times G_2)} \times \sI_{G_3}^{\cU_{\bC}(G_3)} 
\ar[d] \\ \sI_{G_1}^{\cU_{\bC}(G_1)} \times \sI_{G_2 \times
  G_3}^{\cU_{\bC}(G_2 \times G_3)}
\ar[r] & \sI_{G_1 \times G_2 \times G_3}^{\cU_{\bC}(G_1 \times G_2
  \times G_3)}
\end{tikzcd}
\end{equation}
commute. \qed
\end{cor}

We will work with a model of spaces defined in terms of diagrams on
these indexing categories, and construct a model of the cartesian
product that has the property that commutative monoids encode
$E_\infty$ spaces.

\begin{defin}
An $\sI_G^U$-space is a continuous functor from $\sI^U_G$ to
spaces.  We denote the category with objects the $\sI_G^U$-spaces and
morphisms the natural transformations by $\sI_G^U-\Top$.
\end{defin}

The Day convolution product (regarding the category of
spaces as symmetric monoidal under the cartesian product) makes the
category of $\sI^U_G$-spaces into a symmetric monoidal category.  The
monoidal unit is given by the constant $\sI^U_G$-space at a point, and
the symmetric monoidal product can be described by the formula 
\begin{equation}
(X \boxtimes Y)(V) = \colim_{W_1 \oplus W_2 \cong V} X(W_1) \times
Y(W_2).
\end{equation}
By construction, the universal property of this monoidal product is
that a map $X \boxtimes Y \to Z$ of $\sI^U_G$-spaces is specified by
the data of maps 
\begin{equation}
X(W_1) \times Y(W_2) \to Z(W_1 \oplus W_2).
\end{equation} 

There is also an external product on $\sI_G$-spaces, defined as follows.

\begin{defin}
Fix groups $G_1$ and $G_2$ and corresponding universes $U_1$ and
$U_2$.  Let $X$ be an $\sI_{G_1}^{U_1}$-space and $Y$ be an
$\sI_{G_2}^{U_2}$-space.  Then the external product
$\bar{\boxtimes}$ is an $\sI_{G_1 \times G_2}^{U_1 \oplus
  U_2}$-space defined via the assignment 
\begin{equation}
(V_1, V_2) \mapsto X(V_1) \times Y(V_2).
\end{equation}
and the left Kan extension.
\end{defin}

The external multiplication naturally arises from the pairing 
\begin{equation}
\left( \sI^{U_1}_{G_1}-\Top \right)\times \left(\sI^{U_2}_{G_2}-\Top\right) \to
\sI^{U_1 \oplus U_2}_{G_1 \times G_2}-\Top,
\end{equation}
and so Lemma~\ref{lem:vector-space-associativity} implies that the
external product is associative.  Morever, given an isometry $U_1
\oplus U_2 \to U_{12}$, we can take the left Kan extension to produce
a functor  
\begin{equation}
\sI^{U_1 \oplus U_2}_{G_1 \times G_2}-\Top \to \sI^{U_{12}}_{G_1
  \times G_2}.
\end{equation}
The left Kan extension will take a monoidal system of diagrams to a
monoidal external product of diagram categories, and so using
Corollary~\ref{cor:vector-space-associativity}, we can arrange for
such an external product for complete universes when working with the
regular representation universes.

We now turn to the relationship between the category $\sI_G^U-\Top$
and the category of $G$-spaces.

\begin{lem}
The homotopy colimit induces a lax monoidal functor
\begin{equation}
\hocolim_{\sI_G^U} \colon \sI_G^U-\Top \to \cT^G
\end{equation}
\end{lem}

There is a model structure on $\sI_G^U-\Top$ in which the weak
equivalences are detected by equivalences after passage to the
homotopy colimit; this is obtained by localizing the obvious levelwise
model structure, i.e., the projective model structure, at these maps.
We do not require any details about these model structures, so we do
not discuss this further herein.

Note that the homotopy colimit is not a symmetric monoidal functor,
although it does preserve suitable $E_\infty$ operad actions.  It is
however monoidal, which suffices for our applications.  Another key
aspect of the situation is that there is a homotopical right adjoint
to the homotopy colimit functor which has good point-set properties
and in particular is also a monoidal functor, which we will denote by
$R_G^U$ (e.g., see~\cite[\S 4.2]{Schlichtkrull2009} for a discussion
in the symmetric case).

Both the homotopy colimit and the families $R_{(-)}^{(-)}$ of functors
are compatible with the external product.

\begin{prop}
For groups $G_1$ and $G_2$ and universes $U_1$ and $U_2$ respectively,
the diagrams
\begin{equation}
\begin{tikzcd}
(\sI_{G_1}^{U_1}-\Top) \times (\sI_{G_2}^{U_2}-\Top)
  \ar[rrrr,"\hocolim_{\sI_{G_1}^{U_1}} \times \hocolim_{\sI_{G_2}^{U_2}}"] \ar[d] &&&&
\cT^{G_1} \times \cT^{G_2} \ar[d] \\
\sI_{G_1 \times G_2}^{U_1 \oplus U_2}-\Top \ar[rrrr,swap,"\hocolim_{\sI_{G_1
      \times G_2}^{U_1 \oplus U_2}}"] &&&& \cT^{G_1\times G_2}
\end{tikzcd}
\end{equation}
and
\begin{equation}
\begin{tikzcd}
\cT^{G_1} \times \cT^{G_2} \ar[d] \ar[rrr,"R_{G_1}^{U_1} \times
  R_{G_2}^{U_2}"] &&& (\sI_{G_1}^{U_1}-\Top) \times
(\sI_{G_2}^{U_2}-\Top) \ar[d] \\ 
\cT^{G_1\times G_2} \ar[rrr,swap,"R_{G_1 \times G_2}^{U_1
      \oplus U_2}"] &&& \sI_{G_1 \times G_2}^{U_1 \oplus U_2}-\Top 
\end{tikzcd}
\end{equation}
commute.
\end{prop}

Note that the homotopy colimit is also compatible with the left Kan
extension along an isometry $\sI_{G_1}^{U_1} \times \sI_{G_2}^{U_2} \to
\sI_{G_1 \times G_2} U_{12}$; there are induced maps
\begin{equation}
\hocolim_{\sI_{G_1}^{U_1}} X \times \hocolim_{\sI_{G_1}^{U_1}} Y \to 
\hocolim_{\sI_{G_1 \times G_2}^{U_1 \oplus U_2}} X \boxtimes Y \to
\hocolim_{\sI_{G_1 \times G_2}^{U_{12}}} X \boxtimes Y,
\end{equation}
which are associative if the system of isometries has been chosen to
be associative.  In particular, we have the following corollary:

\begin{cor}
For groups $G_1$ and $G_2$ and universes $U_1$ and $U_2$ respectively,
the diagram 
\begin{equation}
\begin{tikzcd}
(\sI_{G_1}^{\cU(G_1)}-\Top) \times (\sI_{G_2}^{\cU(G_2)}-\Top)
  \ar[rrrr,"\hocolim_{\sI_{G_1}^{\cU(G_1)}} \times \hocolim_{\sI_{G_2}^{\cU(G_2)}}"] \ar[d] &&&&
\cT^{G_1} \times \cT^{G_2} \ar[d] \\
\sI_{G_1 \times G_2}^{\cU(G_1 \times G_2)}-\Top \ar[rrrr,swap,"\hocolim_{\sI_{G_1
      \times G_2}^{\cU(G_1 \times G_2)}}"] &&&& \cT^{G_1\times G_2}
\end{tikzcd}
\end{equation}
commutes.
\end{cor}

The point of this framework is that constructions of the Thom
spectrum functor and the theory of orientations that are strictly
multiplicative are easiest to formulate in terms of such diagramatic
models of spaces.  We now give models of the classifying space of
complex vector bundles in this formalism.  Our treatment is derived
from the pioneering approach of~\cite[\S X]{LewisMaySteinberger1986},
and influenced by the modern adaptation in~\cite[\S
  6.1]{Schwede2018book}.

\begin{defin}
Let $V$ be a finite-dimensional complex vector space in a complex
universe $U$.  We denote by $BU_G(q, V \otimes U)$ the $G$-space
of complex $q$-dimensional planes in $V \otimes U$.  Let $EU_{G}(q,V
\otimes U)$ denote the tautological bundle over $BU_G(q, V \otimes U) $, i.e. the $G$-space of pairs $(Z,z)$ with $Z \in BU_G(q,V
\otimes U)$ and $z \in Z$.  There is a natural map
\begin{equation}
\begin{aligned}
EU_G(q,V \otimes U) &\to BU_G(q,V \otimes U) \\ 
(Z,z) &\mapsto Z.
\end{aligned}
\end{equation}
We define $\sI^U_G$-spaces via the formulas
\begin{equation}
BU_G(V) = \coprod_q BU_G(q,V \otimes U) \qquad\textrm{and}\qquad
EU_G(V) = \coprod_q EU_G(q,V \otimes U).
\end{equation}
\end{defin}

\begin{rem}
In the rest of the paper, we have used the notation $EG$ to denote a
specific construction of a free $G$-space whose underlying space is
contractible, namely the two-sided bar construction $B(G,G,\ast)$.
While our notation for the the total space of the tautological bundle
over $BU_G(q,V \otimes U)$ is potentially confusing, it is consistent
with the standard literature on models of the cobordism spectra.
\end{rem}

The standard isomorphisms 
\begin{equation}
(V_1 \otimes U) \oplus (V_2 \otimes U) \cong (V_1
\oplus V_2) \otimes U
\end{equation}
specified by 
\begin{equation}
(v_1 \otimes u), (v_2 \otimes u') \mapsto (v_1, v_2) \otimes (u, u')
\end{equation}
induce a multiplicative structure on these spaces, as follows. 

\begin{lemma}
The $\sI^U_G$-spaces $BU_G(-)$ and $EU_G(-)$ are commutative monoids
and the projection map $EU_G(-) \to BU_G(-)$ is a map of monoids (and
also a vector bundle with varying fibers, i.e., a vector bundle at
each point of $\sI^U_G$).  \qed 
\end{lemma}

The corresponding homotopy colimits are topological monoids
in spaces; we write these as follows.
\begin{equation}
BU_G(U) = \hocolim_{V \in \sI^U_G} BU_G(V) \qquad\textrm{and}\qquad EU_G(U) =
\hocolim_{V \in \sI^U_G} EU_G(V).
\end{equation}

As $G$ varies, the Grassmanians are compatible with the external
multiplication.

\begin{lem}
There is a natural external multiplicative structure on $BU_G{(-)}$
induced by the maps 
\begin{equation}
BU_{G_1}(V_1) \times BU_{G_2}(V_2) \to BU_{G_1 \times G_2}(V_1 \oplus
V_2).
\end{equation}
This external multiplication is associative and unital in the evident
sense.  Moreover, $BU_{G_1 \times G_2}(-)$ is a
$BU_{G_1}(-)$-$BU_{G_2}(-)$ bimodule and the external multiplication
is a map of bimodules.
\end{lem}

\begin{proof}
The external multiplication is induced by the maps
\begin{equation}
\begin{tikzcd}
BU_{G_1}(|V_1|,V_1 \otimes U_1) \times BU_{G_2}(|V_2|, V_2 \otimes U_2)
\ar[d] \\
BU_{G_1 \times G_2}(|V_1|+|V_2|, (V_1 \oplus V_2) \otimes (U_1 \oplus U_2))
\end{tikzcd}
\end{equation}
induced by the direct sum of planes and the inclusions 
\begin{equation}
U_1 \to U_1 \oplus U_2 \quad\textrm{and}\quad U_2 \to U_1 \oplus U_2.
\end{equation}
The bimodule structure on
$BU_{G_1\times G_2}(-)$ is induced by the maps 
\begin{equation}
\begin{tikzcd}
BU_{G_1}(|V_1|,V_1 \otimes U_1) \times BU_{G_1 \times G_2}(|V_{12}|,
V_{12} \otimes (U_1 \oplus U_2))
\ar[d] \\
BU_{G_1 \times G_2}(|V_1|+|V_{12}|, ((V_1 \oplus \{0\}) \oplus V_{12}) \otimes (U_1
\oplus U_2)
\end{tikzcd}
\end{equation}
and the corresponding maps on the other side.  To see that the
external product is a bimodule map, we need to verify that the diagram
\begin{equation}
\begin{tikzpicture}
\node[scale=.7]{
\begin{tikzcd}
BU_{G_1}(|V_1|, V_1 \otimes U_1) \times BU_{G_1}(|V'_1|,V'_1 \otimes
U_1) \times BU_{G_2}(|V_2|, V_{2} \otimes U_2) \ar[r] \ar[d] & 
BU_{G_1}(|V_1|+|V'_1|, (V_1 \oplus V_1') \otimes U_1) \times
BU_{G_2}(|V_2|, V_{2} \otimes U_2) \ar[d] \\
BU_{G_1}(|V_1|, V_1 \otimes U_1) \times BU_{G_1 \times
  G_2}(|V'_1|+|V_2|,(V_1' \oplus V_2) \otimes (U_1 \oplus U_2)) \ar[r] & BU_{G_1
 \times G_2}(|V_1|+|V_1'|+|V_2|, ((V_1 \oplus V_1') \oplus V_2) \otimes (U_1
\oplus U_2)).
\end{tikzcd}
};
\end{tikzpicture}
\end{equation}
commutes.  This amounts to checking that the diagram of vector spaces
\begin{equation}
\begin{tikzcd}
(V_1 \otimes U_1) \times (V'_1 \otimes U_1) \times (V_{2} \otimes U_2) \ar[r] \ar[d] & 
(V_1 \oplus V_1') \otimes U_1) \times (V_{2} \otimes U_2) \ar[d] \\
(V_1 \otimes U_1) \times ((V_1' \oplus V_2) \otimes (U_1 \oplus U_2))
  \ar[r] & ((V_1 \oplus V_1') \oplus V_2) \otimes (U_1 \oplus U_2))
\end{tikzcd}
\end{equation}
commutes, which can be verified by chasing elements.
\end{proof}

Since the homotopy colimit is a monoidal functor, the external
multiplication passes to the associated classifying spaces.

\begin{prop}
On passage to homotopy colimits, we obtain associative and unital
systems of products 
\begin{equation}
BU_{G_1}(U_1) \times BU_{G_2}(U_2) \to BU_{G_1 \times G_2}(U_1 \oplus U_2).
\end{equation}
that are associative in the sense that the diagrams
\begin{equation}
\begin{tikzcd}
BU_{G_1}(U_1) \times BU_{G_2}(U_2) \times BU_{G_3}(U_3) \ar[r] \ar[d]
& BU_{G_1 \times G_2}(U_1 \oplus U_2) \ar[d] \\
BU_{G_1}(U_1) \times BU_{G_2 \times G_3}(U_2 \oplus U_3) \ar[r] &
BU_{G_1 \times G_2 \times G_3}(U_1 \oplus U_2 \oplus U_3)
\end{tikzcd}
\end{equation}
commute.  The monoid $BU_{G_1 \times G_2}(U_1 \oplus U_2)$ is a
$BU_{G_1}(U_1)$-$BU_{G_2}(U_2)$ bimodule, and these external products
are maps of bimodules.\qed
\end{prop}

We now consider the categories of spaces over $BU_G(-)$ and $BU_G(U)$;
these inherit products from the multiplicative structures on $BU_G(-)$.

\begin{defin}
Let $G$ be a finite group and $U$ a universe.  We define $\sI_G^U-\Top
/ BU_G(-)$ to be the category of $\sI_G^U$-spaces
over $BU_G(-)$.  This category inherits an ``internal'' multiplicative
structures from the commutative $\sI$-monoid structure on $BU_G(U)$,
i.e., the product is specified by the maps
\begin{equation}
X(V) \times Y(W) \to BU_{G}(V) \times BU_G(W) \to BU_G(V \oplus
W).
\end{equation}
There is an external multiplicative structures via the evident maps  
\begin{equation}
X \boxtimes Y \to BU_{G_1}(U_1) \boxtimes BU_{G_2} (U_2)\to BU_{G_1 \times
  G_2}(U_1 \oplus U_2),
\end{equation}
where $X$ is an object of $\sI_{G_1}^{U_1}$ and $Y$ is an object of
$\sI_{G_1}^{U_1}$.  (This can also be described levelwise in terms of
the cartesian product, via the universal property of $\boxtimes$.)

Analogously, we have the categories $\cT^G/ BU_G(U)$, which have
analogous internal and external multiplicative structures.  The latter
is specified by the maps  
\begin{equation}
X \times Y \to BU_{G_1}(U_1) \times BU_{G_2}(U_2) \to BU_{G_1 \times
  G_2}(U_1 \oplus U_2).
\end{equation}
\end{defin}

The functors $R_U^G$ induce multiplicative comparisons between these
categories as follows.

\begin{thm}
The functor $R_G^U$ induces a monoidal functor
\begin{equation}
R_G^U \colon \GTop / BU_G(U) \to \sI_G^U-\Top / BU_G(-)
\end{equation}
that is compatible with the external multiplicative structures in the
sense that the diagrams
\begin{equation}
\begin{tikzpicture}
\node[scale=.85]{
\begin{tikzcd}
\cT^{G_1} / BU_{G_1}(U_1) \times \cT^{G_2} / BU_{G_2}(U_2) \ar[r] \ar[d]
& \cT^{G_1 \times G_2} \Top / BU_{G_1 \times G_2}(U_1 \oplus U_2) \ar[d] \\
\sI_{G_1}^{U_1}-\Top / BU_{G_1}(-) \times \sI_{G_2}^{U_2}-\Top /
BU_{G_2}(-) \ar[r] & \sI_{G_1 \times G_2}^{U_1 \oplus U_2}-\Top /
BU_{G_1 \times G_2}(-).  
\end{tikzcd}
};
\end{tikzpicture}
\end{equation}
commute and the evident associativity diagrams also commute. \qed
\end{thm}

This shortcut allows to easily apply multiplicative orientation theory
to space-level data, provided that we are willing to ignore the
symmetric aspect of the monoidal structure.  Since in our current
application we only have an associative orientation, this is no real
limitation.

\subsubsection{Complex orientations and trivializations}

We now turn to discuss the application of complex orientations to the
trivialization of complex bundles.  We begin with a brief discussion
of how trivialization relates to orientations from the perspective of
parametrized spectra.  (See~\cite[20.5]{MaySigurdsson2006} for a more
detailed exposition of this perspective.)  For a space $B$, denote by
$S^n_{B}$ the ex-space over $B$ with total space $S^n \times
B$, fiber $S^n$, and basepoint section given by the canonical
basepoint of $S^n$.  Let $\bbk$ denote an associative ring spectrum.
As explained in~\cite[20.5.5]{MaySigurdsson2006}, a $\bbk$-orientation
of a spherical fibration $f \colon E \to B$ with fiber $S^n$ gives
rise to a trivialization in the form of an equivalence of parametrized
spectra
  \begin{equation}
E \sma \bbk \to S^n_{B} \sma \bbk.     
  \end{equation}
The argument from this perspective is quick and insufficiently known,
so we explain it here.  Let $Mf$ denote the Thom space of $f$.  The
standard notion of a $\bbk$-orientation is given by a Thom class $\mu
\in \bbk^n(Mf)$ that is a unit when restricted to the fibers of $f$.
Now, the Thom class is represented by a map $\Sigma^{\infty} Mf \to
S^n \sma \bbk$, and by adjunction this is the same thing as a map of
parametrized spectra
\begin{equation}
\Sigma_{B}^{\infty} E \to S^n_{B} \sma \bbk.    
\end{equation}
Here we are using the fact that the Thom spectrum functor can be
described as the pushforward of a stable spherical fibration along the
map to $B \to \ast$.  Smashing with $\bbk$ on the right and composing
with the multiplication, we have a map
\begin{equation}
\Sigma_{B}^{\infty} E \sma \bbk \to S^n_{B} \sma \bbk \sma \bbk \to S^n_{B} \sma
\bbk    
\end{equation}
of parametrized spectra over $B$, and the hypothesis that the Thom
class restricts to a unit on the fibers is precisely equivalent to the
statement that this map is an equivalence.

We need a multiplicative model of this trivialization equivalence in
the equivariant context.  We will denote by $MU_G$ the homotopical
complex oriented cobordism spectrum for the group $G$ and $MUP_G$ the
periodic variant.  We work with specific models of these spectra which
we now define.

\begin{defin}
Let $U$ be a complete real universe.  The periodic equivariant complex
bordism orthogonal $G$-spectrum $MUP_G$ has $V$th space (for $V$ an
indexing space in $U$) 
\begin{equation}
MUP_G(V) = \Map(S^{iV}, TU_G(V \otimes_{\bR} \bC)),
\end{equation}
where $TU_G(-)$ denotes the Thom space of the bundle $EU_G(-) \to
BU_G(-)$.  The monoidal structure on $TU_G(-)$ endows $MUP_G$ with the
structure of a commutative ring orthogonal $G$-spectrum.
\end{defin}

\begin{rem}
The looping in the definition of $MUP_G$ ensures that we obtain an
orthogonal spectrum; we could alternatively work with a ``unitary
spectrum'' indexed on complex representations, but it is simpler not
to develop this theory explicitly.
\end{rem}

The spectrum $MUP_G$ is $\bZ$-graded (in the sense of
Section~\ref{sec:graded}), where the $n$th piece $MUP_G^n$ is built
from the the Thom spaces of the $n$th part of the corresponding
grading on $BU_G$; i.e., it is comprised of the Thom spaces associated
to the spaces $BU_G(|V|+n, V \otimes U)$, which we denote by
$TU_G^n(-)$.

\begin{lem}
The orthogonal $G$-spectrum $MUP_G$ is the underlying spectrum of a
$\bZ$-graded orthogonal $G$-spectrum $(MUP_G)^\bullet$ specified by 
\begin{equation}
MUP_G^n(V) = \Map(S^{iV}, TU_G^n(V \otimes_{\bR} \bC).
\end{equation}
The graded orthogonal $G$-spectrum $MUP_G$ is a commutative monoid in
$\bZ$-graded spectra.
\end{lem}

Evaluating at $0$ yields an orthogonal $G$-spectrum $MUP_G^0$ that is
a model for the usual equivariant complex bordism spectrum $MU_G$.

\begin{defin}
The homotopical complex equivariant bordism spectrum $MU_G$ is the
commutative ring orthogonal $G$-spectrum obtained as $MUP_G^0$.
\end{defin}

We have the following periodicity, usually referred to as ``complex
stability''; this is a form of the Thom isomorphism, specialized to
$G$-representations (regarded as trivial bundles over $\ast$).  For
any complex representation $V$, there is a map 
\begin{equation}
MUP_G \to \Sh_V MUP_G
\end{equation}
specified simply by passage to Thom spaces from the evident map 
\begin{equation}
BU_G(|W|,W \otimes U) \to BU_G(|W|+|V|, (W \oplus V) \otimes U)
\end{equation}
induced from the map $W \to W \oplus V$ specified by $w \mapsto
(w,0)$.  Here recall that $\Sh_V$ denotes the $V$-shift functor; see
Definition~\ref{defin:V-shift}.  Observe that this map is independent
of the representation $V$, in the following sense:

\begin{lem}\label{lem:equivariant-isometry}
Let $f \colon V \to V'$ be an isometric isomorphism of complex
representations of $G$.  Then the diagram
\begin{equation}
\begin{tikzcd}
MUP_G \ar[r] \ar[dr] & \Sh_V MUP_G \ar[d] \\
& \Sh_{V'} MUP_G
\end{tikzcd}
\end{equation}
commutes.  Furthermore, for $V \subset W$, the diagram
\begin{equation}
\begin{tikzcd}
MUP_G \ar[r] \ar[drr] & \Sh_{W-V} MUP_G \ar[r] & \Sh_{W-V} \Sh_{V} MUP_G \ar[d]
\\
& & \Sh_W MUP_G 
\end{tikzcd}
\end{equation}
commutes.
\end{lem}

The following proposition can be deduced
from~\cite[X.5.3]{LewisMaySteinberger1986}; in this form, it is close
to the statement~\cite[6.1.14]{Schwede2018book}, although that
discussion is for the real analogue.

\begin{prop}
Let $V$ be a finite-dimensional complex $G$-representation.  There is
a zig-zag of weak equivalences
\begin{equation}\label{eq:complex-stability}
\begin{tikzcd}
MUP_G \ar[r] & \Sh_V MUP_G & S^V \sma MUP_G, \ar[l]  
\end{tikzcd}
\end{equation}
where the lefthand map is the canonical comparison of the suspension
and shift functors.  The  map $MUP_G \to \Sh_V MUP_G$ is degree $-|V|$
with respect to the $\bZ$-grading. \qed
\end{prop}

This property is sometimes referred to as {\em complex stability} in
the literature.  The following easy lemma records a key fact about the
complex stability zig-zag.

\begin{lem}\label{lem:complex-stability-equivariant}
Let $V$ be a finite-dimensional complex $G$-representation.  The
zig-zag of weak equivalences in Equation~\eqref{eq:complex-stability}
is comprised of maps of spectra with $U(V)$ action, where $U(V)$ acts
on $S^V \sma (-)$ and $\Sh_V (-)$ through the action on $V$ and acts
trivially on $MUP_G$.
\end{lem}

\begin{proof}
The map $S^V \sma X \to \Sh_V X$ is equivariant with respect to
isometries of $V$ by construction; it is specified by the structure
map $S^V \sma X(W) \to X(V \oplus W)$, which is $O(V) \times
O(W)$-equivariant.  The map $MUP_G \to \Sh_V MUP_G$ lands at $0$ in
$V$ component, and is therefore equivariant since the $U(V)$-action on
$\Sh_V MUP_G$ fixes $0$; this is the content of
Lemma~\ref{lem:equivariant-isometry}.
\end{proof}

Restricting to $MUP_G^0$, we get a zig-zag which trivializes the
$G$-action (and the $U(V)$ action) on the suspension coordinate.

\begin{cor}\label{cor:thom-zigzag}
Let $V$ be a finite-dimensional complex $G$-representation.  There is
a zig-zag of weak equivalences
\begin{equation}
\begin{tikzpicture}
\node[scale=.8]{
\begin{tikzcd}
S^{|V|} \sma MUP_G^0 \ar[r] & \Sh_{\bR^{|V|}} MUP_G^0 & \ar[l]
MUP_G^{|V|} \ar[r] & \Sh_V MUP_G^0 & S^V \sma MUP_G^0 \ar[l]   
\end{tikzcd}
};
\end{tikzpicture}
\end{equation}
of orthogonal $G$-spectra. \qed
\end{cor}

Lemma~\ref{lem:shift-suspend-monoidal-compare} and the definition of
the map $MUP_G \to \Sh_V MUP_G$  imply that the complex stability
zig-zags of Corollary~\ref{cor:thom-zigzag} are externally
multiplicative.

\begin{prop}\label{prop:thom-multiplicative}
For $G_1$ and $G_2$ finite groups and $V_1$ and $V_2$ representations,
we have commutative diagrams
\begin{equation}
\begin{tikzcd}
S^{V_1} MUP_{G_1}^0 \sma S^{V_2} MUP_{G_2}^0 \ar[r] \ar[d] & S^{V_1 \oplus
  V_2} MUP_{G_1 \times G_2}^0 \ar[d] \\
\Sh_{V_1} MUP_{G_1}^0 \sma \Sh_{V_2} MUP_{G_2}^0 \ar[r] &
\Sh_{V_1 \oplus V_2} MUP_{G_1 \times G_2}^0
\end{tikzcd}
\end{equation}
and
\begin{equation}
\begin{tikzcd}
MUP_{G_1}^{|V_1|} \sma MUP_{G_2}^{|V_2|} \ar[r] \ar[d] & \Sh_{V_1}
MUP^0_{G_1} \sma \Sh_{V_2} MUP^0_{G_2} \ar[d] \\
MUP_{G_1 \times G_2}^{|V_1| + |V_2|} \ar[r] & \Sh_{V_1 \oplus V_2}
MUP^0_{G_1 \times G_2}.
\end{tikzcd}
\end{equation} 
\end{prop}

We now discuss the generalization to arbitrary complex bundles.  To
discuss this, we need to review the construction of the Thom spectrum
functor in the context of $\sI_G^U$-spaces.  This idea goes back
to~\cite{May1977}, which predated the definition of orthogonal
spectra; see~\cite[8.5]{Schlichtkrull2009} for a modern exposition.

Given a map of $\sI_G^U$-spaces $X \to BU_G(-)$, we can form the
pullback  
\begin{equation}
\begin{tikzcd}
Q \ar[r] \ar[d] & EU_G(V) \ar[d] \\
X(V) \ar[r] & BU_G(V)
\end{tikzcd}
\end{equation}
and from this pass to Thom spaces to obtain the $V$th space of a Thom
spectrum functor that lands in orthogonal $G$-spectra.  Because the
Thom space functor is multiplicative, this construction is as well.

\begin{defin}
The Thom spectrum construction yields a functor 
\begin{equation}
M \colon \sI_G^U-\Top / BU_G(-) \to \Sp_G.
\end{equation}
This functor is lax monoidal and externally multiplicative.
\end{defin}

\begin{proof}
The verification that $M$ is a lax monoidal functor to $\Sp_G$ follows
the lines of the arguments given in, for
example,~\cite{SagaveSchlichtkrull2019}.  The external
multiplicativity is a consequence of the fact that given a
$\sI_{G_1}^{U_1}$-space $X_1$ over $BU_{G_1}(-)$ and
$\sI_{G_2}^{U_2}$-space $X_2$ over $BU_{G_2}(-)$, the diagram
\begin{equation}
\begin{tikzcd}
X_1(V_1) \times X_2(V_2) \ar[r] \ar[d] & (X_1 \bar{\boxtimes} X_2)(V_1
\oplus V_2) \ar[d] \\
BU_{G_1}(V_1) \times BU_{G_2}(V_2) \ar[r] & (BU_{G_1}
\bar{\boxtimes} BU_{G_2})(V_1 \oplus V_2)
\end{tikzcd}
\end{equation}
commutes.
\end{proof}

For arbitrary complex $G$-vector bundles with fiber $V$, we have a
multiplicative Thom isomorphism defined as follows.  Such a bundle
with base $B$ is determined up to contractible choice by a classifying map 
\begin{equation}
\xi \colon B \to BU_G(|V|,V \otimes U) \to BU_G(U).
\end{equation}
The Thom diagonal is then classified by the map $B \to B \times B$
over $BU_G(U)$ where the second map is the composite of the projection
and $f$.  The complex orientation determines a map $M\xi \to MU_G$, and
so we obtain the composite 
\begin{equation}
M\xi \sma MU_G \to \Sigma^{\infty}_+ B \sma M\xi \sma
MU_G \to \Sigma^{\infty}_+ B \sma MU_G \sma MU_G \to  
\Sigma^{\infty}_+ B \sma MU_G.
\end{equation}
Since $M\xi$ is equivalent to $\Sigma^{-V} T\xi$, where $T\xi$ denotes
the Thom space of $\xi$, this then yields the following Thom
isomorphism equivalence:

\begin{prop}\label{prop:thom-direct}
For a complex $G$-bundle $\xi$ with fiber $S^V$, there is a natural
equivalence 
\begin{equation}
T\xi \sma MU_G \to \Sigma^V_+ B_{\xi} \sma MU_G
\end{equation}
that is externally multiplicative in the sense that give a
$G_1$-bundle $\xi_1$ and a $G_2$-bundle $\xi_2$, the diagram
\begin{equation}
\begin{tikzcd}
(T\xi_1 \sma MU_{G_1}) \sma (T\xi_2 \sma MU_{G_2}) \ar[r] \ar[d] &
  (\Sigma^{V_1}_+ B_{\xi_1} \sma MU_{G_1}) \sma (\Sigma^{V_2}_+
  B_{\xi_2} \sma MU_{G_2}) \ar[d] \\ 
T(\xi_1 \sma \xi_2) \sma MU_{G_1 \times G_2} \ar[r] &
  \Sigma^{V_1 \oplus V_2}_+ (B_1 \times B_2) \sma MU_{G_1 \times
  G_2} 
\end{tikzcd}
\end{equation}
commutes.
\end{prop}

\begin{proof}
This follows from the external multiplicativity of $MU_G$, the fact
that the Thom spectrum functor is externally multiplicative, and the
compatibility of the complex orientation with the multiplicative
structure (by definition).
\end{proof}

It is also useful to have a version of the Thom isomorphism that
untwists the action on the fiber, using the construction above. 

\begin{prop}
Let $\xi$ be a finite-dimensional complex $G$-bundle with base 
$B_{\xi}$ and fiber $V$.  Then there is a zig-zag of weak equivalences
\begin{equation}
\begin{tikzcd}
MUP_G \sma \Sigma^{\infty}_+ B_{\xi} \ar[r] & \Sh_V MUP_G \sma
  \Sigma^{\infty}_+ B_{\xi} & T\xi \sma MUP_G, \ar[l]  
\end{tikzcd}
\end{equation}
where the lefthand map is the canonical comparison of the suspension
and shift functors and the righthand map is the Thom isomorphism map
induced by the Thom diagonal.  The map $MUP_G \to \Sh_V MUP_G$ is
degree $-|V|$ with respect to the $\bZ$-grading (and the righthand map
has degree $0$). \qed
\end{prop}

We can derive from this the following analogue of
Corollary~\ref{cor:thom-zigzag}.
  
\begin{thm}\label{thm:thom-iso}
Let $\xi$ be a finite-dimensional complex $G$-bundle with base 
$B_{\xi}$ and fiber $V$.  Then there is a natural zig-zag of
externally multiplicative weak equivalences of orthogonal $G$-spectra  
\begin{equation}
T\xi \sma MUP^0_G \htp \Sigma^{|V|}_+ B_{\xi} \sma MUP^0_G,
\end{equation}
where $T\xi$ denotes the Thom space of the spherical fibration
associated to $\xi$.
\end{thm}

\subsection{Rigidifying spheres}\label{sec:rigidifying-spheres}

In this section, we explain how to use the Thom isomorphism to
discretize a topologized category of complex vector spaces.  

\begin{defin}
Let $S \Vect_{\bC}$ be the enriched category of \emph{stable complex
  vector spaces}, with 
\begin{itemize}
\item objects specified by pairs of finite-dimensional complex vector
  spaces $(I,V)$ equipped with an inner product, and 
\item morphisms consisting of a complex embedding $f \co V_0 \to V_1$
  preserving the inner product, and an isomomorphism
\begin{equation}
  I_1 \cong I_0 \oplus V_f^\perp.  
\end{equation}
\end{itemize}
The morphisms are topologized using the topology on the space of
linear isometries.
\end{defin}

The category $S \Vect_{\bC}$ is a monoidal category, with monoidal
structure given by the direct sum   
\begin{equation}
  (I_0,V_0) \oplus (I_1, V_1) \to (I_0 \oplus I_1, V_0 \oplus V_1)  
\end{equation}
and unit the pair $(\{0\}, \{0\})$.  (The monoidal product is also
clearly an enriched functor.)

\begin{lem}
There is an enriched functor
\begin{equation}
\Psi_{\bS} \colon S \Vect_{\bC} \to \Sp
\end{equation}
specified on objects by the assignment $(I,V) \mapsto F(S^V, (S^I)^{\mfib})$
and which on morphisms is given by the map
\begin{equation}
\begin{aligned}
 F(S^{V_0}, (S^{I_0})^{\mfib}) &\to F(S^{V_0} \sma S^{V_f^\perp},
 (S^{I_0})^{\mfib} \sma S^{V_f^\perp}) \\ &\to F(S^{V_0} \sma S^{V_f^\perp},
 (S^{I_0})^{\mfib} \sma (S^{V_f^\perp})^{\mfib}) \\
&\to F(S^{V_0} \sma S^{V_f^\perp},
 (S^{I_0} \sma S^{V_f^\perp})^{\mfib}) \to F(S^{V_1}, (S^{I_1})^{\mfib})
\end{aligned}
\end{equation}
obtained by smashing with the identity on $V_f^\perp$. 
\end{lem}

\begin{rem}
Note that the endomorphisms of any object of $S \Vect_{\bC}$ are given
by a product of unitary groups. This product acts on the target of the
functor side via the inclusions $U(n) \to F(n)$. 
\end{rem}

The functor $\Psi_{\bS}$ is lax monoidal, via the natural map induced
by the smash product
\begin{equation}
\begin{aligned}
 F(S^{V_0}, (S^{I_0})^{\mfib})  \sma  F(S^{V_1}, (S^{I_1})^{\mfib}) &\to  F(S^{V_0}
 \sma S^{V_1}, (S^{I_0})^{\mfib} \sma (S^{I_1})^{\mfib}) \\
 &\to F(S^{V_0 \oplus V_1},  (S^{I_0 \oplus I_1})^{\mfib}) 
\end{aligned}
\end{equation}
and the unit map
\begin{equation}
S^{0} \cong F(S^{0}, S^{0}) \to F(S^{0}, (S^{0})^{\mfib}) \cong (S^{0})^{\mfib}.
\end{equation}

Now let $\bbk$ be a cofibrant associative ring spectrum that equipped
with a multiplicative complex orientation.  Smashing with $\bbk$
yields a lax monoidal functor  
\begin{equation} \label{eq:stable_vector_space_to_spectra}
\Psi_{\bbk} \colon S \Vect_{\bC} \to \kmod 
\end{equation}
specified on objects as $(I,V)\mapsto F(S^V, (S^I)^{\mfib}) \wedge
\bbk$ and with morphisms as above.

On the other hand, we have a functor
\begin{equation}
\Psi_{\disc} \colon S \Vect_{\bC} \to \kmod
\end{equation}
specified on objects as $(I,V) \mapsto \bS[|I| - |V|] \sma
\bbk$, which we describe as follows.

\begin{defin}
Let $\tilde{\bZ}$ be the category with objects pairs $(n,m)$ with
$n,m \in \bN$ and morphisms $(n,m) \to (n',m')$ specified by $m < m'$
and $n' = n + (m' - m)$; i.e., there is a morphism precisely when $n-m
= n' - m'$.  
\end{defin}

The category $\tilde{\bZ}$ is monoidal, with unit $(0,0)$
and product given by 
\begin{equation}
(n,m) \oplus (n',m') = (n+n', m+m'). 
\end{equation}

There is a natural functor
\begin{equation}\label{eq:disc}
S \Vect_{\bC} \to \tilde{\bZ}
\end{equation}
specified on objects by $(I,V) \mapsto (|I|,|V|)$ and on morphisms by
taking a morphism $(I_0,V_0) \to (I_1, V_1)$ to the unique map 
$(|I_0|,|V_0|) \to (|I_1|,|V_1|)$ in $\tilde{\bZ}$.  Moreover, this
functor is clearly lax monoidal.

In addition, there is a functor
\begin{equation}\label{eq:sphere}
\tilde{\bZ} \to \kmod
\end{equation}
specified on objects by $(n,m) \mapsto \bS[n - m] \sma
\bbk$ and on morphisms as the identity map.  This functor is also lax
monoidal, via the associative multiplication map
\begin{equation}
\bS[n-m] \sma \bS[n'-m'] \to \bS[(n-m)+(n'-m')] = \bS[(n+n')-(m+m')]
\end{equation}
and the unit map $\bS \to \bS[0]$.

We now define the lax monoidal functor 
\begin{equation}
\Psi_{\disc} \colon S \Vect_{\bC} \to \tilde{\bZ} \to \kmod 
\end{equation}
as the composite of the functors in Equations~\eqref{eq:disc}
and~\eqref{eq:sphere}; on objects this takes $(I,V)$ to $\bS[|I| -
  |V|] \sma \bbk$.

The goal of the remainder of this section is to construct a zig-zag of
natural equivalences between the functors $\Psi$ and $\Psi_{\disc}$,
which we will do in stages.  First, observe that we can construct a
functor
\begin{equation}
\Psi_{\disc}^0 \colon S \Vect_{\bC} \to \kmod
\end{equation}
that factors through $\tilde{\bZ}$ as follows.  The functor
$\Psi_{\disc}^0$ is specified on objects by the assignment $(I,V)
\mapsto F(S^{|V|}, \bS[|I|]) \sma \bbk$ and on morphisms by assigning
to a morphism $(I_0, V_0) \to (I_1, V_1)$ the map
\begin{equation}
\begin{aligned}
F(S^{|V_0|}, \bS[|I_0|]) \sma \bbk &\to F(S^{|V_0|} \sma S^{|V_f^{\perp}|},
\bS[|I_0|] \sma S^{|V_f^{\perp}|}) \sma \bbk \\
&\to F(S^{|V_1|}, \bS[|I_0|] \sma S^{|V_f^{\perp}|}) \sma \bbk \to F(S^{|V_1|}, \bS[|I_0|] \sma
\bS[|V_f^{\perp}|]) \sma \bbk \\
&\to F(S^{|V_1|}, \bS[|I_1|]) \sma \bbk.
\end{aligned}
\end{equation}
This functor is evidently lax monoidal, using the isomorphisms
$S^{|V|} \sma S^{|W|} \cong S^{|V\oplus W|}$ and the natural maps
$\bS[V] \sma \bS[W] \to \bS[V \oplus W]$.

There is a natural weak equivalence $\Psi_{\disc} \to \Psi^0_{\disc}$ induced by
the natural map
\begin{equation}
\bS[|I|-|V|] \to F(S^{|V|}, \bS[|I|])
\end{equation}
which is the adjoint of the composite
\begin{equation}
\bS[|I|-|V|] \sma S^{|V|} \to \bS[|I|-|V|] \sma \bS[|V|] \to \bS[|I|].
\end{equation}

To check this is a natural transformation, we need to verify that the
diagram
\begin{equation}
\begin{tikzcd}
\bS[|I_0|-|V_0|] \ar[dddd] \ar[r] & F(S^{|V_0|}, \bS[|I_0|]) \ar[d] \\
& F(S^{|V_0|} \sma S^{|V_f^{\perp}|}, \bS[|I_0|]
\sma S^{|V_f^{\perp}|}) \ar[d] \\
& F(S^{|V_1|}, \bS[|I_0|] \sma S^{|V_f^{\perp}|}) \ar[d] \\ 
& F(S^{|V_1|}, \bS[|I_0|] \sma \bS[|V_f^{\perp}|]) \ar[d] \\ 
\bS[|I_1|-|V_1|] \ar[r] & F(S^{|V_1|}, \bS[|I_1|]).
\end{tikzcd}
\end{equation}
commutes.  First, observe that the map 
\begin{equation}
\bS[|I_0|-|V_0|] \to F(S^{|V_0|} \sma S^{|V_f^{\perp}|}, \bS[|I_0|] \sma S^{|V_f^{\perp}|})
\end{equation}
that arises as the adjoint of the map 
\begin{equation}
\bS[|I_0|-|V_0|] \sma S^{|V_0|} \sma S^{|V_f^{\perp}|} \to \bS[|I_0|] \sma S^{|V_f^{\perp}|}
\end{equation}
makes the diagram
\begin{equation}
\begin{tikzcd}
\bS[|I_0|-|V_0|] \ar[dr] \ar[r] & F(S^{|V_0|}, \bS[|I_0|]) \ar[d] \\
& F(S^{|V_0|} \sma S^{|V_f^{\perp}|}, \bS[|I_0|] \sma S^{|V_f^{\perp}|})
\end{tikzcd}
\end{equation}
commute.  This reduces the question to showing that the diagram
\begin{equation}
\begin{tikzcd}
\bS[|I_0|-|V_0|] \ar[dd] \ar[r] & F(S^{|V_1|}, \bS[|I_0|] \sma
S^{|V_f^{\perp}|}) \ar[d] \\  
& F(S^{|V_1|}, \bS[|I_0|] \sma \bS[|V_f^{\perp}|]) \ar[d] \\ 
\bS[|I_1|-|V_1|] \ar[r] & F(S^{|V_1|}, \bS[|I_1|]).
\end{tikzcd}
\end{equation}
commutes.  Next, observe that the diagram
\begin{equation}
\begin{tikzcd}
\bS[|I_0|-|V_0|] \ar[dr] \ar[r] & F(S^{|V_1|}, \bS[|I_0|] \sma
S^{|V_f^{\perp}|}) \ar[d] \\  
& F(S^{|V_1|}, \bS[|I_0|] \sma \bS[|V_f^{\perp}|]) 
\end{tikzcd}
\end{equation}
commutes by consideration of the adjoints.  Therefore, it suffices to
consider the diagram
\begin{equation}
\begin{tikzcd}
\bS[|I_0|-|V_0|] \ar[d] \ar[r] & F(S^{|V_1|}, \bS[|I_0|] \sma \bS[|V_f^{\perp}|]) \ar[d] \\ 
\bS[|I_1|-|V_1|] \ar[r] & F(S^{|V_1|}, \bS[|I_1|]).
\end{tikzcd}
\end{equation}
Next, observe that the diagram
\begin{equation}
\begin{tikzcd}[column sep=small] 
\bS[|I_0|-|V_0|] \ar[d] \ar[r] & F(\bS[|V_0|] \sma \bS[|V_f^{\perp}|],
\bS[|I_0|] \sma \bS[|V_f^{\perp}|]) \ar[d] & \ar[l] F(\bS[|V_1|],
\bS[|I_0|] \sma \bS[|V_f^{\perp}|]) \ar[d] \\
\bS[|I_1|-|V_1|] \ar[r] & F(\bS[|V_0|] \sma \bS[|V_f^{\perp}|],
\bS[|I_1|]) & \ar[l] F(\bS[|V_1|], \bS[|I_1|]).
\end{tikzcd}
\end{equation}
commutes by the associativity of the multplication on $\{\bS[-]\}$.
Thus, we are reduced to consideration of the diagram
\begin{equation}
\begin{tikzcd}
F(\bS[|V_1|], \bS[|I_0|] \sma \bS[|V_f^{\perp}|]) \ar[d] \ar[r] &
F(S^{|V_1|}, \bS[|I_0|] \sma \bS[|V_f^{\perp}|]) \ar[d] \\
F(\bS[|V_1|], \bS[|I_1|]) \ar[r] & F(S^{|V_1|}, \bS[|I_1|]), 
\end{tikzcd}
\end{equation}
which is clear.  It is tedious but straightforward to verify in an
analogous fashion that the natural transformation is compatible with
the monoidal structures.

Next, we have a functor 
\begin{equation}
\Psi_{\disc}^1 \colon S \Vect_{\bC} \to \kmod
\end{equation}
that factors through $\tilde{\bZ}$, specified on objects by the
assignment 
\begin{equation}
(I,V) \mapsto F(S^{|V|}, (S^{|I|})^{\mfib}) \sma \bbk
\end{equation}
and on morphisms by assigning to a morphism $(I_0, V_0) \to (I_1, V_1)$ the
map 
\begin{equation}
\begin{aligned}
F(S^{|V_0|}, (S^{|I_0|})^{\mfib}) \sma \bbk &\to F(S^{|V_0|} \sma S^{|V_f^{\perp}|},
(S^{|I_0|})^{\mfib} \sma S^{|V_f^{\perp}|}) \sma \bbk \\ 
&\to F(S^{|V_1|}, (S^{|I_1|})^{\mfib}) \sma \bbk. 
\end{aligned}
\end{equation}
This functor is evidently lax monoidal.  There is a zigzag of natural
transformations through weak equivalences induced by the evident maps
\begin{equation}
\begin{tikzcd}
F(S^{|V_0|}, (S^{|I_0|})^{\mfib}) \ar[r] & F(S^{|V_0|},
(\bS[|I_0|])^{\mfib}) & F(S^{|V_0|}, \bS[|I_0|]) \ar[l].
\end{tikzcd}
\end{equation}
Again, it is straightforward to check that this zigzag is monoidal.

Finally, we need to compare $\Psi_{\disc}^1$ to $\Psi$.  It is at this
point that we use the complex orientation on $\bbk$ in order to
discretize the topology on the category.  Recall from
Section~\ref{sec:eqMU} that the Thom isomorphism implies that 
for each complex representation $V$, there is a zig-zag of
equivalences 
\begin{equation}
S^V \sma MU_G \htp S^{|V|} \sma MU_G
\end{equation}
which are monoidal and $U(V)$-equivariant, where we give $S^{|V|}$ the
trivial $U(V)$-action.  Since $\bbk$ is complex oriented (see
Section~\ref{sec:ambid-orient}), we obtain analogous equivalences
\begin{equation}
\begin{tikzcd}
S^W \sma \bbk & \ar[l] \ar[r] \ldots & S^{|W|} \sma \bbk.
\end{tikzcd}
\end{equation}
Combining the equivalence
\begin{equation}
\begin{tikzcd}
F(S^W,(S^0)^{\mfib}) \sma \bbk & \ar[l] \ar[r] \ldots &
F(S^{|W|},(S^0)^{\mfib}) \sma \bbk 
\end{tikzcd}
\end{equation}
with the natural equivalences 
\begin{equation}
F(S^W, (S^0)^{\mfib}) \sma S^V \to F(S^W, (S^V)^{\mfib})
\end{equation}
and 
\begin{equation}
F(S^n, (S^0)^{\mfib}) \sma S^m \to F(S^n, (S^m)^{\mfib}).
\end{equation}
to specify an enriched zigzag of monoidal natural transformations
through weak equivalences.  

Putting everything together, we obtain the following comparison
theorem.

\begin{thm}
There is an enriched zigzag of monoidal natural transformations
through weak equivalences connecting $\Psi_{\bbk}$ and $\Psi_{\disc}$.
\end{thm}

\subsection{Morava $K$-theory}
\label{sec:ambid-orient}
\label{sec:ambidextrous-spectra}

Fix a prime $p$.  The chromatic filtration of the $p$-local stable
category mirrors the height filtration on the moduli stack of formal
group laws.  The filtration is controlled by a variety of cohomology
theories; the associated graded is related to certain cohomology
theories known as the Morava $K$-theories, which play an essential
role in our work.

Specifically, for each $n \in \bN$, there exists a periodic
cohomology theory represented by a spectrum $K(n)$; we suppress the
prime $p$ from the notation.  When $n = 0$, we understand $K(n)$ to be
the spectrum $H\bQ$ representing rational ordinary cohomology.  When
$n > 0$, the coefficients of $K(n)$ are specified by the formula
\begin{equation}
  K(n)_* \cong \bF_{p}[v_n^\pm],
\end{equation}
where $|v_n| = 2(p^n-1)$.

The spectra $K(n)$ have remarkable properties, the most salient of
which we summarize here:
\begin{enumerate}
\item The spectra $K(n)$ are associative ring spectra.  Note however
  that they are not $E_2$ (or even homotopy commutative when $p=2$);
  on the other hand, the coefficients are evidently always graded
  commutative.

\item For each $n > 0$, $K(n)$ is a graded field in the sense that any
  $K(n)$-module is free, i.e., is a wedge of shifts of $K(n)$.  (In
  fact, any graded field in spectra must be equivalent to a $K(n)$.)

\item For spaces $X$ and $Y$, there is a Kunneth isomorphism 
\begin{equation}
K(n)_*(X \times Y) \cong K(n)_*(X) \otimes_{K(n)_*} K(n)(Y).
\end{equation}

\item For any finite group $G$, there is an isomorphism 
\begin{equation}
K(n)_*(BG) \cong  K(n)^*(BG)
\end{equation}
of finite-rank modules~\cite{Ravenel82}.
\end{enumerate}
This last property can be interpreted as the statement that the Tate
fixed-point spectrum $K(n)^{tG}$ is trivial for any finite group
$G$~\cite{GreenleesSadofsky1996}.

Fix a compact Lie group $G$.  Let $X$ be an orthogonal $G$-spectrum.
We have the usual cofiber sequence of $G$-spaces 
\[
EG_+ \to S^0 \to \widetilde{EG},
\]
where the first map collapses $EG$ to the non-basepoint in $S^0$.
There is also a natural map $\gamma \colon X \to F(EG_+, X)$ induced
from the collapse map $EG_+ \to S^0$ and the homeomorphism
$F(S^0,X) \cong X$.  Combining $\gamma$ with the cofiber sequence
above, we obtain the diagram
\begin{equation}\label{eq:rawtate}
\xymatrix{
X \sma EG_+ \ar[r] \ar[d]_{\htp}^-{\gamma \sma \id} & X \ar[r] \ar[d]^-{\gamma} &
X \sma \widetilde{EG} \ar[d]^-{\gamma \sma \id} \\
F(EG_+, X) \sma EG_+ \ar[r] & F(EG_+,X) \ar[r] & F(EG_+,X) \sma \widetilde{EG}.\\ 
}
\end{equation}
(Note that since $X \to F(EG_+, X)$ is an underlying equivalence, the lefthand
vertical map is always a weak equivalence.)

Passing to $G$-fixed points on the bottom, we have the cofiber
sequence
\[
(F(EG_+, X) \sma EG_+)^G \to X^{hG} \to (F(EG_+,X) \sma \widetilde{EG})^G,
\]
and when $G$ is finite we can identify
\[
(F(EG_+,X) \sma EG_+)^G \htp X_{hG}.
\]
Therefore, for $G$-finite we have the commutative diagram
\[
\xymatrix{
X_{hG} \ar[d] \ar[r] & X^G \ar[d] \ar[r] & (\widetilde{EG} \sma X)^G \ar[d] \\
X_{hG} \ar[r] & X^{hG} \ar[r] & X^{tG}.
}
\]
Thus the norm map $X_{hG} \to X^{hG}$ is homotopic to the map 
\begin{equation}
X \sma EG_+ \to X \to F(EG_+, X)
\end{equation}
induced by the collapse map $EG_+ \to S^0$.

A useful observation about the homotopical properties of $X^{tG}$ is
that this is an invariant of the ``Borel'' homotopy type.  That is,
$(-)^{tG}$ preserves $G$-maps that are underlying equivalences.  In
particular, we will be interested in taking non-equivariant spectra
regarded as $G$-trivial orthogonal $G$-spectra as input; i.e.,
$\sI_{\bR^{\infty}}^U X$ for cofibrant $X$.  (Classically, this is often
written as $i_* X$.)

Setting $X = K(n)$, the vanishing of $K(n)^{tG}$ implies that the norm
map $X_{hG} \to X^{hG}$ is a weak equivalence, and so we have 
\begin{equation}
K(n) \sma BG_+ \htp X_{hG} \to X^{hG} \cong F(BG_+, K(n)).
\end{equation}
Since the Tate spectrum functor is lax monoidal, another consequence
of the vanishing of $K(n)^{tG}$ is the following more general
vanishing result.

\begin{cor}
Let $G$ be a finite group.  For any $K(n)$-module $M$, $M^{tG} \htp
\ast$.
\end{cor}

\subsection{The existence of $K(n)$ orientations}

We regard as fixed for the body of the paper for each prime $p$ and $n
> 0$ a choice of a point-set model of $K(n)$ that is an associative
ring orthogonal spectrum; we denote this by $\bbk$, as we have in the
body of the paper.  By obstruction theory~\cite{Angeltveit2011}, one
can construct uncountably many $A_\infty$ $MU$-algebra structures on
$\bbk$.  We can parametrize these in terms of formal group law data,
but for all practical purposes it doesn't matter as long as we fix
one.  Thus, we always have an $A_\infty$ orientation $MU \to \bbk$,
i.e., a multiplicative complex orientation of $\bbk$, which we
represent as a map of orthogonal ring spectra $MU \to \bbk$.

\begin{rem}
When $p$ is odd, following Rudyak~\cite{Rudyak1998} we can consider
the composite 
\begin{equation}
MSO \to BP \to MU \to \bbk,
\end{equation}
where the map $MSO \to BP$ is given by the classical decomposition of
$MSO$ into a wedge of shifts of $BP$~\cite{Wilson1982} and the map
$BP \to MU$ is the usual map.  The map $BP \to MU$ is
$E_4$~\cite[]{BasterraMandell2005} and the map $MU \to \bbk$ is $A_\infty$
as observed above.  Direct inspection shows that the projection $MSO
\to BP$ is a map of ring spectra in the stable category, and then by
obstruction theory we can promote this map to an $E_2$
map~\cite{ChadwickMandell2015}.  Therefore, at an odd prime, the
composite $MSO \to \bbk$ is an $A_\infty$ map; i.e., $\bbk$ admits an
$A_\infty$ $MSO$ orientation.

At the prime $2$, the situation is different and this kind of argument
does not appear to apply.  Specifically, in this case, $MSO$ splits as
a wedge of certain Eilenberg-Mac Lane spectra, but these do not
obviously map to $\bbk$.
\end{rem}

We now construct the equivariant complex orientation for $\bbk$.
First, note that $MU_G$ is ``split'' in the sense that there is a
natural map $\sI_{\bR^{\infty}}^U MU \to MU_G$ of ring spectra which
is an equivalence when the $G$-action is forgotten.  This is a shadow
of the global equivariant structure that $MU_G$ possesses.

\begin{warning}
When applied to cofibrant spectra or associative ring orthogonal
spectra, the point-set functor $\sI_{\bR^{\infty}}^U$ computes the
classical derived functor denoted $i_*$.  However, the behavior of
$\sI_{\bR^{\infty}}^U$ on cofibrant commutative ring orthogonal
spectra is different; the functor {\em does not} coincide with
$\sI_{\bR^{\infty}}^U$ on the underlying module.  As a consequence, in
what follows we will tacitly replace $MU$ by a cofibrant associative
ring orthogonal spectrum.
\end{warning}

As a first consequence, if $X$ is a free $G$-space, there is an
equivariant equivalence
\begin{equation}
(\sI_{\bR^{\infty}}^{U} MU) \sma X \to MU_G \sma X.
\end{equation}
In particular, if $R$ is a complex oriented orthogonal ring spectrum,
then for a free $G$-space $X$ there is a zigzag
\begin{equation}
MU_G \sma X \leftarrow (\sI_{\bR^{\infty}}^{U} MU) \sma X \to
(\sI_{\bR^{\infty}}^{U} R) \sma X 
\end{equation}
in which the backward arrow is a weak equivalence.

To arrange for the trivialization of complex representations to apply
to $\bbk$, we consider the $MU_G$-module version of $\bbk$ defined as
the (derived) smash product 
\begin{equation}
\bbk_G = MU_G \sma_{\sI_{\bR^{\infty}} MU} \sI_{\bR^{\infty}} \bbk.
\end{equation}
The spectrum $\bbk_G$ has underlying spectrum equivalent to $\bbk$ and
is again split; there exists a natural map $\sI_{\bR^{\infty}}^{U}
\bbk \to \bbk_G$ of associative ring orthogonal $G$-spectra induced as the map 
\begin{equation}
\sI_{\bR^{\infty}}^{U} \bbk \cong \sI_{\bR^{\infty}}^{U} MU
\sma_{\sI_{\bR^{\infty}}^{U} MU} \sI_{\bR^{\infty}}^{U} \bbk \to MU_G
\sma_{\sI_{\bR^{\infty}}^{U} MU} \sI_{\bR^{\infty}}^{U} \bbk.
\end{equation}
For free $G$-spaces $X$, there is again a natural equivalence
$\sI_{\bR^{\infty}}^{U} \bbk \sma X \to \bbk_G \sma X$.

The equivariant complex orientation data that we require is now
specified by the map of associative orthogonal ring spectra
\begin{equation}
MU_G \to MU_G \sma_{\sI_{\bR^{\infty}}^{U} MU} \sI_{\bR^{\infty}}^{U}
MU \to MU_G \sma_{\sI_{\bR^{\infty}}^{U} MU} \sI_{\bR^{\infty}}^{U}
\bbk,
\end{equation}
obtained by base-change applied to the complex orientation of $\bbk$.

\section{The Adams and Spanier-Whitehead equivalences}
\label{sec:adams-isom-ambid}

The goal of this appendix is to describe models of equivariant
Spanier-Whitehead duality, the Adams isomorphism, and the norm
map that are suitably functorial and multiplicative to carry out the
comparisons described in the body of the paper.  Our ultimate
objective is to prove Proposition~\ref{prop:adams-comparison}, which
asserts there is a natural zig-zag of $\Pi$-equivariant equivalences
connecting the spectral category $C^{*}(B G,
\ssF_{\frac{{\bar{\cX}}|{\cZ}}{N}}({\bar{\cX}}|{\cZ}, \bbk)^{-V-d} )$ 
and the spectral category $B\bar{\cX}|\cZ^{-V-d} \wedge \bbk$ of
virtual cochains.  After assembling the necessary intermediate
results, we summarize the proof of that comparison in
Section~\ref{sec:complete-comparison}.  Because it does not play a
significant role in any of our constructions, we suppress discussion
of the $\Pi$-equivariance throughout the work of this appendix.

\subsection{A review of the Wirthmuller and Adams isomorphisms}
\label{subsec:wirthadams}

We begin by giving an abstract discussion of the Wirthmuller and Adams
isomorphisms.  Our actual implementation (in the remainder of this
section) involves a specific choice of models adapted to our
application.

For a group $G$, a subgroup $H \subset G$, and an $H$-space $X$, there
is a natural map of spaces 
\begin{equation}\label{eq:wirthspace}
G_+ \sma_H X \to F_H(G_+, X)
\end{equation}
specified by the assignment that takes $(g, x)$ to the function
 given as 
\begin{equation}
\tilde{g} \mapsto
\begin{cases}
\tilde{g} g x \qquad \tilde{g}g \in H \\ \ast \qquad \quad
\textrm{otherwise}.
\end{cases}
\end{equation}
This map of spaces gives rise to a corresponding map of orthogonal
$G$-spectra, which is a stable equivalence; this is the point-set
realization of the Wirthmuller isomorphism, which is the name for the
derived equivalence between the left and right adjoints to the
forgetful functor $i_*^H$ from $G$-spaces to $H$-spaces.

We will also use a relative variant of this Wirthmuller map defined as
follows.  For a group $G$ and a subgroup $H \subset G$, we have a map
of $G$-spaces  
\begin{equation}\label{eq:relwirthspace}
G_+ \sma_H F(H_+,S^0) \to F(G_+,S^0)
\end{equation}
specified by the assignment that takes $(g,\gamma)$ to the function
given as
\begin{equation}
\tilde{g} \mapsto
\begin{cases}
\gamma(\tilde{g} g) \qquad \tilde{g}g \in H \\ \ast \qquad \qquad
\textrm{otherwise}.
\end{cases}
\end{equation}
Equivalently, this is specified by the $H$-map from $F(H_+, S^0)$ to
$F(G_+, S^0)$ that is the extension by $0$.
The space-level map of Equation~\eqref{eq:relwirthspace} induces a
spectrum-level map 
\begin{equation}\label{eq:relwirth}
G_+ \sma_H F(H_+, \bS) \to F(G_+, \bS).
\end{equation}

We now turn to a discussion of the Adams isomorphism, which is more
subtle.  The Adams isomorphism for a normal subgroup $H \subset G$
relates two functors defined on the subcategory of $H$-free objects in
$\Sp_G^U$.  Here recall the definition of an $H$-free spectrum.  Let
$E\scrF_H$ denote the classifying space for the family of subgroups
such that the intersection with $H$ is trivial; this is determined up
to homotopy by the requirement that
\begin{equation}
\begin{cases}
(E\scrF_H)^K \htp \ast \qquad K \cap H = \{e\} \\
(E\scrF_H)^K = \emptyset \qquad K \cap H \neq \{e\}.
\end{cases}
\end{equation}
For example, when $H=G$, $E\scrF_H$ is a model of $EG$.  A canonical
model of $E\scrF_H$ can be constructed using Elmendorf's
theorem~\cite{Elmendorf1983}.

\begin{defin}
An orthogonal $G$-spectrum $X$ is $H$-free if the natural map
\begin{equation}
E\scrF_H \sma X \to X
\end{equation}
is an equivalence.
\end{defin}

For example, the category of $G$-free spectra is the Borel equivariant
category.  In analogy with the discussion in Section~\ref{sec:borel},
we can equivalently describe the full subcategory of $H$-free spectra
in $\Sp_G^U$ as the category $\Sp_G^{U^H}$ of orthogonal $G$-spectra
on the universe $U^H$.

For $H \subset G$ a normal subgroup, the orbits $(-)_{H}$
determine a functor   
\begin{equation}
(-)_{H} \colon \Sp_G^{U} \to \Sp_{G/H}^{U^H},
\end{equation}
and the $H$-fixed points $(-)^H$ determine a functor
\begin{equation}
(-)^{H} \colon \Sp_G^U \to \Sp_{G/H}^{U^H}.
\end{equation}

The Adams isomorphism asserts that when restricted to the category of
$H$-free spectra in $\Sp_G^U$, these functors are naturally
equivalent as the next theorem asserts. For its statement, to compute the derived $H$-fixed
points, we fibrantly replace in the complete universe $U$ (rather than universe $U^H$).

\begin{thm}
There is a derived natural transformation 
\begin{equation}
(-)_H \to (-)^H
\end{equation}
that is a weak equivalence when restricted to the full subcategory of
$H$-free spectra in $\Sp_G^U$.  
\end{thm}

\begin{example}
When $H = G$, this implies a natural equivalence 
\begin{equation}
EG_+ \sma_G X \htp (\cQ_U X)^G.
\end{equation}
\end{example}

We will need a point-set model of this equivalence.  As preparation
for constructing this, it is useful to recall the framework for the
classical proof of the Adams isomorphism;
see~\cite[II.7]{LewisMaySteinberger1986}.  Take $N \subset G$ to be a
normal subgroup, and let $\rho \colon G \to G/N$ be the quotient map.
Let $\Gamma = G \rtimes N$ denote the semidirect product with respect
to the conjugation action.  There are natural homomorphisms $\theta,
\epsilon \colon \Gamma \to G$ specified by the assignments
\begin{equation}
\begin{tikzcd}
\theta (g,n) = gn && \epsilon (g,n) \mapsto g.
\end{tikzcd}
\end{equation}
There is an isomorphism $\ker(\epsilon) \cong N$; we will denote this
subgroup of $\Gamma$ by $\tilde{N}$.  In fact, $\Gamma / G \cong N$ as
$\Gamma$-sets, where the action of $\Gamma$ on $N$ is the evident
combination of conjugation by $G$ and the action of $N$ on itself by
left multiplication, i.e., $(g,m) n = g(mn) g^{-1}$.

\begin{rem} \label{rem:why-we-use-semidirect}
The reason to introduce the semidirect product action is that in
general $X \sma_N N$ cannot be given a $G$-action which makes it
homeomorphic to $X$ as a $G$-space or spectrum.  On the other hand,
there is always a homeomorphism
\begin{equation}
\theta^* X \sma_{\tilde{N}} N \cong X
\end{equation}
as $G$-spaces or spectra.
\end{rem}

The core of the Adams isomorphism is the transfer $\tau \colon \bS \to
\Sigma^{\infty_+} N$ in $\Gamma$-spectra associated to the collapse
map $N \to \ast$ (regarded as a map of $\Gamma$-sets).  Specifically,
given an $N$-free $G$-spectrum $X$, we can produce an equivalent
$\tilde{N}$-free $\Gamma$-spectrum on a complete $\Gamma$-universe;
denote this spectrum by $\tilde{X}$.  The Adams isomorphism is then
induced by the transfer
\begin{equation}
\tau \colon \tilde{X} \to N_+ \sma \tilde{X},
\end{equation}
as follows.  Passing to orbit spectra with regard to $\tilde{N}$, we obtain a
map 
\begin{equation}
\tau/\tilde{N} \colon \rho^* X/N \cong \tilde{X}/\tilde{N} \to (N_+
\sma \tilde{X})/\tilde{N} \cong X 
\end{equation}
of $G$-spectra on the complete universe, where we are implicitly
pushing forward from the $N$-fixed universe to the complete universe.
Since $\rho^*$ is the left adjoint of the functor $(-)^N$, this is
equivalent to a map of $G/N$-spectra
\begin{equation}
X/N \to X^N.
\end{equation}
The theorem is now that this map is an equivalence.  Once again, we
emphasize that in order to compute the codomain of this map, we
require a fibrant replacement of $X$ with respect to the complete
universe.  The statement is wildly false if we do not use the derived
fixed points in $\Sp_G^U$, as the example of the spectrum
$\Sigma^{\infty} G_+$ demonstrates.

Understanding the functoriality and multiplicative properties of the
Adams isomorphism therefore boils down to studying the behavior of the
transfer map $\tau$.  It is not hard to verify that the transfer
commutes with the external smash product in the homotopy category
(e.g., see~\cite[III.5.10]{LewisMaySteinberger1986}).  Our work
requires rigidifying this structure.  One way to do this is to enlarge
the category to keep track of the embeddings required for the
Pontryagin-Thom construction; this is the kind of approach that is
taken in Cohen's description of Atiyah duality as an $E_\infty$
map~\cite{Cohen2004}.  Although it would be possible to take the same
tack here, we take a different approach that avoids the need for
tracking this data, in part by working on the other side of the
Wirthmuller isomorphism.

We now describe various technical underpinnings of the approach we
take.  First, we consider the model of the transfer expressed as 
$\bS \to F^{\Gamma}(N_+, \bS)$ to obtain a map
\begin{equation}
\tilde{X} \to \tilde{X} \sma F^{\Gamma}(N_+, \bS).
\end{equation}
Passing to orbits and taking the adjoint again yields the Adams map,
now written as
\begin{equation}
X/N \to (\tilde{X} \sma F^{\Gamma}(N_+, \bS)/ N)^N.
\end{equation}
The transfer $\bS \to F^{\Gamma}(N_+, \bS)$ is evidently strictly
multiplicative in the sense that the diagram
\begin{equation}
\begin{tikzcd}
\bS_{\Gamma_1} \sma \bS_{\Gamma_2} \ar[r] \ar[d] & \bS_{\Gamma_1 \times
  \Gamma_2} \ar [d] \\ 
F^{\Gamma_1}((N_1)_+, \bS_{\Gamma_1}) \sma F^{\Gamma_2}((N_2)_+,
\bS_{\Gamma_2}) \ar[r] & F^{\Gamma_1 \times \Gamma_2}((N_1 \times
N_2)_+, \bS_{G_1 \times G_2}) 
\end{tikzcd}
\end{equation}
commutes, where we are working in point-set categories of orthogonal
$G$-spectra.  We decorate the sphere spectrum in the top row to
indicate the category of spectra involved.

Next, we recall the following convenient point-set model of the Adams
isomorphism, due to Schwede~\cite[\S 8]{Schwede2019}.  We are going to
focus on the case where $N = G$, in which case the semidirect product
$\Gamma$ is simply $G \times G$.  That is, we work with orthogonal
spectra with commuting left and right $G$ actions, which we call
\emph{biequivariant}.  We write $G_r$ when considering the right
action and omit any subscript from the notation when considering the
left action.  

\begin{lem}
For an orthogonal $G$-spectrum $X$, the tensor $G_+ \sma X$ is a
biequivariant spectrum.  Passing to derived $G$-fixed points produces
a $G$-spectrum $(\cQ_{\cU} (G_+ \sma X))^G$.
\end{lem}

There is a natural map
\begin{equation}\label{eq:adams-assemble}
\alpha_1(G,\cU) \colon (\cQ_{\cU} (G_+ \sma X))^G \sma_G EG_+ \to (\cQ_{\cU} (X \sma EG_+))^G
\end{equation}
induced by the action of $G$ on $EG$.  On the other hand, the
Wirthmuller isomorphism yields a comparison
\begin{equation}\label{eq:adams-wirth}
\alpha_2(G,\cU) \colon (\cQ_{\cU} (G_+ \sma X))^G \sma_G EG_+ \to
(\cQ_{\cU} F(G, X))^G \sma_G EG_+ .
\end{equation}
Finally, there is a natural map
\begin{equation}\label{eq:adams-fixorb}
\alpha_3(G,\cU) \colon EG_+ \sma_G X \cong X \sma_G EG_+ \to
(\cQ_{\cU} F(G, X))^G \sma_G EG_+
\end{equation}
induced by the natural map $X \to (F(G,X))^G$, which is a map of
orthogonal $G$-spectra when $X$ is biequivariant. 

Putting this together, we can give the following point-set model of
the Adams isomorphism.

\begin{defin} \label{def:Adams-pointset-biequivariant}
Let $X$ be an orthogonal $G$-spectrum.  The Adams zig-zag is the
comparison of $EG_+ \sma_G X$ and $(\cQ_{\cU} (X \sma EG_+))^G$
constructed as the composite of $\alpha_1$, $\alpha_2$, and
$\alpha_3$.
\end{defin}

The statement of the Adams isomorphism is that the Adams zig-zag is a
weak equivalence.

\begin{thm}
The Adams zig-zag is through weak equivalences.  Thus, when $X$ is a
$G$-free orthogonal $G$-spectrum, the zig-zag represents a natural
isomorphism in the stable category between $X_{hG}$ and $X^G$.
\end{thm}

Next, we discuss the multiplicativity of the maps in the Adams
zig-zag.  All of these facts are essentially consequences of the fact
that $Q_{(-)}(-)$ is externally multiplicative and $(-)^H$ is
externally lax monoidal.

\begin{lem}\label{lem:mult-assem}
Given an orthogonal $G_1$-spectrum $X_1$ and an orthogonal
$G_2$-spectrum $X_2$, there are natural product maps 
\begin{equation}
\begin{tikzcd}
(Q_{\cU_1} (G_{1,+} \sma X_1))^{G_1} \sma_{G_1} EG_{1,+} \sma 
(Q_{\cU_2} (G_{2,+} \sma X_2))^{G_2} \sma_{G_2} EG_{2,+} 
\ar["\mu",d] \\
(Q_{\cU_1 \oplus \cU_2} ((G_{12})_+ \sma (X_1 \sma X_2)))^{G_{12}}
  \sma_{G_{12}} EG_{12,+} 
\end{tikzcd}
\end{equation}
and
\begin{equation}
\begin{tikzcd}
(Q_{\cU_1} (X_1 \sma  EG_{1,+}))^{G_1} \sma (Q_{\cU_2} (X_2 \sma EG_{2,+}))^{G_2} \ar["\mu'",d] \\
(Q_{\cU_1 \oplus \cU_2} ((X_1 \sma X_2) \sma EG_{12,+}))^{G_{12}},
\end{tikzcd}
\end{equation}
where we write $G_{12} = G_1 \times G_2$.  These product maps are
compatible with $\alpha_1(-,-)$ in the sense that
\begin{equation}
\mu' \circ (\alpha_1(G_1, \cU_1) \sma \alpha_1(G_2, \cU_2)) =
\alpha_1(G_1 \times G_2, \cU_1 \oplus \cU_2) \circ \mu.
\end{equation}
Moreover, the evident associativity and unitality diagrams commute as
well. 
\end{lem}

\begin{lem}\label{lem:mult-wirth}
Given an orthogonal $G_1$-spectrum $X$ and an orthogonal
$G_2$-spectrum $X_2$, there are natural product maps 
\begin{equation}
\begin{tikzcd}
\left((Q_{\cU_1} (G_{1,+} \sma X_1))^{G_1} \sma_{G_1} EG_{1,+} \right) \sma  
\left((Q_{\cU_2} (G_{2,+} \sma X_2))^{G_2} \sma_{G_2} EG_{2,+} \right) 
\ar[d] \\
(Q_{\cU_1 \oplus \cU_2} (G_{12,+} \sma (X_1 \sma X_2)))^{G_{12}}
  \sma_{G_12} EG_{12,+} 
\end{tikzcd}
\end{equation}
and
\begin{equation}
\begin{tikzcd}
\left((Q_{\cU_1} F(G_1, X_1))^{G_1} \sma_{G_1} EG_{1,+}\right) \sma \left((Q_{\cU_2} F(G_2, X_2))^{G_2} \sma_{G_2} EG_{2,+}\right) \ar[d] \\
(Q_{\cU_1 \oplus \cU_2} F(G_{12}, X_1 \sma X_2))^{G_{12}} \sma_{G_{12}} EG_{12,+}, 
\end{tikzcd}
\end{equation}
where we write $G_{12} = G_1 \times G_2$.  These products maps are
compatible with $\alpha_2(-,-)$ in the sense that
\begin{equation}
\mu' \circ (\alpha_2(G_1, \cU_1) \sma \alpha_2(G_2, \cU_2)) =
\alpha_2(G_1 \times G_2, \cU_1 \oplus \cU_2) \circ \mu.
\end{equation}
Moreover, the evident
associativity and unitality diagrams commute as well.
\end{lem}

\begin{proof}
This result ultimately depends on the fact that the Wirthmuller map is
externally multiplicative in the sense that the diagram
\begin{equation}
\begin{tikzcd}
((G_1)_+ \sma X_1) \sma ((G_2)_+ \sma X_2) \ar[r] \ar[d] & F((G_1)_+, X_1) \sma
  F((G_2)_+, X_2) \ar[d] \\ 
(G_{12})_+ \sma (X_1 \sma X_2) \ar[r] & F((G_{12})_+, X_1 \sma X_2) \\ 
\end{tikzcd}
\end{equation}
commutes as a diagram of orthogonal $G_{12}$-spectra.
\end{proof}

\begin{lem}
Given an orthogonal $G_1$-spectrum $X_1$ and an orthogonal
$G_2$-spectrum $X_2$, there are natural product maps
\begin{equation}
\begin{tikzcd}
\left(EG_{1,+} \sma_{G_1} X_1\right) \sma \left(EG_{2,+} \sma_{G_2} X_2\right)
\ar[d] \\
EG_{12,+} \sma_{G_{12}} (X_1 \sma X_2) 
\end{tikzcd}
\end{equation}
and
\begin{equation}
\begin{tikzcd}
EG_{1,+} \sma_{G_1} (Q_{\cU_1} F(G_1, X_1))^{G_1} \sma EG_{2,+} \sma_{G_2} (Q_{\cU_2} F(G_2, X_2))^{G_2}
\ar[d] \\
EG_{12,+} \sma_{G_{12}} (Q_{\cU_1 \oplus \cU_2} F(G_{12}, X_1 \sma X_2))^{G_{12}}
\end{tikzcd}
\end{equation}
where we write $G_{12} = G_1 \times G_2$.  These product maps are
compatible with $\alpha_3(-,-)$ in the sense that
\begin{equation}
\mu' \circ (\alpha_3(G_1, \cU_1) \sma \alpha_3(G_2, \cU_2)) =
\alpha_3(G_1 \times G_2, \cU_1 \oplus \cU_2) \circ \mu.
\end{equation}
\end{lem}

Putting this all together, we summarize in the following proposition.

\begin{prop}\label{prop:Adams-external-mult}
The Adams zig-zag is compatible with the external smash product of spectra.
\end{prop}

In the remainder of this appendix, we will establish a pointwise model for the functoriality of the Adams isomorphism with respect to change of groups. This depends in part on the
functoriality of the Wirthmuller isomorphism, and is clearly more
complicated than simple functoriality in the category of equivariant
spectra; given a map $f \colon (G, X) \to (G',Y)$, there is a zig-zag
\begin{equation}
\begin{tikzcd}
F(G_+,X) \ar[r] & F(G_+, f^* Y) & F(G'_+, Y) \ar[l] 
\end{tikzcd}
\end{equation}
of orthogonal $G$-spectra.

\subsection{Fibrewise mappings spaces}
\label{sec:fibr-mapp-spac}

As discussed in Section \ref{sec:funct-span-whit}, functoriality for
Spanier-Whitehead duality isomorphism with respect to covering maps is
conveniently expressed using a notion of fibrewise Spanier-Whitehead
duality; we will employ the same approach to express the functoriality
of Adams isomorphisms, so we introduce the relevant fibrewise spaces
in this section.

\subsubsection{Construction of fibrewise mapping spaces}
\label{sec:constr-fibr-mapp}

Let $G$ be a finite group and $N \subseteq G$ a normal subgroup.  Fix
a complete universe $U$ and let $Y$ be an orthogonal $G$-spectrum.  We
shall consider $Y$ as an orthogonal $\Gamma$-spectrum $\theta^* Y$
indexed on the universe $\theta^* U$.  As in the discussion
surrounding Remark \ref{rem:why-we-use-semidirect}, we write
$\tilde{N}$ for the normal subgroup $\{1,n\} \subset \Gamma$ which
isomorphic to $N$.  There is an action of $\Gamma$ on $N$ specified by
the assignment $(g,n)x \mapsto gnxg^{-1}$.

\begin{defin}\label{defin:fiberwise-dual-total-space}
For each spectrum $\bbk$, the \emph{spectrum of maps from $Y$ to
  $\bbk$ over $Y/N$} is the orthogonal $G$-spectrum
$\ssF_{Y/N}(Y,\bbk)$, indexed on the universe $U$,  given by
\begin{equation}
  \ssF_{Y/N}(Y,\bbk) \equiv \theta^*Y \sma_{\tilde{N}} F(N_+, \bbk^{\mfib}),
\end{equation}
where the $G$-action is induced by the canonical isomorphism $\Gamma /
\tilde{N} \cong G$ and the action of $G$ on $F(N_+, \bbk^{\mfib})$
is induced by the $G$-action on $N$.
\end{defin}

This spectrum is intended as a model of the total space of the
``bundle'' $Y \to Y/N$, where we are working with the
Spanier-Whitehead dual of $N$.  To compare to the usual construction,
we will make use of the duality stable equivalence
\begin{equation}
\eta \colon \Sigma^{\infty} N_+ \to F(N_+, \bS),
\end{equation}
constructed as in Equation~\eqref{eq:wirthspace} but regarded as a map
of $\Gamma$-spectra.

\begin{example}
If $X$ is a based $G$-space, then the spectrum
$\ssF_{\Sigma^{\infty}X/N}(\Sigma^{\infty}X,\bS)$ may be interpreted
as the total space of the \emph{fibrewise Spanier-Whitehead dual of
  $X$ over $X/N$}.  Indeed, if we ignore equivariance, the space
assigned by this spectrum to a finite-dimensional real vector space
$W$ receives a natural map from  
\begin{equation}
  \theta^* X \sma_{\tilde{N}} \Map(N_+, S^W) \cong \cMap_{X/N}(X, S^W),  
\end{equation}
where the right hand side is the space of fibrewise maps from $X$ to
$S^W$ over $X/N$, and $S^W$ denotes the space over $X/N$ with
total space the product $S^W \times X/N$ and structure map the evident
projection.  Recall that this space is constructed by first
considering the space $ \cMap'_{X/N}(X,S^W)$ of fibrewise
unbased maps from $X$ to $S^W \times X/N$.  Since the map $X \to X/N$
is open, this space is weak Hausdorff (see \cite[Proposition
  1.5]{Lewis1985}) and admits a closed section consisting of
basepoint-valued maps, whose quotient is $\cMap_{X/N}(X,S^W)$.  See
e.g.,~\cite[1.3.7]{MaySigurdsson2006} for a more general discussion of
such mapping spaces.
\end{example}

The passage to fiberwise mapping spectra is multiplicative in the
following sense.

\begin{lem}\label{lem:fibmapmult}
Let $Y_1$ be an orthogonal $G_1$-spectrum and let $Y_2$ be an orthogonal
$G_2$-spectrum.  For normal subgroups $N_1 \subseteq G_1$ and $N_2
\subseteq G_2$, there are natural product maps of $G_1 \times
G_2$-spectra   
\begin{equation}
\ssF_{Y_1/N_1}(Y_1,\bbk)  \wedge  \ssF_{Y_2/N_2}(Y_2,\bbk) \to 
\ssF_{(Y_1 \wedge Y_2)/(N_1 \times N_2)}(Y_1 \wedge Y_2,\bbk)
\end{equation}
that are associative and unital.
\end{lem}

\begin{proof}
The pairings in question are induced by the natural maps 
\begin{equation}
Y_1 \wedge_{\tilde{N_1}} F(N_{1,+}, \bbk^{\mfib}) \wedge  Y_2 \wedge_{\tilde{N_2}} F(
N_{2,+}, \bbk^{\mfib}) \to 
Y_1 \wedge Y_2 \wedge_{\widetilde{(N_1 \times N_2)}} F( N_{1,+} \wedge N_{2,+},
\bbk^{\mfib}) 
\end{equation}
where here we are using the fact that there is an evident natural
isomorphism of groups
\begin{equation}
(G_1 \rtimes H_1) \times (G_2 \rtimes H_2) \cong (G_1 \times G_2)
  \rtimes (H_1 \times H_2).  
\end{equation} 
Associativity is straightforward to verify using the associativity of
the external smash product and the pairings on mapping spectra.  For
unitality, observe that $\bS \sma_{\tilde{e}} F(e_+, \bbk^{\mfib})
\cong \bbk^{\mfib}$.
\end{proof}

To construct our comparison zig-zag, we begin with the following observation:

\begin{lem}\label{lem:fiberwise-comparisons}
Let $Y$ be an $N$-free $G$-spectrum.  There are natural $G$-maps
\begin{equation}
\begin{tikzcd}
Y \sma \bbk \ar["\htp",r]  & \ssF_{Y/N}(Y,\bbk) & \ar[l]
\left((\theta^* Y\right) / \tilde{N}) \wedge \bbk, 
\end{tikzcd}
\end{equation}
where the indicated map is an equivalence.
\end{lem}

\begin{proof}
The left arrow (which is an equivalence) is induced by the composition  
\begin{equation}
 Y \sma \bbk \to \left((\theta^* Y) \wedge_{\tilde{N}}  N_+\right)
 \sma \bbk \to (\theta^*
  Y) \wedge_{\tilde{N}} F( N_+, \bbk^{\mfib}), 
\end{equation}
where the map $N_+ \sma \bbk \to F(N_+,\bbk^{\mfib})$ is the duality
map $\eta$.  The right arrow is induced by the map $\bbk \to F(N_+,
\bbk^{\mfib})$ induced by the projection $N_+ \to S^0$ (regarded as a
map of $\Gamma$-spaces).
\end{proof}

The above construction has the correct derived functors when $Y$ is
a cofibrant $N$-free $G$-spectrum.  Next, assume that $V$ is a
finite-dimensional $G$-representation, and define the desuspended
mapping space as follows:
\begin{equation}
\ssF_{Y/N}(Y,\bbk)^{-V} \equiv \Sigma^{-V} \left(\theta^* Y
\wedge_{\tilde{N}} F(N_+, \bbk^{\mfib})\right).
\end{equation}

We can immediately deduce the following refinement of
Lemma~\ref{lem:fiberwise-comparisons}.

\begin{cor}
There are natural $G$-maps
\begin{equation}
\begin{tikzcd}
\Sigma^{-V} Y \wedge \bbk \ar["\htp",r] & \ssF_{Y/N}(Y,\bbk)^{-V}
& \ar[l] \Sigma^{-V} \left((\theta^* Y\right) / \tilde{N}) \wedge \bbk 
\end{tikzcd}
\end{equation}
where the indicated map is an equivalence.\qed
\end{cor}

We now consider the interaction of this construction with the Borel
construction
\begin{equation}
BY^{-V} \equiv EG_+ \wedge_{G} (\Sigma^{-V} Y). 
\end{equation}
Let us denote by $J$ the quotient $G/N$, and consider  a
$J$-representation $V$.  We will also tacitly write $V$ to denote the
pullback $G$-representation $p^* V$ induced by the projection $p
\colon G \to J$.  As above, the map $N_+ \sma \bbk \to F( N_+, \bbk^{\mfib})$ 
induces a natural map
\begin{equation}\label{eq:borel-comparison}
\begin{aligned}
BY^{-V} \sma \bbk &\to   \ssF_{EY/N}(EY,\bbk)^{-V}/G  \\
EG_+ \sma_G (\Sigma^{-V} Y) \sma \bbk &\to 
EG_+ \wedge_G (\Sigma^{-V} (\theta^* Y \wedge_{\tilde{N}} F(N_+, \bbk^{\mfib})))
\\
&\to (\Sigma^{-V} \left(\theta^* (EG_+ \sma Y) \wedge_{\tilde{N}} F(N_+,
\bbk^{\mfib})\right)) / G
\end{aligned}
\end{equation}
where the notation $EY$ denotes $EG_+ \wedge Y$.  On the other hand,
using the identification $((EY)/N)/J \cong BY$, we have a map
\begin{equation}
BY^{-V} \sma \bbk = EG_+ \sma_G Y^{-V} \sma \bbk \to \ssF_{BY}( (EY)/N, \bbk)^{-V}/J,
\end{equation}
constructed as follows:  unpacking the notation, the target of our
desired map is 
\begin{equation}
\ssF_{BY}((EY)/N), \bbk)^{-V}/J = \left(\Sigma^{-V} ((EG_+ \sma_N Y))
\sma_{\tilde{J}} F(J_+, \bbk)\right)/J,
\end{equation}
where $\tilde{J}$ is the subgroup $\{(1,j) \,\mid\, j \in J\}$ of $J
\ltimes J \cong J \times J$.

The map in question is thus simply the map
\begin{equation}
BY^{-V} \sma \bbk = ((EG_+ \sma_N Y)/J)^{-V} \sma \bbk \to \left(\Sigma^{-V} ((EG_+ \sma_N
Y) \sma_{\tilde{J}} F(J_+, \bbk)\right)/J
\end{equation}
induced by the duality map $J_+ \sma \bbk \to F(J_+, \bbk^{\mfib})$
and passage to $J$-orbits.

Next, we consider the map
\begin{equation} \label{eq:commute_E-and_N-quotient}
  (EY)/N \equiv EG_+ \wedge_{N} Y   \to  EJ_+ \wedge_{N} Y/N \equiv E(Y/N)
\end{equation}
induced by the projections $EG \to EJ$ and $Y \to Y/N$, where the
second identification uses the fact that $N$ acts trivially on both $
EJ$ and $Y/N$. This induces the horizontal arrows in the following commuting diagram:
\begin{equation}
  \begin{tikzcd}
   \ssF_{BY}( (EY)/N, \bbk) \ar[r] \ar[d,"="]  & \ssF_{B(Y/N)}(E(Y/N) , \bbk ) \ar[d,"="] \\
 (EY)/N \wedge_{\tilde{J}} F( J_+, \bbk) \ar[r] &  E(Y/N)\wedge_{\tilde{J}} F( J_+, \bbk).
  \end{tikzcd}
\end{equation}
Similarly, the fact that $N$ acts trivially on $F(J_+,\bbk)$ yields an isomorphism
\begin{equation}
 (EY)/N \wedge_{\tilde{J}} F(J_+, \bbk) \to  EY \wedge_{\tilde{G}} F( J_+, \bbk),
\end{equation}
whose composition with the adjoint of the projection $G_+ \to J_+$ yields a map
\begin{equation}
  \begin{tikzcd}
    \ssF_{BY}(EY, \bbk) \ar[d,"="]  & \ar[d,"="]  \ssF_{BY}( (EY)/N, \bbk) \ar[l] \\
    EY \wedge_{\tilde{G}} F( G_+, \bbk) & \ar[l] (EY)/N \wedge_{\tilde{J}} F( J_+, \bbk).
  \end{tikzcd}
\end{equation}
Finally, we have a map 
\begin{equation}\label{eq:proj-rel-wirth}
    \ssF_{EY/N}(EY, \bbk)^{-V} \to  \ssF_{BY}(EY, \bbk)^{-V}.
  \end{equation}
Writing this out, this is the map 
\begin{equation}
EY \sma_{\tilde{N}} F(N_+, \bbk^{\mfib}) \cong EY \sma_G \left( G_+ \sma_{\tilde{N}} F(N_+,
\bbk^{\mfib})\right) \to EY \sma_{\tilde{G}} F(G_+, \bbk^{\mfib}) 
\end{equation}
induced by the relative Wirthmuller map from Equation~\eqref{eq:relwirth}: 
\begin{equation}
G_+ \sma_{\tilde{N}}  F( N_+, (\bbk)^{\mfib})  \to F( G_+, (\bbk)^{\mfib}).  
\end{equation}

Putting this all together, we have the following result summarizing the
situation, which expresses the functoriality of the these fiberwise
mapping spectra and records the fact that they have the correct
homotopy type pointwise.

\begin{lem} \label{lem:comparison-classifying-fibrewise-non-functorial}
The following diagram of orthogonal spectra commutes, and the arrows
between spectra in the first two rows are equivalences.
 \begin{equation} \label{eq:commutative-diagram-compare-BY-to-F_EY/N(EY)}
  \begin{tikzcd}
 BY^{-V} \wedge \bbk \ar[ddd] \ar[r] \ar[dr]   & \ssF_{BY}(EY, \bbk)^{-V}/G   &   BY^{-V}\wedge \bbk  \ar[l] \ar[d] \\
 & \ssF_{BY}( (EY)/N, \bbk)^{-V}/ J  \ar[dd]  \ar[u]  &  \ssF_{(EY)/N}(EY, \bbk)^{-V}/G \ar[ul]   \\
&  &  BY^{-V} \wedge \bbk  \ar[d] \ar[u] \ar[ul] \\
  BY/N^{-V}\wedge \bbk  \ar[r] & \ssF_{B (Y/N)}(E(Y/N) , \bbk )^{-V}/ J   &      BY/N^{-V}\wedge \bbk  \ar[l] 
  \end{tikzcd}
\end{equation} \qed
\end{lem}

We now turn to describing the functorial and multiplicative properties
of these constructions, which will allow us to assemble them into
spectral categories.

\subsubsection{Functoriality of fibrewise mapping spaces over $BY$}
\label{sec:funct-fibr-mapp}

To express the functoriality of the fiberwise mapping spaces and the
coherent comparisons that we require, it is useful to define a
category which encodes the key aspects of the functoriality.

\begin{defin}\label{defn:all-spectra-shifted}
Let $\Sp^{-\Vect}_{\eq}$ denote the category with 
\begin{enumerate}

\item Objects the triples $(G, Y, V)$, where $G$ is a finite
  group, $Y$ a cofibrant orthogonal $G$-spectrum, and $V$ a
  finite-dimensional $G$-inner product space.

\item Morphisms $f \colon (G,Y, V) \to (G', Y', V')$ given by 
\begin{enumerate} 
\item a surjection $f\colon G \to G'$ with kernel $G_f^\perp$ acting freely on
$Y$, 
\item a $G$-equivariant isometric embedding $V \to V'$ with cokernel $V^\perp_f$, 
\item and a $G$-equivariant map 
\begin{equation}
    S^{V^\perp_f} \sma Y \to f^* Y'.
\end{equation}
\end{enumerate}
\end{enumerate}
\end{defin}

\begin{rem}
Our interest in this category arises from the existence of a natural
functor from the category of cubes of Kuranishi charts
\begin{align}
  \square \Kur & \to \Sp^{- \Vect}_{\eq} \\
  \sigma & \mapsto (G_\sigma, X_\sigma|Z_\sigma, V_\sigma),
\end{align}  
as discussed in Equation~\eqref{eq:kur-all-spectra-compare}.
\end{rem}

The category $\Sp^{-\Vect}_{\eq} $ admits a natural monoidal
structure.

\begin{prop}
The category $\Sp^{-\Vect}_{\eq}$ is a symmetric monoidal category
with the product of $(G_0, Y_0, V_0)$ and $(G_1, Y_1, V_1)$ given by
\begin{equation}
(G_0 \times G_1, Y_0 \sma Y_1, V_0 \oplus V_1)
\end{equation}
and unit $(\{e\}, \ast, \{0\})$. \qed
\end{prop} 

We also need a refinement of the category $\Sp^{-\Vect}_{\eq}$.

\begin{defin}\label{defn:all-spectra-shifted-subgroup}
Let $\Sp^{- \Vect}_{\eq, \Sub}$ denote the category with:

\begin{enumerate}

\item Objects the quadruples $(G, Y, V,N)$ where $(G,Y,V)$ is an
  object of $\Sp^{-\Vect}_{\eq}$ and $N \subseteq G$ is a subgroup
  that acts freely on $Y/G_f^\perp$.

\item Morphisms $f \colon (G_0, Y_0, V_0, N_0) \to (G_1, Y_1, V_1,
  N_1)$ specified by a morphism $f\colon \Sp^{-\Vect}_{\eq}$ such that 
  $N_1 \subseteq f(N_0)$.

\end{enumerate}

\end{defin}

The category $\Sp^{-\Vect}_{\eq, \Sub}$ also admits a natural monoidal
structure extending the one on $\Sp^{-\Vect}_{\eq, \Sub}$.

\begin{prop}
The category $\Sp^{-\Vect}_{\eq}$ is a symmetric monoidal category
with the product of $(G_0, Y_0, V_0, N_0)$ and $(G_1, Y_1, V_1, N_0)$ given by
\begin{equation}
(G_0 \times G_1, Y_0 \sma Y_1, V_0 \oplus V_1, N_0 \times N_1)
\end{equation}
and unit $(\{e\}, \ast, \{0\}, \{0\})$.  The evident forgetful functor 
\begin{equation}
\Sp^{-\Vect}_{\eq, \Sub} \to \Sp^{-\Vect}_{\eq}
\end{equation}
is strong symmetric monoidal. \qed
\end{prop} 

We now begin to express the functoriality of our constructions in
terms of these indexing categories.  We begin with the Borel
construction, which plays a central role in the definition of the
virtual cochains.

\begin{defin}
The Borel construction defines a lax monoidal functor 
\begin{align}
  \Sp^{- \Vect}_{\eq} & \to \Sp   \\
  (G, Y, V) & \mapsto BY^{-V},
\end{align}
which corresponds to the first column of Equation
\eqref{eq:commutative-diagram-compare-BY-to-F_EY/N(EY)}, in the sense
that the arrow $  BY^{-V} \to B(Y/N)^{-V}$ is the image under this
functor of the morphism 
\begin{equation}
(G, Y, V) \to  (G/N, Y/N, V)
\end{equation}
specified by the canonical surjections $G \to G/N$ and $Y \to Y/N$
along with the map induced by the identity map $V \to V$.
\end{defin}

Our goal is to assemble the pointwise comparisons
\begin{equation}
BY^{-V} \to \ssF_{BY}(EY, \bbk)^{-V}/G
\end{equation}
into functorial comparisons.  However, consideration of the
functoriality of the assignment 
\begin{equation}
(G,Y,V) \mapsto \ssF_{BY}(EY, \bbk)^{-V}/G
\end{equation}
(as reflected in the second column in
\eqref{eq:commutative-diagram-compare-BY-to-F_EY/N(EY)}) makes it
clear that we need to enlarge the domain category from
$\Sp^{-\Vect}_{\eq}$ in order to capture the functoriality of the
comparison.  Thus, we use the twisted arrow category $\Tw
\Sp^{-\Vect}_{\eq}$.  In the background here is the comparison between
functors out of the twisted arrow category and spans; see
Section~\ref{sec:twisted}, specifically
Theorem~\ref{thm:twisted-to-span}.

First, using the evident projection map
\begin{equation}
   \Tw \Sp^{- \Vect}_{\eq}   \to \Sp^{-\Vect}_{\eq} 
 \end{equation}
which assigns to each arrow its source, we can regard the assignment
\begin{equation}
(f \colon (G, Y, V) \to (G', Y', V')) \mapsto BY^{-V}
\end{equation}
as a functor from $\Tw \Sp^{-\Vect}_{\eq}$ to spectra.  We now turn to
the fiberwise mapping spectra.

\begin{notn}
Given an arrow $f \co (G,Y,V) \to (G', Y',V')$, we write 
\begin{equation}
(EY)_f \equiv EY / G_f^\perp = EG_+ \sma_{G_f^{\perp}} Y.
\end{equation}
Via the isomorphism $G / G_f^\perp \cong G'$, the spectrum $(EY)_f$ admits a
natural action of $G'$, and we write  
\begin{equation}
BY_f = (EY)_f/ G'.
\end{equation}
\end{notn}

Given a $G$-representation $V$, we have an orthogonal $G'$-spectrum 
\begin{equation}
\begin{aligned}
\ssF_{(BY)_f}( (EY)_f, \bbk)^{-V} &\equiv \Sigma^{-V} (EY)_f \sma_{\widetilde{G'}}
F(G'_+, \bbk) \\
&\cong (\Sigma^{-V} EY)_{f} \sma_{\widetilde{G'}} F(G'_+, \bbk),
\end{aligned}
\end{equation}
where here we are regarding $V$ as a $G'$-representation via the
$G$-embedding $V \to V'$.  

We will consider the assignment 
\begin{equation} \label{eq:functor_twisted-arrow-to-fibrewise-over-BY}
  f \mapsto \ssF_{(BY)_f}( (EY)_f, \bbk)^{-V}/ G' \equiv \left((EY)_f
  \sma_{\widetilde{G'}} F(G'_+, \bbk)\right) / G'
\end{equation}
for each arrow $f$ in $\Sp^{-\Vect}_{\eq}$.  For the proof of the next
result, it is convenient to write $Y_f$ for the quotient
$Y/G_f^\perp$.  

\begin{lem}\label{lem:functor-twisted-arrow-objects}
Let $f_0$ and $f_1$ be maps in $\Sp^{-\Vect}_{\eq}$.  Each factorization
\begin{equation}
\begin{tikzcd}
(G_0, Y_0, V_0) \ar[r, "f_0"] \ar[d,"g"] & (G'_0, Y'_0, V'_0) \\
(G_1, Y_1, V_1)  \ar[r,"f_1"] & (G'_1, Y'_1, V'_1) \ar[u,"h"] 
\end{tikzcd}
\end{equation}
induces a map 
\begin{equation}
\ssF_{BY_0}((EY)_{f_0}  , \bbk)^{-V_{0}} \to
\ssF_{BY_1}((EY)_{f_1} , \bbk)^{-V_{1}}
\end{equation}
of orthogonal $G_1'$-spectra.
\end{lem}

\begin{proof} 
The map is a composition of several maps which we now define.  First,
since there is an isomorphism 
\begin{equation}
(EY)_{f_0} \cong (EY)_{f_1 \circ g} / G_h^\perp,
\end{equation}
we have an induced map
\begin{equation}\label{eq:first-map}
(EY)_{f_0} \sma_{\tilde{G'_0}} F(G'_{0,+}, \bbk) \to (EY)_{f_1
      \circ g} \sma_{\tilde{G'_1}} F(G'_{0,+}, \bbk) \to (EY)_{f_1
      \circ g} \sma_{\tilde{G'_1}} F(G'_{1,+}, \bbk), 
\end{equation}
which yields a map
\begin{equation}
  \ssF_{BY_0}((EY)_{f_0}  , \bbk)^{-V_{0}} \to \ssF_{B Y_0}((EY)_{f_1
    \circ g},\bbk)^{-V_{0}}
\end{equation}
after desuspension. Next, observe that the map $G_0 \to G_1$ induces a
surjection $G^\perp_{f_1 \circ g} \to G^{\perp}_{f_1}$; as in Equation
\eqref{eq:commute_E-and_N-quotient}, this yields a map 
\begin{equation}
\begin{aligned}
  (EY)_{f_1 \circ g} &= (EG_{0,+} \sma Y_0) / G_{f_1 \circ g}^{\perp}  \\ 
&\to \left((EG_{0,+} \sma Y_0) /  G_g^{\perp}\right) = ((EY)_g)_{f_1} 
  / G_{f_1}^{\perp}.
\end{aligned}
\end{equation}
Passing to fibrewise mapping spaces gives
\begin{equation}\label{eq:second-map}
 \ssF_{B Y_0}((EY)_{f_1 \circ g}  , \bbk)^{-V_{0}}  \to  \ssF_{B Y_0}((EY_g)_{f_1}, \bbk)^{-V_{0}}.
\end{equation}
Next, the splitting $V_{1} \cong  V_{0} \oplus V^\perp_g$ induces a map
\begin{equation}
  F(S^{V_0},  (EY_g)_{f_1}  \sma_{\tilde{G'_1}} F(G'_{1,+}, \bbk^{\mfib}))
  \to   F(S^{V_1}, S^{V^\perp_g} \sma (EY_g)_{f_1}  \sma_{\tilde{G'_1}}
  F(G'_{1,+}, \bbk^{\mfib})) 
\end{equation}
by smashing with $S^{V^\perp_g}$ on both sides.  Composing with the natural map
\begin{equation}
\begin{aligned}
  S^{V^\perp_g} \sma (EY_g)_{f_1} &= S^{V^\perp_g} \sma \left( (EG_{0,+}
  \sma Y_0) / G_g^{\perp} \right) / G_{f_1}^{\perp} \\
&\to 
(EG_{1,+} \sma S^{V^\perp_g} \sma (Y_0 / G_g^{\perp})) /
  G_{f_1}^{\perp} = (E (S^{V^\perp_g} \sma Y_g))_{f_1},  
\end{aligned}
\end{equation}
yields
\begin{equation}\label{eq:third-map}
   \ssF_{B Y_g}((EY_g)_{f_1}  , \bbk)^{-V_{0}} \to \ssF_{B (S^{V^\perp_g} \wedge Y_g)}(E(S^{V^\perp_g} \wedge Y_g)_{f_1}  , \bbk)^{-V_{1}}.
\end{equation}
Finally, we have a map $S^{V^\perp_g} \wedge Y_g \to Y_1$
induced by the map $S^{V^\perp_g} \sma Y_0 \to Y_1$ and
the fact that $G_g^{\perp}$ acts trivially on $Y_1$.  This map induces
a map
\begin{equation}\label{eq:fourth-map}
   \ssF_{B (S^{V^\perp_g} \wedge Y_g)}(E(S^{V^\perp_g} \wedge
   Y_g)_{f_1}  , \bbk)^{-V_{1}} \to \ssF_{BY_1}((EY)_{f_1} ,
   \bbk)^{-V_{1}}. 
\end{equation}
The desired map is then obtained as the composition of the maps
constructed above.
\end{proof}

The assignment of Lemma~\ref{lem:functor-twisted-arrow-objects} is
functorial; passing to orbits over $G_1'$ yields the following
essential result.

\begin{cor}\label{cor:functorial}
The assignment
\begin{equation}
f \mapsto    \ssF_{(BY)_f}( (EY)_f, \bbk)^{-V}/ G' \equiv \left((EY)_f
  \sma_{G'} F(G'_+, \bbk)\right) / G'
\end{equation}
defines a lax monoidal functor 
\begin{equation}
\Tw \Sp^{- \Vect}_{\eq}   \to \Sp.    
\end{equation}
\end{cor}

\begin{proof}
We begin by considering the composition in the case of the diagram
\begin{equation}\label{eq:twisted-composite}
\begin{tikzcd}
(G_0, Y_0, V_0) \ar[r, "f_0"] \ar[d,swap,"g"] & (G'_0, Y'_0, V'_0) \\
(G_1, Y_1, V_1)  \ar[d,swap,"\id"] \ar[r,"h \circ f_1"] & (G'_0, Y'_0, V'_0) \ar[u,swap,"\id"] \\
(G_1, Y_1, V_1)  \ar[r,"f_1"] & (G'_1, Y'_1, V'_1) \ar[u,swap,"h"], 
\end{tikzcd}
\end{equation}
which we need to compare with
\begin{equation}\label{eq:composite}
\begin{tikzcd}
(G_0, Y_0, V_0) \ar[r, "f_0"] \ar[d,swap,"g"] & (G'_0, Y'_0, V'_0) \\
(G_2, Y_2, V_2)  \ar[r,"f_1"] & (G'_2, Y'_2, V'_2) \ar[u,swap,"h"].
\end{tikzcd}
\end{equation}
The map determing by the bottom square in
Equation~\eqref{eq:twisted-composite} is the map 
\begin{equation}
\ssF_{BY_1}((EY)_{h \circ f_1}, \bbk)^{-V_{1}} \to \ssF_{B
  Y_1}((EY)_{f_1}, \bbk)^{-V_{1}} 
\end{equation}
specified in Equation~\eqref{eq:first-map} above; the other terms in
the composite of Lemma~\ref{lem:functor-twisted-arrow-objects} are the
identity.  In contrast, for the top square in
Equation~\eqref{eq:twisted-composite}, the arrow induced by
Equation~\eqref{eq:first-map} is the identity map.  The associated map
is the composite of the map 
\begin{equation}
\ssF_{B Y_0}((EY)_{f_0}, \bbk)^{-V_{0}} \to \ssF_{B
  Y_0}((EY_g)_{h \circ f_1}, \bbk)^{-V_{0}} 
\end{equation}
of Equation~\eqref{eq:second-map}, the map
\begin{equation}
\ssF_{B Y_g}((EY_g)_{h \circ f_1}  , \bbk)^{-V_{0}} \to \ssF_{B (S^{V^\perp_g}
  \wedge Y_g)}(E(S^{V^\perp_g} \wedge Y_g)_{h \circ f_1}  , \bbk)^{-V_{1}}
\end{equation}
of Equation~\eqref{eq:third-map}, and 
\begin{equation}
   \ssF_{B (S^{V^\perp_g} \wedge Y_g)}(E(S^{V^\perp_g} \wedge
   Y_g)_{h \circ f_1}  , \bbk)^{-V_{1}} \to \ssF_{BY_1}((EY)_{h \circ f_1} ,
   \bbk)^{-V_{1}} 
\end{equation}
of Equation~\eqref{eq:fourth-map}.
Checking that the composite coincides with the map associated to
Equation~\eqref{eq:composite} is now an exercise in commuting orbits
past each other.

Checking the compatibility of composition of a pair of maps of the
form $(\id, h_0)$ and $(\id, h_1)$ is a straightforward commutation of
orbits, as follows.  All maps are induced simply by
Equation~\eqref{eq:first-map} above.  So for the composite, we have
the maps
\begin{equation}
\ssF_{BY_0}((EY)_{f_0}, \bbk)^{-V_{0}} \to \ssF_{B
  Y_0}((EY)_{f_1}, \bbk)^{-V_{0}} 
\end{equation}
and
\begin{equation}
\ssF_{BY_0}((EY)_{f_1}, \bbk)^{-V_{0}} \to \ssF_{B
  Y_0}((EY)_{f_2}, \bbk)^{-V_{0}},
\end{equation}
which we are comparing to the map
\begin{equation}
\ssF_{BY_0}((EY)_{f_0}, \bbk)^{-V_{0}} \to
\ssF_{B Y_0}((EY)_{f_2}, \bbk)^{-V_{0}}.
\end{equation}
These coincide because the composite isomorphism
\begin{equation}
(EY)_{f_0} \cong (EY)_{f_1} / G_{h_0}^\perp \cong 
\left((EY)_{f_2} / G_{h_1}^\perp\right) / G_{h_0}^{\perp}
\end{equation}
is the same as the isomorphism $(EY)_{f_0} \cong (EY)_{f_2} / G_{h_0
  \circ h_1}^\perp$.

To check the composition associated to a pair of maps of the form
$(g_0, \id)$ and $(g_1, \id)$, we use the fact that given a composite
\begin{equation}
\begin{tikzcd}
(G_0, Y_0, V_0) \ar[r,"g_0"] & (G_1, Y_1, V_1) \ar[r,"g_1"] & (G_2, Y_2, V_2),
\end{tikzcd}
\end{equation}
the splittings $V_2 \cong V_0 \oplus V^{\perp}_{g_1 \circ g_0}$ and
$V_2 \cong V_1 \oplus (V_1)^{\perp}_{g_0} \cong V_0 \oplus
(V_0)^{\perp}_{g_1} \oplus (V_1)^{\perp}_{g_0}$
coincide in the sense that the data specifies a canonical
identification 
$V^{\perp}_{g_1 \circ g_0} \cong (V_0)^{\perp}_{g_1} \oplus
(V_1)^{\perp}_{g_0}$.

Checking the compatibility of the composition of the maps associated
to $(g,h)$ and $(\id,h')$ essentially follows
from the work we have already done, and so we can conclude that the
assignment indeed specifies a functor.  

The fact that this is a lax monoidal functor follows from
Lemma~\ref{lem:fibmapmult}, the compatibility of orbits with the
external monoidal product, and the compatibility of the splittings of
the representations with the product.
\end{proof}

We now realize the maps in the first two columns
of Diagram~\eqref{eq:commutative-diagram-compare-BY-to-F_EY/N(EY)} as the
following comparison of functors on the twisted arrow category.

\begin{lem} \label{lem:monoidal_comparison_classifying_fibrewise-1}
There is a zig-zag of lax monoidal natural transformations through
equivalences connecting the functors $BY^{-V} \sma \bbk$ and 
$ \ssF_{BY}( (EY)_f, \bbk)^{-V}/ G'$
on $\Tw \Sp^{- \Vect}_{\eq}$.
\end{lem}

\begin{proof}
First, we construct the pointwise comparison and show it is a weak
equivalence.  For an object $\left(f \colon (G_0, Y_0, V_0) \to (G_0',
Y_0', V_0')\right)$ in $\Tw \Sp^{-\Vect}_{\eq}$, there is a natural
weak equivalence
\begin{equation}
BY_0^{-V} \sma \bbk \to \ssF_{BY_0}( (EY)_f, \bbk)^{-V}/ G_0',
\end{equation}
induced as follows.  Writing this out, this is a weak equivalence
\begin{equation}
(EG_{0,+} \sma (\Sigma^{-V_0} Y_0)) / G_0 \to
\left( \left( \left(EG_{0,+} \sma (\Sigma^{-V_0} Y_0) \right) /
G_f^{\perp}\right) \sma_{G_0'} F(G_{0,+}', \bbk)\right) / G_0'. 
\end{equation}
Recall from Equation~\eqref{eq:borel-comparison} that we can construct
this weak equivalence using the isomorphism $(Y_0 / G_f^\perp) / G_0'
\cong Y_0 / G_0$.

Next, we show that given a morphism
\begin{equation}\label{eq:composite2}
\begin{tikzcd}
(G_0, Y_0, V_0) \ar[r, "f_0"] \ar[d,swap,"g"] & (G'_0, Y'_0, V'_0) \\
(G_1, Y_1, V_1)  \ar[r,"f_1"] & (G'_1, Y'_1, V'_1) \ar[u,swap,"h"].
\end{tikzcd}
\end{equation}
the diagram
\begin{equation}
\begin{tikzcd}
BY_0^{-V_0} \ar[r] \ar[d] & \ssF_{BY_0}( (EY)_{f_0}, \bbk)^{-V_0}/ G_0' \ar[d]
\\
BY_1^{-V_1} \ar[r] & \ssF_{BY_1}( (EY)_{f_1}, \bbk)^{-V_1}/ G_1' 
\end{tikzcd}
\end{equation}
commutes.  To see this, we use the description of natural
transformations from
Lemma~\ref{lem:twisted-arrow-natural-transformation}.  From this 
perspective, the data of the natural transformation is specified by
commutative diagrams
\begin{equation}
\begin{tikzcd}
BY_0^{-V_0} \ar[r] \ar[dr] & \ssF_{(BY)_{f_0}}((EY)_{f_0},
\bbk)^{-V_0}/G_0' \ar[d] \\
& \ssF_{(BY)_{f_1 \circ g}}((EY)_{f_1 \circ g},
\bbk)^{-V_0}/G_1'
\end{tikzcd}
\end{equation}
associated to the factorization $f_0 = g \circ (h\circ f_1) \circ \id$
and
\begin{equation}
\begin{tikzcd}
BY_0^{-V_0} \ar[r] \ar[d] & \ssF_{(BY)_{f_1 \circ g}}((EY)_{f_1 \circ g},
\bbk)^{-V_0}/G_1' \ar[d] \\
BY_1^{-V_1} \ar[r] & \ssF_{(BY)_{f_1}}((EY)_{f_1},
\bbk)^{-V_1}/G_1'
\end{tikzcd}
\end{equation}
associated to the factorization $(f_1 \circ g) = \id \circ f_1 \circ g$.

Writing it out, the top square is the diagram
\begin{equation}
\begin{tikzpicture}
\node[scale=.85]{
\begin{tikzcd}
(EG_{0,+} \sma (\Sigma^{-V_0} Y_0)) / G_0 \ar[r] \ar[dr] & 
\left( \left( \left(EG_{0,+} \sma (\Sigma^{-V_0} Y_0) \right) /
G_f^{\perp}\right) \sma_{G_0'} F(G_{0,+}', \bbk)\right) / G_0' \ar[d]
\\
& \left( \left( \left(EG_{0,+} \sma (\Sigma^{-V_0} Y_0) \right) /
G_{f_1 \circ g}^{\perp}\right) \sma_{G_1'} F(G_{1,+}', \bbk)\right) / G_1',
\end{tikzcd}
};
\end{tikzpicture}
\end{equation}
where the vertical map is induced by Equation~\eqref{eq:first-map}.
To see that this commutes, note that the vertical map is induced by the
isomorphism $(EY)_{f_0} \cong (EY)_{f_1 \circ g} / G_h^\perp$ and the
horizontal maps by the isomorphisms 
\begin{equation}
(Y_0 / G_{f_0}^\perp) / G_0' \cong
Y_0 / G_0 \quad\textrm{and}\quad (Y_0 / G_{f_1 \circ g}^\perp) / G_1'
\cong Y_0 / G_0.
\end{equation}

The bottom square is the diagram 
\begin{equation}
\begin{tikzpicture}
\node[scale=.85]{
\begin{tikzcd}
(EG_{0,+} \sma (\Sigma^{-V_0} Y_0)) / G_0 \ar[r] \ar[d] & 
\left( \left( \left(EG_{0,+} \sma (\Sigma^{-V_0} Y_0) \right) /
G_{f_1 \circ g}^{\perp}\right) \sma_{G_1'} F(G_{1,+}', \bbk)\right) /
G_1' \ar[d]
\\
(EG_{1,+} \sma (\Sigma^{-V_1} Y_1)) / G_1 \ar[r]
& \left( \left( \left(EG_{1,+} \sma (\Sigma^{-V_1} Y_1) \right) /
G_{f_1}^{\perp}\right) \sma_{G_1'} F(G_{1,+}', \bbk)\right) / G_1'.
\end{tikzcd}
};
\end{tikzpicture}
\end{equation}
The horizontal maps are induced by the isomorphisms
\begin{equation}
(Y_0 / G_{f_1 \circ g}^\perp) / G_1'
\cong Y_0 / G_0
\quad\textrm{and}\quad 
(Y_0 / G_{f_1}^\perp) / G_1' \cong
Y_1 / G_1.
\end{equation}
and the vertical maps are clearly compatible with these maps.
\end{proof}

The same arguments also establish the following small extension of
this result.

\begin{lem}\label{lem:monoidal_comparison_classifying_fibrewise-2}
There is a zig-zag of lax monoidal natural transformations through
equivalences
\begin{equation}
BY^{-V_f} \sma \bbk \htp \ssF_{BY}( (EY)_f, \bbk)^{-V_f}/ G'
\end{equation}
on $\Tw \Sp^{- \Vect}_{\eq}$.
\end{lem}

\subsubsection{Functoriality of fibrewise mapping spaces over intermediate quotients}
\label{sec:funct-fibr-mapp-1}

We now expand the discussion above to the setting of the twisted arrow
category $\Tw \Sp^{-\Vect}_{\Sub}$.  
Given an arrow $(G, Y, V, N) \to (G', Y', V', N')$ in
$\Sp^{-\Vect}_{\Sub}$, we consider the orthogonal $G'$-spectrum
\begin{equation} \label{eq:fibrewise-mapping-spectra-functor}
  f \mapsto \ssF_{EY_f/N'}(EY_f, \bbk)^{-V} \equiv (EY)_f
  \sma_{\widetilde{N'}} F(N'_+, \bbk)
\end{equation}
where we again write $Y_f$ for the quotient of $Y$ with respect to the
kernel $G^\perp_f$ of $G \to G'$.

We again will show that this assignment induces a functor from the
twisted arrow category.

\begin{lem} \label{lem:fibrewise-mapping-spectra-functorial}
  Each factorisation
  \begin{equation}
    \begin{tikzcd}
       (G_0, Y_0, V_0,  N_0) \ar[r, "f_0"] \ar[d,"g"] & (G'_0, Y'_0, V'_0,  N'_0) \\
      (G_1, Y_1, V_1, N_1)  \ar[r,"f_1"] & (G'_1, Y'_1, V'_1, N'_1) \ar[u,"h"] 
    \end{tikzcd}
  \end{equation}
  induces a natural map of $G'_1$-spectra:
  \begin{equation} \label{eq:map_internal-mapping-spaces-twisted-arrow-S^0-A}
    \ssF_{EY_{f_0}/N'_0}(EY_{f_0} , \bbk)^{-V_{f_0}} \to  \ssF_{EY_{f_1}/N'_1}(EY_{f_1}, \bbk)^{-V_{f_1}}     
  \end{equation}
\end{lem}
\begin{proof}
We have a natural map $V_{f_0} \to V_{f_1}$, with complement
$V_{g}^{\perp}$. By adjunction, it suffices to construct a map
\begin{equation}
S^{V_g^\perp} \wedge \ssF_{EY_{f_0}/N'_0}(EY_{f_0} , \bbk) \to
\ssF_{EY_{f_1}/N'_1}(EY_{f_1}, \bbk).
  \end{equation}
Recall that $EY_{f_0} \cong EY_{0} / G_{f_0}^{\perp}$. Factoring $f_0 =
h \circ f_1 \circ g$, we obtain the isomorphism $EY_{f_0} \cong
EY_{g}/ G_{h \circ f_1}^{\perp}$. Since the data of the map 
$g$ specifies a map $S^{V^\perp_g} \wedge (EY)_g \to EY_1$, we have a
natural map 
\begin{equation}\label{eq:map_desuspended_fibrewise-spectra-factor-thickening}
\begin{tikzcd}
S^{V_g^\perp} \sma    \ssF_{(EY)_{f_0}/N'_0}((EY)_{f_0} , \bbk) \ar[d] \\
\ssF_{(S^{V_g^\perp} \sma (EY)_{g})/(h \circ
f_1)^{-1}N'_0}((S^{V_g^\perp} \sma  (EY)_{g} )/ G_{h \circ
f_1}^{\perp} , \bbk). 
\end{tikzcd}
\end{equation}
Thus, it suffices to construct a map from the target to
$\ssF_{(EY)_{f_1}/N'_1}((EY)_{f_1}, \bbk)$. 
   
We have a natural map
\begin{equation}
S^{V_g^\perp} \wedge (EY)_{g} \to EY_{1} \to (EY)_{f_1},
\end{equation}
so it suffices to show that this map descends to the quotients. This
follows from the inclusion $h^{-1}(N'_0) \subset N'_1$, and the
isomorphisms:
\begin{equation}
(EY)_{f_1} /N'_1 \equiv EY_{1}/ f_1^{-1}(N'_1) \qquad
  (EY)_{f_0}/N'_0 \equiv (EY)_g/ (g \circ f_1)^{-1}(N'_0).
\end{equation}
The result is a map
\begin{equation} 
(S^{V_g^\perp} \wedge (EY)_{g})/(h \circ f_1)^{-1}N'_0  \to
  (S^{V_g^\perp} \wedge (EY)_{g})/ f_1^{-1} (N'_1)   \to
  (EY)_{f_1}/N'_1.  
\end{equation}
The first map induces
\begin{equation} \label{eq:first-map-functoriality-fibrewise}
\begin{tikzcd}
\ssF_{((S^{V_g^\perp} \wedge (EY)_{g})/ (h \circ
  f_1)^{-1}N'_0}((S^{V_g^\perp} \wedge (EY)_{g})/  G_{h \circ
  f_1}^{\perp} , \bbk) \ar[d] \\ 
 \ssF_{((S^{V_g^\perp} \wedge (EY)_{g})/ f_1^{-1} N'_1 }((S^{V_g^\perp}
 \wedge Y_{g})/  G_{h \circ f_1}^{\perp} , \bbk) , 
\end{tikzcd}
\end{equation}
while the identification $(S^{V_g^\perp} \wedge (EY)_{g})/ f_1^{-1}
(N'_1) = (S^{V_g^\perp} \wedge (EY)_{g})/ G_f^\perp)/ N'_1$, together
with the second map induce 
\begin{equation} \label{eq:second-map-functoriality-fibrewise}
\ssF_{((S^{V_g^\perp} \wedge (EY)_{g})/ f_1^{-1} N'_1}((S^{V_g^\perp}
\wedge (EY)_{g})/ G_{f_1}^\perp , \bbk) \to
\ssF_{(EY)_{f_1}/N'_1}((EY)_{f_1}, \bbk). 
\end{equation}

To complete the comparison, the evident inclusion $G_{f_1}^{\perp} \to
G_{h \circ f_1}^{\perp}$  yields a projection 
\begin{equation}
(S^{V_g^\perp} \sma (EY)_{g})/ G_{f_1}^\perp \to  (S^{V_g^\perp}
  \sma (EY)_{g}) / G_{h \circ f_1}^{\perp}, 
\end{equation}
so that we have a map
\begin{equation}  \label{eq:third-map-functoriality-fibrewise}
\begin{tikzcd}
\ssF_{((S^{V_g^\perp} \sma (EY)_{g})/ f_1^{-1} N'_1}((S^{V_g^\perp}
\sma (EY)_{g})/ G_{h \circ f_1}^{\perp} , \bbk) \ar[d] \\
\ssF_{((S^{V_g^\perp} \sma (EY)_{g})/ f_1^{-1} N'_1}((S^{V_g^\perp}
\wedge (EY)_{g})/ G_{f_1}^\perp , \bbk).
\end{tikzcd} 
\end{equation}
Altogether, the composition of Equations
\eqref{eq:first-map-functoriality-fibrewise},
\eqref{eq:second-map-functoriality-fibrewise}, and
\eqref{eq:third-map-functoriality-fibrewise}, together with Equation
\eqref{eq:map_desuspended_fibrewise-spectra-factor-thickening} yields
the result.
\end{proof}

The relative map constructed in
Lemma~\ref{lem:fibrewise-mapping-spectra-functorial} is compatible
with the absolute map in the following sense.

\begin{lem}\label{lem:relative-composition-compat}
Given a factorization
\begin{equation}
\begin{tikzcd}
(G_0, Y_0, V_0,  N_0) \ar[r, "f_0"] \ar[d,swap,"g"] & (G'_0, Y'_0, V'_0,  N'_0) \\
(G_1, Y_1, V_1, N_1)  \ar[r,"f_1"] & (G'_1, Y'_1, V'_1, N'_1) \ar[u,swap,"h"] 
\end{tikzcd}
\end{equation}
in which $N_0 = N_1 = N_0' = N_1' = \{e\}$, the associated morphism of 
spectra from Lemma~\ref{lem:fibrewise-mapping-spectra-functorial}
coincides with the morphism of
Lemma~\ref{lem:functor-twisted-arrow-objects}.
\end{lem}

Passing to quotients and using essentially the same arguments as in
the proof of Corollary~\ref{cor:functorial}, we obtain the following
proposition recording the functoriality of the relative construction.

\begin{lem}\label{lem:functor-fiberwise-maps-Borel}
The assignment 
\begin{equation} \label{eq:functor-fibrewise-maps-Borel}
\begin{aligned}
  \Tw \Sp^{- \Vect}_{\eq,\Sub}  & \to \Sp \\
   f & \mapsto \ssF_{(EY)_f/N'}((EY)_f, \bbk)^{-V_f}/G'.
\end{aligned}
\end{equation}
specifies a lax monoidal functor. \qed
\end{lem}

Now, using the forgetful functor
\begin{equation}
  \Tw \Sp^{- \Vect}_{\eq,\Sub}  \to  \Tw \Sp^{- \Vect}_{\eq},
\end{equation}
we can pull back the functor $\ssF_{BY}( (EY)_f, \bbk)^{-V}/ G'$ from
Equation~\eqref{eq:functor_twisted-arrow-to-fibrewise-over-BY}.  This
permits the following comparison.

\begin{prop}\label{prop:relative-comparison}
There is a lax monoidal equivalence
\begin{equation}
\ssF_{(EY)_f/N'}((EY)_f, \bbk)^{-V_f}/G'   \Rightarrow  \ssF_{BY}( (EY)_f, \bbk)^{-V_f}/ G' 
\end{equation}
of functors 
\begin{equation}
\Tw \Sp^{- \Vect}_{\eq,\Sub} \to \Sp.
\end{equation}
\end{prop}

\begin{proof}
Working pointwise, we have a map
\begin{equation}
\begin{aligned}
\ssF_{(EY)_f/N'}((EY)_f, \bbk)^{-V_f} &\to \ssF_{BY}( (EY)_f,
\bbk)^{-V_f} \\
\Sigma^{-V_f} (EY)_f \sma_{\tilde{N'}} F(N'_+, \bbk) &\to \Sigma^{-V_f} (EY)_f
\sma_{\tilde{G'}} F(G'_+, \bbk)
\end{aligned}
\end{equation}
induced as in the construction of the map in
Equation~\eqref{eq:proj-rel-wirth}; passing to orbits yields the
transformation of the statement.  The commutative diagram
\begin{equation}
\begin{tikzcd}
\Sigma^{-V_f} (EY)_f \sma_{\tilde{N'}} F(N'_+, \bbk) \ar[r] & \Sigma^{-V_f} (EY)_f
\sma_{\tilde{G'}} F(G'_+, \bbk)  \\
\Sigma^{-V_f} (EY)_f \ar[r,swap,"\id"] \ar[u,"\htp"] & \Sigma^{-V_f} (EY)_f \ar[u,swap,"\htp"]
\end{tikzcd}
\end{equation}
shows that the comparison is a pointwise equivalence.  Moreover, it is
straightforward to check that it is compatible with the external
monoidal structure.  To see that this is a natural transformation, we
need to verify that given a  factorization $f_0 = h \circ f_1 \circ g$ in
$\Sp^{-\Vect}_{\eq, \Sub}$, the associated diagram
\begin{equation}
\begin{tikzcd}
\ssF_{(EY_0)_{f_0}/N'_0}((EY_0)_{f_0}, \bbk)^{-V_{f_0}}/G_0' \ar[r]
\ar[d] & \ar[d] \ssF_{BY_0}( (EY_0)_{f_0}, \bbk)^{-V_{f_0}}/ G_0' 
\\
\ssF_{(EY_1)_{f_1}/N_1'}((EY_1)_{f_1}, \bbk)^{-V_{f_1}}/G_1' \ar[r] & 
\ssF_{BY_1}( (EY_1)_{f_1}, \bbk)^{-V_{f_1}}/ G_1' 
\end{tikzcd}
\end{equation} 
commutes.  When the subgroups $N_i$ and $N_i'$ are trivial, this
reduces to the absolute case.  The result now follows from the
compatibility of the relative Wirthmuller maps with composition.
Specifically, expanding the notation we have the diagram
\begin{equation}
\begin{tikzcd}
\Sigma^{-V_{f_0}} (EY_0)_{f_0} \sma_{\tilde{N_0'}} F(N'_{0,+}, \bbk)/G_0' \ar[r]
\ar[d] & \ar[d] \Sigma^{-V_{f_0}} (EY_0)_{f_0} \sma_{\tilde{G_0'}} F(G'_{0,+}, \bbk) / G_0'
\\
\Sigma^{-V_{f_1}} (EY_1)_{f_1} \sma_{\tilde{N_1'}} F(N'_{1,+}, \bbk)/G_1'
\ar[r] & \Sigma^{-V_{f_1}} (EY_1)_{f_1} \sma_{\tilde{G_1'}} F(G'_{1,+}, \bbk) / G_1',
\end{tikzcd}
\end{equation} 
and one concludes that this commutes since the diagram
\begin{equation}
\begin{tikzcd}
F(N'_{0,+}, \bbk)/G_0' \ar[r] \ar[d] & \ar[d] F(G'_{0,+}, \bbk) / G_0'
\\
F(N'_{1,+}, \bbk)/G_1' \ar[r] & F(G'_{1,+}, \bbk) / G_1'
\end{tikzcd}
\end{equation} 
commutes.
\end{proof}

\begin{rem}
To connect to our high-level outline, note that the combination of the
above result with Lemma~\ref{lem:monoidal_comparison_classifying_fibrewise-1}
effectively yields a lift of
Lemma~\ref{lem:comparison-classifying-fibrewise-non-functorial} to the
level of lax monoidal functors.
\end{rem}

Finally, we combine this functor with the monoidal fibrant replacement
functor.

\begin{lem} \label{lem:genu-fibr-replacement-equivalence}
The assignment
\begin{align} \label{eq:quotient-fibrant-resolution-twisted-arrow-quot-spectra}
   f & \mapsto \cQ_{\cU(G')}\ssF_{(EY)_f/N'}((EY)_f,
   \bbk)^{-V_f}/G'
\end{align}
specifies a lax monoidal functor $\Tw \Sp^{- \Vect}_{\eq,\Sub} \to \Sp$
equipped with a natural equivalence
  \begin{equation}
  \ssF_{(EY)_f/N'}((EY)_f, \bbk)^{-V_f}/G' \to
  \cQ_{\cU(G')}\ssF_{(EY)_f/N'}((EY)_f, \bbk)^{-V_f}/G'
  \end{equation}
\end{lem}

\begin{proof}
The properties of the comparison follow from the argument for the
functoriality of fibrant replacement
(Proposition~\ref{prop:fibrant-twisted-functoriality}) and the natural
transformation $\Id \to \cQ_{\cU(-)}$.  The key point is that in the
construction of
Equation~\eqref{eq:map_internal-mapping-spaces-twisted-arrow-S^0-A}, the
domain inherits a $G_1'$-structure via pullback along the given
projection $G'_1 \to G'_0$.  As a consequence, this comparison is
compatible with the functoriality of the fibrant replacement functor,
and the desired result follows upon passage to orbits. 
\end{proof}

\subsection{The Adams isomorphism via bi-equivariant spectra}
\label{sec:adams-isom-via}

In this section, we continue the comparison to the virtual cochains by
applying a suitably functorial and multiplicative model of the Adams
isomorphism.  The key observation is that the model of the Adams
isomorphism constructed as a zig-zag in
Definition~\ref{def:Adams-pointset-biequivariant} becomes functorial
with respect to the choice of group after passing to the appropriate twisted arrow
category.

\subsubsection{Right fixed points and left quotients}
\label{sec:right-fixed-points}

Recall from the discussion preceding
Definition~\ref{def:Adams-pointset-biequivariant} that a biequivariant
$G$-spectrum is an orthogonal spectrum with commuting left and right
$G$ actions.  We continue to write $X_r$ when
considering the right action and abusively omit any subscript from the
notation when considering the left action.  Note that the 
category of spectra with left action is naturally isomorphic to the
category of spectra with right action, with 
isomorphism provided by the identity on the underlying spectrum, and
the action obtained by precomposing the original action with the
inverse map $G \to G$.  Equivalently, a biequivariant orthogoanl
$G$-spectrum is an orthogonal $G \times G$-spectrum.  

Clearly, the notion of biequivariance is compatible with the external
smash product (recall Section~\ref{sec:ortho}):
\begin{lem}
If $X_1$ is a biequivariant $G_1$-spectrum and $X_2$ is a
biequivariant $G_2$-spectrum, then $X_1 \sma X_2$ is a biequivariant
$G_1 \times G_2$-spectrum.
\end{lem}

Given a $G$-spectrum $X$, the spectrum $G_+ \sma X$ is naturally a
biequivariant spectrum with the diagonal left action and the right
action which is trivial on $X$.  The natural map 
\begin{equation}
G_+ \sma X \to (G_+ \sma X)_r
\end{equation}
specified by $(g,x) \mapsto (g,g^{-1}x)$ is a weak equivalence of
$G$-spectra.  Similarly, the composite 
\begin{equation}
  G_+ \sma X \to X \wedge S^0 \to X
\end{equation}
is a map of bi-equivariant spectra, where we give the target $X$
the trivial right $G$-action. 

To apply this observation in our setting, we introduce the
following definitions.  

\begin{defin}
For an arrow
\begin{equation}
f \co (G, Y, V,N) \to (G', Y', V',N')
\end{equation}
in $\Sp^{- \Vect}_{\eq,\Sub}$, we define the biequivariant orthogonal
$G'$-spectrum 
\begin{equation}
(EY)_f  \sma_{\tilde{N'}} F(N'_+, \bbk^\mfib)  \sma F(G'_{r,+}, \bS^\mfib) 
\end{equation}
where the left $G'$-action is given by the diagonal and the right
$G'$-action is non-trivial only on $F(G'_{r,+}, \bS^\mfib)$.  
\end{defin}

Desuspending as before, we have:

\begin{defin}
For an arrow $f$ in $\Sp^{- \Vect}_{\eq,\Sub}$, we define
\begin{multline}
   \ssF_{(EY)_f/N' \wedge G'/N'}((EY)_f \wedge G'_+, \bbk)^{-V_f}
   \equiv \\ F(S^{V_f}, ((EY)_f  \sma_{\tilde{N'}} F(N'_+, \bbk^\mfib)
   \wedge F(G'_{r,+}, \bS^\mfib))^{\mfib}) 
 \end{multline}
as a biequivariant orthogonal $G'$-spectrum.
\end{defin}

Analogues of the arguments for
Lemma~\ref{lem:functor-fiberwise-maps-Borel} shows that this
assignment of spectra assembles into a functor on the twisted arrow
category $\Tw \Sp^{-\Vect}_{\eq,\Sub}$: 
\begin{lem}
The assignment 
\begin{equation}
f \mapsto \ssF_{(EY)_f/N' \wedge G'/N'}((EY)_f \wedge G'_+,
\bbk)^{-V_f}
\end{equation}
specifies a functor 
\begin{equation}
\Tw\Sp^{-\Vect}_{\eq,\Sub} \to \Sp.
\end{equation} \qed
\end{lem}

We then have the following comparison.

\begin{lem}\label{lem:compare}
There is a natural biequivariant map
\begin{equation}
\begin{tikzcd}
\cQ_{\cU(G')}\ssF_{(EY)_f/N' \wedge G'/N'}((EY)_f \sma G',
        \bbk)^{-V_f} \ar[d] \\
F(G',\cQ_{\cU(G'_r)} \ssF_{(EY)_f/N'}((EY)_f, \bbk)^{-V_f}) 
\end{tikzcd}
\end{equation}
that induces a natural transformation of functors from
$\Tw\Sp^{-\Vect}_{\eq,\Sub}$.
\end{lem}

\begin{proof}
By adjunction, we obtain a map
\begin{equation}
  \begin{aligned}
    \ssF_{(EY)_f/N' \wedge G'/N'} & ((EY)_f \wedge G'_+, \bbk)^{-V_f}
    \\ &= F(S^{V_f}, ( (EY)_f \sma_{\tilde{N'}} F(N'_+, \bbk^\mfib) \wedge F(G'_+, \bS^\mfib))^\mfib) \\
    & \to F(G'_+,  F(S^{V_f},  ((EY)_f  \sma_{\tilde{N'}} F(N'_+, \bbk^\mfib))^\mfib) \\
    & = F(G'_+, \ssF_{(EY)_f/N'}((EY)_f, \bbk)^{-V_f}).
  \end{aligned}
\end{equation}
The desired map now arises by applying the canonical map
  \begin{equation}
    \cQ_{\cU(G')} F(A,B) \to F(A, \cQ_{\cU(G')} B)    
  \end{equation}
described in Equation~\eqref{eq:Qassem}.  This comparison clearly has
the desired functoriality, using the proof of 
Proposition~\ref{prop:fibrant-twisted-functoriality}.
\end{proof}

Given any biequivariant spectrum $X$, we will write $X_{G_r}$ to
denote the quotient with respect to the right $G$-action.  Regarding
$G$ as the subgroup $1 \times G \subset G \times G$ makes it evident
that $X_{G_r}$ has a residual left $G$-action.  Analogously, $X^G$ has
a residual right $G$-action.  The construction 
\begin{equation}
X \mapsto (X^G)_{G_r}
\end{equation}
that first passes to fixed points with respect to the left $G$-action
and then quotients with respect to the residual right $G$-action then
yields a spectrum.  The following proposition records the
functoriality of this construction in our context.

\begin{prop}
For an arrow $f \co (G, Y, V,N) \to (G', Y', V',N')$ in $\Sp^{-
  \Vect}_{\eq,\Sub}$, the assignments 
\begin{equation}
f \mapsto ( \cQ_{\cU(G')}\ssF_{(EY)_f/N' \wedge G'/N'}((EY)_f \wedge G', \bbk)^{-V_f})^{G'})/G'_r
\end{equation}
and
\begin{equation}
f \mapsto \cQ_{\cU(G')}\ssF_{(EY)_f/N'}((EY)_f, \bbk)^{-V_f})/G'
\end{equation}
specify lax monoidal functors $\Tw \Sp^{- \Vect}_{\eq,\Sub} \to \Sp$. \qed
\end{prop}

The functoriality and multiplicativity of the comparison map in
Lemma~\ref{lem:compare} then yields the following.

\begin{cor} \label{cor:right-fixed-points-project}
There is a lax monoidal natural transformation
\begin{align}
  ( \cQ_{\cU(G')}\ssF_{(EY)_f/N' \wedge G'/N'}((EY)_f \wedge G',
  \bbk)^{-V_f})^{G'})/G'_r \\  
\Rightarrow  \cQ_{\cU(G')}\ssF_{(EY)_f/N'}((EY)_f, \bbk)^{-V_f})/G'
\end{align}
of functors from $\Tw \Sp^{- \Vect}_{\eq,\Sub} \to \Sp$.
\end{cor}

\begin{proof}
For any orthogonal $G'$-spectrum $A$, we have a natural level-wise
homeomorphism 
\begin{equation}
\left(F(G',A)^{G'}\right)_r \cong A.
\end{equation}
Applying this pointwise to the functors in question yields the
homeomorphism
\begin{equation}
(F(G', \cQ_{\cU(G')}\ssF_{(EY)_f/N'}((EY)_f, \bbk)^{-V_f})^{G'})_r \cong
  \cQ_{\cU(G')}\ssF_{(EY)_f/N'}((EY)_f, \bbk)^{-V_f}
\end{equation}
which implies the result by passing to quotients; the desired
functoriality is evident.
\end{proof}
 
\subsubsection{Swapping the order of fixed points and quotients}
\label{sec:swapping-order-fixed}

For any biequivariant orthogonal $G$-spectrum $X$, the interchange of
colimits and limits yields a natural comparison map.

\begin{lem}
Let $X$ be a biequivariant $G$-spectrum.  Then there is a natural map
of spectra
\begin{equation}
(X^G)_{G_r} \to (X_{G_r})^G.
\end{equation}
\end{lem}

To construct a derived version of this interchange map, we need to
consider the biequivariant structure of the fibrant replacement
functor.

\begin{lem}
Let $X$ be a biequivariant orthogonal $G$-spectrum.  Let $U$ be a
$G$-universe, considered as a biequivariant universe having trivial
right $G$-action.  Then the fibrant replacement $\cQ_U X$ is a
biequivariant orthogonal $G$-spectrum. 
\qed
\end{lem}

Because the right action on $U$ is trivial, we have the following
commutation result.

\begin{lem}\label{lem:orbits-and-Q}
Let $X$ be a biequivariant orthogonal $G$-spectrum and $U$ a
$G$-universe.  Then there is a natural map 
\begin{equation}
(\cQ_U X) / G_r \to \cQ_U (X/G_r)
\end{equation}
of orthogonal $G$-spectra.  The map is externally multiplicative in
the sense that for a $G_1$-spectrum $X_1$, $G_1$-universe $U_1$,
$G_2$-spectrum $X_2$, and $G_2$-universe $U_2$, the diagram
\begin{equation}
\begin{tikzcd}
(\cQ_{U_1} X_1) / (G_1)_r \sma (\cQ_{U_2} X_2) / (G_2)_r \ar[r] \ar[d] &
  \cQ_{U_1} (X_1/(G_1)_r) \sma \cQ_{U_2} (X_2/(G_2)_r) \ar[d] \\
\cQ_{U_{12}} ((X_1 \sma X_2)) / (G_1 \times G_2)_r) \ar[r] &
  \cQ_{U_{12}} ((X_1 \times X_2)/(G_1 \times G_2)_r) 
\end{tikzcd}
\end{equation}
commutes.
\end{lem}

\begin{proof}
Whenever $A$ and $B$ are $G$-spectra so that $A$ has trivial
$G$-action, there are natural maps
\begin{equation}
F(A,B)/G \to F(A,B/G) \qquad\textrm{and}\qquad A \sma_G B \cong A \sma
(B/G).
\end{equation}
Specializing to the case of $\cQ_U X$, this implies that for each $V$
there are natural maps 
\begin{equation}
\begin{tikzcd}
\displaystyle{\left(\hocolim_{W \in U} \Omega^{W \otimes V} X((W \oplus \bR) \otimes
V) \right) / G_r }\ar[d] \\
\displaystyle{\hocolim_{W \in U} \Omega^{W \otimes V} \left( X((W
\oplus \bR) \otimes V) / G_r\right),}
\end{tikzcd}
\end{equation}
using the fact that orbits commute with homotopy colimits and that the
right action of $G$ is trivial on the representations $V$ and $W$.
These maps are compatible with the structure maps and so assemble into
maps of orthogonal $G$-spectra.  The external multiplicativity now
follows from the monoidality of orbits with respect to the external
product and the lax monoidal structure of the fibrant replacement
functor.
\end{proof}

Next, we consider the quotient
\begin{equation}
\left((EY)_f  \sma_{\tilde{N'}} F(N'_+, \bbk^\mfib)  \sma F(G'_{r,+},
\bS^\mfib)\right) /G'_r. 
\end{equation}
Since $G'_r$ only acts on the smash factor $F(G'_{r,+}, \bS^\mfib)$, 
there is an equivalence
\begin{multline}
\ssF_{(EY)_f/N'}((EY)_f, \bbk) \equiv (EY)_f \sma_{\tilde{N'}} F(N'_{+},
\bbk^\mfib) \\
\to \left((EY)_f \sma_{\tilde{N'}} F(N'_+, \bbk^\mfib) \sma F(G'_{r,+},
\bS^\mfib)\right) /G'_r
\end{multline}
induced by the equivalence 
\begin{equation}
\begin{tikzcd}
\bS \cong G'_{r,+} / G'_r \ar[r,"\htp"] & F(G'_{r,+}, \bS^\mfib)/G'_r.
\end{tikzcd}
\end{equation}

We can now construct the following comparison as a consequence of
Lemma~\ref{lem:orbits-and-Q} and
Corollary~\ref{cor:right-fixed-points-project}.

\begin{lem} \label{lem:swapping-order-fixed-point-quotient}
The interchange of fixed points and quotients defines lax monoidal
transformations 
\begin{equation}
\begin{aligned}
  &\left(\cQ_{\cU(G')} \ssF_{(EY)_f/N' \sma G'/N'}((EY)_f \wedge
  G', \bbk)^{-V_f} \right)^{G'}/G'_r \\
&\Rightarrow \left(\cQ_{\cU(G')}(\ssF_{(EY)_f/N' \wedge G'/N'}((EY)_f \wedge
  G', \bbk)^{-V_f}/G'_r)\right)^{G'} \\
    &\Leftarrow \left( \cQ_{\cU(G')} \ssF_{(EY)_f/N'}((EY)_f, \bbk)^{-V_f}\right)^{G'}
     \end{aligned}
\end{equation}
of functors from $\Tw \Sp^{- \Vect}_{\eq,\Sub}$ to $\Sp$. \qed
\end{lem}

\subsection{The norm map for virtual cochains}
\label{sec:norm-map}

We now complete the zig-zag of comparisons by using the norm map. Fix
a finite group $G$.  For a spectrum $\bbk$, recall that the norm map
is specified as the composite 
\begin{equation}
EG_+ \sma \bbk \to \bbk \to F(EG_+, \bbk),
\end{equation}
where the maps are induced by the projection $EG \to \ast$.  Passing
to derived $G$-fixed points (i.e., via fibrant replacement), we have
the maps  
\begin{equation}
(\cQ_{\cU} (EG_+ \sma \bbk))^G \to (\cQ_{\cU}\bbk)^G \to
  (\cQ_{\cU}F(EG_+, \bbk))^G. 
\end{equation}
The norm map is externally multiplicative. 

\begin{prop}\label{prop:norm-mult}
Let $G_1$ and $G_2$ be finite groups and assume that $\bbk$ is an
associative ring orthogonal spectrum.  Then the following diagram
commutes
\begin{equation}
\begin{tikzcd}
(EG_{1,+} \sma \bbk) \sma (EG_{2,+} \sma \bbk) \ar[r] \ar[d] & \bbk
  \sma \bbk \ar[d] \ar[r] & F(EG_{1,+}, \bbk) \sma F(EG_{2,+}, \bbk)
  \ar[d] \\
(EG_1 \times EG_2)_+ \sma (\bbk \sma \bbk) \ar[r] \ar[d] & \bbk \sma \bbk
  \ar[r] \ar[d] & F((EG_1 \times EG_2)_+, \bbk \sma \bbk) \ar[d] \\
(EG_1 \times EG_2)_+ \sma \bbk \ar[r] & \bbk \ar[r] & F((EG_1 \times
  EG_2)_+, \bbk).
\end{tikzcd}
\end{equation} \qed
\end{prop}

Because the variance of the functors on the different sides of the
norm map is different, to express the functoriality of the norm map
once again requires the twisted arrow category.

\begin{lem}
The assignments 
\begin{equation}
f \mapsto EG_+ \sma Y_f \sma \bbk \qquad\textrm{and}\qquad f \mapsto
F(EG'_+, Y_f \sma \bbk) 
\end{equation}
specify functors $\Tw \Sp^{-\Vect}_{\eq} \to \Sp$.\qed
\end{lem}

In this setting, we have a generalized norm map
\begin{equation}
EG_+ \sma Y_f \to Y_f \to F(EG'_+, Y_f)
\end{equation}
which specifies a natural transformation of functors from $\Tw
\Sp^{-\Vect}_{\eq} \to \Sp$.  The verification that this is a natural
transformation amounts to a check that associated to a factorization
$f_0 = h \circ f_1 \circ g$, the diagram 
\begin{equation}\label{eq:simplified-norm-natural}
\begin{tikzcd}
EG_{0,+} \sma Y_{f_0} \ar[d] \ar[r] & \ar[d] F(EG'_{0,+}, Y_{f_0}) \\
EG_{1,+} \sma Y_{f_0} \ar[d] \ar[r] & \ar[d] F(EG'_{0,+}, Y_{f_0}) \\
EG_{1,+} \sma Y_{f_1} \ar[r] & F(EG'_{1,+}, Y_{f_1})
\end{tikzcd}
\end{equation}
commutes.

We now integrate the norm map into a description of the virtual
cochains.  Recall that given a map $f \colon (G, Y, V) \to (G', Y',
V')$, we have a natural map  
\begin{equation}
(EY)_f \equiv  (EG_+ \wedge Y)/G^\perp_f  \to EG'_+ \sma
(Y/G^\perp_f) \equiv E(Y_f).  
\end{equation}
Since by hypothesis $G^{\perp_f}$ acts freely on $Y$, this map is a
weak equivalence.  Composing this with the norm map 
\begin{equation}
   EG'_+ \sma \bbk \to \bbk \to F( EG'_+, \bbk)
\end{equation}
yields a composite
\begin{equation}\label{eq:norm-pointwise}
\begin{aligned}
\ssF_{(EY)_f/N'}((EY)_f, \bbk)^{-V_f} & \equiv
F(S^{V_f}, ((EG_+ \sma Y)/G^{\perp_f}) \sma_{\tilde{N'}} F(N'_+, \bbk^\mfib))
\\
& \to F(S^{V_f}, (EG'_+ \sma Y_f) \sma_{\tilde{N'}} F(N'_+, \bbk^\mfib)) \\ & \to F(S^{V_f}, Y_f
\sma_{\tilde{N'}} F(N'_+, \bbk^\mfib)) \\  
& \to F(S^{V_f}, F(EG'_+, Y_f \sma_{\tilde{N'}} F(N'_+, \bbk^\mfib))) \\ 
& \cong F( EG'_+,
\ssF_{Y_f/N'}(Y_f, \bbk)^{-V_f}).
\end{aligned}
\end{equation}

This composition induces a weak equivalence on fixed points when the
norm does.

\begin{lem}\label{lem:norm-pointwise-equivalence}
The natural map of Equation~\eqref{eq:norm-pointwise} induces a weak
equivalence of orthogonal spectra on passage to derived $G'$-fixed
points when $\bbk$ is a Morava $K$-theory spectrum.
\end{lem}

\begin{proof}
It suffices to show that the composite
\begin{equation}
\begin{tikzpicture}
\node[scale=.9]{
\begin{tikzcd}
(EG'_+ \sma Y_f) \sma_{\tilde{N'}} F(N'_+, \bbk^\mfib) \to Y_f
\sma_{\tilde{N'}} F(N'_+, \bbk^\mfib) \to F(EG'_+, Y_f \sma_{\tilde{N'}} F(N'_+,
\bbk^\mfib))
\end{tikzcd}
};
\end{tikzpicture}
\end{equation}
induces a weak equivalence on derived $G'$ fixed points, which follows
from consideration of the commutative diagram 
\begin{equation}
\begin{tikzpicture}
\node[scale=.85]{
\begin{tikzcd}
(EG'_+ \sma Y_f) \sma_{\tilde{N'}} F(N'_+, \bbk^\mfib) \ar[r] & Y_f
\sma_{\tilde{N'}} F(N'_+, \bbk^\mfib) \ar[r] & F(EG'_+, Y_f \sma_{\tilde{N'}} F(N'_+,
\bbk^\mfib)) \\
(EG'_+ \wedge Y_f) \sma_{\tilde{N'}} (N'_+ \sma \bbk^\mfib) \ar[u,"\htp"] \ar[r] & Y_f
\sma_{\tilde{N'}} (N'_+ \sma \bbk^\mfib) \ar[r] \ar[u,"\htp"] & F(EG'_+, Y_f
\sma_{\tilde{N'}} (N'_+ \sma \bbk^\mfib) \ar[u,"\htp"] \\
(EG'_+ \wedge Y_f) \sma \bbk^\mfib) \ar[r] \ar[u,"\htp"] & Y_f
\sma \bbk^\mfib \ar[r] \ar[u,"\htp"] & F(EG'_+, Y_f \sma
\bbk^\mfib). \ar[u,"\htp"] 
\end{tikzcd}
};
\end{tikzpicture}
\end{equation}
The bottom horizontal maps are equivalences on derived $G'$-fixed
points by hypothesis and the vertical maps are equivalences by Atiyah
duality.  Therefore, we can conclude that the top horizontal maps are
equivalences.
\end{proof}

The key result of this subsection is that the pointwise comparison of
Lemma~\ref{lem:norm-pointwise-equivalence} is compatible with the
functoriality in the twisted arrow category.  To set this up, we need
to first establish the functoriality of each side of the comparison.
The following proposition records this for the two terminal terms in
the comparison composite; the argument is again essentially the same
as for Lemma~\ref{lem:functor-fiberwise-maps-Borel}.

\begin{prop}
The assignments
\begin{equation}
f \mapsto \left( \cQ_{\cU(G')} \ssF_{(EY)_f/N'}((EY)_f,
\bbk)^{-V_f}\right)^{G'}
\end{equation}
and
\begin{equation}
f \mapsto \left(\cQ_{\cU(G')} F(EG'_+, \ssF_{Y_f/N'}(Y_f, \bbk)^{-V_f})\right)^{G'}
\end{equation}
are functors $\Tw \Sp^{-\Vect}_{\eq, \Sub} \to \Sp$. \qed
\end{prop}

Now we have the following comparison.

\begin{lem}\label{lem:twisted-norm}
The comparison of Equation~\eqref{eq:norm-pointwise} induces a lax
monoidal natural transformation 
\begin{align}
\left( \cQ_{\cU(G')} \ssF_{(EY)_f/N'}((EY)_f,
\bbk)^{-V_f}\right)^{G'}    &\Rightarrow \left( \cQ_{\cU(G')}
F(EG'_+, \ssF_{Y_f/N'}(Y_f, \bbk)^{-V_f})\right)^{G'}  
\end{align}
of functors from $\Tw \Sp^{- \Vect}_{\eq, \Sub} \to \Sp$.
\end{lem}

\begin{proof}
To see that this is a natural transformation, given a factorization
$f_0 = h \circ f_1 \circ g$, we need to show that the diagram 
\begin{equation}
\begin{tikzpicture}
\node[scale=.75]{
\begin{tikzcd}
\cQ_{\cU(G_0')}F(S^{V_{f_0}}, (EY_{f_0} \sma_{\tilde{N_0'}} F(N_{0',+}, \bbk^\mfib))^{G_0'}
\ar[r] \ar[d] & \ar[d] \cQ_{\cU(G_0')} F(S^{V_{f_0}}, F(EG'_{0,+},
  Y_{f_0} \sma_{\tilde{N_0'}} F(N'_{0,+}, \bbk^\mfib)))^{G_0'} \\ 
\cQ_{\cU(G_1')} F(S^{V_{f_1}}, (EY_{f_1} \sma_{\tilde{N_1'}}
F(N'_{1,+}, \bbk^\mfib))^{G_1'} \ar[r] & \cQ_{\cU(G_1')} F(S^{V_{f_1}}, F(EG'_{1,+}, Y_{f_1}
\sma_{\tilde{N'_1}} F(N'_{1,+}, \bbk^\mfib)))^{G_1'}
\end{tikzcd}
};
\end{tikzpicture}
\end{equation}
commutes.  This follows from the discussion surrounding
Equation~\eqref{eq:simplified-norm-natural}; the essential point is that
although the term $F(S^{V_f}, Y_f \sma_{\tilde{N'}} F(N'_+, \bbk^\mfib))$ in
the pointwise comparison is not a functor in $\Tw \Sp^{-\Vect}_{\eq,
  \Sub}$, the generalized norm map nonetheless specifies a natural
transformation because the composite  
\begin{equation}
\begin{tikzcd}
\cQ_{\cU(G')} F(S^{V_f}, (EG'_+ \sma Y_f) \sma_{\tilde{N'}} F(N'_+,
\bbk^\mfib))^{G'} \ar[d] \\
\cQ_{\cU(G')} F(S^{V_f}, F(EG'_+, Y_f \sma_{\tilde{N'}} F(N'_+, \bbk^\mfib)))^{G'} 
\end{tikzcd}
\end{equation}
is compatible with composition.  The fact that the transformation is
monoidal follows from the fact that all of the constituent functors
(the smash product of function spectra, passage to orbits, fibrant
replacement functor) and the norm map (by
Proposition~\ref{prop:norm-mult}) are externally monoidal.
\end{proof}

Lemma~\ref{lem:norm-pointwise-equivalence} now implies that when the
norm map is an equivalence, this comparison transformation is an
equivalence.  Assembling everything, we have the following result
which implements the last piece of the comparison with the virutal
cochains.

\begin{thm} \label{thm:norm-map-equivalence-zig-zag-lax-monoidal}
When $\bbk$ is a Morava $K$-theory spectrum, there are lax monoidal
equivalences 
\begin{align}
\left( \cQ_{\cU(G')} \ssF_{(EY)_f/N'}((EY)_f,
\bbk)^{-V_f}\right)^{G'}    &\Rightarrow \left( \cQ_{\cU(G')}   
F(EG'_+, \ssF_{Y_f/N'}(Y_f, \bbk)^{-V_f})\right)^{G'} \\
&\Leftarrow F(EG', \ssF_{Y_f/N'}(Y_f, \bbk)^{-V_f})^{G'} \\ 
&\equiv C^*(BG', \ssF_{Y_f/N'}(Y_f, \bbk)^{-V_f}) 
\end{align}
of functors on $\Tw \Sp^{- \Vect}_{\eq, \Sub}$. 
\end{thm}

\begin{proof}
The top arrow was constructed above in Lemma~\ref{lem:twisted-norm}
and is an equivalence by Lemma~\ref{lem:norm-pointwise-equivalence}.
The middle (left-pointing) arrow is induced by the natural
transformation $\id \to \cQ_{\cU(G')}$; this is lax monoidal by
Corollary~\ref{cor:unit-fibrant-funct} and is a weak equivalence
because the homotopy fixed points only depend on the Borel homotopy
type.
\end{proof}

\subsection{The comparison of the virtual cochains}\label{sec:complete-comparison}

Finally, we can assemble the comparisons of the preceding subsections
into zig-zags of lax monoidal natural transformations of functors from
$\Tw \Sp^{-\Vect}_{\eq, \Sub}$ to the category of orthogonal spectra.
Note that although we have ignored the $\Pi$-action throughout this
section, all of the transformations are evidently $\Pi$-equivariant
because of the rigidity of the $\Pi$-action.  By
Lemma~\ref{lem:monoidal_comparison_classifying_fibrewise-2}, there is
a zig-zag of lax monoidal natural transformations through equivalences
connecting the functors
\begin{equation}
BY^{-V_f} \sma \bbk \htp \ssF_{BY}( (EY)_f, \bbk)^{-V_f}/ G'.
\end{equation}
(Here we are tacitly pulling back along the projection 
$\Tw \Sp^{-\Vect}_{\eq, \Sub}$ to $\Tw \Sp^{-\Vect}_{\eq}$.) 
Next, by Proposition~\ref{prop:relative-comparison} there is a lax
monoidal equivalence 
\begin{equation}
\ssF_{(EY)_f/N'}((EY)_f, \bbk)^{-V_f}/G' \Rightarrow \ssF_{BY}(
(EY)_f, \bbk)^{-V_f}/ G'.
\end{equation}
Lemma~\ref{lem:genu-fibr-replacement-equivalence} establishes a lax
monoidal equivalence 
\begin{equation}
\ssF_{(EY)_f/N'}((EY)_f, \bbk)^{-V_f}/G' \to
\cQ_{\cU(G')}\ssF_{(EY)_f/N'}((EY)_f, \bbk)^{-V_f}/G',
\end{equation}
and Corollary~\ref{cor:right-fixed-points-project} produces a lax
monoidal equivalence
\begin{align}
  ( \cQ_{\cU(G')}\ssF_{(EY)_f/N' \wedge G'/N'}((EY)_f \wedge G',
  \bbk)^{-V_f})^{G'})/G'_r \\  
\Rightarrow  \cQ_{\cU(G')}\ssF_{(EY)_f/N'}((EY)_f, \bbk)^{-V_f})/G'
\end{align}
Next, Lemma~\ref{lem:swapping-order-fixed-point-quotient} implies that
we have lax monoidal equivalences
\begin{equation}
\begin{aligned}
  &\left(\cQ_{\cU(G')} \ssF_{(EY)_f/N' \wedge G'/N'}((EY)_f \wedge
  G', \bbk)^{-V_f} \right)^{G'}/G'_r \\
&\Rightarrow \left(\cQ_{\cU(G')}(\ssF_{(EY)_f/N' \wedge G'/N'}((EY)_f \wedge
  G', \bbk)^{-V_f}/G'_r)\right)^{G'} \\
    &\Leftarrow \left( \cQ_{\cU(G')} \ssF_{(EY)_f/N'}((EY)_f, \bbk)^{-V_f}\right)^{G'}
     \end{aligned}
\end{equation}
Finally, Theorem~\ref{thm:norm-map-equivalence-zig-zag-lax-monoidal}
establishes a zig-zag of equivalences
\begin{align}
\left( \cQ_{\cU(G')} \ssF_{(EY)_f/N'}((EY)_f,
\bbk)^{-V_f}\right)^{G'}    &\Rightarrow \left( \cQ_{\cU(G')}   
F(EG'_+, \ssF_{Y_f/N'}(Y_f, \bbk)^{-V_f})\right)^{G'} \\
&\Leftarrow F(EG', \ssF_{Y_f/N'}(Y_f, \bbk)^{-V_f})^{G'} \\ 
&\equiv C^*(BG', \ssF_{Y_f/N'}(Y_f, \bbk)^{-V_f}).
\end{align}

\section{Morse theory and homotopy types}
\label{sec:morse-theory-and-homotopy-types}

Our goal in this section is to reprove from our perspective the result
of Cohen, Jones, and Segal~\cite[Section 5]{CohenJonesSegal1995},
which asserts that the homotopy type associated to natural (relative)
framings of the moduli spaces of gradient flow trajectories of a
Morse function on a closed smooth manifold $M$ agrees with the (stable)
homotopy type of $M$.

\subsection{Morse-theoretic setup}
\label{sec:morse-theor-setup}

Let $(f,g)$ be a a Morse-Smale pair consisting of a Morse function $f$ and a metric $g$ on a closed smooth manifold $M$.  For each pair $(x,y)$ of critical points, let  $\barT(x,y)$ denote the compactified moduli space of gradient flow lines of $f$, converging at $-\infty$ to $x$ and at $+\infty$ to $y$.
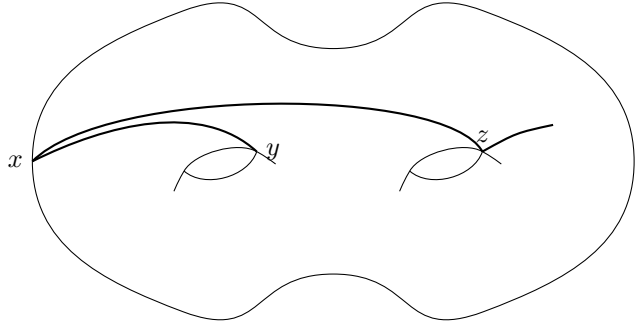
\begin{figure}[h]
  \centering
  \begin{tikzpicture}
    \begin{scope}
    \draw plot [smooth cycle, tension=1] coordinates {(-4,0) (-2,-2) (0,-1.5) (2, -2) (4,0) (2,2) (0,1.5) (-2,2) };
    \begin{scope}[shift = {(1,0)}, rotate around = {15:(-2.5,0)}]
      \coordinate[label= right:$y$] (y) at (-2,0);
     \draw (-3.2,-.225) .. controls (-3.1,-.1) and (-3.1, -.1) .. (-3,0) .. controls (-2.8, .2) and (-2.2, .2) .. (-2,0) .. controls (-1.9,-.1) and (-1.9, -.1) .. (-1.8,-.225);
     \draw  (-3,0) .. controls (-2.8, -.3) and (-2.2,- .3) .. (-2,0) ;
   \end{scope}
   \coordinate[label= left:$x$] (x) at (-4,0);
     \draw [thick] (-4,0) .. controls (-3,.5) and (-1.75, .8) .. (y);
     \begin{scope}[shift = {(4,0)}, rotate around = {15:(-2.5,0)}]
        \coordinate[label= above :$z$] (z) at (-2,0);
      \draw (-3.2,-.225) .. controls (-3.1,-.1) and (-3.1, -.1) .. (-3,0) .. controls (-2.8, .2) and (-2.2, .2) .. (-2,0) .. controls (-1.9,-.1) and (-1.9, -.1) .. (-1.8,-.225);
      \draw  (-3,0) .. controls (-2.8, -.3) and (-2.2,- .3) .. (-2,0) ;
      \draw [thick] (z) .. controls (-1.5,.125) .. (-1,.1);
    \end{scope}
     \draw [thick] (-4,0) .. controls (-3,1) and (1.5, 1) .. (z);
    \end{scope}
  \end{tikzpicture}
  \caption{Elements of the moduli spaces of flow lines $\barT(x,y)$ and $\barT(x,M)$ for the function given by projection to the horizontal line. The element of $ \barT(x,M)$ lies in the image of the stratum $\barT(x,z) \times \barT(z,M)$.}
  \label{fig:moduli-spaces-flow-lines}
\end{figure}
There are natural inclusions of codimension $1$ boundary strata
\begin{align} \label{eq:breaking-gradient-trajectory}
  \barT(x,y)  \times \barT(y,z) & \to     \barT(x,z),
 \end{align}
 which yields a flow category $\barT(f)$ with objects the set of
 critical points. The Morse-Smale assumption implies that the morphism
 spaces $\barT(x,y)$ admit global Kuranishi charts with trivial group
 $G$ and trivial obstruction space $V$.   

For expository purposes, it is convenience to arrange for the moduli
spaces to be smooth. We thus work in the setup introduced by Burghelea
and Haller \cite{Burghelea2001} and revisited by Wehrheim in
\cite{Wehrheim2012}: assume that the metric $g$ is flat near each
critical point and that the function $f$ takes the standard form $\sum
x_i^2 - \sum y_j^2$ in local flat coordinates. As proved in
\cite[Theorem 1]{Burghelea2001}  the moduli spaces $\barT(x,y)$
acquire natural smooth structures as manifolds with corners. Under
these assumptions, we will prove the following result using the
framework of this paper. 

\begin{prop} \label{prop:computation-Morse}
The category $\barT(f)$ lifts to a complex oriented Kuranishi flow
category. The corresponding homotopy type $CM_*(f; \bbk)$ is
well-defined over any complex oriented spectrum $\bbk$, and is
equivalent to the spectrum of ordinary chains $C_*(M; \bbk)$. 
\end{prop}

The proof occupies the remainder of this Appendix.

\begin{rem}
Although it is not relevant for this paper, the following may be of
interest to the reader: in the next section we shall in fact reprove
that $\barT(f)$ is a framed flow category in the sense
of~\cite{CohenJonesSegal1995}, which implies that the condition that
$\bbk$ be complex oriented can be dropped from the above statement. 
\end{rem}

\subsection{Stable framings of moduli spaces}
\label{sec:stable-fram-moduli}

Let $\barT(x,M)$ and $\barT(M,y)$ denote the moduli spaces of negative
(i.e., with domain $(-\infty,0])$) and positive (with domain
  $[0,+\infty)$) half-gradient flow lines; in the first case, this is
    a compactification of the ascending manifold of $x$, and in the
    second case, this is a compactification of the descending
    manifold. These are smooth contractible manifolds with corners
    equipped with evaluation maps
\begin{equation}
\barT(x,M) \to M \leftarrow \barT(M,y).
\end{equation}
The tangent space of $\barT(x,M) $ at the constant flow lines is
canonically isomorphic to the positive-definite subspace of the
Hessian matrix of $f$ and $x$, which we denote $V^+_x$.  We shall fix: 
\begin{equation} \label{eq:fix-trivialisation}
  \parbox{31em}{An identification of $ T\barT(x,M)$ with $V_x^+$ which is the identity at the constant flow line.}
\end{equation}

For each critical point $x$, consider the embedding
\begin{equation} \label{eq:breaking-gradient-negative-flow-line}
  \barT(x,y) \times \barT(y,M) \to \barT(x,M)
\end{equation}
of a codimension $1$ boundary stratum. Given a sufficiently small constant $\epsilon$, we define a function to $\bR$ on a neighbourhood of this stratum which assigns to every flow line $\gamma$ the quantity
\begin{equation}
e^{-T_{y}(\gamma)} \in [0,\infty)    
\end{equation}
where $T_y(\gamma)$ is the length of unique interval in the domain of $\gamma$, whose endpoints map to the level sets $f(y) \pm \epsilon$. The key point of the flatness assumption of the metric near the critical points is that it ensures that the moduli space can be equipped with a smooth structure so that this map is the projection to the collar direction near the stratum $\barT(x,y) \times \barT(y,M)$. In this way, we obtain a short exact sequence
\begin{equation} 
 T \barT(x,y) \oplus T \barT(y,M) \to T \barT(x,M) \to \ell_y,
\end{equation}
where $\ell_y$ is a real line canonically isomorphic to $\bR$. A choice of Riemannian metric on $ \barT(x,M)$ then determines an isomorphism
\begin{equation}
     T \barT(x,M) \cong T \barT (x,y) \oplus \ell_y \oplus T \barT(y,M).
  \end{equation}

Using the choice fixed in Equation \eqref{eq:fix-trivialisation}, we
obtain a stable trivialisation
\begin{equation} \label{eq:stable_trivialisation-moduli-flow-lines}
  T \barT (x,y) \oplus \ell_y \oplus V^+_y \cong V_x^+
\end{equation}
of the tangent bundle of $ T \barT (x,y) $ (relative the vector spaces
$V_x^+$ and $V_y^+ \oplus \ell_y$).

We now consider the restriction of this bundle to the boundary strata
of the moduli spaces of gradient flow lines: given a triple of
critical points $x$, $y$, and $z$, we have a commutative diagram 
\begin{equation}
  \begin{tikzcd}
    \barT(x,y)  \times \barT(y,z) \times \barT(z,M) \ar[r] \ar[d] & \barT(x,y)  \times \barT(y,M) \ar[d] \\
    \barT(x,z) \times \barT(z,M) \ar[r] & \barT(x,M),
  \end{tikzcd}
\end{equation}
where each map is an inclusion of a codimension $1$ stratum. Choosing a Riemannian metric on these manifolds so that the inclusions are isometric embeddings yields a commutative diagram
\begin{equation}
  \begin{tikzcd}
   T \barT(x,y)  \oplus \ell_y \oplus  T \barT(y,z) \oplus \ell_z \oplus  T \barT(z,M) \ar[r] \ar[d] & T \barT(x,y)  \oplus \ell_y \oplus T \barT(y,M) \ar[d] \\
   T \barT(x,z) \oplus \ell_z \oplus T \barT(z,M) \ar[r] & T \barT(x,M),
  \end{tikzcd}
\end{equation}
of isomorphisms of tangent spaces. 

Using our given stable trivializations, we conclude:

\begin{cor}
The stable trivializations from Equation
\eqref{eq:stable_trivialisation-moduli-flow-lines} fit in the
following commutative diagram 
\begin{equation}
\begin{tikzcd}
T \barT (x,y) \oplus \ell_y \oplus T \barT (y,z) \oplus \ell_z \oplus V_{z}^{+}  \ar[r]  \ar[d] &  \ar[d]  T \barT (x,y) \oplus \ell_y \oplus V^+_y  \ar[d] \\
  T \barT (x,z) \oplus \ell_z \oplus V_{z}^{+}  \ar[r]  &  V_x^+.
\end{tikzcd}
\end{equation} \qed
\end{cor}

Comparing with Definition~\ref{def:complex-oriented-flow-category}, we
conclude that we have constructed a complex oriented Kuranishi
presentation of $\barT(f)$.

\begin{prop}
The constructions above specify a complex oriented Kuranishi
presentation of $\barT(f)$, where $\Pi$ is the trivial group.  The data
consists entirely of global charts with trivial isotropy, and where
the complex vector bundle (and the obstruction bundle) appearing in
the definitions are both trivial.
\end{prop}

It should be clear at this stage that our construction amounts to
producing compatible (relative) framings of the moduli spaces in the
sense of \cite{CohenJonesSegal1995}, which is the data required to
build a stable homotopy type over the sphere spectrum $\bS$.

\subsection{Computation of the Morse homotopy type}
\label{sec:comp-morse-homot}

We begin by extending the category $\barT(f)$ to a slightly larger
category $\barT(f,M)$ and producing a corresponding Kuranishi flow
category.

\begin{defin}
The category $\barT(f,M)$ has:
\begin{itemize}
\item Objects consisting of the critical
points of $f$, together with a new terminal object which we denote by $M$.
\item Morphisms specified by setting the endomorphisms of $M$ to be a
  point and morphisms from $x$ to $M$ are given by the moduli spaces
  $\barT(x,M)$, with the new compositions given by Equation
\eqref{eq:breaking-gradient-negative-flow-line}.
\end{itemize}
\end{defin}

In Equation~\eqref{eq:fix-trivialisation}, we have already chosen
trivialisations of the moduli spaces of descending gradient flow
lines, so that the category $\barT(f,M)$ also lifts to a Kuranishi
flow category with complex orientation and trivial group $\Pi$, and
with trivial isotropy and obstruction bundle; to use precisely the
same context, we add a copy of a real line $V^+_M \cong \bR$ to the
two sides, to get
\begin{equation}
V^+_M \oplus T\barT(x,M) \cong \bR \oplus V_x^+.
\end{equation}

We shall consider a homotopy type associated to this category which
corresponds to a functor whose value at $M$ is $M_+ \wedge \bbk$.  
To formulate this, consider the following consequence of the results
of this paper: 

\begin{lem} \label{lem:bimodule_equivalence_Morse_theory}
There is a spectral bimodule representing an equivalence between
$C^*_{\rel \partial}(\barT(f,M), \Omega \bbk)$ and the category
$\barT(f,M)^{-d} \wedge \bbk $, whose objects are those of
$\barT(f,M)$, and whose morphisms are 
\begin{align}\label{eq:morsecat}
\barT(x,y)^{-d} \wedge \bbk & \equiv  S^{V_y^+} \wedge S^{-V_x^+} \wedge \barT(x,y)_+ \wedge \bbk\\
\barT(x,M)^{-d} \wedge \bbk & \equiv S^{V^+_M} \wedge S^{-V_x^+} \wedge \barT(x,M)_+ \wedge \bbk.
\end{align}
\end{lem}
\begin{proof}
The construction follows the procedure of Sections
\ref{sec:local-glob-constr} and \ref{sec:compl-orient-flow-1}, and
results in a comparison with the category of virtual cochains
constructed in Section \ref{sec:impl-atlas-with}, whose morphisms are
given by Definition \ref{def:virtual_cochains_chart}. In this case,
these morphisms are: 
\begin{align}\label{eq:rawmorsecat}
&  S^{V_y^+} \wedge S^{-V_x^+} \wedge (\barT(x,y)_+)^{\mfib} \wedge \bbk \\
& S^{V^+_M} \wedge S^{-V_x^+} \wedge (\barT(x,M)_+)^{\mfib} \wedge \bbk,
\end{align}
   where we use the fact that the reduced degree vanishes because we
   are working with a trivial group $\Pi$.  Now the monoidal natural
   transformation $\id \to (-)^{\mfib}$ yields a $DK$-equivalence
   between the spectral category in Equation~\eqref{eq:morsecat} and
   the spectral category in Equation~\eqref{eq:rawmorsecat}.
\end{proof}

We denote by $\barT(f)^{-d}$ the full spectral subcategory of
$\barT(f,M)^{-d}$ with objects the critical points of $f$.  The
virtual fundamental chain of the Morse flow category $\barT(f)$ is
obtained from the Spanier-Whitehead duality equivalence of the
relative cochains with $\barT(f)^{-d} $, together with the functor 
\begin{equation}
\delta^{f} \colon \barT(f)^{-d} \wedge \bbk \to \bbk-\mod
\end{equation}
which assigns $S^{V_x^+} \wedge \bbk$ to each critical point $x$, and is given at the level of morphisms by the projection from $\barT(x,y)$ to a point.  Specifically, the functor is given by the map
\begin{equation}
\begin{tikzcd}
S^{V_y^+} \wedge S^{-V_x^+} \wedge \barT(x,y)_+ \wedge \bbk \ar[d] \\
S^{V_y^+} \wedge S^{-V_x^+} \wedge \bbk \ar[d] \\
F(S^{V_x^+}, S^{V_y^+} \sma \bbk) \ar[d] \\
F_{\bbk}(S^{V_x^+} \sma \bbk, S^{V_y^+} \sma \bbk).
\end{tikzcd}
\end{equation}
It is straightforward to check that this composite is compatible with the composition.

We extend this to a functor
\begin{equation}
\delta^{f,M} \colon \barT(f,M)^{-d} \wedge \bbk \to  \bbk-\mod,
\end{equation}
which assigns to the terminal object $M$ the chains $S^{V_M^+} \wedge
M_+ \wedge \bbk$. The action of morphisms is given by the map 
\begin{align}
 S^{V_x^+} \wedge S^{V_M^+} \wedge S^{-V_x^+} \wedge \barT(x,M)_+ \wedge \bbk \to S^{V^+_M} \wedge M_+ \wedge \bbk,
\end{align}
obtained by pairing positive and negative spheres against each other, together with the evaluation map from $\barT(x,M)$ to $M$.

\begin{proof}[Proof of Proposition \ref{prop:computation-Morse}]
The restriction $\delta^M$ of $\delta^{f,M}$ to the full spectral
subcategory of $\barT(f,M)^{-d}$ consisting of the terminal object
induces a comparison
\begin{equation}\label{eq:Minclude}
|\delta^M| \to |\delta^{f,M}|,
\end{equation}
since we can regard $\delta^{f,M}$ as specifying categorical
continuation data in the sense of
Definition~\ref{defn:categorical-continuation}.

Next, observe that we have an equivalence
\begin{equation}
|\delta^M| \htp \bS^{V_M^+} \sma M_+ \sma \bbk,
\end{equation}
as $|\delta^M|$ is specified over a poset with one object plus
$\infty$ (recall part (2) of
Example~\ref{ex:compute_filtered_homotopy_type}).  Since by Theorem 
\ref{thm:cofiber} the cofiber of the map in
Equation~\eqref{eq:Minclude} can be identified with the Kan suspension
of $|\delta^{f}|$, we have a cofiber sequence
\begin{equation}
\bS^{V_M^+} \sma M_+ \sma \bbk \to |\delta^{f,M}| \to \Sigma
|\delta^{f}|.
\end{equation}
As a consequence, it will suffice to show that $|\delta^{f,M}| \htp
\ast$; then there is an induced equivalence
\begin{equation}
|\delta^{f}| \to \bS^{V_M^+} \sma M_+ \sma \bbk.
\end{equation}

We will construct a different filtration on $|\delta^{f,M}|$ to show
that it is contractible.  For each pair $(a,b)$ of real numbers,
consider the spectral category  
\begin{equation}
\barT(f,M)_{[a,b]}
\end{equation}
consisting of all critical points $x$ such that $a \leq f(x) \leq b$,
together with the terminal object $M$, with morphisms between critical
points given as before, and morphisms from $x$ to $M$ given by the
quotient 
\begin{equation}
\barT(f,M)_{[a,b]}(x,M) \equiv  \barT(x,M)/ f^{-1}((-\infty,a]),
\end{equation}
where $f$ is defined on the space of flow lines starting at $x$ via
the evaluation map to $M$.  The composition 
\begin{equation}
\barT(f,M)_{[a,b]}(x,y) \sma \barT(f,M)_{[a,b]}(y,M) \to 
\barT(f,M)_{[a,b]}(x,M)  
\end{equation}
is defined as the composition in $\barT(f,M)$ composed with the
projection onto the quotient.

We now consider the functor
\begin{equation}
\delta^{f,M}_{[a,b]} \colon \barT(f,M)_{[a,b]} \wedge \bbk \to  \bbk-\mod,
\end{equation}
defined on objects as 
\begin{equation}
\delta^{f,M}_{[a,b]}(x) = 
\begin{cases}
\delta^{f,M}(x) & x \in [a,b] \\
\left(M_{\leq b}/M_{\leq a}\right)_+ \sma \bbk \qquad & x = M.
\end{cases}
\end{equation}
Here to specify the action on morphisms we are using the fact that the
evaluation map from the descending manifolds of critical points of
value $\leq b$ is well-defined in this space.

For $a \leq b$ there exist coherent zig-zags representing natural
transformations
\begin{equation}
\delta^{f,M}_{[-\infty,b]}(x) \to \delta^{f,M}_{[-\infty,a]}(x)
\end{equation}
and so we obtain a filtration on $|\delta^{f,M}|$ with associated
graded consisting of terms equivalent to $|\delta^{f,M}_{[a,b]}|$.
(The construction of the filtration here in terms of zig-zags is
precisely analogous to the construction of the filtration by action
discussed in Section~\ref{sec:from-flow-categories}.) 

We will show that for a sufficiently fine filtration the associated
graded of $|\delta^{f,M}|$ with respect to this filtration is
contractible and therefore by induction so is $|\delta^{f,M}|$.  Since
we can assume without loss of generality that all critical points have
distinct values of $f$, we can reduce to considering the associated
graded homotopy type $|\delta^{f,M}_{[a,b]}|$ for a window $[a,b]$
containing a single critical point $x$.

This homotopy type is computed over a poset with
two objects $M$ and $x$ plus an additional terminal object, where the
virtual fundamental chain maps $M$ to
\begin{equation}
M_{[a,b]} \equiv f^{-1}((-\infty,a])/ f^{-1}((-\infty,b])
\end{equation}
and the critical point $x$ to $S^{V_x}$.  At the level of morphisms,
Morse theory implies that the map
\begin{equation}
\barT(f,M)_{[a,b]} \to  M_{[a,b]}
\end{equation}
is an equivalence, from which we deduce that the morphism from $x$ to
$M$ in $C^*_{\rel \partial}(\barT(f,M), \Omega \bbk)$ induces an
equivalence 
\begin{equation}
C^*_{\rel \partial}(\barT(f,M), \Omega \bbk)(x,M) \sma
\delta^{f,M}_{[a,b]}(x) \to \delta^{f,M}_{[a,b]}(M).
\end{equation}
As a consequence, the computation of part (iii) of
Example~\ref{ex:compute_filtered_homotopy_type} shows that 
$|\delta^{f,M}_{[a,b]}|$ is contractible.  
\end{proof}

\bibliographystyle{halpha}
\bibliography{large-bib}

\end{document}